# Lectures on Linear Algebra over Division Ring

Aleks Kleyn

Aleks_Kleyn@MailAPS.org
http://AleksKleyn.dyndns-home.com:4080/
http://sites.google.com/site/AleksKleyn/
http://arxiv.org/a/kleyn_a_1
http://AleksKleyn.blogspot.com/

Abstract. In this book I treat linear algebra over division ring. A system of linear equations over a division ring has properties similar to properties of a system of linear equations over a field. However, noncommutativity of a product creates a new picture.

Matrices allow two products linked by transpose. Biring is algebra which defines on the set two correlated structures of the ring.

As in the commutative case, solutions of a system of linear equations build up right or left vector space depending on type of system. We study vector spaces together with the system of linear equations because their properties have a close relationship. As in a commutative case, the group of automorphisms of a vector space has a single transitive representation on a frame manifold. This gives us an opportunity to introduce passive and active representations.

Studying a vector space over a division ring uncovers new details in the relationship between passive and active transformations, makes this picture clearer.

Considering of twin representations of division ring in Abelian group leads to the concept of $D$-vector space and their linear map. Based on polylinear map I considered definition of tensor product of rings and tensor product of $D$-vector spaces.

# Contents

CHAPTER 1

# Preface

## 1.1. Preface to Version 1

Starting a new journey you do not know in the beginning what to expect down the road. I started studying noncommutative algebra just out of curiosity. In the case of a module over a ring it is impossible to give the definition of basis in the way we do it in the case of a vector space over a field. I wanted to understand how the picture changes when I use division ring instead of a field.

I started from the study of systems of linear equations. I started from system of two equations in two unknowns. Even I solved system without problem, it was impossible to express this solution as the ratio of two determinants.

I understood that this problem might be interesting not only for me. I started to search for the mathematicians, who are interested in the same problem. Professor Retakh introduced me papers [6, 7] dedicated to theory of quasideterminants. It was beginning of my research in theory of vector space.

I dedicated chapter 2 to biring of matrices. There are two reasons why I study this algebra.

Suppose we are given basis in vector space. Then we can describe transformation of vector space using a matrix. Product of matrices corresponds to the product of transformations. In contrast to commutative case it is not each time possible to represent this product as the product of rows of first matrix over columns of second one. The same time representation of matrix as

$$A = \begin{pmatrix} A_1^1 & ... & A_n^1 \\ ... & ... & ... \\ A_1^n & ... & A_n^n \end{pmatrix}$$

depends on convention. Without loss of generality we can represent matrix as

$$A = \begin{pmatrix} A_1^1 & ... & A_1^n \\ ... & ... & ... \\ A_n^1 & ... & A_n^n \end{pmatrix}$$

However in this case the algorithm to find product changes.

Structures like biring are known in algebra. In lattice we define operations $a \vee b$ and $a \wedge b$ which exchange places when order on the set reverses. The symmetry between ${}_*{}^*$-product and ${}^*{}_*$-product is phrased as duality principle. Subsequently I extend the duality principle to representation theory and theory of vector spaces. Without considering the duality principle the book would be four times its size, let alone the fact that endless reiterations would make the text hard to read.

According to each product we can extend the definition of a quasideterminant given in [6, 7] and introduce two different types of a quasideterminant.

The chapter 3 is review of representation theory and is basis for following chapters. We extend to representation theory convention described in remark 2.2.15. Theorem about an existence of twin representations in homogeneous space finishes the chapter.

In the chapter 5 I study a few concepts of linear algebra over division ring $D$. I recall definitions of a vector space and a basis in the beginning.[1.1] Linear algebra over a division ring is more diverse than linear algebra over a field. In contrast to vector space over a field, we can define left and right vector space over an arbitrary division ring $D$. To make the definition of a basis stronger I turn to theory of an arbitrary system of linear equations (section 5.5) for each type of vector space. However, in spite of this diversity statements from linear algebra over a division ring are very similar to statements from linear algebra over a field.

Because I apply statements of this chapter in geometry, I will follow the same notation as we use in geometry. To write coordinates of a vector and elements of a matrix we follow the convention described in section 2.1.

We identify the vector and the set of its coordinates relative selected basis. However, they are different objects. To underscore this difference, I returned to traditional notation of vector as $\overline{a}$, when vector and its coordinates are present in the same equation; at the same time we will use notation $a^b$ for coordinates of vector $\overline{a}$. We use the same root letter for notation of a basis and vectors which form it. To distinguish a vector and a basis, we will use notation $\overline{\overline{e}}$ for the basis and notation $\overline{e}_a$ for vectors which form the basis $\overline{\overline{e}}$. Studying theory of vector space we use convention described in remark 2.2.15 and in section 4.1.

Chapter 6 is dedicated to the theory of linear representation. The studying of a homogenous space of a group of symmetry of a $D_*{}^*$-vector space leads us to the definition of a basis of this $D_*{}^*$-vector space and a basis manifold. We introduce two types of transformation of a basis manifold: active and passive transformations. The difference between them is that the active transformation can be expressed as a transformation of an original space. As it is shown in [3] passive transformation gives ability to define concepts of invariance and of geometric object. Based on this theory I study basis manifold, passive and active transformations in section 6.2. I study geometric object of the $D_*{}^*$-vector space in section 6.3.

We have two opposite points of view about a geometric object. On the one hand we determine coordinates of the geometric object relative to a given basis and introduce the law of transition of coordinates during transformation of the basis. At the same time we study the set of coordinates of the geometric object relative to different bases as a single whole. This gives us an opportunity to study the geometric object without using coordinates.

Chapter 7 is dedicate to linear map of ${}_*{}^*D$-vector spaces over division ring. This statement is basis for development of theory of twin representations of ring $D$. The notation that I introduced for the theory of vector space seems ponderous. However this notation is guided by the existence of twin representations.

Linear algebra over a division ring is more diverse than linear algebra over a field. Noncommutativity compels us to take care about a proper order of factors

[1.1] You can see definitions also in [8].

in an expression. This leads not only to diversity of concepts, but also helps more clear see statements of commutative algebra.

January, 2007

## 1.2. Preface to Version 2

> Any criminal once a while makes mistake.
>
> When detective discovered my cap on the table, he became furious. Next time there was rain, and I left a wet track on the floor. Needless to say, I had to start from the scratch again. And I started again and again. The Tunguska Cosmic Body, an unplanned solar eclipse, ... However, finally detective did not find defects in the design.
>
> - Tic-tac-toe, - detective sought sorrowfully. - You were not tired to code this game? It is time to start something more deep.
>
> Author is unknown. The art of programming.

Science research sometimes reminds programming. Complex and large project gradually becomes surrounded by code. Everything seems to be simple and clear. However, suddenly the bug appears. Fix of this bug demands to alter nearly half of the code. As the wise man observed in any software there is at least one bug, though.

When I started the research in the field of vector spaces over division ring, I saw that I stepped to terra incognita. The concept of quasideterminant extremely lightened my task. The opportunity to solve a system of linear equations made this theory very similar to theory of vector spaces over field. This gave ability to define basis and extent the concept of geometric object.

Tensors of order 2 appeared on the scene as dark cloud. Non commutability seemed to be putting on the structure of tensor insurmountable constrains. However at this time I discovered light in the tunnel.

The considering of $D_*{}^*$-linear maps leads to the concept of twin representations of division ring in Abelian group. Abelian group, in which we define twin representation of division ring $D$, is called $D$-vector space. Since we define multiplication over elements of division ring $D$ from left, as well as right, then homomorphisms of $D$-vector spaces cannot keep this operation. This leads to concept of additive map which is morphism of $D$-vector space.[1.2]

$D$-vector space is the best candidate for building of tensor algebra. However one more important structure was required here.

When we study rings, we study not only left or right module over ring, but bimodule also. Division ring is no exception. However if we add up the concept of the biring, it becomes obvious that it does not matter from what side we multiply.

[1.2] Design of additive mapping which I made proving theorem 12.1.4, is instructive methodically. From the time when I started to learn differential geometry, I accepted the fact that we can everything express using tensor, that we enumerate components of tensor by indexes and this relationship is permanent. Once I saw that it was not the case, I found myself unprepared to recognize it It got month before I was able to write down expression that I saw by my eyes. The problem was because I tried to write down this expression in tensor form.

As consequence it appears that bivector space is particular case of direct product of $D_*{}^*$-vector spaces. We cannot define polylinear map. However consideration of polyadditive maps of direct product of $D_*{}^*$-vector spaces leads to definition of tensor product. Such way we can build tensors of arbitrary structure.

I designed this paper such way that it allows reproduce my not simple and interesting adventure in the world of tensors over division ring.

March, 2009

## 1.3. Preface to Version 3

Quaternion algebra is the simplest example of division ring. This is why I verify theorems in quaternion algebra to see how they work. Quaternion algebra is similar to complex field and so it is natural to search some parallel.

In real field any additive map automatically turns out to be linear over real field. This statement is caused by statement that real field is completion of rational field and is corollary of theorem 10.1.3.

However the statement changes for complex field. Not every additive map of complex field is linear over complex field. Conjugation is the simplest example of such map. As soon as I discovered this extremely interesting statement, I returned to question about analytical representation of additive map.

Exploring additive maps, I realized that I too abruptly expanded the set of linear maps while switching from field to division ring. The reason for this was not a very clear understanding of how to overcome the noncommutativity of the product. However during the process of research it became more evident that any additive map is linear over some field. For the first time I used this concept when constructing the tensor product and the concept took shape in a subsequent study.

During constructing the tensor product of division rings $D_1$, ..., $D_n$ I assumed existence of field $F$ such that additive map of division ring $D_i$ is linear over this field for any $i$. If all division rings have characteristic 0, then according to theorem 10.1.3 such field always exists. However dependence of tensor product on selected field $F$ arises here. To get rid of this dependence, I assume that field $F$ is maximal field that possesses stated property.

If $D_1 = ... = D_n = D$, then such field is center $Z(D)$ of division ring $D$. If the product in division ring $D$ is commutative[1.3], then $Z(D) = D$. So, starting with additive map, I arrive at concept of linear map, which is generalization of linear map over field.

Research in area of complex numbers and quaternions revealed one more interesting phenomenon. In spite of the fact that complex field is extension of real field, the structure of linear map over complex field is different from the structure of linear map over real field. This difference leads to the statement that conjugation of complex numbers is additive map but not linear map over complex field.

Similarly, the structure of linear map over division ring of quaternions is different from the structure of linear map over complex field. The source of difference is statement that the center of quaternion algebra has more simple structure than complex field. This difference leads to the statement that conjugation of quaternion satisfies to equation

$$\overline{p} = -\frac{1}{2}(p + ipi + jpj + kpk)$$

[1.3] In other words, division ring is field.

Consequently, the problem to find a mapping satisfying to a theorem similar to the Riemann theorem (theorem 11.1.1), is a nontrivial task for quaternions.

August, 2009

## 1.4. Preface to Version 5

Passenger: ... We will be late at the train station!

Little Engine: ... But if we do not see the first lily of the valley, then we will be late for the entire spring!

...

Passenger: ... We will be too late!

Little Engine: Yes. But if we do not hear the first nightingales, then we will be late for the entire summer!

...

Little Engine: Dawn! ... every sunrise is the only in our life! ... It is time to ride. After all we will be late.

Passenger: Yes. But if we do not see the dawn, we may be late for the whole life!

Gennadi Tsyferov, Little Engine from Romashkovo.

Young man John departed from point A to point B. Two hours later his friend Peter went on a bicycle from point $B$ in the direction to the point $A$. How fast did Peter ride, if they met in an hour?..

This is the familiar problem of school problem book. But the problem is not so simple as it is formulated. When Peter was near the pond he saw swans and stopped to make picture. John stopped in the grove to hear the nightingales. There they met.

Very often I tempted to take a look at the end of the problem book. The solution is almost there, just a jiffy away. However some uncertainty remains when you see the printed answer. I already see in general the theory that I am looking for. However I understand that I have a long way ahead of me before I will be able to tell about this to others. Besides, my experience tells ne that when you stroll along unfamiliar terrain, suddenly see such an amazing landscape that unintentionally stay there longer, to explore this place better.

However, once I could not resist and, in the paper [9], I tried to analyze how it may look like geometry over division ring. Although the paper was premature, the results formed the basis of my subsequent research.

---

In the paper [12], I explored the system of additive equations in Ω-algebra where the operation of addition is defined. The format of notation of equations in [12] is determined by the format of operations in Ω-algebra. This format is different from format accepted in this book by default. The theory of additive equations is generalization of the theory of ${}_*{}^*D$-linear equations and uses the same methods. Different order of variable and coefficient in equation made it difficult to compare results in both books. However it is not important for me whether I consider a system of $D^*{}_*$-linear equations or a system of ${}_*{}^*D$-linear equations. This stimulated me to edit text of this book and to consider the system of ${}_*{}^*D$-linear equations by default.

Besides, I realized that the notion of a set of additive mappings is artificial in the study of linear algebra over division ring. That is why I removed corresponding definitions and theorems from this version. Now I am directly considering a linear mapping over a given field.

The main difference of linear mappings over division ring is that I cannot represent a linear mapping as product of element of division ring over variable. In papers [11, 12], I found way to functionally separate variable and mapping. This allows me to represent a notation of linear mapping in the familiar form. However in this book I keep the previous notation. It has several reasons for this.

First of all, both forms of notation are correct. We can use any of these forms of notation choosing the one that seems more appropriate. It also enables the reader to compare both forms of notation. In addition, in [12], I consider linear mapping of algebras in parallel with their tensor product, since these two topics are inseparable. In this book, I decided to follow the same path that I covered when I worked on the second version of this book and present the tensor product as distinctive peak to ascend at the end of the book.

However, the road covered after the last version was reflected on last chapters of the book. For instance, since division ring is algebra over field, I decided to use the standard representation of indexes. The key changes are attributed with my perception of $D$-vector spaces. Let us take a closer look at this subject in greater detail since this subject remains outside of the scope of the book.

$D$-vector space appeared as corollary of theorem about twin representations of division ring in Abelian group (theorems 7.3.2, 7.3.3). Initially I considered $D$-vector space as consolidation of structures of $D^*{}_*$-vector space and ${}_*{}^*D$-vector space. Correspondingly I considered $D^*{}_*$-basis and ${}_*{}^*D$-basis.

Let $\overline{\overline{p}}$ be basis of ${}_*{}^*D$-vector space $V$ and basis $\overline{\overline{r}}$ be basis of $D^*{}_*$-vector space $V$. Vector ${}_ir$ of basis $\overline{\overline{r}}$ has expansion

$$ {}_ir = p_j\,{}_i^jR \qquad {}_ir = p_*{}^*{}_iR \tag{1.4.1} $$

relative to basis $\overline{\overline{p}}$. Vector $p_j$ of basis $\overline{\overline{p}}$ has expansion

$$ p_j = P_j^i\,{}_ir \qquad p_j = P_j{}^*{}_*r \tag{1.4.2} $$

relative to basis $\overline{\overline{r}}$.

It is easy to see from design that $R$ is coordinate matrix of basis $\overline{\overline{r}}$ relative to basis $\overline{\overline{p}}$. Columns of matrix $R$ are $D^*{}_*$-linearly independent.

In the same way, $P$ is coordinate matrix of basis $\overline{\overline{p}}$ relative to basis $\overline{\overline{p}}$. Columns of matrix $P$ are ${}_*{}^*D$-linearly independent.

From equations (1.4.1) and (1.4.2) it follows

$$ \begin{aligned} {}_ir &= p_j\,{}_i^jR = P_j^k\,{}_kr\,{}_i^jR \\ {}_ir &= p_*{}^*{}_iR = (P^*{}_*r)_*{}^*{}_iR \end{aligned} \tag{1.4.3} $$

From equation (1.4.3) we see that order of brackets is important.

Though matrices $P$ and $R$ are not inter inverse, we see that equation (1.4.3) represents identical transformation of $D$-vector space. It is possible to write this

transformation as

$$p_j = P_j^i{}_i r = P_j^i (p_k{}^k_i R)$$
$$p_j = P_j{}^*{}_* r = P_j{}^*{}_* (p_*{}^* R) \tag{1.4.4}$$

From comparison of equations (1.4.3) and (1.4.4) it follows that change of order of brackets changes order of summation. These equations express the symmetry in choice of bases $\overline{\overline{p}}$ and $\overline{\overline{r}}$.

A non-trivial presentation of identity map like expression (1.4.3) or (1.4.4) was a complete surprise to me.

There is another interesting statement that I saved for later analysis, namely the theorem 7.3.4, which states existence of $D^*{}_*$-basis which is not ${}_*{}^*D$-basis. We will make a pause here.

In [10], I explored the theory of representations of $\Omega$-algebra. According to this theory, $D$-vector space is representation of algebra $D \otimes D$ in Abelian group. Therefore, choice of basis is determined neither by $D\star$-linear dependence or by $\star D$-linear dependence, but by dependence like in expression

$$a^i_{s\cdot 0}\overline{e}_i a^i_{s\cdot 1}$$

Therefore, $D^*{}_*$-basis is not also basis of $D$-vector space if it is not ${}_*{}^*D$-bais. However, it is easy to see that the set of transformations of basis forms a group. This structure of basis leads to the fact that the difference between right and left bases disappears. This is why for $D$-vector space, I use standard representation of index.

However, even more surprising statement follows from this. Coordinates of vector of $D$-vector space belong to algebra $D \otimes D$. Few questions arise.

- Can a basis of $D$-vector space be ${}_*{}^*D$-basis?
- What is relation between dimensions of $D$-vector space and ${}_*{}^*D$-vector space?
- What is the structure of 1-dimensional $D$-vector space?

I intend to explore these questions for free modules over arbitrary algebra. This is the main reason why this topic is outside the scope of this book.

But we can already say that the set of ${}_*{}^*D$-linearly dependent vectors cannot be the basis of $D$-vector space. Therefore dimension of $D$-vector space does not exceed dimension of corresponding $D_*{}^*$-vector space. Therefore, I gave in this book proper description of $D$-vector space.

Well, I could not resist, and once again looked at the end of the problem book.

August, 2010

## 1.5. Conventions

Convention 1.5.1. *In any expression where we use index I assume that this index may have internal structure. For instance, considering the algebra $A$ we enumerate coordinates of $a \in A$ relative to basis $\overline{\overline{e}}$ by an index $i$. This means that $a$ is a vector. However, if $a$ is matrix, then we need two indexes, one enumerates rows, another enumerates columns. In the case, when index has structure, we begin the index from symbol $\cdot$ in the corresponding position. For instance, if I consider the matrix $a^i_j$ as an element of a vector space, then I can write the element of matrix as $a^{\cdot i}_{\cdot j}$.* □

Convention 1.5.2. *I assume sum over index $s$ in expression like*

$$a_{s\cdot 0}xa_{s\cdot 1}$$

□

Convention 1.5.3. *We can consider division ring $D$ as $D$-vector space of dimension 1. According to this statement, we can explore not only homomorphisms of division ring $D_1$ into division ring $D_2$, but also linear maps of division rings.* □

Convention 1.5.4. *In spite of noncommutativity of product a lot of statements remain to be true if we substitute, for instance, right representation by left representation or right vector space by left vector space. To keep this symmetry in statements of theorems I use symmetric notation. For instance, I consider $D\star$-vector space and $\star D$-vector space. We can read notation $D\star$-vector space as either $D$-star-vector space or left vector space.* □

Convention 1.5.5. *Let $A$ be free algebra with finite or countable basis. Considering expansion of element of algebra $A$ relative basis $\overline{\overline{e}}$ we use the same root letter to denote this element and its coordinates. In expression $a^2$, it is not clear whether this is component of expansion of element $a$ relative basis, or this is operation $a^2 = aa$. To make text clearer we use separate color for index of element of algebra. For instance,*

$$a = a^{i}e_{i}$$

□

Convention 1.5.6. *It is very difficult to draw the line between the module and the algebra. Especially since sometimes in the process of constructing, we must first prove that the set $A$ is a module, and then we prove that this set is an algebra. Therefore, to write the element of the module, we will also use the convention 1.5.5.* □

Convention 1.5.7. *The identification of the vector and matrix of its coordinates can lead to ambiguity in the equation*

$$a = a^{*}{}_{*}e \tag{1.5.1}$$

*where $\overline{\overline{e}}$ is a basis of vector space. Therefore, we write the equation (1.5.1) in the following form*

$$\overline{a} = a^{*}{}_{*}e$$

*in order to see where we wrote vector.* □

Convention 1.5.8. *If free finite dimensional algebra has unit, then we identify the vector of basis $\overline{e}_{0}$ with unit of algebra.* □

Without a doubt, the reader may have questions, comments, objections. I will appreciate any response.

CHAPTER 2

# Biring of Matrices

## 2.1. Concept of Generalized Index

Studying tensor calculus we start from studying univalent covariant and contravariant tensors. In spite on difference of properties both these objects are elements of respective vector spaces. Suppose we introduce a generalized index according to the rule $a^{\boldsymbol{i}} = a^i$, $b^{\boldsymbol{i}} = b^{\cdot -}_{\cdot i}$. Then we see that these tensors have the similar behavior. For instance, the transformation of a covariant tensor gets form

$$b'^{\boldsymbol{i}} = b'^{\cdot -}_{\cdot i} = f^{\cdot - \cdot j}_{\cdot i \cdot -} b^{\cdot -}_{\cdot j} = f^{\boldsymbol{i}}_{\boldsymbol{j}} b^{\boldsymbol{j}}$$

This similarity goes as far as we need because tensors also form vector space.

These observations of the similarity between properties of covariant and contravariant tensors lead us to the concept of generalized index. We will use the symbol $\cdot$ in front of a generalized index when we need to describe its structure. I put the sign $'-'$ in place of the index whose position was changed. For instance, if an original term was $a_{ij}$ I will use notation $a_i{}^j{}_-$ instead of notation $a_i^j$.

Even though the structure of a generalized index is arbitrary we assume that there exists a one-to-one map of the interval of positive integers 1, ..., $n$ to the range of index. Let $\mathcal{i}$ be the range of the index $\boldsymbol{i}$. We denote the power of this set by symbol $|\mathcal{i}|$ and assume that $|\mathcal{i}| = n$. If we want to enumerate elements $a_{\boldsymbol{i}}$ we use notation $a_1$, ..., $a_n$.

Representation of coordinates of a vector as a matrix allows making a notation more compact. The question of the presentation of vector as a row or a column of the matrix is just a question of convention. We extend the concept of generalized index to entries of the matrix. A matrix is a two dimensional table, the rows and columns of which are enumerated by generalized indexes. To represent a matrix we will use one of the following forms:

**Standard representation:** in this case we write entries of matrix $\mathcal{A}$ as $\mathcal{A}^{\boldsymbol{a}}_{\boldsymbol{b}}$.

**Alternative representation:** in this case we write entries of matrix $\mathcal{A}$ as ${}^{\boldsymbol{a}}\mathcal{A}_{\boldsymbol{b}}$ or ${}_{\boldsymbol{b}}\mathcal{A}^{\boldsymbol{a}}$.

Since we use generalized index, we cannot tell whether index $\boldsymbol{a}$ of matrix enumerates rows or columns until we know the structure of index.

We could use notation $*$-column and $*$-row which is more close to our custom. However as we can see bellow the form of presentation of matrix is not important for us. To make sure that notation offered below is consistent with the traditional we will assume that the matrix is presented in the form

$$A=\begin{pmatrix}{}^1A_1 & ... & {}^1A_n\\ ... & ... & ...\\ {}^mA_1 & ... & {}^mA_n\end{pmatrix}$$

DEFINITION 2.1.1. I use the following names and notation for different **minor matrices** of the matrix $A$

$A_a$ : $^*$**-row** with the index $a$ is generalization of a column of a matrix. The upper index enumerates entries of $^*$-rows and the lower index enumerates $^*$-rows.

$A_T$ : the minor matrix obtained from $A$ by selecting $^*$-rows with an index from the set $T$

$A_{[a]}$ : the minor matrix obtained from $A$ by deleting $^*$-row $A_a$

$A_{[T]}$ : the minor matrix obtained from $A$ by deleting $^*$-rows with an index from the set $T$

${}^bA$ : $_*$**-row** with the index $b$ is generalization of a row of a matrix. The lower index enumerates entries of $_*$-rows and the upper index enumerates $_*$-rows.

${}^SA$ : the minor matrix obtained from $A$ by selecting $_*$-rows with an index from the set $S$

${}^{[b]}A$ : the minor matrix obtained from $A$ by deleting $_*$-row ${}^bA$

${}^{[S]}A$ : the minor matrix obtained from $A$ by deleting $_*$-rows with an index from the set $S$

□

REMARK 2.1.2. We will combine the notation of indexes. Thus ${}^bA_a$ is $1\times 1$ minor matrix. The same time this is the notation for a matrix entry. This allows an identifying of $1\times 1$ matrix and its entry. The index $a$ is number of $^*$-row of matrix and the index $b$ is number of $_*$-rows of matrix. □

Each form of the notation of a matrix has its own advantages. The standard notation is more natural when we study matrix theory. The alternative form of the notation makes expressions in the theory of vector spaces more clear. Extending the alternative notation of indexes to arbitrary tensors we can better understand an interaction of different geometric objects. Using the duality principle (theorem 2.2.14) improves our expressivity.

REMARK 2.1.3. We can read symbol $^*$- as $c$- and symbol $_*$- as $r$- creating this way names $c$**-row** and $r$**-row**. Further we extend this rule to other objects of linear algebra. I will use this convention designing index. □

Since transpose of the matrix exchanges $_*$-rows and $^*$-rows we get equation

$$ {}_i(A^T)^j={}^iA_j \tag{2.1.1}$$

REMARK 2.1.4. As we can see from the equation (2.1.1), it is not important for us the choice of a side to place a number of $_*$-row and and the choice of a side to place a number of $^*$-row. This is due to the fact that we can enumerate the entries of matrix in different ways. If we want to show the numbers of $_*$-row and $^*$-row according to the definition 2.1.1, then the equation (2.1.1) has form

$${}^j(A^T)_i={}^iA_j$$

In standard representation, the equation (2.1.1) has form

$$(A^T)_i^j = A_j^i$$

□

We call matrix[2.1]

$$\mathcal{H}A = ({}^j(\mathcal{H}A)_i) = (({}^{\cdot}{}_i^{-}A_{\cdot}{}_{-}^{j})^{-1}) \tag{2.1.2}$$

**Hadamard inverse of matrix** $A = ({}_bA^a)$ ([7]-page 4).

I will use the Einstein convention about sums. This means that when an index is present in an expression twice and a set of index is known, I have the sum over this index. If needed to clearly show set of index, I will do it. Also, in this paper I will use the same root letter for a matrix and its entries.

We will study matrices entries of which belong to division ring $D$. We will also keep in mind that instead of division ring $D$ we may write in text field $F$. We will clearly write field $F$ in case when commutativity creates new details. We will denote by 1 identity element of division ring $D$.

Let $I$, $|I| = n$ be a set of indexes. We introduce the **Kronecker symbol**

$$\delta_j^i = \begin{cases} 1 & i = j \\ 0 & i \neq j \end{cases} \qquad i, j \in I$$

## 2.2. Biring

We consider matrices whose entries belong to division ring $D$.

The product of matrices is associated with the product of homomorphisms of vector spaces over field. According to the custom the product of matrices $A$ and $B$ is defined as product of ${}_*$-rows of the matrix $A$ and ${}^*$-rows of the matrix $B$. Conventional character of this definition becomes evident when we put attention that ${}_*$-row of the matrix $A$ may be a column of this matrix. In such case we multiply columns of the matrix $A$ over rows of the matrix $B$. Thus we can define two products of matrices. To distinguish between these products we introduced a new notation.[2.2]

DEFINITION 2.2.1. **${}_*{}^*$-product of matrices** $A$ and $B$ has form

$$\begin{cases} A_*{}^*B = ({}^aA_c\,{}^cB_b) \\ {}^a(A_*{}^*B)_b = {}^aA_c\,{}^cB_b \end{cases} \tag{2.2.1}$$

[2.1] The notation $({}^{\cdot}{}_i^{-}A_{\cdot}{}_{-}^{j})^{-1}$ means that we exchange rows and columns in Hadamard inverse. We can formally write this expression in following form

$$({}^{\cdot}{}_i^{-}A_{\cdot}{}_{-}^{j})^{-1} = \frac{1}{{}^iA_j}$$

[2.2] In order to keep this notation consistent with the existing one we assume that we have in mind ${}_*{}^*$-product when no clear notation is present.

and can be expressed as product of a ${}_*$-row of matrix $A$ over a ${}^*$-row of matrix $B$.[2.3] □

Definition 2.2.2. ${}^*{}_*$**-product of matrices** $A$ and $B$ has form

$$\begin{cases} A^*{}_*B = ({}_aA^c\ {}_cB^b) \\ {}_a(A^*{}_*B)^b = {}_aA^c\ {}_cB^b \end{cases} \tag{2.2.2}$$

and can be expressed as product of a ${}^*$-row of matrix $A$ over a ${}_*$-row of matrix $B$.[2.4] □

Remark 2.2.3. We will use symbol ${}_*{}^*$- or ${}^*{}_*$- in name of properties of each product and in the notation. According to remark 2.1.3 we can read symbols ${}_*{}^*$ and ${}^*{}_*$ as $rc$-product and $cr$-product. This rule we extend to following terminology. □

Remark 2.2.4. Just as in remark 2.1.4, I want to draw attention to the fact that I change the numbering of entries of the matrix. If we want to show the numbers of ${}_*$-row and ${}^*$-row according to the definition 2.1.1, then the equation (2.2.2) has form

$$ {}^b(A^*{}_*B)_a = {}^cA_a\ {}^bB_c \tag{2.2.3}$$

However the format of the equation (2.2.3) is unusual. □

Set of $n\times n$ matrices is closed relative ${}_*{}^*$-product and ${}^*{}_*$-product as well relative sum which is defined by rule

$$(A+B)^b_a = A^b_a + B^b_a$$

Theorem 2.2.5.

$$(A_*{}^*B)^T = A^T{}^*{}_*B^T \tag{2.2.4}$$

---

[2.3] In alternative form operation consists from two symbols $*$ which we put in the place of index which participate in sum. In standard notation we write operation as

$$\begin{cases} A_*{}^*B = (A^a_c B^c_b) \\ (A_*{}^*B)^a_b = A^a_c B^c_b \end{cases}$$

and can be construed as symbolic notation

$$A_*{}^*B = A_*B^*$$

where we write symbol $*$ on place of index which participate in sum.

[2.4] In alternative form operation consists from two symbols $*$ which we put in the place of index which participate in sum. In standard notation we write operation as

$$\begin{cases} A^*{}_*B = (A^c_a\ B^b_c) \\ (A^*{}_*B)^a_b = A^c_a\ B^b_c \end{cases}$$

and can be construed as symbolic notation

$$A^*{}_*B = A^*B_*$$

where we write symbol $*$ on place of index which participate in sum.

PROOF. The chain of equations

$$
\begin{aligned}
{}_a((A_*{}^*B)^T)^b &= {}^a(A_*{}^*B)_b \\
&= {}^aA_c\,{}^cB_b \\
&= {}_a(A^T)^c{}_c(B^T)^b \\
&= {}_a((A^T)^*{}_*(B^T))^b
\end{aligned}
\tag{2.2.5}
$$

follows from (2.1.1), (2.2.1) and (2.2.2). The equation (2.2.4) follows from (2.2.5). □

Matrix $\delta = (\delta^c_a)$ is identity for both products.

DEFINITION 2.2.6. $\mathcal{A}$ is a **biring** if we defined on $\mathcal{A}$ an unary operation, say transpose, and three binary operations, say ${}_*{}^*$-product, ${}^*{}_*$-product and sum, such that

- ${}_*{}^*$-product and sum define structure of ring on $\mathcal{A}$
- ${}^*{}_*$-product and sum define structure of ring on $\mathcal{A}$
- both products have common identity $\delta$
- products satisfy equation

$$(A_*{}^*B)^T = A^T{}^*{}_*B^T \tag{2.2.6}$$

- transpose of identity is identity

$$\delta^T = \delta \tag{2.2.7}$$

- double transpose is original element

$$(A^T)^T = A \tag{2.2.8}$$

□

THEOREM 2.2.7.

$$(A^*{}_*B)^T = (A^T)_*{}^*(B^T) \tag{2.2.9}$$

PROOF. We can prove (2.2.9) in case of matrices the same way as we proved (2.2.6). However it is more important for us to show that (2.2.9) follows directly from (2.2.6).

Applying (2.2.8) to each term in left side of (2.2.9) we get

$$(A^*{}_*B)^T = ((A^T)^T{}^*{}_*(B^T)^T)^T \tag{2.2.10}$$

From (2.2.10) and (2.2.6) it follows that

$$(A^*{}_*B)^T = ((A^T{}_*{}^*B^T)^T)^T \tag{2.2.11}$$

(2.2.9) follows from (2.2.11) and (2.2.8). □

DEFINITION 2.2.8. We introduce ${}_*{}^*$**-power** of element $A$ of biring $\mathcal{A}$ using recursive definition

$$A^{0_*{}^*} = \delta \tag{2.2.12}$$

$$A^{n_*{}^*} = A^{n-1_*{}^*}{}_*{}^*A \tag{2.2.13}$$

□

Definition 2.2.9. We introduce ${}^*{}_*$-**power** of element $A$ of biring $\mathcal{A}$ using recursive definition

$$A^{0^*{}_*} = \delta \tag{2.2.14}$$

$$A^{n^*{}_*} = A^{n-1^*{}_*}{}^*{}_*A \tag{2.2.15}$$

□

Theorem 2.2.10.

$$(A^T)^{n_*{}^*} = (A^{n^*{}_*})^T \tag{2.2.16}$$

$$(A^T)^{n^*{}_*} = (A^{n_*{}^*})^T \tag{2.2.17}$$

Proof. We proceed by induction on $n$.

For $n = 0$ the statement immediately follows from equations (2.2.12), (2.2.14), and (2.2.7).

Suppose the statement of theorem holds when $n = k - 1$

$$(A^T)^{n-1_*{}^*} = (A^{n-1^*{}_*})^T \tag{2.2.18}$$

It follows from (2.2.13) that

$$(A^T)^{k_*{}^*} = (A^T)^{k-1_*{}^*}{}_*{}^*A^T \tag{2.2.19}$$

It follows from (2.2.19) and (2.2.18) that

$$(A^T)^{k_*{}^*} = (A^{k-1^*{}_*})^T{}_*{}^*A^T \tag{2.2.20}$$

It follows from (2.2.20) and (2.2.9) that

$$(A^T)^{k_*{}^*} = (A^{k-1^*{}_*}{}^*{}_*A)^T \tag{2.2.21}$$

(2.2.16) follows from (2.2.19) and (2.2.15).

We can prove (2.2.17) by similar way. □

Definition 2.2.11. Element $A^{-1_*{}^*}$ of biring $\mathcal{A}$ is ${}_*{}^*$-**inverse element** of element $A$ if

$$A{}_*{}^*A^{-1_*{}^*} = \delta \tag{2.2.22}$$

Element $A^{-1^*{}_*}$ of biring $\mathcal{A}$ is ${}^*{}_*$-**inverse element** of element $A$ if

$$A{}^*{}_*A^{-1^*{}_*} = \delta \tag{2.2.23}$$

□

Theorem 2.2.12. *Suppose element $A \in \mathcal{A}$ has ${}_*{}^*$-inverse element. Then transpose element $A^T$ has ${}^*{}_*$-inverse element and these elements satisfy equation*

$$(A^T)^{-1^*{}_*} = (A^{-1_*{}^*})^T \tag{2.2.24}$$

*Suppose element $A \in \mathcal{A}$ has ${}^*{}_*$-inverse element. Then transpose element $A^T$ has ${}_*{}^*$-inverse element and these elements satisfy equation*

$$(A^T)^{-1_*{}^*} = (A^{-1^*{}_*})^T \tag{2.2.25}$$

PROOF. If we get transpose of both side (2.2.22) and apply (2.2.7) we get

$$(A_*{}^*A^{-1_*{}^*})^T = \delta^T = \delta$$

Applying (2.2.6) we get

$$\delta = A^T{}^*{}_*(A^{-1_*{}^*})^T \tag{2.2.26}$$

(2.2.24) follows from comparison (2.2.23) and (2.2.26).

We can prove (2.2.25) similar way. □

Theorems 2.2.5, 2.2.7, 2.2.10, and 2.2.12 show that some kind of duality exists between ${}_*{}^*$-product and ${}^*{}_*$-product. We can combine these statements.

THEOREM 2.2.13 (**duality principle for biring**). *Let $\mathfrak{A}$ be true statement about biring $\mathcal{A}$. If we exchange the same time*

- *$A \in \mathcal{A}$ and $A^T$*
- *${}_*{}^*$-product and ${}^*{}_*$-product*

*then we soon get true statement.*

THEOREM 2.2.14 (**duality principle for biring of matrices**). *Let $\mathcal{A}$ be biring of matrices. Let $\mathfrak{A}$ be true statement about matrices. If we exchange the same time*

- *${}^*$-rows and ${}_*$-rows of all matrices*
- *${}_*{}^*$-product and ${}^*{}_*$-product*

*then we soon get true statement.*

PROOF. This is the immediate consequence of the theorem 2.2.13. □

REMARK 2.2.15. We execute operations in expression

$$A_*{}^*B_*{}^*C$$

from left to right. However we can execute product from right to left. In custom notation this expression is

$$C^*{}_*B^*{}_*A$$

We follow the rule that to write power from right of expression. If we use standard representation, then we write indexes from right of expression. If we use alternative representation, then we read indexes in the same order as symbols of operation and root letters. For instance, let original expression be like

$$A^{-1_*{}^*}{}_*{}^*B_a$$

Then expression which we reed from right to left is like

$$B_a{}^*{}_*A^{-1^*{}_*}$$

in standard representation or

$${}_aB^*{}_*A^{-1^*{}_*}$$

in alternative representation.

Suppose we established the order in which we write indexes. Then we state that we read an expression from top to bottom reading first upper indexes, then lower ones. We assume that this is standard form of reading. We can read this expression from bottom to top. We extend this rule stating that we read symbols of operation in the same order as indexes. For instance, if we read expression

$$A^a{}_*{}^*B^{-1_*{}^*} = C^a$$

from bottom to top, then we can write this expression in standard form

$$A_a{}^*{}_*B^{-1^*{}_*} = C_a$$

According to the duality principle if we can prove one statement then we can prove other as well. □

THEOREM 2.2.16. *Let matrix $A$ have ${}_*{}^*$-inverse matrix. Then for any matrices $B$ and $C$ equation*

$$B = C \tag{2.2.27}$$

*follows from the equation*

$$B{}_*{}^*A = C{}_*{}^*A \tag{2.2.28}$$

PROOF. Equation (2.2.27) follows from the equation (2.2.28) if we multiply both parts of the equation (2.2.28) over $A^{-1_*{}^*}$. □

## 2.3. Quasideterminant

THEOREM 2.3.1. *Suppose $n\times n$ matrix $A$ has ${}_*{}^*$-inverse matrix.*[2.5] *Then $k\times k$ minor matrix of ${}_*{}^*$-inverse matrix satisfy to the equation*

$$\left({}^I(A^{-1_*{}^*})_J\right)^{-1_*{}^*} = {}^JA_I - {}^JA_{[I]*}{}^*\left({}^{[J]}A_{[I]}\right)^{-1_*{}^*}{}_*{}^{*[J]}A_I \tag{2.3.1}$$

PROOF. Definition (2.2.22) of ${}_*{}^*$-inverse matrix leads to the system of linear equations

$${}^{[J]}A_{[I]*}{}^{*[I]}(A^{-1_*{}^*})_J + {}^{[J]}A_{I*}{}^{*I}(A^{-1_*{}^*})_J = 0 \tag{2.3.2}$$

$${}^JA_{[I]*}{}^{*[I]}(A^{-1_*{}^*})_J + {}^JA_{I*}{}^{*I}(A^{-1_*{}^*})_J = \delta \tag{2.3.3}$$

We multiply (2.3.2) by $\left({}^{[J]}A_{[I]}\right)^{-1_*{}^*}$

$${}^{[I]}(A^{-1_*{}^*})_J + \left({}^{[J]}A_{[I]}\right)^{-1_*{}^*}{}_*{}^{*[J]}A_{I*}{}^{*I}(A^{-1_*{}^*})_J = 0 \tag{2.3.4}$$

Now we can substitute (2.3.4) into (2.3.3)

$$-{}^JA_{[I]*}{}^*\left({}^{[J]}A_{[I]}\right)^{-1_*{}^*}{}_*{}^{*[J]}A_{I*}{}^{*I}(A^{-1_*{}^*})_J + {}^JA_{I*}{}^{*I}(A^{-1_*{}^*})_J = \delta \tag{2.3.5}$$

(2.3.1) follows from (2.3.5). □

COROLLARY 2.3.2. *Suppose $n\times n$ matrix $A$ has ${}_*{}^*$-inverse matrix. Then entries of ${}_*{}^*$-inverse matrix satisfy to the equation*[2.1]

$${}^i(A^{-1_*{}^*})_j = \left({}^jA_i - {}^jA_{[i]*}{}^*\left({}^{[j]}A_{[i]}\right)^{-1_*{}^*}{}_*{}^{*[j]}A_i\right)^{-1} \tag{2.3.6}$$

$${}^j\left(\mathcal{H}A^{-1_*{}^*}\right)_i = {}^jA_i - {}^jA_{[i]*}{}^*\left({}^{[j]}A_{[i]}\right)^{-1_*{}^*}{}_*{}^{*[j]}A_i \tag{2.3.7}$$

□

[2.5] This statement and its proof is based on statement 1.2.1 from [6] (page 8) for matrix over free division ring.

EXAMPLE 2.3.3. Consider matrix

$$\begin{pmatrix} {}^1A_1 & {}^1A_2 \\ {}^2A_1 & {}^2A_2 \end{pmatrix}$$

According to (2.3.6)

$$ {}^1(A^{-1_*{}^*})_1 = ({}^1A_1 - {}^1A_2({}^2A_2)^{-1}\ {}^2A_1)^{-1} \tag{2.3.8}$$

$$ {}^2(A^{-1_*{}^*})_1 = ({}^1A_2 - {}^1A_1({}^2A_1)^{-1}\ {}^2A_2)^{-1} \tag{2.3.9}$$

$$ {}^1(A^{-1_*{}^*})_2 = ({}^2A_1 - {}^2A_2({}^1A_2)^{-1}\ {}^1A_1)^{-1} \tag{2.3.10}$$

$$ {}^2(A^{-1_*{}^*})_2 = ({}^2A_2 - {}^2A_1({}^1A_1)^{-1}\ {}^1A_2)^{-1} \tag{2.3.11}$$

$$A^{-1^*{}_*} = \begin{pmatrix} ({}_1A^1 - {}_1A^2({}_2A^2)^{-1}\ {}_2A^1)^{-1} & ({}_1A^2 - {}_1A^1({}_2A^1)^{-1}\ {}_2A^2)^{-1} \\ ({}_2A^1 - {}_2A^2({}_1A^2)^{-1}\ {}_1A^1)^{-1} & ({}_2A^2 - {}_2A^1({}_1A^1)^{-1}\ {}_1A^2)^{-1} \end{pmatrix}$$

□

According to [6], page 3 we do not have an appropriate definition of a determinant for a division ring. However, we can define a quasideterminant which finally gives a similar picture. In definition below we follow definition [6]-1.2.2.

DEFINITION 2.3.4. $\binom{j}{i}$-${}_*{}^*$**-quasideterminant** of $n\times n$ matrix $A$ is formal expression[2.1]

$$ {}^j\det({}_*{}^*)_i\,A = {}^j\left(\mathcal{H}A^{-1_*{}^*}\right)_i \tag{2.3.12}$$

According to the remark 2.1.2 we consider $\binom{j}{i}$-${}_*{}^*$-quasideterminant as an entry of the matrix $\det({}_*{}^*)\,A$ which is called ${}_*{}^*$**-quasideterminant**. □

THEOREM 2.3.5. *Expression for ${}_*{}^*$-inverse matrix has form*

$$A^{-1_*{}^*} = \mathcal{H}\det({}_*{}^*)A \tag{2.3.13}$$

PROOF. (2.3.13) follows from (2.3.12). □

THEOREM 2.3.6. *Expression for $\binom{j}{i}$-${}_*{}^*$-quasideterminant can be evaluated by either form*[2.6]

$$ {}^j\det({}_*{}^*)_i\,A = {}^jA_i - {}^jA_{[i]}{}_*{}^*\left({}^{[j]}A_{[i]}\right)^{-1_*{}^*}{}_*{}^*{}^{[j]}A_i \tag{2.3.14}$$

$$ {}^j\det({}_*{}^*)_i\,A = {}^jA_i - {}^jA_{[i]}{}_*{}^*\mathcal{H}\det({}_*{}^*)^{[j]}A_{[i]}{}_*{}^*{}^{[j]}A_i \tag{2.3.15}$$

PROOF. Statement follows from (2.3.7) and (2.3.12). □

[2.6]We can provide similar proof for $\binom{j}{i}$-${}^*{}_*$**-quasideterminant**. However we can write corresponding statement using the duality principle. Thus, if we read equation (2.3.14) from right to left, we get equation

$$ {}^j\det({}^*{}_*)A_i = {}^jA_i - {}^{[j]}A_i{}^*{}_*\left({}^{[j]}A_{[i]}\right)^{-1^*{}_*}{}^*{}_*{}^jA_{[i]}$$

$$ {}^j\det({}^*{}_*)A_i = {}^jA_i - {}^{[j]}A_i{}^*{}_*\mathcal{H}\det({}^*{}_*)^{[j]}A_{[i]}{}^*{}_*{}^jA_{[i]}$$

Example 2.3.7. Consider matrix

$$\begin{pmatrix} {}^1A_1 & {}^1A_2 \\ {}^2A_1 & {}^2A_2 \end{pmatrix}$$

According to (2.3.14)

$$\det({}_*{}^*)A = \begin{pmatrix} {}^1A_1 - {}^1A_2({}^2A_2)^{-1}\,{}^2A_1 & {}^1A_2 - {}^1A_1({}^2A_1)^{-1}\,{}^2A_2 \\ {}^2A_1 - {}^2A_2({}^1A_2)^{-1}\,{}^1A_1 & {}^2A_2 - {}^2A_1({}^1A_1)^{-1}\,{}^1A_2 \end{pmatrix}$$

$$\det({}^*{}_*)A = \begin{pmatrix} {}_1A^1 - {}_1A^2({}_2A^2)^{-1}\,{}_2A^1 & {}_2A^1 - {}_2A^2({}_1A^2)^{-1}\,{}_1A^1 \\ {}_1A^2 - {}_1A^1({}_2A^1)^{-1}\,{}_2A^2 & {}_2A^2 - {}_2A^1({}_1A^1)^{-1}\,{}_1A^2 \end{pmatrix}$$

□

Theorem 2.3.8.

$$(2.3.16)\qquad {}_j\det({}_*{}^*)^i\,A^T = {}^j\det({}^*{}_*)_i\,A$$

Proof. According to (2.3.12) and (2.1.2)

$${}_j\det({}_*{}^*)^i\,A^T = ({}^{\ j}_{\cdot\,-}((A^T)^{-1_*{}^*})^{\cdot\,-}_{\ \ i})^{-1}$$

Using theorem 2.2.12 we get

$${}_j\det({}_*{}^*)^i\,A^T = ({}^{\ j}_{\cdot\,-}((A^{-1^*{}_*})^T)^{\cdot\,-}_{\ \ i})^{-1}$$

Using (2.1.1) we get

$$(2.3.17)\qquad {}_j\det({}_*{}^*)^i\,A^T = ({}^{\cdot\,-}_{\ \ j}(A^{-1^*{}_*}){}^{\ i}_{\cdot\,-})^{-1}$$

Using (2.3.17), (2.1.2), (2.3.12) we get (2.3.16). □

The theorem 2.3.8 extends the duality principle stated in the theorem 2.2.14 to statements on quasideterminants and tells us that the same expression is ${}_*{}^*$-quasideterminant of matrix $A$ and ${}^*{}_*$-quasideterminant of matrix $A^T$. Using this theorem, we can write any statement for ${}^*{}_*$-matrix on the basis of similar statement for ${}_*{}^*$-matrix.

Theorem 2.3.9 (duality principle). *Let $\mathfrak{A}$ be true statement about matrix biring. If we exchange the same time*

- ${}_*$*-row and* ${}^*$*-row*
- ${}_*{}^*$*-quasideterminant and* ${}^*{}_*$*-quasideterminant*

*then we soon get true statement.*

Theorem 2.3.10.

$$(2.3.18)\qquad (mA)^{-1_*{}^*} = A^{-1_*{}^*}m^{-1}$$

$$(2.3.19)\qquad (Am)^{-1_*{}^*} = m^{-1}A^{-1_*{}^*}$$

Proof. To prove equation (2.3.18) we proceed by induction on size of the matrix.

Since

$$(mA)^{-1_*{}^*} = \left((mA)^{-1}\right) = \left(A^{-1}m^{-1}\right) = \left(A^{-1}\right)m^{-1} = A^{-1_*{}^*}m^{-1}$$

the statement is evident for $1\times 1$ matrix.

Let the statement holds for $(n-1)\times(n-1)$ matrix. Then from equation (2.3.1) it follows that

$$\begin{aligned}({}^I((mA)^{-1_*{}^*})_J)^{-1_*{}^*} &= {}^J(mA)_I - {}^J(mA)_{[I]*}{}^* \left({}^{[J]}(mA)_{[I]}\right)^{-1_*{}^*}{}_*{}^{*[J]}(mA)_I \\ &= m{}^JA_I - m\ {}^JA_{[I]*}{}^*\left({}^{[J]}A_{[I]}\right)^{-1_*{}^*}\ m^{-1}{}_*{}^*m\ {}^{[J]}A_I \\ &= m\ {}^JA_I - m\ {}^JA_{[I]*}{}^*\left({}^{[J]}A_{[I]}\right)^{-1_*{}^*}{}_*{}^{*[J]}A_I\end{aligned}$$

$$({}^I((mA)^{-1_*{}^*})_J)^{-1_*{}^*} = m\ {}^I(A^{-1_*{}^*})_J \tag{2.3.20}$$

The equation (2.3.18) follows from the equation (2.3.20). In the same manner we prove the equation (2.3.19). □

THEOREM 2.3.11. *Let*

$$A=\begin{pmatrix}1&0\\0&1\end{pmatrix} \tag{2.3.21}$$

*Then*

$$A^{-1_*{}^*}=\begin{pmatrix}1&0\\0&1\end{pmatrix} \tag{2.3.22}$$

$$A^{-1^*{}_*}=\begin{pmatrix}1&0\\0&1\end{pmatrix} \tag{2.3.23}$$

PROOF. It is clear from (2.3.8) and (2.3.11) that ${}^1(A^{-1_*{}^*})_1=1$ and ${}^2(A^{-1_*{}^*})_2=1$. However expression for ${}^2(A^{-1_*{}^*})_1$ and ${}^1(A^{-1_*{}^*})_2$ cannot be defined from (2.3.9) and (2.3.10) since ${}^2A_1={}^1A_2=0$. We can transform these expressions. For instance

$$\begin{aligned}{}^2(A^{-1_*{}^*})_1 &= ({}^1A_2 - {}^1A_1({}^2A_1)^{-1}\ {}^2A_2)^{-1} \\ &= ({}^1A_1(({}^1A_1)^{-1}\ {}^1A_2 - ({}^2A_1)^{-1}\ {}^2A_2))^{-1} \\ &= (({}^2A_1)^{-1}\ {}^1A_1({}^2A_1({}^1A_1)^{-1}\ {}^1A_2 - {}^2A_2))^{-1} \\ &= ({}^1A_1({}^2A_1({}^1A_1)^{-1}\ {}^1A_2 - {}^2A_2))^{-1}\ {}^2A_1\end{aligned}$$

It follows immediately that ${}^2(A^{-1_*{}^*})_1=0$. In the same manner we can find that ${}^1(A^{-1_*{}^*})_2=0$. This completes the proof of (2.3.22).

Equation (2.3.23) follows from (2.3.22), theorem 2.3.8 and symmetry of matrix (2.3.21). □

CHAPTER 3

# Representation of Universal Algebra

## 3.1. Representation of Universal Algebra

DEFINITION 3.1.1. Let the structure of $\Omega_2$-algebra be defined on the set $M$ ([2, 13]). Endomorphism of $\Omega_2$-algebra

$$t : M \to M$$

is called **transformation of universal algebra** $M$.[3.1] □

We denote $\delta$ identical transformation.

DEFINITION 3.1.2. Let ${}^*M$ be the set of **left-side transformations**

$$u' = tu$$

of $\Omega_2$-algebra $M$. Let ${}^*M$ be $\Omega_1$-algebra. The homomorphism

$$f : A \to {}^*M \tag{3.1.1}$$

of $\Omega_1$-algebra $A$ into $\Omega_1$-algebra ${}^\star M$ is called **left-side representation of $\Omega_1$-algebra** $A$ or $A*$**-representation** in $\Omega_2$-algebra $M$. □

DEFINITION 3.1.3. Let $M^*$ be the set of **right-side transformations**

$$u' = ut$$

of $\Omega_2$-algebra $M$. Let $M^*$ be $\Omega_1$-algebra. The homomorphism

$$f : A \to M^*$$

of $\Omega_1$-algebra $A$ into $\Omega_1$-algebra $M^*$ is called **right-side representation of $\Omega_1$-algebra** $A$ or $*A$**-representation** in $\Omega_2$-algebra $M$. □

We extend to representation theory convention described in remark 2.2.15. We can write duality principle in the following form

THEOREM 3.1.4 (duality principle). *Any statement which holds for left-side representation of $\Omega_1$-algebra $A$ holds also for right-side representation of $\Omega_1$-algebra $A$.*

REMARK 3.1.5. There exist two forms of notation for transformation of $\Omega_2$-algebra $M$. In operational notation, we write the transformation $A$ as either $Aa$ which corresponds to the left-side transformation or $aA$ which corresponds to the right-side transformation. In functional notation, we write the transformation $A$

[3.1] If the set of operations of $\Omega_2$-algebra is empty, then

$$t : M \to M$$

is a map.

as $A(a)$ regardless of the fact whether this is left-side or right-side transformation. This notation is in agreement with duality principle.

This remark serves as a basis for the following convention. When we use functional notation we do not make a distinction whether this is left-side or right-side transformation. We denote ${}^*M$ the set of transformations of $\Omega_2$-algebra $M$. Suppose we defined the structure of $\Omega_1$-algebra on the set ${}^*M$. Let $A$ be $\Omega_1$-algebra. We call homomorphism

$$f : A \rightarrow {}^*M \tag{3.1.2}$$

**representation of $\Omega_1$-algebra $A$ in $\Omega_2$-algebra $M$**. We also use record

$$f : A \xrightarrow{*} M$$

to denote the representation of $\Omega_1$-algebra $A$ in $\Omega_2$-algebra $M$.

Correspondence between operational notation and functional notation is unambiguous. We can select any form of notation which is convenient for presentation of particular subject. □

There are several ways to describe the representation. We can define the map $f$ keeping in mind that the domain is $\Omega_1$-algebra $A$ and range is $\Omega_1$-algebra ${}^*M$. Either we can specify $\Omega_1$-algebra $A$ and $\Omega_2$-algebra $M$ keeping in mind that we know the structure of the map $f$.[3.2]

Diagram

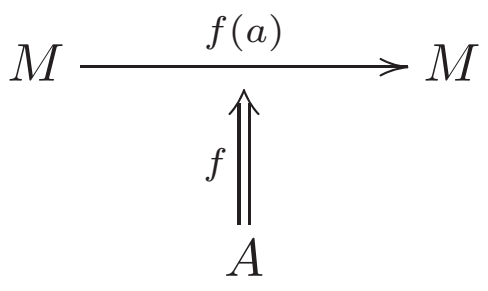

means that we consider the representation of $\Omega_1$-algebra $A$. The map $f(a)$ is image of $a \in A$.

Definition 3.1.6. Let the map (3.1.2) be an isomorphism of the $\Omega_1$-algebra $A$ into ${}^*M$. Then the representation of the $\Omega_1$-algebra $A$ is called **effective**. □

Remark 3.1.7. If the left-side representation of $\Omega_1$-algebra is effective, then we identify an element of $\Omega_1$-algebra and its image and write left-side transformation caused by element $a \in A$ as

$$v' = av$$

If the right-side representation of $\Omega_1$-algebra is effective, then we identify an element of $\Omega_1$-algebra and its image and write right-side transformation caused by element $a \in A$ as

$$v' = va$$

□

Definition 3.1.8. We call a representation of $\Omega_1$-algebra **transitive** if for any $a, b \in V$ exists such $g$ that

$$a = f(g)(b)$$

We call a representation of $\Omega_1$-algebra **single transitive** if it is transitive and effective. □

[3.2] For instance, we consider vector space $\overline{V}$ over field $D$ ( definition 5.1.4 ).

Theorem 3.1.9. *Representation is single transitive iff for any $a, b \in M$ exists one and only one $g \in A$ such that $a = f(g)(b)$*

Proof. Corollary of definitions 3.1.6 and 3.1.8. □

## 3.2. Morphism of Representations of Universal Algebra

Theorem 3.2.1. *Let $A$ and $B$ be $\Omega_1$-algebras. Representation of $\Omega_1$-algebra $B$*

$$g : B \rightarrow {}^*M$$

*and homomorphism of $\Omega_1$-algebra*

$$h : A \rightarrow B \tag{3.2.1}$$

*define representation $f$ of $\Omega_1$-algebra $A$*

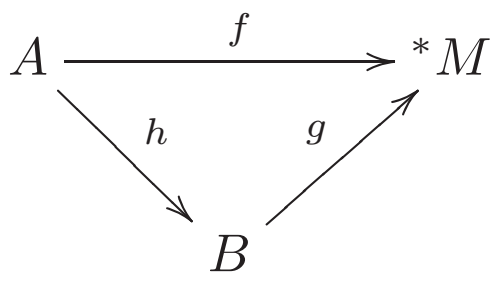

Proof. Since map $g$ is homomorphism of $\Omega_1$-algebra $B$ into $\Omega_1$-algebra ${}^*M$, the map $f$ is homomorphism of $\Omega_1$-algebra $A$ into $\Omega_1$-algebra ${}^*M$. □

Considering representations of $\Omega_1$-algebra in $\Omega_2$-algebras $M$ and $N$, we are interested in a map that preserves the structure of representation.

Definition 3.2.2. Let

$$f : A \rightarrow {}^*M$$

be representation of $\Omega_1$-algebra $A$ in $\Omega_2$-algebra $M$ and

$$g : B \rightarrow {}^*N$$

be representation of $\Omega_1$-algebra $B$ in $\Omega_2$-algebra $N$. Tuple of maps

$$(r : A \rightarrow B, R : M \rightarrow N) \tag{3.2.2}$$

such, that

- $r$ is homomorphism of $\Omega_1$-algebra
- $R$ is homomorphism of $\Omega_2$-algebra
- $$R \circ f(a) = g(r(a)) \circ R \tag{3.2.3}$$

is called **morphism of representations from $f$ into $g$**. We also say that **morphism of representations of $\Omega_1$-algebra in $\Omega_2$-algebra** is defined. □

Remark 3.2.3. We may consider a pair of mappins $r$, $R$ as map

$$F : A \cup M \rightarrow B \cup N$$

such that

$$F(A) = B \quad F(M) = N$$

Therefore, hereinafter we will say that we have the map $(r, R)$. □

Definition 3.2.4. If representation $f$ and $g$ coincide, then morphism of representations $(r, R)$ is called **morphism of representation $f$**. □

For any $m \in M$ equation (3.2.3) has form

$$R(f(a)(m)) = g(r(a))(R(m)) \tag{3.2.4}$$

Remark 3.2.5. Consider morphism of representations (3.2.2). We denote elements of the set $B$ by letter using pattern $b \in B$. However if we want to show that $b$ is image of element $a \in A$, we use notation $r(a)$. Thus equation

$$r(a) = r(a)$$

means that $r(a)$ (in left part of equation) is image $a \in A$ (in right part of equation). Using such considerations, we denote element of set $N$ as $R(m)$. We will follow this convention when we consider correspondences between homomorphisms of $\Omega_1$-algebra and maps between sets where we defined corresponding representations. □

Remark 3.2.6. There are two ways to interpret (3.2.4)

- Let transformation $f(a)$ map $m \in M$ into $f(a)(m)$. Then transformation $g(r(a))$ maps $R(m) \in N$ into $R(f(a)(m))$.
- We represent morphism of representations from $f$ into $g$ using diagram

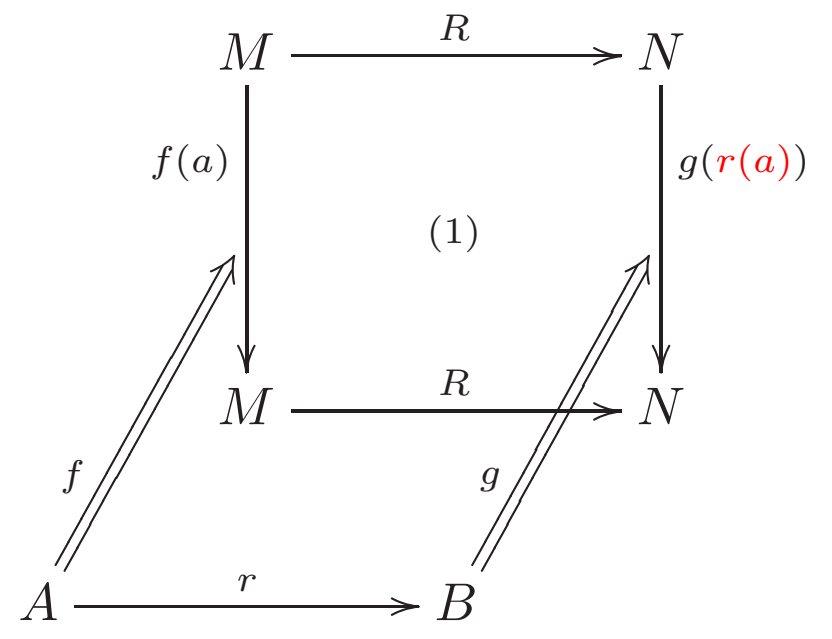

From (3.2.3), it follows that diagram (1) is commutative.

□

Theorem 3.2.7. *Consider representation*

$$f : A \rightarrow {}^*M$$

*of $\Omega_1$-algebra $A$ and representation*

$$g : B \rightarrow {}^*N$$

*of $\Omega_1$-algebra $B$. Morphism*

$$h : A \longrightarrow B \qquad H : M \longrightarrow N$$

*of representations from $f$ into $g$ satisfies equation*

$$H \circ \omega(f(a_1), ..., f(a_n)) = \omega(g(h(a_1)), ..., g(h(a_n))) \circ H \tag{3.2.5}$$

*for any $n$-ary operation $\omega$ of $\Omega_1$-algebra.*

Proof. Since $f$ is homomorphism, we have

$$H \circ \omega(f(a_1), ..., f(a_n)) = H \circ f(\omega(a_1, ..., a_n)) \tag{3.2.6}$$

From (3.2.3) and (3.2.6) it follows that

$$H \circ \omega(f(a_1), ..., f(a_n)) = g(h(\omega(a_1, ..., a_n))) \circ H \tag{3.2.7}$$

Since $h$ is homomorphism, from (3.2.7) it follows that

$$H\circ\omega(f(a_1),...,f(a_n))=g(\omega(h(a_1),...,h(a_n)))\circ H \tag{3.2.8}$$

Since $g$ is homomorphism, (3.2.5) follows from (3.2.8). □

THEOREM 3.2.8. *Let the map*

$$h:A\longrightarrow B\qquad H:M\longrightarrow N$$

*be morphism from representation*

$$f:A\rightarrow {}^*M$$

*of $\Omega_1$-algebra $A$ into representation*

$$g:B\rightarrow {}^*N$$

*of $\Omega_1$-algebra $B$. If representation $f$ is effective, then the map*

$${}^*H:{}^*M\rightarrow {}^*N$$

*defined by equation*

$${}^*H(f(a))=g(h(a)) \tag{3.2.9}$$

*is homomorphism of $\Omega_1$-algebra.*

PROOF. Because representation $f$ is effective, then for given transformation $f(a)$ element $a$ is determined uniquely. Therefore, transformation $g(h(a))$ is properly defined in equation (3.2.9).

Since $f$ is homomorphism, we have

$${}^*H(\omega(f(a_1),...,f(a_n)))={}^*H(f(\omega(a_1,...,a_n))) \tag{3.2.10}$$

From (3.2.9) and (3.2.10) it follows that

$${}^*H(\omega(f(a_1),...,f(a_n)))=g(h(\omega(a_1,...,a_n))) \tag{3.2.11}$$

Since $h$ is homomorphism, from (3.2.11) it follows that

$${}^*H(\omega(f(a_1),...,f(a_n)))=g(\omega(h(a_1),...,h(a_n))) \tag{3.2.12}$$

Since $g$ is homomorphism,

$${}^*H(\omega(f(a_1),...,f(a_n)))=\omega(g(h(a_1)),...,g(h(a_n)))=\omega({}^*H(f(a_1)),...,{}^*H(f(a_n)))$$

follows from (3.2.12). Therefore, the map ${}^*H$ is homomorphism of $\Omega_1$-algebra. □

THEOREM 3.2.9. *Given single transitive representation*

$$f:A\rightarrow {}^*M$$

*of $\Omega_1$-algebra $A$ and single transitive representation*

$$g:B\rightarrow {}^*N$$

*of $\Omega_1$-algebra $B$, there exists morphism*

$$h:A\rightarrow B\quad H:M\rightarrow N$$

*of representations from $f$ into $g$.*

Proof. Let us choose homomorphism $h$. Let us choose element $m \in M$ and element $n \in N$. To define map $H$, consider following diagram

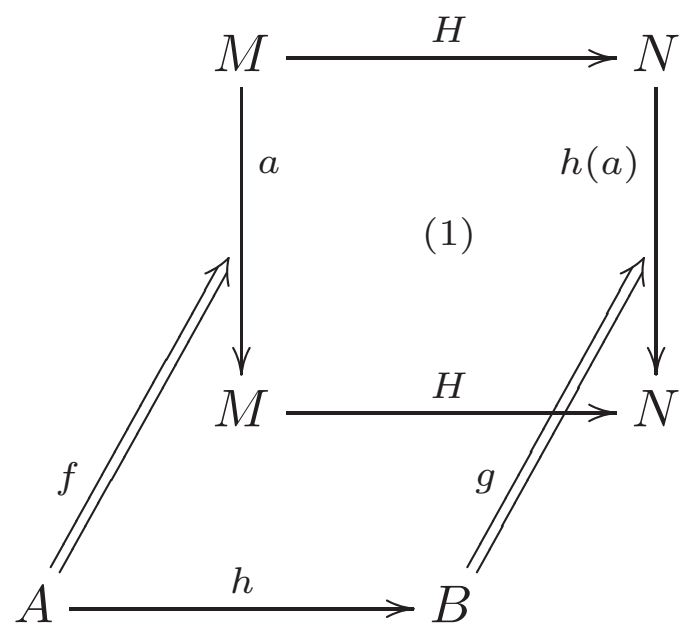

From commutativity of diagram (1), it follows that

$$H(am) = h(a)H(m)$$

For arbitrary $m' \in M$, we defined unambiguously $a \in A$ such that $m' = am$. Therefore, we defined map $H$ which satisfies to equation (3.2.3). □

Theorem 3.2.10. *Let*

$$f : A \to {}^*M$$

*be single transitive representation of* $\Omega_1$*-algebra* $A$ *and*

$$g : B \to {}^*N$$

*be single transitive representation of* $\Omega_1$*-algebra* $B$*. Given homomorphism of* $\Omega_1$*-algebra*

$$h : A \longrightarrow B$$

*consider a homomorphism of* $\Omega_2$*-algebra*

$$H : M \longrightarrow N$$

*such that* $(h, H)$ *is morphism of representations from* $f$ *into* $g$*. This map is unique up to choice of image* $n = H(m) \in N$ *of given element* $m \in M$*.*

Proof. From proof of theorem 3.2.9, it follows that choice of homomorphism $h$ and elements $m \in M$, $n \in N$ uniquely defines the map $H$. □

Theorem 3.2.11. *Given single transitive representation*

$$f : A \to {}^*M$$

*of* $\Omega_1$*-algebra* $A$*, for any endomorphism* $h$ *of* $\Omega_1$*-algebra* $A$ *there exists morphism of representation* $f$

$$h : A \to B \quad H : M \to N$$

PROOF. Consider following diagram

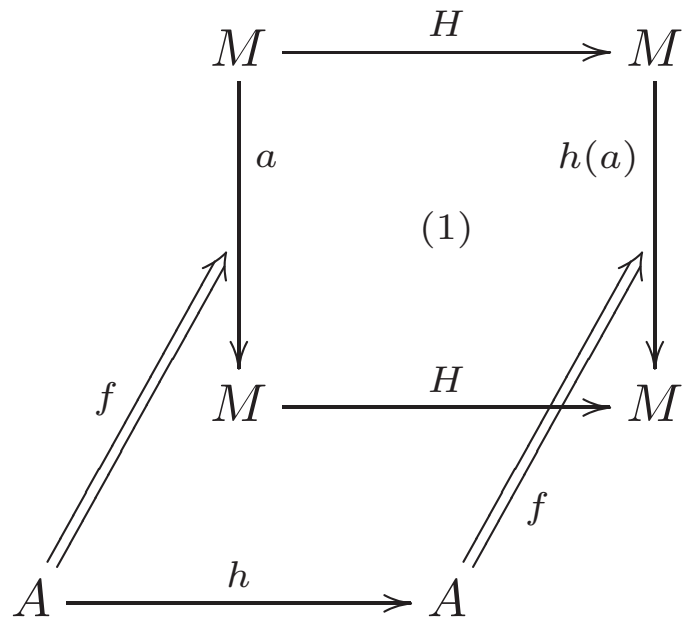

Statement of theorem is corollary of the theorem 3.2.9. □

THEOREM 3.2.12. *Let*

$$f : A \to {}^*M$$

*be representation of $\Omega_1$-algebra $A$,*

$$g : B \to {}^*N$$

*be representation of $\Omega_1$-algebra $B$,*

$$h : C \to {}^*L$$

*be representation of $\Omega_1$-algebra $C$. Given morphisms of representations of $\Omega_1$-algebra*

$$p : A \longrightarrow B \qquad P : M \longrightarrow N$$

$$q : B \longrightarrow C \qquad Q : N \longrightarrow L$$

*There exists morphism of representations of $\Omega_1$-algebra*

$$r : A \longrightarrow C \qquad R : M \longrightarrow L$$

*where $r = qp$, $R = QP$. We call morphism $(r, R)$ of representations from $f$ into $h$* **product of morphisms** $(p, P)$ **and** $(q, Q)$ **of representations of universal algebra**.

PROOF. We represent statement of theorem using diagram

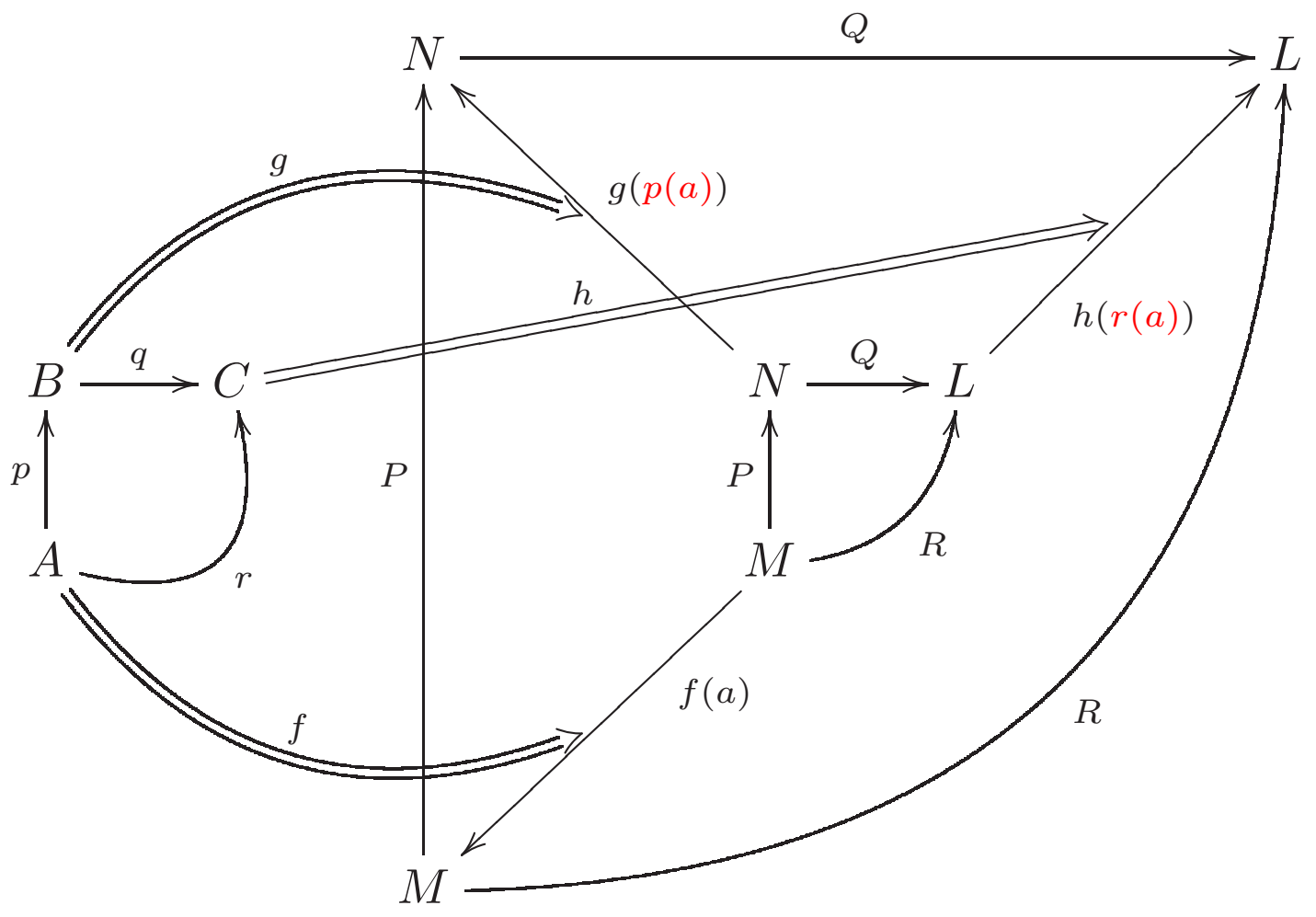

Map $r$ is homomorphism of $\Omega_1$-algebra $A$ into $\Omega_1$-algebra $C$. We need to show that tuple of maps $(r, R)$ satisfies to (3.2.3):

$$\begin{aligned}R(f(a)m) &= QP(f(a)m)\\ &= Q(g(p(a))P(m))\\ &= h(qp(a))QP(m))\\ &= h(r(a))R(m)\end{aligned}$$

□

Definition 3.2.13. Let $\mathcal{A}$ be category of $\Omega_1$-algebras. We define **category $\mathcal{A}*$ of left-side representations of $\Omega_1$-algebra from category $\mathcal{A}$**. Left-side representations of $\Omega_1$-algebra are objects of this category. Morphisms of left-side representations of $\Omega_1$-algebra are morphisms of this category. □

Definition 3.2.14. Let us define equivalence $S$ on the set $M$. Transformation $f$ is called **coordinated with equivalence** $S$, when $f(m_1) \equiv f(m_2)(\mathrm{mod}\ S)$ follows from condition $m_1 \equiv m_2(\mathrm{mod} S)$ . □

Theorem 3.2.15. *Consider equivalence $S$ on set $M$. Consider $\Omega_1$-algebra on set ${}^*M$. If any transformation $f \in {}^*M$ is coordinated with equivalence $S$, then we can define the structure of $\Omega_1$-algebra on the set ${}^*(M/S)$.*

Proof. Let $h = \mathrm{nat}\ S$. If $m_1 \equiv m_2(\mathrm{mod} S)$ , then $h(m_1) = h(m_2)$. Since $f \in {}^*M$ is coordinated with equivalence $S$, then $h(f(m_1)) = h(f(m_2))$. This allows us to define transformation $F$ according to rule

$$F([m]) = h(f(m)) \tag{3.2.13}$$

Let $\omega$ be n-ary operation of $\Omega_1$-algebra. Suppose $f_1$, ..., $f_n \in {}^\star M$ and

$$F_1([m]) = h(f_1(m)) \quad ... \quad F_n([m]) = h(f_n(m))$$

According to condition of theorem, the transformation

$$f = \omega(f_1, ..., f_n) \in {}^*M$$

is coordinated with equivalence $S$. Therefore,

$$\begin{gathered}f(m_1) \equiv f(m_2)(\mathrm{mod} S)\\ \omega(f_1, ..., f_n)(m_1) \equiv \omega(f_1, ..., f_n)(m_2)(\mathrm{mod} S)\end{gathered} \tag{3.2.14}$$

follows from condition $m_1 \equiv m_2(\mathrm{mod} S)$ and the definition 3.2.14. Therefore, we can define operation $\omega$ on the set ${}^\star(M/S)$ according to rule

$$\omega(F_1, ..., F_n)[m] = h(\omega(f_1, ..., f_n)(m)) \tag{3.2.15}$$

From the definition (3.2.13) and equation (3.2.14), it follows that we properly defined operation $\omega$ on the set ${}^\star(M/S)$. □

Theorem 3.2.16. *Let*

$$f : A \to {}^*M$$

*be representation of $\Omega_1$-algebra $A$,*

$$g : B \to {}^*N$$

*be representation of* $\Omega_1$*-algebra* $B$*. Let*

$$r : A \longrightarrow B \qquad R : M \longrightarrow N$$

*be morphism of representations from* $f$ *into* $g$*. Suppose*

$$s = rr^{-1} \qquad S = RR^{-1}$$

*Then there exist decompositions of* $r$ *and* $R$*, which we describe using diagram*

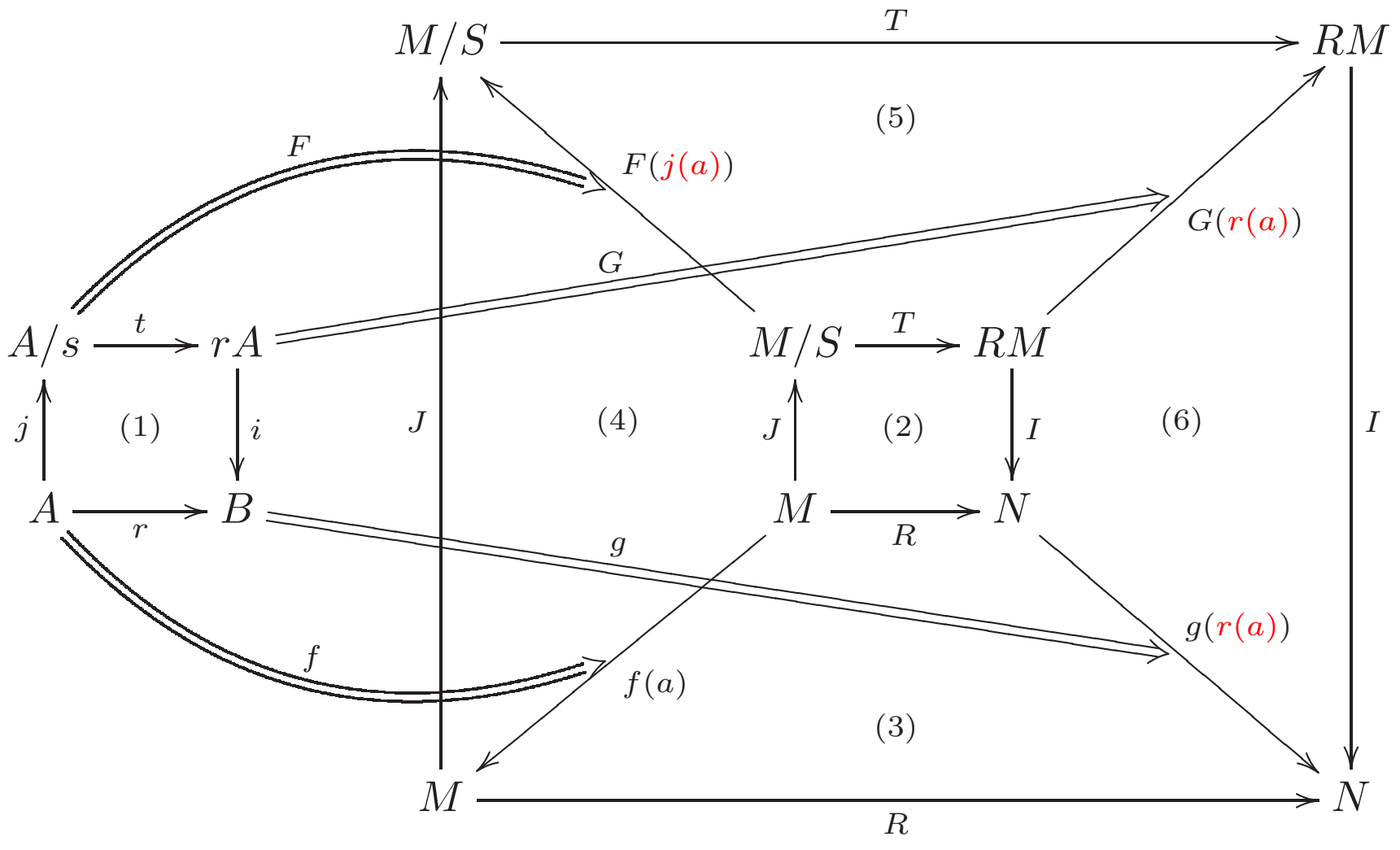

(1) $s = \ker r$ *is a congruence on* $A$*. There exists decompositions of homomorphism* $r$

$$r = itj \tag{3.2.16}$$

$j = \mathrm{nat}\ s$ *is the natural homomorphism*

$$j(a) = j(a)$$

$t$ *is isomorphism*

$$r(a) = t(j(a)) \tag{3.2.17}$$

$i$ *is the inclusion map*

$$r(a) = i(r(a)) \tag{3.2.18}$$

(2) $S = \ker R$ *is an equivalence on* $M$*. There exists decompositions of homomorphism* $R$

$$R = ITJ \tag{3.2.19}$$

$J = \mathrm{nat}\ S$ *is surjection*

$$J(m) = J(m)$$

$T$ *is bijection*

$$R(m) = T(J(m)) \tag{3.2.20}$$

$I$ *is the inclusion map*

$$R(m) = I(R(m)) \tag{3.2.21}$$

(3) $F$ *is left-side representation of* $\Omega_1$*-algebra* $A/s$ *in* $M/S$

(4) *$G$ is left-side representation of $\Omega_1$-algebra $rA$ in $RM$*
(5) *$(j,J)$ is morphism of representations $f$ and $F$*
(6) *$(t,T)$ is morphism of representations $F$ and $G$*
(7) *$(t^{-1},T^{-1})$ is morphism of representations $G$ and $F$*
(8) *$(i,I)$ is morphism of representations $G$ and $g$*
(9) *There exists decompositions of morphism of representations*

$$(r,R)=(i,I)(t,T)(j,J) \tag{3.2.22}$$

Proof. Existence of diagrams (1) and (2) follows from theorem II.3.7 ([13], p. 60).

We start from diagram (4).

Let $m_1\equiv m_2(\text{mod } S)$. Then

$$R(m_1)=R(m_2) \tag{3.2.23}$$

Since $a_1\equiv a_2(\text{mod } s)$, then

$$r(a_1)=r(a_2) \tag{3.2.24}$$

Therefore, $j(a_1)=j(a_2)$. Since $(r,R)$ is morphism of representations, then

$$R(f(a_1)(m_1))=g(r(a_1))(R(m_1)) \tag{3.2.25}$$

$$R(f(a_2)(m_2))=g(r(a_2))(R(m_2)) \tag{3.2.26}$$

From (3.2.23), (3.2.24), (3.2.25), (3.2.26), it follows that

$$R(f(a_1)(m_1))=R(f(a_2)(m_2)) \tag{3.2.27}$$

From (3.2.27) it follows

$$f(a_1)(m_1)\equiv f(a_2)(m_2)(\text{mod } S) \tag{3.2.28}$$

and, therefore,

$$J(f(a_1)(m_1))=J(f(a_2)(m_2)) \tag{3.2.29}$$

From (3.2.29) it follows that map

$$F(j(a))(J(m))=J(f(a)(m))) \tag{3.2.30}$$

is well defined and this map is transformation of set $M/S$.

From equation (3.2.28) (in case $a_1=a_2$) it follows that for any $a$ transformation is coordinated with equivalence $S$. From theorem 3.2.15 it follows that we defined structure of $\Omega_1$-algebra on the set ${}^*(M/S)$. Consider $n$-ary operation $\omega$ and $n$ transformations

$$F(j(a_i))(J(m))=J(f(a_i)(m))) \quad i=1,...,n$$

of the set $M/S$. We assume

$$\omega(F(j(a_1)),...,F(j(a_n)))(J(m))=J(\omega(f(a_1),...,f(a_n)))(m))$$

Therefore, map $F$ is representations of $\Omega_1$-algebra $A/s$.

From (3.2.30) it follows that $(j,J)$ is morphism of representations $f$ and $F$ (the statement (5) of the theorem).

Consider diagram (5).

Since $T$ is bijection, then we identify elements of the set $M/S$ and the set $MR$, and this identification has form

$$T(J(m)) = R(m) \tag{3.2.31}$$

We can write transformation $F(j(a))$ of the set $M/S$ as

$$F(j(a)) : J(m) \to F(j(a))(J(m)) \tag{3.2.32}$$

Since $T$ is bijection, we define transformation

$$T(J(m)) \to T(F(j(a))(J(m))) \tag{3.2.33}$$

of the set $RM$. Transformation (3.2.33) depends on $j(a) \in A/s$. Since $t$ is bijection, we identify elements of the set $A/s$ and the set $rA$, and this identification has form

$$t(j(a)) = r(a) \tag{3.2.34}$$

Therefore, we defined map

$$G : rA \to {}^\star RM$$

according to equation

$$G(t(j(a)))(T(J(m))) = T(F(j(a))(J(m))) \tag{3.2.35}$$

Consider $n$-ary operation $\omega$ and $n$ transformations

$$G(r(a_i))(R(m)) = T(F(j(a_i))(J(m))) \quad i = 1, ..., n$$

of space $RM$. We assume

$$\omega(G(r(a_1)), ..., G(r(a_n)))(R(m)) = T(\omega(F(j(a_1), ..., F(j(a_n)))(J(m))) \tag{3.2.36}$$

According to (3.2.35) operation $\omega$ is well defined on the set ${}^\star RM$. Therefore, the map $G$ is representations of $\Omega_1$-algebra.

From (3.2.35), it follows that $(t, T)$ is morphism of representations $F$ and $G$ (the statement (6) of the theorem).

Since $T$ is bijection, then from equation (3.2.31), it follows that

$$J(m) = T^{-1}(R(m)) \tag{3.2.37}$$

We can write transformation $G(r(a))$ of the set $RM$ as

$$G(r(a)) : R(m) \to G(r(a))(R(m)) \tag{3.2.38}$$

Since $T$ is bijection, we define transformation

$$T^{-1}(R(m)) \to T^{-1}(G(r(a))(R(m))) \tag{3.2.39}$$

of the set $M/S$. Transformation (3.2.39) depends on $r(a) \in rA$. Since $t$ is bijection, then from equation (3.2.34) it follows that

$$j(a) = t^{-1}(r(a)) \tag{3.2.40}$$

Since, by construction, diagram (5) is commutative, then transformation (3.2.39) coincides with transformation (3.2.32). We can write the equation (3.2.36) as

$$T^{-1}(\omega(G(r(a_1)), ..., G(r(a_n)))(R(m))) = \omega(F(j(a_1), ..., F(j(a_n)))(J(m)) \tag{3.2.41}$$

Therefore $(t^{-1}, T^{-1})$ is morphism of representations $G$ and $F$ (the statement (7) of the theorem).

Diagram (6) is the simplest case in our proof. Since map $I$ is immersion and diagram (2) is commutative, we identify $n \in N$ and $R(m)$ when $n \in \mathrm{Im}R$. Similarly, we identify corresponding transformations.

$$g'(i(r(a)))(I(R(m))) = I(G(r(a))(R(m))) \tag{3.2.42}$$
$$\omega(g'(r(a_1)), ..., g'(r(a_n)))(R(m)) = I(\omega(G(r(a_1), ..., G(r(a_n)))(R(m)))$$

Therefore, $(i, I)$ is morphism of representations $G$ and $g$ (the statement (8) of the theorem).

To prove the statement (9) of the theorem we need to show that defined in the proof representation $g'$ is congruent with representation $g$, and operations over transformations are congruent with corresponding operations over ${}^*N$.

$$\begin{aligned}
g'(i(r(a)))(I(R(m))) &= I(G(r(a))(R(m))) && \text{by (3.2.42)}\\
&= I(G(t(j(a)))(T(J(m)))) && \text{by (3.2.17), (3.2.20)}\\
&= IT(F(j(a))(J(m))) && \text{by (3.2.35)}\\
&= ITJ(f(a)(m)) && \text{by (3.2.30)}\\
&= R(f(a)(m)) && \text{by (3.2.19)}\\
&= g(r(a))(R(m)) && \text{by (3.2.3)}
\end{aligned}$$

$$\begin{aligned}
\omega(G(r(a_1)), ..., G(r(a_n)))(R(m)) &= T(\omega(F(j(a_1), ..., F(j(a_n)))(J(m)))\\
&= T(F(\omega(j(a_1), ..., j(a_n)))(J(m)))\\
&= T(F(j(\omega(a_1, ..., a_n)))(J(m)))\\
&= T(J(f(\omega(a_1, ..., a_n))(m)))
\end{aligned}$$

□

Definition 3.2.17. Let

$$f : A \to {}^*M$$

be representation of $\Omega_1$-algebra $A$,

$$g : B \to {}^*N$$

be representation of $\Omega_1$-algebra $B$. Let

$$r : A \longrightarrow B \qquad R : M \longrightarrow N$$

be morphism of representations from $f$ into $g$ such that $f$ is isomorphism of $\Omega_1$-algebra and $g$ is isomorphism of $\Omega_2$-algebra. Then map $(r, R)$ is called **isomorphism of repesentations**. □

Theorem 3.2.18. *In the decomposition (3.2.22), the map $(t, T)$ is isomorphism of representations $F$ and $G$.*

Proof. The statement of the theorem is corollary of definition 3.2.17 and statements (6) and (7) of the theorem 3.2.16. □

From theorem 3.2.16 it follows that we can reduce the problem of studying of morphism of representations of $\Omega_1$-algebra to the case described by diagram

(3.2.43)

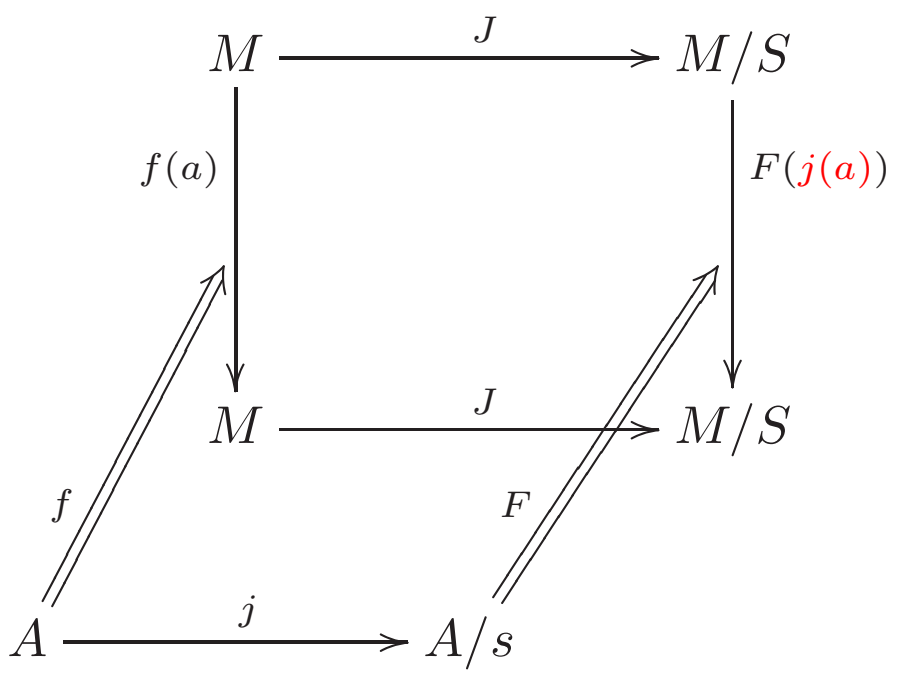

THEOREM 3.2.19. *We can supplement diagram* (3.2.43) *with representation* $F_1$ *of* $\Omega_1$*-algebra* $A$ *into set* $M/S$ *such that diagram*

(3.2.44)

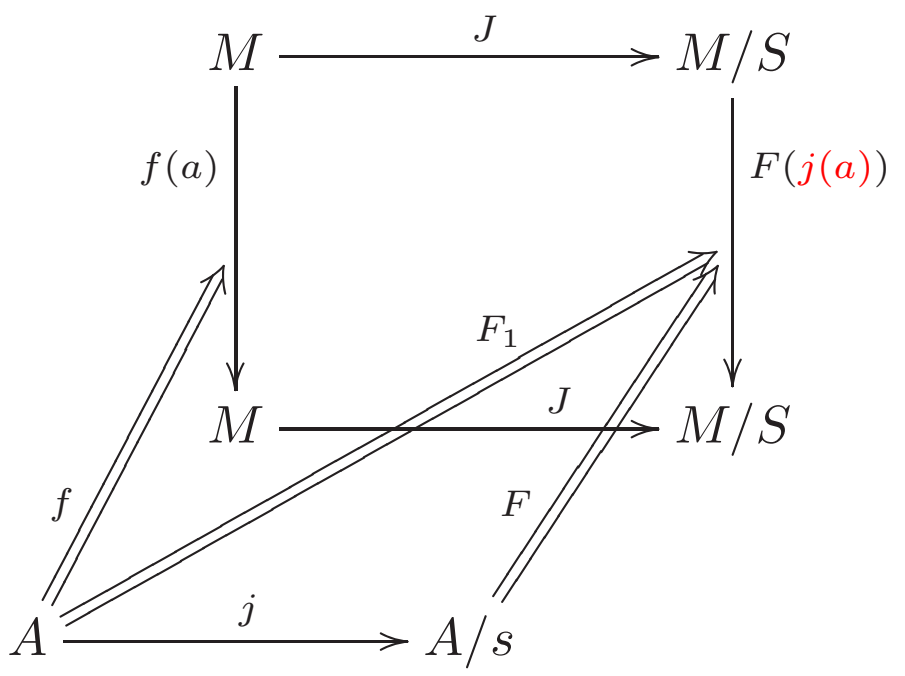

*is commutative. The set of transformations of representation* $F$ *and the set of transformations of representation* $F_1$ *coincide.*

PROOF. To prove theorem it is enough to assume

$$F_1(a) = F(j(a))$$

Since map $j$ is surjection, then $\mathrm{Im}F_1 = \mathrm{Im}F$. Since $j$ and $F$ are homomorphisms of $\Omega_1$-algebra, then $F_1$ is also homomorphism of $\Omega_1$-algebra. □

Theorem 3.2.19 completes the series of theorems dedicated to the structure of morphism of representations $\Omega_1$-algebra. From these theorems it follows that we can simplify task of studying of morphism of representations $\Omega_1$-algebra and not go beyond morphism of representations of form

$$id : A \longrightarrow A \qquad R : M \longrightarrow N$$

In this case we identify morphism of $(\mathrm{id}, R)$ representations of $\Omega_1$-algebra and corresponding homomorphism $R$ of $\Omega_2$-algebra and use the same letter $R$ to denote these maps. We will use diagram

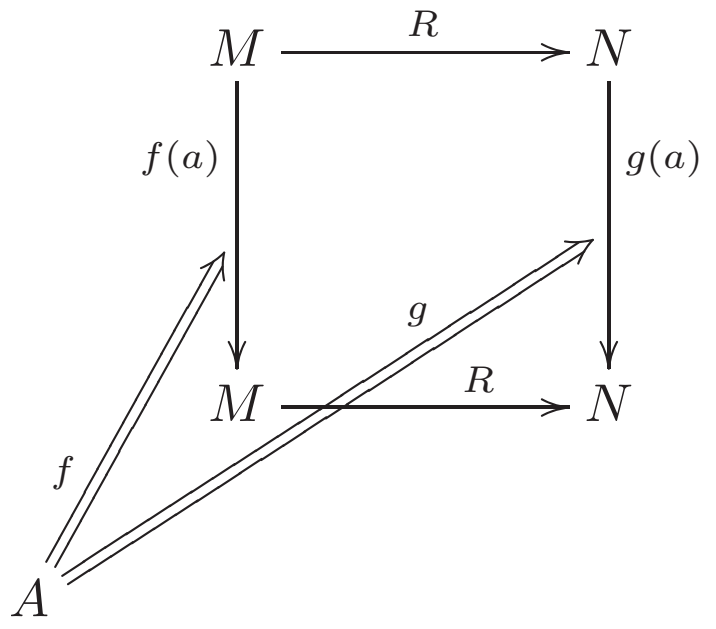

to represent morphism $(\mathrm{id}, R)$ of representations of $\Omega_1$-algebra. From diagram it follows

$$R \circ f(a) = g(a) \circ R \tag{3.2.45}$$

By analogy with definition 3.2.13. we give following definition.

DEFINITION 3.2.20. We define **category** $A*$ **of left-side representations of** $\Omega_1$**-algebra** $A$. Left-side representations of $\Omega_1$-algebra $A$ are objects of this category. Morphisms $(id, R)$ of left-side representations of $\Omega_1$-algebra $A$ are morphisms of this category. □

## 3.3. Automorphism of Representation of Universal Algebra

DEFINITION 3.3.1. Let

$$f : A \to {}^*M$$

be representation of $\Omega_1$-algebra $A$ in $\Omega_2$-algebra $M$. The morphism of representations of $\Omega_1$-algebra

$$(\mathrm{id} : A \to A, R : M \to M)$$

such, that $R$ is endomorphism of $\Omega_2$-algebra is called **endomorphism of representation** $f$. □

THEOREM 3.3.2. *Given single transitive representation*

$$f : A \to {}^*M$$

*of* $\Omega_1$*-algebra* $A$*, for any* $p$*,* $q \in M$ *there exists unique endomorphism*

$$H : M \to M$$

*of representation* $f$ *such that* $H(p) = q$*.*

PROOF. Consider following diagram

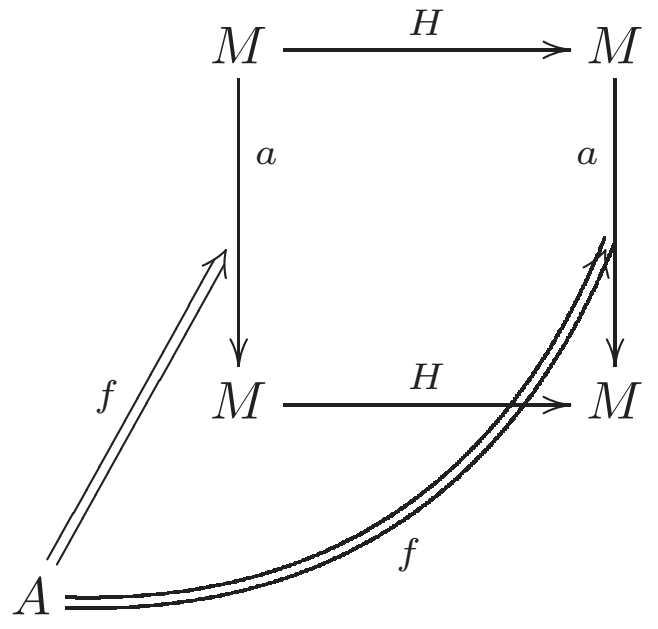

Existence of endomorphism is corollary of the theorem 3.2.9. For given $p, q \in M$, uniqueness of endomorphism follows from the theorem 3.2.10 when $h = \mathrm{id}$. □

THEOREM 3.3.3. *Endomorphisms of representation $f$ form semigroup.*

PROOF. From theorem 3.2.12, it follows that the product of endomorphisms $(p, P)$, $(r, R)$ of the representation $f$ is endomorphism $(pr, PR)$ of the representation $f$. □

DEFINITION 3.3.4. Let

$$f : A \to {}^*M$$

be representation of $\Omega_1$-algebra $A$ in $\Omega_2$-algebra $M$. The morphism of representations of $\Omega_1$-algebra

$$(\mathrm{id} : A \to A, R : M \to M)$$

such, that $R$ is automorphism of $\Omega_2$-algebra is called **automorphism of representation** $f$. □

THEOREM 3.3.5. *Let*

$$f : A \to {}^*M$$

*be representation of $\Omega_1$-algebra $A$ in $\Omega_2$-algebra $M$. The set of automorphisms of the representation $f$ forms* **group** $GA(f)$ .

PROOF. Let $R$, $P$ be automorphisms of the representation $f$. According to definition 3.3.4, maps $R$, $P$ are automorphisms of $\Omega_2$-algebra $M$. According to theorem II.3.2, ([13], p. 57), the map $R \circ P$ is automorphism of $\Omega_2$-algebra $M$. From the theorem 3.2.12 and the definition 3.3.4, it follows that product of automorphisms $R \circ P$ of the representation $f$ is automorphism of the representation $f$.

Let $R$, $P$, $Q$ be automorphisms of the representation $f$. The associativity of product of maps $R$, $P$, $Q$ follows from the chain of equations[3.3]

$$\begin{aligned}((R \circ P) \circ Q)(a) &= (R \circ P)(Q(a)) = R(P(Q(a))) \\ &= R((P \circ Q)(a)) = (R \circ (P \circ Q))(a)\end{aligned}$$

Let $R$ be an automorphism of the representation $f$. According to definition 3.3.4 the map $R$ is automorphism of $\Omega_2$-algebra $M$. Therefore, the map $R^{-1}$ is automorphism of $\Omega_2$-algebra $M$. The equation (3.2.4) is true for automorphism $R$ of representation. Assume $m' = R(m)$. Since $R$ is automorphism of $\Omega_2$-algebra, then $m = R^{-1}(m')$ and we can write (3.2.4) in the form

$$R(f(a')(R^{-1}(m'))) = f(a')(m') \tag{3.3.1}$$

Since the map $R$ is automorphism of $\Omega_2$-algebra $M$, then from the equation (3.3.1) it follows that

$$f(a')(R^{-1}(m')) = R^{-1}(f(a')(m')) \tag{3.3.2}$$

The equation (3.3.2) corresponds to the equation (3.2.4) for the map $R^{-1}$. Therefore, map $R^{-1}$ of the representation $f$. □

[3.3] To prove the associativity of product I follow to the example of the semigroup from [4], p. 20, 21.

CHAPTER 4

# Representation of Group

## 4.1. Representation of Group

Group is among few algebras that allow somebody to consider the product of transformations of the $\Omega$-algebra $M$ in such a way that, if transformations belong to the representation, then their product also belongs to the representation. We should remember that order of maps in product depends on order of maps on diagram and how these maps act over elements of the set (from left or from right).

DEFINITION 4.1.1. Let ${}^*M$ be a group with product

$$(f \circ g)x = f(gx)$$

and $\delta$ be unit of group ${}^\star M$. Let $G$ be group. We call a homomorphism of group

$$f : G \to {}^\star M \tag{4.1.1}$$

**left-side representation of group** $G$ or $G*$**-representation** in $\Omega$-algebra $M$ if the map $f$ holds

$$f(ab)u = f(a)(f(b)u) \tag{4.1.2}$$

□

REMARK 4.1.2. Since the map (4.1.1) is homomorhism, then

$$f(ab)u = (f(a)f(b))u \tag{4.1.3}$$

We use here convention

$$f(a)f(b) = f(a) \circ f(b)$$

Thus, the idea of representation of group is that we multiply elements of group in the same order as we multiply transformations of representation. From equations (4.1.2) and (4.1.3) it follows

$$(f(a)f(b))u = f(a)(f(b)u) \tag{4.1.4}$$

Equation (4.1.4) together with associativity of product of transformations expresses **associative law** for $G*$-representation. This allows writing of equation (4.1.4) without using of brackets

$$f(ab)u = f(a)f(b)u$$

□

DEFINITION 4.1.3. Let $M^\star$ be a group with product

$$x(f \circ g) = (xf)g$$

and $\delta$ be unit of group $M^\star$. Let $G$ be group. We call a homomorphism of group

$$f : G \to M^* \tag{4.1.5}$$

**right-side representation of group** $G$ or **$*G$-representation** in $\Omega$-algebra $M$ if the map $f$ holds

$$uf(ab) = (uf(a))f(b) \tag{4.1.6}$$

□

REMARK 4.1.4. Since the map (4.1.5) is homomorhism, then

$$uf(ab) = u(f(a)f(b)) \tag{4.1.7}$$

From equations (4.1.6) and (4.1.7) it follows

$$u(f(a)f(b)) = (uf(a))f(b) \tag{4.1.8}$$

Equation (4.1.8) together with associativity of product of transformations expresses **associative law** for $*G$-representation. This allows writing of equation (4.1.8) without using of brackets

$$uf(ab) = uf(a)f(b)$$

□

DEFINITION 4.1.5. We call the transformation

$$t : M \to M$$

**nonsingular transformation**, if there exists inverse map. □

THEOREM 4.1.6. *For any $g \in G$ transformation is nonsingular and satisfies equation*

$$f(g^{-1}) = f(g)^{-1} \tag{4.1.9}$$

PROOF. Since (4.1.2) and

$$f(e) = \delta$$

we have

$$u = \delta(u) = f(gg^{-1})(u) = f(g)(f(g^{-1})(u))$$

This completes the proof. □

THEOREM 4.1.7. *The group operation determines two different representations on the group:*

- *The* **left shift** $t_\star$

$$\begin{aligned} b' &= t_\star(a)b = ab \\ b' &= t_\star(a)(b) = ab \end{aligned} \tag{4.1.10}$$

*is $G*$-representation on the set*[4.1] $G$

$$t_\star(ab) = t_\star(a) \circ t_\star(b) \tag{4.1.11}$$

- *The* **right shift** ${}_\star t$

$$\begin{aligned} b' &= b\ {}_\star t(a) = ba \\ b' &= {}_\star t(a)(b) = ba \end{aligned} \tag{4.1.12}$$

*is $*G$-representation on the set $G$*

$${}_\star t(ab) = {}_\star t(a) \circ {}_\star t(b) \tag{4.1.13}$$

[4.1] Left shift is not a representation of group in a group, because the transformation $t_\star$ is not a homomorphism of group. Similar remark is true for right shift.

Proof. Equation (4.1.11) follows from associativity of product

$$t_\star(ab)c = (ab)c = a(bc) = t_\star(a)(t_\star(b)c) = (t_\star(a) \circ t_\star(b))c$$

In a similar manner we prove the equation (4.1.13). □

Definition 4.1.8. Let $G$ be group. Let $f$ be $G*$-representation in $\Omega$-algebra $M$. For any $v \in M$ we define **orbit of representation of the group** $G$ as set

$$f(G)v = \{w = f(g)v : g \in G\}$$

□

Since $f(e) = \delta$ we have $v \in f(G)v$.

Theorem 4.1.9. *Suppose*

$$v \in f(G)u \tag{4.1.14}$$

*Then*

$$f(G)u = f(G)v$$

Proof. From (4.1.14) it follows that there exists $a \in G$ such that

$$v = f(a)u \tag{4.1.15}$$

Suppose $w \in f(G)v$. Then there exists $b \in G$ such that

$$w = f(b)v \tag{4.1.16}$$

If we substitute (4.1.15) into (4.1.16) we get

$$w = f(b)(f(a)u) \tag{4.1.17}$$

Since (4.1.2), we see that from (4.1.17) it follows that $w \in f(G)u$. Thus

$$f(G)v \subseteq f(G)u$$

Since (4.1.9), we see that from (4.1.15) it follows that

$$u = f(a)^{-1}v = f(a^{-1})v \tag{4.1.18}$$

From (4.1.18) it follows that $u \in f(G)v$ and therefore

$$f(G)u \subseteq f(G)v$$

This completes the proof. □

Thus, $G*$-representation $f$ in $\Omega$-algebra $M$ forms equivalence $S$ and the orbit $f(G)u$ is equivalence class. We will use notation $M/f(G)$ for quotient set $M/S$ and this set is called **space of orbits of $G*$-representation** $f$.

Theorem 4.1.10. *Suppose $f_1$ is $G*$-representation in $\Omega$-algebra $M_1$ and $f_2$ is $G*$-representation in $\Omega$-algebra $M_2$. Then we introduce* **direct product of $G*$-representations** *$f_1$ and $f_2$*

$$f = f_1 \times f_2 : G \to M_1 \otimes M_2$$
$$f(g) = (f_1(g), f_2(g))$$

Proof. To show that $f$ is a representation, it is enough to prove that $f$ satisfies the definition 4.1.1.

$$f(e) = (f_1(e), f_2(e)) = (\delta_1, \delta_2) = \delta$$

$$\begin{aligned} f(ab)u &= (f_1(ab)u_1, f_2(ab)u_2) \\ &= (f_1(a)(f_1(b)u_1), f_2(a)(f_2(b)u_2)) \\ &= f(a)(f_1(b)u_1, f_2(b)u_2) \\ &= f(a)(f(b)u) \end{aligned}$$

□

## 4.2. Single Transitive Right-Side Representation of Group

Definition 4.2.1. We call **kernel of inefficiency of $G*$-representation** a set

$$K_f = \{g \in G : f(g) = \delta\}$$

□

Theorem 4.2.2. *A kernel of inefficiency of $G*$-representation is a subgroup of the group $G$.*

Proof. Assume $f(a_1) = \delta$ and $f(a_2) = \delta$. Then

$$f(a_1a_2)u = f(a_1)(f(a_2)u) = u$$

$$f(a^{-1}) = f^{-1}(a) = \delta$$

□

Theorem 4.2.3. *$G*$-representation is* **effective** *iff kernel of inefficiency $K_f = \{e\}$.*

Proof. Statement is corollary of definitions 3.1.6 and 4.2.1 and of the theorem 4.2.2. □

If an action is not effective we can switch to an effective one by changing group $G_1 = G|K_f$ using factorization by the kernel of inefficiency. This means that we can study only an effective action.

Definition 4.2.4. Consider $G*$-representation $f$ in $\Omega$-algebra $M$. A **little group** or **stability group** of $x \in M$ is the set

$$= \{g \in G : f(g)x = x\}$$

$G*$-representation $f$ is said to be **free**, if for any $x \in M$ stability group $G_x = \{e\}$. □

Theorem 4.2.5. *Given free $G*$-representation $f$ in the $\Omega$-algebra $A$, there exist $1-1$ correspondence between orbits of representation, as well between orbit of representation and group $G$.*

Proof. Given $a \in A$ there exist $g_1$, $g_2 \in G$

$$f(g_1)a = f(g_2)a \tag{4.2.1}$$

We multiply both parts of equation (4.2.1) by $f(g_1^{-1})$

$$a = f(g_1^{-1})f(g_2)a$$

Since the representation is free, $g_1 = g_2$. Since we established $1-1$ correspondence between orbit and group $G$, we proved the statement of the theorem. □

DEFINITION 4.2.6. We call a space $V$ **homogeneous space of group** $G$ if we have single transitive $G*$-representation on $V$. □

THEOREM 4.2.7. *If we define a single transitive representation $f$ of the group $G$ on the $\Omega$-algebra $A$ then we can uniquely define coordinates on $A$ using coordinates on the group $G$.*

*If $f$ is left-side representation then $f(a)$ is equivalent to the left shift $t_\star(a)$ on the group $G$. If $f$ is right-side representation then $f(a)$ is equivalent to the right shift ${}_\star t(a)$ on the group $G$.*

PROOF. We select a point $v \in A$ and define coordinates of a point $w \in A$ as coordinates of $a \in G$ such that $w = f(a)v$. Coordinates defined this way are unique up to choice of an initial point $v \in A$ because the action is effective.

If $f$ is left-side representation, we will use the notation

$$f(a)v = av$$

Because the notation

$$f(a)(f(b)v) = a(bv) = (ab)v = f(ab)v$$

is compatible with the group structure we see that left-side representation $f$ is equivalent to the left shift.

If $f$ is right-side representation, we will use the notation

$$vf(a) = va$$

Because the notation

$$(vf(b))f(a) = (vb)a = v(ba) = vf(ba)$$

is compatible with the group structure we see that right-side representation $f$ is equivalent to the right shift. □

REMARK 4.2.8. We will write effective $G*$-representation as

$$v' = t_\star(a)v = av$$

Orbit of this representation is

$$Gv = t_\star(G)v$$

We will use notation $M/t_\star(G)$ for the space of orbits of effective $G*$-representation. □

REMARK 4.2.9. We will write effective $*G$-representation as

$$v' = v\ {}_\star t(a) = va$$

Orbit of this representation is

$$vG = v\ {}_\star t(G)$$

We will use notation $M/{}_\star t(G)$ for the space of orbits of effective $*G$-representation. □

THEOREM 4.2.10. *Free $G**$-representation is effective. Free $G**$-representation $f$ in $\Omega$-algebra $M$ is single transitive representation on orbit.*

PROOF. The statement of theorem is the corollary of definition 4.2.4. □

Theorem 4.2.11. *Left and right shifts on group $G$ are commuting.*

Proof. This is the consequence of the associativity on the group $G$

$$(t_\star(a)\circ\ {}_\star t(b))c = a(cb) = (ac)b = ({}_\star t(b)\circ t_\star(a))c$$

□

Theorem 4.2.11 can be phrased n the following way.

Theorem 4.2.12. *Let $G$ be group. For any $a\in G$, the map $t_\star(a)$ is automorphism of representation ${}_\star t$.*

Proof. According to theorem 4.2.11

$$t_\star(a)\circ {}_\star t(b) = {}_\star t(b)\circ t_\star(a) \tag{4.2.2}$$

Equation (4.2.2) coincides with equation (3.2.3) from definition 3.2.2 when $r=id$, $R=t_\star(a)$. □

Theorem 4.2.13. *Let $G*$-representation $f$ on $\Omega$-algebra $M$ be single transitive. Then we can uniquely define a single transitive $*G$-representation $h$ on $\Omega$-algebra $M$ such that diagram*

$$\begin{array}{ccc} M & \xrightarrow{\ \ h(a)\ \ } & M \\ \Big\downarrow{\scriptstyle f(b)} & & \Big\downarrow{\scriptstyle f(b)} \\ M & \xrightarrow[\ \ h(a)\ \ ]{} & M \end{array}$$

*is commutative for any $a$, $b\in G$.*[4.2]

Proof. We use group coordinates for points $v\in M$. Then according to theorem 4.2.7 we can write the left shift $t_\star(a)$ instead of the transformation $f(a)$.

Let $v_0, v\in M$. Then we can find one and only one $a\in G$ such that

$$v = v_0a = v_0\ {}_\star t(a)$$

We assume

$$h(a) =\ {}_\star t(a)$$

For some $b\in G$ we have

$$w_0 = f(b)v_0 = t_\star(b)v_0 \qquad w = f(b)v = t_\star(b)v$$

According to the theorem 4.2.11, the diagram

$$\begin{array}{ccc} v_0 & \xrightarrow{\ h(a)=\ {}_\star t(a)\ } & v \\ \Big\downarrow{\scriptstyle f(b)=t_\star(b)} & & \Big\downarrow{\scriptstyle f(b)=t_\star(b)} \\ w_0 & \xrightarrow[\ h(a)=\ {}_\star t(a)\ ]{} & w \end{array} \tag{4.2.3}$$

is commutative.

Changing $b$ we get that $w_0$ is an arbitrary point of $M$.

We see from the diagram that if $v_0=v$ then $w_0=w$ and therefore $h(e)=\delta$. On other hand if $v_0\ne v$ then $w_0\ne w$ because the $G*$-representation $f$ is single transitive. Therefore the $*G$-representation $h$ is effective.

[4.2] You can see this statement in [3].

In the same way we can show that for given $w_0$ we can find $a$ such that $w = h(a)w_0$. Therefore the $*G$-representation $h$ is single transitive.

In general the product of transformations of the $G*$-representation $f$ is not commutative and therefore the $*G$-representation $h$ is different from the $G*$-representation $f$. In the same way we can create a $G*$-representation $f$ using the $*G$-representation $h$. □

Representations $f$ and $h$ are called **twin representations of the group** $G$.

REMARK 4.2.14. It is clear that transformations $t_\star(a)$ and ${}_\star t(a)$ are different until the group $G$ is nonabelian. However they both are maps onto. Theorem 4.2.13 states that if both right and left shift presentations exist on the set $M$, then we can define two commuting representations on the set $M$. The right shift or the left shift only cannot represent both types of representation. To understand why it is so let us change diagram (4.2.3) and assume $h(a)v_0 = t_\star(a)v_0 = v$ instead of $h(a)v_0 = v_{0\star}t(a) = v$ and let us see what expression $h(a)$ has at the point $w_0$. The diagram

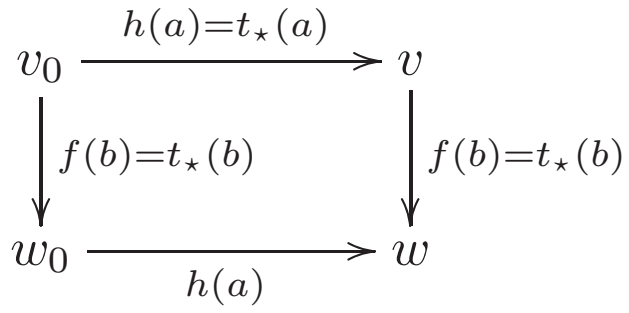

is equivalent to the diagram

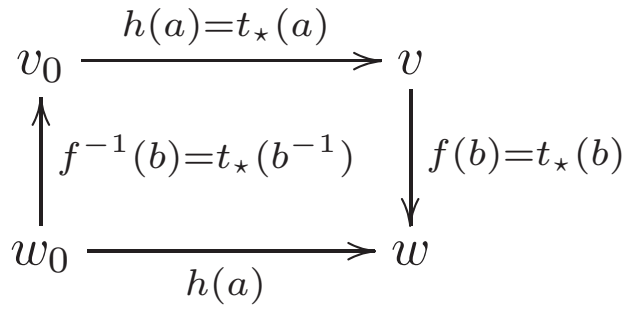

and we have $w = bv = bav_0 = bab^{-1}w_0$. Therefore

$$h(a)w_0 = (bab^{-1})w_0$$

We see that the representation of $h$ depends on its argument. □

THEOREM 4.2.15. *Let $f$ and $h$ be twin representations of the group $G$. For any $a \in G$ the map $h(a)$ is automorphism of representation $f$.*

PROOF. The statement of theorem is corollary of theorems 4.2.12 and 4.2.13. □

REMARK 4.2.16. Is there a morphism of representations from $t_\star$ to $t_\star$ different from automorphism $(\mathrm{id}, {}_\star t(a))$? If we assume

$$r(g) = cgc^{-1}$$

$$R(a)(m) = cmac^{-1}$$

then it is easy to see that the map $(r, R(a))$ is morphism of the representations from $t_\star$ to $t_\star$. However this map is not automorphism of the representation $t_\star$, because $r \neq \mathrm{id}$. □

CHAPTER 5

# Vector Space over Division Ring

## 5.1. Vector Space

To define left-side representation

$$f : D \xrightarrow{*} M \qquad f(d) : v \to d\,v$$

of ring $D$ in the $\Omega$-algebra $M$ we need to define the structure of the ring on the set ${}^{\star}M$. [5.1]

THEOREM 5.1.1. *Left-side representation $f$ of the ring $D$ in the $\Omega$-algebra $M$ is defined iff left-side representations of multiplicative and additive groups of the ring $D$ are defined and these representations hold relationship*

$$f(a(b+c)) = f(a)f(b) + f(a)f(c)$$

PROOF. Theorem follows from definition 3.1.2. □

DEFINITION 5.1.2. An Abelian group $M$ is a $D*$**-module** if there exists $D*$-representation

$$f : D \xrightarrow{*} M \qquad f(d) : v \to d\,v$$

□

According to our notation $D*$-module is **left module over a ring** $D$ and $*D$**-module** is **right module over a ring** $D$.

---

[5.1] Is it possible to define an addition on the set ${}^{\star}M$, if this operation is not defined on the set $M$. The answer on this question is positive.

Let $M = B \cup C$ and let $F : B \to C$ be one to one map. We define the set ${}^{\star}M$ of left-side transformations of the set $M$ according to the following rule. Let $V \subseteq B$. Let the left-side transformation $F_V$ be given by

$$F_V x = \begin{cases} x & x \in B \backslash V \\ Fx & x \in V \\ x & x \in C \backslash F(V) \\ F^{-1}x & x \in F(V) \end{cases}$$

We define sum of left-side transformations according rule

$$F_V + F_W = F_{V \triangle W}$$
$$V \triangle W = (V \cup W) \backslash (V \cap W)$$

It is evident that

$$F_\emptyset + F_V = F_V$$
$$F_V + F_V = F_\emptyset$$

Therefore, the map $F_\emptyset$ is zero of the addition, and the set ${}^{\star}M$ is the Abelian group.

Since a field is a special case of a ring, vector space over the field has more properties then module over the ring. It is very hard, if possible at all, to extend definitions, which work in a vector space, to a module over an arbitrary ring. A definition of a basis and dimension of vector space are closely linked with the possibility of finding a solution of a linear equation in a ring. Properties of the linear equation in division ring are close to properties of the linear equation in field. This is why we hope that properties of vector space over division ring are close to properties of the vector space over the field.

THEOREM 5.1.3. *Left-side representation of the division ring $D$ is* **effective** *iff left-side representation of its multiplicative group is effective.*

PROOF. Suppose

$$f : D \xrightarrow{*} M \qquad f(d) : v \to d\,v$$

is left-side representation of the division ring $D$. Suppose elements $a$, $b$ of of the multiplicative group cause the same left-side transformation. Then

$$f(a)m = f(b)m \tag{5.1.1}$$

for any $m \in M$. Performing transformation $f(a^{-1})$ on both sides of the equation (5.1.1), we obtain

$$m = f(a^{-1})(f(b)m) = f(a^{-1}b)m$$

□

According to the remark 3.1.7, since the representation of the division ring is effective, we identify an element of the division ring and left-side transformation corresponding to this element.

DEFINITION 5.1.4. Let $D$ be division ring. Abelian group $V$ is a **$D*$-vector space** if there exists effective $D*$-representation

$$f : D \xrightarrow{*} M \qquad f(d) : v \to d\,v \tag{5.1.2}$$

Abelian group $V$ is a **$*D$-vector space** if there exists effective $*D$-representation

$$f : D \xrightarrow{*} M \qquad f(d) : v \to v\,d \tag{5.1.3}$$

□

$D*$-vector space is also called **left $D$-vector space** or left vector space over a division ring $D$. $*D$-vector space is also called **right $D$-vector space** or right vector space over a division ring $D$.

THEOREM 5.1.5. *Following conditions hold for $D*$-vector space:*

- **associative law**

$$(ab)m = a(bm) \tag{5.1.4}$$

- **distributive law**

$$a(m+n) = am + an \tag{5.1.5}$$

$$(a+b)m = am + bm \tag{5.1.6}$$

- **unitarity law**

$$1m = m \tag{5.1.7}$$

*for any* $a,\ b \in D,\ m,\ n \in V$.

Proof. Since left-side transformation $a$ is endomorphism of the Abelian group, we obtain the equation (5.1.5). Since representation (5.1.2) is homomorphism of the aditive group of division ring $D$, we obtain the equation (5.1.6). Since representation (5.1.2) is left-side representation of the multiplicative group of division ring $D$, we obtain the equations (5.1.4) and (5.1.7). □

According to our notation $D*$-vector space is **left $D$-vector space**. Map

$$(d, v) \in D \times V \to dv \in V$$

generated by $D*$-representation (5.1.2) is called **left-side product of vector over scalar**.

According to our notation **$*D$-vector space** is **right $D$-vector space**. Map

$$(v, d) \in D \times V \to vd \in V$$

generated by $*D$-representation (5.1.3), is called **right-side product of vector over scalar**.

Any statement that is valid for left $D$-vector space, is valid for right $D$-vector space, if we substitude left-side product of vector over scalar by right-side product of vector over scalar.

Definition 5.1.6. Let $V$ be a $D*$-vector space over a division ring $D$. Set of vectors $N$ is a **subspace of $D*$-vector space** $V$ if

$$a + b \in N$$

$$ka \in N$$

$$a, b \in N \quad k \in D$$

□

Example 5.1.7. Let $D^m_n$ be set of $m \times n$ matrices over division ring $D$. We define addition

$$a + b = \left(a^j_i\right) + \left(b^j_i\right) = \left(a^j_i + b^j_i\right)$$

and product over scalar

$$da = d\left(a^j_i\right) = \left(da^j_i\right)$$

$a = 0$ iff $a^j_i = 0$ for any $i$, $j$. We can verify directly that $D^m_n$ is a $D*$-vector space. when product is defined from left. Otherwise $D^m_n$ is $*D$-vector space. Vector space $D^m_n$ is called **$D*$-matrices vector space**. □

## 5.2. Vector Space Type

The product of vector over scalar is asymmetric because the product is defined for objects of different sets. However we see difference between $D*$-and $*D$-vector space only when we work with coordinate representation. When we speak vector space is $D*$- or $*D$- we point out how we multiply coordinates of vector over elements of division ring: from left or right.

Definition 5.2.1. Let $u$, $v$ be vectors of $D*$-vector space $V$. Vector $w$ is called **linear composition of vectors** $u$ and $v$ when we can write

$$w = au + bv$$

where $a$ and $b$ are scalars. □

We can extend definition of the linear composition on any finite set of vectors. Using generalized indexes to enumerate vectors we can represent set of vectors as one dimensional matrix. We use the convention that we represent any set of vectors of the vector space or as ${}_*$-row (row matrix) either as ${}^*$-row (column matrix). This representation defines type of notation of linear composition. Getting this representation in $D*$-or $*D$-vector space we get four different models of vector space considered in exemples 5.2.2, 5.2.3, 5.2.4, 5.2.5.

For an opportunity to show without change of the notation what kind of vector space ($D*$- or $*D$-) we study we introduce new notation. The symbol $D*$ is called **vector space type** and this symbol means that we study $D*$-vector space. The symbol of product in the type of vector space points to matrix operation used in the linear composition.

Example 5.2.2. Let ${}^*$-row (column matrix)

$$a = \begin{pmatrix} {}^1a \\ ... \\ {}^na \end{pmatrix}$$

represent the set of vectors ${}^ia$, $i \in I$, of $D*$-vector space $V$ and ${}_*$-row (row matrix)

$$c = \begin{pmatrix} c_1 & ... & c_n \end{pmatrix}$$

represent the set of scalars $c_i$, $i \in I$. Then we can write the linear composition of vectors ${}^ia$ as

$$c_i\ {}^ia = c_*{}^*a$$

Such implementation of $D*$-vector space is called **$D_*{}^*$-vector space** or **left $D$-vector space of rows**. □

Example 5.2.3. Let ${}_*$-row (row matrix)

$$a = \begin{pmatrix} {}_1a & ... & {}_na \end{pmatrix}$$

represent the set of vectors ${}_ia$, $i \in I$, of $D*$-vector space $V$ and ${}^*$-row (column matrix)

$$c = \begin{pmatrix} c^1 \\ ... \\ c^n \end{pmatrix}$$

represent the set of scalars $c^i$, $i \in I$. Then we can write the linear composition of vectors ${}_ia$ as

$$c^i\ {}_ia = c^*{}_*a$$

Such implementation of $D*$-vector space is called **$D^*{}_*$-vector space** or **left $D$-vector space of columns**. □

Example 5.2.4. Let ${}^*$-row (column matrix)

$$a=\begin{pmatrix}a^1\\ ...\\ a^n\end{pmatrix}$$

represent the set of vectors $a^i$, $i\in I$, of $*D$-vector space $V$ and ${}_*$-row (row matrix)

$$c=\begin{pmatrix}{}_1c & ... & {}_nc\end{pmatrix}$$

represent the set of scalars ${}_ic$, $i\in I$. Then we can write the linear composition of vectors $a^i$ as

$$a^i\ {}_ic=a^*{}_*c$$

Such implementation of $*D$-vector space is called ${}^*{}_*D$**-vector space** or **right** $D$**-vector space of rows**. □

Example 5.2.5. Let ${}_*$-row (row matrix)

$$a=\begin{pmatrix}a_1 & ... & a_n\end{pmatrix}$$

represent the set of vectors $a_i$, $i\in I$, of $*D$-vector space $V$ and ${}^*$-row (column matrix)

$$c=\begin{pmatrix}{}^1c\\ ...\\ {}^nc\end{pmatrix}$$

represent the set of scalars ${}^ic$, $i\in I$. Then we can write the linear composition of vectors $a_i$ as

$$a_i\ {}^ic=a_*{}^*c$$

Such implementation of $*D$-vector space is called ${}_*{}^*D$**-vector space** or **right** $D$**-vector space of columns**. □

Remark 5.2.6. We extend to vector space and its type convention described in remark 2.2.15. For instance, we execute operations in expression

$$A_*{}^*B_*{}^*v\lambda$$

from left to right. This corresponds to the ${}_*{}^*D$-vector space. However we can execute product from right to left. In custom notation this expression is

$$\lambda v^*{}_*B^*{}_*A$$

and corresponds to $D^*{}_*$-vector space. Similarly, reading this expression from down up we get expression

$$A^*{}_*B^*{}_*v\lambda$$

corresponding to ${}^*{}_*D$-vector space. □

## 5.3. Basis of ${}_*{}^*D$-Vector Space

Definition 5.3.1. Vectors $a_i$, $i\in I$, of ${}_*{}^*D$-vector space $V$ are **linearly independent** if $c=0$ follows from the equation

$$a_*{}^*c=0$$

Otherwise vectors $a_i$ are **linearly dependent**. □

Definition 5.3.2. We call set of vectors $\overline{\overline{e}}=(e_i, i\in I)$ a **basis for ${}_*{}^*D$-vector space** if vectors $e_i$ are linearly independent and adding to this system any other vector we get a new system which is linearly dependent. □

Theorem 5.3.3. *If $\overline{\overline{e}}$ is a basis ${}_*{}^*D$-of vector space $V$ then any vector $\overline{v}\in V$ has one and only one expansion*

$$\overline{v}=e_*{}^*v \tag{5.3.1}$$

*relative to this basis.*

Proof. Because system of vectors $e_i$ is a maximal set of linearly independent vectors the system of vectors $\overline{v}$, $e_i$ is linearly dependent and in equation

$$\overline{v}b+e_*{}^*c=0 \tag{5.3.2}$$

at least $b$ is different from 0. Then equation

$$\overline{v}=e_*{}^*(-cb^{-1}) \tag{5.3.3}$$

follows from (5.3.2). (5.3.1) follows from (5.3.3).

Assume we get another expansion

$$\overline{v}=e_*{}^*v' \tag{5.3.4}$$

We subtract (5.3.1) from (5.3.4) and get

$$0=e_*{}^*(v'-v)$$

Because vectors $e_i$ are linearly independent we get

$$v'-v=0$$

□

Definition 5.3.4. We call the matrix $v$ in expansion (5.3.1) **coordinate matrix of vector $\overline{v}$** in basis $\overline{\overline{e}}$ and we call its entries **coordinates of vector $\overline{v}$** relative to basis $\overline{\overline{e}}$. □

Theorem 5.3.5. *Set of coordinates $a$ of vector $\overline{a}$ relative basis $\overline{\overline{e}}$ of ${}_*{}^*D$-vector space forms ${}_*{}^*D$-vector space $D^n$ isomorphic ${}_*{}^*D$-vector space $V$. This ${}_*{}^*D$-vector space is called* **coordinate ${}_*{}^*D$-vector space**. *This isomorphism is called* **coordinate isomorphism**.

Proof. Suppose vectors $\overline{a}$ and $\overline{b}\in V$ have expansion

$$\overline{a}=e_*{}^*a$$
$$\overline{b}=e_*{}^*b$$

relative basis $\overline{\overline{e}}$. Then

$$\overline{a}+\overline{b}=e_*{}^*a+e_*{}^*b=e_*{}^*(a+b)$$
$$\overline{a}m=(e_*{}^*a)m=e_*{}^*(am)$$

for any $m\in D$. Thus, operations in a vector space are defined by coordinates

$$ {}^i(a+b)={}^ia+{}^ib $$
$$ {}^i(am)={}^iam $$

This completes the proof. □

Example 5.3.6. Let $\overline{\overline{e}}=({}^je,j\in J,|J|=n)$ be a **basis for $D_*{}^*$-vector space** $V$. According to the example 5.2.2, we can present the basis $\overline{\overline{e}}$ as ${}^*$-row (column matrix)

$$ \overline{\overline{e}}=\begin{pmatrix}{}^1e\\ \dots\\ {}^ne\end{pmatrix} $$

The coordinate matrix

$$ a=\begin{pmatrix}a_1&\dots&a_n\end{pmatrix}=(a_j,j\in J) $$

of vector $\overline{a}$ in basis $\overline{\overline{e}}$ is called **$D_*{}^*$-vector**[5.2] or **row $D*$-vector**.

Let ${}^*$-row

(5.3.5)
$$ \overline{A}=\begin{pmatrix}{}^1\overline{A}\\ \dots\\ {}^m\overline{A}\end{pmatrix}=({}^i\overline{A},i\in I) $$

be set of vectors. Vectors ${}^i\overline{A}$ have expansion

$$ {}^i\overline{A}={}^iA_*{}^*e $$

If we substitute coordinate matrices of vectors ${}^i\overline{A}$ in the matrix (5.3.5) we get matrix

$$ A=\begin{pmatrix}\begin{pmatrix}{}^1A_1&\dots&{}^1A_n\end{pmatrix}\\ \dots\\ \begin{pmatrix}{}^mA_1&\dots&{}^mA_n\end{pmatrix}\end{pmatrix}=\begin{pmatrix}{}^1A_1&\dots&{}^1A_n\\ \dots&\dots&\dots\\ {}^mA_1&\dots&{}^mA_n\end{pmatrix}=({}^iA_j) $$

We call the matrix $A$ **coordinate matrix of set of vectors** $({}^i\overline{A},i\in I)$ in basis $\overline{\overline{e}}$ and we call its elements **coordinates of set of vectors** $({}^i\overline{A},i\in I)$ in basis $\overline{\overline{e}}$.

Let ${}^*$-row

$$ \overline{\overline{f}}=\begin{pmatrix}{}^1f\\ \dots\\ {}^nf\end{pmatrix}=({}^if,i\in J) $$

be the basis for $D_*{}^*$-vector space $V$. We tell that coordinate matrix $f$ of set of vectors $({}^if,i\in J)$ defines **coordinates** ${}^if_j$ **of basis** $\overline{\overline{f}}$ relative basis $\overline{\overline{e}}$. □

[5.2] $D_*{}^*$-vector is an analogue of row vector.

EXAMPLE 5.3.7. Let $\overline{\overline{e}} = ({}_je, j \in J, |J| = n)$ be a **basis for $D^*{}_*$-vector space $V$**. According to the example 5.2.3, we can present the basis $\overline{\overline{e}}$ as ${}_*$-row (row matrix)

$$\overline{\overline{e}} = \begin{pmatrix} {}_1e & ... & {}_ne \end{pmatrix}$$

The coordinate matrix

$$a = \begin{pmatrix} a^1 \\ ... \\ a^n \end{pmatrix} = (a_j, j \in J)$$

of vector $\overline{a}$ in basis $\overline{\overline{e}}$ is called **$D^*{}_*$-vector**[5.3] or **column $D*$-vector**.

Let ${}_*$-row

$$\overline{A} = \begin{pmatrix} {}_1\overline{A} & ... & {}_m\overline{A} \end{pmatrix} = ({}^i\overline{A}, i \in I) \tag{5.3.6}$$

be set of vectors. Vectors ${}^i\overline{A}$ have expansion

$$ {}_i\overline{A} = {}_iA^*{}_*e$$

If we substitute coordinate matrices of vectors ${}_i\overline{A}$ in the matrix (5.3.6) we get matrix

$$A = \begin{pmatrix} \begin{pmatrix} {}_1A^1 \\ ... \\ {}_1A^n \end{pmatrix} & ... & \begin{pmatrix} {}_mA^1 \\ ... \\ {}_mA^n \end{pmatrix} \end{pmatrix} = \begin{pmatrix} {}_1A^1 & ... & {}_1A^n \\ ... & ... & ... \\ {}_mA^1 & ... & {}_mA^n \end{pmatrix} = ({}_iA^j)$$

We call the matrix $A$ **coordinate matrix of set of vectors** $({}_i\overline{A}, i \in I)$ in basis $\overline{\overline{e}}$ and we call its elements **coordinates of set of vectors** $({}_i\overline{A}, i \in I)$ in basis $\overline{\overline{e}}$.

Let ${}_*$-row

$$\overline{\overline{f}} = \begin{pmatrix} {}_1f & ... & {}_nf \end{pmatrix} = ({}_if, i \in J)$$

be the basis for $D^*{}_*$-vector space $V$. We tell that coordinate matrix $f$ of set of vectors $({}_if, i \in J)$ defines **coordinates** ${}_if^j$ **of basis $\overline{\overline{f}}$** relative basis $\overline{\overline{e}}$. □

EXAMPLE 5.3.8. Let $\overline{\overline{e}} = (e^i, i \in I, |I| = n)$ be a **basis for ${}^*{}_*D$-vector space** $V$. According to the example 5.2.4, we can present the basis $\overline{\overline{e}}$ as ${}^*$-row (column matrix)

$$\overline{\overline{e}} = \begin{pmatrix} e^1 \\ ... \\ e^n \end{pmatrix}$$

The coordinate matrix

$$a = \begin{pmatrix} {}_1a & ... & {}_na \end{pmatrix} = ({}_ia, i \in I)$$

of vector $\overline{a}$ in basis $\overline{\overline{e}}$ is called **${}^*{}_*D$-vector**[5.4] or **row $*D$-vector**.

[5.3] $D_*{}^*$-vector is an analogue of row vector.

[5.4] ${}_*{}^*D$-vector is an analogue of column vector.

Let ${}^*$-row

$$\overline{A} = \begin{pmatrix} \overline{A}^1 \\ ... \\ \overline{A}^m \end{pmatrix} = (\overline{A}_j, j \in J) \tag{5.3.7}$$

be set of vectors. Vectors $\overline{A}_j$ have expansion

$$\overline{A}^j = e^*{}_* A^j$$

If we substitute coordinate matrices of vectors $\overline{A}_j$ in the matrix (5.3.7) we get matrix

$$A = \begin{pmatrix} \begin{pmatrix} {}_1A^1 & ... & {}_nA^1 \end{pmatrix} \\ ... \\ \begin{pmatrix} {}_1A^m & ... & {}_nA^m \end{pmatrix} \end{pmatrix} = \begin{pmatrix} {}_1A^1 & ... & {}_nA^1 \\ ... & ... & ... \\ {}_1A^m & ... & {}_nA^m \end{pmatrix} = ({}^iA_j)$$

We call the matrix $A$ **coordinate matrix of set of vectors** $(\overline{A}_j, j \in J)$ in basis $\overline{\overline{e}}$ and we call its entries **coordinates of set of vectors** $(\overline{A}_j, j \in J)$ in basis $\overline{\overline{e}}$.

Let ${}^*$-row

$$\overline{\overline{f}} = \begin{pmatrix} f^1 \\ ... \\ f^n \end{pmatrix} = (f^j, j \in I)$$

be the basis of ${}^*{}_*D$-vector space $V$. We tell that coordinate matrix $f$ of set of vectors $(f^j, j \in I)$ defines **coordinates** ${}_if^j$ **of basis** $\overline{\overline{f}}$ relative basis $\overline{\overline{e}}$. □

Example 5.3.9. Let $\overline{\overline{e}} = (e_i, i \in I, |I| = n)$ be a **basis for** ${}_*{}^*D$**-vector space** $V$. According to the example 5.2.5, we can present the basis $\overline{\overline{e}}$ as ${}_*$-row (row matrix)

$$\overline{\overline{e}} = \begin{pmatrix} e_1 & ... & e_n \end{pmatrix}$$

The coordinate matrix

$$a = \begin{pmatrix} {}^1a \\ ... \\ {}^na \end{pmatrix} = ({}^ia, i \in I)$$

of vector $\overline{a}$ in basis $\overline{\overline{e}}$ is called ${}_*{}^*D$**-vector**[5.5] or **column** $D*$**-vector**.

Let ${}_*$-row

$$\overline{A} = \begin{pmatrix} \overline{A}_1 & ... & \overline{A}_m \end{pmatrix} = (\overline{A}_j, j \in J) \tag{5.3.8}$$

be set of vectors. Vectors $\overline{A}_j$ have expansion

$$\overline{A}_j = e_*{}^* A_j$$

---

[5.5] ${}_*{}^*D$-vector is an analogue of column vector.

If we substitute coordinate matrices of vectors $\overline{A}_j$ in the matrix (5.3.8) we get matrix

$$A=\left(\begin{pmatrix}{}^1A_1\\ ...\\ {}^nA_1\end{pmatrix} \ ... \ \begin{pmatrix}{}^1A_m\\ ...\\ {}^nA_m\end{pmatrix}\right)=\begin{pmatrix}{}^1A_1 & ... & {}^1A_m\\ ... & ... & ...\\ {}^nA_1 & ... & {}^nA_m\end{pmatrix}=({}^iA_j)$$

We call the matrix $A$ **coordinate matrix of set of vectors** $(\overline{A}_j, j\in J)$ in basis $\overline{\overline{e}}$ and we call its entries **coordinates of set of vectors** $(\overline{A}_j, j\in J)$ in basis $\overline{\overline{e}}$.

Let ${}_*$-row

$$\overline{\overline{f}}=\begin{pmatrix} f_1 & ... & f_n\end{pmatrix}=(f_j, j\in I)$$

be the basis of ${}_*{}^*D$-vector space $V$. We tell that coordinate matrix $f$ of set of vectors $(f_j, j\in I)$ defines **coordinates** ${}^if_j$ **of basis** $\overline{\overline{f}}$ relative basis $\overline{\overline{e}}$. □

Since we express linear composition using matrices we can extend the duality principle to the vector space theory. We can write duality principle in one of the following forms

THEOREM 5.3.10 (duality principle). *Let $\mathfrak{A}$ be true statement about vector spaces. If we exchange the same time*

- *$D_*{}^*$-vector and $D^*{}_*$-vector*
- *${}_*{}^*D$-vector and ${}^*{}_*D$-vector*
- *${}_*{}^*$-product and ${}^*{}_*$-product*

*then we soon get true statement.*

THEOREM 5.3.11 (duality principle). *Let $\mathfrak{A}$ be true statement about vector spaces. If we exchange the same time*

- *$D_*{}^*$-vector and ${}_*{}^*D$-vector or $D^*{}_*$-vector and ${}^*{}_*D$-vector*
- *${}_*{}^*$-quasideterminant and ${}^*{}_*$-quasideterminant*

*then we soon get true statement.*

## 5.4. Linear Map of ${}_*{}^*D$-Vector Spaces

DEFINITION 5.4.1. Suppose $V$ is ${}_*{}^*S$-vector space. Suppose $U$ is ${}_*{}^*T$-vector space. Morphism

$$f:S\longrightarrow T \qquad A:V\longrightarrow U$$

of right representations of division ring in Abelian group is called **linear map** of ${}_*{}^*S$-vector space $V$ into ${}_*{}^*T$-vector space $U$. □

By theorem 3.2.16, studying linear map we can consider case $S=T$.

DEFINITION 5.4.2. Suppose $V$ and $W$ are ${}_*{}^*D$-vector spaces. We call map

$$A:V\to W$$

**linear map** of ${}_*{}^*D$-vector space if[5.6]

(5.4.1) $$A(m_*{}^*a)=A(m)_*{}^*a$$

for any ${}^ia\in D$, $m_i\in V$. □

[5.6] Expression $A(m)_*{}^*a$ means expression $A(m_i)\ {}^ia$

Theorem 5.4.3. *Let*

$$\overline{\overline{f}}=(f_i,i\in I)$$

*be a basis of* ${}_*{}^*D$*-vector space* $V$ *and*

$$\overline{\overline{e}}=(e_j,j\in J)$$

*be a basis of* ${}_*{}^*D$*-vector space* $U$*. Then linear map*

$$A:V\to W \tag{5.4.2}$$

*of* ${}_*{}^*D$*-vector spaces has presentation*

$$b=A_*{}^*a \tag{5.4.3}$$

*relative to selected bases. Here*

- $a$ *is coordinate matrix of vector* $\overline{a}$ *relative the basis* $\overline{\overline{f}}$*.*
- $b$ *is coordinate matrix of vector*

$$\overline{b}=\overline{A}(\overline{a})$$

  *relative the basis* $\overline{\overline{e}}$*.*
- $A$ *is coordinate matrix of set of vectors* $(\overline{A}(\overline{f}_i))$ *relative the basis* $\overline{\overline{e}}$*. The matrix* $A$ *is called* **matrix of linear map** *relative bases* $\overline{\overline{f}}$ *and* $\overline{\overline{e}}$*.*

Proof. Vector $\overline{a}\in V$ has expansion

$$\overline{a}=f_*{}^*a$$

relative to the basis $\overline{\overline{f}}$. Vector $\overline{b}\in U$ has expansion

$$\overline{b}=e_*{}^*b \tag{5.4.4}$$

relative to the basis $\overline{\overline{e}}$.

Since $\overline{A}$ is a linear map, from (5.4.1) it follows that

$$\overline{b}=\overline{A}(\overline{a})=\overline{A}(f_*{}^*a)=\overline{A}(f)_*{}^*a \tag{5.4.5}$$

$\overline{A}(f_i)$ is a vector of ${}_*{}^*D$-vector space $U$ and has expansion

$$\overline{A}(\overline{f}_i)=\overline{e}_*{}^*A_i=\overline{e}_j{}^jA_i \tag{5.4.6}$$

relative to basis $\overline{\overline{e}}$. Combining (5.4.5) and (5.4.6) we get

$$\overline{b}=e_*{}^*A_*{}^*a \tag{5.4.7}$$

(5.4.3) follows from comparison of (5.4.4) and (5.4.7) and theorem 5.3.3. □

On the basis of theorem 5.4.3 we identify the linear map (5.4.2) of ${}_*{}^*D$-vector spaces and the matrix of its presentation (5.4.3).

Theorem 5.4.4. *Let*

$$\overline{\overline{f}}=(f_i,i\in I)$$

*be a basis of* ${}_*{}^*D$*-vector space* $V$*,*

$$\overline{\overline{e}}=(e_j,j\in J)$$

*be a basis of* ${}_*{}^*D$*-vector space* $U$*, and*

$$\overline{\overline{g}}=(g_l,l\in L)$$

*be a basis of ${}_*{}^*D$-vector space $W$. Suppose diagram of maps*

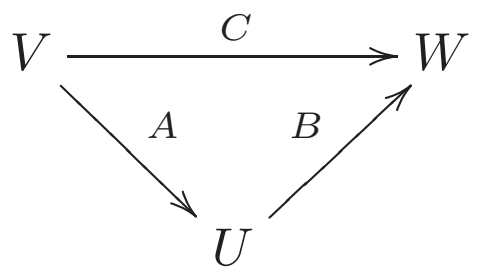

*is commutative diagram where linear map $A$ has presentation*

$$b = A_*{}^*a \tag{5.4.8}$$

*relative to selected bases and linear map $B$ has presentation*

$$c = B_*{}^*b \tag{5.4.9}$$

*relative to selected bases. Then map $C$ is linear map and has presentation*

$$c = B_*{}^*A_*{}^*a \tag{5.4.10}$$

*relative to selected bases.*

Proof. The map $C$ is linear map because

$$\begin{aligned} C_*{}^*(f_*{}^*a) &= (A_*{}^*B)_*{}^*(f_*{}^*a) = B_*{}^*(A_*{}^*(f_*{}^*a)) \\ &= B_*{}^*(e_*{}^*(A_*{}^*a)) = g_*{}^*(B_*{}^*(A_*{}^*a)) \\ &= g_*{}^*((B_*{}^*A)_*{}^*a) = g_*{}^*(C_*{}^*a) \end{aligned}$$

Equation (5.4.10) follows from substituting (5.4.8) into (5.4.9). □

Presenting linear map as ${}_*{}^*$-product we can rewrite (5.4.1) as

$$\overline{A}_*{}^*(\overline{a}k) = (\overline{A}_*{}^*\overline{a})k \tag{5.4.11}$$

We can express the statement of the theorem 5.4.4 in the next form

$$\overline{B}_*{}^*(\overline{A}_*{}^*\overline{a}) = (\overline{B}_*{}^*\overline{A})_*{}^*\overline{a} \tag{5.4.12}$$

Equations (5.4.11) and (5.4.12) represent the **associative law for linear maps of ${}_*{}^*D$-vector spaces**. This allows us writing of such expressions without using of brackets.

Equation (5.4.3) is coordinate notation for linear map. Based theorem 5.4.3, non coordinate notation also can be expressed using ${}_*{}^*$-product

$$\overline{b} = \overline{A}_*{}^*\overline{a} = \overline{A}_*{}^*\overline{f}_*{}^*a = \overline{e}_*{}^*A_*{}^*a \tag{5.4.13}$$

If we substitute equation (5.4.13) into theorem 5.4.4, then we get chain of equations

$$\overline{c} = \overline{B}_*{}^*\overline{b} = \overline{B}_*{}^*\overline{e}_*{}^*b = \overline{g}_*{}^*B_*{}^*b$$

$$\overline{c} = \overline{B}_*{}^*\overline{A}_*{}^*\overline{a} = \overline{B}_*{}^*\overline{A}_*{}^*\overline{f}_*{}^*a = \overline{g}_*{}^*B_*{}^*A_*{}^*a$$

Remark 5.4.5. One can easily see from the example of linear map how theorem 3.2.16 makes our reasoning simpler in study of the morphism of representations of $\Omega$-algebra. In the framework of this remark, we agree to call the theory of linear maps reduced theory, and theory stated in this remark is called enhanced theory.

Suppose $V$ is ${}_*{}^*S$-vector space. Suppose $U$ is ${}_*{}^*T$-vector space. Suppose

$$r : S \longrightarrow T \qquad \overline{A} : V \longrightarrow U$$

is linear map of ${}_*{}^*S$-vector space $V$ into ${}_*{}^*T$-vector space $U$. Let

$$\overline{\overline{f}} = (f_i, i \in I)$$

be a basis of ${}_*{}^*S$-vector space $V$ and

$$\overline{\overline{e}}=(e_j,j\in J)$$

be a basis of ${}_*{}^*T$-vector space $U$.

From definitions 5.4.1 and 3.2.2 it follows

$$\overline{b}=\overline{A}(\overline{a})=\overline{A}(f_*{}^*a)=\overline{A}(f)_*{}^*r(a) \tag{5.4.14}$$

$\overline{A}(f_i)$ is also a vector of $U$ and has expansion

$$\overline{A}(\overline{f}_i)=\overline{e}_*{}^*A_i=\overline{e}_j\ {}^jA_i \tag{5.4.15}$$

relative to basis $\overline{\overline{e}}$. Combining (5.4.14) and (5.4.15), we get

$$\overline{b}=e_*{}^*A_*{}^*r(a) \tag{5.4.16}$$

Suppose $W$ is ${}_*{}^*D$-vector space. Suppose

$$p:T\longrightarrow D \qquad \overline{B}:U\longrightarrow W$$

is linear map of ${}_*{}^*T$-vector space $U$ into ${}_*{}^*D$-vector space $W$. Let

$$\overline{\overline{g}}=(g_l,l\in L)$$

be basis of ${}_*{}^*D$-vector space $W$. Then, according to (5.4.16), the product of linear map $(r,\overline{A})$ and linear map $(p,\overline{B})$ has form

$$\overline{c}=h_*{}^*B_*{}^*p(A)_*{}^*pr(a) \tag{5.4.17}$$

Comparison of equations (5.4.10) and (5.4.17) shows that extended theory of linear maps is more complicated then reduced theory.

If we need we can use extended theory, however we will not get new results comparing with reduced theory. At the same time plenty of details makes picture less clear and demands permanent attention. □

## 5.5. System of Linear Equations

Definition 5.5.1. Let $V$ be a ${}_*{}^*D$-vector space and $\{A_i\in V,i\in I\}$ be set of vectors. **Linear span in ${}_*{}^*D$-vector space** is set $\mathrm{span}(A_i,i\in I)$ of vectors linearly dependent on vectors $A_i$. □

Theorem 5.5.2. *Let $span(A_i,i\in I)$ be linear span in ${}_*{}^*D$-vector space $V$. Then $span(A_i,i\in I)$ is subspace of ${}_*{}^*D$-vector space $V$.*

Proof. Suppose

$$\overline{b}\in\mathrm{span}(A_i,i\in I)$$
$$\overline{c}\in\mathrm{span}(A_i,i\in I)$$

According to definition 5.5.1

$$\overline{b}=A_*{}^*b$$
$$\overline{c}=A_*{}^*c$$

Then

$$\overline{b}+\overline{c}=A_*{}^*b+A_*{}^*c=A_*{}^*(b+c)\in\mathrm{span}(A_i,i\in I)$$
$$\overline{b}k=(A_*{}^*b)k=A_*{}^*(bk)\in\mathrm{span}(A_i,i\in I)$$

This proves the statement. □

EXAMPLE 5.5.3. Let $V$ be a ${}_*{}^*D$-vector space and ${}_*$-row

$$\overline{\overline{A}} = \begin{pmatrix} \overline{A}_1 & ... & \overline{A}_n \end{pmatrix} = (\overline{A}_i, i \in I)$$

be set of vectors. To answer the question of whether vector $\overline{b} \in \mathrm{span}(\overline{A}_i, i \in I)$ we write linear equation

$$\overline{b} = \overline{A}_*{}^*x \tag{5.5.1}$$

where

$$x = \begin{pmatrix} {}^1x \\ ... \\ {}^nx \end{pmatrix}$$

is ${}^*$-row of unknown coefficients of expansion. $\overline{b} \in \mathrm{span}(\overline{A}_i, i \in I)$ if equation (5.5.1) has a solution. Suppose $\overline{\overline{f}} = (f_j, j \in J)$ is a basis. Then vectors $\overline{b}$, $\overline{A}_i$ have expansion

$$\overline{b} = f_*{}^*b \tag{5.5.2}$$

$$\overline{A}_i = f_*{}^*A_i \tag{5.5.3}$$

If we substitute (5.5.2) and (5.5.3) into (5.5.1) we get

$$f_*{}^*b = f_*{}^*A_*{}^*x \tag{5.5.4}$$

Applying theorem 5.3.3 to (5.5.4) we get **system of linear equations**

$$A_*{}^*x = b \tag{5.5.5}$$

We can write system of linear equations (5.5.5) in one of the next forms

$$\begin{pmatrix} {}^1A_1 & ... & {}^1A_n \\ ... & ... & ... \\ {}^mA_1 & ... & {}^mA_n \end{pmatrix} {}_*{}^* \begin{pmatrix} {}^1x \\ ... \\ {}^nx \end{pmatrix} = \begin{pmatrix} {}^1b \\ ... \\ {}^mb \end{pmatrix}$$

$${}^jA_i\ {}^ix = {}^jb \tag{5.5.6}$$

$$\begin{matrix} {}^1A_1\ {}^1x & +... & +{}^1A_n\ {}^nx & = {}^1b \\ ... & ... & ... & ... \\ {}^1A_m\ {}^1x & +... & +{}^mA_n\ {}^nx & = {}^mb \end{matrix}$$

□

EXAMPLE 5.5.4. Let $V$ be a $D_*{}^*$-vector space and ${}^*$-row

$$\overline{\overline{A}} = \begin{pmatrix} {}^1\overline{A} \\ ... \\ {}^m\overline{A} \end{pmatrix} = ({}^j\overline{A}, j \in J)$$

be set of vectors. To answer the question of whether vector $\overline{b} \in \mathrm{span}({}^jA, j \in J)$ we write linear equation

$$\overline{b} = x_*{}^*\overline{A} \tag{5.5.7}$$

where

$$x=\begin{pmatrix} x_1 & ... & x_m \end{pmatrix}$$

is ${}_*$-row of unknown coefficients of expansion. $\overline{b}\in \mathrm{span}({}^jA,j\in J)$ if equation (5.5.7) has a solution. Suppose $\overline{\overline{f}}=({}^i\overline{f},i\in I)$ is a basis. Then vectors $\overline{b}$, ${}^jA$ have expansion

$$\overline{b}=b_*{}^*f \tag{5.5.8}$$

$${}^jA={}^jA_*{}^*f \tag{5.5.9}$$

If we substitute (5.5.8) and (5.5.9) into (5.5.7) we get

$$b_*{}^*f=x_*{}^*A_*{}^*f \tag{5.5.10}$$

Applying theorem 5.3.3 to (5.5.10) we get **system of linear equations**[5.7]

$$x_*{}^*A=b \tag{5.5.11}$$

We can write system of linear equations (5.5.11) in one of the next forms

$$\begin{pmatrix} x_1 & ... & x_m \end{pmatrix}{}_*{}^*\begin{pmatrix} {}^1A_1 & ... & {}^1A_n \\ ... & ... & ... \\ {}^mA_1 & ... & {}^mA_n \end{pmatrix}=\begin{pmatrix} b_1 & ... & b_n \end{pmatrix}$$

$$x_j\ {}^jA_i=b_i \tag{5.5.12}$$

$$\begin{matrix} x_1\ {}^1A_1 & ... & x_1\ {}^1A_n \\ + & ... & + \\ ... & ... & ... \\ + & ... & + \\ x_m\ {}^mA_1 & ... & x_m\ {}^mA_n \\ = & ... & = \\ b_1 & ... & b_n \end{matrix}$$

□

To find a solution of system of linear equations, we need a matrix of this system. From examples 5.3.7, 5.3.9, we see that we may consider a column of matrix as vector of left or right vector space. To make statements more clear, we will use type of vector space before word linear. For instance, the statement

> Columns of matrix are $D^*{}_*$-linear dependent.

means that

> Columns of matrix are vectors of $D^*{}_*$-vector space, and corresponding vectors are linear dependent.

In particular, system of linear equations (5.5.5) in ${}_*{}^*D$-vector space is called **system of ${}_*{}^*D$-linear equations**. and system of linear equations (5.5.11) in $D_*{}^*$-vector space is called **system of $D_*{}^*$-linear equations**.

[5.7] Reading system of linear equations (5.5.5) in ${}_*{}^*D$-vector space from bottom up and from left to right we get system of linear equations (5.5.11) in $D_*{}^*$-vector space.

Definition 5.5.5. If $n\times n$ matrix $A$ has ${}_*{}^*$-inverse matrix we call such matrix ${}_*{}^*$**-nonsingular matrix**. Otherwise, we call such matrix ${}_*{}^*$**-singular matrix**. □

Definition 5.5.6. Suppose $A$ is ${}_*{}^*$-nonsingular matrix. We call appropriate system of ${}_*{}^*D$-linear equations

$$A_*{}^*x=b \tag{5.5.13}$$

**nonsingular system of** ${}_*{}^*D$**-linear equations**. □

Theorem 5.5.7. *Solution of nonsingular system of ${}_*{}^*D$-linear equations* (5.5.13) *is determined uniquely and can be presented in either form*[5.8]

$$x=A^{-1_*{}^*}{}_*{}^*b \tag{5.5.14}$$

$$x=\mathcal{H}\det({}_*{}^*)A_*{}^*b \tag{5.5.15}$$

Proof. Multiplying both sides of equation (5.5.13) from left by $A^{-1_*{}^*}$ we get (5.5.14). Using definition (2.3.12) we get (5.5.15). Since theorem 2.2.16 the solution is unique. □

## 5.6. Rank of Matrix

Definition 5.6.1. Matrix[5.9] ${}^SA_T$ is a minor matrix of an order $k$. □

Definition 5.6.2. If minor matrix ${}^SA_T$ is ${}_*{}^*$-nonsingular matrix then we say that ${}_*{}^*$-rank of matrix $A$ is not less then $k$. ${}_*{}^*$**-rank of matrix** $A$

$$\mathrm{rank}_{{}_*{}^*}A$$

is the maximal value of $k$. We call an appropriate minor matrix the ${}_*{}^*$**-major minor matrix**. □

Theorem 5.6.3. *Let matrix $A$ be ${}_*{}^*$-singular matrix and minor matrix ${}^SA_T$ be major minor matrix. Then*

$${}^p\det({}_*{}^*)_r\,{}^{S\cup\{p\}}A_{T\cup\{r\}}=0 \tag{5.6.1}$$

Proof. To understand why minor matrix

$$B={}^{S\cup\{p\}}A_{T\cup\{r\}} \tag{5.6.2}$$

does not have ${}_*{}^*$-inverse matrix,[5.10] we assume that there exists ${}_*{}^*$-inverse matrix $B^{-1_*{}^*}$. We write down the system of linear equations (2.3.2), (2.3.3) in case

[5.8] We can see a solution of system (5.5.13) in theorem [6]-1.6.1. I repeat this statement because I slightly changed the notation.

[5.9] In this section, we will make the following assumption.

- $i\in M$, $|M|=m$, $j\in N$, $|N|=n$.
- $A=({}^iA_j)$ is an arbitrary matrix.
- $k,\ s\in S\supseteq M$, $l,\ t\in T\supseteq N$, $k=|S|=|T|$.
- $p\in M\setminus S$, $r\in N\setminus T$.

[5.10] It is natural to expect relationship between ${}_*{}^*$-singularity of the matrix and its ${}_*{}^*$-quasideterminant similar to relationship which is known in commutative case. However ${}_*{}^*$-quasideterminant is defined not always. For instance, it is not defined when ${}_*{}^*$-inverse matrix has too much elements equal 0. As it follows from this theorem, the ${}_*{}^*$-quasideterminant is undefined also in case when ${}_*{}^*$-rank of the matrix is less then $n-1$.

$I = \{r\}$, $J = \{p\}$ (then $[I] = T$, $[J] = S$ )

$$ {}^{S}B_{T*}{}^{*T}B^{-1_*{}^*}{}_p + {}^{S}B_r\ {}^{r}B^{-1_*{}^*}{}_p = 0 \tag{5.6.3}$$

$$ {}^{p}B_{T*}{}^{*T}B^{-1_*{}^*}{}_p + {}^{p}B_r\ {}^{r}B^{-1_*{}^*}{}_p = 1 \tag{5.6.4}$$

and will try to solve this system. We multiply (5.6.3) by $({}^{S}B_T)^{-1_*{}^*}$

$$ {}^{T}B^{-1_*{}^*}{}_p + ({}^{S}B_T)^{-1_*{}^*}{}_*{}^{*S}B_r\ {}^{r}B^{-1_*{}^*}{}_p = 0 \tag{5.6.5}$$

Now we can substitute (5.6.5) into (5.6.4)

$$ -{}^{p}B_{T*}{}^*({}^{S}B_T)^{-1_*{}^*}{}_*{}^{*S}B_r\ {}^{r}B^{-1_*{}^*}{}_p + {}^{p}B_r\ {}^{r}B^{-1_*{}^*}{}_p = 1 \tag{5.6.6}$$

From (5.6.6) it follows that

$$ ({}^{p}B_r - {}^{p}B_{T*}{}^*({}^{S}B_T)^{-1_*{}^*}{}_*{}^{*S}B_r)\ {}^{r}B^{-1_*{}^*}{}_p = 1 \tag{5.6.7}$$

Expression in brackets is quasideterminant ${}^{p}\det({}_*{}^*)_r B$. Substituting this expression into (5.6.7), we get

$$ {}^{p}\det({}_*{}^*)_r B\ {}^{r}B^{-1_*{}^*}{}_p = 1 \tag{5.6.8}$$

Thus we proved that quasideterminant ${}^{p}\det({}_*{}^*)_r B$ is defined and its equation to 0 is necessary and sufficient condition that the matrix $B$ is singular. Since (5.6.2) the statement of theorem is proved. □

THEOREM 5.6.4. *Suppose $A$ is a matrix*[5.9],

$$\operatorname{rank}_{{}_*{}^*} A = k < m$$

*and ${}^{S}A_T$ is ${}_*{}^*$-major minor matrix. Then ${}_*$-row ${}^{p}A$ is a $D_*{}^*$-linear composition of ${}_*$-rows ${}^{S}A$.*

$$ {}^{M\setminus S}A = R_*{}^{*S}A \tag{5.6.9}$$

$$ {}^{p}A = {}^{p}R_*{}^{*S}A \tag{5.6.10}$$

$$ {}^{p}A_b = {}^{p}R_s\ {}^{s}A_b \tag{5.6.11}$$

PROOF. If the matrix $A$ has $k$ ${}^*$-rows, then assuming that ${}_*$-row ${}^{p}A$ is a $D_*{}^*$-linear combination (5.6.10) of ${}_*$-rows ${}_sA$ with coefficients ${}^{p}R_s$ we get system of $D_*{}^*$-linear equations (5.6.11). According to theorem 5.5.7 the system of $D_*{}^*$-linear equations (5.6.11) has a unique solution[5.11] and this solution is nontrivial because all ${}_*{}^*$-quasideterminants are different from 0.

It remains to prove this statement in case when a number of ${}^*$-rows of the matrix $A$ is more then $k$. I get ${}_*$-row ${}^{p}A$ and ${}^*$-row $A_r$. According to assumption, minor matrix ${}^{S\cup\{p\}}A_{T\cup\{r\}}$ is a ${}_*{}^*$-singular matrix and its ${}_*{}^*$-quasideterminant

$$ {}^{p}\det({}_*{}^*)_r\ {}^{S\cup\{p\}}A_{T\cup\{r\}} = 0 \tag{5.6.12}$$

According to (2.3.14) the equation (5.6.12) has form

$$ {}^{p}A_r - {}^{p}A_{T*}{}^*(({}^{S}A_T)^{-1_*{}^*})_*{}^{*S}A_r = 0$$

Matrix

$$ {}^{p}R = {}^{p}A_{T*}{}^*(({}^{S}A_T)^{-1_*{}^*}) \tag{5.6.13}$$

does not depend on $r$, Therefore, for any $r \in N \setminus T$

$$ {}^{p}A_r = {}^{p}R_*{}^{*S}A_r \tag{5.6.14}$$

---

[5.11]We assume that unknown variables are $x_s = {}^{p}R_s$

From equation

$$(({}^S A_T)^{-1_*{}^*})_*{}^*{}^S A_l = {}^T\delta_l$$

it follows that

$$ {}^p A_l = {}^p A_T{}_*{}^*{}^T\delta_l = {}^p A_T{}_*{}^*(({}^S A_T)^{-1_*{}^*})_*{}^*{}^S A_l \tag{5.6.15}$$

Substituting (5.6.13) into (5.6.15) we get

$$ {}^p A_l = {}^p R_*{}^*{}^S A_l \tag{5.6.16}$$

(5.6.14) and (5.6.16) finish the proof. □

COROLLARY 5.6.5. *Suppose $A$ is a matrix, $\mathrm{rank}_{{}_*{}^*} A = k < m$. Then ${}_*$-rows of the matrix are $D_*{}^*$-linearly dependent.*

$$\lambda_*{}^*A = 0$$

PROOF. Suppose ${}_*$-row ${}^pA$ is a $D_*{}^*$-linear composition (5.6.10). We assume $\lambda_p = -1$, $\lambda_s = {}^pR_s$ and the rest $\lambda_c = 0$. □

THEOREM 5.6.6. *Let $({}^i\overline{A}, i \in M, |M| = m)$ be set of $D_*{}^*$-linearly independent vectors. Then ${}_*{}^*$-rank of their coordinate matrix equal $m$.*

PROOF. Let $\overline{\overline{e}}$ be the basis of $D_*{}^*$-vector space. According to model built in example 5.3.6, the coordinate matrix of set of vectors $({}^i\overline{A})$ relative basis $\overline{\overline{e}}$ consists from ${}_*$-rows which are coordinate matrices of vectors ${}^i\overline{A}$ relative the basis $\overline{\overline{e}}$. Therefore ${}_*{}^*$-rank of this matrix cannot be more then $m$.

Let ${}_*{}^*$-rank of the coordinate matrix be less then $m$. According to corollary 5.6.5, ${}_*$-rows of matrix are $D_*{}^*$-linear dependent

$$\lambda_*{}^*A = 0 \tag{5.6.17}$$

Suppose $c = \lambda_*{}^*A$. From the equation (5.6.17) it follows that $D_*{}^*$-linear composition

$$c_*{}^*\overline{e} = 0$$

of vectors of basis equal 0. This contradicts to statement that vectors $\overline{e}$ form basis. We proved statement of theorem. □

THEOREM 5.6.7. *Suppose $A$ is a matrix*[5.9]*,*

$$\mathrm{rank}_{{}_*{}^*} A = k < n$$

*and ${}^S A_T$ is ${}_*{}^*$-major minor matrix. Then ${}^*$-row $A_r$ is a ${}_*{}^*D$-linear composition of ${}^*$-rows $A_t$*

$$A_{N\setminus T} = A_T{}_*{}^*R \tag{5.6.18}$$

$$A_r = A_T{}_*{}^*R_r \tag{5.6.19}$$

$${}^aA_r = {}^aA_t\ {}^tR_r \tag{5.6.20}$$

PROOF. If the matrix $A$ has $k$ ${}_*$-rows, then assuming that ${}^*$-row $A_r$ is a ${}_*{}^*D$-linear combination (5.6.19) of ${}^*$-rows $A_t$ with coefficients ${}^tR_r$ we get system of ${}_*{}^*D$-linear equations (5.6.20). According to theorem 5.5.7 the system of ${}_*{}^*D$-linear equations (5.6.20) has a unique solution[5.12] and this solution is nontrivial because all ${}_*{}^*$-quasideterminants are different from 0.

[5.12] We assume that unknown variables are ${}_tx = {}^tR_r$

It remains to prove this statement in case when a number of ${}_*$-rows of the matrix $A$ is more then $k$. I get ${}^*$-row $A_r$ and ${}_*$-row ${}^pA$. According to assumption, minor matrix ${}^{S\cup\{p\}}A_{T\cup\{r\}}$ is a ${}_*{}^*$-singular matrix and its $\overline{RC}$-quasideterminant

$$ {}^p\det({}_*{}^*)_r\ {}^{S\cup\{p\}}A_{T\cup\{r\}} = 0 \tag{5.6.21}$$

According to (2.3.14) (5.6.21) has form

$$ {}^pA_r - {}^pA_T{}_*{}^*(({}^SA_T)^{-1_*{}^*})_*{}^*{}^SA_r = 0$$

Matrix

$$ R_r = (({}^SA_T)^{-1_*{}^*})_*{}^*{}^SA_r \tag{5.6.22}$$

does not depend on $p$, Therefore, for any $p \in M\setminus S$

$$ {}^pA_r = {}^pA_T{}_*{}^*R_r \tag{5.6.23}$$

From equation

$$ {}^kA_T{}_*{}^*(({}^SA_T)^{-1_*{}^*})_s = {}^k\delta_s$$

it follows that

$$ {}^kA_r = {}^k\delta_S{}_*{}^*{}^SA_r = {}^kA_T{}_*{}^*(({}^SA_T)^{-1_*{}^*})^s{}_*{}^*{}^SA_r \tag{5.6.24}$$

Substituting (5.6.22) into (5.6.24) we get

$$ {}^kA_r = {}^kA_T{}_*{}^*R_r \tag{5.6.25}$$

(5.6.23) and (5.6.25) finish the proof. □

Corollary 5.6.8. *Suppose $A$ is a matrix,* $\mathrm{rank}_{*}{}^{*}\,A = k < m$. *Then ${}^*$-rows of the matrix are ${}_*{}^*D$-linearly dependent.*

$$ A_*{}^*\lambda = 0$$

Proof. Suppose ${}^*$-row $A_r$ is a right linear composition (5.6.19). We assume ${}^r\lambda = -1$, ${}^t\lambda = {}^tR_r$ and the rest ${}^c\lambda = 0$. □

Base on theorem 2.3.9 we can write similar statements for ${}^*{}_*$-rank of matrix.

Theorem 5.6.9. *Suppose $A$ is a matrix,*

$$ \mathrm{rank}^{*}{}_{*}\,A = k < m$$

*and ${}_TA^S$ is ${}^*{}_*$-major minor matrix. Then ${}_*$-row $A^p$ is a ${}^*{}_*D$-linear composition of ${}_*$-rows $A^s$.*

$$ A^{M\setminus S} = A^S{}^*{}_*R \tag{5.6.26}$$

$$ A^p = A^S{}^*{}_*R^p \tag{5.6.27}$$

$$ {}_bA^p = {}_bA^s\ {}_sR^p \tag{5.6.28}$$

Corollary 5.6.10. *Suppose $A$ is a matrix,* $\mathrm{rank}^{*}{}_{*}\,A = k < m$. *Then ${}_*$-rows of matrix are ${}^*{}_*D$-linearly dependent.*

$$ A^*{}_*\lambda = 0$$

THEOREM 5.6.11. *Suppose $A$ is a matrix,*

$$\mathrm{rank}_{*_*} A = k < n$$

*and ${}_TA^S$ is ${}^*{}_*$-major minor matrix. Then ${}^*$-row ${}_rA$ is a $D^*{}_*$-linear composition of ${}^*$-rows ${}_tA$*

$$ {}_{N\setminus T}A = R^*{}_*{}_TA \tag{5.6.29}$$

$$ {}_rA = {}_rR^*{}_*{}_TA \tag{5.6.30}$$

$$ {}_rA^a = {}_rR^t\ {}_tA^a \tag{5.6.31}$$

COROLLARY 5.6.12. *Suppose $A$ is a matrix,* $\mathrm{rank}_{*_*}\ a = k < m$*. Then ${}^*$-rows of matrix are $D^*{}_*$-linearly dependent.*

$$\lambda^*{}_*A = 0$$

## 5.7. System of ${}_*{}^*D$-Linear Equations

DEFINITION 5.7.1. Suppose[5.9] $A$ is a matrix of system of $D_*{}^*$-linear equations (5.5.12). We call matrix

$$\begin{pmatrix} {}^jA_i \\ b_i \end{pmatrix} = \begin{pmatrix} {}^1A_1 & ... & {}^1A_n \\ ... & ... & ... \\ {}^mA_1 & ... & {}^mA_n \\ b_1 & ... & b_n \end{pmatrix} \tag{5.7.1}$$

an **extended matrix** of this system. □

DEFINITION 5.7.2. Suppose[5.9] $A$ is a matrix of system of ${}_*{}^*D$-linear equations (5.5.6). We call matrix

$$\begin{pmatrix} {}^jA_i & {}^jb \end{pmatrix} = \begin{pmatrix} {}^1A_1 & ... & {}^1A_n & {}^1b \\ ... & ... & ... & ... \\ {}^mA_1 & ... & {}^mA_n & {}^mb \end{pmatrix} \tag{5.7.2}$$

an **extended matrix** of this system. □

THEOREM 5.7.3. *System of ${}_*{}^*D$-linear equations* (5.5.6) *has a solution iff*

$$\mathrm{rank}_{*^*}({}^jA_i) = \mathrm{rank}_{*^*}\begin{pmatrix} {}^jA_i & {}^jb \end{pmatrix} \tag{5.7.3}$$

PROOF. Let ${}^SA_T$ be ${}_*{}^*$-major minor matrix of matrix $A$.

Let a system of ${}_*{}^*D$-linear equations (5.5.6) have solution ${}^ix = {}^id$. Then

$$A_*{}^*d = b \tag{5.7.4}$$

Equation (5.7.4) can be rewritten in form

$$A_{T*}{}^{*T}d + A_{N\setminus T*}{}^{*N\setminus T}d = b \tag{5.7.5}$$

Substituting (5.6.18) into (5.7.5) we get

$$A_{T*}{}^{*T}d + A_{T*}{}^*R_*{}^{*N\setminus T}d = b \tag{5.7.6}$$

From (5.7.6) it follows that ${}_*$-row $b$ is a ${}_*{}^*D$-linear combination of ${}_*$-rows $A_T$

$$A_{T*}{}^*({}^Td + R_*{}^{*N\setminus T}d) = b$$

This holds equation (5.7.3).

It remains to prove that an existence of solution of system of ${}_*{}^*D$-linear equations (5.5.6) follows from (5.7.3). Holding (5.7.3) means that ${}^SA_T$ is ${}_*{}^*$-major minor matrix of extended matrix as well. From theorem 5.6.7 it follows that ${}_*$-row $b$ is a ${}_*{}^*D$-linear composition of ${}_*$-rows $A_T$

$$b = A_T{}_*{}^*{}^TR$$

Assigning ${}^rR = 0$ we get

$$b = A_*{}^*R$$

Therefore, we found at least one solution of system of ${}_*{}^*D$-linear equations (5.5.6). □

Theorem 5.7.4. *Suppose* (5.5.6) *is a system of ${}_*{}^*D$-linear equations satisfying* (5.7.3). *If* $\operatorname{rank}_{{}_*{}^*} A = k \le m$ *then solution of the system depends on arbitrary values of $m-k$ variables not included into ${}_*{}^*$-major minor matrix.*

Proof. Let ${}^SA_T$ be ${}_*{}^*$-major minor matrix of matrix $a$. Suppose

$$ {}^pA_*{}^*x = {}^pb \tag{5.7.7}$$

is an equation with number $p$. Applying theorem 5.6.4 to extended matrix (5.7.1) we get

$$ {}^pA = {}^pR_*{}^*{}^SA \tag{5.7.8}$$

$$ {}^pb = {}^pR_*{}^*{}^Sb \tag{5.7.9}$$

Substituting (5.7.8) and (5.7.9) into (5.7.7) we get

$$ {}^pR_*{}^*{}^SA_*{}^*x = {}^pR_*{}^*{}^Sb \tag{5.7.10}$$

(5.7.10) means that we can exclude equation (5.7.7) from system (5.5.6) and the new system is equivalent to the old one. Therefore, a number of equations can be reduced to $k$.

At this point, we have two choices. If the number of variables is also $k$ then according to theorem 5.5.7 the system has unique solution (5.5.15). If the number of variables $m > k$ then we can move $m-k$ variables that are not included into ${}_*{}^*$-major minor matrix in right side. Giving arbitrary values to these variables, we determine value of the right side and for this value we get a unique solution according to theorem 5.5.7. □

Corollary 5.7.5. *System of ${}_*{}^*D$-linear equations* (5.5.6) *has a unique solution iff its matrix is nonsingular.* □

Theorem 5.7.6. *Solutions of a homogenous system of ${}_*{}^*D$-linear equations*

$$A_*{}^*x = 0 \tag{5.7.11}$$

*form a ${}_*{}^*D$-vector space.*

Proof. Let $\overline{X}$ be set of solutions of system of ${}_*{}^*D$-linear equations (5.7.11). Suppose $x = ({}^ax) \in \overline{X}$ and $y = ({}^ay) \in \overline{X}$. Then

$$x^a\ {}_aa^b = 0$$

$$y^a\ {}_aa^b = 0$$

Therefore

$${}^iA_j\ ({}^jx + {}^jy) = {}^iA_j\ {}^jx + {}^iA_j\ {}^jy = 0$$

$$x+y=({}^jx+{}^jy)\in\overline{X}$$

The same way we see

$$ {}^iA_j\ ({}^jxb)=({}^iA_j\ {}^jx)b=0$$
$$xb=({}^jxb)\in\overline{X}$$

According to definition 5.1.4 $\overline{X}$ is a ${}_*{}^*D$-vector space. □

## 5.8. Nonsingular Matrix

Suppose $A$ is $n\times n$ matrix. Corollaries 5.6.5 and 5.6.8 tell us that if $\mathrm{rank}_{{}_*{}^*}\ A<n$ then ${}_*$-rows are $D_*{}^*$-linearly dependent and ${}^*$-rows are ${}_*{}^*D$-linearly dependent.[5.13]

Theorem 5.8.1. *Let $A$ be $n\times n$ matrix and ${}^*$-row $A_r$ be a ${}_*{}^*D$-linear combination of other ${}^*$-rows. Then* $\mathrm{rank}_{{}_*{}^*}\ A<n$.

Proof. The statement that ${}^*$-row $A_r$ is a ${}_*{}^*D$-linear combination of other ${}^*$-rows means that system of ${}_*{}^*D$-linear equations

$$A_r=A_{[r]}{}_*{}^*\lambda$$

has at least one solution. According theorem 5.7.3

$$\mathrm{rank}_{{}_*{}^*}\ A=\mathrm{rank}_{{}_*{}^*}\ A_{[r]}$$

Since a number of ${}_*$-rows is less then $n$ we get $\mathrm{rank}_{{}_*{}^*}\ A[r]<n$. □

Theorem 5.8.2. *Let $A$ be $n\times n$ matrix and ${}_*$-row ${}^pA$ be a $D_*{}^*$-linear combination of other ${}_*$-rows. Then* $\mathrm{rank}_{{}_*{}^*}\ A<n$.

Proof. Proof of statement is similar to proof of theorem 5.8.1 □

Theorem 5.8.3. *Suppose $A$ and $B$ are $n\times n$ matrices and*

$$C=A_*{}^*B \tag{5.8.1}$$

*$C$ is ${}_*{}^*$-singular matrix iff either matrix $A$ or matrix $B$ is ${}_*{}^*$-singular matrix.*

Proof. Suppose matrix $B$ is ${}_*{}^*$-singular. According theorem 5.6.7 ${}^*$-rows of matrix $B$ are ${}_*{}^*D$-linearly dependent. Therefore

$$0=B_*{}^*\lambda \tag{5.8.2}$$

where $\lambda\ne 0$. From (5.8.1) and (5.8.2) it follows that

$$C_*{}^*\lambda=A_*{}^*B_*{}^*\lambda=0$$

According theorem 5.8.1 matrix $C$ is ${}_*{}^*$-singular.

Suppose matrix $B$ is not ${}_*{}^*$-singular, but matrix $A$ is ${}_*{}^*$-singular According to theorem 5.6.4 ${}^*$-rows of matrix $A$ are ${}_*{}^*D$-linearly dependent. Therefore

$$0=A_*{}^*\mu \tag{5.8.3}$$

where $\mu\ne 0$. According to theorem 5.5.7 the system

$$B_*{}^*\lambda=\mu$$

has only solution where $\lambda\ne=0$. Therefore

$$C_*{}^*\lambda=A_*{}^*B_*{}^*\lambda=A_*{}^*\mu=0$$

According to theorem 5.8.1 matrix $C$ is ${}_*{}^*$-singular.

[5.13] This statement is similar to proposition [7]-1.2.5.

Suppose matrix $C$ is ${}_*{}^*$-singular matrix. According to the theorem 5.6.4 ${}_*$-rows of matrix $C$ are ${}_*{}^*D$-linearly dependent. Therefore

$$0 = C_*{}^*\lambda \tag{5.8.4}$$

where $\lambda \neq 0$. From (5.8.1) and (5.8.4) it follows that

$$0 = A_*{}^*B_*{}^*\lambda$$

If

$$0 = B_*{}^*\lambda$$

satisfied then matrix $B$ is ${}_*{}^*$-singular. Suppose that matrix $B$ is not ${}_*{}^*$-singular. Let us introduce

$$\mu = B_*{}^*\lambda$$

where $\mu \neq 0$. Then

$$0 = A_*{}^*\mu \tag{5.8.5}$$

From (5.8.5) it follows that matrix $A$ is ${}_*{}^*$-singular. □

Basing theorem 2.3.9 we can write similar statements for $D^*{}_*$-linear combination of ${}_*$-rows or ${}^*{}_*D$-linear combination of ${}^*$-rows and ${}^*{}_*$-quasideterminant.

Theorem 5.8.4. *Let $A$ be $n\times n$ matrix and ${}^*$-row ${}_rA$ be a ${}^*{}_*D$-linear combination of other ${}^*$-rows. Then* $\mathrm{rank}_{{}^*{}_*} A < n$.

Theorem 5.8.5. *Let $A$ be $n \times n$ matrix and ${}_*$-row $A^p$ be a $D^*{}_*$-linear combination of other ${}_*$-rows. Then* $\mathrm{rank}_{{}^*{}_*} A < n$.

Theorem 5.8.6. *Suppose $A$ and $B$ are $n \times n$ matrices and $C = A^*{}_*B$. $C$ is ${}^*{}_*$-singular matrix iff either matrix $A$ or matrix $B$ is ${}^*{}_*$-singular matrix.*

Definition 5.8.7. **${}_*{}^*$-matrix group** $GL(n, {}_*{}^*, D)$ is a group of ${}_*{}^*$-nonsingular matrices where we define ${}_*{}^*$-product of matrices (2.2.1) and ${}_*{}^*$-inverse matrix $A^{-1_*{}^*}$. □

Definition 5.8.8. **${}^*{}_*$-matrix group** $GL(n, {}^*{}_*, D)$ is a group of ${}^*{}_*$-nonsingular matrices where we define ${}^*{}_*$-product of matrices (2.2.2) and ${}^*{}_*$-inverse matrix $A^{-1^*{}_*}$. □

Theorem 5.8.9.

$$GL(n, {}_*{}^*, D) \neq GL(n, {}^*{}_*, D)$$

Remark 5.8.10. From theorem 2.3.11 it follows that there are matrices which are ${}^*{}_*$-nonsingular and ${}_*{}^*$-nonsingular. Theorem 5.8.9 implies that sets of ${}^*{}_*$-nonsingular matrices and ${}_*{}^*$-nonsingular matrices are not identical. For instance, there exists such ${}_*{}^*$-nonsingular matrix which is ${}^*{}_*$-singular matrix. □

Proof. It is enough to prove this statement for $n = 2$. Assume every ${}^*{}_*$-singular matrix

$$A = \begin{pmatrix} {}^1A_1 & {}^2A_1 \\ {}^1A_2 & {}^2A_2 \end{pmatrix} \tag{5.8.6}$$

is ${}_*{}^*$-singular matrix. It follows from theorem 5.6.4 and theorem 5.6.7 that ${}_*{}^*$-singular matrix satisfies to condition

$$ {}^1A_2 = b\ {}^1A_1 \tag{5.8.7}$$

$$ {}^2A_2 = b\ {}^2A_1 \tag{5.8.8}$$

$$ {}^2A_1 = {}^1A_1 c \tag{5.8.9}$$

$$ {}^2A_2 = {}^1A_2 c \tag{5.8.10}$$

If we substitute (5.8.9) into (5.8.8) we get

$$ {}^2A_2 = b\ {}^1A_1\ c$$

$b$ and $c$ are arbitrary elements of division ring $D$ and ${}_*{}^*$-singular matrix matrix (5.8.6) has form ($d = {}^1A_1$)

$$A = \begin{pmatrix} d & dc \\ bd & bdc \end{pmatrix} \tag{5.8.11}$$

The similar way we can show that ${}^*{}_*$-singular matrix has form

$$A = \begin{pmatrix} d & c'd \\ db' & c'db' \end{pmatrix} \tag{5.8.12}$$

From assumption it follows that (5.8.12) and (5.8.11) represent the same matrix. Comparing (5.8.12) and (5.8.11) we get that for every $d, c \in D$ exists such $c' \in D$ which does not depend on $d$ and satisfies equation

$$dc = c'd$$

This contradicts the fact that $D$ is division ring. □

Example 5.8.11. Since we get division ring of quaternions we assume $b = 1+k$, $c = j$, $d = k$. Then we get

$$A = \begin{pmatrix} k & kj \\ (1+k)k & (1+k)kj \end{pmatrix} = \begin{pmatrix} k & -i \\ k-1 & -i-j \end{pmatrix}$$

$$\begin{aligned} {}^2\det({}_*{}^*)_2 A &= {}^2A_2 - {}^2A_1({}^1A_1)^{-1}\ {}^1A_2 \\ &= -i-j-(k-1)(k)^{-1}(-i) = -i-j-(k-1)(-k)(-i) \\ &= -i-j-kki+ki = -i-j+i+j = 0 \end{aligned}$$

$$\begin{aligned} {}_1\det({}^*{}_*)^1 A &= {}_1A^1 - {}_1A^2({}_2A^2)^{-1}\ {}_2A^1 \\ &= k-(k-1)(-i-j)^{-1}(-i) = k-(k-1)\frac{1}{2}(i+j)(-i) \\ &= k+\frac{1}{2}((k-1)i+(k-1)j)i = k+\frac{1}{2}(ki-i+kj-j)i \\ &= k+\frac{1}{2}(j-i-i-j)i = k-ii = k+1 \end{aligned}$$

$$\begin{aligned}{}_1\det({}^*{}_*)^2 A &= {}_1A^2 - {}_1A^1({}_2A^1)^{-1}\ {}_2A^2\\ &= k-1-k(-i)^{-1}(-i-j) = k-1+ki(i+j)\\ &= k-1+j(i+j) = k-1+ji+jj\\ &= k-1-k-1 = -2\end{aligned}$$

$$\begin{aligned}{}_2\det({}^*{}_*)^1 A &= {}_2A^1 - {}_2A^2({}_1A^2)^{-1}\ {}_1A^1\\ &= (-i)-(-i-j)(k-1)^{-1}k = -i+(i+j)\frac{1}{2}(-k-1)k\\ &= -i-\frac{1}{2}(i+j)(k+1)k = -i-\frac{1}{2}(ik+i+jk+j)k\\ &= -i-\frac{1}{2}(-j+i+i+j)k = -i-ik = -i+j\end{aligned}$$

$$\begin{aligned}{}_2\det({}^*{}_*)^2 A &= {}_2A^2 - {}_2A^1({}_1A^1)^{-1}\ {}_1A^2\\ &= -i-j-(-i)(k)^{-1}(k-1) = -i-j+i(-k)(k-1)\\ &= -i-j+j(k-1) = -i-j+jk-j = -i-j+i-j = -2j\end{aligned}$$

The system of ${}_*{}^*D$-linear equations

$$\begin{pmatrix} k & -i\\ k-1 & -i-j\end{pmatrix} {}_*{}^* \begin{pmatrix}{}^1x\\ {}^2x\end{pmatrix} = \begin{pmatrix}{}^1b\\ {}^2b\end{pmatrix} \tag{5.8.13}$$

has ${}_*{}^*$-singular matrix. We can write the system of ${}_*{}^*D$-linear equations (5.8.13) in the form

$$\left\{\begin{array}{rl} k\ {}^1x - & i\ {}^2x = {}^1b\\ (k-1)\ {}^1x - & (i+j)\ {}^2x = {}^2b\end{array}\right.$$

The system of ${}^*{}_*D$-linear equations

$$\begin{pmatrix} k & -i\\ k-1 & -i-j\end{pmatrix} {}^*{}_* \begin{pmatrix}{}_1x & {}_2x\end{pmatrix} = \begin{pmatrix}{}_1b & {}_2b\end{pmatrix} \tag{5.8.14}$$

has ${}^*{}_*$-nonsingular matrix. We can write the system of ${}^*{}_*D$-linear equations (5.8.14) in the form

$$\left\{\begin{array}{rr} k\ {}_1x & -\quad i\ {}_1x\\ +\ (k-1)\ {}_2x & -\ (i+j)\ {}_2x\\ =\quad {}_1b & =\quad {}_2b\end{array}\right. \qquad \left\{\begin{array}{l} k\ {}_1x + (k-1)\ {}_2x = {}_1b\\ -i\ {}_1x - (i+j)\ {}_2x = {}_2b\end{array}\right.$$

The system of $D_*{}^*$-linear equations

$$\begin{pmatrix}x_1 & x_2\end{pmatrix} {}_*{}^* \begin{pmatrix} k & -i\\ k-1 & -i-j\end{pmatrix} = \begin{pmatrix}b_1 & b_2\end{pmatrix} \tag{5.8.15}$$

has ${}_*{}^*$-singular matrix. We can write the system of $D_*{}^*$-linear equations (5.8.15) in the form

$$\left\{\begin{array}{ll} x_1k & -x_1i \\ +\ x_2(k-1) & -x_2(i+j) \\ =b_1 & =b_2 \end{array}\right. \qquad \left\{\begin{array}{lll} x_1k & +x_2(k-1) & =b_1 \\ -x_1i & -x_2(i+j) & =b_2 \end{array}\right.$$

The system of $D^*{}_*$-linear equations

$$\begin{pmatrix} x^1 \\ x^2 \end{pmatrix} {}^*{}_* \begin{pmatrix} k & -i \\ k-1 & -i-j \end{pmatrix} = \begin{pmatrix} {}^1b \\ {}^2b \end{pmatrix} \tag{5.8.16}$$

has ${}^*{}_*$-nonsingular matrix. We can write the system of $D^*{}_*$-linear equations (5.8.16) in the form

$$\left\{\begin{array}{l} k\ {}^1x - i\ {}^2x = {}^1b \\ (k-1)\ {}^1x - (i+j)\ {}^2x = {}^2b \end{array}\right.$$

□

## 5.9. Dimension of ${}_*{}^*D$-Vector Space

Theorem 5.9.1. *Let $V$ be a ${}_*{}^*D$-vector space. Suppose $V$ has bases $\overline{\overline{e}} = (e_i, i \in I)$ and $\overline{\overline{g}} = (g_j, j \in J)$. If $|I|$ and $|J|$ are finite numbers then $|I| = |J|$.*

Proof. Suppose $|I| = m$ and $|J| = n$. Suppose

$$m < n \tag{5.9.1}$$

Because $\overline{\overline{e}}$ is a basis any vector $g_j,\ j \in J$ has expansion

$$g_j = e_*{}^*A_j$$

Because $\overline{\overline{g}}$ is a basis,

$$\lambda = 0 \tag{5.9.2}$$

should follow from

$$g_*{}^*\lambda = e_*{}^*A_*{}^*\lambda = 0$$

Because $\overline{\overline{e}}$ is a basis we get

$$A_*{}^*\lambda = 0 \tag{5.9.3}$$

According to (5.9.1) $\operatorname{rank}_{*}{}^{*} A \le m$ and system (5.9.3) has more variables then equations. According to the theorem 5.7.4, $\lambda \ne 0$. This contradicts statement (5.9.2). Therefore, statement $m < n$ is not valid.

In the same manner we can prove that the statement $n < m$ is not valid. This completes the proof of the theorem. □

Definition 5.9.2. We call **dimension of ${}_*{}^*D$-vector space** the number of vectors in a basis □

Theorem 5.9.3. *The coordinate matrix of ${}_*{}^*D$-basis $\overline{\overline{g}}$ relative ${}_*{}^*D$-basis $\overline{\overline{e}}$ of vector space $V$ is ${}_*{}^*$-nonsingular matrix.*

Proof. According to theorem 5.6.6 ${}_*{}^*D$-rank of the coordinate matrix of basis $\overline{\overline{g}}$ relative basis $\overline{\overline{e}}$ equal to the dimension of vector space. This proves the statement of the theorem. □

Definition 5.9.4. We call one-to-one map

$$A : V \to W$$

**isomorphism of ${}_*{}^*D$-vector spaces** if this map is a linear map of ${}_*{}^*D$-vector spaces. □

Definition 5.9.5. **Automorphism of ${}_*{}^*D$-vector space** $V$ is isomorphism $A : V \to V$. □

Theorem 5.9.6. *Suppose that $\overline{\overline{f}}$ is a ${}_*{}^*D$-basis of vector space $V$. Then any automorphism $\overline{A}$ of ${}_*{}^*D$-vector space $V$ has form*

$$v' = A_*{}^*v \tag{5.9.4}$$

*where $A$ is a ${}_*{}^*$-nonsingular matrix.*

Proof. (5.9.4) follows from theorem 5.4.3. Because $\overline{A}$ is an isomorphism for each vector $v'$ exist one and only one vector $v$ such that $v' = v_*{}^*\overline{A}$. Therefore, system of ${}_*{}^*D$-linear equations (5.9.4) has a unique solution. According to corollary 5.7.5 matrix $A$ is a nonsingular matrix. □

Theorem 5.9.7. *Automorphisms of ${}_*{}^*D$-vector space form a group $GL(n, {}_*{}^*, D)$.*

Proof. If we have two automorphisms $\overline{A}$ and $\overline{B}$ then we can write

$$v' = A_*{}^*v$$

$$v'' = B_*{}^*v' = B_*{}^*A_*{}^*v$$

Therefore, the resulting automorphism has matrix $B_*{}^*A$. □

CHAPTER 6

# Basis Manifold

## 6.1. Linear Representation of Group

Let $V$ be ${}_*{}^*D$-vector space. We proved in the theorem 5.9.6 that we can identify any automorphism of a ${}_*{}^*D$-vector space $V$ with certain matrix. We proved in the theorem 5.9.7 that automorphisms of ${}_*{}^*D$-vector space form a group. When we study a representation in the ${}_*{}^*D$-vector space linear maps of the ${}_*{}^*D$-vector space are important for us.

Definition 6.1.1. Let $V$ be ${}_*{}^*D$-vector space.[6.1] The left-side representation[6.2] $f$ of group $G$ in ${}_*{}^*D$-vector space $V$ is called **linear $G*$-representation**. □

Theorem 6.1.2. *Automorphisms of ${}_*{}^*D$-vector space form a linear effective $GL(n, {}_*{}^*, D)*$-representation.*

Proof. If we have two automorphisms $\overline{A}$ and $\overline{B}$ then we can write

$$v' = A_*{}^*v$$
$$v'' = B_*{}^*v' = B_*{}^*A_*{}^*v$$

Therefore, the resulting automorphism has matrix $B_*{}^*A$.

It remains to prove that the kernel of inefficiency consists only of identity. Identity transformation satisfies to equation

$${}^iv = {}^iA_j\, {}^jv$$

Choosing values of coordinates as ${}^iv = {}^i\delta_k$ where we selected $k$ we get

$${}^i\delta_k = {}^iA_j\, {}^j\delta_k \tag{6.1.1}$$

From (6.1.1) it follows

$${}^i\delta_k = {}^iA_k$$

Since $k$ is arbitrary, we get the conclusion $A = \delta$. □

Theorem 6.1.3. *Let $\overline{\overline{f}}$, $\overline{\overline{g}}$ be bases of ${}_*{}^*D$-vector space $\overline{V}$. Let $\overline{\overline{e}}$, $\overline{\overline{h}}$ be bases of ${}_*{}^*D$-vector space $\overline{W}$. Let $A_1$ be matrix of linear map*

$$A : V \to W \tag{6.1.2}$$

*relative to bases $\overline{\overline{f}}$ and $\overline{\overline{e}}$ and $A_2$ be matrix of linear map (6.1.2) relative to bases $\overline{\overline{g}}$ and $\overline{\overline{h}}$. Suppose the basis $\overline{\overline{f}}$ has coordinate matrix $B$ relative the basis $\overline{\overline{g}}$*

$$\overline{\overline{f}} = \overline{\overline{g}}_*{}^*B$$

[6.1] Studying representation of the group in the ${}_*{}^*D$-vector space we follow to the agreement described in the remark 2.2.15.

[6.2] According ton the theorem 5.9.6, linear map of the right vector space is left-side transformation.

*and $\overline{\overline{e}}$ has coordinate matrix $C$ relative the basis $\overline{\overline{h}}$*

$$\overline{\overline{e}} = \overline{\overline{h}}_*{}^*C \tag{6.1.3}$$

*Then there is relationship between matrices $A_1$ and $A_2$*

$$A_1 = C^{-1_*{}^*}{}_*{}^*A_2{}_*{}^*B \tag{6.1.4}$$

Proof. Vector $\overline{a} \in V$ has expansion

$$\overline{a} = f_*{}^*a = g_*{}^*B_*{}^*a$$

relative to bases $\overline{\overline{f}}$ and $\overline{\overline{g}}$. Since $A$ is linear map, we can write it as

$$\overline{b} = e_*{}^*A_1{}_*{}^*a \tag{6.1.5}$$

relative to bases $\overline{\overline{f}}$ and $\overline{\overline{e}}$ and as

$$\overline{b} = h_*{}^*A_2{}_*{}^*B_*{}^*a \tag{6.1.6}$$

relative to bases $\overline{\overline{g}}$ and $\overline{\overline{h}}$. Since theorem 5.9.3 matrix $C$ has ${}_*{}^*$-inverse and from equation (6.1.3) it follows that

$$\overline{\overline{h}} = \overline{\overline{e}}_*{}^*C^{-1_*{}^*} \tag{6.1.7}$$

Substituting (6.1.7) into equation (6.1.6) we get

$$\overline{b} = e_*{}^*C^{-1_*{}^*}{}_*{}^*A_2{}_*{}^*B_*{}^*a \tag{6.1.8}$$

From theorem 5.3.3 and comparison of equations (6.1.5) and (6.1.8) it follows that

$$A_1{}_*{}^*a = C^{-1_*{}^*}{}_*{}^*A_2{}_*{}^*B_*{}^*a \tag{6.1.9}$$

Since vector $a$ is arbitrary vector, from theorem 6.1.2 and equation (6.1.9) statement of theorem follows. □

Theorem 6.1.4. *Let $\overline{A}$ be automorphism of ${}_*{}^*D$-vector space. Let $A_1$ be matrix of this automorphism defined relative to basis $\overline{\overline{f}}$ and $A_2$ be matrix of the same automorphism defined relative to basis $\overline{\overline{g}}$. Suppose the basis $\overline{\overline{f}}$ has coordinate matrix $B$ relative the basis $\overline{\overline{g}}$*

$$\overline{\overline{f}} = \overline{\overline{g}}_*{}^*B$$

*Then there is relationship between matrices $A_1$ and $A_2$*

$$A_1 = B^{-1_*{}^*}{}_*{}^*A_2{}_*{}^*B$$

Proof. Statement follows from theorem 6.1.3, because in this case $C = B$. □

## 6.2. Basis Manifold for ${}_*{}^*D$-Vector Space

Theorem 6.2.1. *Automorphism $A$ acting on each vector of basis of ${}_*{}^*D$-vector space maps a basis into another basis.*

Proof. Let $\overline{\overline{e}}$ be basis of ${}_*{}^*D$-vector space $V$. According to theorem 5.9.6, vector $e_a$ maps into a vector $e'_a$

$$e'_a = A_*{}^*e_a \tag{6.2.1}$$

Suppose vectors $e'_a$ are linearly dependent. Then $\lambda \neq 0$ in equation

$$e'{}_*{}^*\lambda = 0 \tag{6.2.2}$$

From equations (6.2.1) and (6.2.2) it follows that

$$A^{-1_*{}^*}{}_*{}^*e'{}_*{}^*\lambda = e_*{}^*\lambda = 0$$

and $\lambda \neq 0$. This contradicts to the statement that vectors $e_a$ are linearly independent. Therefore vectors $e'_a$ are linearly independent and form basis. □

Thus we can extend a linear $GL(n, {}_*{}^*, D)*$-representation to the set of bases. Transformation of this left-side representation is called **active transformation** because the linear map of the vector space induced this transformation ([5]). Since a linear operation is not defined on the basis manifold the active transformation is not a linear transformation. According to definition we write the action of the transformation $A \in GL(n, {}_*{}^*, D)$ on the basis $\overline{\overline{e}}$ as $A_*{}^*\overline{\overline{e}}$. Homomorphism of the group $G$ into the group $GL(n, {}_*{}^*, D)$ of active transformations is called **active $*G$-representation**.

THEOREM 6.2.2. *Active $GL(n, {}_*{}^*, D)*$-representation on the set of bases is single transitive representation.*

PROOF. To prove this theorem it is sufficient to show that at least one transformation of left-side representation is defined for any two bases and this transformation is unique. Homomorphism $A$ operating on basis $\overline{\overline{e}}$ has form

$$g_i = A_*{}^*e_i$$

where $g_i$ is coordinate matrix of vector $\overline{g}_i$ and $e_i$ is coordinate matrix of vector $\overline{e}_i$ relative basis $\overline{\overline{h}}$. Therefore, coordinate matrix of image of basis equal to ${}_*{}^*$-product of coordinate matrix of original basis over matrix of automorphism

$$g = A_*{}^*e$$

Since the theorem 5.9.3, matrices $g$ and $e$ are nonsingular. Therefore, matrix

$$A = g_*{}^*e^{-1_*{}^*}$$

is the matrix of automorphism mapping basis $\overline{\overline{e}}$ to basis $\overline{\overline{g}}$.

Suppose elements $g_1$, $g_2$ of group $G$ and basis $\overline{\overline{e}}$ satisfy equation

$$g_{1*}{}^*e = g_{2*}{}^*e \tag{6.2.3}$$

Since theorems 5.9.3 and 2.2.16 we get $g_1 = g_2$. This proves statement of theorem. □

Let us define an additional structure on vector space $V$. Then not every linear map keeps properties of the selected structure. In this case we need subgroup $G$ of the group $GL(n, {}_*{}^*, D)$ such that subgroup $G$ generates linear maps which hold properties of the selected structure. We usually call group $G$ **symmetry group**. Without loss of generality we identify element $g$ of group $G$ with corresponding transformation of representation and write its action on vector $v \in V$ as $g_*{}^*v$.

Not every two bases can be mapped by a transformation from the symmetry group because not every nonsingular linear transformation belongs to the representation of group $G$. Therefore, we can represent the set of bases as a union of orbits of group $G$.

DEFINITION 6.2.3. We call orbit $G_*{}^*\overline{\overline{e}}$ of the selected basis $\overline{\overline{e}}$ the **basis manifold $\mathcal{B}(\overline{V}, G)$ of ${}_*{}^*D$-vector space $V$**. □

THEOREM 6.2.4. *Active $G*$-representation on basis manifold is single transitive representation.*

PROOF. This is corollary of theorem 6.2.2 and definition 6.2.3. □

Theorem 6.2.4 means that the basis manifold $\mathcal{B}(V,G)$ is a homogenous space of group $G$. According to theorem 4.2.13 $*G$-representation, commuting with active, exists on the basis manifold. As we see from remark 4.2.14 transformation of $*G$-representation is different from an active transformation and cannot be reduced to transformation of space $V$. To emphasize the difference this transformation is called a **passive transformation** of basis manifold $\mathcal{B}(\overline{V},G)$ and the $*G$-representation is called **passive $*G$-representation**. According to the definition we write the passive transformation, defined by element $A\in G$, of basis $\overline{\overline{e}}$ as $\overline{\overline{e}}{}_*{}^*A$.

REMARK 6.2.5. I show samples of active and passive representations in table 6.2.1. □

TABLE 6.2.1. Active and Passive Representations

| Vector space | Group of Representation | Active Representation | Passive Representation |
|---|---|---|---|
| ${}_*{}^*D$-vector space | $GL(n,{}_*{}^*,D)$ | left-side | right-side |
| ${}^*{}_*D$-vector space | $GL(n,{}^*{}_*,D)$ | left-side | right-side |
| $D_*{}^*$-vector space | $GL(n,{}_*{}^*,D)$ | right-side | left-side |
| $D^*{}_*$-vector space | $GL(n,{}^*{}_*,D)$ | right-side | left-side |

According to the theorem 4.2.7 we can introduce on $\mathcal{B}(V,G)$ two types of coordinates defined on group $G$. Since we defined two representations of group $G$ on $\mathcal{B}(V,G)$, we use passive $*G$-representation to define coordinates. Our choice is based on the following theorem.

THEOREM 6.2.6. *The coordinate matrix of basis $\overline{\overline{g}}$ relative basis $\overline{\overline{e}}$ of ${}_*{}^*D$-vector space $V$ is identical with the matrix of passive transformation mapping basis $\overline{\overline{e}}$ to basis $\overline{\overline{g}}$.*

PROOF. According to model described in the example 5.3.9, the coordinate matrix of basis $\overline{\overline{g}}$ relative basis $\overline{\overline{e}}$ consist from ${}^*$-rows which are coordinate matrices of vectors $g_i$ relative the basis $\overline{\overline{e}}$. Therefore,

$$\overline{g}_i=\overline{e}_*{}^*g_i \tag{6.2.4}$$

At the same time the passive transformation $A$ mapping one basis to another has a form

$$\overline{g}_i=\overline{e}_*{}^*A_i \tag{6.2.5}$$

Since theorem 5.3.3,

$$g_i=A_i$$

for any $i$. This proves the theorem. □

Coordinates of representation are called **standard coordinates of basis**. This point of view allows introduction of two types of coordinates for element $g$ of group $G$. We can either use coordinates defined on the group, or introduce coordinates as elements of the matrix of the corresponding transformation. The former type of coordinates is more effective when we study properties of group $G$. The latter type of coordinates contains redundant data; however, it may be more convenient when we study representation of group $G$. The latter type of coordinates is called **coordinates of representation**.

## 6.3. Geometric Object of ${}_*{}^*D$-Vector Space

An active transformation changes bases and vectors uniformly and coordinates of vector relative basis do not change. A passive transformation changes only the basis and it leads to change of coordinates of vector relative to basis.

Let passive transformation $A \in G$ maps basis $\overline{\overline{e}} \in \mathcal{B}(V,G)$ into basis $\overline{\overline{e}}' \in \mathcal{B}(V,G)$

$$e' = e_*{}^*A \tag{6.3.1}$$

Let vector $v \in V$ has expansion

$$\overline{v} = e_*{}^*v \tag{6.3.2}$$

relative to basis $\overline{\overline{e}}$ and has expansion

$$\overline{v} = e'{}_*{}^*v' \tag{6.3.3}$$

relative to basis $\overline{\overline{e}}'$. From (6.3.1) and (6.3.3), it follows that

$$\overline{v} = e_*{}^*A_*{}^*v' \tag{6.3.4}$$

Comparing (6.3.2) and (6.3.4) we get

$$v = A_*{}^*v' \tag{6.3.5}$$

Because $A$ is ${}_*{}^*$-nonsingular matrix we get from (6.3.5)

$$v' = A^{-1_*{}^*}{}_*{}^*v \tag{6.3.6}$$

Coordinate transformation (6.3.6) does not depend on vector $\overline{v}$ or basis $\overline{\overline{e}}$, but is defined only by coordinates of vector $\overline{v}$ relative to basis $\overline{\overline{e}}$.

THEOREM 6.3.1. *Coordinate transformations (6.3.6) form effective linear $G*$-representation which is called* **coordinate representation in ${}_*{}^*D$-vector space**.

PROOF. Suppose we have two consecutive passive transformations $A$ and $B$. Coordinate transformation (6.3.6) corresponds to passive transformation $A$. Coordinate transformation

$$v'' = B^{-1_*{}^*}{}_*{}^*v' \tag{6.3.7}$$

corresponds to passive transformation $B$. Product of coordinate transformations (6.3.6) and (6.3.7) has form

$$v'' = B^{-1_*{}^*}{}_*{}^*A^{-1_*{}^*}{}_*{}^*v = (A_*{}^*B)^{-1_*{}^*}{}_*{}^*v \tag{6.3.8}$$

and is coordinate transformation corresponding to passive transformation $A_*{}^*B$. It proves that coordinate transformations form linear $G*$-representation.

Suppose coordinate transformation does not change vectors $\delta_k$. Then unit of group $G$ corresponds to it because representation is single transitive. Therefore, coordinate representation is effective. □

Let map of group $G$ to the group of passive transformations of ${}_*{}^*D$-vector space $N$ be coordinated with symmetry group of ${}_*{}^*D$-vector space $V$.[6.3] This means that passive transformation $A(a)$ of ${}_*{}^*D$-vector space $N$ corresponds to passive transformation $a$ of ${}_*{}^*D$-vector space $V$.

$$e'_N = e_N{}_*{}^*A(a) \tag{6.3.9}$$

Then coordinate transformation in $N$ gets form

$$w' = A(a^{-1_*{}^*}){}_*{}^*w = A(a)^{-1_*{}^*}{}_*{}^*w \tag{6.3.10}$$

Definition 6.3.2. Orbit

$$(A(G)^{-1_*{}^*}{}_*{}^*w, \overline{\overline{e}}_V{}_*{}^*G)$$

is called **geometric object in coordinate representation defined in ${}_*{}^*D$-vector space $V$**. For any basis $\overline{\overline{e}}'_V = \overline{\overline{e}}_V{}_*{}^*A$ corresponding point (6.3.10) of orbit defines **coordinates of geometric object in coordinate ${}_*{}^*D$-vector space** relative basis $\overline{\overline{e}}'_V$. □

Definition 6.3.3. Orbit

$$(A(G)^{-1_*{}^*}{}_*{}^*w, \overline{\overline{e}}_N{}_*{}^*A(G), \overline{\overline{e}}_V{}_*{}^*G)$$

is called **geometric object defined in ${}_*{}^*D$-vector space** $V$. For any basis $\overline{\overline{e}}'_V = a_*{}^*\overline{\overline{e}}_V$ corresponding point (6.3.10) of orbit defines **coordinates of a geometric object in ${}_*{}^*D$-vector space** relative to basis $\overline{\overline{e}}'_V$ and the corresponding vector

$$\overline{w} = e'_N{}_*{}^*w'$$

is called **representative of geometric object in $D_*{}^*$-vector space** $V$ in basis $\overline{\overline{e}}'_V$. □

Since a geometric object is an orbit of representation, we see that according to theorem 4.1.9 the definition of the geometric object is a proper definition.

We also say that $\overline{w}$ is a **geometric object of type** $A$

Definition 6.3.2 introduces a geometric object in coordinate space. We assume in definition 6.3.3 that we selected a basis of vector space $W$. This allows using a representative of the geometric object instead of its coordinates.

The question how large a diversity of geometric objects is well studied in case of vector spaces. However it is not such obvious in case of ${}_*{}^*D$-vector spaces. As we can see from table 6.2.1 $D_*{}^*$-vector space and ${}_*{}^*D$-vector space have common symmetry group $GL(n, {}_*{}^*, D)$. This allows study a geometric object in $D_*{}^*$-vector space when we study passive representation in ${}_*{}^*D$-vector space. Can we study the same time a geometric object in ${}^*{}_*D$-vector space. At first glance, the answer is negative based theorem 5.8.9. However, equation 2.2.6 determines required linear map between $GL(n, {}_*{}^*, D)$ and $GL(n, {}^*{}_*, D)$.

Theorem 6.3.4 (invariance principle). *Representative of geometric object does not depend on selection of basis $\overline{\overline{e}}'_V$.*

Proof. To define representative of geometric object, we need to select basis $\overline{\overline{e}}_V$, basis $\overline{\overline{e}}_W$ and coordinates of geometric object $w^\alpha$. Corresponding representative of geometric object has form

$$\overline{w} = e_W{}_*{}^*w$$

[6.3] We use the same notation for type of vector space for vector spaces $N$ and $V$. However their type may be different.

Suppose we map basis $\overline{\overline{e}}_V$ to basis $\overline{\overline{e}}'_V$ by passive transformation

$$e'_V = e_V{}_*{}^*A$$

According construction this generates passive transformation (6.3.9) and coordinate transformation (6.3.10). Corresponding representative of geometric object has form

$$\overline{w}' = e'_W{}_*{}^*w' = e_W{}_*{}^*A(a)_*{}^*A(a)^{-1}{}_*{}^*w = e_W{}_*{}^*w = \overline{w}$$

Therefore representative of geometric object is invariant relative selection of basis. □

Definition 6.3.5. Let

$$\overline{w}_1 = e_W{}_*{}^*w_1$$
$$\overline{w}_2 = e_W{}_*{}^*w_2$$

be geometric objects of the same type defined in ${}_*{}^*D$-vector space $V$. Geometric object

$$\overline{w} = e_W{}_*{}^*(w_1 + w_2)$$

is called **sum**

$$\overline{w} = \overline{w}_1 + \overline{w}_2$$

**of geometric objects** $\overline{w}_1$ and $\overline{w}_2$. □

Definition 6.3.6. Let

$$\overline{w}_1 = e_W{}_*{}^*w_1$$

be geometric object defined in ${}_*{}^*D$-vector space $V$. Geometric object

$$\overline{w}_2 = e_W{}_*{}^*(kw_1)$$

is called **product**

$$\overline{w}_2 = k\overline{w}_1$$

**of geometric object** $\overline{w}_1$ **and constant** $k \in D$. □

Theorem 6.3.7. *Geometric objects of type $A$ defined in ${}_*{}^*D$-vector space $V$ form ${}_*{}^*D$-vector space.*

Proof. The statement of the theorem follows from immediate verification of the properties of vector space. □

CHAPTER 7

# Linear Map of ${}_*{}^*D$-Vector Space

## 7.1. Linear Map of ${}_*{}^*D$-Vector Space

In this subsection we assume $V$, $W$ are ${}_*{}^*D$-vector spaces.

Definition 7.1.1. Let us denote by $\mathcal{L}({}_*{}^*D;V;W)$ set of linear maps

$$A:V\rightarrow W$$

of ${}_*{}^*D$-vector space $V$ into ${}_*{}^*D$-vector space $W$. Let us denote by $\mathcal{L}(D_*{}^*;V;W)$ set of linear maps

$$A:V\rightarrow W$$

of $D_*{}^*$-vector space $V$ into $D_*{}^*$-vector space $W$. □

We can consider division ring $D$ as ${}_*{}^*D$-vector space of dimension 1. Correspondingly we can consider set $\mathcal{L}(D_*{}^*;D;W)$ and $\mathcal{L}({}_*{}^*D;V;D)$.

Theorem 7.1.2. *Suppose $V$, $W$ are ${}_*{}^*D$-vector spaces. Then set $\mathcal{L}({}_*{}^*D;V;W)$ is an Abelian group relative composition law*

$$(A+B)_*{}^*x=A_*{}^*x+B_*{}^*x \tag{7.1.1}$$

*Linear map $A+B$ is called* **sum of maps** *$A$ and $B$.*

Proof. We need to show that map

$$A+B:V\rightarrow W$$

defined by equation (7.1.1) is linear map of ${}_*{}^*D$-vector spaces. According to the definition 5.4.2

$$A_*{}^*(x_*{}^*a)=(A_*{}^*x)_*{}^*a$$
$$B_*{}^*(x_*{}^*a)=(B_*{}^*x)_*{}^*a$$

We see that

$$\begin{aligned}(A+B)_*{}^*(x_*{}^*a)&=A_*{}^*(x_*{}^*a)+B_*{}^*(x_*{}^*a)\\&=(A_*{}^*x)_*{}^*a+(B_*{}^*x)_*{}^*a\\&=(A_*{}^*x+B_*{}^*x)_*{}^*a\\&=((A+B)_*{}^*x)_*{}^*a\end{aligned}$$

We need to show also that this operation is commutative.

$$\begin{aligned}(A+B)_*{}^*x&=A_*{}^*x+B_*{}^*x\\&=B_*{}^*x+A_*{}^*x\\&=(B+A)_*{}^*x\end{aligned}$$

□

THEOREM 7.1.3. *Let* $\overline{\overline{e}}_V=(e_{V\cdot i}, i\in I)$ *be a basis in* ${}_*{}^*D$*-vector space* $V$ *and* $\overline{\overline{e}}_W=(e_{W\cdot j}, j\in J)$ *be a basis in* ${}_*{}^*D$*-vector space* $W$*. Let* $A=({}^jA_i)$, $i\in I$, $j\in J$, *be arbitrary matrix. Then map*

$$\overline{A}:V\to W$$

*defined by equation*

$$\overline{a}=e_V{}_*{}^*a\rightarrow\overline{A}_*{}^*\overline{a}=e_W{}_*{}^*A_*{}^*a \tag{7.1.2}$$

*relative to selected bases is a linear map of* ${}_*{}^*D$*-vector spaces.*

PROOF. Theorem 7.1.3 is inverse statement to theorem 5.4.3. From the equation (7.1.2), it follows that

$$\overline{A}_*{}^*(\overline{v}_*{}^*a)=e_W{}_*{}^*A_*{}^*v_*{}^*a=(\overline{A}_*{}^*\overline{v})_*{}^*a$$

□

THEOREM 7.1.4. *Let* $\overline{\overline{e}}_V=(e_{V\cdot i}, i\in I)$ *be a basis in* ${}_*{}^*D$*-vector space* $V$ *and* $\overline{\overline{e}}_W=(e_{W\cdot j}, j\in J)$ *be a basis in* ${}_*{}^*D$*-vector space* $W$*. Suppose linear map* $\overline{A}$ *has matrix* $A=({}^jA_i)$, $i\in I$, $j\in J$, *relative to selected bases. Let* $m\in D$*. Then matrix*

$${}^j(mA)_i=m\ {}^jA_i$$

*defines linear map*

$$m\overline{A}:V\to W$$

*which we call* **left-side product of map** $A$ **over scalar**.

PROOF. The statement of the theorem is corollary of the theorem 7.1.3. □

THEOREM 7.1.5. *Set* $\mathcal{L}({}_*{}^*D;V;W)$ *is* $D_*{}^*$*-vector space.*

PROOF. Theorem 7.1.2 states that $\mathcal{L}({}_*{}^*D;V;W)$ is an Abelian group. It follows from theorem 7.1.4 that element of division ring $D$ defines left-side transformation on the Abelian group $\mathcal{L}({}_*{}^*D;V;W)$. From theorems 7.1.3, 5.1.1, and 5.1.3 it follows that set $\mathcal{L}({}_*{}^*D;V;W)$ is $D\star$-vector space.

Writing elements of basis $D\star$-vector space $\mathcal{L}({}_*{}^*D;V;W)$ as ${}_*$-rows or ${}^*$-rows, we represent $D\star$-vector space $\mathcal{L}({}_*{}^*D;V;W)$ as $D^*{}_*$- or ${}_*{}^*D$-vector space. I want to stress that choice between $D_*{}^*$- and $D^*{}_*$-linear combination in $D\star$-vector space $\mathcal{L}({}_*{}^*D;V;W)$ does not depend on type of vector spaces $V$ and $W$.

To select the type of vector space $\mathcal{L}({}_*{}^*D;V;W)$ I draw attention to the following observation. let $V$ and $W$ be ${}_*{}^*D$-vector spaces. Suppose $\mathcal{L}({}_*{}^*D;V;W)$ is $D_*{}^*$-vector space. Then we can represent the operation of ${}^*$-row of ${}_*{}^*D$-linear maps ${}^jA$ on ${}_*$-row of vectors $f_i$ as matrix

$$\begin{pmatrix}{}^1A\\ \dots\\ {}^mA\end{pmatrix}{}_*{}^*\begin{pmatrix}f_1 & \dots & f_n\end{pmatrix}=\begin{pmatrix}{}^1A_*{}^*f_1 & \dots & {}^1A_*{}^*f_n\\ \dots & \dots & \dots\\ {}^mA_*{}^*f_1 & \dots & {}^mA_*{}^*f_n\end{pmatrix}$$

This notation is coordinated with matrix notation of action of $D_*{}^*$-linear combination $a_*{}^*A$ of ${}_*{}^*D$-linear maps $A$. □

We can also define right-side product of ${}_*{}^*D$-linear map $A$ over scalar. However in general we cannot carry this left-side representation of division ring $D$ in $D\star$-vector space $\mathcal{L}({}_*{}^*D;V;W)$ into ${}_*{}^*D$-vector space $W$. Indeed, in case of right-side product we get

$$(Am){}_*{}^*v=(Am){}_*{}^*f{}_*{}^*v=e{}_*{}^*(Am){}_*{}^*v$$

Since product in division ring is noncommutative, we cannot express this expression as product of $A{}_*{}^*v$ over $m$.

Ambiguity of notation

$$m{}_*{}^*A{}_*{}^*v$$

is corollary of theorem 7.1.5. We may assume that meaning of this notation is clear from the text. To make notation more clear we will use brackets. Expression

$$w=[m{}_*{}^*A]{}_*{}^*v$$

means that $D{}_*{}^*$-linear composition of ${}_*{}^*D$-linear maps ${}^iA$ maps the vector $v$ to the vector $w$. Expression

$$w=B{}_*{}^*A{}_*{}^*v$$

means that ${}_*{}^*$-product of linear maps $A$ and $B$ maps the vector $v$ to the vector $w$.

## 7.2. 1-Form on ${}_*{}^*D$-Vector Space

Definition 7.2.1. **1-form** on ${}_*{}^*D$-vector space $V$ is linear map

$$b:V\to D \tag{7.2.1}$$

□

We can write value of 1-form $b$, defined for vector $a$, as

$$b(a)=<b,a>$$

Theorem 7.2.2. *Set $\mathcal{L}({}_*{}^*D;V;D)$ is $D\star$-vector space.*

Proof. ${}_*{}^*D$-vector space of dimension 1 is equivalent to division ring $D$ □

Theorem 7.2.3. *Let $\overline{\overline{e}}$ be a basis in ${}_*{}^*D$-vector space $V$. 1-form $\overline{b}$ has presentation*

$$<\overline{b},\overline{a}>=b{}_*{}^*a \tag{7.2.2}$$

*relative to selected basis, where vector $a$ has expansion*

$$\overline{a}=e{}_*{}^*a \tag{7.2.3}$$

*and*

$$b_i=<\overline{b},e_i> \tag{7.2.4}$$

Proof. Because $\overline{b}$ is 1-form, it follows from (7.2.3) that

$$<\overline{b},\overline{a}>=<\overline{b},e{}_*{}^*a>=<\overline{b},e_i>\ {}^ia \tag{7.2.5}$$

(7.2.2) follows from (7.2.5) and (7.2.4). □

Theorem 7.2.4. *Let $\overline{\overline{e}}$ be a basis in ${}_*{}^*D$-vector space $V$. 1-form*

$$\overline{b}:V\to D$$

*is uniquely defined by values (7.2.4) into which 1-form $\overline{b}$ maps vectors of basis.*

Proof. Statement follows from theorems 7.2.3 and 5.3.3. □

THEOREM 7.2.5. *Let $\overline{\overline{e}}$ be a basis in ${}_*{}^*D$-vector space $V$. The set of 1-forms ${}^id$ such that*

$$< {}^jd, e_i >= {}^j\delta_i \tag{7.2.6}$$

*is basis $\overline{\overline{d}}$ of $D_*{}^*$-vector space $\mathcal{L}({}_*{}^*D; V; D)$.*

PROOF. Since we assume $b_i = {}^j\delta_i$ , then according to theorem 7.2.4 there exists the 1-form ${}^jd = \overline{b}$ for given $j$. If we assume that there exists 1-form

$$\overline{c} = c_*{}^*d = 0$$

then

$$c_j < {}^jd, e_i >= 0$$

According to equation (7.2.6),

$$c_i = c_j\ {}^j\delta_i = 0$$

Therefore, 1-forms ${}^jd$ are linearly independent. □

DEFINITION 7.2.6. Let $V$ be ${}_*{}^*D$-vector space. $D_*{}^*$-vector space

$$V^* = \mathcal{L}({}_*{}^*D; V; D)$$

is called **dual space of ${}_*{}^*D$-vector space** $V$. Let $\overline{\overline{e}}$ be a basis in ${}_*{}^*D$-vector space $V$. Basis $\overline{\overline{d}}$ of $D_*{}^*$-vector space $V^*$, satisfying to equation (7.2.6), is called **basis dual to basis** $\overline{\overline{e}}$. □

THEOREM 7.2.7. *Let $A$ be passive transformation of basis manifold $\mathcal{B}(V, GL(n, {}_*{}^*, D))$. Let basis*

$$\overline{\overline{e}}' = \overline{\overline{e}}_*{}^*A \tag{7.2.7}$$

*be image of basis $\overline{\overline{e}}$. Let $B$ be passive transformation of basis manifold $\mathcal{B}(V^*, GL(n, {}_*{}^*, D))$ such, that basis*

$$\overline{\overline{d}}' = B_*{}^*\overline{\overline{d}} \tag{7.2.8}$$

*is dual to basis. Then*

$$B = A^{-1_*{}^*} \tag{7.2.9}$$

PROOF. From equations (7.2.6), (7.2.7), (7.2.8) it follows

$$\begin{aligned} {}^j\delta_i &= < {}^jd', e'_i > \\ &= {}^jB_l\ < {}^ld, e_k >\ {}^kA_i \\ &= {}^jB_l\ {}^l\delta_k\ {}^kA_i \\ &= {}^jB_k\ {}^kA_i \end{aligned} \tag{7.2.10}$$

Equation (7.2.9) follows from equation (7.2.10). □

## 7.3. Twin Representations of Division Ring

Theorem 7.3.1. *We can introduce structure of $D*$-vector space in any ${}_*{}^*D$-vector space defining right-side product of vector over scalar using equation*

$$mv = m\delta_*{}^*v$$

Proof. We verify directly that the map

$$f : D \to V^*$$

defined by equation

$$f(m) = m\delta$$

defines $D*$-representation. □

We can formulate the theorem 7.3.1 by other way.

Theorem 7.3.2. *Suppose we defined an effective right-side representation $f$ of division ring $D$ on the Abelian group $V$. Then we can uniquely define an effective left-side representation $h$ of division ring $D$ on the Abelian group $V$ such that diagram*

$$\begin{array}{ccc} V & \xrightarrow{h(a)} & V \\ \downarrow{\scriptstyle f(b)} & & \downarrow{\scriptstyle f(b)} \\ V & \xrightarrow[h(a)]{} & V \end{array}$$

*is commutative for any $a$, $b \in D$.* □

We call representations $f$ and $h$ **twin representations of the division ring** $D$.

Theorem 7.3.3. *In vector space $V$ over division ring $D$ we can define $D\star$-product and $\star D$-product of vector over scalar. According to theorem 7.3.2 these operations satisfy equation*

$$(am)b = a(mb) \tag{7.3.1}$$

Equation (7.3.1) represents **associative law for twin representations**. This allows us writing of such expressions without using of brackets.

Proof. In section 7.1 there is definition of right-side product of linear map $A$ of ${}_*{}^*D$-vector space over scalar. According to theorem 7.3.1, $D*$-representation in Abelian group $\mathcal{L}({}_*{}^*D; V; V)$ can be carried into ${}_*{}^*D$-vector space $V$ according to rule

$$[mA]_*{}^*v = [m\delta]_*{}^*(A_*{}^*v) = m(A_*{}^*v)$$

□

The analogy of the vector space over field goes so far that we can assume an existence of the concept of basis which serves for left-side and right-side product of vector over skalar.

Theorem 7.3.4. *Basis manifold of ${}_*{}^*D$-vector space $V$ and basis manifold of ${}^*{}_*D$-vector space $V$ are different*

$$\mathcal{B}(V, D_*{}^*) \neq \mathcal{B}(V, {}^*{}_*D)$$

PROOF. To prove this theorem we use the standard representation of a matrix. Without loss of generality we prove theorem in coordinate vector space $D^n$.

Let $\overline{\overline{e}} = (e_i = (\delta_i^j),\ i,\ j \in I,\ |I| = n)$ be the set of vectors of vector space $D^n$. $\overline{\overline{e}}$ is evidently basis of ${}_*{}^*D$-vector space and basis of ${}^*{}_*D$-vector space. For arbitrary set of vectors $(f_i, i \in I, |I| = n)$ ${}_*{}^*D$-coordinate matrix

$$f = \begin{pmatrix} f_1^1 & ... & f_1^n \\ ... & ... & ... \\ f_n^1 & ... & f_n^n \end{pmatrix} \tag{7.3.2}$$

relative to basis $\overline{\overline{e}}$ coincide with ${}^*{}_*D$-coordinate matrix relative to basis $\overline{\overline{e}}$.

Let the set of vectors $(f_i, i \in I, |I| = n)$ be basis of ${}_*{}^*D$-vector space. According to theorem 5.9.3 matrix (7.3.2) is ${}_*{}^*$-nonsingular matrix.

Let the set of vectors $(f_i, i \in I, |I| = n)$ be basis of ${}^*{}_*D$-vector space. According to theorem 5.9.3 matrix (7.3.2) is ${}^*{}_*$-nonsingular matrix.

Therefore, if the set of vectors $(f_i, i \in I, |I| = n)$ is basis ${}_*{}^*D$-of vector space and basis of ${}^*{}_*D$-vector space, their coordinate matrix (7.3.2) is ${}_*{}^*$-nonsingular and ${}^*{}_*$-nonsingular matrix. The statement follows from theorem 5.8.9. □

From theorem 7.3.4, it follows that there exists basis $\overline{\overline{e}}$ of $D_*{}^*$-vector space $V$ which is not basis of ${}^*{}_*D$-vector space $V$.

## 7.4. $D$-Vector Space

For many problems we may confine ourselves to considering of $D\star$-vector space or $\star D$-vector space. However there are problems where we forced to reject simple model and consider both structures of vector space at the same time.

Twin representations of division ring $D$ in Abelian group $V$ determine the structure of **$D$-vector space**. From this definition, it follows that if $v$, $w \in V$, then for any $a$, $b$, $c$, $d \in D$

$$avb + cwd \in V \tag{7.4.1}$$

The expression (7.4.1) is called **linear combination of vectors** of $D$-vector space $V$.

DEFINITION 7.4.1. Vectors $a_i$, $i \in I$, of $D$-vector space $V$ are **linearly independent** if $b^ic^i = 0$, $i = 1, ..., n$, follows from the equation

$$b^i a_i c^i = 0$$

Otherwise vectors $a_i$ are **linearly dependent**. □

DEFINITION 7.4.2. We call set of vectors

$$\overline{\overline{e}} = (e_i, i \in I)$$

**basis for $D$-vector space** if vectors $e_i$ are linearly independent and adding to this system any other vector we get a new system which is linearly dependent. □

THEOREM 7.4.3. *A basis of $D$-vector space $V$ is basis of $D^*{}_*$-vector space $V$.*

PROOF. Let vectors $a_i$, $i \in I$ are linear dependent in $D^*{}_*$-vector space. Then there exists ${}_*$-row

$$\begin{pmatrix} c^i & i \in I \end{pmatrix} \neq 0 \tag{7.4.2}$$

such that

$$c^{*}{}_{*}a = 0 \tag{7.4.3}$$

From the equation (7.4.3), it folows that

$$c^{i}a_{i}1 = 0 \tag{7.4.4}$$

According to the definition 7.4.1, vectors $a_i$ are linear dependent in $D$-vector space.

Let $D^{*}{}_{*}$-dimension of vector space $V$ is equal to $\boldsymbol{n}$. Let $\overline{\overline{e}}$ be basis of $D^{*}{}_{*}$-vector space. Let

$$\begin{pmatrix} 1 & ... & 0 \\ ... & ... & ... \\ 0 & ... & 1 \end{pmatrix}$$

be coordinate matrix of basis $\overline{\overline{e}}$. It is evident that we cannot present vector $e_{\boldsymbol{n}}$ as $D$-linear combination of vectors $e_{\boldsymbol{1}}$, ..., $e_{\boldsymbol{n-1}}$. □

Whether the set of vectors $(e_{\boldsymbol{i}}, \boldsymbol{i} \in \boldsymbol{I})$ is a basis of ${}_{*}{}^{*}D$-vector space or a basis of $D^{*}{}_{*}$-vector space is of no concern to us. So we will use index in standard format to represent coordinates of vector. Based on the theorem 7.4.3, considering a basis of $D$-vector space $V$ we will consider it as a basis of $D^{*}{}_{*}$-vector space $V$.

CHAPTER 8

# Product of Representations

## 8.1. Bimodule

DEFINITION 8.1.1. $\overline{V}$ is a **$(S\star, \star T)$-bimodule** if we define structures of a $S\star$-vector space and a $\star T$-vector space on the set $\overline{V}$.[8.1] □

We also use notation ${}_S\overline{V}_T$ when we want to tell that $\overline{V}$ is $(S, T)$-bimodule.

EXAMPLE 8.1.2. In section 7.4 we considered $D$-vector space where we can define the structure of $(D^*{}_*, {}_*{}^*D)$-bimodule. □

EXAMPLE 8.1.3. Set of $n\times m$ matrices builds up $(D\star, \star D)$-bimodule. We represent $D\star$-basis as set of matrices ${}_i{}^je = (\delta_i^l \delta_k^j)$. □

EXAMPLE 8.1.4. To see example of $(S_*{}^*, {}_*{}^*T)$-bimodule we can use set of matrices

$$\begin{pmatrix} (A_1, {}^1B) & ... & (A_1, {}^mB) \\ ... & ... & ... \\ (A_n, {}^1B) & ... & (A_n, {}^mB) \end{pmatrix} \tag{8.1.1}$$

Formally we can represent this matrix in form

$$\begin{pmatrix} (A_1, {}^1B) & ... & (A_1, {}^mB) \\ ... & ... & ... \\ (A_n, {}^1B) & ... & (A_n, {}^mB) \end{pmatrix} = \begin{pmatrix} A_1 \\ ... \\ A_n \end{pmatrix} {}_*{}^* \begin{pmatrix} {}^1B & ... & {}^mB \end{pmatrix}$$

This presentation clearly shows that $S_*{}^*$-dimension of $(S_*{}^*, {}_*{}^*T)$-bimodule is $n$ and ${}_*{}^*T$-dimension is $m$. However in contrast to example 8.1.2 we cannot build up $S_*{}^*$-basis or ${}_*{}^*T$-basis which generate $(S_*{}^*, {}_*{}^*T)$-bimodule. As a matter of fact, basis of $(S_*{}^*, {}_*{}^*T)$-bimodule has form $({}^a\overline{e}, \overline{f}_b)$; and any vector of $(S_*{}^*, {}_*{}^*T)$-bimodule has expansion

$$(\overline{A}, \overline{B}) = A_a({}^a\overline{e}, \overline{f}_b)^bB \tag{8.1.2}$$

Coefficients of this expansion generate matrix (8.1.1). □

We can write equation (8.1.2) as

$$(\overline{A}, \overline{B}) = (A_*{}^*\overline{e}, \overline{f}_*{}^*B) \tag{8.1.3}$$

Equation (8.1.3) looks unusual, however it becomes clear, if we assume

$$(\overline{A}, \overline{B}) = (\overline{A}, \overline{0}) + (\overline{0}, \overline{B})$$

[8.1] I recall, that we use notation $S\star$-vector space for left $S$-vector space and $\star T$-vector space for right $T$-vector space.

Convention 5.2.6 holds for representation of vectors of ${}_*{}^*T$-vector space is independent from representation of vectors of $S_*{}^*$-vector space. We write down this equation as

$$(\overline{A},\overline{B}) = (A^*{}_*\overline{e}, B_*{}^*\overline{f}) \tag{8.1.4}$$

Coordinates of vector $(\overline{A},\overline{B})$ do not generate matrix. However we do not care about this now.

We want to write components $A$ and $B$ in equation (8.1.2) on the one hand relative to vectors $({}^a\overline{e},\overline{f}_b)$. At first glance it appears impossible. However we accept convention to write this equation as

$$(\overline{A},\overline{B}) = (A_a, B^b)({}^a e, {}_b f) = (A_*{}^*, B^*{}_*)(e, f)$$

This structure can be generalized.

## 8.2. Direct Product of Division Rings

Definition 8.2.1. Let $\mathcal{A}$ be a category. Let $\{B_i, i\in I\}$ be the set of objects of $\mathcal{A}$. Object

$$P = \prod_{i\in I} B_i$$

and set of morphisms

$$\{\ f_i : P \longrightarrow B_i\ , i\in I\}$$

is called a **product of objects** $\{B_i, i\in I\}$ **in category** $\mathcal{A}$[8.2] if for any object $R$ and set of morphisms

$$\{\ g_i : R \longrightarrow B_i\ , i\in I\}$$

there exists a unique morphism

$$h : R \longrightarrow P$$

such that diagram

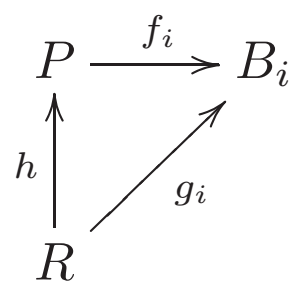

$$f_i\circ h = g_i$$

is commutative for all $i\in I$.

If $|I| = n$, then we also will use notation

$$G = \prod_{i=1}^{n} B_i = B_1\times ...\times B_n$$

for product of objects $\{B_i, i\in I\}$ in $\mathcal{A}$. □

Example 8.2.2. Let $\mathcal{S}$ be the category of sets.[8.3] According to the definition 8.2.1, Cartesian product

$$A = \prod_{i\in I} A_i$$

of family of sets $(A_i, i\in I)$ and family of projections on the $i$-th factor

$$p_i : A\to A_i$$

[8.2] I made definition according to [1], page 58.
[8.3] See also the example in [1], page 59.

are product in the category $\mathcal{S}$. □

THEOREM 8.2.3. *Let product exist in the category $\mathcal{A}$ of $\Omega$-algebras. Let $\Omega$-algebra $A$ and family of morphisms*

$$p_i : A \to A_i \quad i \in I$$

*be product in the category $\mathcal{A}$. Then*

8.2.3.1: *The set $A$ is Cartesian product of family of sets $(A_i, i \in I)$*

8.2.3.2: *The homomorphism of $\Omega$-algebra*

$$p_i : A \to A_i$$

*is projection on $i$-th factor.*

8.2.3.3: *We can represent any $A$-number $a$ as tuple $(p_i(a), i \in I)$ of $A_i$-numbers.*

8.2.3.4: *Let $\omega \in \Omega$ be $n$-ary operation. Then operation $\omega$ is defined componentwise*

$$a_1 ... a_n\omega = (a_{1i} ... a_{ni}\omega, i \in I) \tag{8.2.1}$$

*where $a_1 = (a_{1i}, i \in I)$, ..., $a_n = (a_{ni}, i \in I)$ .*

PROOF. Let

$$A = \prod_{i \in I} A_i$$

be Cartesian product of family of sets $(A_i, i \in I)$ and, for each $i \in I$, the map

$$p_i : A \to A_i$$

be projection on the $i$-th factor. Consider the diagram of morphisms in category of sets $\mathcal{S}$

(8.2.2)

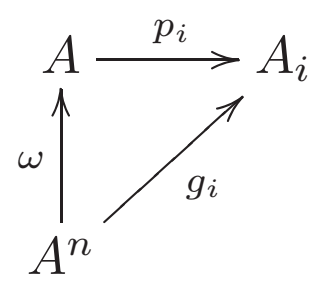

$$p_i \circ \omega = g_i$$

where the map $g_i$ is defined by the equation

$$g_i(a_1, ..., a_n) = p_i(a_1) ... p_i(a_n)\omega$$

According to the definition 8.2.1, the map $\omega$ is defined uniquely from the set of diagrams (8.2.2)

$$a_1 ... a_n\omega = (p_i(a_1) ... p_i(a_n)\omega, i \in I) \tag{8.2.3}$$

The equation (8.2.1) follows from the equation (8.2.3). □

DEFINITION 8.2.4. Let product exist in the category of $\Omega$-algebras. If $\Omega$-algebra $A$ and family of morphisms

$$p_i : A \to A_i \quad i \in I$$

is product in the category $\mathcal{A}$, then $\Omega$-algebra $A$ is called **direct** or **Cartesian product of $\Omega$-algebras** $(A_i, i \in I)$ . □

EXAMPLE 8.2.5. Let $\{G_i, i \in I\}$ be the set of groups. Let

$$G = \prod_{i \in I} G_i$$

be Cartesian product of the sets $G_i$, $i \in I$. We can define a group structure on $G$ by componentwise multiplication. If $x = (x_i, i \in I) \in D$ and $y = (y_i, i \in I) \in D$, we define their product

$$xy = (x_i y_i, i \in I)$$

If $x = (x_i, i \in I) \in D$, we define the inverse

$$x^{-1} = (x_i^{-1}, i \in I)$$

The set $G$ is called **Cartesian product of groups**[8.4] $G_i$, $i \in I$.

If $|I| = n$, then we also will use notation

$$G = \prod_{i=1}^{n} G_i = G_1 \times ... \times G_n$$

for Cartesian product of groups $G_1$, ..., $G_n$. □

THEOREM 8.2.6. *Cartesian product of Abelian groups is Abelian group.*

PROOF. Let $x = (x_i, i \in I) \in G$ and $y = (y_i, i \in I) \in G$. Then

$$x + y = (x_i + y_i, i \in I) = (y_i + x_i, i \in I) = y + x$$

□

EXAMPLE 8.2.7. Let $\{D_i, i \in I\}$ be the set of division rings. Let

$$D = \prod_{i \in I} D_i$$

be product of additive groups of division rings $D_i$, $i \in I$. If $x = (x_i, i \in I) \in D$ and $y = (y_i, i \in I) \in D$, we define their product componentwise

$$xy = (x_i y_i, i \in I)$$

The multiplication unit is $e = (e_i, i \in I) \in D$ where $e_i$, $i \in I$ is multiplication unit of $D_i$. The set $D$ is called **direct product of division rings**[8.5] $D_i$, $i \in I$.

If $|I| = n$, then we also will use notation

$$D = \prod_{i=1}^{n} D_i = D_1 \times ... \times D_n$$

for product of division rings $D_1$, ..., $D_n$. □

Direct product of division rings $D_i$, $i \in I$, in general, is not division ring. For instance, let $x_1 \in D_1$, $x_1 \neq 0$, $x_2 \in D_2$, $x_2 \neq 0$. Then

$$(x_1, 0)(0, x_2) = (0, 0)$$

However, direct product of division rings is ring. Therefore, direct product is not defined in category of division rings, however direct product of division rings is defined in category of rings.

[8.4] I made definition according to example from [1], page 9.

[8.5] I made definition according to proposition 1.1 from [1], p. 91

Theorem 8.2.8. *Let set $A$ be Cartesian product of sets $(A_i, i\in I)$ and set $B$ be Cartesian product of sets $(B_i, i\in I)$ . For each $i\in I$, let*

$$f_i:A_i\rightarrow B_i$$

*be the map from the set $A_i$ into the set $B_i$. For each $i\in I$, consider commutative diagram*

(8.2.4)

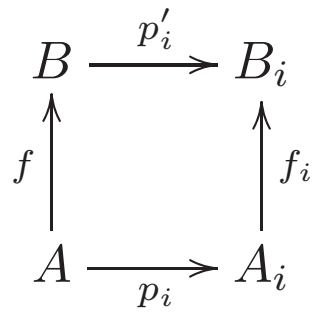

*where maps $p_i$, $p'_i$ are projection on the $i$-th factor. The set of commutative diagrams (8.2.4) uniquely defines map*

$$f:A\rightarrow B$$
$$f(a_i, i\in I)=(f_i(a_i), i\in I)$$

Proof. For each $i\in I$, consider commutative diagram

(8.2.5)

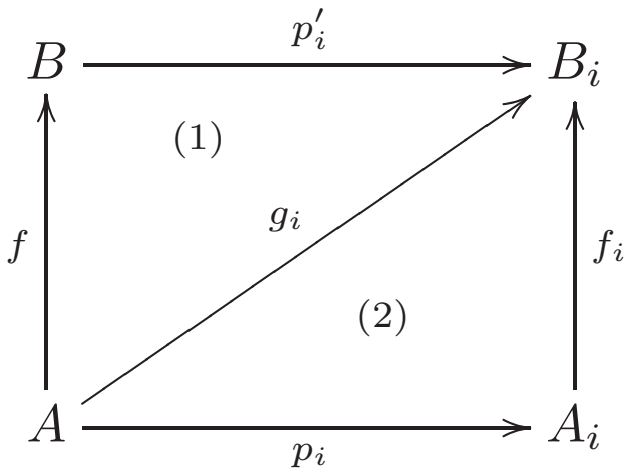

Let $a\in A$. According to the statement 8.2.3.3, we can represent $A$-number $a$ as tuple of $A_i$-numbers

$$a=(a_i, i\in I)\quad a_i=p_i(a)\in A_i \tag{8.2.6}$$

Let

$$b=f(a)\in B \tag{8.2.7}$$

According to the statement 8.2.3.3, we can represent $B$-number $b$ as tuple of $B_i$-numbers

$$b=(b_i, i\in I)\quad b_i=p'_i(b)\in B_i \tag{8.2.8}$$

From commutativity of diagram (1) and from equations (8.2.7), (8.2.8), it follows that

$$b_i=g_i(b) \tag{8.2.9}$$

From commutativity of diagram (2) and from the equation (8.2.6), it follows that

$$b_i=f_i(a_i)$$

□

Theorem 8.2.9. *Let $\Omega$-algebra $A$ be Cartesian product of $\Omega$-algebras $(A_i, i \in I)$ and $\Omega$-algebra $B$ be Cartesian product of $\Omega$-algebras $(B_i, i \in I)$ . For each $i \in I$, let the map*

$$f_i : A_i \to B_i$$

*be homomorphism of $\Omega$-algebra. Then the map*

$$f : A \to B$$

*defined by the equation*

$$f(a_i, i \in I) = (f_i(a_i), i \in I) \tag{8.2.10}$$

*is homomorphism of $\Omega$-algebra.*

Proof. Let $\omega \in \Omega$ be n-ary operation. Let $a_1 = (a_{1i}, i \in I)$, ..., $a_n = (a_{ni}, i \in I)$ , $b_1 = (b_{1i}, i \in I)$, ..., $b_n = (b_{ni}, i \in I)$ . From equations (8.2.1), (8.2.10), it follows that

$$\begin{aligned}
f(a_1...a_n\omega) &= f(a_{1i}...a_{ni}\omega, i \in I)\\
&= (f_i(a_{1i}...a_{ni}\omega), i \in I)\\
&= ((f_i(a_{1i}))...(f_i(a_{ni})), i \in I)\\
&= (b_{1i}...b_{ni}\omega, i \in I)
\end{aligned}$$

$$f(a_1)...f(a_n)\omega = b_1...b_n\omega \qquad = (b_{1i}...b_{ni}\omega, i \in I)$$

□

Remark 8.2.10. A theorem which is converse to theorems 8.2.8, 8.2.9 in general is not true. For instance, consider $R$-vector space $V$ of dimension 2 and linear map

$$f : V \to V$$

$$x'^1 = x^1$$
$$x'^2 = x^1 + x^2$$

Vector space $V$ is Abelian group which is Cartesian product of Abelian groups

$$V = R \times R$$

The map $f$ is homomorphism of Abelian group $V$. However, the corespondence

$$x^2 \to x'^2$$

depends on the value of $x^1$; therefore this corespondence is not a map. □

## 8.3. Direct Product of $D_*{}^*$-Vector Spaces

Lemma 8.3.1. Let

$$A = \prod_{i \in I} A_i$$

be Cartesian product of family of $\Omega_2$-algebras $(A_i, i \in I)$ . For each $i \in I$, let the set ${}^*A_i$ be $\Omega_2$-algebra. Then the set

$$^\circ A = \{f \in {}^*A : f(a_i, i \in I) = (f_i(a_i), i \in I)\} \tag{8.3.1}$$

is Cartesian product of $\Omega_2$-algebras ${}^*A_i$.

Proof. According to the definition (8.3.1), we can represent a map $f \in {}^\circ A$ as tuple

$$f = (f_i, i \in I)$$

of maps $f_i \in {}^* A_i$. According to the definition (8.3.1),

$$(f_i, i \in I)(a_i, i \in I) = (f_i(a_i), i \in I)$$

Let $\omega \in \Omega$ be n-ary operation. We define operation $\omega$ on the set ${}^\circ A$ using equation

$$((f_{1i}, i \in I)...(f_{ni}, i \in I)\omega)(a_i, i \in I) = ((f_{1i}(a_i))...(f_{ni}(a_i))\omega, i \in I)$$

□

Theorem 8.3.2. *Let category $\mathcal{A}_1$ of $\Omega_1$-algebras have product. Let category $\mathcal{A}_2$ of $\Omega_2$-algebras have product. Then in category*[8.6] $(\mathcal{A}_1 *)\mathcal{A}_2$ *there exists product of single transitive left-side representations of $\Omega_1$-algebra in $\Omega_2$-algebra.*

Proof. For $j = 1, 2$, let

$$P_j = \prod_{i \in I} B_{ji}$$

be product of family of $\Omega_j$-algebras $\{B_{ji}, i \in I\}$ and for any $i \in I$ the map

$$t_{ji} : P_j \longrightarrow B_{ji}$$

be projection onto factor $i$. For each $i \in I$, let

$$h_i : B_{1i} \xrightarrow{*} B_{2i}$$

be single transitive $B_{1i}*$-representation in $\Omega_2$-algebra $B_{2i}$.

Let $b_1 \in P_1$. According to the statement 8.2.3.3, we can represent $P_1$-number $b_1$ as tuple of $B_{1i}$-numbers

(8.3.2) $$b_1 = (b_{1i}, i \in I) \quad b_{1i} = t_{1i}(b_1) \in B_{1i}$$

Let $b_2 \in P_2$. According to the statement 8.2.3.3, we can represent $P_2$-number $b_2$ as tuple of $B_{2i}$-numbers

(8.3.3) $$b_2 = (b_{2i}, i \in I) \quad b_{2i} = t_{2i}(b_2) \in B_{2i}$$

Lemma 8.3.3. For each $i \in I$, consider diagram of maps

(8.3.4)

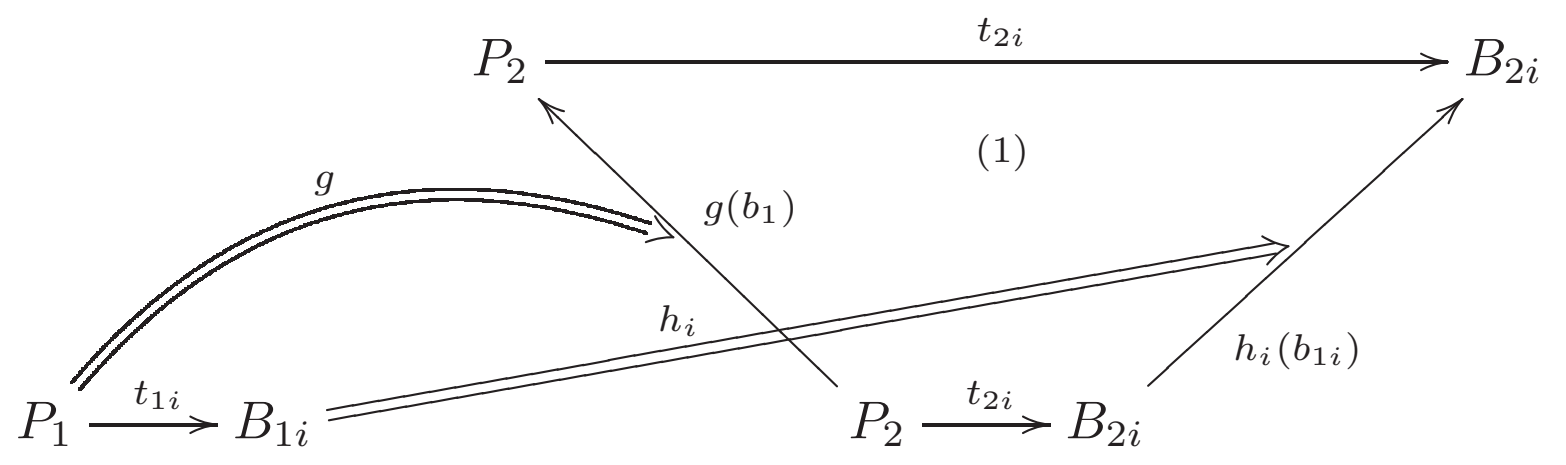

Let map

$$g : P_1 \rightarrow {}^* P_2$$

be defined by the equation

(8.3.5) $$g(b_1)(b_2) = (h_i(b_{1i})(b_{2i}), i \in I)$$

[8.6] See definition 3.2.13

Then the map $g$ is single transitive $P_1*$-representation in $\Omega_2$-algebra $P_2$

$$g : P_1 \xrightarrow{*} P_2$$

The map $(t_{1i}, t_{2i})$ is morphism of representation $g$ into representation $h_i$.

Proof.

8.3.3.1: According to definitions 3.1.1, 3.1.2, the map $h_i(b_{1i})$ is homomorphism of $\Omega_2$-algebra $B_{2i}$. According to the theorem 8.2.9, from commutativity of the diagram (1) for each $i \in I$, it follows that the map

$$g(b_1) : P_2 \to P_2$$

defined by the equation (8.3.5) is homomorphism of $\Omega_2$-algebra $P_2$.

8.3.3.2: According to the definition 3.1.2, the set ${}^*B_{2i}$ is $\Omega_1$-algebra. According to the lemma 8.3.1, the set ${}^\circ P_2 \subseteq {}^*P_2$ is $\Omega_1$-algebra.

8.3.3.3: According to the definition 3.1.2, the map

$$h_i : B_{1i} \to {}^*B_{2i}$$

is homomorphism of $\Omega_1$-algebra. According to the theorem 8.2.9, the map

$$g : P_1 \to {}^*P_2$$

defined by the equation

$$g(b_1) = (h_i(b_{1i}), i \in I) \tag{8.3.6}$$

is homomorphism of $\Omega_1$-algebra.

According to statements 8.3.3.1, 8.3.3.3 and to the definition 3.1.2, the map $g$ is $P_1*$-representation in $\Omega_2$-algebra $P_2$.

Let $b_{21}, b_{22} \in P_2$. According to the statement 8.2.3.3, we can represent $P_2$-numbers $b_{21}, b_{22}$ as tuples of $B_{2i}$-numbers

$$\begin{array}{ll} b_{21} = (b_{21i}, i \in I) & b_{21i} = t_{2i}(b_{21}) \in B_{2i} \\ b_{22} = (b_{22i}, i \in I) & b_{22i} = t_{2i}(b_{22}) \in B_{2i} \end{array} \tag{8.3.7}$$

According to the theorem 3.1.9, since the representation $h_i$ is single transitive, then there exists unique $B_{1i}$-number $b_{1i}$ such that

$$b_{22i} = h_i(b_{1i})(b_{21i})$$

According to definitions (8.3.2), (8.3.5), (8.3.7), there exists unique $P_1$-number $b_1$ such that

$$b_{22} = g(b_1)(b_{21})$$

According to the theorem 3.1.9, the representation $g$ is single transitive.

From commutativity of diagram (1) and from the definition 3.2.2, it follows that map $(t_{1i}, t_{2i})$ is morphism of representation $g$ into representation $h_i$. ⊙

Let

$$d_2 = g(b_1)(b_2) \quad d_2 = (d_{2i}, i \in I) \tag{8.3.8}$$

From equations (8.3.5), (8.3.8), it follows that

$$d_{2i} = h_i(b_{1i})(b_{2i}) \tag{8.3.9}$$

For $j = 1, 2$, let $R_j$ be other object of category $\mathcal{A}_j$. For any $i \in I$, let

$$r_{1i} : R_1 \longrightarrow B_{1i}$$

be morphism from $\Omega_1$-algebra $R_1$ into $\Omega_1$-algebra $B_{1i}$. According to definition 8.2.1, there exists a unique morphism of $\Omega_1$-algebra

$$s_1 : R_1 \longrightarrow P_1$$

such that following diagram is commutative

(8.3.10) 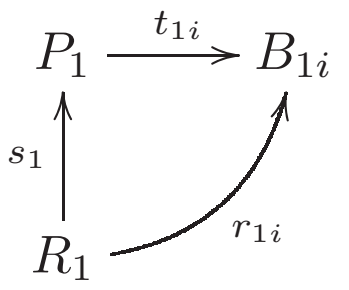

$$t_{1i} \circ s_1 = r_{1i}$$

Let $a_1 \in R_1$. Let

$$b_1 = s_1(a_1) \in P_1 \tag{8.3.11}$$

From commutativity of the diagram (8.3.10) and statements (8.3.11), (8.3.2), it follows that

$$b_{1i} = r_{1i}(a_1) \tag{8.3.12}$$

Let

$$f : R_1 \overset{*}{\longrightarrow} R_2$$

be single transitive $R_1*$-representation in $\Omega_2$-algebra $R_2$. According to the theorem 3.2.10, a morphism of $\Omega_2$-algebra

$$r_{2i} : R_2 \longrightarrow B_{2i}$$

such that map $(r_{1i}, r_{2i})$ is morphism of representations from $f$ into $h_i$ is unique up to choice of image of $R_2$-number $a_2$. According to the remark 3.2.6, in diagram of maps

(8.3.13)

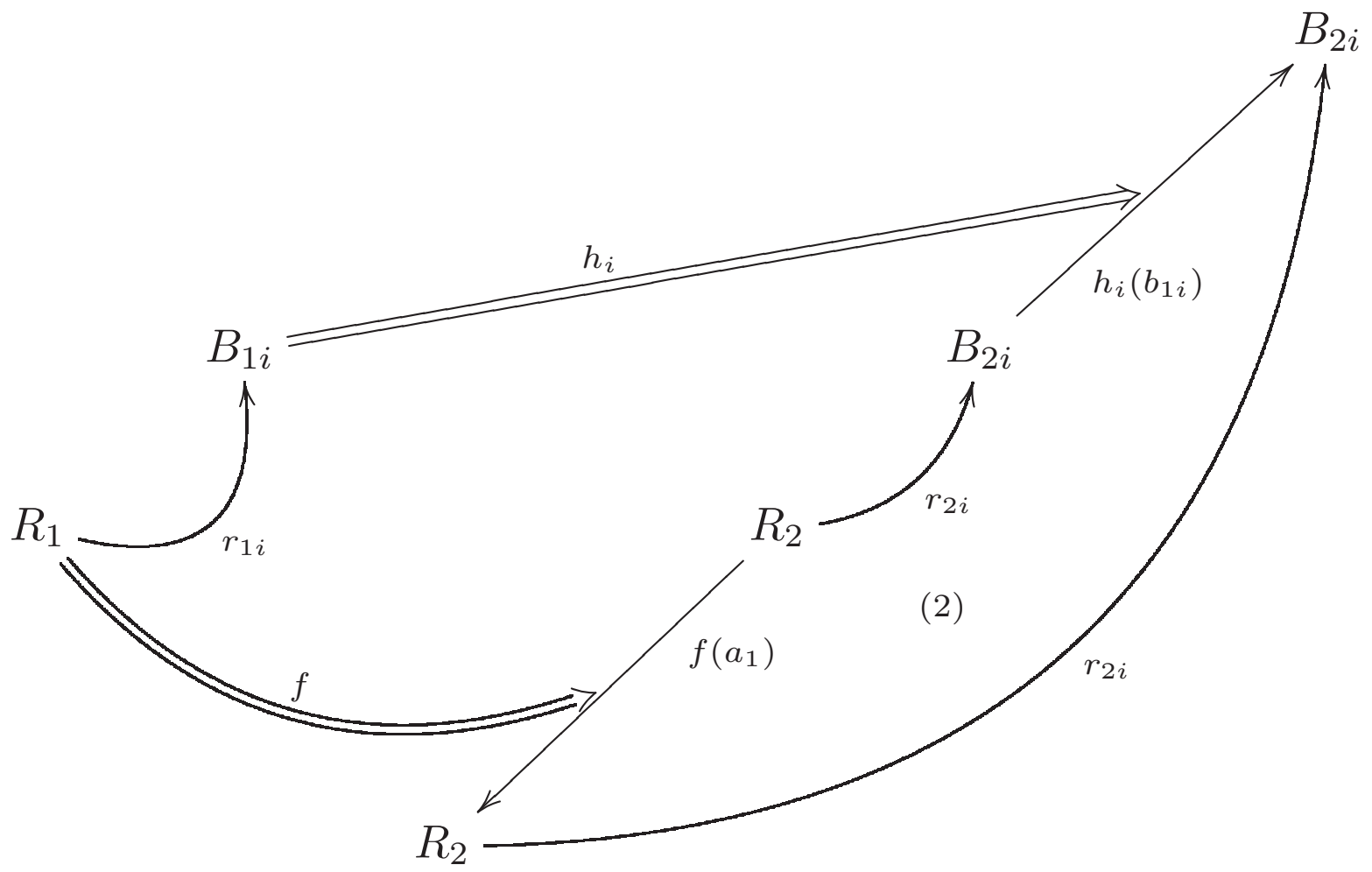

diagram (2) is commutative. According to definition 8.2.1, there exists a unique morphism of $\Omega_2$-algebra

$$s_2 : R_2 \longrightarrow P_2$$

such that following diagram is commutative

(8.3.14)

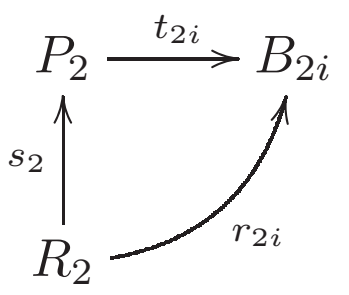

$t_{2i} \circ s_2 = r_{2i}$

Let $a_2 \in R_2$. Let

(8.3.15) $$b_2 = s_2(a_2) \in P_2$$

From commutativity of the diagram (8.3.14) and statements (8.3.15), (8.3.3), it follows that

(8.3.16) $$b_{2i} = r_{2i}(a_2)$$

Let

(8.3.17) $$c_2 = f(a_1)(a_2)$$

From commutativity of the diagram (2) and equations (8.3.9), (8.3.16), (8.3.17), it follows that

(8.3.18) $$d_{2i} = r_{2i}(c_2)$$

From equations (8.3.9), (8.3.18), it follows that

(8.3.19) $$d_2 = s_2(c_2)$$

and this is consistent with commutativity of the diagram (8.3.14).

For each $i \in I$,we join diagrams of maps (8.3.4), (8.3.10), (8.3.14), (8.3.13)

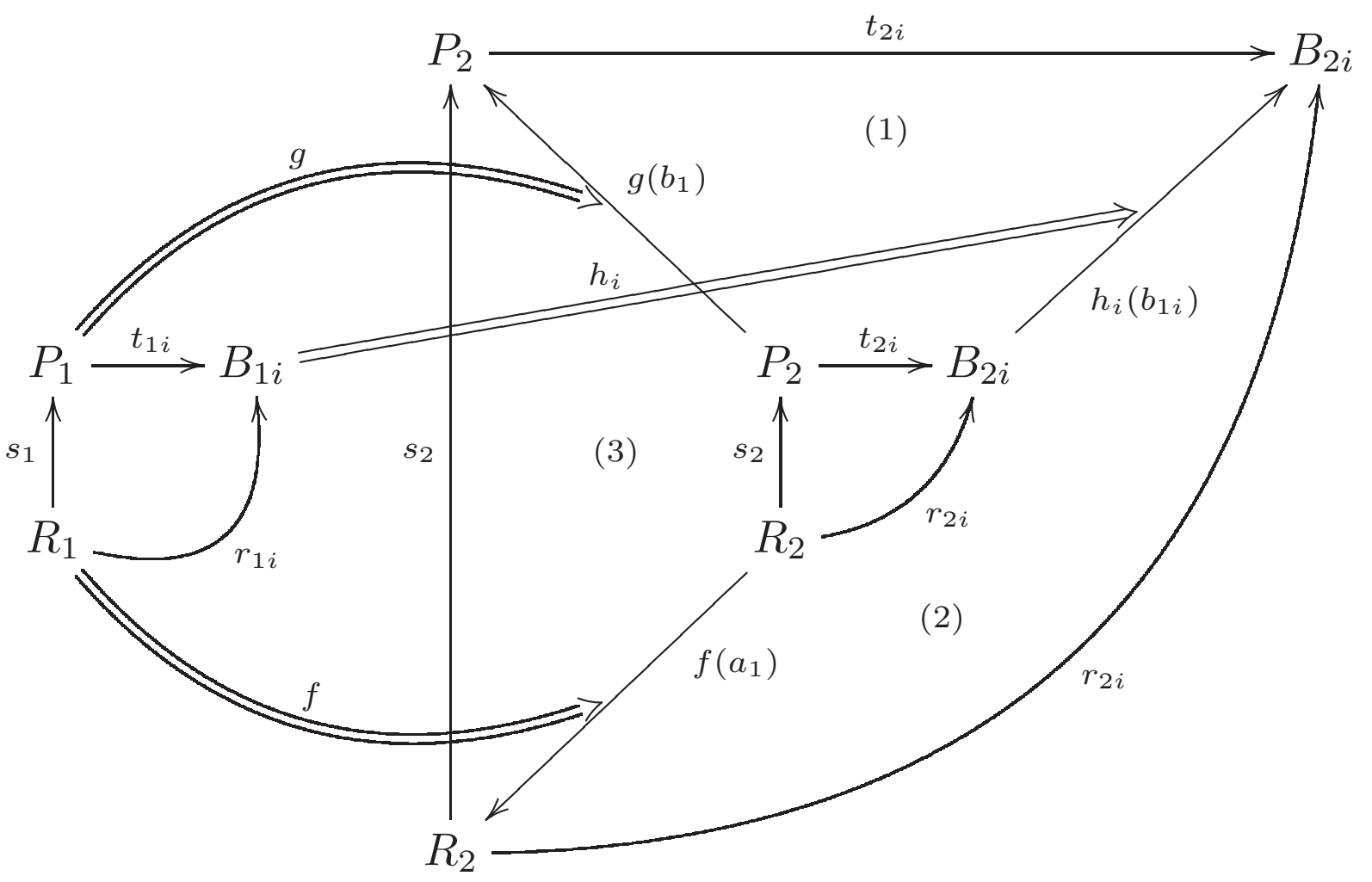

From equations (8.3.8), (8.3.15) and from equations (8.3.17), (8.3.19), commutativity of the diagram (3) follows. Therefore, the map $(s_1, s_2)$ is morphism of representations from $f$ into $g$. According to the theorem 3.2.10, the morphism $(s_1, s_2)$ is defined unambiguously, since we require (8.3.19).

According to the definition 8.2.1, the representation $g$ and family of morphisms of representation $((t_{1i}, t_{2i}), i \in I)$ is product in the category $(\mathcal{A}_1 *)\mathcal{A}_2$ . □

THEOREM 8.3.4. *In category $\mathcal{G}*$ of left-side representations of group there exists product of effective left-side representations of group.*

PROOF. Statement of theorem is corollary of theorems 8.3.2 and 4.2.10 $\square$

DEFINITION 8.3.5. Let $\{V_i, i \in I\}$ be the set of ${}_*{}^*D_i$-vector spaces. Left-side representation of ring $D = \prod_{i \in I} D_i$ in the Abelian group $V = \prod_{i \in I} V_i$ is called **direct product**

$$V = \prod_{i \in I} V_i$$

**of ${}_*{}^*D_i$-vector spaces**[8.7] $V_i$, $i \in I$. $\square$

If $|I| = n$, then we also will use notation

$$V = \prod_{i=1}^{n} V_i = V_1 \times \dots \times V_n$$

for direct product of ${}_*{}^*D_i$-vector spaces $V_1$, ..., $V_n$

According to definition we define representation of ring $D$ in Abelian group $V$ componentwise. If $a = (a_i, i \in I) \in D$ and $v = (v_i, i \in I) \in V$, we define right-side representation corresponding to element $a$

$$va = (v_i a_i, i \in I)$$

We consider convention described in remark 2.2.15 componentwise. We will write linear combination of ${}_*{}^*D$-vectors $\overline{v}_k$ as

$$\overline{v}_*{}^*a = (\overline{v}_{i*}{}^*a_i, i \in I) = (v_{i \cdot k}\ {}^k a_i, i \in I)$$

We also will use notation

$$\overline{v}_*{}^*a = (v_{i \cdot *}, i \in I)\ ({}^*a_i, i \in I) = (v_{i \cdot k}, i \in I)\ ({}^k a_i, i \in I) \tag{8.3.20}$$

If $|I| = n$, then we can write equation (8.3.20) as

$$\overline{v}_*{}^*a = (v_{1 \cdot *}, ..., v_{n \cdot *})\ ({}^*a_1, ..., {}^*a_n) = (v_{1 \cdot k_1}, ..., v_{n \cdot k_n})\ ({}^{k_1}a_1, ..., {}^{k_n}a_n)$$

Since direct product of division rings is not division ring, then direct product of ${}_*{}^*D_i$-vector spaces, in general, is module. However structure of this module is close to structure of ${}_*{}^*D$-vector space.

THEOREM 8.3.6. *Let $V_1$, ..., $V_n$ be the set of ${}_*{}^*D_i$-vector spaces. Let $\overline{\overline{e}}_i$ be the ${}_*{}^*D_i$-basis of vector space $V_i$. Set of vectors $(e_{1 \cdot i_1}, ..., e_{n \cdot i_n})$ builds up a basis of direct product*

$$V = V_1 \times \dots \times V_n$$

PROOF. We prove the statement of theorem by induction on $n$.

When $n = 1$ the statement is evident.

Let statement be true when $n = k - 1$. We represent space $V$ as[8.8]

$$V = V_1 \times \dots \times V_k = (V_1 \times \dots \times V_{k-1}) \times V_k$$

---

[8.7] I made definition according to [1], p. 127

[8.8] In general, we need consider more wide class of modules like

$$V_1 \times \dots \times V_k$$

including associativity of product.

Correspondingly, we can represent arbitrary vector $(\overline{v}_1, ..., \overline{v}_k) \in V$ as

$$(\overline{v}_1, ..., \overline{v}_k) = ((\overline{v}_1, ..., \overline{v}_{k-1}), \overline{v}_k) = ((\overline{v}_1, ..., \overline{v}_{k-1}), 0) + (0, \overline{v}_k) \tag{8.3.21}$$

$(\overline{v}_1, ..., \overline{v}_{k-1}) \in V_1 \times ... \times V_{k-1}$ and according to the assumption of induction

$$(\overline{v}_1, ..., \overline{v}_{k-1}) = (e_{1\cdot*}, ..., e_{k-1\cdot*})\ ({}^*v_1, ..., {}^*v_{k-1}) \tag{8.3.22}$$

From equations (8.3.21) and (8.3.22) it follows

$$\begin{aligned}(\overline{v}_1, ..., \overline{v}_k) &= ((e_{1\cdot*}, ..., e_{k-1\cdot*})\ ({}^*v_1, ..., {}^*v_{k-1}), 0) + (0, e_{k*}{}^*v_k)\\ &= (e_{1\cdot*}, ..., e_{k-1\cdot*}, e_k)\ ({}^*v_1, ..., {}^*v_{k-1}, 0) + (e_1, ..., e_{k-1}, e_{k\cdot*})\ (0, ..., {}^*v_k)\\ &= (e_{1\cdot*}, ..., e_{k-1\cdot*}, e_{k\cdot*})\ ({}^*v_1, ..., {}^*v_{k-1}, {}^*v_k)\end{aligned}$$

Therefore the statement is true when $n = k$. □

When $V_1$, ..., $V_n$ are $D$-vector spaces, direct product is called **direct product of $D$-vector spaces**.

THEOREM 8.3.7. *Let category $\mathcal{A}$ of $\Omega$-algebras have product. Then in category $A*$*[8.9] *where $A$ is $\Omega$-algebra, there exists product of single transitive $A*$-representations.*

PROOF. Proof of theorem is similar to proof of theorem 8.3.2 with the only difference that we use diagram

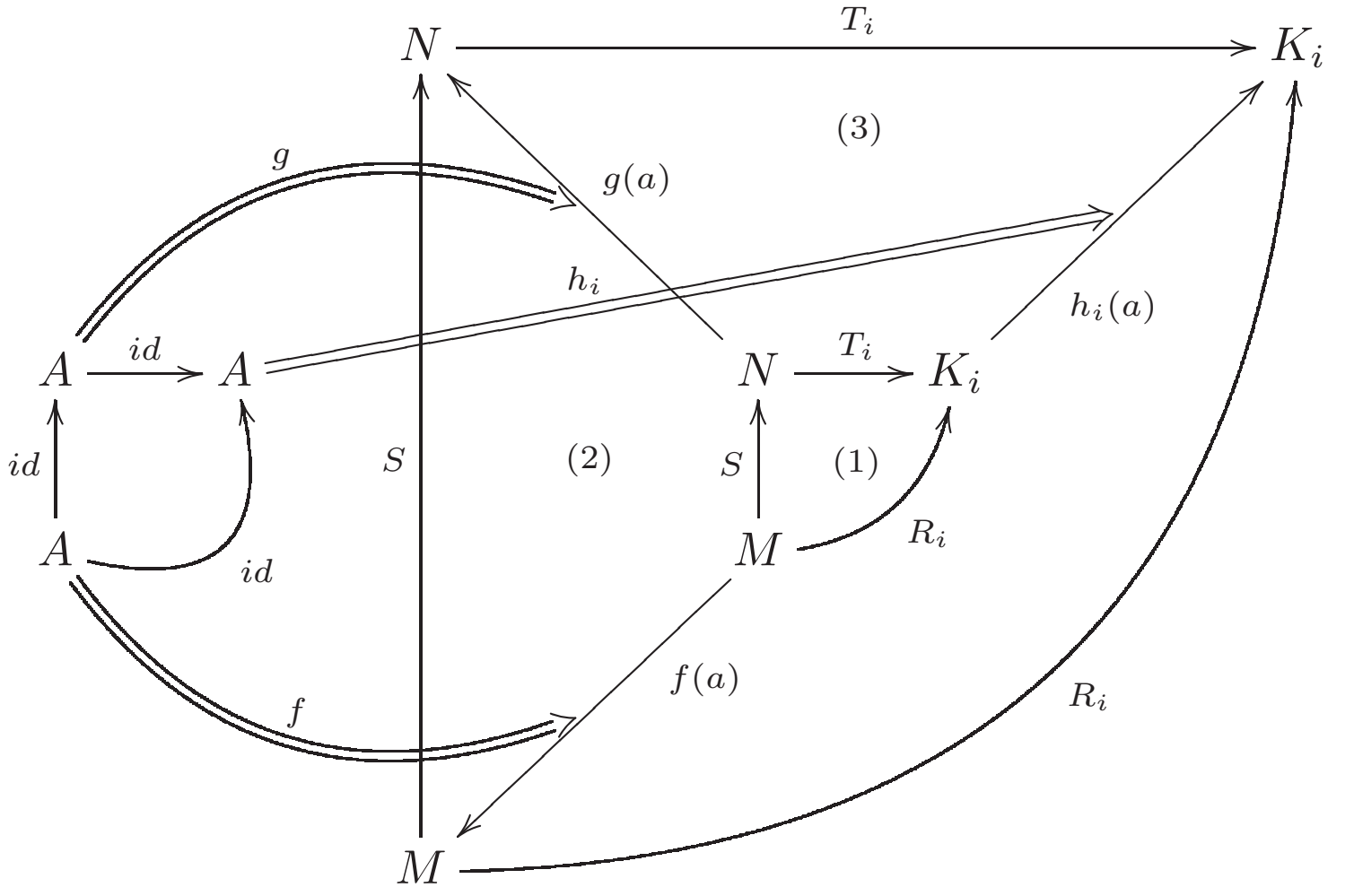

□

Let us consider relationship between product of single transitive left-side representations of $\Omega$-algebra $A$ in category $\mathcal{A}*$ and product of single transitive left-side representations of $\Omega$-algebra $A$ in category $A*$. Let for any $i \in I$ representation $\Omega$-of algebra $A$ be defined on set $K_i$. For every $i \in I$ let us consider diagram

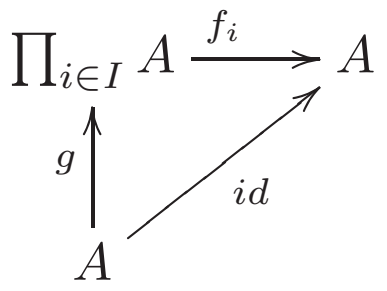

[8.9]See definition 3.2.20

where $g$ is map of algebra $A$ onto diagonal in $\prod_{i\in I} A$ and $f_i$ is projection onto factor $i$. From diagram, it follows that maps $g$ and $f_i$ are injective. Therefore both products are equivalent.

## 8.4. Morphisms of Direct Product of ${}_*{}^*D$-Vector Spaces

Suppose $V_i$, $i\in I$, $W_j$, $j\in J$ are ${}_*{}^*D$-vector spaces. Assume

$$V=\prod_{i\in I} V_i$$

$$W=\prod_{j\in J} W_j$$

Linear map

$$f:V\longrightarrow W \tag{8.4.1}$$

holds linear operations. Since we defined operations in module $W$ componentwise, we can represent map $f$ as

$$f=(\ f_j:V\longrightarrow W_j, j\in J\ )$$

Let

$$I_i:V_i\longrightarrow V \qquad i\in I$$

be injection $V_i$ into $V$. Then for any vector $v=(v_i,i\in I)\in V$, we get

$$f_j(v)=\sum_{i\in I} f_j{}_*{}^*I_i{}_*{}^*v_i$$

Therefore, we can represent the mapping (8.4.1) as matrix of ${}_*{}^*D$-linear maps

$$f=(f_{ij}=f_j{}_*{}^*I_i:\ V_i\longrightarrow W_j\ , i\in I, j\in J)$$

CHAPTER 9

# Geometry of Division Ring

## 9.1. Center of Division Ring

DEFINITION 9.1.1. Let $D$ be a ring.[9.1] The set $Z(D)$ of elements $a \in D$ such that

$$ax = xa \tag{9.1.1}$$

for all $x \in D$, is called **center of ring** $D$. □

THEOREM 9.1.2. *The center $Z(D)$ of ring $D$ is subring of ring $D$.*

PROOF. The statement follows immediately from definition 9.1.1. □

DEFINITION 9.1.3. Let $D$ be a ring with unit element $e$.[9.2] The map

$$l : Z \to D$$

such that $l(n) = ne$ is a homomorphism of rings, and its kernel is an ideal $(n)$, generated by integer $n \geq 0$. We have canonical injective homomorphism

$$Z/nZ \to D$$

which is an isomorphism between $Z/nZ$ and subring of $D$. If $nZ$ is prime ideal, then we have two cases.

- $n = 0$. $D$ contains as subring a ring which isomorphic to $Z$, and which is often identified with $Z$. In that case, we say that $D$ has **characteristic** 0.
- $n = p$ for some prime number $p$. $D$ has **characteristic** $p$, and $D$ contains an isomorphic image of $F_p = Z/pZ$.

□

THEOREM 9.1.4. *Let $D$ be ring of characteristic 0 and let $d \in D$. Then every integer $n \in Z$ commutes with $d$.*

PROOF. We prove statement by induction. The statement is evident for $n = 0$ and $n = 1$. Let statement be true for $n = k$. From chain of equation

$$(k + 1)d = kd + d = dk + d = d(k + 1)$$

Evidence of statement for $n = k + 1$ follows. □

THEOREM 9.1.5. *Let $D$ be ring of characteristic 0. Then ring of integers $Z$ is subring of center $Z(D)$ of ring $D$.*

PROOF. Corollary of theorem 9.1.4. □

[9.1] [1], page 89.
[9.2] I made definition according to definition from [1], pages 89, 90.

Let $D$ be division ring. If $D$ has characteristic 0, $D$ contains as subfield an isomorphic image of the field $Q$ of rational numbers. If $D$ has characteristic $p$, $D$ contains as subfield an isomorphic image of $F_p$. In either case, this subfield will be called the prime field. Since the prime field is the smallest subfield of $D$ containing 1 and has no automorphism except identity, it is customary to identify it with $Q$ or $F_p$ as the case may be.

THEOREM 9.1.6. *The center $Z(D)$ of division ring $D$ is subfield of division ring $D$.*

PROOF. According to theorem 9.1.2 it is enough to verify that $a^{-1} \in Z(D)$ if $a \in Z(D)$. Let $a \in Z(D)$. Repeatedly using the equation (9.1.1) we get chain of equations

$$aa^{-1}x = x = xaa^{-1} = axa^{-1} \tag{9.1.2}$$

From (9.1.2) it follows

$$a^{-1}x = xa^{-1}$$

Therefore, $a^{-1} \in Z(D)$. □

THEOREM 9.1.7. *Let $D$ be division ring of characteristic 0 and let $d \in D$. Then for any integer $n \in Z$*

$$n^{-1}d = dn^{-1} \tag{9.1.3}$$

PROOF. According to theorem 9.1.4 following chain of equation is true

$$n^{-1}dn = nn^{-1}d = d \tag{9.1.4}$$

Let us multiply right and left sides of equation (9.1.4) by $n^{-1}$. We get

$$n^{-1}d = n^{-1}dnn^{-1} = dn^{-1} \tag{9.1.5}$$

(9.1.3) follows from (9.1.5). □

THEOREM 9.1.8. *Let $D$ be division ring of characteristic 0 and let $d \in D$. Then every rational number $p \in Q$ commutes with $d$.*

PROOF. Let us represent rational number $p \in Q$ as $p = mn^{-1}$, $m$, $n \in Z$. Statement of theorem follows from chain of equations

$$pd = mn^{-1}d = n^{-1}dm = dmn^{-1} = dp$$

based on the statement of theorem 9.1.4 and equation (9.1.3). □

THEOREM 9.1.9. *Let $D$ be division ring of characteristic 0. Then field of rational numbers $Q$ is subfield of center $Z(D)$ of division ring $D$.*

PROOF. Corollary of theorem 9.1.8. □

## 9.2. Geometry of Division Ring over Field

We may consider division ring $D$ as vector space over field $F \subset Z(D)$. Because $F$ is field, we can write all indexes on right side of root letter. We will use convention 1.5.5.

Remark 9.2.1. Let $\overline{\overline{e}}$ be basis of division ring $D$ over field $F$. Then we may present any element $a \in D$ as

$$a = e_i a^i \quad a^i \in F \tag{9.2.1}$$

When dimension of division ring $D$ over field $F$ infinite, then basis may be either countable, or its power is not less than power of continuum. If basis is countable, then we put constraints on coefficients $a^i$ of expansion (9.2.1). If power of the set $I$ is continuum, then we assume that there is measure on the set $I$ and sum in expansion (9.2.1) is integral over this measure. □

Remark 9.2.2. Since we defined product in the division ring $D$, we consider the division ring as algebra over field $F \subset Z(D)$. For elements of basis we assume

$$e_i e_j = e_k C^k_{ij} \tag{9.2.2}$$

Coefficients $C^k_{ij}$ of expansion (9.2.2) are called **structural constants** of division ring $D$ over field $F$.

From equations (9.2.1), (9.2.2), it follows

$$ab = e_k C^k_{ij} a^i b^j \tag{9.2.3}$$

From equation (9.2.3) it follows that

$$(ab)c = e_k C^k_{ij} (ab)^i c^j = e_k C^k_{ij} C^i_{mn} a^m b^n c^j \tag{9.2.4}$$

$$a(bc) = e_k C^k_{ij} a^i (bc)^j = e_k C^k_{ij} a^i C^j_{mn} b^m c^n \tag{9.2.5}$$

From associativity of product

$$(ab)c = a(bc)$$

and equations (9.2.4) and (9.2.5) it follows that

$$e_k C^k_{ij} C^i_{mn} a^m b^n c^j = e_k C^k_{ij} a^i C^j_{mn} b^m c^n \tag{9.2.6}$$

Because vectors $a$, $b$, $c$ are arbitrary, and vectors $e_k$ are linearly independent, then from equation (9.2.6) it follows that

$$C^k_{jn} C^j_{im} = C^k_{ij} C^j_{mn} \tag{9.2.7}$$

Theorem 9.2.3. *Coordinates $a^j$ of vector $a$ are tensor*

$$a^j = A^j_i a'^i \tag{9.2.8}$$

Proof. Let $\overline{\overline{e}}'$ be another basis. Let

$$e'_i = e_j A^j_i \tag{9.2.9}$$

be transformation, mapping basis $\overline{\overline{e}}$ into basis $\overline{\overline{e}}'$. Because vector $a$ does not change, then

$$a = e'_i a'^i = e_j a^j \tag{9.2.10}$$

From equations (9.2.9) and (9.2.10) it follows that

$$e_j a^j = e'_i a'^i = e_j A^j_i a'^i \tag{9.2.11}$$

Because vectors $e_j$ are linearly independent, then equation (9.2.8) follows from equation (9.2.11). Therefore, coordinates of vector are tensor. □

Theorem 9.2.4. *Structural constants of division ring $D$ over field $F$ are tensor*

$$A_{\boldsymbol{k}}^{\boldsymbol{l}} C'^{\boldsymbol{k}}_{\boldsymbol{ij}} A^{-1\cdot}{}^{\boldsymbol{i}}_{\boldsymbol{n}} A^{-1\cdot}{}^{\boldsymbol{j}}_{\boldsymbol{m}} = C^{\boldsymbol{l}}_{\boldsymbol{nm}} \tag{9.2.12}$$

Proof. Consider similarly the transformation of product. Equation (9.2.3) has form

$$ab = e'_{\boldsymbol{k}} C'^{\boldsymbol{k}}_{\boldsymbol{ij}} a'^{\boldsymbol{i}} b'^{\boldsymbol{j}} \tag{9.2.13}$$

relative to basis $\overline{\overline{e}}'$. Let us substitute (9.2.8) and (9.2.9) into (9.2.13). We get

$$ab = e_{\boldsymbol{l}} A_{\boldsymbol{k}}^{\boldsymbol{l}} C'^{\boldsymbol{k}}_{\boldsymbol{ij}} a^{\boldsymbol{n}} A^{-1\cdot}{}^{\boldsymbol{i}}_{\boldsymbol{n}} b^{\boldsymbol{m}} A^{-1\cdot}{}^{\boldsymbol{j}}_{\boldsymbol{m}} \tag{9.2.14}$$

From (9.2.3) and (9.2.14) it follows that

$$e_{\boldsymbol{l}} A_{\boldsymbol{k}}^{\boldsymbol{l}} C'^{\boldsymbol{k}}_{\boldsymbol{ij}} A^{-1\cdot}{}^{\boldsymbol{i}}_{\boldsymbol{n}} a^{\boldsymbol{n}} A^{-1\cdot}{}^{\boldsymbol{j}}_{\boldsymbol{m}} b^{\boldsymbol{m}} = e_{\boldsymbol{l}} C^{\boldsymbol{l}}_{\boldsymbol{nm}} a^{\boldsymbol{n}} b^{\boldsymbol{m}} \tag{9.2.15}$$

Because vectors $a$ and $b$ are arbitrary, and vectors $e_{\boldsymbol{l}}$ are linearly independent, then equation (9.2.12) follows from equation (9.2.15). Therefore, structural constants are tensor. □

CHAPTER 10

# Linear Map of Division Ring

## 10.1. Linear Map of Division Ring

According to the remark 9.2.2, we consider division ring $D$ as algebra over field $F \subset Z(D)$.

Definition 10.1.1. Let $D_1$, $D_2$ be division rings. Let $F$ be field such that $F \subset Z(D_1)$, $F \subset Z(D_2)$. Linear map

$$f : D_1 \to D_2$$

of $F$-vector space $D_1$ into $F$-vector space $D_2$ is called **linear map of division ring $D_1$ into division ring $D_2$**. □

According to definition 10.1.1, linear map $f$ of division ring $D_1$ into division ring $D_2$ holds

$$f(a+b) = f(a) + f(b) \quad a, b \in D$$

$$f(pa) = pf(a) \quad p \in F$$

Theorem 10.1.2. *Let map*

$$f : D_1 \to D_2$$

*is linear map of division ring $D_1$ of characteristic 0 into division ring $D_2$ of characteristic 0. Then*[10.1]

$$f(nx) = nf(x)$$

*for any integer $n$.*

Proof. We prove the theorem by induction on $n$. Statement is obvious for $n = 1$ because

$$f(1x) = f(x) = 1f(x)$$

Let statement is true for $n = k$. Then

$$f((k+1)x) = f(kx + x) = f(kx) + f(x) = kf(x) + f(x) = (k+1)f(x)$$

□

[10.1] Let

$$f : D_1 \to D_2$$

be a linear map of division ring $D_1$ of characteristic 2 into division ring $D_2$ of characteristic 3. Then for any $a \in D_1$, $2a = 0$, although $2f(a) \neq 0$. Therefore, if we assume that characteristic of division ring $D_1$ is greater than 0, then must demand that characteristic of division ring $D_1$ equal to characteristic of division ring $D_2$.

Theorem 10.1.3. *Let map*

$$f : D_1 \to D_2$$

*be linear map of division ring $D_1$ of characteristic 0 into division ring $D_2$ of characteristic 0. Then*

$$f(ax) = af(x)$$

*for any rational $a$.*

Proof. Let $a = \frac{p}{q}$. Assume $y = \frac{1}{q}x$. Then according to the theorem 10.1.2

$$f(x) = f(qy) = qf(y) = qf\left(\frac{1}{q}x\right) \tag{10.1.1}$$

From equation (10.1.1) it follows

$$\frac{1}{q}f(x) = f\left(\frac{1}{q}x\right) \tag{10.1.2}$$

From equation (10.1.2) it follows

$$f\left(\frac{p}{q}x\right) = pf\left(\frac{1}{q}x\right) = \frac{p}{q}f(x)$$

□

We cannot extend the statement of theorem 10.1.3 for arbitrary subfield of center $Z(D)$ of division ring $D$.

Theorem 10.1.4. *Let division ring $D$ is algebra over field $F \subset Z(D)$. If $F \neq Z(D)$, then there exists linear map*

$$f : D \to D$$

*which is not linear over field $Z(D)$.*

Proof. To prove the theorem it is enough to consider the complex field $C$ because $C = Z(C)$. Because the complex field is algebra over real field, then the function

$$z \to \overline{z}$$

is linear. However the equation

$$\overline{az} = a\overline{z}$$

is not true. □

Based on theorem 10.1.4, the question arises. Why do we consider linear maps over field $F \neq Z(D)$, if this leads us to sharp expansion of the set of linear maps? The answer to this question is a rich experience of the theory of complex function.

Theorem 10.1.5. *Let maps*

$$f : D_1 \to D_2$$

$$g : D_1 \to D_2$$

*be linear maps of division ring $D_1$ into division ring $D_2$. Then map $f + g$ is linear.*

Proof. Statement of theorem follows from chain of equations

$$\begin{aligned}(f+g)(x+y)=&f(x+y)+g(x+y)=f(x)+f(y)+g(x)+g(y)\\=&(f+g)(x)+(f+g)(y)\\(f+g)(px)=&f(px)+g(px)=pf(x)+pg(x)=p(f(x)+g(x))\\=&p(f+g)(x)\end{aligned}$$

□

Theorem 10.1.6. *Let map*

$$f:D_1\to D_2$$

*be linear map of division ring $D_1$ into division ring $D_2$. Then maps $af$, $fb$, $a$, $b\in R_2$ are linear.*

Proof. Statement of theorem follows from chain of equations

$$\begin{aligned}(af)(x+y)=&a(f(x+y))=a(f(x)+f(y))=af(x)+af(y)\\=&(af)(x)+(af)(y)\\(af)(px)=&a(f(px))=a(pf(x))=p(af(x))\\=&p(af)(x)\\(fb)(x+y)=&(f(x+y))b=(f(x)+f(y))b=f(x)b+f(y)b\\=&(fb)(x)+(fb)(y)\\(fb)(px)=&(f(px))b=(pf(x))b=p(f(x)b)\\=&p(fb)(x)\end{aligned}$$

□

Definition 10.1.7. Denote $\mathcal{L}(D_1;D_2)$ set of linear maps

$$f:D_1\to D_2$$

of division ring $D_1$ into division ring $D_2$. □

Theorem 10.1.8. *We may represent linear map*

$$f:D_1\to D_2$$

*of division ring $D_1$ into division ring $D_2$ as*

$$f(x)=f_{k\cdot s_k\cdot 0}\ G_k(x)\ f_{k\cdot s_k\cdot 1}\tag{10.1.3}$$

*where $(G_k,k\in K)$ is set of additive maps of division ring $D_1$ into division ring $D_2$.*[10.2] *Expression $f_{k\cdot s_k\cdot p}$, $p=0,\ 1,$ in equation (10.1.3) is called* **component of linear map** *$f$.*

Proof. The statement of theorem follows from theorems 10.1.5 and 10.1.6. □

If in the theorem 10.1.8 $|K|=1$, then the equation (10.1.3) has form

$$f(x)=f_{s\cdot 0}\ G(x)\ f_{s\cdot 1}\tag{10.1.4}$$

and map $f$ is called **linear map generated by map** $G$. Map $G$ is called **generator of linear map**.

[10.2] Here and in the following text we assume sum over index that is used in product few times. Equation (10.1.3) is recursive definition and there is hope that it is possible to simplify it.

THEOREM 10.1.9. *Let $D_1$, $D_2$ be division rings of characteristic* 0. *Let $F$, $F \subset Z(D_1)$, $F \subset Z(D_2)$, be field. Let $G$ be linear map. Let $\overline{\overline{e}}$ be basis of division ring $D_2$ over field $F$.* **Standard representation of linear map** (10.1.4) *has form*[10.3]

$$f(x) = f_G^{ij}\ e_i\ G(x)\ e_j \tag{10.1.5}$$

*Expression $f_G^{ij}$ in equation* (10.1.5) *is called* **standard component of linear map** *$f$.*

PROOF. Components of linear map $f$ have expansion

$$f_{s\cdot p} = f_{s\cdot p}^{i}\ e_i \tag{10.1.6}$$

relative to basis $\overline{\overline{e}}$. If we substitute (10.1.6) into (10.1.4), we get

$$f(x) = f_{s\cdot 0}^{i}\ e_i\ G(x)\ f_{s\cdot 1}^{j}\ e_j \tag{10.1.7}$$

If we substitute expression

$$f_G^{ij} = f_{s\cdot 0}^{i}\ f_{s\cdot 1}^{j}$$

into equation (10.1.7) we get equation (10.1.5). □

THEOREM 10.1.10. *Let $D_1$, $D_2$ be division rings of characteristic* 0. *Let $F$, $F \subset Z(D_1)$, $F \subset Z(D_2)$, be field Let $G$ be linear map. Let $\overline{\overline{e}}_1$ be basis of division ring $D_1$ over field $F$. Let $\overline{\overline{e}}_2$ be basis of division ring $D_2$ over field $F$. Let $C_{2\cdot}{}_{kl}^{p}$ be structural constants of division ring $D_2$. Then it is possible to represent linear map* (10.1.4) *generated by linear map $G$ as*

$$\begin{aligned} f(a) &= e_{2\cdot j}\ f_i^j a^i & f_k^j &\in F \\ a &= e_{1\cdot i}\ a^i & a^i &\in F \quad a \in D_1 \end{aligned} \tag{10.1.8}$$

$$f_i^j = G_i^l f_G^{kr} C_{2\cdot}{}_{kl}^{p} C_{2\cdot}{}_{pr}^{j} \tag{10.1.9}$$

PROOF. Consider map

$$\begin{aligned} G : D_1 \to D_2 \quad a = e_{1\cdot i} a^i \to G(a) = e_{2\cdot j} G_i^j a^i \\ a^i \in F \quad G_i^j \in F \end{aligned} \tag{10.1.10}$$

According to the theorem 5.4.3, linear map $f(a)$ relative to bases $\overline{\overline{e}}_1$ and $\overline{\overline{e}}_2$ has form (10.1.8). From equations (10.1.5) and (10.1.10), it follows

$$f(a) = G_i^l a^i f_G^{kj} e_{2\cdot k} e_{2\cdot l} e_{2\cdot j} \tag{10.1.11}$$

From equations (10.1.8) and (10.1.11), it follows

$$e_{2\cdot j}\ f_i^j\ a^i = G_i^l a^i f_G^{kr} e_{2\cdot k} e_{2\cdot l} e_{2\cdot r} = G_i^l a^i f_G^{kr} C_{2\cdot}{}_{kl}^{p} C_{2\cdot}{}_{pr}^{j}\ e_{2\cdot j} \tag{10.1.12}$$

Since vectors $e_{2\cdot r}$ are linear independent over field $F$ and values $a^k$ are arbitrary, then equation (10.1.9) follows from equation (10.1.12). □

[10.3]Representation of linear map of of division ring using components of linear map is ambiguous. We can increase or decrease number of summands using algebraic operations. Since dimension of division ring $D_2$ over field $F$ is finite, standard representation of linear map guarantees finiteness of set of items in the representation of map.

Considering map

$$f : D \to D \tag{10.1.13}$$

we assume $G(x) = x$.

THEOREM 10.1.11. *Let $D$ be division ring of characteristic 0. Linear map (10.1.13) has form*

$$f(x) = f_{s\cdot 0}\ x\ f_{s\cdot 1} \tag{10.1.14}$$

THEOREM 10.1.12. *Let $D$ be division ring of characteristic 0. Let $\overline{\overline{e}}$ be the basis of division ring $D$ over field $F \subset Z(D)$. Standard representation of linear map (10.1.14) of division ring has form*

$$f(x) = f^{ij}\ e_i x e_j \tag{10.1.15}$$

THEOREM 10.1.13. *Let $D$ be division ring of characteristic 0. Let $\overline{\overline{e}}$ be basis of division ring $D$ over field $F \subset Z(D)$. Then it is possible to represent linear map (10.1.13) as*

$$\begin{aligned} f(a) &= e_j f^j_i a^i & f^j_k &\in F \\ a &= e_i a^i & a^i &\in F \quad a \in D \end{aligned} \tag{10.1.16}$$

$$f^j_i = f^{kr} C^p_{ki} C^j_{pr} \tag{10.1.17}$$

THEOREM 10.1.14. *Consider matrix*

$$\mathcal{C} = \left(\mathcal{C}^{\cdot j}_{i \cdot kr}\right) = \left(C^p_{ki} C^j_{pr}\right) \tag{10.1.18}$$

*whose rows are indexed by $^{\cdot j}_{i}$ and columns are indexed by $\cdot kr$ . If $\det \mathcal{C} \neq 0$, then, for given coordinates of linear transformation $f^j_i$, the system of linear equations (10.1.17) with standard components of this transformation $f^{kr}$ has the unique solution. If $\det \mathcal{C} = 0$, then the equation*

$$\operatorname{rank} \left(\mathcal{C}^{\cdot j}_{i \cdot kr} \quad f^j_i\right) = \operatorname{rank} \mathcal{C} \tag{10.1.19}$$

*is the condition for the existence of solutions of the system of linear equations (10.1.17). In such case the system of linear equations (10.1.17) has infinitely many solutions and there exists linear dependence between values $f^j_i$.*

PROOF. Equation (10.1.14) is special case of equation (10.1.4) when $G(x) = x$. Theorem 10.1.12 is special case of theorem 10.1.9 when $G(x) = x$. Theorem 10.1.13 is special case of theorem 10.1.10 when $G(x) = x$. The statement of the theorem 10.1.14 is corollary of the theory of linear equations over field. □

THEOREM 10.1.15. *Standard components of the identity map have the form*

$$f^{kr} = \delta^k_0 \delta^r_0 \tag{10.1.20}$$

PROOF. The equation (10.1.20) is corollary of the equation

$$x = e_0\ x\ e_0$$

Let us show that the standard components (10.1.20) of a linear transformation satisfy to the equation

$$\delta^j_i = f^{kr} C^p_{ki} C^j_{pr} \tag{10.1.21}$$

which follows from the equation (10.1.17) if $f = \delta$. From equations (10.1.20), (10.1.21), it follows that

$$\delta_i^j = C_{0i}^p C_{p0}^j \tag{10.1.22}$$

The equation (10.1.22) is true, because, from equations

$$e_j e_0 = e_0\ e_j = e_j$$

it follows that

$$C_{0r}^j = \delta_r^j \quad C_{r0}^j = \delta_r^j$$

If $\det \mathcal{C} \neq 0$, then the solution (10.1.20) is unique. If $\det \mathcal{C} = 0$, then the system of linear equations (10.1.21) has infinitely many solutions. However, we are looking for at least one solution. □

THEOREM 10.1.16. *If* $\det \mathcal{C} \neq 0$, *then standard components of the zero map*

$$z : A \to A \quad z(x) = 0$$

*are defined uniquely and have form* $z^{ij} = 0$. *If* $\det \mathcal{C} = 0$, *then the set of standard components of the zero map forms a vector space.*

PROOF. The theorem is true because standard components $z^{ij}$ are solution of homogeneous system of linear equations

$$0 = z^{kr} C_{ki}^p C_{pr}^j$$

□

REMARK 10.1.17. Consider equation

$$a^{kr}\ e_k\ x\ e_r = b^{kr}\ e_k\ x\ e_r \tag{10.1.23}$$

From the theorem 10.1.16, it follows that only when condition $\det \mathcal{C} \neq 0$ is true, from the equation (10.1.23), it follows that

$$a^{kr} = b^{kr} \tag{10.1.24}$$

Otherwise, we must assume equality

$$a^{kr} = b^{kr} + z^{kr} \tag{10.1.25}$$

Despite this, in case $\det \mathcal{C} = 0$, we also use standard representation because in general it is very hard to show the set of linear independent vectors. If we want to define operation over linear maps in standard representation, then as well as in the case of the theorem 10.1.15 we choose one element from the set of possible representations. □

THEOREM 10.1.18. *Expression*

$$f_k^r = f^{ij} C_{ik}^p C_{pj}^r$$

*is tensor over field* $F$

$$f'^j_i = A_i^k f_k^l A^{-1\cdot j}_{\ \ l} \tag{10.1.26}$$

Proof. $D$-linear map has form (10.1.16) relative to basis $\overline{\overline{e}}$. Let $\overline{\overline{e}}'$ be another basis. Let

$$e'_i = e_j A^j_i \tag{10.1.27}$$

be transformation map basis $\overline{\overline{e}}$ to basis $\overline{\overline{e}}'$. Since linear map $f$ is the same, then

$$f(x) = e'_l f'^l_k x'^k \tag{10.1.28}$$

Let us substitute (9.2.8), (10.1.27) into equation (10.1.28)

$$f(x) = e_j A^j_l f'^l_k A^{-1\cdot k}_i x^i \tag{10.1.29}$$

Because vectors $e_j$ are linear independent and components of vector $x^i$ are arbitrary, the equation (10.1.26) follows from equation (10.1.29). Therefore, expression $f^r_k$ is tensor over field $F$. □

Definition 10.1.19. The set

$$\ker f = \{x \in D_1 : f(x) = 0\}$$

is called **kernel of linear map**

$$f : D_1 \to D_2$$

of division ring $D_1$ into division ring $D_2$. □

Theorem 10.1.20. *Kernel of linear map*

$$f : D_1 \to D_2$$

*is subgroup of additive group of division ring $D_1$.*

Proof. Let $a, b \in \ker f$. Then

$$f(a) = 0$$
$$f(b) = 0$$
$$f(a+b) = f(a) + f(b) = 0$$

Therefore, $a + b \in \ker f$. □

Definition 10.1.21. The linear map

$$f : D_1 \to D_2$$

of division ring $D_1$ into division ring $D_2$ is called **singular**, when

$$\ker f \neq \{0\}$$

□

Theorem 10.1.22. *Let $D$ be division ring of characteristic 0. Let $\overline{\overline{e}}$ be basis of division ring $D$ over center $Z(D)$ of division ring $D$. Let*

$$f : D \to D \quad f(x) = f_{s\cdot 0}\ x\ f_{s\cdot 1} \tag{10.1.30}$$

$$= f^{ij}\ e_i x e_j \tag{10.1.31}$$

$$g : D \to D \quad g(x) = g_{t\cdot 0}\ x\ g_{t\cdot 1} \tag{10.1.32}$$

$$= g^{ij}\ e_i x e_j \tag{10.1.33}$$

*be linear maps of division ring $D$. Map*

$$h(x) = gf(x) = g(f(x)) \tag{10.1.34}$$

*is linear map*

$$h(x) = h_{ts\cdot 0}\ x\ h_{ts\cdot 1} \tag{10.1.35}$$

$$= h^{pr}\ e_p\ x\ e_r \tag{10.1.36}$$

*where*

$$h_{ts\cdot 0} = g_{t\cdot 0}\ f_{s\cdot 0} \tag{10.1.37}$$

$$h_{ts\cdot 1} = f_{s\cdot 1}\ g_{t\cdot 1} \tag{10.1.38}$$

$$h^{pr} = g^{ij} f^{kl} C^p_{ik} C^r_{lj} \tag{10.1.39}$$

Proof. Map (10.1.34) is linear because

$$h(x+y) = g(f(x+y)) = g(f(x)+f(y)) = g(f(x)) + g(f(y)) = h(x) + h(y)$$

$$h(ax) = g(f(ax)) = g(af(x)) = ag(f(x)) = ah(x)$$

If we substitute (10.1.30) and (10.1.32) into (10.1.34), we get

$$h(x) = g_{t\cdot 0}\ f(x)\ g_{t\cdot 1} = g_{t\cdot 0}\ f_{s\cdot 0}\ x\ f_{s\cdot 1}\ g_{t\cdot 1} \tag{10.1.40}$$

Comparing (10.1.40) and (10.1.35), we get (10.1.37), (10.1.38).

If we substitute (10.1.31) and (10.1.33) into (10.1.34), we get

$$\begin{aligned} h(x) =& g^{ij}\ e_i\ f(x)\ e_j \\ =& g^{ij}\ e_i\ f^{kl}\ e_k\ x\ e_l e_j \\ =& g^{ij} f^{kl} C^p_{ik} C^r_{lj}\ e_p\ x\ e_r \end{aligned} \tag{10.1.41}$$

Comparing (10.1.41) and (10.1.36), we get (10.1.39). □

## 10.2. Polylinear Map of Division Ring

Definition 10.2.1. Let $R_1$, ..., $R_n$, $P$ be rings of characteristic 0. Let $S$ be module over ring $P$. Let $F$ be commutative ring which is for any $i$ is subring of center of ring $R_i$. Map

$$f : R_1 \times ... \times R_n \to S$$

is called **polylinear over commutative ring** $F$, if

$$f(p_1, ..., p_i + q_i, ..., p_n) = f(p_1, ..., p_i, ..., p_n) + f(p_1, ..., q_i, ..., p_n)$$

$$f(a_1, ..., ba_i, ..., a_n) = bf(a_1, ..., a_i, ..., a_n)$$

for any $i$, $1 \le i \le n$, and any $p_i$, $q_i \in R_i$, $b \in F$. Let us denote $\mathcal{L}(R_1, ..., R_n; S)$ set of polylinear maps of rings $R_1$, ..., $R_n$ into module $S$. □

Theorem 10.2.2. *Let $D$ be division ring of characteristic 0. Polylinear map*

$$f : D^n \to D, d = f(d_1, ..., d_n) \tag{10.2.1}$$

*has form*

$$d = f^n_{s\cdot 0}\ \sigma_s(d_1)\ f^n_{s\cdot 1}\ ...\ \sigma_s(d_n)\ f^n_{s\cdot n} \tag{10.2.2}$$

*$\sigma_s$ is a transposition of set of variables* $\{d_1, ..., d_n\}$

$$\sigma_s = \begin{pmatrix} d_1 & ... & d_n \\ \sigma_s(d_1) & ... & \sigma_s(d_n) \end{pmatrix}$$

Proof. We prove statement by induction on $n$.

When $n=1$ the statement of theorem is corollary of theorem 10.1.11. In such case we may identify[10.4]

$$f^1_{s\cdot p}=f_{s\cdot p}\quad p=0,1$$

Let statement of theorem be true for $n=k-1$. Then it is possible to represent map (10.2.1) as

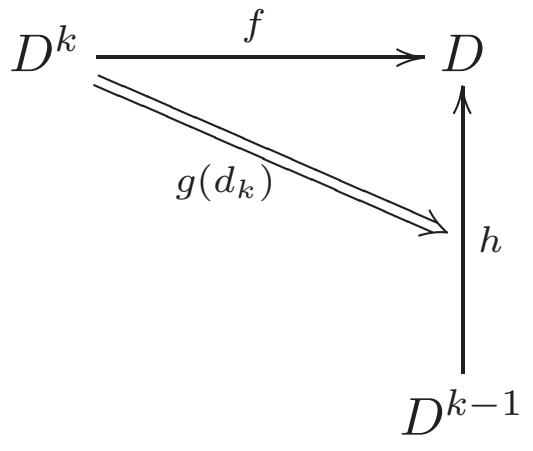

$$d=f(d_1,...,d_k)=g(d_k)(d_1,...,d_{k-1})$$

According to statement of induction pollinear map $h$ has form

$$d=h^{k-1}_{t\cdot 0}\ \sigma_t(d_1)\ h^{k-1}_{t\cdot 1}\ ...\ \sigma_t(d_{k-1})\ h^{k-1}_{t\cdot k-1}$$

According to construction $h=g(d_k)$. Therefore, expressions $h_{t\cdot p}$ are functions of $d_k$. Since $g(d_k)$ is linear map of $d_k$, then only one expression $h_{t\cdot p}$ is linear map of $d_k$, and rest expressions ${}_{t\cdot q}h$ do not depend on $d_k$.

Without loss of generality, assume $p=0$. According to equation (10.1.14) for given $t$

$$h^{k-1}_{t\cdot 0}=g_{tr\cdot 0}\ d_k\ g_{tr\cdot 1}$$

Assume $s=tr$. Let us define transposition $\sigma_s$ according to rule

$$\sigma_s=\sigma_{tr}=\begin{pmatrix}d_k & d_1 & ... & d_{k-1}\\ d_k & \sigma_t(d_1) & ... & \sigma_t(d_{k-1})\end{pmatrix}$$

Suppose

$$f^k_{tr\cdot q+1}=h^{k-1}_{t\cdot q}\quad q=1,...,k-1$$
$$f^k_{tr\cdot q}=g_{tr\cdot q}\quad q=0,1$$

We proved step of induction. □

Definition 10.2.3. Expression $f^n_{s\cdot p}$ in equation (10.2.2) is called **component of polylinear map** $f$. □

Theorem 10.2.4. *Let $D$ be division ring of characteristic* 0*. Let $\overline{\overline{e}}$ be basis in division ring $D$ over field $F\subset Z(D)$.* **Standard representation of polylinear map** *of division ring has form*

$$f(d_1,...,d_n)=f^{i_0...i_n}_t\ e_{i_0}\ \sigma_t(d_1)\ e_{i_1}\ ...\ \sigma_t(d_n)\ e_{i_n}\tag{10.2.3}$$

[10.4] In representation (10.2.2) we will use following rules.

- If range of any index is set consisting of one element, then we will omit corresponding index.
- If $n=1$, then $\sigma_s$ is identical transformation. We will not show such transformation in the expression.

*Index $t$ enumerates every possible transpositions $\sigma_t$ of the set of variables $\{d_1, ..., d_n\}$. Expression $f_t^{i_0...i_n}$ in equation (10.2.3) is called* **standard component of polylinear map** $f$.

Proof. Components of polylinear map $f$ have expansion

$$f^n_{s\cdot p} = e_i f^{ni}_{s\cdot p} \tag{10.2.4}$$

relative to basis $\overline{\overline{e}}$. If we substitute (10.2.4) into (10.2.2), we get

$$d = f^{nj_1}_{s\cdot 0}\ e_{j_1}\ \sigma_s(d_1)\ f^{nj_2}_{s\cdot 1}\ e_{j_2}\ ...\ \sigma_s(d_n)\ f^{nj_n}_{s\cdot n}\ e_{j_n} \tag{10.2.5}$$

Consider expression

$$f_t^{j_0...j_n} = f^{nj_1}_{s\cdot 0}\ ...f^{nj_n}_{s\cdot n} \tag{10.2.6}$$

The right-hand side is supposed to be the sum of the terms with the index $s$, for which the transposition $\sigma_s$ is the same. Each such sum has a unique index $t$. If we substitute expression (10.2.6) into equation (10.2.5) we get equation (10.2.3). □

Theorem 10.2.5. *Let $\overline{\overline{e}}$ be basis of division ring $D$ over field $F\subset Z(D)$. Polylinear map (10.2.1) can be represented as $D$-valued form of degree $n$ over field $F\subset Z(D)$*[10.5]

$$f(a_1, ..., a_n) = a_1^{i_1} ...a_n^{i_n} f_{i_1...i_n} \tag{10.2.7}$$

*where*

$$a_j = e_i a_j^i$$

$$f_{i_1...i_n} = f(e_{i_1}, ..., e_{i_n}) \tag{10.2.8}$$

*and values ${}_{i_1...i_n}f$ are coordinates of $D$-valued covariant tensor over field $F$.*

Proof. According to the definition 10.2.1, the equation (10.2.7) follows from the chain of equations

$$f(a_1, ..., a_n) = f(e_{i_1}a_1^{i_1}, ..., e_{i_n}a_n^{i_n}) = a_1^{i_1} ...a_n^{i_n} f(e_{i_1}, ..., e_{i_n})$$

Let $\overline{\overline{e}}'$ be another basis. Let

$$e'_i = e_j A^j_i \tag{10.2.9}$$

be transformation, mapping basis $\overline{\overline{e}}$ into basis $\overline{\overline{e}}'$. From equations (10.2.9) and (10.2.8) it follows

$$\begin{aligned} f'_{i_1...i_n} &= f(e'_{i_1}, ..., e'_{i_n}) \\ &= f(e_{j_1}A^{j_1}_{i_1}, ..., e'_{j_n}A^{j_n}_{i_n}) \\ &= A^{j_1}_{i_1}...A^{j_n}_{i_n} f(e_{j_1}, ..., e_{j_n}) \\ &= A^{j_1}_{i_1}...A^{j_n}_{i_n} f_{j_1...j_n} \end{aligned} \tag{10.2.10}$$

From equation (10.2.10) the tensor law of transformation of coordinates of polylinear map follows. From equation (10.2.10) and theorem 9.2.3 it follows that value of the map $f(a_1, ..., a_n)$ does not depend from choice of basis. □

Polylinear map (10.2.1) is **symmetric**, if

$$f(d_1, ..., d_n) = f(\sigma(d_1), ..., \sigma(d_n))$$

for any transposition $\sigma$ of set $\{d_1, ..., d_n\}$.

[10.5] We proved the theorem by analogy with theorem in [3], p. 107, 108

THEOREM 10.2.6. *If polyadditive map $f$ is symmetric, then*

$$f_{i_1,...,i_n}=f_{\sigma(i_1),...,\sigma(i_n)} \tag{10.2.11}$$

PROOF. Equation (10.2.11) follows from equation

$$\begin{aligned}a_1^{i_1}...a_n^{i_n}f_{i_1...i_n}=&f(a_1,...,a_n)\\=&f(\sigma(a_1),...,\sigma(a_n))\\=&a_1^{i_1}...a_n^{i_n}f_{\sigma(i_1)...\sigma(i_n)}\end{aligned}$$

□

Polylinear map (10.2.1) is **skew symmetric**, if

$$f(d_1,...,d_n)=|\sigma|f(\sigma(d_1),...,\sigma(d_n))$$

for any transposition $\sigma$ of set $\{d_1,...,d_n\}$. Here

$$|\sigma|=\begin{cases}1 & \text{transposition } \sigma \text{ even}\\-1 & \text{transposition } \sigma \text{ odd}\end{cases}$$

THEOREM 10.2.7. *If polylinear map $f$ is skew symmetric, then*

$$f_{i_1,...,i_n}=|\sigma|f_{\sigma(i_1),...,\sigma(i_n)} \tag{10.2.12}$$

PROOF. Equation (10.2.12) follows from equation

$$\begin{aligned}a_1^{i_1}...a_n^{i_n}f_{i_1...i_n}=&f(a_1,...,a_n)\\=&|\sigma|f(\sigma(a_1),...,\sigma(a_n))\\=&a_1^{i_1}...a_n^{i_n}|\sigma|f_{\sigma(i_1)...\sigma(i_n)}\end{aligned}$$

□

THEOREM 10.2.8. *Coordinates of the polylinear over field $F$ map (10.2.1) and its components relative basis $\overline{\overline{e}}$ satisfy to the equation*

$$f_{j_1...j_n}=f_t^{i_0...i_n}C^{k_1}_{i_0\sigma_t(j_1)}C^{l_1}_{k_1i_1}...B^{k_n}_{l_{n-1}\sigma_t(j_n)}C^{l_n}_{k_ni_n}e_{l_n} \tag{10.2.13}$$

$$f^p_{j_1...j_n}=f_t^{i_0...i_n}C^{k_1}_{i_0\sigma_t(j_1)}C^{l_1}_{k_1i_1}...C^{k_n}_{l_{n-1}\sigma_t(j_n)}C^{p}_{k_ni_n} \tag{10.2.14}$$

PROOF. In equation (10.2.3), we assume

$$d_i=e_{j_i}d_i^{j_i}$$

Then equation (10.2.3) gets form

$$\begin{aligned}f(d_1,...,d_n)&=f_t^{i_0...i_n}\ e_{i_0}\sigma_t(d_1^{j_1}e_{j_1})e_{i_1}...\sigma_t(d_n^{j_n}e_{j_n})e_{i_n}\\&=d_1^{j_1}...d_n^{j_n}f_t^{i_0...i_n}e_{i_0}\sigma_t(e_{j_1})e_{i_1}...\sigma_t(e_{j_n})e_{i_n}\\&=d_1^{j_1}...d_n^{j_n}f_t^{i_0...i_n}C^{k_1}_{i_0\sigma_t(j_1)}C^{l_1}_{k_1i_1}\\&\quad...C^{k_n}_{l_{n-1}\sigma_t(j_n)}C^{l_n}_{k_ni_n}e_{l_n}\end{aligned} \tag{10.2.15}$$

From equation (10.2.7) it follows that

$$f(a_1,...,a_n)=e_pf^p_{i_1...i_n}a_1^{i_1}...a_n^{i_n} \tag{10.2.16}$$

Equation (10.2.13) follows from comparison of equations (10.2.15) and (10.2.7). Equation (10.2.14) follows from comparison of equations (10.2.15) and (10.2.16).

□

CHAPTER 11

# Quaternion Algebra

## 11.1. Linear Function of Complex Field

THEOREM 11.1.1 (the Cauchy-Riemann equations). *Let us consider complex field $C$ as two-dimensional algebra over real field. Let*

$$e_{C\cdot 0} = 1 \quad e_{C\cdot 1} = i \tag{11.1.1}$$

*be the basis of algebra $C$. Then in this basis product has form*

$$e_{C\cdot 1}^2 = -e_{C\cdot 0} \tag{11.1.2}$$

*and structural constants have form*

$$\begin{aligned} C_{C\cdot}{}^{0}_{00} &= 1 & C_{C\cdot}{}^{1}_{01} &= 1 \\ C_{C\cdot}{}^{1}_{10} &= 1 & C_{C\cdot}{}^{0}_{11} &= -1 \end{aligned} \tag{11.1.3}$$

*Matrix of linear function*

$$y^i = x^j f^i_j$$

*of complex field over real field satisfies relationship*

$$f^0_0 = f^1_1 \tag{11.1.4}$$

$$f^1_0 = -f^0_1 \tag{11.1.5}$$

PROOF. Equations (11.1.2) and (11.1.3) follow from equation $i^2 = -1$. Using equation (10.1.17) we get relationships

$$f^0_0 = f^{kr} C_{C\cdot}{}^{p}_{k0} C_{C\cdot}{}^{0}_{pr} = f^{0r} C_{C\cdot}{}^{0}_{00} C_{C\cdot}{}^{0}_{0r} + f^{1r} C_{C\cdot}{}^{1}_{10} C_{C\cdot}{}^{0}_{1r} = f^{00} - f^{11} \tag{11.1.6}$$

$$f^1_0 = f^{kr} C_{C\cdot}{}^{p}_{k0} C_{C\cdot}{}^{1}_{pr} = f^{0r} C_{C\cdot}{}^{0}_{00} C_{C\cdot}{}^{1}_{0r} + f^{1r} C_{C\cdot}{}^{1}_{10} C_{C\cdot}{}^{1}_{1r} = f^{01} + f^{10} \tag{11.1.7}$$

$$f^0_1 = f^{kr} C_{C\cdot}{}^{p}_{k1} C_{C\cdot}{}^{0}_{pr} = f^{0r} C_{C\cdot}{}^{1}_{01} C_{C\cdot}{}^{0}_{1r} + f^{1r} C_{C\cdot}{}^{0}_{11} C_{C\cdot}{}^{0}_{0r} = -f^{01} - f^{10} \tag{11.1.8}$$

$$f^1_1 = f^{kr} C_{C\cdot}{}^{p}_{k1} C_{C\cdot}{}^{1}_{pr} = f^{0r} C_{C\cdot}{}^{1}_{01} C_{C\cdot}{}^{1}_{1r} + f^{1r} C_{C\cdot}{}^{0}_{11} C_{C\cdot}{}^{1}_{0r} = f^{00} - f^{11} \tag{11.1.9}$$

(11.1.4) follows from equations (11.1.6) and (11.1.9). (11.1.5) follows from equations (11.1.7) and (11.1.8). □

## 11.2. Quaternion Algebra

In this paper I explore the set of quaternion algebras defined in [14].

Definition 11.2.1. Let $F$ be field. Extension field $F(i,j,k)$ is called **the quaternion algebra $E(F,a,b)$ over the field $F$**[11.1] if multiplication in algebra $E$ is defined according to rule

(11.2.1)

| | $i$ | $j$ | $k$ |
|---|---|---|---|
| $i$ | $a$ | $k$ | $aj$ |
| $j$ | $-k$ | $b$ | $-bi$ |
| $k$ | $-aj$ | $bi$ | $-ab$ |

where $a$, $b\in F$, $ab\neq 0$. □

Elements of the algebra $E(F,a,b)$ have form

$$x=x^0+x^1i+x^2j+x^3k$$

where $x^i\in F$, $i=0$, $1$, $2$, $3$. Quaternion

$$\overline{x}=x^0-x^1i-x^2j-x^3k$$

is called conjugate to the quaternion $x$. We define **the norm of the quaternion** $x$ using equation

$$|x|^2=x\overline{x}=(x^0)^2-a(x^1)^2-b(x^2)^2+ab(x^3)^2 \tag{11.2.2}$$

From equation (11.2.2), it follows that $E(F,a,b)$ is algebra with division only when $a<0$, $b<0$. In this case we can renorm basis such that $a=-1$, $b=-1$.

We use symbol $E(F)$ to denote the quaternion division algebra $E(F,-1,-1)$ over the field $F$. We will use notation $E(R,-1,-1)$. Multiplication in quaternion algebra $E(F)$ is defined according to rule

(11.2.3)

| | $i$ | $j$ | $k$ |
|---|---|---|---|
| $i$ | $-1$ | $k$ | $-j$ |
| $j$ | $-k$ | $-1$ | $i$ |
| $k$ | $j$ | $-i$ | $-1$ |

In algebra $E(F)$, the norm of the quaternion has form

$$|x|^2=x\overline{x}=(x^0)^2+(x^1)^2+(x^2)^2+(x^3)^2 \tag{11.2.4}$$

In this case inverse element has form

$$x^{-1}=|x|^{-2}\overline{x} \tag{11.2.5}$$

The inner automorphism of quaternion algebra $H$[11.2]

$$\begin{gathered}p\rightarrow qpq^{-1}\\ q(ix+jy+kz)q^{-1}=ix'+jy'+kz'\end{gathered} \tag{11.2.6}$$

[11.1] I follow definition from [14].
[11.2] See [15], p. 643.

describes the rotation of the vector with coordinates $x$, $y$, $z$. The norm of quaternion $q$ is irrelevant, although usually we assume $|q| = 1$. If $q$ is written as sum of scalar and vector

$$q = \cos\alpha + (ia + jb + kc)\sin\alpha \quad a^2 + b^2 + c^2 = 1$$

then (11.2.6) is a rotation of the vector $(x, y, z)$ about the vector $(a, b, c)$ through an angle $2\alpha$.

## 11.3. Linear Function of Quaternion Algebra

THEOREM 11.3.1. *Let*

$$e_0 = 1 \quad e_1 = i \quad e_2 = j \quad e_3 = k \tag{11.3.1}$$

*be basis of quaternion algebra $H$. Then in the basis* (11.3.1)*, structural constants have form*

$$\begin{matrix} C^0_{00} = 1 & C^1_{01} = 1 & C^2_{02} = 1 & C^3_{03} = 1 \\ C^1_{10} = 1 & C^0_{11} = -1 & C^3_{12} = 1 & C^2_{13} = -1 \\ C^2_{20} = 1 & C^3_{21} = -1 & C^0_{22} = -1 & C^1_{23} = 1 \\ C^3_{30} = 1 & C^2_{31} = 1 & C^1_{32} = -1 & C^0_{33} = -1 \end{matrix}$$

PROOF. Value of structural constants follows from multiplication table (11.2.3). □

Since calculations in this section get a lot of space, I put in one place references to theorems in this section.

**Theorem 11.3.2, page 125:** the definition of coordinates of linear map of quaternion algebra $H$ using standard components of this map.

**Equation** (11.3.22)**, page 128:** matrix form of dependence of coordinates of linear map of quaternion algebra $H$ from standard components of this map.

**Equation** (11.3.23)**, page 129:** matrix form of dependence of standard components of linear map of quaternion algebra $H$ from coordinates of this map.

**Theorem 11.3.4, page 131:** dependence of standard components of a linear map of quaternion algebra $H$ from coordinates of this map.

THEOREM 11.3.2. *Standard components of linear function of quaternion algebra $H$ relative to basis* (11.3.1) *and coordinates of corresponding linear map satisfy relationship*

$$\begin{cases} f^0_0 = f^{00} - f^{11} - f^{22} - f^{33} \\ f^1_1 = f^{00} - f^{11} + f^{22} + f^{33} \\ f^2_2 = f^{00} + f^{11} - f^{22} + f^{33} \\ f^3_3 = f^{00} + f^{11} + f^{22} - f^{33} \end{cases} \tag{11.3.2}$$

$$\left\{\begin{array}{l} f_0^1 = \ \ f^{01} + f^{10} + f^{23} - f^{32} \\ f_1^0 = -f^{01} - f^{10} + f^{23} - f^{32} \\ f_2^3 = -f^{01} + f^{10} - f^{23} - f^{32} \\ f_3^2 = \ \ f^{01} - f^{10} - f^{23} - f^{32} \end{array}\right. \tag{11.3.3}$$

$$\left\{\begin{array}{l} f_0^2 = \ \ f^{02} - f^{13} + f^{20} + f^{31} \\ f_1^3 = \ \ f^{02} - f^{13} - f^{20} - f^{31} \\ f_2^0 = -f^{02} - f^{13} - f^{20} + f^{31} \\ f_3^1 = -f^{02} - f^{13} + f^{20} - f^{31} \end{array}\right. \tag{11.3.4}$$

$$\left\{\begin{array}{l} f_0^3 = \ \ f^{03} + f^{12} - f^{21} + f^{30} \\ f_1^2 = -f^{03} - f^{12} - f^{21} + f^{30} \\ f_2^1 = \ \ f^{03} - f^{12} - f^{21} - f^{30} \\ f_3^0 = -f^{03} + f^{12} - f^{21} - f^{30} \end{array}\right. \tag{11.3.5}$$

Proof. Using equation (10.1.17) we get relationships

$$\begin{aligned} f_0^0 &= f^{kr} C_{k0}^p C_{pr}^0 \\ &= f^{00} C_{00}^0 C_{00}^0 + f^{11} C_{10}^1 C_{11}^0 + f^{22} C_{20}^2 C_{22}^0 + f^{33} C_{30}^3 C_{33}^0 \\ &= f^{00} - f^{11} - f^{22} - f^{33} \end{aligned} \tag{11.3.6}$$

$$\begin{aligned} f_0^1 &= f^{kr} C_{k0}^p C_{pr}^1 \\ &= f^{01} C_{00}^0 C_{01}^1 + f^{10} C_{10}^1 C_{10}^1 + f^{23} C_{20}^2 C_{23}^1 + f^{32} C_{30}^3 C_{32}^1 \\ &= f^{01} + f^{10} + f^{23} - f^{32} \end{aligned} \tag{11.3.7}$$

$$\begin{aligned} f_0^2 &= f^{kr} C_{k0}^p C_{pr}^2 \\ &= f^{02} C_{00}^0 C_{02}^2 + f^{13} C_{10}^1 C_{13}^2 + f^{20} C_{20}^2 C_{20}^2 + f^{31} C_{30}^3 C_{31}^2 \\ &= f^{02} - f^{13} + f^{20} + f^{31} \end{aligned} \tag{11.3.8}$$

$$\begin{aligned} f_0^3 &= f^{kr} C_{k0}^p C_{pr}^3 \\ &= f^{03} C_{00}^0 C_{03}^3 + f^{12} C_{10}^1 C_{12}^3 + f^{21} C_{20}^2 C_{21}^3 + f^{30} C_{30}^3 C_{30}^3 \\ &= f^{03} + f^{12} - f^{21} + f^{30} \end{aligned} \tag{11.3.9}$$

$$\begin{aligned} f_{1}^{0} &= f^{kr} C_{k1}^{p} C_{pr}^{0} \\ &= f^{01} C_{01}^{1} C_{11}^{0} + f^{10} C_{11}^{0} C_{00}^{0} + f^{23} C_{21}^{3} C_{33}^{0} + f^{32} C_{31}^{2} C_{22}^{0} \\ &= -f^{01} - f^{10} + f^{23} - f^{32} \end{aligned} \tag{11.3.10}$$

$$\begin{aligned} f_{1}^{1} &= f^{kr} C_{k1}^{p} C_{pr}^{1} \\ &= f^{00} C_{01}^{1} C_{10}^{1} + f^{11} C_{11}^{0} C_{01}^{1} + f^{22} C_{21}^{3} C_{32}^{1} + f^{33} C_{31}^{2} C_{23}^{1} \\ &= f^{00} - f^{11} + f^{22} + f^{33} \end{aligned} \tag{11.3.11}$$

$$\begin{aligned} f_{1}^{2} &= f^{kr} C_{k1}^{p} C_{pr}^{2} \\ &= f^{03} C_{01}^{1} C_{13}^{2} + f^{12} C_{11}^{0} C_{02}^{2} + f^{21} C_{21}^{3} C_{31}^{2} + f^{30} C_{31}^{2} C_{20}^{2} \\ &= -f^{03} - f^{12} - f^{21} + f^{30} \end{aligned} \tag{11.3.12}$$

$$\begin{aligned} f_{1}^{3} &= f^{kr} C_{k1}^{p} C_{pr}^{3} \\ &= f^{02} C_{01}^{1} C_{12}^{3} + f^{13} C_{11}^{0} C_{03}^{3} + f^{20} C_{21}^{3} C_{30}^{3} + f^{31} C_{31}^{2} C_{21}^{3} \\ &= f^{02} - f^{13} - f^{20} - f^{31} \end{aligned} \tag{11.3.13}$$

$$\begin{aligned} f_{2}^{0} &= f^{kr} C_{k2}^{p} C_{pr}^{0} \\ &= f^{02} C_{02}^{2} C_{22}^{0} + f^{13} C_{12}^{3} C_{33}^{0} + f^{20} C_{22}^{0} C_{00}^{0} + f^{31} C_{32}^{1} C_{11}^{0} \\ &= -f^{02} - f^{13} - f^{20} + f^{31} \end{aligned} \tag{11.3.14}$$

$$\begin{aligned} f_{2}^{1} &= f^{kr} C_{k2}^{p} C_{pr}^{1} \\ &= f^{03} C_{02}^{2} C_{23}^{1} + f^{12} C_{12}^{3} C_{32}^{1} + f^{21} C_{22}^{0} C_{01}^{1} + f^{30} C_{32}^{1} C_{10}^{1} \\ &= f^{03} - f^{12} - f^{21} - f^{30} \end{aligned} \tag{11.3.15}$$

$$\begin{aligned} f_{2}^{2} &= f^{kr} C_{k2}^{p} C_{pr}^{2} \\ &= f^{00} C_{02}^{2} C_{20}^{2} + f^{11} C_{12}^{3} C_{31}^{2} + f^{22} C_{22}^{0} C_{02}^{2} + f^{33} C_{32}^{1} C_{13}^{2} \\ &= f^{00} + f^{11} - f^{22} + f^{33} \end{aligned} \tag{11.3.16}$$

$$\begin{aligned} f_{2}^{3} &= f^{kr} C_{k2}^{p} C_{pr}^{3} \\ &= f^{01} C_{02}^{2} C_{21}^{3} + f^{10} C_{12}^{3} C_{30}^{3} + f^{23} C_{22}^{0} C_{03}^{3} + f^{32} C_{32}^{1} C_{12}^{3} \\ &= -f^{01} + f^{10} - f^{23} - f^{32} \end{aligned} \tag{11.3.17}$$

$$\begin{aligned} f_3^0 &= f^{kr}C_{k3}^pC_{pr}^0 \\ &= f^{03}C_{03}^3C_{33}^0 + f^{12}C_{13}^2C_{22}^0 + f^{21}C_{23}^1C_{11}^0 + f^{30}C_{33}^0C_{00}^0 \\ &= -f^{03} + f^{12} - f^{21} - f^{30} \end{aligned} \tag{11.3.18}$$

$$\begin{aligned} f_3^1 &= f^{kr}C_{k3}^pC_{pr}^1 \\ &= f^{02}C_{03}^3C_{32}^1 + f^{13}C_{13}^2C_{23}^1 + f^{20}C_{23}^1C_{10}^1 + f^{31}C_{33}^0C_{01}^1 \\ &= -f^{02} - f^{13} + f^{20} - f^{31} \end{aligned} \tag{11.3.19}$$

$$\begin{aligned} f_3^2 &= f^{kr}C_{k3}^pC_{pr}^2 \\ &= f^{01}C_{03}^3C_{31}^2 + f^{10}C_{13}^2C_{20}^2 + f^{23}C_{23}^1C_{13}^2 + f^{32}C_{33}^0C_{02}^2 \\ &= f^{01} - f^{10} - f^{23} - f^{32} \end{aligned} \tag{11.3.20}$$

$$\begin{aligned} f_3^3 &= f^{kr}C_{k3}^pC_{pr}^3 \\ &= f^{00}C_{03}^3C_{30}^3 + f^{11}C_{13}^2C_{21}^3 + f^{22}C_{23}^1C_{12}^3 + f^{33}C_{33}^0C_{03}^3 \\ &= f^{00} + f^{11} + f^{22} - f^{33} \end{aligned} \tag{11.3.21}$$

Equations (11.3.6), (11.3.11), (11.3.16), (11.3.21) form the system of linear equations (11.3.2).

Equations (11.3.7), (11.3.10), (11.3.17), (11.3.20) form the system of linear equations (11.3.3).

Equations (11.3.8), (11.3.13), (11.3.14), (11.3.19) form the system of linear equations (11.3.4).

Equations (11.3.9), (11.3.12), (11.3.15), (11.3.18) form the system of linear equations (11.3.5). □

Theorem 11.3.3. *Consider quaternion algebra $H$ with the basis (11.3.1). Standard components of linear function over field $F$ and coordinates of this function over field $F$ satisfy relationship*

$$\begin{pmatrix} f_0^0 & f_1^0 & f_2^0 & f_3^0 \\ f_1^1 & -f_1^0 & f_3^1 & -f_2^1 \\ f_2^2 & -f_3^2 & -f_0^2 & f_1^2 \\ f_3^3 & f_2^3 & -f_1^3 & -f_0^3 \end{pmatrix} = \begin{pmatrix} 1 & -1 & -1 & -1 \\ 1 & -1 & 1 & 1 \\ 1 & 1 & -1 & 1 \\ 1 & 1 & 1 & -1 \end{pmatrix} \begin{pmatrix} f^{00} & -f^{01} & -f^{02} & -f^{03} \\ f^{11} & f^{10} & f^{13} & -f^{12} \\ f^{22} & -f^{23} & f^{20} & f^{21} \\ f^{33} & f^{32} & -f^{31} & f^{30} \end{pmatrix} \tag{11.3.22}$$

$$\begin{pmatrix} f^{00} & -f^{01} & -f^{02} & -f^{03} \\ f^{11} & f^{10} & f^{13} & -f^{12} \\ f^{22} & -f^{23} & f^{20} & f^{21} \\ f^{33} & f^{32} & -f^{31} & f^{30} \end{pmatrix} = \frac{1}{4}\begin{pmatrix} 1 & 1 & 1 & 1 \\ -1 & -1 & 1 & 1 \\ -1 & 1 & -1 & 1 \\ -1 & 1 & 1 & -1 \end{pmatrix}\begin{pmatrix} f_0^0 & f_1^0 & f_2^0 & f_3^0 \\ f_1^1 & -f_1^0 & f_3^1 & -f_2^1 \\ f_2^2 & -f_3^2 & -f_0^2 & f_1^2 \\ f_3^3 & f_2^3 & -f_1^3 & -f_0^3 \end{pmatrix} \tag{11.3.23}$$

*where*

$$\begin{pmatrix} 1 & -1 & -1 & -1 \\ 1 & -1 & 1 & 1 \\ 1 & 1 & -1 & 1 \\ 1 & 1 & 1 & -1 \end{pmatrix}^{-1} = \frac{1}{4}\begin{pmatrix} 1 & 1 & 1 & 1 \\ -1 & -1 & 1 & 1 \\ -1 & 1 & -1 & 1 \\ -1 & 1 & 1 & -1 \end{pmatrix}$$

PROOF. Let us write the system of linear equations (11.3.2) as product of matrices

$$\begin{pmatrix} f_0^0 \\ f_1^1 \\ f_2^2 \\ f_3^3 \end{pmatrix} = \begin{pmatrix} 1 & -1 & -1 & -1 \\ 1 & -1 & 1 & 1 \\ 1 & 1 & -1 & 1 \\ 1 & 1 & 1 & -1 \end{pmatrix}\begin{pmatrix} f^{00} \\ f^{11} \\ f^{22} \\ f^{33} \end{pmatrix} \tag{11.3.24}$$

Let us write the system of linear equations (11.3.3) as product of matrices

$$\begin{pmatrix} f_0^1 \\ f_1^0 \\ f_2^3 \\ f_3^2 \end{pmatrix} = \begin{pmatrix} 1 & 1 & 1 & -1 \\ -1 & -1 & 1 & -1 \\ -1 & 1 & -1 & -1 \\ 1 & -1 & -1 & -1 \end{pmatrix}\begin{pmatrix} f^{01} \\ f^{10} \\ f^{23} \\ f^{32} \end{pmatrix} \tag{11.3.25}$$

From the equation (11.3.25), it follows that

$$\begin{pmatrix} -f_{\mathbf{1}}^{\mathbf{0}} \\ f_{\mathbf{1}}^{\mathbf{0}} \\ f_{\mathbf{2}}^{\mathbf{3}} \\ -f_{\mathbf{3}}^{\mathbf{2}} \end{pmatrix} = \begin{pmatrix} -1 & -1 & -1 & 1 \\ -1 & -1 & 1 & -1 \\ -1 & 1 & -1 & -1 \\ -1 & 1 & 1 & 1 \end{pmatrix} \begin{pmatrix} f^{\mathbf{01}} \\ f^{\mathbf{10}} \\ f^{\mathbf{23}} \\ f^{\mathbf{32}} \end{pmatrix}$$

$$= \begin{pmatrix} 1 & -1 & 1 & 1 \\ 1 & -1 & -1 & -1 \\ 1 & 1 & 1 & -1 \\ 1 & 1 & -1 & 1 \end{pmatrix} \begin{pmatrix} -f^{\mathbf{01}} \\ f^{\mathbf{10}} \\ -f^{\mathbf{23}} \\ f^{\mathbf{32}} \end{pmatrix}$$

$$\begin{pmatrix} f_{\mathbf{1}}^{\mathbf{0}} \\ -f_{\mathbf{1}}^{\mathbf{0}} \\ -f_{\mathbf{3}}^{\mathbf{2}} \\ f_{\mathbf{2}}^{\mathbf{3}} \end{pmatrix} = \begin{pmatrix} 1 & -1 & -1 & -1 \\ 1 & -1 & 1 & 1 \\ 1 & 1 & -1 & 1 \\ 1 & 1 & 1 & -1 \end{pmatrix} \begin{pmatrix} -f^{\mathbf{01}} \\ f^{\mathbf{10}} \\ -f^{\mathbf{23}} \\ f^{\mathbf{32}} \end{pmatrix} \tag{11.3.26}$$

Let us write the system of linear equations (11.3.4) as product of matrices

$$\begin{pmatrix} f_{\mathbf{0}}^{\mathbf{2}} \\ f_{\mathbf{1}}^{\mathbf{3}} \\ f_{\mathbf{2}}^{\mathbf{0}} \\ f_{\mathbf{3}}^{\mathbf{1}} \end{pmatrix} = \begin{pmatrix} 1 & -1 & 1 & 1 \\ 1 & -1 & -1 & -1 \\ -1 & -1 & -1 & 1 \\ -1 & -1 & 1 & -1 \end{pmatrix} \begin{pmatrix} f^{\mathbf{02}} \\ f^{\mathbf{13}} \\ f^{\mathbf{20}} \\ f^{\mathbf{31}} \end{pmatrix} \tag{11.3.27}$$

From the equation (11.3.27), it follows that

$$\begin{pmatrix} -f_{\mathbf{0}}^{\mathbf{2}} \\ -f_{\mathbf{1}}^{\mathbf{3}} \\ f_{\mathbf{2}}^{\mathbf{0}} \\ f_{\mathbf{3}}^{\mathbf{1}} \end{pmatrix} = \begin{pmatrix} -1 & 1 & -1 & -1 \\ -1 & 1 & 1 & 1 \\ -1 & -1 & -1 & 1 \\ -1 & -1 & 1 & -1 \end{pmatrix} \begin{pmatrix} f^{\mathbf{02}} \\ f^{\mathbf{13}} \\ f^{\mathbf{20}} \\ f^{\mathbf{31}} \end{pmatrix}$$

$$= \begin{pmatrix} 1 & 1 & -1 & 1 \\ 1 & 1 & 1 & -1 \\ 1 & -1 & -1 & -1 \\ 1 & -1 & 1 & 1 \end{pmatrix} \begin{pmatrix} -f^{\mathbf{02}} \\ f^{\mathbf{13}} \\ f^{\mathbf{20}} \\ -f^{\mathbf{31}} \end{pmatrix}$$

$$\begin{pmatrix} f^{0}_{2} \\ f^{1}_{3} \\ -f^{2}_{0} \\ -f^{3}_{1} \end{pmatrix} = \begin{pmatrix} 1 & -1 & -1 & -1 \\ 1 & -1 & 1 & 1 \\ 1 & 1 & -1 & 1 \\ 1 & 1 & 1 & -1 \end{pmatrix} \begin{pmatrix} -f^{02} \\ f^{13} \\ f^{20} \\ -f^{31} \end{pmatrix} \tag{11.3.28}$$

Let us write the system of linear equations (11.3.5) as product of matrices

$$\begin{pmatrix} f^{3}_{0} \\ f^{2}_{1} \\ f^{1}_{2} \\ f^{0}_{3} \end{pmatrix} = \begin{pmatrix} 1 & 1 & -1 & 1 \\ -1 & -1 & -1 & 1 \\ 1 & -1 & -1 & -1 \\ -1 & 1 & -1 & -1 \end{pmatrix} \begin{pmatrix} f^{03} \\ f^{12} \\ f^{21} \\ f^{30} \end{pmatrix} \tag{11.3.29}$$

From the equation (11.3.29), it follows that

$$\begin{pmatrix} -f^{3}_{0} \\ f^{2}_{1} \\ -f^{1}_{2} \\ f^{0}_{3} \end{pmatrix} = \begin{pmatrix} -1 & -1 & 1 & -1 \\ -1 & -1 & -1 & 1 \\ -1 & 1 & 1 & 1 \\ -1 & 1 & -1 & -1 \end{pmatrix} \begin{pmatrix} f^{03} \\ f^{12} \\ f^{21} \\ f^{30} \end{pmatrix}$$

$$= \begin{pmatrix} 1 & 1 & 1 & -1 \\ 1 & 1 & -1 & 1 \\ 1 & -1 & 1 & 1 \\ 1 & -1 & -1 & -1 \end{pmatrix} \begin{pmatrix} -f^{03} \\ -f^{12} \\ f^{21} \\ f^{30} \end{pmatrix}$$

$$\begin{pmatrix} f^{0}_{3} \\ -f^{1}_{2} \\ f^{2}_{1} \\ -f^{3}_{0} \end{pmatrix} = \begin{pmatrix} 1 & -1 & -1 & -1 \\ 1 & -1 & 1 & 1 \\ 1 & 1 & -1 & 1 \\ 1 & 1 & 1 & -1 \end{pmatrix} \begin{pmatrix} -f^{03} \\ -f^{12} \\ f^{21} \\ f^{30} \end{pmatrix} \tag{11.3.30}$$

We join equations (11.3.24), (11.3.26), (11.3.28), (11.3.30) into equation (11.3.22). □

Theorem 11.3.4. *Standard components of linear function of quaternion algebra $H$ relative to basis (11.3.1) and coordinates of corresponding linear map satisfy relationship*

$$\begin{cases} 4f^{00} = f^{0}_{0} + f^{1}_{1} + f^{2}_{2} + f^{3}_{3} \\ 4f^{11} = -f^{0}_{0} - f^{1}_{1} + f^{2}_{2} + f^{3}_{3} \\ 4f^{22} = -f^{0}_{0} + f^{1}_{1} - f^{2}_{2} + f^{3}_{3} \\ 4f^{33} = -f^{0}_{0} + f^{1}_{1} + f^{2}_{2} - f^{3}_{3} \end{cases} \tag{11.3.31}$$

$$\begin{cases} 4f^{10} = -f^{0}_{1} + f^{0}_{1} - f^{2}_{3} + f^{3}_{2} \\ 4f^{01} = -f^{0}_{1} + f^{0}_{1} + f^{2}_{3} - f^{3}_{2} \\ 4f^{32} = -f^{0}_{1} - f^{0}_{1} - f^{2}_{3} - f^{3}_{2} \\ 4f^{23} = f^{0}_{1} + f^{0}_{1} - f^{2}_{3} - f^{3}_{2} \end{cases} \tag{11.3.32}$$

$$\begin{cases} 4f^{20} = -f^{0}_{2} + f^{1}_{3} + f^{2}_{0} - f^{3}_{1} \\ 4f^{31} = f^{0}_{2} - f^{1}_{3} + f^{2}_{0} - f^{3}_{1} \\ 4f^{02} = -f^{0}_{2} - f^{1}_{3} + f^{2}_{0} + f^{3}_{1} \\ 4f^{13} = -f^{0}_{2} - f^{1}_{3} - f^{2}_{0} - f^{3}_{1} \end{cases} \tag{11.3.33}$$

$$\begin{cases} 4f^{30} = -f^{0}_{3} - f^{1}_{2} + f^{2}_{1} + f^{3}_{0} \\ 4f^{21} = -f^{0}_{3} - f^{1}_{2} - f^{2}_{1} - f^{3}_{0} \\ 4f^{12} = f^{0}_{3} - f^{1}_{2} - f^{2}_{1} + f^{3}_{0} \\ 4f^{03} = -f^{0}_{3} + f^{1}_{2} - f^{2}_{1} + f^{3}_{0} \end{cases} \tag{11.3.34}$$

Proof. We get systems of linear equations (11.3.31), (11.3.32), (11.3.33), (11.3.34) as the product of matrices in equation (11.3.23). □

Theorem 11.3.5. *We can identify quaternion*

$$a = a^{0} + a^{1}i + a^{2}j + a^{3}k \tag{11.3.35}$$

*and matrix*

$$J_a = \begin{pmatrix} a^{0} & -a^{1} & -a^{2} & -a^{3} \\ a^{1} & a^{0} & -a^{3} & a^{2} \\ a^{2} & a^{3} & a^{0} & -a^{1} \\ a^{3} & -a^{2} & a^{1} & a^{0} \end{pmatrix} \tag{11.3.36}$$

Proof. The product of quaternions (11.3.35) and

$$x = x^{0} + x^{1}i + x^{2}j + x^{3}k$$

has form

$$\begin{aligned} ax = {} & a^{0}x^{0} - a^{1}x^{1} - a^{2}x^{2} - a^{3}x^{3} + (a^{0}x^{1} + a^{1}x^{0} + a^{2}x^{3} - a^{3}x^{2})i \\ & + (a^{0}x^{2} + a^{2}x^{0} + a^{3}x^{1} - a^{1}x^{3})j + (a^{0}x^{3} + a^{3}x^{0} + a^{1}x^{2} - a^{2}x^{1})k \end{aligned}$$

Therefore, function $f_a(x) = ax$ has Jacobian matrix (11.3.36). It is evident that $f_a \circ f_b = f_{ab}$. Similar equation is true for matrices

$$\begin{pmatrix} a^0 & -a^1 & -a^2 & -a^3 \\ a^1 & a^0 & -a^3 & a^2 \\ a^2 & a^3 & a^0 & -a^1 \\ a^3 & -a^2 & a^1 & a^0 \end{pmatrix} \begin{pmatrix} b^0 & -b^1 & -b^2 & -b^3 \\ b^1 & b^0 & -b^3 & b^2 \\ b^2 & b^3 & b^0 & -b^1 \\ b^3 & -b^2 & b^1 & b^0 \end{pmatrix}$$

$$= \begin{pmatrix} \begin{matrix} a^0b^0 - a^1b^1 \\ -a^2b^2 - a^3b^3 \end{matrix} & \begin{matrix} -a^0b^1 - a^1b^0 \\ -a^2b^3 + a^3b^2 \end{matrix} & \begin{matrix} -a^0b^2 + a^1b^3 \\ -a^2b^0 - a^3b^1 \end{matrix} & \begin{matrix} -a^0b^3 - a^1b^2 \\ +a^2b^1 - a^3b^0 \end{matrix} \\ \\ \begin{matrix} a^0b^1 + a^1b^0 \\ +a^2b^3 - a^3b^2 \end{matrix} & \begin{matrix} a^0b^0 - a^1b^1 \\ -a^2b^2 - a^3b^3 \end{matrix} & \begin{matrix} -a^0b^3 - a^1b^2 \\ +a^2b^1 - a^3b^0 \end{matrix} & \begin{matrix} a^0b^2 - a^1b^3 \\ +a^2b^0 + a^3b^1 \end{matrix} \\ \\ \begin{matrix} a^0b^2 - a^1b^3 \\ +a^2b^0 + a^3b^1 \end{matrix} & \begin{matrix} a^0b^3 + a^1b^2 \\ -a^2b^1 + a^3b^0 \end{matrix} & \begin{matrix} a^0b^0 - a^1b^1 \\ -a^2b^2 - a^3b^3 \end{matrix} & \begin{matrix} -a^0b^1 - a^1b^0 \\ -a^2b^3 - a^3b^2 \end{matrix} \\ \\ \begin{matrix} a^0b^3 + a^1b^2 \\ -a^2b^1 + a^3b^0 \end{matrix} & \begin{matrix} -a^0b^2 + a^1b^3 \\ -a^2b^0 - a^3b^1 \end{matrix} & \begin{matrix} a^0b^1 + a^1b^0 \\ +a^2b^3 - a^3b^2 \end{matrix} & \begin{matrix} a^0b^0 - a^1b^1 \\ -a^2b^2 - a^3b^3 \end{matrix} \end{pmatrix}$$

Therefore, we can identify the quaternion $a$ and the matrix $J_a$. □

CHAPTER 12

# Linear Map of $D$-Vector Spaces

## 12.1. Linear Map of $D$-Vector Spaces

Considering linear map of $D$-vector spaces we assume that division ring $D$ is finite dimensional algebra over field $F$.

Definition 12.1.1. Let field $F$ be subring of center $Z(D)$ of division ring $D$. Suppose $V$ and $W$ are $D$-vector spaces. We call map

$$A: V \rightarrow W$$

of $D$-vector space $V$ into $D$-vector space $W$ **linear map** if

$$A(x+y) = A(x) + A(y) \quad x, y \in V$$
$$A(px) = pA(x) \qquad p \in F$$

Let us denote $\mathcal{L}(D; V; W)$ set of linear maps

$$A: V \rightarrow W$$

of $D$-vector space $V$ into $D$-vector space $W$. □

It is evident that linear map of $D_*{}^*$-vector space as well linear map of ${}^*{}_*D$-vector space are linear maps. Set of morphisms of $D$-vector space is wider then set of morphisms of $D_*{}^*$-vector space. To consider linear map of vector space, we will follow method used in section 5.4.

Theorem 12.1.2. *Let $D$ be division ring of characteristic 0. Let $V$, $W$ be $D$-vector spaces. Linear map*

$$A: V \rightarrow W$$

*relative to basis $\overline{\overline{e}}_V$ of ${}_*{}^*D$-vector space $V$ and basis $\overline{\overline{e}}_W$ of ${}_*{}^*D$-vector space $W$ has form*

$$A(v) = e_{W \cdot j}\ A^j_i(v^i) \qquad v = e_{V*}{}^*v \tag{12.1.1}$$

*where $A^j_i(v^i)$ linearly depends on one variable $v^i$ and does not depends on the rest of coordinates of vector $v$.*

Proof. According to definition 12.1.1

$$A(v) = A(e_{V*}{}^*v) = A\left(\sum_i e_{V \cdot i}\ v^i\right) = \sum_i A(e_{V \cdot i}\ v^i) \tag{12.1.2}$$

For any given $i$ vector $A(e_{V \cdot i}\ v^i) \in W$ has only expansion

$$A(e_{V \cdot i}\ v^i) = e_{W \cdot j}\ A^j_i(v^i) \qquad A(e_{V \cdot i}\ v^i) = e_{W*}{}^*A_i(v^i) \tag{12.1.3}$$

relative to basis $\overline{\overline{e}}_W$ of ${}_*{}^*D$-vector space. Let us substitute (12.1.3) into (12.1.2). We get (12.1.1). □

DEFINITION 12.1.3. The linear map

$$A_i^j : D \to D$$

is called **partial linear map** of variable $v^i$. □

We can write linear map as product of matrices

$$A(v) = \begin{pmatrix} e_{W\cdot 1} & \dots & e_{W\cdot m} \end{pmatrix} {}_*{}^* \begin{pmatrix} A_i^1(v^i) \\ \dots \\ A_i^m(v^i) \end{pmatrix} \tag{12.1.4}$$

Let us define product of matrices

$$\begin{pmatrix} A_1^1 & \dots & A_n^1 \\ \dots & \dots & \dots \\ A_1^m & \dots & A_n^m \end{pmatrix} {}_\circ{}^\circ \begin{pmatrix} v^1 \\ \dots \\ v^n \end{pmatrix} = \begin{pmatrix} A_i^1(v^i) \\ \dots \\ A_i^m(v^i) \end{pmatrix} \tag{12.1.5}$$

where $A = \begin{pmatrix} A_i^j \end{pmatrix}$ is matrix of partial linear maps. Using the equation (12.1.5), we can write the equation (12.1.4) in the form

$$A(v) = \begin{pmatrix} e_{W\cdot 1} & \dots & e_{W\cdot m} \end{pmatrix} {}_*{}^* \begin{pmatrix} A_1^1 & \dots & A_n^1 \\ \dots & \dots & \dots \\ A_1^m & \dots & A_n^m \end{pmatrix} {}_\circ{}^\circ \begin{pmatrix} v^1 \\ \dots \\ v^n \end{pmatrix} \tag{12.1.6}$$

THEOREM 12.1.4. *Let $D$ be division ring of characteristic 0. Linear map*

$$A : v \in V \to w \in W \quad w = A(v) \tag{12.1.7}$$

*relative to basis $\overline{\overline{e}}_V$ of $D$-vector space $V$ and basis $\overline{\overline{e}}_W$ of $D$-vector space $W$*[12.1] *has form*

$$v = e_V {}_*{}^* v$$

$$w = e_W {}_*{}^* w \tag{12.1.8}$$

$$w^j = A_i^j(v^i) = A_{s\cdot 0\cdot}{}_i^j \; v^i \; A_{s\cdot 1\cdot}{}_i^j \tag{12.1.9}$$

PROOF. According to theorem 12.1.2 we can write linear map $A(v)$ as (12.1.1). Because for given indexes $i$, $j$ partial linear map $A_i^j(v^i)$ is linear with respect to variable $v^i$, then according to (10.1.14) it is possible to represent expression $A_i^j(v^i)$ as

$$A_i^j(v^i) = A_{s\cdot 0\cdot}{}_i^j \; v^i \; A_{s\cdot 1\cdot}{}_i^j \tag{12.1.10}$$

[12.1] Coordinate representation of map (12.1.7) depends on choice of basis. Equations change form if, for instance, we choose ${}^*{}_*D$-basis $\overline{\overline{r}}_*$ in $D$-vector space $W$.

where index $s$ is numbering items. Range of index $s$ depends on indexes $i$ and $j$. Combining equations (12.1.2) and (12.1.10), we get

$$
\begin{aligned}
A(v) &= e_{W\cdot j} A^{j}_{i}(v^{i}) &&= e_{W\cdot j} A_{s\cdot 0\cdot}{}^{j}_{i}\ v^{i}\ A_{s\cdot 1\cdot}{}^{j}_{i} \\
A(v) &= e_{W*}{}^{*} A_{i}(v^{i}) &&= e_{W*}{}^{*}(A_{s\cdot 0\cdot i}\ v^{i}\ A_{s\cdot 1\cdot i})
\end{aligned}
\tag{12.1.11}
$$

In equation (12.1.11), we also summarize on the index $i$. Equation (12.1.9) follows from comparison of equations (12.1.8) and (12.1.11). □

Definition 12.1.5. Expression $A_{s\cdot p\cdot}{}^{j}_{i}$ in equation (12.1.3) is called **component of linear map** $A$. □

Theorem 12.1.6. *Let $D$ be division ring of characteristic 0. Let $\overline{\overline{e}}_V$ be a basis of $D$-vector space $V$, $\overline{\overline{e}}_U$ be a basis of $D$-vector space $U$, and $\overline{\overline{e}}_W$ be a basis of $D$-vector space $W$. Suppose diagram of maps*

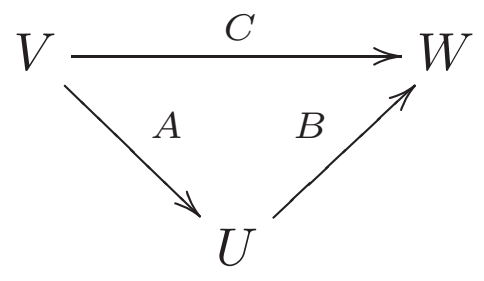

*is commutative diagram where linear map $A$ has presentation*

$$
u = A(v) = e_{U\cdot j} A^{j}_{i}(v^{i}) = e_{U\cdot j} A_{s\cdot 0\cdot}{}^{j}_{i}\ v^{i}\ A_{s\cdot 1\cdot}{}^{j}_{i}
\tag{12.1.12}
$$

*relative to selected bases and linear map $B$ has presentation*

$$
w = B(u) = e_{W\cdot k} B^{k}_{j}(u^{j}) = e_{W\cdot k} B_{t\cdot 0\cdot}{}^{k}_{j}\ u^{j}\ B_{t\cdot 1\cdot}{}^{k}_{j}
\tag{12.1.13}
$$

*relative to selected bases. Then map $C$ is linear map and has presentation*

$$
w = C(v) = e_{W\cdot k} C^{k}_{i}(v^{i}) = e_{W\cdot k} C_{u\cdot 0\cdot}{}^{k}_{i}\ v^{i}\ C_{u\cdot 1\cdot}{}^{k}_{i}
\tag{12.1.14}
$$

*relative to selected bases, where*[12.2]

$$
\begin{aligned}
& C^{k}_{i}(v^{i}) && = B^{k}_{j}(A^{j}_{i}(v^{i})) \\
C_{u\cdot 0\cdot}{}^{k}_{i} &= C_{ts\cdot 0\cdot}{}^{k}_{i} && = B_{t\cdot 0\cdot}{}^{k}_{j}\ A_{s\cdot 0\cdot}{}^{j}_{i} \\
C_{u\cdot 1\cdot}{}^{k}_{i} &= C_{ts\cdot 1\cdot}{}^{k}_{i} && = A_{s\cdot 1\cdot}{}^{j}_{i}\ B_{t\cdot 1\cdot}{}^{k}_{j}
\end{aligned}
\tag{12.1.15}
$$

Proof. The map $C$ is linear map because

$$
\begin{aligned}
C(a+b) &= B(A(a+b)) \\
&= B(A(a)+A(b)) \\
&= B(A(a)) + B(A(b)) \\
&= C(a) + C(b) \\
C(ab) &= B(A(ab)) = B(aA(b)) \\
&= aB(A(b)) = aC(b) \\
& a, b \in V \\
& a \in F
\end{aligned}
$$

Equation (12.1.14) follows from substituting (12.1.12) into (12.1.13). □

[12.2] Index $u$ appeared composite index, $u = st$. However it is possible that some items in (12.1.15) may be summed together.

Theorem 12.1.7. *For linear map $A$ there exists linear map $B$ such, that*

$$A(axb)=B(x)$$
$$B_{s\cdot 0\cdot}{}^{j}_{i}=A_{s\cdot 0\cdot}{}^{j}_{i}\ a$$
$$B_{s\cdot 1\cdot}{}^{j}_{i}=b\ A_{s\cdot 1\cdot}{}^{j}_{i}$$

Proof. The map $B$ is linear map because

$$B(x+y)=A(a(x+y)b)=A(axb+ayb)=A(axb)+A(ayb)=B(x)+B(y)$$
$$B(cx)=A(acxb)=A(caxb)=cA(axb)=cB(x)$$
$$x,y\in V,c\in F$$

According to equation (12.1.9)

$$B_{s\cdot 0\cdot}{}^{j}_{i}\ v^{i}\ B_{s\cdot 1\cdot}{}^{j}_{i}=A_{s\cdot 0\cdot}{}^{j}_{i}\ (a\ v^{i}\ b)A_{s\cdot 1\cdot}{}^{j}_{i}=(A_{s\cdot 0\cdot}{}^{j}_{i}\ a)v^{i}\ (b\ A_{s\cdot 1\cdot}{}^{j}_{i})$$

□

Theorem 12.1.8. *Let $D$ be division ring of characteristic 0. Let*

$$A:V\rightarrow W$$

*linear map of $D$-vector space $V$ into $D$-vector space $W$. Then* $A(0)=0$.

Proof. Corollary of equation

$$A(a+0)=A(a)+A(0)$$

□

Definition 12.1.9. The set

$$\ker f=\{x\in V:f(x)=0\}$$

is called **kernel of linear map**

$$A:V\rightarrow W$$

of $D$-vector space $V$ into $D$-vector space $W$. □

Definition 12.1.10. The linear map

$$A:V\rightarrow W$$

of $D$-vector space $V$ into $D$-vector space $W$ is called **singular**, if

$$\ker A=V$$

□

## 12.2. Polylinear Map of $D$-Vector Spaces

Definition 12.2.1. Let field $F$ be subring of center $Z(D)$ of division ring $D$ of characteristic 0. Suppose $V_1$, ..., $V_n$, $W_1$, ..., $W_m$ are $D$-vector spaces. We call map

$$A:V_1\times...\times V_n\rightarrow W_1\times...\times W_m$$
$$w_1\times...\times w_m=A(v_1,...,v_n) \tag{12.2.1}$$

**polylinear map** of $\times$-$D$-vector space $V_1\times ...\times V_n$ into $\times$-$D$-vector space $W_1\times ...\times W_m$, if

$$A(p_1,...,p_i+q_i,...,p_n)=A(p_1,...,p_i,...,p_n)+A(p_1,...,q_i,...,p_n)$$
$$A(p_1,...,ap_i,...,p_n)=aA(p_1,...,p_i,...,p_n)$$
$$1\le i\le n\quad p_i,q_i\in V_i\quad a\in F$$

□

Definition 12.2.2. Let us denote $\mathcal{L}(D;V_1,...,V_n;W_1,...,W_m)$ set of polylinear maps of $\times$-$D$-vector space $V_1\times ...\times V_n$ into $\times$-$D$-vector space $W_1\times ...\times W_m$. □

Theorem 12.2.3. *Let $D$ be division ring of characteristic $0$. For each $k\in K=[1,n]$ let $\overline{\overline{e}}_{V_k}$ be basis in $D$-vector space $V_k$ and*

$$v_k={v_k}^*{}_*e_{V_k}\quad v_k\in V_k$$

*For each $l$, $1\le l\le m$, let $\overline{\overline{e}}_{W_l}$ be basis in $D$-vector space $W_l$ and*

$$w_l={w_l}^*{}_*e_{W_l}\quad w_l\in W_l$$

*Polylinear map (12.2.1) relative to bases $\overline{\overline{e}}_{V_1}$, ..., $\overline{\overline{e}}_{V_n}$, $\overline{\overline{e}}_{W_1}$, ..., $\overline{\overline{e}}_{W_m}$ has form*

$$\begin{aligned}w_l^{j}&=A_{l\cdot}{}^{j}_{i_1...i_n}(v_1^{i_1},...,v_n^{i_n})\\&=A^n_{s\cdot 0\cdot l\cdot}{}^{j}_{i_1...i_n}\ \sigma_s(v_1^{i_1})\ A^n_{s\cdot 1.l\cdot}{}^{j}_{i_1...i_n}\ ...\ \sigma_s(v_n^{i_n})\ A^n_{s\cdot n\cdot l\cdot}{}^{j}_{i_1...i_n}\end{aligned}\tag{12.2.2}$$

*Range $S$ of index $s$ depends on values of indexes $i_1$, ..., $i_n$. $\sigma_s$ is a transposition of set of variables $\{v_1^{i_1},...,v_n^{i_n}\}$.*

Proof. Since we may consider map $A$ into $\times$-$D$-vector space $W_1\times ...\times W_m$ componentwise, then we may confine to considering of map

$$A_l:V_1\times ...\times V_n\to W_l\quad w_l=A_l(v_1,...,v_n)\tag{12.2.3}$$

We prove statement by induction on $n$.

When $n=1$ the statement of theorem is statement of theorem 12.1.4. In such case we may identify[12.3]

$$A^1_{s\cdot p\cdot}{}^{j}_{i}=A_{s\cdot p\cdot}{}^{j}_{i}\quad p=0,1$$

[12.3] In representation (12.2.2) we will use following rules.

- If range of any index is set consisting of one element, then we will omit corresponding index.
- If $n=1$, then $\sigma_s$ is identical transformation. We will not show such transformation in the expression.

Let statement of theorem be true for $n = k-1$. Then it is possible to represent map (12.2.3) as

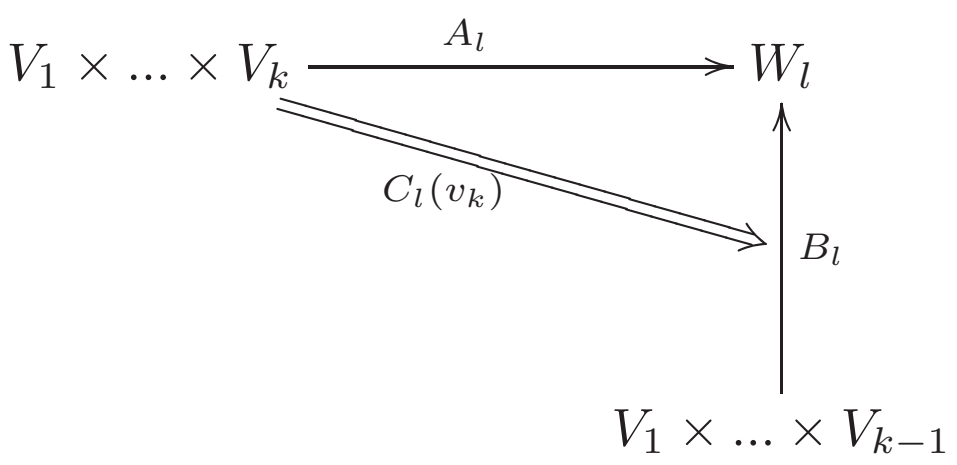

$$w_l = A_l(v_1, ..., v_k) = C_l(v_k)(v_1, ..., v_{k-1})$$

According to statement of induction polylinear map $B_l$ has form

$$w_l^{j} = B^{k-1}{}_{t\cdot 0\cdot l\cdot}{}^{j}_{i_1...i_{k-1}}\ \sigma_t(v_1^{i_1})\ B^{k-1}{}_{t\cdot 1.l\cdot}{}^{j}_{i_1...i_{k-1}}\ ...\ \sigma_t(v_{k-1}^{i_{k-1}})\ B^{k-1}{}_{t\cdot k-1.l\cdot}{}^{j}_{i_1...i_{k-1}}$$

According to construction $B_l = C_l(v_k)$. Therefore, expressions $B^{k-1}{}_{t\cdot p\cdot l\cdot}{}^{j}_{i_1...i_{k-1}}$ are functions of $v_k$. Since $C_l(v_k)$ is linear map of $v_k$, then only one expression $B^{k-1}{}_{t\cdot p\cdot l\cdot}{}^{j}_{i_1...i_{k-1}}$ is linear map $v_k$, and rest expressions $B^{k-1}{}_{t\cdot q\cdot l\cdot}{}^{j}_{i_1...i_{k-1}}$ do not depend on $v_k$.

Without loss of generality, assume $p = 0$. According to theorem 12.1.4

$$B^{k-1}{}_{t\cdot 0\cdot l\cdot}{}^{j}_{i_1...i_{k-1}} = C^{k}{}_{tr\cdot 0\cdot l\cdot}{}^{j}_{i_k i_1...i_{k-1}}\ v_k^{i_k}\ C^{k}{}_{tr\cdot 1.l\cdot}{}^{j}_{i_k i_1...i_{k-1}}$$

Assume $s = tr$. Let us define transposition $\sigma_s$ according to rule

$$\sigma_s = \sigma(tr) = \begin{pmatrix} v_k^{i_k} & v_1^{i_1} & ... & v_{k-1}^{i_{k-1}} \\ v_k^{i_k} & \sigma_t(v_1^{i_1}) & ... & \sigma_t(v_{k-1}^{i_{k-1}}) \end{pmatrix}$$

Suppose

$$A^{k}{}_{tr\cdot q+1.l\cdot}{}^{j}_{i_k i_1...i_{k-1}} = B^{k-1}{}_{t\cdot q\cdot l\cdot}{}^{j}_{i_1...i_{k-1}} \qquad q = 1, ..., k-1$$
$$A^{k}{}_{tr\cdot q\cdot l\cdot}{}^{j}_{i_k i_1...i_{k-1}} = C^{k}{}_{tr\cdot q\cdot l\cdot}{}^{j}_{i_k i_1...i_{k-1}} \qquad q = 0, 1$$

We proved step of induction. □

Definition 12.2.4. Expression $A^{n}{}_{s\cdot p\cdot l\cdot}{}^{j}_{i_1...i_n}$ in equation (12.2.2) is called **component of polylinear map** $A$. □

CHAPTER 13

# Tensor Product

## 13.1. Tensor Product of Rings

THEOREM 13.1.1. *Let $G$ be semigroup. Let $M$ be free module with basis $G$ over commutative ring $F$. Then the structure of ring is defined on $M$ such that $F$ is subring of center $Z(M)$ of ring $M$.*

PROOF. Arbitrary vectors $a$, $b$, $c \in M$ have unique expansion

$$a = a^i g_i \qquad b = b^j g_j \qquad c = c^k g_k$$

relative to basis $G$. $i \in I$, $j \in J$, $k \in K$ where $I$, $J$, $K$ are finite sets. Without loss of generality, we assume $i$, $j$, $k \in I \cup J \cup K$.

Let us define product $a$ and $b$ by equation

$$ab = (a^i b^j)(g_i g_j) \tag{13.1.1}$$

From chain of equations

$$\begin{aligned}(a+b)c &= ((a^i + b^i)g_i)(c^j g_j)\\ &= ((a^i + b^i)c^j)(g_i g_j)\\ &= (a^i c^j + b^i c^j)(g_i g_j)\\ &= (a^i c^j)(g_i g_j) + (b^i c^j)(g_i g_j)\\ &= (a^i g_i)(c^j g_j) + (b^i g_i)(c^j g_j)\\ &= ac + bc\end{aligned}$$

it follows that this product is right-distributive relative to sum. The same way we prove that this product is left-distributive relative to sum. From chain of equations

$$\begin{aligned}(ab)c &= ((a^i g_i)(b^j g_j))(c^k g_k)\\ &= ((a^i b^j)(g_i g_j))(c^k g_k)\\ &= ((a^i b^j)c^k)((g_i g_j)g_k)\\ &= (a^i (b^j c^k))(g_i (g_j g_k))\\ &= (a^i g_i)((b^j c^k)(g_j g_k))\\ &= (a^i g_i)((b^j g_j)(c^k g_k))\\ &= a(bc)\end{aligned}$$

it follows that this product is associative.

Let $e \in G$ be unit of semigroup $G$. From equation

$$ae = (a^i g_i)e = a^i (g_i e) = a^i g_i = a$$

it follows that $e$ is unit of ring $M$.

Map

$$f : p \in F \to pe \in M$$

is homomorphic injection of ring $F$ into ring $M$. Therefore, center $Z(M)$ of ring $M$ contains ring $F$. □

Let $R_1$, ..., $R_n$ be rings of characteristic 0.[13.1] Let $F$ be maximum commutative ring which for any $i$, $i = 1,...,n$, is subring of center $Z(R_i)$. Consider category $\mathcal{A}$ whose objects are polylinear over commutative ring $F$ maps

$$f : R_1 \times ... \times R_n \longrightarrow S_1 \qquad g : R_1 \times ... \times R_n \longrightarrow S_2$$

where $S_1$, $S_2$ are modules over ring $F$. We define morphism $f \to g$ to be linear over commutative ring $F$ map $h : S_1 \to S_2$ making commutative following diagram

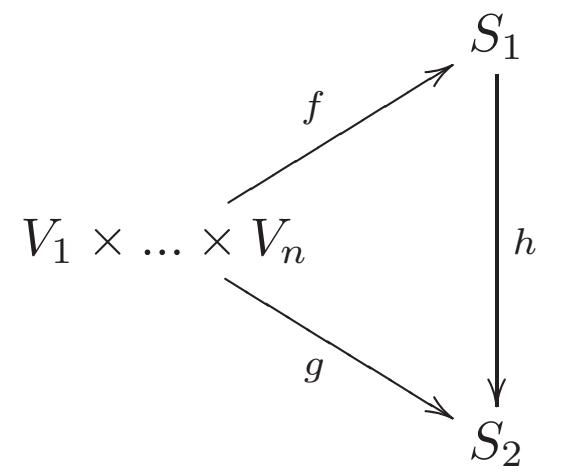

Universal object $R_1 \otimes ... \otimes R_n$ of category $\mathcal{A}$ is called **tensor product of rings $R_1$, ..., $R_n$ over commutative ring $F$**.

Theorem 13.1.2. *There exists tensor product of rings.*

Proof. Let $M$ be module over ring $F$ generated by product $R_1 \times ... \times R_n$ of multiplicative semigroups of rings $R_1$, ..., $R_n$. Injection

$$i : R_1 \times ... \times R_n \longrightarrow M$$

is defined according to rule

$$i(d_1, ..., d_n) = (d_1, ..., d_n) \tag{13.1.2}$$

and is homomorphism of multiplicative semigroup of ring $R_1 \times ... \times R_n$ onto basis of module $M$. Therefore, product of vectors of basis is defined componentwise

$$(d_1, ..., d_n)(c_1, ..., c_n) = (d_1c_1, ..., d_nc_n) \tag{13.1.3}$$

Let $N \subset M$ be submodule generated by elements of the following type

$$\begin{aligned} &(d_1, ..., d_i + c_i, ..., d_n) - (d_1, ..., d_i, ..., d_n) - (d_1, ..., c_i, ..., d_n) \\ &(d_1, ..., ad_i, ..., d_n) - a(d_1, ..., d_i, ..., d_n) \end{aligned} \tag{13.1.4}$$

where $d_i, c_i \in R_i$, $a \in F$. Let

$$j : M \to M/N$$

be canonical map on factor module. Consider commutative diagram

(13.1.5)

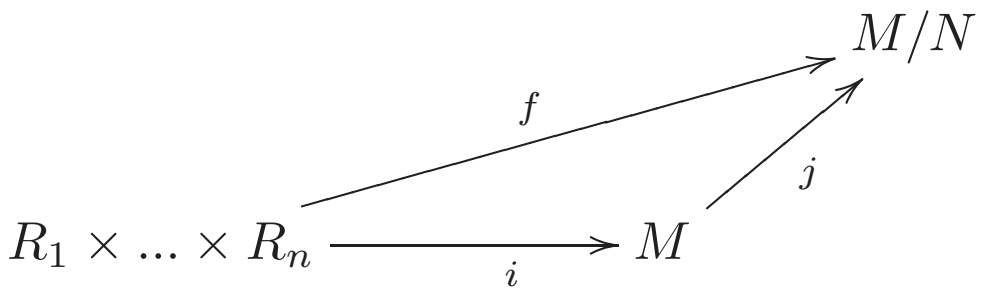

[13.1] I give definition of tensor product of rings following to definition in [1], p. 601 - 603.

Since elements (13.1.4) belong to kernel of linear map $j$, then from equation (13.1.2) it follows

(13.1.6)
$$\begin{aligned} f(d_1, ..., d_i + c_i, ..., d_n) =& f(d_1, ..., d_i, ..., d_n) + f(d_1, ..., c_i, ..., d_n) \\ f(d_1, ..., ad_i, ..., d_n) =& af(d_1, ..., d_i, ..., d_n) \end{aligned}$$

From equations (13.1.6), it follows that map $f$ is polylinear over field $F$. Since $M$ is module with basis $R_1 \times ... \times R_n$, then according to theorem [1]-4.1, on p. 135 for any module $V$ and any polylinear over $F$ map

$$g : R_1 \times ... \times R_n \longrightarrow V$$

there exists a unique homomorphism $k : M \to V$, for which following diagram is commutative

(13.1.7)

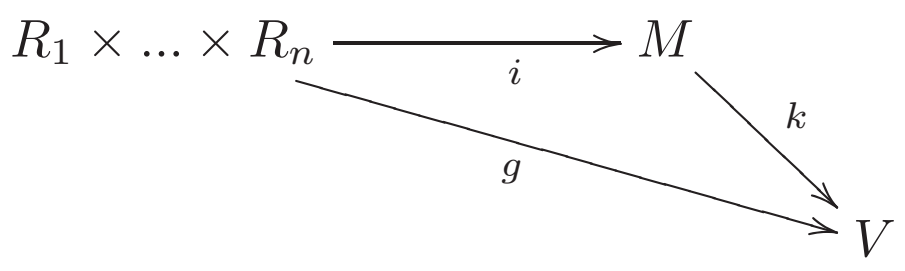

Since $g$ is polylinear over $F$, then $\ker k \subseteq N$. According to statement on p. [1]-119 map $j$ is universal in the category of homomorphisms of vector space $M$ whose kernel contains $N$. Therefore, we have homomorphism

$$h : M/N \to V$$

which makes the following diagram commutative

(13.1.8)

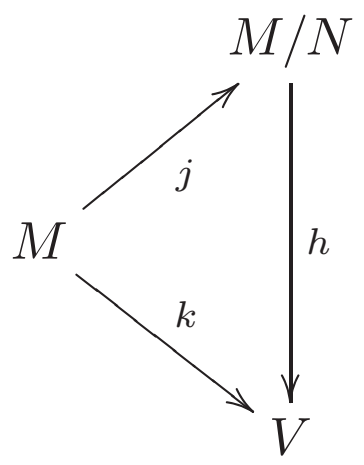

We join diagrams (13.1.5), (13.1.7), (13.1.8), and get commutative diagram

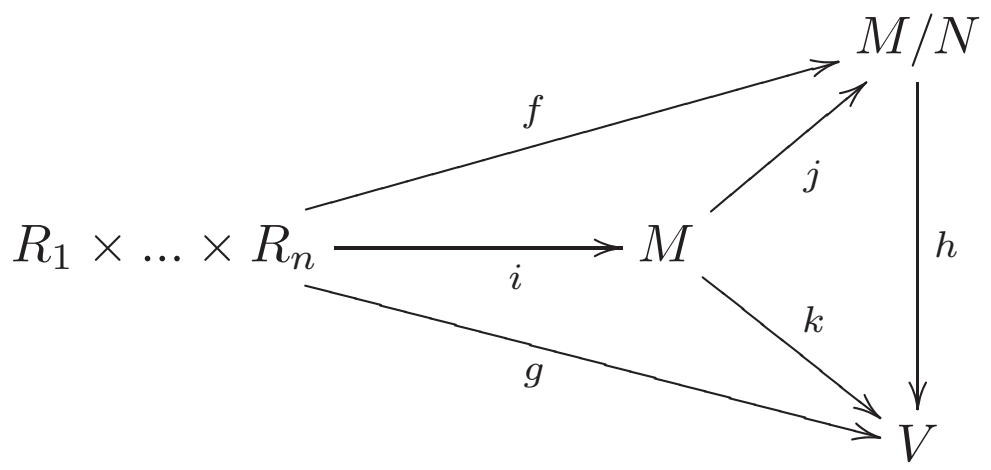

Since $\mathrm{Im} f$ generates $M/N$, then map $h$ is uniquely determined. □

According to proof of theorem 13.1.2

$$R_1 \otimes ... \otimes R_n = M/N$$

If $d_i \in R_i$, we write

(13.1.9)
$$d_1 \otimes ... \otimes d_n = j(d_1, ..., d_n)$$

From equations (13.1.2) and (13.1.9) it follows

$$f(d_1, ..., d_n) = d_1 \otimes ... \otimes d_n$$

We can write equations (13.1.6) as

$$
\begin{aligned}
(13.1.10)\qquad & d_1 \otimes ... \otimes (d_i + c_i) \otimes ... \otimes d_n \\
&= d_1 \otimes ... \otimes d_i \otimes ... \otimes d_n + d_1 \otimes ... \otimes c_i \otimes ... \otimes d_n \\
& d_1 \otimes ... \otimes (ad_i) \otimes ... \otimes d_n = a(d_1 \otimes ... \otimes d_i \otimes ... \otimes d_n)
\end{aligned}
$$

THEOREM 13.1.3. *Tensor product $R_1 \otimes ... \otimes R_n$ of rings $R_1$, ..., $R_n$ of characteristic 0 over ring $F$ is ring.*

PROOF. According to proof of theorem 13.1.2 basis of module $M$ over ring $F$ is semigroup. According to theorem 13.1.1 structure of ring is defined in module $M$. Submodule $N$ is two-sided ideal of ring $M$. Canonical map on factor module

$$j : M \to M/N$$

is also canonical map on factor ring. Therefore, there is the structure of ring in module $M/N = R_1 \otimes ... \otimes R_n$. From equation (13.1.3) it follows that

$$(d_1 \otimes ... \otimes d_n)(c_1 \otimes ... \otimes c_n) = (d_1c_1) \otimes ... \otimes (d_nc_n)$$

From equation (13.1.10), it follows that

$$
\begin{aligned}
& d_1 \otimes ... \otimes d_i \otimes ... \otimes d_n \\
&= d_1 \otimes ... \otimes (d_i + 0) \otimes ... \otimes d_n \\
&= d_1 \otimes ... \otimes d_i \otimes ... \otimes d_n + d_1 \otimes ... \otimes 0 \otimes ... \otimes d_n
\end{aligned}
$$

Therefore, we can identify tensor $d_1 \otimes ... \otimes 0 \otimes ... \otimes d_n$ with zero $0 \otimes ... \otimes 0$. According to statement of theorem product equal $0 \otimes ... \otimes 0$ if only one of factors equal $0 \otimes ... \otimes 0$. Tensor $e_1 \otimes ... \otimes e_n$ is unit of product. □

## 13.2. Tensor Product of Division Rings

Let $D_1$, ..., $D_n$ be division rings of characteristic 0. Let $F$ be field which for any $i$, $i = 1,...,n$ is subring of center $Z(D_i)$. Tensor

$$(a_1 \otimes ... \otimes a_n)^{-1} = (a_1)^{-1} \otimes ... \otimes (a_n)^{-1}$$

is inverse tensor to tensor

$$a_1 \otimes ... \otimes a_n \in D_1 \otimes ... \otimes D_n$$

However **tensor product $D_1 \otimes ... \otimes D_n$ of division rings $D_1$, ..., $D_n$**, in general, is not division ring, because we cannot answer on following question. Does element

$$p(a_1 \otimes ... \otimes a_n) + q(b_1 \otimes ... \otimes b_n)$$

have inverse one?

REMARK 13.2.1. Representation of tensor as

$$a^s\ d_{s\cdot 1} \otimes ... \otimes d_{s\cdot n} \tag{13.2.1}$$

is ambiguous. Since $r \in F$, then

$$(d_1r) \otimes d_2 = d_1 \otimes (rd_2)$$

We can increase or decrease number of summands using algebraic operations. Following transformation

$$\begin{aligned}&d_1\otimes d_2+c_1\otimes c_2\\=&d_1\otimes(d_2-c_2+c_2)+c_1\otimes c_2\\=&d_1\otimes(d_2-c_2)+d_1\otimes c_2+c_1\otimes c_2\\=&d_1\otimes(d_2-c_2)+(d_1+c_1)\otimes c_2\end{aligned}$$

is between possible transformations. □

Consider following theorem for the purposes of illustration of remark 13.2.1.

THEOREM 13.2.2. *Since in division ring $D$ of characteristic 0 there exist elements $a$, $b$, $c$ product of which is not commutative, then there exists nontrivial representation of zero tensor.*

PROOF. Proof of theorem follows from chain of equations

$$\begin{aligned}0\otimes 0=&a\otimes a-a\otimes a+b\otimes b-b\otimes b\\=&(a+c)\otimes a-c\otimes a-a\otimes(a-c)-a\otimes c\\&+(b-c)\otimes b+c\otimes b-b\otimes(b+c)+b\otimes c\\=&(a+c)\otimes a-a\otimes(a-c)+(b-c)\otimes b\\&+c\otimes(b-a)-b\otimes(b+c)+(b-a)\otimes c\end{aligned}$$

□

From ambiguity of representation of tensor follows that we must to find canonical representation of tensor. We can find answer of this problem in case of tensor product of division rings. Let division ring $D_i$ be vector space over field $F\subset Z(D_i)$. Let this vector space have finite basis $e_{i\cdot\boldsymbol{s}_i}$, $\boldsymbol{s}_i\in\boldsymbol{S}_i$, $|\boldsymbol{S}_i|=m_i$. Therefore, any element $f_i$ of division ring $D_i$ has expansion

$$f_i=e_{i\cdot\boldsymbol{s}_i}f_i^{\boldsymbol{s}_i}$$

In this case we can write tensor (13.2.1) in form

$$a^s({}_{s\cdot}f_1^{\boldsymbol{s}_1}e_{1\cdot\boldsymbol{s}_1})\otimes...\otimes({}_{s\cdot}f_n^{\boldsymbol{s}_n}e_{n\cdot\boldsymbol{s}_n})\tag{13.2.2}$$

where $a^s$, ${}_sf_1^{s_1}\in F$. We can reduce expression (13.2.2)

$$a^s\,{}_{s\cdot}f_1^{\boldsymbol{s}_1}...{}_{s\cdot}f_n^{\boldsymbol{s}_n}e_{1\cdot\boldsymbol{s}_1}\otimes...\otimes e_{n\cdot\boldsymbol{s}_n}\tag{13.2.3}$$

Let

$$a^s\,{}_{s\cdot}f_1^{\boldsymbol{s}_1}...{}_{s\cdot}f_n^{\boldsymbol{s}_n}=f^{s_1...s_n}$$

Then equation (13.2.3) has form

$$f^{s_1...s_n}e_{1\cdot\boldsymbol{s}_1}\otimes...\otimes e_{n\cdot\boldsymbol{s}_n}\tag{13.2.4}$$

Expression (13.2.4) is defined uniquely up to selected basis. Expression $f^{s_1...s_n}$ is called **standard component of tensor**.

## 13.3. Tensor Product of $D_*{}^*$-Vector Spaces

Let $D_1$, ..., $D_n$ be division rings of characteristic 0.[13.2] Let $V_i$ be $D_{i*}{}^*$-vector space, $i = 1$, ..., $n$. Consider category $\mathcal{A}$ whose objects are polylinear maps

$$f : V_1 \times ... \times V_n \longrightarrow W_1 \qquad g : V_1 \times ... \times V_n \longrightarrow W_2$$

where $W_1$, $W_2$ are $D\star$-modules over ring $D_1 \otimes ... \otimes D_n$. We define morphism $f \to g$ to be additive map

$$h : W_1 \to W_2$$

making commutative following diagram

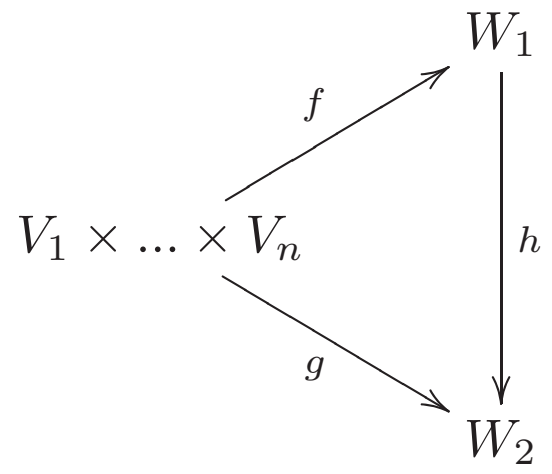

Universal object $V_1 \otimes ... \otimes V_n$ of category $\mathcal{A}$ is called **tensor product of $D\star$-vector spaces** $V_1$, ..., $V_n$.

Theorem 13.3.1. *There exists tensor product of $D\star$-vector spaces.*

Proof. Let $F$ be field which for any $i$, $i = 1,...,n$, is subring of center $Z(D_i)$.

Let $D$ be free vector space over ring $F$ generated by product $D_1 \times ... \times D_n$ of multiplicative semigroups of division rings $D_1$, ..., $D_n$. Injection

$$i' : D_1 \times ... \times D_n \longrightarrow D$$

is defined according to rule

$$i'(d_1, ..., d_n) = (d_1, ..., d_n) \tag{13.3.1}$$

and is homomorphism of multiplicative semigroup of ring $D_1 \times ... \times D_n$ onto basis of module $D$. Therefore, product of vectors of basis is defined componentwise

$$(d_1, ..., d_n)(c_1, ..., c_n) = (d_1c_1, ..., d_nc_n) \tag{13.3.2}$$

According to theorem 13.1.1 the structure of ring is defined on vector space $D$.

Consider direct product $F \times D$ of field $F$ and ring $D$. We will identify element $(f, e)$ with element $f \in F$ and element $(1, d)$ with element $d \in D$.

Let $M$ be free module over ring $F \times D$ generated by Cartesian product $V_1 \times ... \times V_n$ . Since $v_1 \in V_1$, ..., $v_n \in V_n$, then we denote corresponding vector from $M$ as $(v_1, ..., v_n)$. Let

$$i : V_1 \times ... \times V_n \longrightarrow M$$

be injection defined according to rule

$$i(v_1, ..., v_n) = (v_1, ..., v_n) \tag{13.3.3}$$

Let $N \subset M$ be vector subspacee generated by elements of the following type

$$\begin{array}{l} (v_1, ..., v_i + w_i, ..., v_n) - (v_1, ..., v_i, ..., v_n) - (v_1, ..., w_i, ..., v_n) \\ (v_1, ..., av_i, ..., v_n) - a(v_1, ..., v_i, ..., v_n) \end{array} \tag{13.3.4}$$

[13.2] I give definition of tensor product of $D_*{}^*$-vector spaces following to definition in [1], p. 601 - 603.

where $v_i,\ w_i \in V_i,\ a \in F.$ Let

$$j : M \to M/N$$

be canonical map on factor module. Consider commutative diagram

(13.3.5)

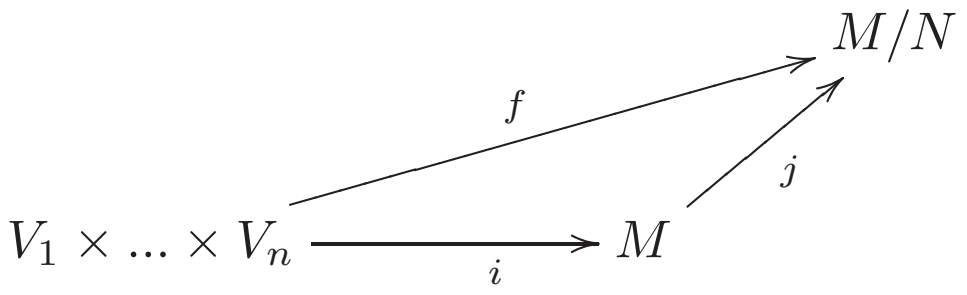

Since elements (13.3.4) belong to kernel of linear map $j$, then from equation (13.3.3) it follows

$$\begin{aligned} f(v_1, ..., v_i + w_i, ..., v_n) &= f(v_1, ..., v_i, ..., v_n) + f(v_1, ..., w_i, ..., v_n) \\ f(v_1, ..., av_i, ..., v_n) &= af(v_1, ..., v_i, ..., v_n) \end{aligned} \tag{13.3.6}$$

From equations (13.3.6), it follows that map $f$ is polylinear over field $F$. Since $M$ is module with basis $V_1 \times ... \times V_n$, then according to theorem [1]-4.1, on p. 135 for any module $V$ and any polylinear over $F$ map

$$g : V_1 \times ... \times V_n \longrightarrow V$$

there exists a unique homomorphism $k : M \to V$, for which following diagram is commutative

(13.3.7)

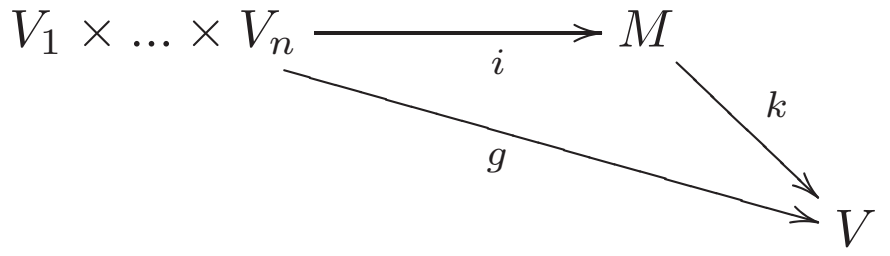

Since $g$ is polylinear over $F$, then $\ker k \subseteq N$. According to statement on p. [1]-119 map $j$ is universal in the category of homomorphisms of module $M$ whose kernel contains $N$. Therefore, we have homomorphism

$$h : M/N \to V$$

which makes the following diagram commutative

(13.3.8)

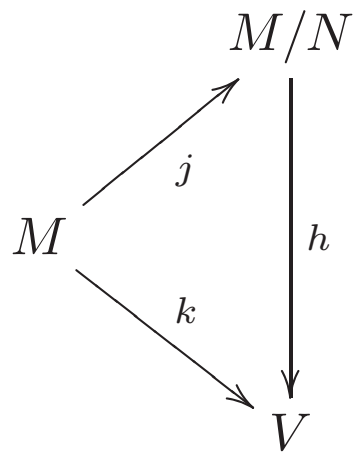

We join diagrams (13.3.5), (13.3.7), (13.3.8), and get commutative diagram

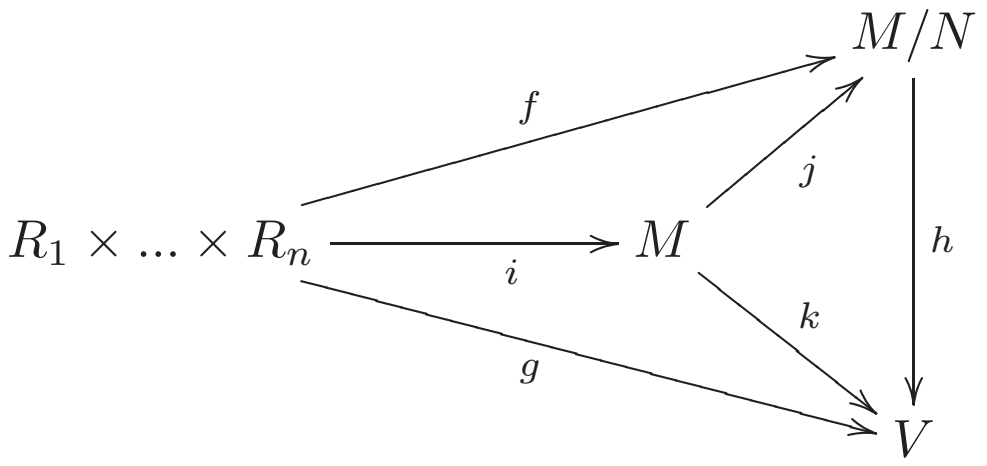

Since Im$f$ generates $M/N$, then map $h$ is uniquely determined. □

According to proof of theorem 13.3.1

$$V_1 \otimes ... \otimes V_n = M/N$$

If $v_i \in V_i$, we write

(13.3.9) $$v_1 \otimes ... \otimes v_n = j(v_1, ..., v_n)$$

From equations (13.3.3) and (13.3.9) it follows

$$f(v_1, ..., v_n) = v_1 \otimes ... \otimes v_n$$

We can write equations (13.3.6), as

$$\begin{aligned}
&v_1 \otimes ... \otimes (v_i + w_i) \otimes ... \otimes v_n\\
=&v_1 \otimes ... \otimes v_i \otimes ... \otimes v_n + v_1 \otimes ... \otimes w_i \otimes ... \otimes v_n\\
&v_1 \otimes ... \otimes (av_i) \otimes ... \otimes v_n = a(v_1 \otimes ... \otimes v_i \otimes ... \otimes v_n)
\end{aligned}$$

THEOREM 13.3.2. *Tensor product $V_1 \otimes ... \otimes V_n$ is module over tensor product $D_1 \otimes ... \otimes D_n$.*

PROOF. To prove statement of theorem we prove that representation

$$(d_1 \otimes ... \otimes d_n)(v_1 \otimes ... \otimes v_n) = (d_1, ..., d_n)(v_1 \otimes ... \otimes v_n)$$

of ring $D_1 \otimes ... \otimes D_n$ in module $V_1 \otimes ... \otimes V_n$ is defined properly. This follows from chain of equations

$$\begin{aligned}
(d_1 \otimes ... \otimes (ad_i) \otimes ... \otimes d_n)(v_1 \otimes ... \otimes v_n) &= (d_1, ..., ad_i, ..., d_n)(v_1 \otimes ... \otimes v_n)\\
&= (d_1 v_1) \otimes ... \otimes (ad_i v_i) \otimes ... \otimes (d_n v_n)\\
&= a(d_1 v_1) \otimes ... \otimes (d_n v_n)\\
&= a((d_1, ..., d_n)(v_1 \otimes ... \otimes v_n))\\
&= a((d_1 \otimes ... \otimes d_n)(v_1 \otimes ... \otimes v_n))\\
&= (a(d_1 \otimes ... \otimes d_n))(v_1 \otimes ... \otimes v_n)
\end{aligned}$$

$$\begin{aligned}
&(d_1 \otimes ... \otimes (d_i + c_i) \otimes ... \otimes d_n)(v_1 \otimes ... \otimes v_n)\\
=&(d_1, ..., d_i + c_i, ..., d_n)(v_1 \otimes ... \otimes v_n)\\
=&(d_1 v_1) \otimes ... \otimes ((d_i + c_i)v_i) \otimes ... \otimes (d_n v_n)\\
=&(d_1 v_1) \otimes ... \otimes (d_i v_i + c_i v_i) \otimes ... \otimes (d_n v_n)\\
=&(d_1 v_1) \otimes ... \otimes (d_i v_i) \otimes ... \otimes (d_n v_n) + (d_1 v_1) \otimes ... \otimes (c_i v_i) \otimes ... \otimes (d_n v_n)\\
=&(d_1 v_1) \otimes ... \otimes (d_i v_i) \otimes ... \otimes (d_n v_n) + (d_1 v_1) \otimes ... \otimes (c_i v_i) \otimes ... \otimes (d_n v_n)\\
=&(d_1, ..., d_i, ..., d_n)(v_1 \otimes ... \otimes v_n) + (d_1, ..., c_i, ..., d_n)(v_1 \otimes ... \otimes v_n)\\
=&(d_1 \otimes ... \otimes d_i \otimes ... \otimes d_n)(v_1 \otimes ... \otimes v_n) + (d_1 \otimes ... \otimes c_i \otimes ... \otimes d_n)(v_1 \otimes ... \otimes v_n)\\
=&((d_1 \otimes ... \otimes d_i \otimes ... \otimes d_n) + (d_1 \otimes ... \otimes c_i \otimes ... \otimes d_n))(v_1 \otimes ... \otimes v_n)
\end{aligned}$$

□

## 13.4. Tensor Product of $D$-Vector Spaces

Let $D_1$, ..., $D_n$ be division rings of characteristic 0.[13.3] Let $V_i$ be $D_i$-vector space, $i = 1$, ..., $n$. Consider category $\mathcal{A}$ whose objects are polylinear maps

$$f : V_1 \times ... \times V_n \longrightarrow W_1 \qquad g : V_1 \times ... \times V_n \longrightarrow W_2$$

where $W_1$, $W_2$ are modules over ring $D_1 \otimes ... \otimes D_n$. We define morphism $f \to g$ to be linear map

$$h : W_1 \to W_2$$

making commutative following diagram

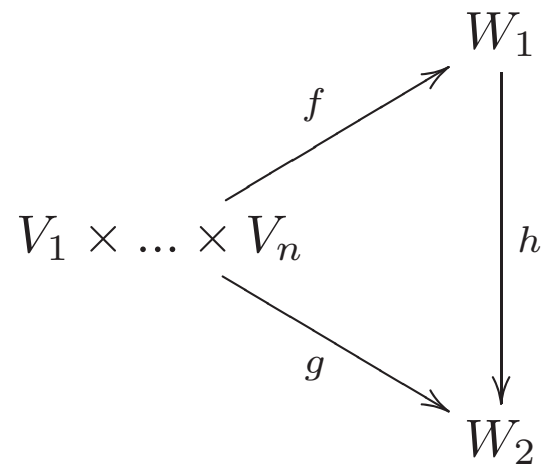

Universal object $V_1 \otimes ... \otimes V_n$ of category $\mathcal{A}$ is called **tensor product of $D$-vector spaces** $V_1$, ..., $V_n$.

Theorem 13.4.1. *There exists tensor product of $D$-vector spaces.*

Proof. Let $F$ be field which for any $i$, $i = 1,...,n$, is subring of center $Z(D_i)$.

Let $D$ be free vector space over ring $F$ generated by product $D_1 \times ... \times D_n$ of multiplicative semigroups of division rings $D_1$, ..., $D_n$. Injection

$$i' : D_1 \times ... \times D_n \longrightarrow D$$

is defined according to rule

$$i'(d_1, ..., d_n) = (d_1, ..., d_n) \tag{13.4.1}$$

and is homomorphism of multiplicative semigroup of ring $D_1 \times ... \times D_n$ onto basis of module $D$. Therefore, product of vectors of basis is defined componentwise

$$(d_1, ..., d_n)(c_1, ..., c_n) = (d_1c_1, ..., d_nc_n) \tag{13.4.2}$$

According to theorem 13.1.1 the structure of ring is defined on vector space $D$.

Consider direct product $F \times D$ of field $F$ and ring $D$. We will identify element $(f, e)$ with element $f \in F$ and element $(1, d)$ with element $d \in D$.

Let $M$ be free module over ring $F \times D$ generated by Cartesian product $V_1 \times ... \times V_n$. Since $v_1 \in V_1$, ..., $v_n \in V_n$, then we denote corresponding vector from $M$ as $(v_1, ..., v_n)$. Let

$$i : V_1 \times ... \times V_n \longrightarrow M$$

be injection defined according to rule

$$i(v_1, ..., v_n) = (v_1, ..., v_n) \tag{13.4.3}$$

Let $N \subset M$ be vector subspace generated by elements of the following type

$$\begin{array}{l} (v_1, ..., v_i + w_i, ..., v_n) - (v_1, ..., v_i, ..., v_n) - (v_1, ..., w_i, ..., v_n) \\ (v_1, ..., av_i, ..., v_n) - a(v_1, ..., v_i, ..., v_n) \end{array} \tag{13.4.4}$$

[13.3] I give definition of tensor product of $D$-vector spaces following to definition in [1], p. 601 - 603.

where $v_i,\ w_i \in V_i,\ a \in F$. Let

$$j : M \to M/N$$

be canonical map on factor module. Consider commutative diagram

(13.4.5)

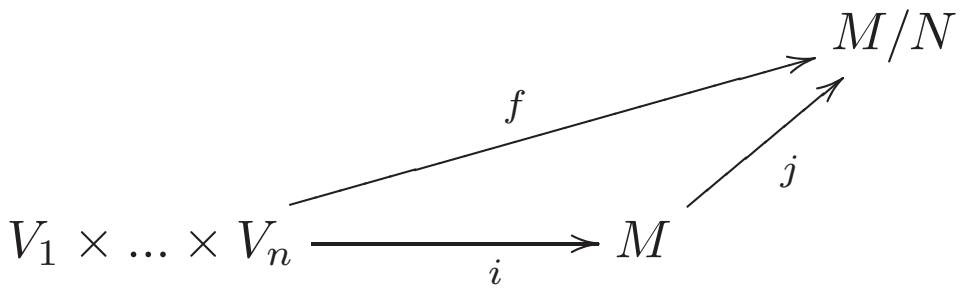

Since elements (13.4.4) belong to kernel of linear map $j$, then from equation (13.4.3) it follows

$$\begin{aligned} f(v_1, ..., v_i + w_i, ..., v_n) &= f(v_1, ..., v_i, ..., v_n) + f(v_1, ..., w_i, ..., v_n) \\ f(v_1, ..., av_i, ..., v_n) &= af(v_1, ..., v_i, ..., v_n) \end{aligned} \tag{13.4.6}$$

From equations (13.4.6), it follows that map $f$ is polylinear over field $F$. Since $M$ is module with basis $V_1 \times ... \times V_n$, then according to theorem [1]-4.1, on p. 135 for any module $V$ and any polylinear over $F$ map

$$g : V_1 \times ... \times V_n \longrightarrow V$$

there exists a unique homomorphism $k : M \to V$, for which following diagram is commutative

(13.4.7)

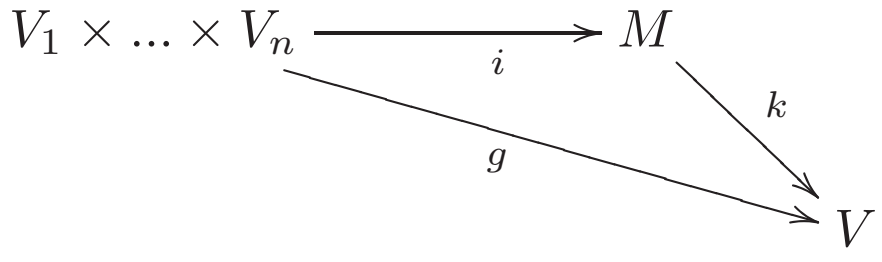

Since $g$ is polylinear over $F$, then $\ker k \subseteq N$. According to statement on p. [1]-119 map $j$ is universal in the category of homomorphisms of module $M$ whose kernel contains $N$. Therefore, we have homomorphism

$$h : M/N \to V$$

which makes the following diagram commutative

(13.4.8)

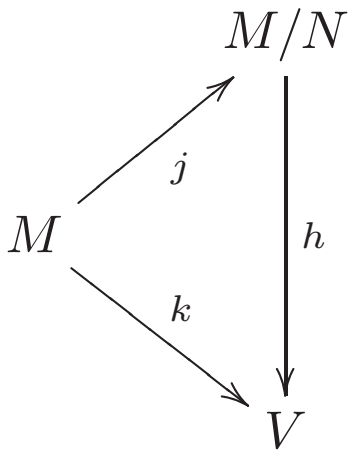

We join diagrams (13.4.5), (13.4.7), (13.4.8), and get commutative diagram

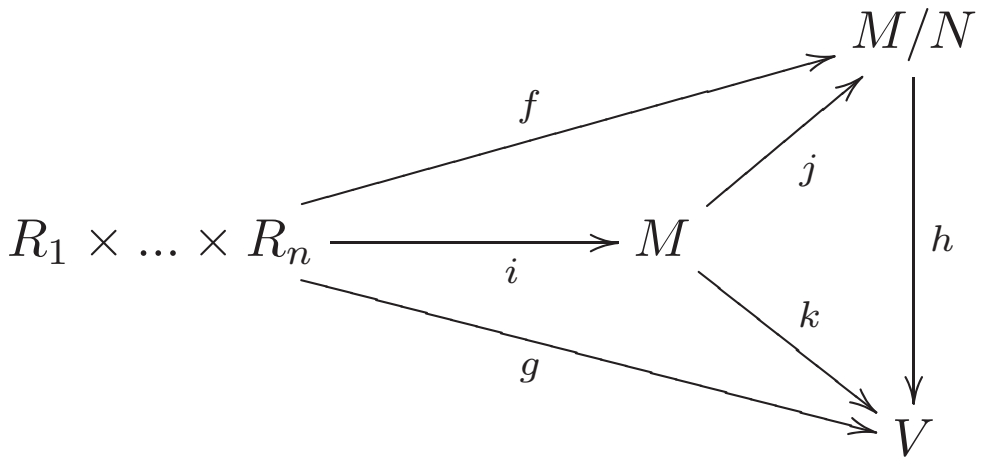

Since Im$f$ generates $M/N$, then map $h$ is uniquely determined. □

According to proof of theorem 13.4.1

$$V_1 \otimes ... \otimes V_n = M/N$$

If $v_i \in V_i$, we write

(13.4.9) $$v_1 \otimes ... \otimes v_n = j(v_1, ..., v_n)$$

From equations (13.4.3) and (13.4.9) it follows

$$f(v_1, ..., v_n) = v_1 \otimes ... \otimes v_n$$

We can write equations (13.4.6) as

$$\begin{aligned} & v_1 \otimes ... \otimes (v_i + w_i) \otimes ... \otimes v_n \\ = & v_1 \otimes ... \otimes v_i \otimes ... \otimes v_n + v_1 \otimes ... \otimes w_i \otimes ... \otimes v_n \\ & v_1 \otimes ... \otimes (av_i) \otimes ... \otimes v_n = a(v_1 \otimes ... \otimes v_i \otimes ... \otimes v_n) \end{aligned}$$

THEOREM 13.4.2. *Tensor product $V_1 \otimes ... \otimes V_n$ is bimodule over tensor product $D_1 \otimes ... \otimes D_n$.*

PROOF. Proof of statement of theorem is similar to proof of statement of theorem 13.3.2. □

CHAPTER 14

CHAPTER 15

# Index

CHAPTER 16

# Special Symbols and Notations

# Лекции по линейной алгебре над телом

Александр Клейн

Aleks_Kleyn@MailAPS.org
http://AleksKleyn.dyndns-home.com:4080/
http://sites.google.com/site/AleksKleyn/
http://arxiv.org/a/kleyn_a_1
http://AleksKleyn.blogspot.com/

Аннотация. В книге рассматриваются вопросы линейной алгебры над телом. Система линейных уравнений над телом имеет свойства, похожие на свойства систем линейных уравнений над полем. Тем не менее, некоммутативность произведения порождает новую картину.

Матрицы допускают две операции произведения, связанные операцией транспонирования. Бикольцо - это алгебра, определяющая на множестве две взаимосвязанные структуры кольца.

Подобно коммутативному случаю, решения системы линейных уравнений порождают правый или левое векторное пространство в зависимости от вида системы. Мы изучаем векторные пространства совместно с системами линейных уравнений потому, что существует тесная связь между их свойствами. Также как и в коммутативном случае, группа автоморфизмов векторного пространства имеет одно транзитивное представление на многообразии базисов, что даёт нам возможность определить пассивное и активное представления.

Изучение векторного пространства над телом раскрывает новые детали во взаимоотношении между пассивными и активными преобразованиями, делая картину более ясной.

Изучение парных представлений тела в абелевой группе приводит к концепции $D$-векторных пространств и их линейных отображений. На основе полилинейных отображений рассмотрено определение тензорного произведения колец и тензорное произведение $D$-векторных пространств.

# Оглавление

Глава 1

# Предисловие

## 1.1. Предисловие к изданию 1

Отправляясь в новое путешествие, не знаешь вначале, что ждёт тебя в пути. Знакомство с некоммутативной алгеброй началось с простого любопытства. В случае модуля над кольцом нельзя дать определение базиса подобно тому, как мы это делаем в случае векторного пространства над полем. Я хотел понять, как изменится картина, если вместо поля я буду рассматривать тело.

Прежде всего мне надо было научиться решать системы линейных уравнений. Я начал с системы двух уравнений с двумя неизвестными. Хотя решение было найдено легко, это решение нельзя выразить как отношение двух определителей.

Я понимал, что эта проблема интересует не только меня. Я начал искать математиков, интересующихся подобными задачами. Профессор Ретах познакомил меня со статьями [6, 7] по теории квазидетерминантов. Это было началом моего исследования в теории векторных пространств.

Глава 2 посвящена бикольцу матриц. Существует две причины, почему я рассматриваю эту алгебру.

Если в векторном пространстве задан базис, то преобразование векторного пространства можно описать с помощью матрицы. Произведению преобразований соответствует произведение матриц. В отличии от коммутативного случая мы не всегда можем представить это произведение в виде произведения строк первой матрицы на столбцы второй. В то же время запись элементов матрицы в виде

$$A=\begin{pmatrix}A_1^1 & ... & A_n^1\\ ... & ... & ...\\ A_1^n & ... & A_n^n\end{pmatrix}$$

является вопросом соглашения. Не нарушая общности, мы можем записать матрицу в виде

$$A=\begin{pmatrix}A_1^1 & ... & A_1^n\\ ... & ... & ...\\ A_n^1 & ... & A_n^n\end{pmatrix}$$

Но тогда изменяется правило умножения матриц.

Конструкции подобные бикольцу известны в алгебре. В структуре определены операции $a\vee b$ и $a\wedge b$, которые меняются местами, если порядок на множестве меняется на противоположный. Симметрия между ${}_*{}^*$-произведением и ${}^*{}_*$-произведением сформулирована в форме принципа двойственности.

Впоследствии я распространяю принцип двойственности на теорию представлений и теорию векторных пространств. Без учёта принципа двойственности эта книга была бы в четыре раза больше, не говоря о том, что бесконечные повторения текста сделали бы чтение текста невозможным.

Согласно каждому виду произведения мы можем расширить определение квазидетерминанта, данное в [6, 7], и определить два разных вида квазидетерминанта.

Глава 3 является обзором теории представлений и является базой для последующих глав. Мы распространим на представление алгебры соглашение, описанное в замечании 2.2.15. Теорема о существовании парных представлений в однородном пространстве завершает главу.

В главе 5 я изучаю несколько концепций линейной алгебры над телом $D$. Сперва я напоминаю определения векторного пространства и базиса.[1.1] Линейная алгебра над телом более богата фактами по сравнению с линейной алгеброй над полем. В отличии от векторных пространств над полем, мы можем определить левое и правое векторное пространство над произвольным телом $D$. Чтобы дать более строгое определение базиса, я излагаю теорию произвольной системы линейных уравнений (раздел 5.5) для каждого типа векторного пространства. Тем не менее, не смотря на это разнообразие, утверждения линейной алгебры над телом очень похожи на утверждения из линейной алгебры над полем.

Так как я пользуюсь утверждениями этой главы в геометрии, я следую той же форме записи, что мы используем в геометрии. При записи координат вектора и элементов матрицы мы будем следовать соглашению, описанному в разделе 2.1.

Мы отождествляем вектор с набором его координат относительно некоторого базиса. Тем не менее, это различные объекты. Чтобы подчеркнуть это различие, я вернулся к традиционному обозначению вектора в форме $\overline{a}$, когда вектор и его координаты присутствуют в одном и том же равенстве; в тоже время мы будем пользоваться записью $a^b$ для координат вектора $\overline{a}$. Мы используем одну и туже корневую букву для обозначения базиса и составляющих его векторов. Чтобы отличить вектор и базис, мы будем пользоваться обозначением $\overline{\overline{e}}$ для базиса и обозначением $\overline{e}_a$ для векторов, составляющих базис $\overline{\overline{e}}$. Изучая теорию векторных пространств мы будем пользоваться соглашением, описанным в замечании 2.2.15 и в разделе 4.1.

Глава 6 посвящена теории линейных представлений. Изучение однородного пространства группы симметрий $D_*{}^*$-векторного пространства ведёт нас к определению базиса этого $D_*{}^*$-векторного пространства и многообразия базисов. Мы вводим два типа преобразований многообразия базисов: активные и пассивные преобразования. Различие между ними состоит в том, что активное преобразование может быть выражено как преобразование исходного пространства. Как показано в [3], пассивное преобразование даёт возможность определить понятия инвариантности и геометрического объекта. Опираясь на эту теорию, я изучаю многообразие базисов, пассивные и активные преобразования в разделе 6.2. Я рассматриваю геометрический объект в $D_*{}^*$-векторном пространстве в разделе 6.3.

---

[1.1] Определения можно найти также в [8].

Мы имеем две противоположные точки зрения на геометрический объект. С одной стороны мы фиксируем координаты геометрического объекта относительно заданного базиса и указываем закон преобразования координат при замене базиса. В то же время мы рассматриваем всю совокупность координат геометрического объекта относительно различных базисов как единое целое. Это даёт нам возможность бескоординатного изучения геометрического объекта.

Глава 7 посвящена линейным отображениям ${}_*{}^*D$-векторных пространств над телом. Это явилось основой построения парных представлений кольца D. Несмотря на кажущуюся громоздкость, обозначения, которыми я пользуюсь в теории векторных пространств, опираются на существование парных представлений.

Линейная алгебра над телом более разнообразна по сравнению с линейной алгеброй над полем. Некоммутативность заставляет нас позаботиться о правильном порядке множителей в выражении. Это порождает не только разнообразие понятий, но также помогает лучше понять утверждения коммутативной алгебры.

Январь, 2007

## 1.2. Предисловие к изданию 2

> Любой преступник рано или поздно совершает ошибку.
>
> Когда детектив нашёл мою кепку на столе, он пришёл в ярость. В следующий раз был дождь, и я оставил мокрые следы на полу. Нечего говорить, опять пришлось начать сначала. И я начинал снова и снова. То тунгуский метеорит, то незапланированное затмение Солнца... Но в конце концов детектив не нашёл изъянов в дизайне.
>
> - Крестики и нолики, - печально вздохнул детектив. - Вам не надоело кодироать эту игру? Пора уже заняться чем-то более серьёзным.
>
> Автор неизвестен. Искуство программирования.

Научный поиск иногда напоминает программирование. Сложный и большой проект постепенно обрастает кодами. Всё кажется просто и понятно. Но вдруг появляется баг, исправление которого требует переделать почти половину кода. Впрочем, как заметил мудрец, в любой программе есть по крайней мере один баг.

Когда я начал исследование в области векторных пространств над телами, я понимал, что я ступаю на terra incognita. Концепция квазидетерминанта весьма облегчила мою задачу. Возможность решить систему линейных уравнений сделала эту теорию весьма похожей на теорию векторных пространств над полями. Это дало возможность построить базис и расширить концепцию геометрического объекта.

Тёмным облаком встали на горизонте двухвалентные тензоры. Казалось, некоммутативность накладывает на их структуру непреодолимые ограничения. Но именно в эту минуту я увидел свет в тунеле.

Изучение $D_*{}^*$-линейных отображений приводит к концепции парных представлений тела в абелевой группе. Абелева группа, в которой определено парное представление тела $D$, называется $D$-векторным пространством. Поскольку умножение на элементы тела $D$ определено и слева, и справа, гомоморфизмы $D$-векторных пространств не могут сохранять эту операцию. Это приводит к концепции аддитивного отображения, которое является морфизмом $D$-векторных пространств.[1.2]

$D$-векторное пространство является наиболее удачным претендентом для построения тензорной алгебры. Однако здесь потребовалась ещё одна важная конструкция.

Когда мы имеем дело с кольцами, мы изучаем не только модули, левые или правые, но и бимодули. Тело не является исключением. Но если добавить сюда концепцию бикольца, становится очевидным, что мне без разницы, с какой стороны умножать. Как следствие, оказывается, что бивекторное пространство - это частный случай прямого произведения $D_*{}^*$-векторных пространств. Мы не можем определить полилинейное отображение. Однако рассмотрение полиадитивных отображений прямого произведения $D_*{}^*$-векторных пространств приводит к определению тензорного произведения. Таким образом, можно построить тензора произвольной природы.

Статья эта построена таким образом, чтобы воспроизвести моё непростое и интересное путешествие в мир тензоров над телом.

Март, 2009

## 1.3. Предисловие к изданию 3

Алгебра кватернионов - самый простой пример тела. Поэтому я проверяю полученные мной теоремы в алгебре кватернионов, чтобы увидеть как они работают. Алгебра кватернионов похожа на поле комплексных чисел и поэтому естественно искать некоторые параллели.

В поле действительных чисел любое аддитивное отображение автоматически оказывается линейным над полем действительных чисел. Это связано с тем, что поле действительных чисел является пополнением поля рациональных чисел и является следствием теоремы 10.1.3.

Однако в случае комплексных чисел ситуация меняется. Не всякое аддитивное отображение поля комплексных чисел оказывается линейным над полем комплексных чисел. Операция сопряжения - простейший пример такого отображения. Когда я обнаружил этот крайне интересный факт, я вернулся к вопросу об аналитическом представлении аддитивного отображения.

[1.2] Построение аддитивного отображения, выполненное при доказательстве теоремы 12.1.4, весьма поучительно методически. С того самого момента, когда я начал изучать дифференциальную геометрию, я признал факт, что всё сводится к тензорам, что компоненты тензора нумеруются с помощью индексов и эта связь неразрывна. В тот момент, когда это оказалось не так, я оказался не готов признать это. Понадобился месяц прежде, чем я смог записать выражение, которое я видел глазами, но которое я не мог записать до этого, так как пытался выразить его в тензорной или операторной форме.

Изучая аддитивные отображения, я понял, что я слишком резко расширил множество линейных отображений при переходе от поля к телу. Причиной этому было не вполне ясное понимание как преодолеть некоммутативность произведения. Однако в процессе построений становилось всё более очевидным, что любое аддитивное отображение линейно над некоторым полем. Впервые эта концепция появилась при построении тензорных произведений и отчётливо проявилась в последующем исследовании.

При построении тензорного произведения тел $D_1$, ..., $D_n$ я предполагаю существование поля $F$, над которым аддитивное преобразование тела $D_i$ линейно для любого $i$. Если все тела имеют характеристику 0, то согласно теореме 10.1.3 такое поле всегда существует. Однако здесь возникает зависимость тензорного произведения от выбранного поля $F$. Чтобы избавиться от этой зависимости, я предполагаю, что поле $F$ - максимальное поле, обладающее указанным свойством.

Если $D_1 = ... = D_n = D$, то такое поле является центром $Z(D)$ тела $D$. Если произведение в теле $D$ коммутативно[1.3], то $Z(D) = D$. Следовательно, начав с аддитивного отображения, я пришёл к концепции линейного отображения, которая является обобщением линейного отображения над полем.

Исследование в области комплексных чисел и кватернионов проявили ещё одно интересное явление. Несмотря на то, что поле комплексных чисел является расширением поля действительных чисел, структура линейного отображения над полем комплексных чисел отличается от структуры линейного отображения над полем действительных чисел. Это различие приводит к тому, что операция сопряжения комплексных чисел является аддитивным отображением, но не является линейным отображением над полем комплексных чисел.

Аналогично, структура линейного отображения над телом кватернионов отличается от структуры линейного отображения над полем комплексных чисел. Причина различия в том, что центр алгебры кватернионов имеет более простую структуру, чем поле комплексных чисел. Это отличие приводит к тому, что операция сопряжения кватерниона удовлетворяет равенству

$$\overline{p} = -\frac{1}{2}(p + ipi + jpj + kpk)$$

Вследствие этого задача найти отображения, удовлетворяющие теореме подобной теореме Римана (теорема 11.1.1), является нетривиальной задачей для кватернионов.

Август, 2009

[1.3] Иными словами, тело является полем.

## 1.4. Предисловие к изданию 5

Пассажир: ... Мы же опоздаем на станцию!

Паровозик: ... Но если мы не увидим первых ландышей, то мы опоздаем на всю весну!

...
Пассажир: ... Мы же совсем опоздаем!
Паровозик: Да. Но если мы не услышим первых соловьёв, то мы опоздаем на всё лето!

...
Паровозик: Рассвет! ... каждый рассвет единственный в жизни! ... Ехать пора. Ведь мы опоздаем.
Пассажир: Да. Но если мы не увидим рассвет, мы можем опоздать на всю жизнь!

Геннадий Цыферов, Паровозик из Ромашково.

Молодой человек по имени Женя вышел из пункта $A$ в пункт $B$. Через два часа его товарищ Петя выехал из пункта $B$ на велосипеде по направлению к пункту $A$. С какой скоростью ехал Петя, если они встретились через час?..

Знакомая задача из школьного задачника. Но задача не столь проста как она сформулирована. Петя проезжал мимо пруда, увидел лебедей и остановился, чтобы их сфотографировать. А Женя остановился в роще послушать соловьёв. Там они и встретились.

Очень часто хочется заглянуть в конец задачника. Вот оно решение, уже близко. Однако что-то непознанное останется за рамками готового ответа. Я уже сейчас вижу в общих чертах теорию, которую я ищу. Но я понимаю, что прежде я смогу рассказать об этом другим, мне предстоит пройти длинный путь. Кроме того, мой опыт мне подсказывает, что когда бредёшь по незнакомой местности, внезапно увидишь столь удивительный ландшафт, что невольно задержишься, чтобы познакомиться с ним получше.

Впрочем, один раз я не удержался и в статье [9] попытался проанализировать, как может выглядеть геометрия над телом. И хотя статья выглядела незавершенной, её результаты легли в основу моих последующих исследований.

---

В статье [12] я рассматриваю систему аддитивных уравнений в $\Omega$-алгебре, в которой определена операция сложения. Формат записи уравнений в [12] определён форматом операций в $\Omega$-алгебре. Этот формат отличался от формата принятому в этой книге по умолчанию. Теория аддитивных уравнений является обобщением теории ${}_*{}^*D$-линейных уравнений и пользуется теми же методами. Другой порядок переменной и коэффициента в уравнении затруднял сравнение результатов в обеих книгах. Однако для меня не имеет значения, рассматриваю ли я систему $D^*{}_*$-линейных уравнений или систему ${}_*{}^*D$-линейных уравнений. Это побудило меня отредактировать текст этой книги и рассматривать систему ${}_*{}^*D$-линейных уравнений по умолчанию.

Кроме того я понял, что множество аддитивных отображений является понятием искусственным при изучении линейной алгебры над телом. Поэтому я при переработке этого издания удалил соответствующие определения и теоремы и рассматриваю непосредственно линейное отображение над заданным полем.

Основное отличие линейных отображений над телом состоит в том, что я не могу представить линейное отображение в виде произведения элемента тела на переменную. В статьях [11, 12] я нашёл способ функционального разделения переменной и отображения. Это позволяет представить запись линейного отображения в знакомой форме. Однако в этой книге я сохранил прежнюю форму записи. Причин этому несколько.

Прежде всего, обе формы записи верны. Мы можем пользоваться любой из этих форм записи, выбирая ту, которая кажется более адекватной. Это даёт также возможность читателю сопоставлять обе формы записи. Кроме того, в [12] я рассматриваю линейное отображение алгебр параллельно с их тензорным произведением, так как эти две темы неразделимы. В этой книге я решил сохранить тот путь, который я прошёл при написании второго издания этой книги и сохранить тензорное произведение как своеобразную вершину, к которой я подымаюсь в конце книги.

Тем не менее, путь, пройденный после последнего издания, отразился на последних главах книги. Например, поскольку тело является алгеброй над полем, я решил перейти к стандартному представлению индексов. Основные изменения коснулись моего понимания $D$-векторных пространств. Я остановлюсь на этом вопросе немного подробнее, поскольку эта тема осталась за пределами книги.

$D$-векторное пространство возникло как следствие теоремы о парных представлениях тела в абелевой группе (теоремы 7.3.2, 7.3.3). Вначале я рассматривал $D$-векторное пространство как объединение структур $D^*{}_*$-векторного пространства и ${}_*{}^*D$-векторного пространства. Соответственно я рассматриваю $D^*{}_*$-базис и ${}_*{}^*D$-базис.

Допустим, что $\overline{\overline{p}}$ - базис ${}_*{}^*D$-векторного пространства $V$ и $\overline{\overline{r}}$ - базис $D^*{}_*$-векторного пространства $V$. Вектор ${}_ir$ базиса $\overline{\overline{r}}$ имеет разложение

$$ {}_ir = p_j\,{}_i^jR \quad {}_ir = p_*{}^*{}_iR \tag{1.4.1} $$

относительно базиса $\overline{\overline{p}}$. Вектор $p_j$ базиса $\overline{\overline{p}}$ имеет разложение

$$ p_j = P_j^i\,{}_ir \quad p_j = P_j{}^*{}_*r \tag{1.4.2} $$

относительно базиса $\overline{\overline{r}}$.

Нетрудно видеть из построения, что $R$ - координатная матрица базиса $\overline{\overline{r}}$ относительно базиса $\overline{\overline{p}}$. Столбцы матрицы $R$ $D^*{}_*$-линейно независимы.

Аналогично, $P$ - координатная матрица базиса $\overline{\overline{p}}$ относительно базиса $\overline{\overline{r}}$. Столбцы матрицы $P$ ${}_*{}^*D$-линейно независимы.

Из равенств (1.4.1) и (1.4.2) следует

$$ \begin{aligned} {}_ir &= p_j\,{}_i^jR = P_j^k\,{}_kr\,{}_i^jR \\ {}_ir &= p_*{}^*{}_iR = (P^*{}_*r)_*{}^*{}_iR \end{aligned} \tag{1.4.3} $$

Из равенства (1.4.3) видно, что порядок скобок существенен.

Хотя матрицы $P$ и $R$ не являются взаимно обратными, мы можем сказать, что равенство (1.4.3) описывает тождественное преобразование $D$-векторного пространства. Это преобразование можно записать также в виде

$$(1.4.4) \qquad \begin{aligned} p_j &= P_j^i\ {}_i r = P_j^i\ (p_k\ {}^k_i R) \\ p_j &= P_j{}^*{}_* r = P_j{}^*{}_*(p_*{}^* R) \end{aligned}$$

Из сравнения равенств (1.4.3) и (1.4.4) следует, что изменение порядка скобок приводят к изменению порядка суммирования. Эти равенства выражают симметрию в выборе базисов $\overline{\overline{p}}$ и $\overline{\overline{r}}$.

Нетривиальная запись тождественного отображения в виде (1.4.3) либо (1.4.4) явилась для меня полнейшей неожиданностью.

Другое интересное утверждение, которое я сохранил для последующего анализа - это теорема 7.3.4, которая утверждает существование $D^*{}_*$-базиса, который не является ${}_*{}^*D$-базисом. Здесь мы сделаем остановку.

Согласно теории представлений $\Omega$-алгебры, рассмотренной в [10], $D$-векторное пространство является представлением алгебры $D \otimes D$ в абелевой группе. Следовательно, выбор базиса определяется не $D\star$-линейной зависимостью или $\star D$-линейной зависимостью, а зависимостью вида

$$a^i_{s\cdot 0}\overline{e}_i a^i_{s\cdot 1}$$

Следовательно, $D^*{}_*$-базис, который не является ${}_*{}^*D$-базисом, не является также базисом $D$-векторного пространства. Тем не менее не трудно убедиться, что множество преобразований базиса порождает группу. Эта структура базиса приводит к тому, что исчезает разница между правым и левым базисами. Поэтому я для $D$-векторного пространства пользуюсь стандартным представлением индекса.

Однако из этого следует ещё более удивительное утверждение. Координаты вектора $D$-векторного пространства принадлежат алгебре $D \otimes D$. Сразу появляются вопросы.

- Может ли базис $D$-векторного пространства быть ${}_*{}^*D$-базисом?
- Каково соотношение размерностей $D$-векторного пространства и ${}_*{}^*D$-векторного пространства?
- Какова структура 1-мерного $D$-векторного пространства?

Я планирую рассмотреть эти вопросы для свободных модулей над произвольной алгеброй. Это основная причина почему эта тема лежит вне рамок этой книги.

Однако уже сейчас можно сказать, что множество векторов, которое ${}_*{}^*D$-линейно зависимо, не может быть базисом $D$-векторного пространства. Поэтому размерность $D$-векторного пространства не превышает размерность соответствующего $D_*{}^*$-векторного пространства. А следовательно описание $D$-векторного пространства, предложенное в данной книге остаётся верным.

Что ж. Я не удержался, и заглянул в конец задачника.

Август, 2010

## 1.5. Соглашения

Соглашение 1.5.1. *В любом выражении, где появляется индекс, я предполагаю, что этот индекс может иметь внутреннюю структуру. Например,*

*при рассмотрении алгебры $A$ координаты $a \in A$ относительно базиса $\overline{\overline{e}}$ пронумерованы индексом $i$. Это означает, что $a$ является вектором. Однако, если $a$ является матрицей, нам необходимо два индекса, один нумерует строки, другой - столбцы. В том случае, когда мы уточняем структуру индекса, мы будем начинать индекс с символа $\cdot$ в соответствующей позиции. Например, если я рассматриваю матрицу $a^i_j$ как элемент векторного пространства, то я могу записать элемент матрицы в виде $a^{\cdot i}_{\ j}$.* □

Соглашение 1.5.2. *В выражении вида*

$$a_{s\cdot 0} x a_{s\cdot 1}$$

*предполагается сумма по индексу $s$.* □

Соглашение 1.5.3. *Тело $D$ можно рассматривать как $D$-векторное пространство размерности 1. Соответственно этому, мы можем изучать не только гомоморфизм тела $D_1$ в тело $D_2$, но и линейное отображение тел.* □

Соглашение 1.5.4. *Несмотря на некоммутативность произведения многие утверждения сохраняются, если заменить например правое представление на левое представление или правое векторное пространство на левое векторное пространство. Чтобы сохранить эту симметрию в формулировках теорем я пользуюсь симметричными обозначениями. Например, я рассматриваю $D\star$-векторное пространство и $\star D$-векторное пространство. Запись $D\star$-векторное пространство можно прочесть как $D$-star-векторное пространство либо как левое векторное пространство.* □

Соглашение 1.5.5. *Пусть $A$ - свободная алгебра с конечным или счётным базисом. При разложении элемента алгебры $A$ относительно базиса $\overline{\overline{e}}$ мы пользуемся одной и той же корневой буквой для обозначения этого элемента и его координат. В выражении $a^2$ не ясно - это компонента разложения элемента $a$ относительно базиса или это операция возведения в степень. Для облегчения чтения текста мы будем индекс элемента алгебры выделять цветом. Например,*

$$a = a^{i} e_{i}$$

□

Соглашение 1.5.6. *Очень трудно провести границу между модулем и алгеброй. Тем более, иногда в процессе построения мы должны сперва доказать, что множество $A$ является модулем, а потом мы доказываем, что это множество является алгеброй. Поэтому для записи координат элемента модуля мы также будем пользоваться соглашением 1.5.5.* □

Соглашение 1.5.7. *Отождествление вектора и матрицы его координат может привести к неоднозначности в равенстве*

$$a = a^*{}_* e \tag{1.5.1}$$

*где $\overline{\overline{e}}$ - базис векторного пространства. Поэтому мы будем записывать равенство (1.5.1) в виде*

$$\overline{a} = a^*{}_* e$$

*чтобы видеть, где записан вектор.* □

Соглашение 1.5.8. *Если свободная конечномерная алгебра имеет единицу, то мы будем отождествлять вектор базиса* $\overline{e}_0$ *с единицей алгебры.* □

Без сомнения, у читателя могут быть вопросы, замечания, возражения. Я буду признателен любому отзыву.

Глава 2

# Бикольцо матриц

## 2.1. Концепция обобщённого индекса

Изучая тензорное исчисление, мы начинаем с изучения одновалентных ковариантного и контравариантного тензоров. Несмотря на различие свойств, оба эти объекта являются элементами соответствующих векторных пространств. Если мы введём обобщённый индекс по правилу $a^i = a^i$, $b^i = b^{\cdot -}_{\ \ i}$, то мы видим, что эти тензоры ведут себя одинаково. Например, преобразование ковариантного тензора принимает форму

$$b'^i = b'^{\cdot -}_{\ \ i} = f^{\cdot -}_{\ \ i}{}^{\ j}_{\cdot -} b^{\cdot -}_{\ \ j} = f^i_j b^j$$

Это сходство идёт сколь угодно далеко, так как тензоры также порождают векторное пространство.

Эти наблюдения сходства свойств ковариантного и контравариантного тензоров приводят нас к концепции обобщённого индекса. Я пользуюсь символом $\cdot$ перед обобщённым индексом, когда мне необходимо описать его структуру. Я помещаю символ $'-'$ на месте индекса, позиция которого изменилась. Например, если исходное выражение было $a_{ij}$, я пользуюсь записью $a_i{}^j_-$ вместо записи $a^j_i$.

Хотя структура обобщённого индекса произвольна, мы будем предполагать, что существует взаимно однозначное отображение отрезка натуральных чисел 1, ..., $n$ на множество значений индекса. Пусть $i$ - множество значений индекса $i$. Мы будем обозначать мощность этого множества символом $|i|$ и будем полагать $|i| = n$. Если нам надо перечислить элементы $a_i$, мы будем пользоваться обозначением $a_1$, ..., $a_n$.

Представление координат вектора в форме матрицы позволяет сделать запись более компактной. Вопрос о представлении вектора как строка или столбец матрицы является вопросом соглашения. Мы можем распространить концепцию обобщённого индекса на элементы матрицы. Матрица - это двумерная таблица, строки и столбцы которой занумерованы обобщёнными индексами. Для представления матрицы мы будем пользоваться одним из следующих представлений:

**Стандартное представление:** в этом случае мы представляем элементы матрицы $A$ в виде $A^a_b$.

**Альтернативное представление:** в этом случае мы представляем элементы матрицы $A$ в виде ${}^aA_b$ или ${}_bA^a$.

Так как мы пользуемся обобщёнными индексами, мы не можем сказать, нумерует ли индекс $a$ строки матрицы или столбцы, до тех пор, пока мы не знаем структуры индекса.

Мы могли бы пользоваться терминами $*$-столбец и $*$-строка, которые более близки традиционным. Однако как мы увидим ниже для нас несущественна форма представления матрицы. Для того, чтобы обозначения, предлагаемые ниже, были согласованы с традиционными, мы будем предполагать, что матрица представлена в виде

$$A = \begin{pmatrix} {}^1A_1 & ... & {}^1A_n \\ ... & ... & ... \\ {}^mA_1 & ... & {}^mA_n \end{pmatrix}$$

Определение 2.1.1. Я использую следующие имена и обозначения различных **минорных матриц** матрицы $A$

$A_a$ **:** ${}^*$**-строка** с индексом $a$ является обобщением столбца матрицы. Верхний индекс перечисляет элементы ${}^*$-строки, нижний индекс перечисляет ${}^*$-строки.

$A_T$ **:** минорная матрица, полученная из $A$ выбором ${}^*$-строк с индексом из множества $T$

$A_{[a]}$ **:** минорная матрица, полученный из $A$ удалением ${}^*$-строки $A_a$

$A_{[T]}$ **:** минорная матрица, полученный из $A$ удалением ${}^*$-строк с индексом из множества $T$

${}^bA$ **:** ${}_*$**-строка** с индексом $b$ является обобщением строки матрицы. Нижний индекс перечисляет элементы ${}_*$-строки, верхний индекс перечисляет ${}_*$-строки.

${}^SA$ **:** минорная матрица, полученный из $A$ выбором ${}_*$-строк с индексом из множества $S$

${}^{[b]}A$ **:** минорная матрица, полученный из $A$ удалением ${}_*$-строки ${}^bA$

${}^{[S]}A$ **:** минорная матрица, полученный из $A$ удалением ${}_*$-строк с индексом из множества $S$

□

Замечание 2.1.2. Мы будем комбинировать запись индексов. Так ${}^bA_a$ является $1\times 1$ минорной матрицей. Одновременно, это обозначение элемента матрицы. Это позволяет отождествить $1\times 1$ матрицу и её элемент. Индекс $a$ является номером ${}^*$-строки матрицы и индекс $b$ является номером ${}_*$-строки матрицы. □

Каждая форма записи матрицы имеет свои преимущества. Стандартная форма более естественна, когда мы изучаем теорию матриц. Альтернативная форма записи делает выражения более ясными в теории векторных пространств. Распространив альтернативную запись индексов на произвольные тензоры, мы сможем лучше понять взаимодействие различных геометрических объектов. Опираясь на принцип двойственности (теорема 2.2.14), мы можем расширить наши выразительные возможности.

Замечание 2.1.3. Мы можем договориться, что при чтении мы произносим символ ${}^*$- как $c$- и символ ${}_*$- как $r$-, формируя тем самым названия $c$**-строка** и $r$ **-строка**. В последующем мы распространим это соглашение на другие элементы линейной алгебры. Я буду пользоваться этим соглашением при составлении индекса. □

Так как транспонирование матрицы меняет местами ${}_*$-строки и ${}^*$-строки, то мы получаем равенство

$$ {}_i(A^T)^j = {}^iA_j \tag{2.1.1} $$

Замечание 2.1.4. Как видно из равенства (2.1.1), для нас несущественно, с какой стороны мы указываем номер ${}_*$-строки и с какой стороны мы указываем номер ${}^*$-строки. Это связано с тем, что мы можем нумеровать элементы матрицы различными способами. Если мы хотим указывать номера ${}_*$-строки и ${}^*$-строки согласно определению 2.1.1, то равенство (2.1.1) примет вид

$$ {}^j(A^T)_i = {}^iA_j $$

В стандартном представлении равенство (2.1.1) примет вид

$$ (A^T)^j_i = A^i_j $$

□

Мы называем матрицу[2.1]

$$ \mathcal{H}A = ({}^j(\mathcal{H}A)_i) = (({}^{\cdot}{}_i^{-}A_{\cdot}{}^j_{-})^{-1}) \tag{2.1.2} $$

**обращением Адамара матрицы** $A = ({}_bA^a)$ ([7]-page 4).

Я пользуюсь Эйнштейновским соглашением о суммировании. Это означает, что, когда индекс присутствует в выражении дважды и множество индексов известно, у меня есть сумма по этому индексу. Я буду явно указывать множество индексов, если это необходимо. Кроме того, в этой статье я использую туже корневую букву для матрицы и её элементов.

Мы будем изучать матрицы, элементы которых принадлежат телу $D$. Мы также будем иметь в виду, что вместо тела $D$ мы можем в тексте писать поле $F$. Мы будем явно писать поле $F$ в тех случаях, когда коммутативность будет порождать новые детали. Обозначим через 1 единичный элемент тела $D$.

Пусть $I$, $|I| = n$ - множество индексов. **Символ Кронекера** определён равенством

$$ \delta^i_j = \begin{cases} 1 & i = j \\ 0 & i \neq j \end{cases} \qquad i, j \in I $$

## 2.2. Бикольцо

Мы будем рассматривать матрицы, элементы которых принадлежат телу $D$.

Произведение матриц связано с произведением гомоморфизмов векторных пространств над полем. Согласно традиции произведение матриц $A$ и $B$ определено как произведение ${}_*$-строк матрицы $A$ и ${}^*$-строк матрицы $B$. Условность этого определения становится очевидной, если мы обратим внимание, что ${}_*$-строка матрицы $A$ может быть столбцом этой матрицы. В этом случае мы

[2.1] Запись $({}^{\cdot}{}_i^{-}A_{\cdot}{}^j_{-})^{-1}$ означает, что при обращении Адамара столбцы и строки меняются местами. Мы можем формально записать это выражение следующим образом

$$ ({}^{\cdot}{}_i^{-}A_{\cdot}{}^j_{-})^{-1} = \frac{1}{{}^iA_j} $$

умножаем столбцы матрицы $A$ на строки матрицы $B$. Таким образом, мы можем определить два вида произведения матриц. Чтобы различать эти произведения, мы вводим новые обозначения.[2.2]

Определение 2.2.1. ${}_*{}^*$**-произведение матриц** $A$ и $B$ имеет вид

$$\begin{cases} A_*{}^*B = ({}^aA_c\ {}^cB_b) \\ {}^a(A_*{}^*B)_b = {}^aA_c\ {}^cB_b \end{cases} \tag{2.2.1}$$

и может быть выражено как произведение ${}_*$-строки матрицы $A$ и ${}^*$-строки матрицы $B$.[2.3] □

Определение 2.2.2. ${}^*{}_*$**-произведение матриц** $A$ и $B$ имеет вид

$$\begin{cases} A^*{}_*B = ({}_aA^c\ {}_cB^b) \\ {}_a(A^*{}_*B)^b = {}_aA^c\ {}_cB^b \end{cases} \tag{2.2.2}$$

и может быть выражено как произведение ${}^*$-строки матрицы $A$ на ${}_*$-строку матрицы $B$.[2.4] □

Замечание 2.2.3. Мы будем пользоваться символом ${}_*{}^*$- или ${}^*{}_*$- в имени свойств каждого произведения и в обозначениях. Согласно замечанию 2.1.3 мы можем читать символы ${}_*{}^*$ и ${}^*{}_*$ как $rc$-произведение и $cr$-произведение. Это правило мы распространим на последующую терминологию. □

Замечание 2.2.4. Также как и в замечании 2.1.4, я хочу обратить внимание на то, что я меняю нумерацию элементов матрицы. Если мы хотим указывать номера ${}_*$-строки и ${}^*$-строки согласно определению 2.1.1, то равенство (2.2.2) примет вид

$$ {}^b(A^*{}_*B)_a = {}^cA_a\ {}^bB_c \tag{2.2.3}$$

Однако формат равенства (2.2.3) несколько необычен. □

[2.2]Для совместимости обозначений с существующими мы будем иметь в виду ${}_*{}^*$-произведение, когда нет явных обозначений.

[2.3]В альтернативной форме операция состоит из двух символов $*$, которые записываются на месте индекса суммирования. В стандартной форме операция имеет вид

$$\begin{cases} A_*{}^*B = (A_c^aB_b^c) \\ (A_*{}^*B)_b^a = A_c^aB_b^c \end{cases}$$

и может быть интерпретирована как символическая запись

$$A_*{}^*B = A_*B^*$$

где мы записываем символ $*$ на месте индекса, по которому предполагается суммирование.

[2.4]В альтернативной форме операция состоит из двух символов $*$, которые записываются на месте индекса суммирования. В стандартной форме операция имеет вид

$$\begin{cases} A^*{}_*B = (A_a^c\ B_c^b) \\ (A^*{}_*B)_b^a = A_a^c\ B_c^b \end{cases}$$

и может быть интерпретирована как символическую запись

$$A^*{}_*B = A^*B_*$$

где мы записываем символ $*$ на месте индекса, по которому предполагается суммирование.

Множество $n\times n$ матриц замкнуто относительно ${}_*{}^*$-произведения и ${}^*{}_*$-произведения, а также относительно суммы, определённой согласно правилу

$$(A+B)_a^b=A_a^b+B_a^b$$

Теорема 2.2.5.

$$(A_*{}^*B)^T=A^T{}^*{}_*B^T \tag{2.2.4}$$

Доказательство. Цепочка равенств

$$\begin{aligned}{}_a((A_*{}^*B)^T)^b&={}^a(A_*{}^*B)_b\\&={}^aA_c\ {}^cB_b\\&={}_a(A^T)^c{}_c(B^T)^b\\&={}_a((A^T)^*{}_*(B^T))^b\end{aligned} \tag{2.2.5}$$

следует из (2.1.1), (2.2.1) и (2.2.2). Равенство (2.2.4) следует из (2.2.5). □

Матрица $\delta=(\delta_a^c)$ является единицей для обоих произведений.

Определение 2.2.6. **Бикольцо** $\mathcal{A}$ - это множество, на котором мы определили унарную операцию, называемую транспозицией, и три бинарных операции, называемые ${}_*{}^*$-произведение, ${}^*{}_*$-произведение и сумма, такие что

- ${}_*{}^*$-произведение и сумма определяют структуру кольца на $\mathcal{A}$
- ${}^*{}_*$-произведение и сумма определяют структуру кольца на $\mathcal{A}$
- оба произведения имеют общую единицу $\delta$
- произведения удовлетворяют равенству

$$(A_*{}^*B)^T=A^T{}^*{}_*B^T \tag{2.2.6}$$

- транспозиция единицы есть единица

$$\delta^T=\delta \tag{2.2.7}$$

- двойная транспозиция есть исходный элемент

$$(A^T)^T=A \tag{2.2.8}$$

□

Теорема 2.2.7.

$$(A^*{}_*B)^T=(A^T)_*{}^*(B^T) \tag{2.2.9}$$

Доказательство. Мы можем доказать (2.2.9) в случае матриц тем же образом, что мы доказали (2.2.6). Тем не менее для нас более важно показать, что (2.2.9) следует непосредственно из (2.2.6).

Применяя (2.2.8) к каждому слагаемому в левой части (2.2.9), мы получим

$$(A^*{}_*B)^T=((A^T)^T{}^*{}_*(B^T)^T)^T \tag{2.2.10}$$

Из (2.2.10) и (2.2.6) следует, что

$$(A^*{}_*B)^T=((A^T{}_*{}^*B^T)^T)^T \tag{2.2.11}$$

(2.2.9) следует из (2.2.11) и (2.2.8). □

Определение 2.2.8. Мы определим ${}_*{}^*$**-степень** элемента $A$ бикольца $\mathcal{A}$, пользуясь рекурсивным правилом

$$A^{0_*{}^*}=\delta \tag{2.2.12}$$

$$A^{n_*{}^*}=A^{n-1_*{}^*}{}_*{}^*A \tag{2.2.13}$$

□

Определение 2.2.9. Мы определим ${}^*{}_*$**-степень** элемента $A$ бикольца $\mathcal{A}$, пользуясь рекурсивным правилом

$$A^{0^*{}_*}=\delta \tag{2.2.14}$$

$$A^{n^*{}_*}=A^{n-1^*{}_*}{}^*{}_*A \tag{2.2.15}$$

□

Теорема 2.2.10.

$$(A^T)^{n_*{}^*}=(A^{n^*{}_*})^T \tag{2.2.16}$$

$$(A^T)^{n^*{}_*}=(A^{n_*{}^*})^T \tag{2.2.17}$$

Доказательство. Мы проведём доказательство индукцией по $n$.

При $n=0$ утверждение непосредственно следует из равенств (2.2.12), (2.2.14) и (2.2.7).

Допустим утверждение справедливо при $n=k-1$

$$(A^T)^{n-1_*{}^*}=(A^{n-1^*{}_*})^T \tag{2.2.18}$$

Из (2.2.13) следует

$$(A^T)^{k_*{}^*}=(A^T)^{k-1_*{}^*}{}_*{}^*A^T \tag{2.2.19}$$

Из (2.2.19) и (2.2.18) следует

$$(A^T)^{k_*{}^*}=(A^{k-1^*{}_*})^T{}_*{}^*A^T \tag{2.2.20}$$

Из (2.2.20) и (2.2.9) следует

$$(A^T)^{k_*{}^*}=(A^{k-1^*{}_*}{}^*{}_*A)^T \tag{2.2.21}$$

Из (2.2.19) и (2.2.15) следует (2.2.16).

Мы можем доказать (2.2.17) подобным образом. □

Определение 2.2.11. Элемент $A^{-1_*{}^*}$ бикольца $\mathcal{A}$ - это ${}_*{}^*$**-обратный элемент** элемента $A$, если

$$A{}_*{}^*A^{-1_*{}^*}=\delta \tag{2.2.22}$$

Элемент $A^{-1^*{}_*}$ бикольца $\mathcal{A}$ - это ${}^*{}_*$**-обратный элемент** элемента $A$, если

$$A{}^*{}_*A^{-1^*{}_*}=\delta \tag{2.2.23}$$

□

Теорема 2.2.12. *Предположим, что элемент $A\in\mathcal{A}$ имеет ${}_*{}^*$-обратный элемент. Тогда транспонированный элемент $A^T$ имеет ${}^*{}_*$-обратный элемент и эти элементы удовлетворяют равенству*

$$(A^T)^{-1^*{}_*}=(A^{-1_*{}^*})^T \tag{2.2.24}$$

*Предположим, что элемент $A \in \mathcal{A}$ имеет ${}^*{}_*$-обратный элемент. Тогда транспонированный элемент $A^T$ имеет ${}_*{}^*$-обратный элемент и эти элементы удовлетворяет равенству*

$$(A^T)^{-1_*{}^*} = (A^{-1^*{}_*})^T \tag{2.2.25}$$

Доказательство. Если мы возьмём транспонирование обеих частей (2.2.22) и применим (2.2.7), мы получим

$$(A_*{}^*A^{-1_*{}^*})^T = \delta^T = \delta$$

Применяя (2.2.6), мы получим

$$\delta = A^T{}^*{}_*(A^{-1_*{}^*})^T \tag{2.2.26}$$

(2.2.24) следует из сравнения (2.2.23) и (2.2.26).

Мы можем доказать (2.2.25) подобным образом. □

Теоремы 2.2.5, 2.2.7, 2.2.10 и 2.2.12 показывают, что существует двойственность между ${}_*{}^*$-произведением и ${}^*{}_*$-произведением. Мы можем объединить эти утверждения.

Теорема 2.2.13 (**принцип двойственности для бикольца**). *Пусть $\mathfrak{A}$ - истинное утверждение о бикольце $\mathcal{A}$. Если мы заменим одновременно*

- *$A \in \mathcal{A}$ и $A^T$*
- *${}_*{}^*$-произведение и ${}^*{}_*$-произведение*

*то мы снова получим истинное утверждение.*

Теорема 2.2.14 (**принцип двойственности для бикольца матриц**). *Пусть $\mathcal{A}$ является бикольцом матриц. Пусть $\mathfrak{A}$ - истинное утверждение о матрицах. Если мы заменим одновременно*

- *${}^*$-строки и ${}_*$-строки всех матриц*
- *${}_*{}^*$-произведение и ${}^*{}_*$-произведение*

*то мы снова получим истинное утверждение.*

Доказательство. Непосредственное следствие теоремы 2.2.13. □

Замечание 2.2.15. В выражении

$$A_*{}^*B_*{}^*C$$

мы выполняем операцию умножения слева направо. Однако мы можем выполнять операцию умножения справа налево. В традиционной записи это выражение примет вид

$$C^*{}_*B^*{}_*A$$

Мы сохраним правило, что показатель степени записывается справа от выражения. Если мы пользуемся стандартным представлением, то индексы также записываются справа от выражения. Если мы пользуемся альтернативным представлением, то индексы читаются в том же порядке, что и символы операции и корневые буквы. Например, если исходное выражение имеет вид

$$A^{-1_*{}^*}{}_*{}^*B_a$$

то выражение, читаемое справа налево, примет вид

$$B_a{}^*{}_*A^{-1^*{}_*}$$

в стандартном представлении либо примет вид

$$ {}_aB^*{}_*A^{-1^*{}_*} $$

в альтернативном представлении.

Если задать порядок, в котором мы записываем индексы, то мы можем утверждать, что мы читаем выражение сверху вниз, читая сперва верхние индексы, потом нижние. Договорившись, что это стандартная форма чтения, мы можем прочесть выражение снизу вверх. При этом мы дополним правило, что символы операции также читаются в том же направлении, что и индексы. Например, выражение

$$ A^a{}_*{}^*B^{-1_*{}^*} = C^a $$

прочтённое снизу вверх, в стандартной форме имеет вид

$$ A_a{}^*{}_*B^{-1^*{}_*} = C_a $$

Согласно принципу двойственности, если верно одно утверждение, то верно и другое. □

Теорема 2.2.16. *Если матрица $A$ имеет ${}_*{}^*$-обратную матрицу, то для любых матриц $B$ и $C$ из равенства*

$$ B_*{}^*A = C_*{}^*A \tag{2.2.27} $$

*следует равенство*

$$ B = C \tag{2.2.28} $$

Доказательство. Равенство (2.2.28) следует из (2.2.27), если обе части равенства (2.2.27) умножить на $A^{-1_*{}^*}$. □

## 2.3. Квазидетерминант

Теорема 2.3.1. *Предположим, что $n \times n$ матрица $A$ имеет ${}_*{}^*$-обратную матрицу.*[2.5] *Тогда $k \times k$ минорная матрица ${}_*{}^*$-обратной матрицы удовлетворяет равенству*

$$ \left({}^I(A^{-1_*{}^*})_J\right)^{-1_*{}^*} = {}^JA_I - {}^JA_{[I]*}{}^*\left({}^{[J]}A_{[I]}\right)^{-1_*{}^*}{}_*{}^{*[J]}A_I \tag{2.3.1} $$

Доказательство. Определение (2.2.22) ${}_*{}^*$-обратной матрицы приводит к системе линейных уравнений

$$ {}^{[J]}A_{[I]*}{}^{*[I]}(A^{-1_*{}^*})_J + {}^{[J]}A_{I*}{}^{*I}(A^{-1_*{}^*})_J = 0 \tag{2.3.2} $$

$$ {}^JA_{[I]*}{}^{*[I]}(A^{-1_*{}^*})_J + {}^JA_{I*}{}^{*I}(A^{-1_*{}^*})_J = \delta \tag{2.3.3} $$

Мы умножим (2.3.2) на $\left({}^{[J]}A_{[I]}\right)^{-1_*{}^*}$

$$ {}^{[I]}(A^{-1_*{}^*})_J + \left({}^{[J]}A_{[I]}\right)^{-1_*{}^*}{}_*{}^{*[J]}A_{I*}{}^{*I}(A^{-1_*{}^*})_J = 0 \tag{2.3.4} $$

Теперь мы можем подставить (2.3.4) в (2.3.3)

$$ -{}^JA_{[I]*}{}^*\left({}^{[J]}A_{[I]}\right)^{-1_*{}^*}{}_*{}^{*[J]}A_{I*}{}^{*I}(A^{-1_*{}^*})_J + {}^JA_{I*}{}^{*I}(A^{-1_*{}^*})_J = \delta \tag{2.3.5} $$

(2.3.1) следует из (2.3.5). □

[2.5] Это утверждение и его доказательство основаны на утверждении 1.2.1 из [6] (page 8) для матриц над свободным кольцом с делением.

Следствие 2.3.2. *Предположим, что $n\times n$ матрица $A$ имеет ${}_*{}^*$-обратную матрицу. Тогда элементы ${}_*{}^*$-обратной матрицы удовлетворяют равенству*[2.1]

$$ {}^i(A^{-1_*{}^*})_j=\left({}^jA_i-{}^jA_{[i]}{}_*{}^*\left({}^{[j]}A_{[i]}\right)^{-1_*{}^*}{}_*{}^*{}^{[j]}A_i\right)^{-1} \tag{2.3.6}$$

$$ {}^j\left(\mathcal{H}A^{-1_*{}^*}\right)_i={}^jA_i-{}^jA_{[i]}{}_*{}^*\left({}^{[j]}A_{[i]}\right)^{-1_*{}^*}{}_*{}^*{}^{[j]}A_i \tag{2.3.7}$$

□

Пример 2.3.3. Рассмотрим матрицу

$$\begin{pmatrix}{}^1A_1&{}^1A_2\\{}^2A_1&{}^2A_2\end{pmatrix}$$

Согласно (2.3.6)

$$ {}^1(A^{-1_*{}^*})_1=({}^1A_1-{}^1A_2({}^2A_2)^{-1}\ {}^2A_1)^{-1} \tag{2.3.8}$$

$$ {}^2(A^{-1_*{}^*})_1=({}^1A_2-{}^1A_1({}^2A_1)^{-1}\ {}^2A_2)^{-1} \tag{2.3.9}$$

$$ {}^1(A^{-1_*{}^*})_2=({}^2A_1-{}^2A_2({}^1A_2)^{-1}\ {}^1A_1)^{-1} \tag{2.3.10}$$

$$ {}^2(A^{-1_*{}^*})_2=({}^2A_2-{}^2A_1({}^1A_1)^{-1}\ {}^1A_2)^{-1} \tag{2.3.11}$$

$$A^{-1^*{}_*}=\begin{pmatrix}({}_1A^1-{}_1A^2({}_2A^2)^{-1}\ {}_2A^1)^{-1}&({}_1A^2-{}_1A^1({}_2A^1)^{-1}\ {}_2A^2)^{-1}\\({}_2A^1-{}_2A^2({}_1A^2)^{-1}\ {}_1A^1)^{-1}&({}_2A^2-{}_2A^1({}_1A^1)^{-1}\ {}_1A^2)^{-1}\end{pmatrix}$$

□

Согласно [6], page 3 у нас нет определения детерминанта в случае тела. Тем не менее, мы можем определить квазидетерминант, который в конечном итоге даёт похожую картину. В определении, данном ниже, мы следуем определению [6]-1.2.2.

Определение 2.3.4. $\binom{j}{i}$**-**${}_*{}^*$**-квазидетерминант** $n\times n$ матрицы $A$ - это формальное выражение[2.1]

$$ {}^j\det({}_*{}^*)_i\,A={}^j\left(\mathcal{H}A^{-1_*{}^*}\right)_i \tag{2.3.12}$$

Согласно замечанию 2.1.2 мы можем рассматривать $\binom{j}{i}$-${}_*{}^*$-квазидетерминант как элемент матрицы $\det({}_*{}^*)\,A$ , которую мы будем называть ${}_*{}^*$**-квазидетерминантом**. □

Теорема 2.3.5. *Выражение для ${}_*{}^*$-обратной матрицы имеет вид*

$$A^{-1_*{}^*}=\mathcal{H}\det({}_*{}^*)A \tag{2.3.13}$$

Доказательство. (2.3.13) следует из (2.3.12). □

Теорема 2.3.6. *Выражение для* $\binom{j}{i}$-${}_*{}^*$*-квазидетерминанта имеет любую из следующих форм*[2.6]

$$ {}^j\det({}_*{}^*)_i\, A = {}^jA_i - {}^jA_{[i]}{}_*{}^* \left({}^{[j]}A_{[i]}\right)^{-1_*{}^*} {}_*{}^*{}^{[j]}A_i \tag{2.3.14}$$

$$ {}^j\det({}_*{}^*)_i\, A = {}^jA_i - {}^jA_{[i]}{}_*{}^*\mathcal{H}\det({}_*{}^*)^{[j]}A_{[i]}{}_*{}^*{}^{[j]}A_i \tag{2.3.15}$$

Доказательство. Утверждение следует из (2.3.7) и (2.3.12). □

Пример 2.3.7. Рассмотрим матрицу

$$\begin{pmatrix} {}^1A_1 & {}^1A_2 \\ {}^2A_1 & {}^2A_2 \end{pmatrix}$$

Согласно (2.3.14)

$$\det({}_*{}^*)A = \begin{pmatrix} {}^1A_1 - {}^1A_2({}^2A_2)^{-1}\ {}^2A_1 & {}^1A_2 - {}^1A_1({}^2A_1)^{-1}\ {}^2A_2 \\ {}^2A_1 - {}^2A_2({}^1A_2)^{-1}\ {}^1A_1 & {}^2A_2 - {}^2A_1({}^1A_1)^{-1}\ {}^1A_2 \end{pmatrix}$$

$$\det({}^*{}_*)A = \begin{pmatrix} {}_1A^1 - {}_1A^2({}_2A^2)^{-1}\ {}_2A^1 & {}_2A^1 - {}_2A^2({}_1A^2)^{-1}\ {}_1A^1 \\ {}_1A^2 - {}_1A^1({}_2A^1)^{-1}\ {}_2A^2 & {}_2A^2 - {}_2A^1({}_1A^1)^{-1}\ {}_1A^2 \end{pmatrix}$$

□

Теорема 2.3.8.

$$ {}_j\det({}_*{}^*)^i\, A^T = {}^j\det({}^*{}_*)_i\, A \tag{2.3.16}$$

Доказательство. Согласно (2.3.12) и (2.1.2)

$$ {}_j\det({}_*{}^*)^i\, A^T = ({}_{\cdot}{}^j{}_{-}((A^T)^{-1_*{}^*}){}^{\cdot}{}_{i}{}^{-})^{-1}$$

Пользуясь теоремой 2.2.12, мы получим

$$ {}_j\det({}_*{}^*)^i\, A^T = ({}_{\cdot}{}^j{}_{-}((A^{-1^*{}_*})^T){}^{\cdot}{}_{i}{}^{-})^{-1}$$

Пользуясь (2.1.1), мы имеем

$$ {}_j\det({}_*{}^*)^i\, A^T = ({}^{\cdot}{}_{j}{}^{-}(A^{-1^*{}_*}){}_{\cdot}{}^{i}{}_{-})^{-1} \tag{2.3.17}$$

Пользуясь (2.3.17), (2.1.2), (2.3.12), мы получим (2.3.16). □

Теорема 2.3.8 расширяет принцип двойственности, теорема 2.2.14, на утверждения о квазидетерминантах и утверждает, что одно и тоже выражение является ${}_*{}^*$-квазидетерминантом матрицы $A$ и ${}^*{}_*$-квазидетерминантом матрицы $A^T$. Пользуясь этой теоремой, мы можем записать любое утверждение о ${}^*{}_*$-матрице, опираясь на подобное утверждение о ${}_*{}^*$-матрице.

Теорема 2.3.9 (принцип двойственности). *Пусть* $\mathfrak{A}$ *- истинное утверждение о бикольце матриц. Если мы одновременно заменим*

[2.6] Мы можем дать подобное доказательство для $\binom{j}{i}$-${}^*{}_*$**-квазидетерминанта**. Однако мы можем записать соответствующие утверждения, опираясь на принцип двойственности. Так, если прочесть равенство (2.3.14) справа налево, то мы получим равенство

$$ {}^j\det({}^*{}_*)A_i = {}^jA_i - {}^{[j]}A_i{}^*{}_* \left({}^{[j]}A_{[i]}\right)^{-1^*{}_*} {}^*{}_*{}^jA_{[i]}$$

$$ {}^j\det({}^*{}_*)A_i = {}^jA_i - {}^{[j]}A_i{}^*{}_*\mathcal{H}\det({}^*{}_*)^{[j]}A_{[i]}{}^*{}_*{}^jA_{[i]}$$

- *${}_*$-строку и ${}^*$-строку*
- *${}_*{}^*$-квазидетерминант и ${}^*{}_*$-квазидетерминант*

*то мы снова получим истинное утверждение.*

Теорема 2.3.10.

$$(mA)^{-1_*{}^*}=A^{-1_*{}^*}m^{-1}\tag{2.3.18}$$

$$(Am)^{-1_*{}^*}=m^{-1}A^{-1_*{}^*}\tag{2.3.19}$$

Доказательство. Мы докажем равенство (2.3.18) индукцией по размеру матрицы.

Для $1\times 1$ матрицы утверждение очевидно, так как

$$(mA)^{-1_*{}^*}=\left((mA)^{-1}\right)=\left(A^{-1}m^{-1}\right)=\left(A^{-1}\right)m^{-1}=A^{-1_*{}^*}m^{-1}$$

Допустим утверждение справедливо для $(n-1)\times(n-1)$ матрицы. Тогда из равенства (2.3.1) следует

$$\begin{aligned}({}^I((mA)^{-1_*{}^*})_J)^{-1_*{}^*}=&{}^J(mA)_I-{}^J(mA)_{[I]}{}_*{}^*\left({}^{[J]}(mA)_{[I]}\right)^{-1_*{}^*}{}_*{}^*{}^{[J]}(mA)_I\\=&m{}^JA_I-m\ {}^JA_{[I]}{}_*{}^*\left({}^{[J]}A_{[I]}\right)^{-1_*{}^*}m^{-1}{}_*{}^*m\ {}^{[J]}A_I\\=&m\ {}^JA_I-m\ {}^JA_{[I]}{}_*{}^*\left({}^{[J]}A_{[I]}\right)^{-1_*{}^*}{}_*{}^*{}^{[J]}A_I\end{aligned}$$

$$({}^I((mA)^{-1_*{}^*})_J)^{-1_*{}^*}=m\ {}^I(A^{-1_*{}^*})_J\tag{2.3.20}$$

Из равенства (2.3.20) следует равенство (2.3.18). Аналогично доказывается равенство (2.3.19). □

Теорема 2.3.11. *Пусть*

$$A=\begin{pmatrix}1&0\\0&1\end{pmatrix}\tag{2.3.21}$$

*Тогда*

$$A^{-1_*{}^*}=\begin{pmatrix}1&0\\0&1\end{pmatrix}\tag{2.3.22}$$

$$A^{-1^*{}_*}=\begin{pmatrix}1&0\\0&1\end{pmatrix}\tag{2.3.23}$$

Доказательство. Из (2.3.8) и (2.3.11) очевидно, что ${}^1(A^{-1_*{}^*})_1=1$ и ${}^2(A^{-1_*{}^*})_2=1$. Тем не менее выражение для ${}^2(A^{-1_*{}^*})_1$ и ${}^1(A^{-1_*{}^*})_2$ не может быть определено из (2.3.9) и (2.3.10) так как ${}^2A_1={}^1A_2=0$. Мы можем преобразовать эти выражения. Например

$$\begin{aligned}{}^2(A^{-1_*{}^*})_1&=({}^1A_2-{}^1A_1({}^2A_1)^{-1}\,{}^2A_2)^{-1}\\&=({}^1A_1(({}^1A_1)^{-1}\,{}^1A_2-({}^2A_1)^{-1}\,{}^2A_2))^{-1}\\&=(({}^2A_1)^{-1}\,{}^1A_1({}^2A_1({}^1A_1)^{-1}\,{}^1A_2-{}^2A_2))^{-1}\\&=({}^1A_1({}^2A_1({}^1A_1)^{-1}\,{}^1A_2-{}^2A_2))^{-1}\,{}^2A_1\end{aligned}$$

Мы непосредственно видим, что ${}^2(A^{-1_*{}^*})_1 = 0$. Таким же образом мы можем найти, что ${}^1(A^{-1_*{}^*})_2 = 0$. Это завершает доказательство (2.3.22).

Равенство (2.3.23) следует из (2.3.22), теоремы 2.3.8 и симметрии матрицы (2.3.21). □

Глава 3

# Представление универсальной алгебры

## 3.1. Представление универсальной алгебры

Определение 3.1.1. Пусть на множестве $M$ определена структура $\Omega_2$-алгебры ([2, 13]). Эндоморфизм $\Omega_2$-алгебры

$$t : M \to M$$

называется **преобразованием универсальной алгебры** $M$.[3.1] □

Мы будем обозначать $\delta$ тождественное преобразование.

Определение 3.1.2. Пусть ${}^*M$ - множество **левосторонних преобразований**

$$u' = tu$$

$\Omega_2$-алгебры $M$. Пусть на множестве ${}^*M$ определена структура $\Omega_1$-алгебры. Гомоморфизм

$$f : A \to {}^*M \tag{3.1.1}$$

$\Omega_1$-алгебры $A$ в $\Omega_1$-алгебру ${}^*M$ называется **левосторонним представлением $\Omega_1$-алгебры** $A$ или $A*$**-представлением** в $\Omega_2$-алгебре $M$. □

Определение 3.1.3. Пусть $M^*$ - множество **правосторонних преобразований**

$$u' = ut$$

$\Omega_2$-алгебры $M$. Пусть на множестве $M^*$ определена структура $\Omega_1$-алгебры. Гомоморфизм

$$f : A \to M^*$$

$\Omega_1$-алгебры $A$ в $\Omega_1$-алгебру $M^*$ называется **правосторонним представлением $\Omega_1$-алгебры** $A$ или $*A$**-представлением** в $\Omega_2$-алгебре $M$. □

Мы распространим на теорию представлений соглашение, описанное в замечании 2.2.15. Мы можем записать принцип двойственности в следующей форме

Теорема 3.1.4 (принцип двойственности). *Любое утверждение, справедливое для левостороннего представления $\Omega_1$-алгебры $A$, будет справедливо для правостороннего представления $\Omega_1$-алгебры $A$.*

---

[3.1] Если множество операций $\Omega_2$-алгебры пусто, то

$$t : M \to M$$

является отображением.

Замечание 3.1.5. Существует две формы записи преобразования $\Omega_2$-алгебры $M$. Если мы пользуемся операторной записью, то преобразование $A$ записывается в виде $Aa$ или $aA$, что соответствует левостороннему преобразованию или правостороннему преобразованию. Если мы пользуемся функциональной записью, то преобразование $A$ записывается в виде $A(a)$ независимо от того, это левостороннее или правостороннее преобразование. Эта запись согласована с принципом двойственности.

Это замечание является основой следующего соглашения. Когда мы пользуемся функциональной записью, мы не различаем левостороннее и правостороннее преобразование. Мы будем обозначать ${}^*M$ множество преобразований $\Omega_2$-алгебры $M$. Пусть на множестве ${}^*M$ определена структура $\Omega_1$-алгебры. Пусть $A$ является $\Omega_1$-алгеброй. Мы будем называть гомоморфизм

$$f : A \to {}^*M \tag{3.1.2}$$

**представлением $\Omega_1$-алгебры $A$ в $\Omega_2$-алгебре $M$**. Мы будем также пользоваться записью

$$f : A \xrightarrow{*} M$$

для обозначения представления $\Omega_1$-алгебры $A$ в $\Omega_2$-алгебре $M$.

Соответствие между операторной записью и функциональной записью однозначно. Мы можем выбирать любую форму записи, которая удобна для изложения конкретной темы. □

Существует несколько способов описать представление. Мы можем указать отображение $f$, имея в виду что область определения - это $\Omega_1$-алгебра $A$ и область значений - это $\Omega_1$-алгебра ${}^*M$. Либо мы можем указать $\Omega_1$-алгебру $A$ и $\Omega_2$-алгебру $M$, имея в виду что нам известна структура отображения $f$.[3.2]

Диаграмма

$$\begin{array}{ccc} M & \xrightarrow{f(a)} & M \\ & \Uparrow{\scriptstyle f} & \\ & A & \end{array}$$

означает, что мы рассматриваем представление $\Omega_1$-алгебры $A$. Отображение $f(a)$ является образом $a \in A$.

Определение 3.1.6. Мы будем называть представление $\Omega_1$-алгебры $A$ **эффективным**, если отображение (3.1.2) - изоморфизм $\Omega_1$-алгебры $A$ в ${}^*M$. □

Замечание 3.1.7. Если левостороннее представление $\Omega_1$-алгебры эффективно, мы можем отождествлять элемент $\Omega_1$-алгебры с его образом и записывать левостороннее преобразование, порождённое элементом $a \in A$, в форме

$$v' = av$$

Если правостороннее представление $\Omega_1$-алгебры эффективно, мы можем отождествлять элемент $\Omega_1$-алгебры с его образом и записывать правостороннее преобразование, порождённое элементом $a \in A$, в форме

$$v' = va$$

□

[3.2]Например, мы рассматриваем векторное пространство $\overline{V}$ над полем $D$ ( definition 5.1.4 ).

Определение 3.1.8. Мы будем называть представление $\Omega_1$-алгебры **транзитивным**, если для любых $a, b \in V$ существует такое $g$, что

$$a = f(g)(b)$$

Мы будем называть представление $\Omega_1$-алгебры **однотранзитивным**, если оно транзитивно и эффективно. □

Теорема 3.1.9. *Представление однотранзитивно тогда и только тогда, когда для любых $a, b \in M$ существует одно и только одно $g \in A$ такое, что $a = f(g)(b)$*

Доказательство. Следствие определений 3.1.6 и 3.1.8. □

## 3.2. Морфизм представлений универсальной алгебры

Теорема 3.2.1. *Пусть $A$ и $B$ - $\Omega_1$-алгебры. Представление $\Omega_1$-алгебры $B$*

$$g : B \to {}^\star M$$

*и гомоморфизм $\Omega_1$-алгебры*

$$h : A \to B \tag{3.2.1}$$

*определяют представление $f$ $\Omega_1$-алгебры $A$*

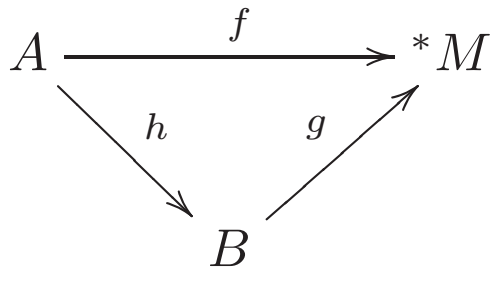

Доказательство. Отображение $f$ является гомоморфизмом $\Omega_1$-алгебры $A$ в $\Omega_1$-алгебру ${}^*M$, так как отображение $g$ является гомоморфизмом $\Omega_1$-алгебры $B$ в $\Omega_1$-алгебру ${}^*M$. □

Если мы изучаем представление $\Omega_1$-алгебры в $\Omega_2$-алгебрах $M$ и $N$, то нас интересуют отображения из $M$ в $N$, сохраняющие структуру представления.

Определение 3.2.2. Пусть

$$f : A \to {}^*M$$

представление $\Omega_1$-алгебры $A$ в $\Omega_2$-алгебре $M$ и

$$g : B \to {}^*N$$

представление $\Omega_1$-алгебры $B$ в $\Omega_2$-алгебре $N$. Пара отображений

$$(r : A \to B, R : M \to N) \tag{3.2.2}$$

таких, что

- $r$ - гомоморфизм $\Omega_1$-алгебры
- $R$ - гомоморфизм $\Omega_2$-алгебры
-

$$R \circ f(a) = g(r(a)) \circ R \tag{3.2.3}$$

называется **морфизмом представлений из $f$ в $g$**. Мы также будем говорить, что определён **морфизм представлений $\Omega_1$-алгебры в $\Omega_2$-алгебре**. □

Замечание 3.2.3. Мы можем рассматривать пару отображений $r$, $R$ как отображение

$$F : A \cup M \to B \cup N$$

такое, что

$$F(A) = B \qquad F(M) = N$$

Поэтому в дальнейшем мы будем говорить, что дано отображение $(r, R)$. □

Определение 3.2.4. Если представления $f$ и $g$ совпадают, то морфизм представлений $(r, R)$ называется **морфизмом представления** $f$. □

Для произвольного $m \in M$ равенство (3.2.3) имеет вид

$$R(f(a)(m)) = g(r(a))(R(m)) \tag{3.2.4}$$

Замечание 3.2.5. Рассмотрим морфизм представлений (3.2.2). Мы можем обозначать элементы множества $B$, пользуясь буквой по образцу $b \in B$. Но если мы хотим показать, что $b$ является образом элемента $a \in A$, мы будем пользоваться обозначением $r(a)$. Таким образом, равенство

$$r(a) = r(a)$$

означает, что $r(a)$ (в левой части равенства) является образом $a \in A$ (в правой части равенства). Пользуясь подобными соображениями, мы будем обозначать элемент множества $N$ в виде $R(m)$. Мы будем следовать этому соглашению, изучая соотношения между гомоморфизмами $\Omega_1$-алгебр и отображениями между множествами, где определены соответствующие представления. □

Замечание 3.2.6. Мы можем интерпретировать (3.2.4) двумя способами

- Пусть преобразование $f(a)$ отображает $m \in M$ в $f(a)(m)$. Тогда преобразование $g(r(a))$ отображает $R(m) \in N$ в $R(f(a)(m))$.
- Мы можем представить морфизм представлений из $f$ в $g$, пользуясь диаграммой

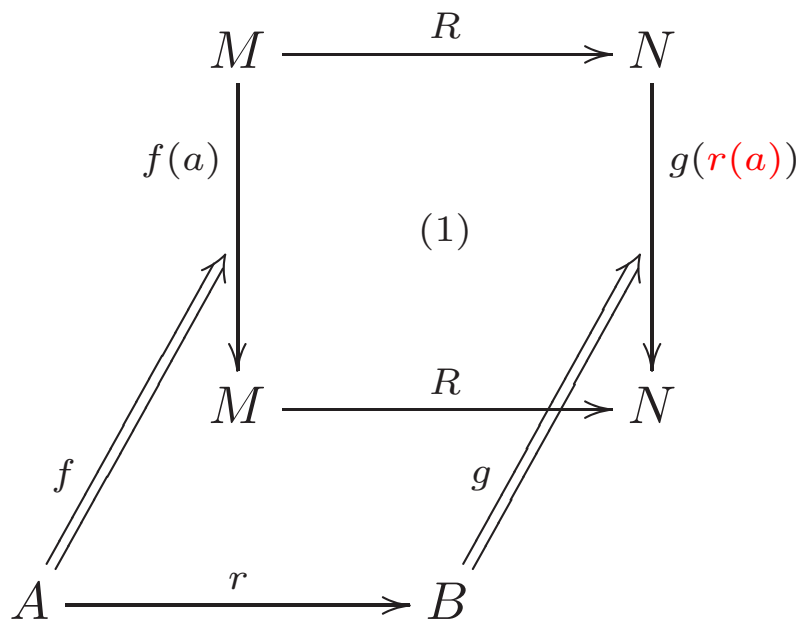

Из (3.2.3) следует, что диаграмма (1) коммутативна.

□

Теорема 3.2.7. *Рассмотрим представление*

$$f : A \to {}^*M$$

$\Omega_1$*-алгебры* $A$ *и представление*

$$g : B \to {}^*N$$

$\Omega_1$-*алгебры* $B$. *Морфизм*

$$h:A\longrightarrow B\qquad H:M\longrightarrow N$$

*представлений из* $f$ *в* $g$ *удовлетворяет соотношению*

$$H\circ\omega(f(a_1),...,f(a_n))=\omega(g(h(a_1)),...,g(h(a_n)))\circ H\tag{3.2.5}$$

*для произвольной* $n$*-арной операции* $\omega$ $\Omega_1$*-алгебры.*

Доказательство. Так как $f$ - гомоморфизм, мы имеем

$$H\circ\omega(f(a_1),...,f(a_n))=H\circ f(\omega(a_1,...,a_n))\tag{3.2.6}$$

Из (3.2.3) и (3.2.6) следует

$$H\circ\omega(f(a_1),...,f(a_n))=g(h(\omega(a_1,...,a_n)))\circ H\tag{3.2.7}$$

Так как $h$ - гомоморфизм, из (3.2.7) следует

$$H\circ\omega(f(a_1),...,f(a_n))=g(\omega(h(a_1),...,h(a_n)))\circ H\tag{3.2.8}$$

Так как $g$ - гомоморфизм, из (3.2.8) следует (3.2.5). □

Теорема 3.2.8. *Пусть отображение*

$$h:A\longrightarrow B\qquad H:M\longrightarrow N$$

*является морфизмом из представления*

$$f:A\to{}^*M$$

$\Omega_1$*-алгебры* $A$ *в представление*

$$g:B\to{}^*N$$

$\Omega_1$*-алгебры* $B$. *Если представление* $f$ *эффективно, то отображение*

$${}^*H:{}^*M\to{}^*N$$

*определённое равенством*

$${}^*H(f(a))=g(h(a))\tag{3.2.9}$$

*является гомоморфизмом* $\Omega_1$*-алгебры.*

Доказательство. Так как представление $f$ эффективно, то для выбранного преобразования $f(a)$ выбор элемента $a$ определён однозначно. Следовательно, преобразование $g(h(a))$ в равенстве (3.2.9) определено корректно.

Так как $f$ - гомоморфизм, мы имеем

$${}^*H(\omega(f(a_1),...,f(a_n)))={}^*H(f(\omega(a_1,...,a_n)))\tag{3.2.10}$$

Из (3.2.9) и (3.2.10) следует

$${}^*H(\omega(f(a_1),...,f(a_n)))=g(h(\omega(a_1,...,a_n)))\tag{3.2.11}$$

Так как $h$ - гомоморфизм, из (3.2.11) следует

$${}^*H(\omega(f(a_1),...,f(a_n)))=g(\omega(h(a_1),...,h(a_n)))\tag{3.2.12}$$

Так как $g$ - гомоморфизм,

$${}^*H(\omega(f(a_1),...,f(a_n)))=\omega(g(h(a_1)),...,g(h(a_n)))=\omega({}^*H(f(a_1)),...,{}^*H(f(a_n)))$$

следует из (3.2.12). Следовательно, отображение ${}^*H$ является гомоморфизмом $\Omega_1$-алгебры. □

Теорема 3.2.9. *Если представление*

$$f : A \to {}^*M$$

*$\Omega_1$-алгебры $A$ однотранзитивно и представление*

$$g : B \to {}^*N$$

*$\Omega_1$-алгебры $B$ однотранзитивно, то существует морфизм*

$$h : A \to B \quad H : M \to N$$

*представлений из $f$ в $g$.*

Доказательство. Выберем гомоморфизм $h$. Выберем элемент $m \in M$ и элемент $n \in N$. Чтобы построить отображение $H$, рассмотрим следующую диаграмму

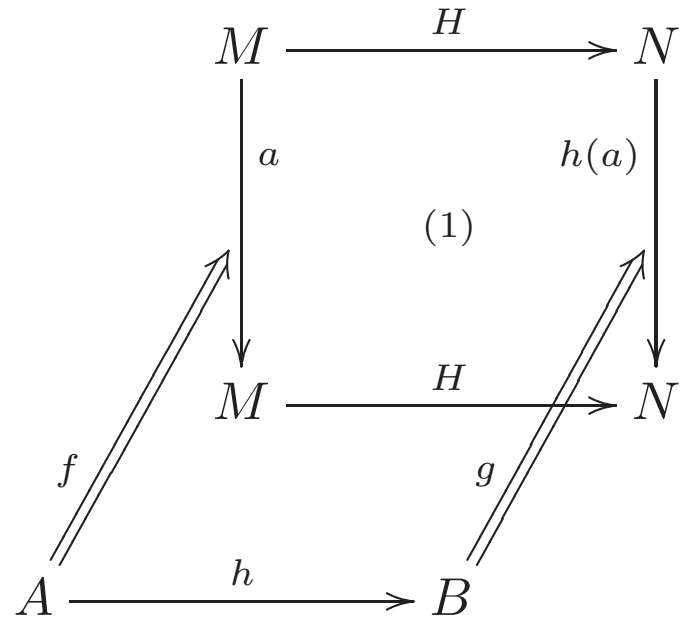

Из коммутативности диаграммы (1) следует

$$H(am) = h(a)H(m)$$

Для произвольного $m' \in M$ однозначно определён $a \in A$ такой, что $m' = am$. Следовательно, мы построили отображении $H$, которое удовлетворяет равенству (3.2.3). □

Теорема 3.2.10. *Если представление*

$$f : A \to {}^*M$$

*$\Omega_1$-алгебры $A$ однотранзитивно и представление*

$$g : B \to {}^*N$$

*$\Omega_1$-алгебры $B$ однотранзитивно, то для заданного гомоморфизма $\Omega_1$-алгебры*

$$h : A \longrightarrow B$$

*гомоморфизм $\Omega_2$-алгебры*

$$H : M \longrightarrow N$$

*такой, что $(h, H)$ является морфизмом представлений из $f$ в $g$, определён однозначно с точностью до выбора образа $n = H(m) \in N$ заданного элемента $m \in M$.*

Доказательство. Из доказательства теоремы 3.2.9 следует, что выбор гомоморфизма $h$ и элементов $m \in M$, $n \in N$ однозначно определяет отображение $H$. □

Теорема 3.2.11. *Если представление*

$$f : A \to {}^*M$$

*$\Omega_1$-алгебры $A$ однотранзитивно, то для любого эндоморфизма $h$ $\Omega_1$-алгебры $A$ существует морфизм представления $f$*

$$h : A \to B \quad H : M \to N$$

Доказательство. Рассмотрим следующую диаграмму

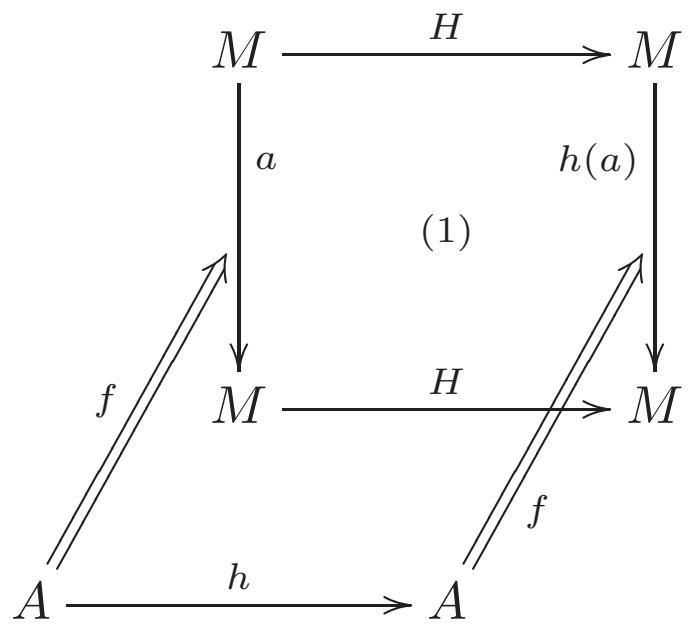

Утверждение теоремы является следствием теоремы 3.2.9. □

Теорема 3.2.12. *Пусть*

$$f : A \to {}^*M$$

*представление $\Omega_1$-алгебры $A$,*

$$g : B \to {}^*N$$

*представление $\Omega_1$-алгебры $B$,*

$$h : C \to {}^*L$$

*представление $\Omega_1$-алгебры $C$. Пусть определены морфизмы представлений $\Omega_1$-алгебры*

$$p : A \longrightarrow B \qquad P : M \longrightarrow N$$

$$q : B \longrightarrow C \qquad Q : N \longrightarrow L$$

*Тогда определён морфизм представлений $\Omega_1$-алгебры*

$$r : A \longrightarrow C \qquad R : M \longrightarrow L$$

*где $r = qp$, $R = QP$. Мы будем называть морфизм $(r, R)$ представлений из $f$ в $h$* **произведением морфизмов** $(p, P)$ **и** $(q, Q)$ **представлений универсальной алгебры**.

Доказательство. Мы можем представить утверждение теоремы, пользуясь диаграммой

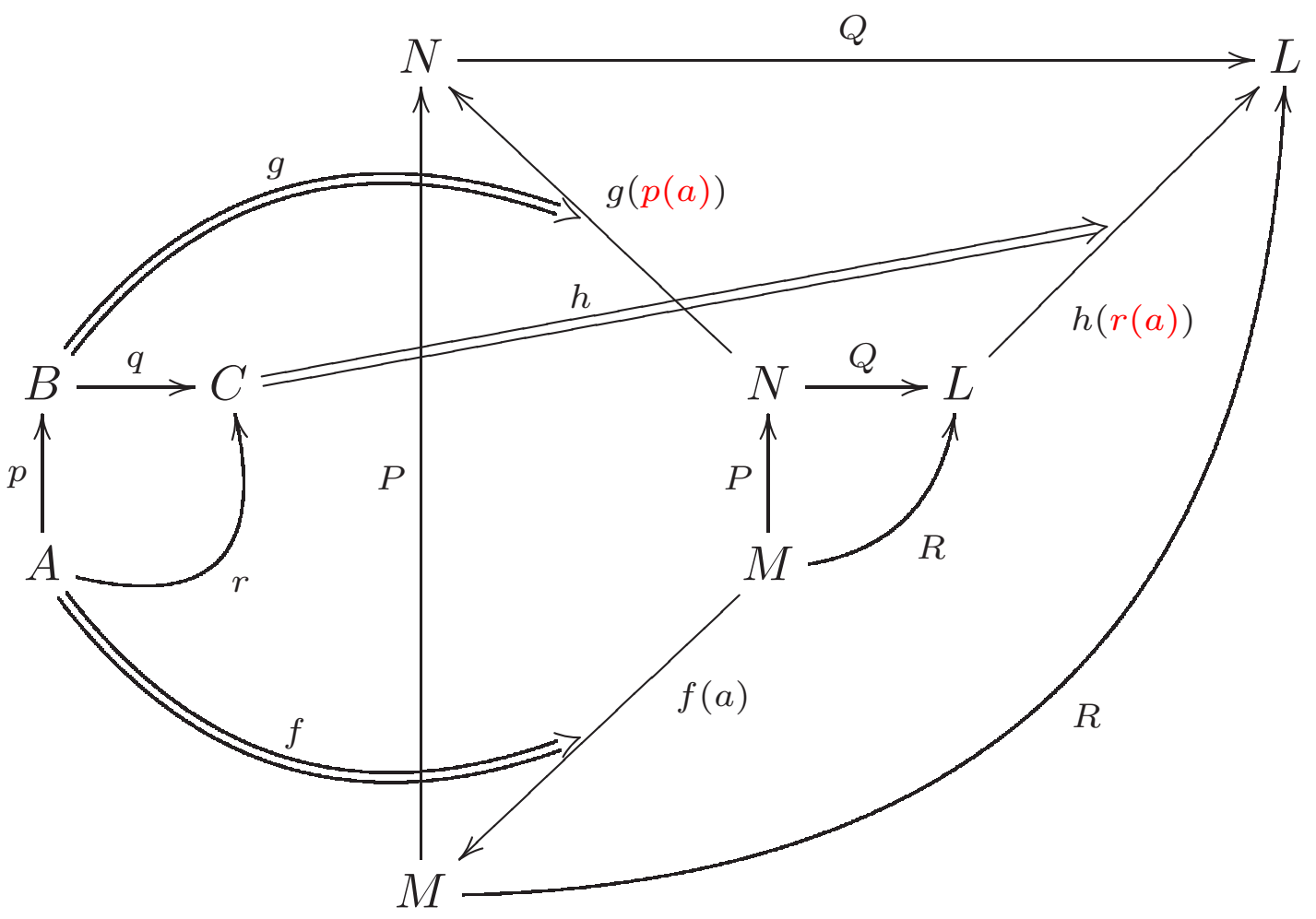

Отображение $r$ является гомоморфизмом $\Omega_1$-алгебры $A$ в $\Omega_1$-алгебру $C$. Нам надо показать, что пара отображений $(r, R)$ удовлетворяет (3.2.3):

$$\begin{aligned}R(f(a)m) &= QP(f(a)m)\\ &= Q(g(p(a))P(m))\\ &= h(qp(a))QP(m))\\ &= h(r(a))R(m)\end{aligned}$$

□

Определение 3.2.13. Допустим $\mathcal{A}$ категория $\Omega_1$-алгебр. Мы определим **категорию $\mathcal{A}*$ левосторонних представлений универсальной алгебры из категории** $\mathcal{A}$. Объектами этой категории являются левосторонние представлениями $\Omega_1$-алгебры. Морфизмами этой категории являются морфизмы левосторонних представлений $\Omega_1$-алгебры. □

Определение 3.2.14. Пусть на множестве $M$ определена эквивалентность $S$. Преобразование $f$ называется **согласованным с эквивалентностью** $S$, если из условия $m_1 \equiv m_2(\mathrm{mod}S)$ следует $f(m_1) \equiv f(m_2)(\mathrm{mod}\ S)$. □

Теорема 3.2.15. *Пусть на множестве $M$ определена эквивалентность $S$. Пусть на множестве ${}^*M$ определена $\Omega_1$-алгебра. Если любое преобразование $f \in {}^*M$ согласованно с эквивалентностью $S$, то мы можем определить структуру $\Omega_1$-алгебры на множестве ${}^*(M/S)$.*

Доказательство. Пусть $h = \mathrm{nat}\ S$. Если $m_1 \equiv m_2(\mathrm{mod}S)$ , то $h(m_1) = h(m_2)$. Поскольку $f \in {}^*M$ согласованно с эквивалентностью $S$, то $h(f(m_1)) = h(f(m_2))$. Это позволяет определить преобразование $F$ согласно правилу

$$F([m]) = h(f(m)) \tag{3.2.13}$$

Пусть $\omega$ - n-арная операция $\Omega_1$-алгебры. Пусть $f_1$, ..., $f_n \in {}^\star M$ и

$$F_1([m]) = h(f_1(m)) \quad ... \quad F_n([m]) = h(f_n(m))$$

Согласно условию теоремы, преобразование

$$f = \omega(f_1, ..., f_n) \in {}^*M$$

согласованно с эквивалентностью $S$. Следовательно, из условия $m_1 \equiv m_2 (\mathrm{mod} S)$ и определения 3.2.14 следует

$$\begin{gathered} f(m_1) \equiv f(m_2)(\mathrm{mod} S) \\ \omega(f_1, ..., f_n)(m_1) \equiv \omega(f_1, ..., f_n)(m_2)(\mathrm{mod} S) \end{gathered} \tag{3.2.14}$$

Следовательно, мы можем определить операцию $\omega$ на множестве ${}^\star(M/S)$ по правилу

$$\omega(F_1, ..., F_n)[m] = h(\omega(f_1, ..., f_n)(m)) \tag{3.2.15}$$

Из определения (3.2.13) и равенства (3.2.14) следует, что мы корректно определили операцию $\omega$ на множестве ${}^\star(M/S)$. □

Теорема 3.2.16. *Пусть*

$$f : A \to {}^*M$$

*представление $\Omega_1$-алгебры $A$,*

$$g : B \to {}^*N$$

*представление $\Omega_1$-алгебры $B$. Пусть*

$$r : A \longrightarrow B \qquad R : M \longrightarrow N$$

*морфизм представлений из $f$ в $g$. Положим*

$$s = rr^{-1} \qquad S = RR^{-1}$$

*Тогда для отображений $r$, $R$ существуют разложения, которые можно описать диаграммой*

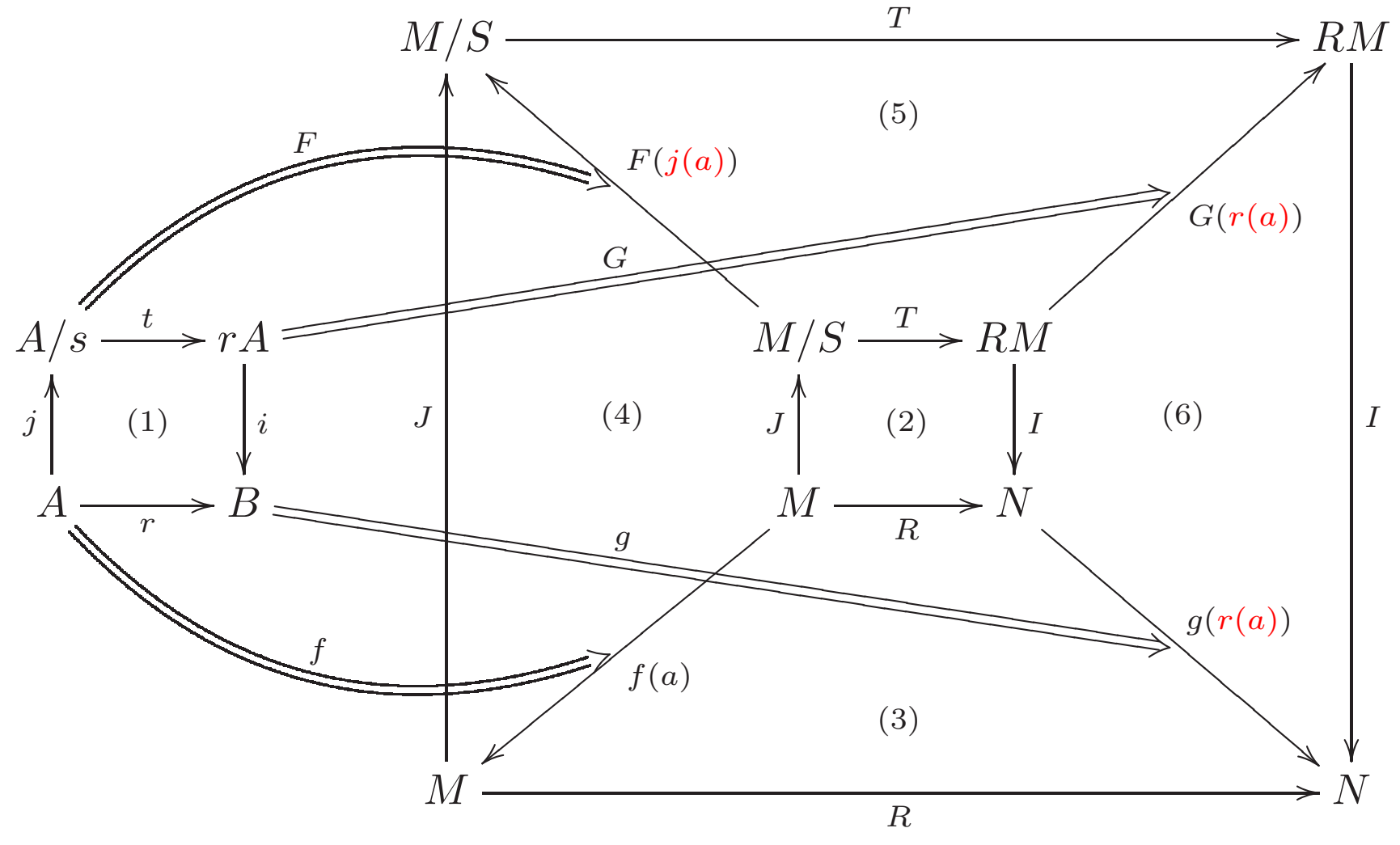

(1) *$s = \ker r$ является конгруэнцией на $A$. Существует разложение гомоморфизма $r$*

$$r = itj \tag{3.2.16}$$

$j = \mathrm{nat}\ s$ - *естественный гомоморфизм*

$$j(a) = j(a)$$

$t$ - *изоморфизм*

$$r(a) = t(j(a)) \tag{3.2.17}$$

$i$ - *вложение*

$$r(a) = i(r(a)) \tag{3.2.18}$$

(2) $S = \ker R$ *является эквивалентностью на* $M$. *Существует разложение отображения* $R$

$$R = ITJ \tag{3.2.19}$$

$J = \mathrm{nat}\ S$ - *сюръекция*

$$J(m) = J(m)$$

$T$ - *биекция*

$$R(m) = T(J(m)) \tag{3.2.20}$$

$I$ - *вложение*

$$R(m) = I(R(m)) \tag{3.2.21}$$

(3) $F$ - *левостороннее представление* $\Omega_1$*-алгебры* $A/s$ *в* $M/S$
(4) $G$ - *левостороннее представление* $\Omega_1$*-алгебры* $rA$ *в* $RM$
(5) $(j, J)$ - *морфизм представлений* $f$ *и* $F$
(6) $(t, T)$ - *морфизм представлений* $F$ *и* $G$
(7) $(t^{-1}, T^{-1})$ - *морфизм представлений* $G$ *и* $F$
(8) $(i, I)$ - *морфизм представлений* $G$ *и* $g$
(9) *Существует разложение морфизма представлений*

$$(r, R) = (i, I)(t, T)(j, J) \tag{3.2.22}$$

Доказательство. Существование диаграмм (1) и (2) следует из теоремы II.3.7 ([13], с. 74).

Мы начнём с диаграммы (4).

Пусть $m_1 \equiv m_2 (\mathrm{mod}\ S)$. Следовательно,

$$R(m_1) = R(m_2) \tag{3.2.23}$$

Если $a_1 \equiv a_2 (\mathrm{mod}\ s)$, то

$$r(a_1) = r(a_2) \tag{3.2.24}$$

Следовательно, $j(a_1) = j(a_2)$. Так как $(r, R)$ - морфизм представлений, то

$$R(f(a_1)(m_1)) = g(r(a_1))(R(m_1)) \tag{3.2.25}$$

$$R(f(a_2)(m_2)) = g(r(a_2))(R(m_2)) \tag{3.2.26}$$

Из (3.2.23), (3.2.24), (3.2.25), (3.2.26) следует

$$R(f(a_1)(m_1)) = R(f(a_2)(m_2)) \tag{3.2.27}$$

Из (3.2.27) следует

$$f(a_1)(m_1) \equiv f(a_2)(m_2) (\mathrm{mod}\ S) \tag{3.2.28}$$

и, следовательно,

(3.2.29) $$J(f(a_1)(m_1)) = J(f(a_2)(m_2))$$

Из (3.2.29) следует, что отображение

(3.2.30) $$F(j(a))(J(m)) = J(f(a)(m)))$$

определено корректно и является преобразованием множества $M/S$.

Из равенства (3.2.28) (в случае $a_1 = a_2$) следует, что для любого $a$ преобразование согласованно с эквивалентностью $S$. Из теоремы 3.2.15 следует, что на множестве ${}^*(M/S)$ определена структура $\Omega_1$-алгебры. Рассмотрим $n$-арную операцию $\omega$ и $n$ преобразований

$$F(j(a_i))(J(m)) = J(f(a_i)(m))) \quad i = 1, ..., n$$

пространства $M/S$. Мы положим

$$\omega(F(j(a_1)), ..., F(j(a_n)))(J(m)) = J(\omega(f(a_1), ..., f(a_n)))(m))$$

Следовательно, отображение $F$ является представлением $\Omega_1$-алгебры $A/s$.

Из (3.2.30) следует, что $(j, J)$ является морфизмом представлений $f$ и $F$ (утверждение (5) теоремы).

Рассмотрим диаграмму (5).

Так как $T$ - биекция, то мы можем отождествить элементы множества $M/S$ и множества $MR$, причём это отождествление имеет вид

(3.2.31) $$T(J(m)) = R(m)$$

Мы можем записать преобразование $F(j(a))$ множества $M/S$ в виде

(3.2.32) $$F(j(a)) : J(m) \to F(j(a))(J(m))$$

Так как $T$ - биекция, то мы можем определить преобразование

(3.2.33) $$T(J(m)) \to T(F(j(a))(J(m)))$$

множества $RM$. Преобразование (3.2.33) зависит от $j(a) \in A/s$. Так как $t$ - биекция, то мы можем отождествить элементы множества $A/s$ и множества $rA$, причём это отождествление имеет вид

(3.2.34) $$t(j(a)) = r(a)$$

Следовательно, мы определили отображение

$$G : rA \to {}^\star RM$$

согласно равенству

(3.2.35) $$G(t(j(a)))(T(J(m))) = T(F(j(a))(J(m)))$$

Рассмотрим $n$-арную операцию $\omega$ и $n$ преобразований

$$G(r(a_i))(R(m)) = T(F(j(a_i))(J(m))) \quad i = 1, ..., n$$

пространства $RM$. Мы положим

(3.2.36) $$\omega(G(r(a_1)), ..., G(r(a_n)))(R(m)) = T(\omega(F(j(a_1), ..., F(j(a_n)))(J(m)))$$

Согласно (3.2.35) операция $\omega$ корректно определена на множестве ${}^\star RM$. Следовательно, отображение $G$ является представлением $\Omega_1$-алгебры.

Из (3.2.35) следует, что $(t, T)$ является морфизмом представлений $F$ и $G$ (утверждение (6) теоремы).

Так как $T$ - биекция, то из равенства (3.2.31) следует

$$J(m) = T^{-1}(R(m)) \tag{3.2.37}$$

Мы можем записать преобразование $G(r(a))$ множества $RM$ в виде

$$G(r(a)) : R(m) \to G(r(a))(R(m)) \tag{3.2.38}$$

Так как $T$ - биекция, то мы можем определить преобразование

$$T^{-1}(R(m)) \to T^{-1}(G(r(a))(R(m))) \tag{3.2.39}$$

множества $M/S$. Преобразование (3.2.39) зависит от $r(a) \in rA$. Так как $t$ - биекция, то из равенства (3.2.34) следует

$$j(a) = t^{-1}(r(a)) \tag{3.2.40}$$

Так как по построению диаграмма (5) коммутативна, то преобразование (3.2.39) совпадает с преобразованием (3.2.32). Равенство (3.2.36) можно записать в виде

$$T^{-1}(\omega(G(r(a_1)), ..., G(r(a_n)))(R(m))) = \omega(F(j(a_1), ..., F(j(a_n)))(J(m)) \tag{3.2.41}$$

Следовательно, $(t^{-1}, T^{-1})$ является морфизмом представлений $G$ и $F$ (утверждение (7) теоремы).

Диаграмма (6) является самым простым случаем в нашем доказательстве. Поскольку отображение $I$ является вложением и диаграмма (2) коммутативна, мы можем отождествить $n \in N$ и $R(m)$, если $n \in \mathrm{Im}R$. Аналогично, мы можем отождествить соответствующие преобразования.

$$\begin{gathered} g'(i(r(a)))(I(R(m))) = I(G(r(a))(R(m))) \\ \omega(g'(r(a_1)), ..., g'(r(a_n)))(R(m)) = I(\omega(G(r(a_1), ..., G(r(a_n)))(R(m))) \end{gathered} \tag{3.2.42}$$

Следовательно, $(i, I)$ является морфизмом представлений $G$ и $g$ (утверждение (8) теоремы).

Для доказательства утверждения (9) теоремы осталось показать, что определённое в процессе доказательства представление $g'$ совпадает с представлением $g$, а операции над преобразованиями совпадают с соответствующими операциями на ${}^*N$.

$$\begin{aligned} g'(i(r(a)))(I(R(m))) &= I(G(r(a))(R(m))) && \text{by (3.2.42)} \\ &= I(G(t(j(a)))(T(J(m)))) && \text{by (3.2.17), (3.2.20)} \\ &= IT(F(j(a))(J(m))) && \text{by (3.2.35)} \\ &= ITJ(f(a)(m)) && \text{by (3.2.30)} \\ &= R(f(a)(m)) && \text{by (3.2.19)} \\ &= g(r(a))(R(m)) && \text{by (3.2.3)} \end{aligned}$$

$$\begin{aligned} \omega(G(r(a_1)), ..., G(r(a_n)))(R(m)) &= T(\omega(F(j(a_1), ..., F(j(a_n)))(J(m))) \\ &= T(F(\omega(j(a_1), ..., j(a_n)))(J(m))) \\ &= T(F(j(\omega(a_1, ..., a_n)))(J(m))) \\ &= T(J(f(\omega(a_1, ..., a_n))(m))) \end{aligned}$$

□

Определение 3.2.17. Пусть

$$f : A \to {}^*M$$

представление $\Omega_1$-алгебры $A$,

$$g : B \to {}^*N$$

представление $\Omega_1$-алгебры $B$. Пусть

$$r : A \longrightarrow B \qquad R : M \longrightarrow N$$

морфизм представлений из $r$ в $R$ такой, что $f$ - изоморфизм $\Omega_1$-алгебры и $g$ - изоморфизм $\Omega_2$-алгебры. Тогда отображение $(r, R)$ называется **изоморфизмом представлений**. □

Теорема 3.2.18. *В разложении* (3.2.22) *отображение* $(t, T)$ *является изоморфизмом представлений* $F$ *и* $G$.

Доказательство. Следствие определения 3.2.17 и утверждений (6) и (7) теоремы 3.2.16. □

Из теоремы 3.2.16 следует, что мы можем свести задачу изучения морфизма представлений $\Omega_1$-алгебры к случаю, описываемому диаграммой

(3.2.43)

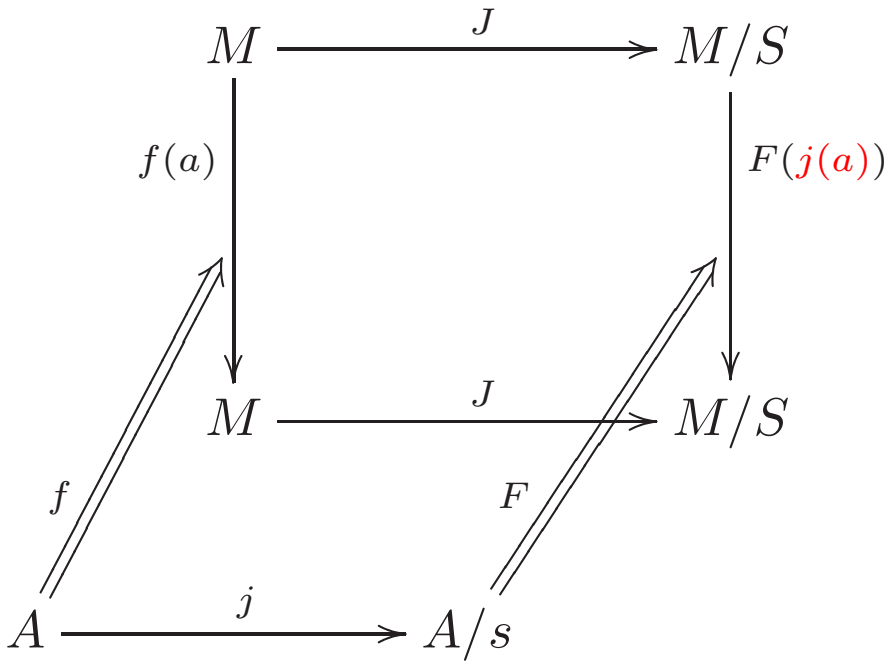

Теорема 3.2.19. *Диаграмма* (3.2.43) *может быть дополнена представлением* $F_1$ $\Omega_1$*-алгебры* $A$ *в множестве* $M/S$ *так, что диаграмма*

(3.2.44)

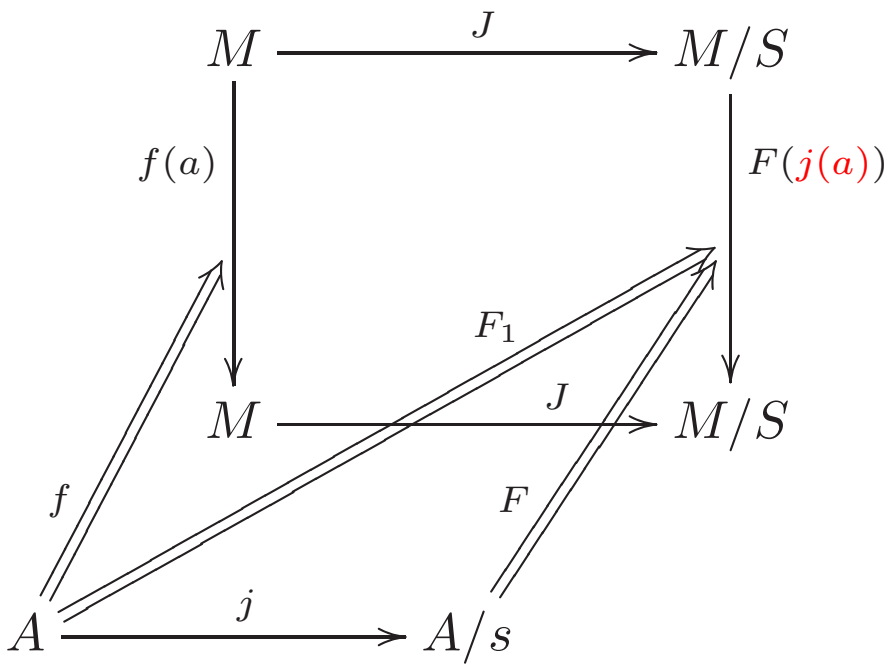

*коммутативна. При этом множество преобразований представления* $F$ *и множество преобразований представления* $F_1$ *совпадают.*

Доказательство. Для доказательства теоремы достаточно положить

$$F_1(a) = F(j(a))$$

Так как отображение $j$ - сюрьекция, то $\mathrm{Im}F_1 = \mathrm{Im}F$. Так как $j$ и $F$ - гомоморфизмы $\Omega_1$-алгебры, то $F_1$ - также гомоморфизм $\Omega_1$-алгебры. □

Теорема 3.2.19 завершает цикл теорем, посвящённых структуре морфизма представлений $\Omega_1$-алгебры. Из этих теорем следует, что мы можем упростить задачу изучения морфизма представлений $\Omega_1$-алгебры и ограничиться морфизмом представлений вида

$$id : A \longrightarrow A \qquad R : M \longrightarrow N$$

В этом случае мы можем отождествить морфизм $(\mathrm{id}, R)$ представлений $\Omega_1$-алгебры и соответствующий гомоморфизм $R$ $\Omega_2$-алгебры и пользоваться одной и тойже буквой $R$ для обозначения этих отображений. Мы будем пользоваться диаграммой

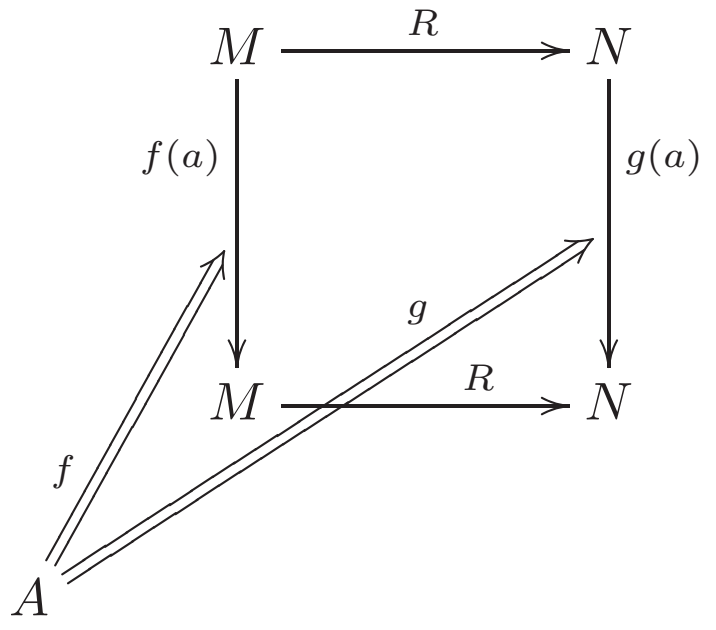

для представления морфизма $(\mathrm{id}, R)$ представлений $\Omega_1$-алгебры. Из диаграммы следует

$$R \circ f(a) = g(a) \circ R \tag{3.2.45}$$

Мы дадим следующее определение по аналогии с определением 3.2.13.

Определение 3.2.20. Мы определим **категорию** $A*$ **левосторонних представлений** $\Omega_1$**-алгебры** $A$. Объектами этой категории являются левосторонними представлениями $\Omega_1$-алгебры $A$. Морфизмами этой категории являются морфизмы $(id, R)$ левосторонних представлений $\Omega_1$-алгебры $A$. □

## 3.3. Автоморфизм представления универсальной алгебры

Определение 3.3.1. Пусть

$$f : A \to {}^*M$$

представление $\Omega_1$-алгебры $A$ в $\Omega_2$-алгебре $M$. Морфизм представлений $\Omega_1$-алгебры

$$(\mathrm{id} : A \to A, R : M \to M)$$

такой, что $R$ - эндоморфизм $\Omega_2$-алгебры называется **эндоморфизмом представления** $f$. □

Теорема 3.3.2. *Если представление*

$$f : A \to {}^*M$$

*$\Omega_1$-алгебры $A$ однотранзитивно, то для любых $p$, $q \in M$ существует единственный эндоморфизм*

$$H : M \to M$$

*представления $f$ такой, что $H(p) = q$.*

Доказательство. Рассмотрим следующую диаграмму

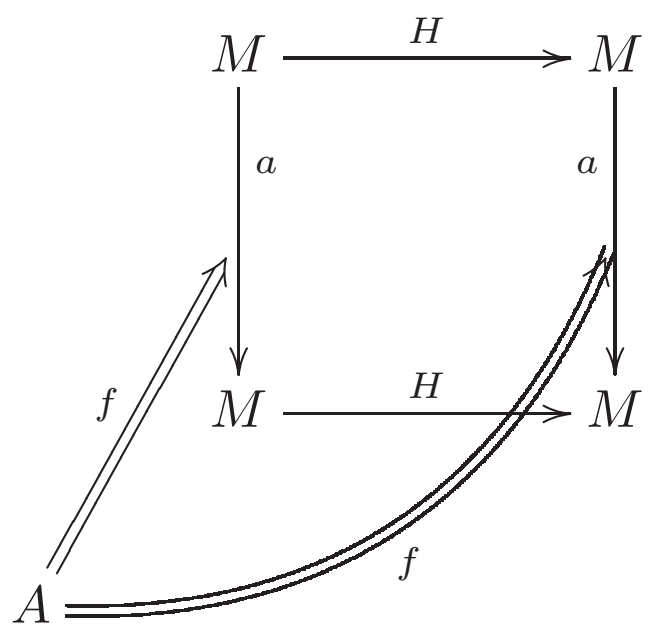

Существование эндоморфизма является следствием теоремы 3.2.9. Единственность эндоморфизма для заданных $p$, $q \in M$ является следствием теоремы 3.2.10, когда $h = \mathrm{id}$. □

Теорема 3.3.3. *Эндоморфизмы представления $f$ порождают полугруппу.*

Доказательство. Из теоремы 3.2.12 следует, что произведение эндоморфизмов $(p, P)$, $(r, R)$ представления $f$ является эндоморфизмом $(pr, PR)$ представления $f$. □

Определение 3.3.4. Пусть

$$f : A \to {}^*M$$

представление $\Omega_1$-алгебры $A$ в $\Omega_2$-алгебре $M$. Морфизм представлений $\Omega_1$-алгебры

$$(\mathrm{id} : A \to A, R : M \to M)$$

такой, что $R$ - автоморфизм $\Omega_2$-алгебры называется **автоморфизмом представления** $f$. □

Теорема 3.3.5. *Пусть*

$$f : A \to {}^*M$$

*представление $\Omega_1$-алгебры $A$ в $\Omega_2$-алгебре $M$. Множество автоморфизмов представления $f$ порождает* **группу** *$GA(f)$ .*

Доказательство. Пусть $R$, $P$ - автоморфизмы представления $f$. Согласно определению 3.3.4, отображения $R$, $P$ являются автоморфизмами $\Omega_2$-алгебры $M$. Согласно теореме II.3.2, ([13], с. 60), отображение $R \circ P$ является автоморфизмом $\Omega_2$-алгебры $M$. Из теоремы 3.2.12 и определения 3.3.4 следует, что произведение автоморфизмов $R \circ P$ представления $f$ является автоморфизмом представления $f$.

Пусть $R$, $P$, $Q$ - автоморфизмы представления $f$. Из цепочки равенств

$$\begin{aligned}((R\circ P)\circ Q)(a) &= (R\circ P)(Q(a)) = R(P(Q(a)))\\ &= R((P\circ Q)(a)) = (R\circ(P\circ Q))(a)\end{aligned}$$

следует ассоциативность произведения для отображений $R$, $P$, $Q$.[3.3]

Пусть $R$ - автоморфизм представления $f$. Согласно определению 3.3.4 отображение $R$ является автоморфизмом $\Omega_2$-алгебры $M$. Следовательно, отображение $R^{-1}$ является автоморфизмом $\Omega_2$-алгебры $M$. Для автоморфизма $R$ представления справедливо равенство (3.2.4). Положим $m' = R(m)$. Так как $R$ - автоморфизм $\Omega_2$-алгебры, то $m = R^{-1}(m')$ и равенство (3.2.4) можно записать в виде

$$R(f(a')(R^{-1}(m'))) = f(a')(m') \tag{3.3.1}$$

Так как отображение $R$ является автоморфизмом $\Omega_2$-алгебры $M$, то из равенства (3.3.1) следует

$$f(a')(R^{-1}(m')) = R^{-1}(f(a')(m')) \tag{3.3.2}$$

Равенство (3.3.2) соответствует равенству (3.2.4) для отображения $R^{-1}$. Следовательно, отображение $R^{-1}$ является автоморфизмом представления $f$. □

[3.3] При доказательстве ассоциативности произведения я следую примеру полугруппы из [4], с. 20, 21.

Глава 4

# Представление группы

## 4.1. Представление группы

Группа - одна из немногих алгебр, которая позволяет рассматривать произведение преобразований $\Omega$-алгебры $M$ таким образом, что если преобразования принадлежат представлению, то их произведение также принадлежит представлению. При этом следует помнить, что порядок отображений при суперпозиции зависит от порядка отображений на диаграмме и с какой стороны отображения действуют на элементы множества.

Определение 4.1.1. Пусть ${}^*M$ - группа с произведением

$$(f \circ g)x = f(gx)$$

и $\delta$ - единица группы ${}^\star M$. Пусть $G$ - группа. Мы будем называть гомоморфизм групп

$$f : G \to {}^\star M \tag{4.1.1}$$

**левосторонним представлением группы** $G$ или $G*$**-представлением группы** в $\Omega$-алгебре $M$, если отображение $f$ удовлетворяет условиям

$$f(ab)u = f(a)(f(b)u) \tag{4.1.2}$$

□

Замечание 4.1.2. Поскольку отображение (4.1.1) - гомоморфизм, то

$$f(ab)u = (f(a)f(b))u \tag{4.1.3}$$

Мы здесь пользуемся соглашением

$$f(a)f(b) = f(a) \circ f(b)$$

Таким образом, концепция представления групп состоит в том, что в каком порядке мы перемножаем элементы группы, в том же порядке перемножаются соответствующие преобразования представления. Из равенств (4.1.2) и (4.1.3) следует

$$(f(a)f(b))u = f(a)(f(b)u) \tag{4.1.4}$$

Равенство (4.1.4) совместно с ассоциативностью произведения преобразований представляет собой **закон ассоциативности** для $G*$-представления. Это позволяет записывать равенство (4.1.4) без использования скобок

$$f(ab)u = f(a)f(b)u$$

□

Определение 4.1.3. Пусть $M^\star$ - группа с произведением

$$x(f\circ g)=(xf)g$$

и $\delta$ - единица группы $M^\star$. Пусть $G$ - группа. Мы будем называть гомоморфизм групп

$$f:G\to M^* \tag{4.1.5}$$

**правосторонним представлением группы** $G$ или $*G$**-представлением** в $\Omega$-алгебре $M$, если отображение $f$ удовлетворяет условиям

$$uf(ab)=(uf(a))f(b) \tag{4.1.6}$$

□

Замечание 4.1.4. Поскольку отображение (4.1.5) - гомоморфизм, то

$$uf(ab)=u(f(a)f(b)) \tag{4.1.7}$$

Из равенств (4.1.6) и (4.1.7) следует

$$u(f(a)f(b))=(uf(a))f(b) \tag{4.1.8}$$

Равенство (4.1.8) совместно с ассоциативностью произведения преобразований представляет собой **закон ассоциативности** для $*G$-представления. Это позволяет записывать равенство (4.1.8) без использования скобок

$$uf(ab)=uf(a)f(b)$$

□

Определение 4.1.5. Мы будем называть преобразование

$$t:M\to M$$

**невырожденным преобразованием**, если существует обратное отображение. □

Теорема 4.1.6. *Для любого $g\in G$ преобразование $f(g)$ является невырожденным и удовлетворяет равенству*

$$f(g^{-1})=f(g)^{-1} \tag{4.1.9}$$

Доказательство. На основании (4.1.2) и

$$f(e)=\delta$$

мы можем записать

$$u=\delta(u)=f(gg^{-1})(u)=f(g)(f(g^{-1})(u))$$

Это завершает доказательство. □

Теорема 4.1.7. *Групповая операция определяет два различных представления на группе:*

- **Левый сдвиг** $t_\star$

$$\begin{aligned}b'&=t_\star(a)b=ab\\ b'&=t_\star(a)(b)=ab\end{aligned} \tag{4.1.10}$$

*является $G*$-представлением на множестве*[4.1] $G$

$$t_\star(ab) = t_\star(a) \circ t_\star(b) \tag{4.1.11}$$

- **Правый сдвиг** $_\star t$

$$\begin{aligned} b' &= b\ {}_\star t(a) = ba \\ b' &= {}_\star t(a)(b) = ba \end{aligned} \tag{4.1.12}$$

*является $*G$-представлением на множестве $G$*

$${}_\star t(ab) = {}_\star t(a) \circ {}_\star t(b) \tag{4.1.13}$$

Доказательство. Равенство (4.1.11) следует из ассоциативности произведения

$$t_\star(ab)c = (ab)c = a(bc) = t_\star(a)(t_\star(b)c) = (t_\star(a) \circ t_\star(b))c$$

Аналогично доказывается равенство (4.1.13). □

Определение 4.1.8. Пусть $G$ - группа. Пусть $f$ - $G*$-представление в $\Omega$-алгебре $M$. Для любого $v \in M$ мы определим **орбиту представления группы** $G$ как множество

$$f(G)v = \{w = f(g)v : g \in G\}$$

□

Так как $f(e) = \delta$, то $v \in f(G)v$.

Теорема 4.1.9. *Если*

$$v \in f(G)u \tag{4.1.14}$$

*то*

$$f(G)u = f(G)v$$

Доказательство. Из (4.1.14) следует существование $a \in G$ такого, что

$$v = f(a)u \tag{4.1.15}$$

Если $w \in f(G)v$, то существует $b \in G$ такой, что

$$w = f(b)v \tag{4.1.16}$$

Подставив (4.1.15) в (4.1.16), мы получим

$$w = f(b)(f(a)u) \tag{4.1.17}$$

На основании (4.1.2) из (4.1.17) следует, что $w \in f(G)u$. Таким образом,

$$f(G)v \subseteq f(G)u$$

На основании (4.1.9) из (4.1.15) следует, что

$$u = f(a)^{-1}v = f(a^{-1})v \tag{4.1.18}$$

Равенство (4.1.18) означает, что $u \in f(G)v$ и, следовательно,

$$f(G)u \subseteq f(G)v$$

Это завершает доказательство. □

[4.1] Левый сдвиг не является представлением группы в группе, так как преобразование $t_\star$ не является гомоморфизмом группы. Аналогичное замечание верно для правого сдвига.

Таким образом, $G*$-представление $f$ в $\Omega$-алгебре $M$ порождает отношение эквивалентности $S$ и орбита $f(G)u$ является классом эквивалентности. Мы будем пользоваться обозначением $M/f(G)$ для фактор множества $M/S$ и мы будем называть это множество **пространством орбит** $G*$**-представления** $f$.

Теорема 4.1.10. *Если определены $G*$-представление $f_1$ в $\Omega$-алгебре $M_1$ и $G*$-представление $f_2$ в $\Omega$-алгебре $M_2$, то мы можем определить* **прямое произведение** $G*$**-представлений** *$f_1$ и $f_2$ группы*

$$f = f_1 \times f_2 : G \to M_1 \otimes M_2$$
$$f(g) = (f_1(g), f_2(g))$$

Доказательство. Чтобы показать, что $f$ является представлением, достаточно показать, что $f$ удовлетворяет определению 4.1.1.

$$f(e) = (f_1(e), f_2(e)) = (\delta_1, \delta_2) = \delta$$

$$\begin{aligned} f(ab)u &= (f_1(ab)u_1, f_2(ab)u_2) \\ &= (f_1(a)(f_1(b)u_1), f_2(a)(f_2(b)u_2)) \\ &= f(a)(f_1(b)u_1, f_2(b)u_2) \\ &= f(a)(f(b)u) \end{aligned}$$

□

## 4.2. Однотранзитивное правостороннее представление группы

Определение 4.2.1. Мы будем называть **ядром неэффективности** $G*$**-представления** множество

$$K_f = \{g \in G : f(g) = \delta\}$$

□

Теорема 4.2.2. *Ядро неэффективности $G*$-представления - это подгруппа группы $G$.*

Доказательство. Допустим $f(a_1) = \delta$ и $f(a_2) = \delta$. Тогда

$$f(a_1a_2)u = f(a_1)(f(a_2)u) = u$$
$$f(a^{-1}) = f^{-1}(a) = \delta$$

□

Теорема 4.2.3. *$G*$-представление* **эффективно** *тогда и только тогда, когда ядро неэффективности $K_f = \{e\}$.*

Доказательство. Утверждение является следствием определений 3.1.6 и 4.2.1 и теоремы 4.2.2. □

Если действие не эффективно, мы можем перейти к эффективному заменив группой $G_1 = G|K_f$, пользуясь факторизацией по ядру неэффективности. Это означает, что мы можем изучать только эффективное действие.

Определение 4.2.4. Рассмотрим $G*$-представление $f$ в $\Omega$-алгебре $M$. **Малая группа** или **группа стабилизации** элемента $x \in M$ - это множество

$$= \{g \in G : f(g)x = x\}$$

Мы будем называть $G*$-представление $f$ **свободным**, если для любого $x \in M$ группа стабилизации $G_x = \{e\}$. □

Теорема 4.2.5. *Если определено свободное $G*$-представление $f$ группы $G$ на $\Omega$-алгебре $A$, то определено взаимно однозначное соответствие между орбитами представления, а также между орбитой представления и группой $G$.*

Доказательство. Допустим для точки $a \in A$ существуют $g_1, g_2 \in G$

$$f(g_1)a = f(g_2)a \tag{4.2.1}$$

Умножим обе части равенства (4.2.1) на $f(g_1^{-1})$

$$a = f(g_1^{-1})f(g_2)a$$

Поскольку представление свободное, $g_1 = g_2$. Теорема доказана, так как мы установили взаимно однозначное соответствие между орбитой и группой $G$. □

Определение 4.2.6. Мы будем называть пространство $V$ **однородным пространством группы** $G$, если мы имеем однотранзитивное $G*$-представление на $V$. □

Теорема 4.2.7. *Если мы определим однотранзитивное представление $f$ группы $G$ на $\Omega$-алгебре $A$, то мы можем однозначно определить координаты на $A$, пользуясь координатами на группе $G$.*

*Если $f$ - левостороннее представление, то $f(a)$ эквивалентно левому сдвигу $t_\star(a)$ на группе $G$. Если $f$ - правостороннее представление, то $f(a)$ эквивалентно правому сдвигу ${}_\star t(a)$ на группе $G$.*

Доказательство. Мы выберем точку $v \in A$ и определим координаты точки $w \in A$ как координаты $a \in G$ такого, что $w = f(a)v$. Координаты, определённые таким образом, однозначны с точностью до выбора начальной точки $v \in A$, так как действие эффективно.

Если $f$ - левостороннее представление, мы будем пользоваться записью

$$f(a)v = av$$

Так как запись

$$f(a)(f(b)v) = a(bv) = (ab)v = f(ab)v$$

совместима с групповой структурой, мы видим, что левостороннее представление $f$ эквивалентно левому сдвигу.

Если $f$ - правостороннее представление, мы будем пользоваться записью

$$vf(a) = va$$

Так как запись

$$(vf(b))f(a) = (vb)a = v(ba) = vf(ba)$$

совместима с групповой структурой, мы видим, что правостороннее представление $f$ эквивалентно правому сдвигу. □

Замечание 4.2.8. Мы будем записывать эффективное $G*$-представление в форме

$$v' = t_\star(a)v = av$$

Орбита этого представления имеет вид

$$Gv = t_\star(G)v$$

Мы будем пользоваться обозначением $M/t_\star(G)$ для пространства орбит эффективного $G*$-представления. □

Замечание 4.2.9. Мы будем записывать эффективное $*G$-представление в форме

$$v' = v\ {}_\star t(a) = va$$

Орбита этого представления имеет вид

$$vG = v\ {}_\star t(G)$$

Мы будем пользоваться обозначением $M/{}_\star t(G)$ для пространства орбит эффективного $*G$-представления. □

Теорема 4.2.10. *Свободное $G*$-представление эффективно. Свободное $G*$-представление $f$ в $\Omega$-алгебре $M$ однотранзитивно на орбите.*

Доказательство. Следствие определения 4.2.4. □

Теорема 4.2.11. *Правый и левый сдвиги на группе $G$ перестановочны.*

Доказательство. Это следствие ассоциативности группы $G$

$$(t_\star(a) \circ\ {}_\star t(b))c = a(cb) = (ac)b = ({}_\star t(b) \circ t_\star(a))c$$

□

Теорема 4.2.11 может быть сформулирована следующим образом.

Теорема 4.2.12. *Пусть $G$ - группа. Для любого $a \in G$ отображение $t_\star(a)$ является автоморфизмом представления ${}_\star t$.*

Доказательство. Согласно теореме 4.2.11

$$t_\star(a) \circ {}_\star t(b) = {}_\star t(b) \circ t_\star(a) \tag{4.2.2}$$

Равенство (4.2.2) совпадает с равенством (3.2.3) из определения 3.2.2 при условии $r = id$, $R = t_\star(a)$. □

Теорема 4.2.13. *Пусть $G*$-представление $f$ на $\Omega$-алгебре $M$ однотранзитивно. Тогда мы можем однозначно определить однотранзитивное $*G$-представление $h$ на $\Omega$-алгебре $M$ такое, что диаграмма*

$$\begin{array}{ccc}
M & \xrightarrow{h(a)} & M \\
\Big\downarrow{\scriptstyle f(b)} & & \Big\downarrow{\scriptstyle f(b)} \\
M & \xrightarrow[h(a)]{} & M
\end{array}$$

*коммутативна для любых $a, b \in G$.*[4.2]

[4.2] Это утверждение можно также найти в [3].

Доказательство. Мы будем пользоваться групповыми координатами для точек $v \in M$. Тогда согласно теореме 4.2.7 мы можем записать левый сдвиг $t_\star(a)$ вместо преобразования $f(a)$.

Пусть $v_0, v \in M$. Тогда мы можем найти одно и только одно $a \in G$ такое, что

$$v = v_0 a = v_0\ {}_\star t(a)$$

Мы предположим

$$h(a) = {}_\star t(a)$$

Существует $b \in G$ такое, что

$$w_0 = f(b)v_0 = t_\star(b)v_0 \qquad w = f(b)v = t_\star(b)v$$

Согласно теореме 4.2.11 диаграмма

$$\begin{array}{ccc} v_0 & \xrightarrow{h(a)=\ {}_\star t(a)} & v \\ \Big\downarrow{\scriptstyle f(b)=t_\star(b)} & & \Big\downarrow{\scriptstyle f(b)=t_\star(b)} \\ w_0 & \xrightarrow[h(a)=\ {}_\star t(a)]{} & w \end{array} \tag{4.2.3}$$

коммутативна.

Изменяя $b$ мы получим, что $w_0$ - это произвольная точка, принадлежащая $M$.

Мы видим из диаграммы, что, если $v_0 = v$, то $w_0 = w$ и следовательно $h(e) = \delta$. С другой стороны, если $v_0 \neq v$, то $w_0 \neq w$ потому, что $G*$-представление $f$ однотранзитивно. Следовательно $*G$-представление $h$ эффективно.

Таким же образам мы можем показать, что для данного $w_0$ мы можем найти $a$ такое, что $w = h(a)w_0$. Следовательно $*G$-представление $h$ однотранзитивно.

В общем случае, произведение преобразований $G*$-представления $f$ не коммутативно и следовательно $*G$-представление $h$ отлично от $G*$-представления $f$. Таким же образом мы можем создать $G*$-представление $f$, пользуясь $*G$-представлением $h$. □

Мы будем называть представления $f$ и $h$ **парными представлениями группы** $G$.

Замечание 4.2.14. Очевидно, что преобразования $t_\star(a)$ и ${}_\star t(a)$ отличаются, если группа $G$ неабелева. Тем не менее, они являются отображениями на. Теорема 4.2.13 утверждает, что, если оба представления правого и левого сдвига существуют на множестве $M$, то мы можем определить два перестановочных представления на множестве $M$. Только правый или левый сдвиг не может представлять оба типа представления. Чтобы понять почему это так, мы можем изменить диаграмму (4.2.3) и предположить $h(a)v_0 = t_\star(a)v_0 = v$ вместо $h(a)v_0 = v_{0\star}t(a) = v$ и проанализировать, какое выражение $h(a)$ имеет в точке $w_0$. Диаграмма

$$\begin{array}{ccc} v_0 & \xrightarrow{h(a)=t_\star(a)} & v \\ \Big\downarrow{\scriptstyle f(b)=t_\star(b)} & & \Big\downarrow{\scriptstyle f(b)=t_\star(b)} \\ w_0 & \xrightarrow[h(a)]{} & w \end{array}$$

эквивалентна диаграмме

$$\begin{array}{ccc}
v_0 & \xrightarrow{h(a)=t_\star(a)} & v \\
\Big\uparrow{\scriptstyle f^{-1}(b)=t_\star(b^{-1})} & & \Big\downarrow{\scriptstyle f(b)=t_\star(b)} \\
w_0 & \xrightarrow[h(a)]{} & w
\end{array}$$

и мы имеем $w = bv = bav_0 = bab^{-1}w_0$. Следовательно

$$h(a)w_0 = (bab^{-1})w_0$$

Мы видим, что представление $h$ зависит от его аргумента. □

Теорема 4.2.15. *Пусть $f$ и $h$ - парные преставления группы $G$. Для любого $a \in G$ отображение $h(a)$ является автоморфизмом представления $f$.*

Доказательство. Следствие теорем 4.2.12 и 4.2.13. □

Замечание 4.2.16. Существует ли морфизм представлений из $t_\star$ в $t_\star$, отличный от автоморфизма $(\mathrm{id}, {}_\star t(a))$? Если мы положим

$$r(g) = cgc^{-1}$$
$$R(a)(m) = cmac^{-1}$$

то нетрудно убедиться, что отображение $(r, R(a))$ является морфизмом представлений из $t_\star$ в $t_\star$. Но это отображение не является автоморфизмом представления $t_\star$, так как $r \neq \mathrm{id}$. □

Глава 5

# Векторное пространство над телом

## 5.1. Векторное пространство

Чтобы определить левостороннее представление

$$f : D \xrightarrow{*} M \qquad f(d) : v \to d\,v$$

кольца $D$ в $\Omega$-алгебре $M$, мы должны определить структуру кольца на множестве ${}^{\star}M$. [5.1]

Теорема 5.1.1. *Левостороннее представление $f$ кольца $D$ в $\Omega$-алгебре $M$ определенно тогда и только тогда, когда определены левосторонние представления мультипликативной и аддитивной групп кольца $D$ и эти представления удовлетворяют соотношению*

$$f(a(b+c)) = f(a)f(b) + f(a)f(c)$$

Доказательство. Теорема следует из определения 3.1.2. □

Определение 5.1.2. Абелева группа $M$ является $D*$**-модулем**, если определено $D*$-представление

$$f : D \xrightarrow{*} M \qquad f(d) : v \to d\,v$$

□

[5.1]Можно ли определить операцию сложения на множестве ${}^{\star}M$, если эта операция не определена на множестве $M$. Ответ на этот вопрос положительный.

Допустим $M = B \cup C$ и $F : B \to C$ - взаимно однозначное отображение. Мы определим множество ${}^{\star}M$ левосторонних преобразований множества $M$ согласно следующему правилу. Пусть $V \subseteq B$. левостороннее преобразование $F_V$ имеет вид

$$F_V x = \begin{cases} x & x \in B\backslash V \\ Fx & x \in V \\ x & x \in C\backslash F(V) \\ F^{-1}x & x \in F(V) \end{cases}$$

Мы определим сложение левосторонних преобразований по правилу

$$F_V + F_W = F_{V \triangle W}$$
$$V \triangle W = (V \cup W)\backslash(V \cap W)$$

Очевидно, что

$$F_\emptyset + F_V = F_V$$
$$F_V + F_V = F_\emptyset$$

Следовательно, отображение $F_\emptyset$ является нулём относительно сложения, а множество ${}^{\star}M$ является абелевой группой.

Согласно нашим обозначениям $D*$-модуль - это **левый модуль над кольцом** $D$ и $*D$**-модуль** - это **правый модуль над кольцом** $D$.

Поле является частным случаем кольца. Поэтому векторное пространство над полем имеет больше свойств чем модуль над кольцом. Очень трудно, если вообще возможно, распространить определения, работающие в векторном пространстве, на модуль над произвольным кольцом. Определение базиса и размерности векторного пространства тесно связаны с возможностью найти решение линейного уравнения в кольце. Свойства линейного уравнения в теле близки к свойствам линейного уравнения в поле. Поэтому мы надеемся что свойства векторного пространства над телом близки к свойствам векторного пространства над полем.

Теорема 5.1.3. *Левостороннее представление тела $D$* **эффективно**, *если эффективно левостороннее представление мультипликативной группы тела $D$.*

Доказательство. Пусть

$$f : D \xrightarrow{*} M \qquad f(d) : v \to d\,v$$

левостороннее представление тела $D$. Если элементы $a$, $b$ мультипликативной группы порождают одно и то же левостороннее преобразование, то

$$f(a)m = f(b)m \tag{5.1.1}$$

для любого $m \in M$. Выполняя преобразование $f(a^{-1})$ над обеими частями равенства (5.1.1), мы получим

$$m = f(a^{-1})(f(b)m) = f(a^{-1}b)m$$

□

Согласно замечанию 3.1.7, если представление тела эффективно, мы отождествляем элемент тела и соответствующее ему левостороннее преобразование.

Определение 5.1.4. Пусть $D$ является телом. Абелева группа $V$ является $D*$**-векторным пространством**, если определено эффективное $D*$-представление

$$f : D \xrightarrow{*} M \qquad f(d) : v \to d\,v \tag{5.1.2}$$

Абелева группа $V$ является $*D$**-векторным пространством**, если определено эффективное $*D$-представление

$$f : D \xrightarrow{*} M \qquad f(d) : v \to v\,d \tag{5.1.3}$$

□

$D*$-векторное пространство называется также **левым $D$-векторным пространством** или левым векторным пространством над телом $D$. $*D$-векторное пространство называется также **правым $D$-векторным пространством** или правым векторным пространством над телом $D$.

Теорема 5.1.5. *Элементы $D*$-векторного пространства $V$ удовлетворяют соотношениям*

- **закону ассоциативности**

$$(ab)m = a(bm) \tag{5.1.4}$$

- **закону дистрибутивности**

$$a(m+n) = am + an \tag{5.1.5}$$

$$(a+b)m = am + bm \tag{5.1.6}$$

- **закону унитарности**

$$1m = m \tag{5.1.7}$$

*для любых* $a,\ b \in D,\ m,\ n \in V$.

Доказательство. Равенство (5.1.5) следует из утверждения, что левостороннее преобразование $a$ является эндоморфизмом абелевой группы. Равенство (5.1.6) следует из утверждения, что представление (5.1.2) является гомоморфизмом аддитивной группы тела $D$. Равенства (5.1.4) и (5.1.7) следуют из утверждения, что представление (5.1.2) является левосторонним представлением мультипликативной группы тела $D$. □

Согласно нашим обозначениям $D*$-векторное пространство - это **левое $D$-векторное пространство**. Отображение

$$(d, v) \in D \times V \to dv \in V$$

порождённое $D*$-представлением (5.1.2), называется **левосторонним произведением вектора на скаляр**.

Согласно нашим обозначениям $*D$**-векторное пространство** - это **правое $D$-векторное пространство**. Отображение

$$(v, d) \in D \times V \to vd \in V$$

порождённое $*D$-представлением (5.1.3), называется **правосторонним произведением вектора на скаляр**.

Любое утверждение, справедливое для левого $D$-векторного пространства, справедливо для правого $D$-векторного пространства, если левостороннее произведение вектора на скаляр заменить правосторонним произведением вектора на скаляр.

Определение 5.1.6. Пусть $V$ - $D*$-векторное пространство над телом $D$. Множество векторов $N$ - **подпространство $D*$-векторного пространства** $V$, если

$$a + b \in N$$

$$ka \in N$$

$$a, b \in N \quad k \in D$$

□

Пример 5.1.7. Определим на множестве $D^m_n$ $m \times n$ матриц над телом $D$ операцию сложения

$$a + b = \left(a^j_i\right) + \left(b^j_i\right) = \left(a^j_i + b^j_i\right)$$

и умножения на скаляр

$$da = d\left(a^j_i\right) = \left(da^j_i\right)$$

$a = 0$ тогда и только тогда, когда $a^j_i = 0$ для любых $i$, $j$. Непосредственная проверка показывает, что $D^m_n$ является $D*$-векторным пространством, если произведение действует слева. В противном случае $D_n^m$ - $*D$-векторное пространство. Мы будем называть векторное пространство $D_n^m$ $D*$**-векторным пространством матриц**. □

## 5.2. Тип векторного пространства

Произведение вектора на скаляр асимметрично, так как произведение определено для объектов разных множеств. Однако различие между $D*$-и $*D$ -векторным пространством появляется только тогда, когда мы переходим к координатному представлению. Говоря, векторное пространство является $D*$ - или $*D$-, мы указываем, с какой стороны, слева или справа, мы умножаем координаты вектора на элементы тела.

Определение 5.2.1. Допустим $u$, $v$ - векторы $D*$-векторного пространства $V$. Мы будем говорить, что вектор $w$ является **линейной комбинацией векторов** $u$ и $v$, если мы можем записать

$$w = au + bv$$

где $a$ и $b$ - скаляры. □

Мы можем распространить понятие линейной комбинации на любое конечное семейство векторов. Пользуясь обобщённым индексом для нумерации векторов, мы можем представить семейство векторов в виде одномерной матрицы. Мы пользуемся соглашением, что в заданном векторном пространстве мы представляем любое семейство векторов либо в виде $_*$-строки (матрицы строки), либо $^*$-строки (матрицы столбца). Это представление определяет характер записи линейной комбинации. Рассматривая это представление в $D*$-или $*D$-векторном пространстве, мы получаем четыре разных модели векторного пространства, рассмотренные в примерах 5.2.2, 5.2.3, 5.2.4, 5.2.5.

Чтобы иметь возможность, не меняя записи указать, является векторное пространство $D*$- или $*D$-векторным пространством, мы вводим новое обозначение. Символ $D*$ называется **типом векторного пространства** и означает, что мы изучаем $D*$-векторное пространство. Операция умножения в типе векторного пространства указывает на матричную операцию, используемую в линейной комбинации.

Пример 5.2.2. Представим множество векторов ${}^ia$, $i \in I$, $D*$-векторного пространства $V$ в виде $^*$-строки (матрицы столбца)

$$a = \begin{pmatrix} {}^1a \\ \dots \\ {}^na \end{pmatrix}$$

и множество скаляров $c_i$, $i \in I$, в виде $_*$-строки (матрицы строки)

$$c = \begin{pmatrix} c_1 & \dots & c_n \end{pmatrix}$$

Тогда мы можем записать линейную комбинацию векторов ${}^ia$ в виде

$$c_i\,{}^ia = c_*{}^*a$$

Соответствующая реализация $D*$-векторного пространства называется $D_*{}^*$-**векторным пространством** или **левым $D$-векторным пространством строк**. □

Пример 5.2.3. Представим множество векторов ${}_ia$, $i \in I$, $D*$-векторного пространства $V$ в виде subs строки (матрицы строки)

$$a = \begin{pmatrix} {}_1a & ... & {}_na \end{pmatrix}$$

и множество скаляров $c^i$, $i \in I$, в виде ${}^*$-строки (матрицы столбца)

$$c = \begin{pmatrix} c^1 \\ ... \\ c^n \end{pmatrix}$$

Тогда мы можем записать линейную комбинацию векторов ${}_ia$ в виде

$$c^i\ {}_ia = c^*{}_*a$$

Соответствующая реализация $D*$-векторного пространства называется $D^*{}_*$-**векторным пространством** или **левым $D$-векторным пространством столбцов**. □

Пример 5.2.4. Представим множество векторов $a^i$, $i \in I$, $*D$-векторного пространства $V$ в виде ${}^*$-строки (матрицы столбца)

$$a = \begin{pmatrix} a^1 \\ ... \\ a^n \end{pmatrix}$$

и множество скаляров ${}_ic$, $i \in I$, в виде ${}_*$-строки (матрицы строки)

$$c = \begin{pmatrix} {}_1c & ... & {}_nc \end{pmatrix}$$

Тогда мы можем записать линейную комбинацию векторов $a^i$ в виде

$$a^i\ {}_ic = a^*{}_*c$$

Соответствующая реализация $*D$-векторного пространства называется ${}^*{}_*D$-**векторным пространством** или **правым $D$-векторным пространством строк**. □

Пример 5.2.5. Представим множество векторов $a_i$, $i \in I$, $*D$-векторного пространства $V$ в виде ${}_*$-строки (матрицы строки)

$$a = \begin{pmatrix} a_1 & ... & a_n \end{pmatrix}$$

и множество скаляров ${}^ic$, $i \in I$, в виде ${}^*$-строки (матрицы столбца)

$$c = \begin{pmatrix} {}^1c \\ ... \\ {}^nc \end{pmatrix}$$

Тогда мы можем записать линейную комбинацию векторов $a_i$ в виде

$$a_i\,{}^i c = a_*{}^* c$$

Соответствующая реализация $*D$-векторного пространства называется ${}_*{}^*D$-**векторным пространством** или **правым** $D$**-векторным пространством столбцов**. □

Замечание 5.2.6. Мы распространим на векторное пространство и его тип соглашение, описанное в замечании 2.2.15. Например, в выражении

$$A_*{}^* B_*{}^* v\lambda$$

мы выполняем операцию умножения слева направо. Это соответствует ${}_*{}^*D$-векторному пространству. Однако мы можем выполнять операцию умножения справа налево. В традиционной записи это выражение примет вид

$$\lambda v^*{}_* B^*{}_* A$$

и будет соответствовать $D^*{}_*$-векторному пространству. Аналогично, читая выражение снизу вверх мы получим выражение

$$A^*{}_* B^*{}_* v\lambda$$

соответствующего ${}^*{}_*D$-векторному пространству. □

## 5.3. Базис ${}_*{}^*D$-векторного пространства

Определение 5.3.1. Векторы $a_i$, $i \in I$, ${}_*{}^*D$-векторного пространства $V$ **линейно независимы**, если $c = 0$ следует из уравнения

$$a_*{}^* c = 0$$

В противном случае, векторы $a_i$ **линейно зависимы**. □

Определение 5.3.2. Множество векторов $\overline{\overline{e}} = (e_i, i \in I)$ - **базис** ${}_*{}^*D$**-векторного пространства**, если векторы $e_i$ линейно независимы и добавление любого вектора к этой системе делает эту систему линейно зависимой. □

Теорема 5.3.3. *Если $\overline{\overline{e}}$ - базис ${}_*{}^*D$-векторного пространства $V$, то любой вектор $\overline{v} \in V$ имеет одно и только одно разложение*

$$\overline{v} = e_*{}^* v \tag{5.3.1}$$

*относительно этого базиса.*

Доказательство. Так как система векторов $e_i$ является максимальным множеством линейно независимых векторов, система векторов $\overline{v}$, $e_i$ - линейно зависима и в уравнении

$$\overline{v}b + e_*{}^* c = 0 \tag{5.3.2}$$

по крайней мере $b$ отлично от 0. Тогда равенство

$$\overline{v} = e_*{}^*(-cb^{-1}) \tag{5.3.3}$$

следует из (5.3.2). (5.3.1) следует из (5.3.3).

Допустим мы имеем другое разложение

$$\overline{v} = e_*{}^* v' \tag{5.3.4}$$

Вычтя (5.3.1) из (5.3.4), мы получим

$$0 = e_*{}^*(v' - v)$$

Так как векторы $e_i$ линейно независимы, мы имеем

$$v' - v = 0$$

□

Определение 5.3.4. Мы будем называть матрицу $v$ разложения (5.3.1) **координатной матрицей вектора** $\overline{v}$ в базисе $\overline{\overline{e}}$ и её элементы **координатами вектора** $\overline{v}$ относительно базиса $\overline{\overline{e}}$. □

Теорема 5.3.5. *Множество координат $a$ вектора $\overline{a}$ в базисе $\overline{\overline{e}}$ ${}_*{}^*D$-векторного пространства порождают ${}_*{}^*D$-векторное пространство $D^n$, изоморфное ${}_*{}^*D$-векторному пространству $V$. Это ${}_*{}^*D$-векторное пространство называется* **координатным ${}_*{}^*D$-векторным пространством**, *а изоморфизм* **координатным изоморфизмом**.

Доказательство. Допустим векторы $\overline{a}$ и $\overline{b} \in V$ имеют разложение

$$\overline{a} = e_*{}^*a$$
$$\overline{b} = e_*{}^*b$$

в базисе $\overline{\overline{e}}$. Тогда

$$\overline{a} + \overline{b} = e_*{}^*a + e_*{}^*b = e_*{}^*(a+b)$$
$$\overline{a}m = (e_*{}^*a)m = e_*{}^*(am)$$

для любого $m \in D$. Таким образом, операции в векторном пространстве определены по координатно

$${}^i(a+b) = {}^ia + {}^ib$$
$${}^i(am) = {}^iam$$

Это доказывает теорему. □

Пример 5.3.6. Пусть $\overline{\overline{e}} = ({}^je, j \in J, |J| = n)$ - **базис $D_*{}^*$-векторного пространства** $V$. Согласно примеру 5.2.2, мы можем представить базис $\overline{\overline{e}}$ в виде ${}^*$-строки (матрицы столбца)

$$\overline{\overline{e}} = \begin{pmatrix} {}^1e \\ ... \\ {}^ne \end{pmatrix}$$

Координатная матрица

$$a = \begin{pmatrix} a_1 & ... & a_n \end{pmatrix} = (a_j, j \in J)$$

вектора $\overline{a}$ в базисе $\overline{\overline{e}}$ называется **$D_*{}^*$-вектором**[5.2] или **$D*$-вектор строкой**.

Пусть ${}^*$-строка

$$\overline{A} = \begin{pmatrix} {}^1\overline{A} \\ ... \\ {}^m\overline{A} \end{pmatrix} = ({}^i\overline{A}, i \in I) \tag{5.3.5}$$

задаёт множество векторов. Векторы ${}^i\overline{A}$ имеют разложение

$${}^i\overline{A} = {}^iA_*{}^*e$$

[5.2] $D_*{}^*$-вектор является аналогом вектор-строки.

Если мы подставим координатные матрицы вектора ${}^i\overline{A}$ в матрицу (5.3.5), мы получим матрицу

$$A=\begin{pmatrix}\begin{pmatrix}{}^1A_1 & ... & {}^1A_n\end{pmatrix}\\ ...\\ \begin{pmatrix}{}^mA_1 & ... & {}^mA_n\end{pmatrix}\end{pmatrix}=\begin{pmatrix}{}^1A_1 & ... & {}^1A_n\\ ... & ... & ...\\ {}^mA_1 & ... & {}^mA_n\end{pmatrix}=({}^iA_j)$$

Мы будем называть матрицу $A$ **координатной матрицей множества векторов** $({}^i\overline{A}, i\in I)$ в базисе $\overline{\overline{e}}$ и её элементы **координатами множества векторов** $({}^i\overline{A}, i\in I)$ в базисе $\overline{\overline{e}}$.

Пусть ${}^*$-строка

$$\overline{\overline{f}}=\begin{pmatrix}{}^1f\\ ...\\ {}^nf\end{pmatrix}=({}^if, i\in J)$$

является базисом $D_*{}^*$-векторного пространства $V$. Мы будем говорить, что координатная матрица $f$ множества векторов $({}^if, i\in J)$ определяет **координаты** ${}^if_j$ **базиса** $\overline{\overline{f}}$ относительно базиса $\overline{\overline{e}}$. □

Пример 5.3.7. Пусть $\overline{\overline{e}}=({}_je, j\in J, |J|=n)$ - **базис** $D^*{}_*$**-векторного пространства** $V$. Согласно примеру 5.2.3, мы можем представить базис $\overline{\overline{e}}$ в виде ${}_*$-строки (матрицы строки)

$$\overline{\overline{e}}=\begin{pmatrix}{}_1e & ... & {}_ne\end{pmatrix}$$

Координатная матрица

$$a=\begin{pmatrix}a^1\\ ...\\ a^n\end{pmatrix}=(a_j, j\in J)$$

вектора $\overline{a}$ в базисе $\overline{\overline{e}}$ называется $D_*{}^*$**-вектором**[5.3] или $D*$**-вектор столбцом**.

Пусть ${}_*$-строка

$$\overline{A}=\begin{pmatrix}{}_1\overline{A} & ... & {}_m\overline{A}\end{pmatrix}=({}^i\overline{A}, i\in I) \tag{5.3.6}$$

задаёт множество векторов. Векторы ${}^i\overline{A}$ имеют разложение

$${}_i\overline{A}={}_iA^*{}_*e$$

Если мы подставим координатные матрицы вектора ${}_i\overline{A}$ в матрицу (5.3.6), мы получим матрицу

$$A=\begin{pmatrix}\begin{pmatrix}{}_1A^1\\ ...\\ {}_1A^n\end{pmatrix} & ... & \begin{pmatrix}{}_mA^1\\ ...\\ {}_mA^n\end{pmatrix}\end{pmatrix}=\begin{pmatrix}{}_1A^1 & ... & {}_1A^n\\ ... & ... & ...\\ {}_mA^1 & ... & {}_mA^n\end{pmatrix}=({}_iA^j)$$

[5.3] $D_*{}^*$-вектор является аналогом вектор-строки.

Мы будем называть матрицу $A$ **координатной матрицей множества векторов** $({}_i\overline{A}, i \in I)$ в базисе $\overline{\overline{e}}$ и её элементы **координатами множества векторов** $({}_i\overline{A}, i \in I)$ в базисе $\overline{\overline{e}}$.

Пусть ${}_*$-строка

$$\overline{\overline{f}} = \begin{pmatrix} {}_1f & ... & {}_nf \end{pmatrix} = ({}_if, i \in J)$$

является базисом $D^*{}_*$-векторного пространства $V$. Мы будем говорить, что координатная матрица $f$ множества векторов $({}_if, i \in J)$ определяет **координаты** ${}_if^j$ **базиса** $\overline{\overline{f}}$ относительно базиса $\overline{\overline{e}}$. □

Пример 5.3.8. Пусть $\overline{\overline{e}} = (e^i, i \in I, |I| = n)$ - **базис** ${}^*{}_*D$**-векторного пространства** $V$. Согласно примеру 5.2.4, мы можем представить базис $\overline{\overline{e}}$ в виде ${}^*$-строки (матрицы столбца)

$$\overline{\overline{e}} = \begin{pmatrix} e^1 \\ ... \\ e^n \end{pmatrix}$$

Координатная матрица

$$a = \begin{pmatrix} {}_1a & ... & {}_na \end{pmatrix} = ({}_ia, i \in I)$$

вектора $\overline{a}$ в базисе $\overline{\overline{e}}$ называется ${}^*{}_*D$**-вектором**[5.4] или $*D$**-вектор строкой**.

Пусть ${}^*$-строка

$$\overline{A} = \begin{pmatrix} \overline{A}^1 \\ ... \\ \overline{A}^m \end{pmatrix} = (\overline{A}_j, j \in J) \tag{5.3.7}$$

задаёт множество векторов. Векторы $\overline{A}_j$ имеют разложение

$$\overline{A}^j = e^*{}_*A^j$$

Если мы подставим координатные матрицы вектора $\overline{A}_j$ в матрицу (5.3.7), мы получим матрицу

$$A = \begin{pmatrix} \begin{pmatrix} {}_1A^1 & ... & {}_nA^1 \end{pmatrix} \\ ... \\ \begin{pmatrix} {}_1A^m & ... & {}_nA^m \end{pmatrix} \end{pmatrix} = \begin{pmatrix} {}_1A^1 & ... & {}_nA^1 \\ ... & ... & ... \\ {}_1A^m & ... & {}_nA^m \end{pmatrix} = ({}^iA_j)$$

Мы будем называть матрицу $A$ **координатной матрицей множества векторов** $(\overline{A}_j, j \in J)$ в базисе $\overline{\overline{e}}$ и её элементы **координатами множества векторов** $(\overline{A}_j, j \in J)$ в базисе $\overline{\overline{e}}$.

Пусть ${}^*$-строка

$$\overline{\overline{f}} = \begin{pmatrix} f^1 \\ ... \\ f^n \end{pmatrix} = (f^j, j \in I)$$

[5.4] ${}_*{}^*D$-вектор является аналогом вектор-столбца.

является базисом ${}^*{}_*D$-векторного пространства $V$. Мы будем говорить, что координатная матрица $f$ множества векторов $(f^j, j \in I)$ определяет **координаты** ${}_i f^j$ **базиса** $\overline{\overline{f}}$ относительно базиса $\overline{\overline{e}}$. □

Пример 5.3.9. Пусть $\overline{\overline{e}} = (e_i, i \in I, |I| = n)$ - **базис** ${}_*{}^*D$**-векторного пространства** $V$. Согласно примеру 5.2.5, мы можем представить базис $\overline{\overline{e}}$ в виде ${}_*$-строки (матрицы строки)

$$\overline{\overline{e}} = \begin{pmatrix} e_1 & ... & e_n \end{pmatrix}$$

Координатная матрица

$$a = \begin{pmatrix} {}^1a \\ ... \\ {}^na \end{pmatrix} = ({}^ia, i \in I)$$

вектора $\overline{a}$ в базисе $\overline{\overline{e}}$ называется ${}_*{}^*D$**-вектором**[5.5] или $D*$**-вектор столбцом**.

Пусть ${}_*$-строка

$$\overline{A} = \begin{pmatrix} \overline{A}_1 & ... & \overline{A}_m \end{pmatrix} = (\overline{A}_j, j \in J) \tag{5.3.8}$$

задаёт множество векторов. Векторы $\overline{A}_j$ имеют разложение

$$\overline{A}_j = e_*{}^*A_j$$

Если мы подставим координатные матрицы вектора $\overline{A}_j$ в матрицу (5.3.8), мы получим матрицу

$$A = \begin{pmatrix} \begin{pmatrix} {}^1A_1 \\ ... \\ {}^nA_1 \end{pmatrix} & ... & \begin{pmatrix} {}^1A_m \\ ... \\ {}^nA_m \end{pmatrix} \end{pmatrix} = \begin{pmatrix} {}^1A_1 & ... & {}^1A_m \\ ... & ... & ... \\ {}^nA_1 & ... & {}^nA_m \end{pmatrix} = ({}^iA_j)$$

Мы будем называть матрицу $A$ **координатной матрицей множества векторов** $(\overline{A}_j, j \in J)$ в базисе $\overline{\overline{e}}$ и её элементы **координатами множества векторов** $(\overline{A}_j, j \in J)$ в базисе $\overline{\overline{e}}$.

Пусть ${}_*$-строка

$$\overline{\overline{f}} = \begin{pmatrix} f_1 & ... & f_n \end{pmatrix} = (f_j, j \in I)$$

является базисом ${}_*{}^*D$-векторного пространства $V$. Мы будем говорить, что координатная матрица $f$ множества векторов $(f_j, j \in I)$ определяет **координаты** ${}^if_j$ **базиса** $\overline{\overline{f}}$ относительно базиса $\overline{\overline{e}}$. □

Так как мы линейную комбинацию выражаем с помощью матриц, мы можем распространить принцип двойственности на теорию векторных пространств. Мы можем записать принцип двойственности в одной из следующих форм

Теорема 5.3.10 (принцип двойственности). *Пусть* $\mathfrak{A}$ *- истинное утверждение о векторных пространствах. Если мы заменим одновременно*

- $D_*{}^*$*-вектор и* $D^*{}_*$*-вектор*
- ${}_*{}^*D$*-вектор и* ${}^*{}_*D$*-вектор*
- ${}_*{}^*$*-произведение и* ${}^*{}_*$*-произведение*

[5.5] ${}_*{}^*D$-вектор является аналогом вектор-столбца.

*то мы снова получим истинное утверждение.*

Теорема 5.3.11 (принцип двойственности). *Пусть $\mathfrak{A}$ - истинное утверждение о векторных пространствах. Если мы одновременно заменим*

- *$D_*{}^*$-вектор и ${}_*{}^*D$-вектор или $D^*{}_*$-вектор и ${}^*{}_*D$-вектор*
- *${}_*{}^*$-квазидетерминант и ${}^*{}_*$-квазидетерминант*

*то мы снова получим истинное утверждение.*

## 5.4. Линейное отображение ${}_*{}^*D$-векторных пространств

Определение 5.4.1. Пусть $V$ - ${}_*{}^*S$-векторное пространство. Пусть $U$ - ${}_*{}^*T$-векторное пространство. Мы будем называть морфизм

$$f : S \longrightarrow T \qquad A : V \longrightarrow U$$

правых представлений тела в абелевой группе **линейным отображением** ${}_*{}^*S$-векторного пространства $V$ в ${}_*{}^*T$-векторное пространство $U$. □

Согласно теореме 3.2.16 при изучении линейного отображения мы можем ограничиться случаем $S = T$.

Определение 5.4.2. Пусть $V$ и $W$ - ${}_*{}^*D$-векторные пространства. Мы будем называть отображение

$$A : V \to W$$

**линейным отображением** ${}_*{}^*D$-векторного пространства, если[5.6]

(5.4.1) $$A(m_*{}^*a) = A(m)_*{}^*a$$

для любых ${}^ia \in D$, $m_i \in V$. □

Теорема 5.4.3. *Пусть*

$$\overline{\overline{f}} = (f_i, i \in I)$$

*базис в ${}_*{}^*D$-векторном пространстве $V$ и*

$$\overline{\overline{e}} = (e_j, j \in J)$$

*базис в ${}_*{}^*D$-векторном пространстве $U$. Тогда линейное отображение*

(5.4.2) $$A : V \to W$$

*${}_*{}^*D$-векторных пространств имеет представление*

(5.4.3) $$b = A_*{}^*a$$

*относительно заданных базисов. Здесь*

- *$a$ - координатная матрица вектора $\overline{a}$ относительно базиса $\overline{\overline{f}}$.*
- *$b$ - координатная матрица вектора*
  $$\overline{b} = \overline{A}(\overline{a})$$
  *относительно базиса $\overline{\overline{e}}$.*
- *$A$ - координатная матрица множества векторов $(\overline{A}(\overline{f}_i))$ относительно базиса $\overline{\overline{e}}$. Мы будем называть матрицу $A$* **матрицей линейного отображения** *относительно базисов $\overline{\overline{f}}$ и $\overline{\overline{e}}$.*

[5.6] Выражение $A(m)_*{}^*a$ означает выражение $A(m_i)\ {}^ia$

Доказательство. Вектор $\overline{a} \in V$ имеет разложение

$$\overline{a} = f_*{}^*a$$

относительно базиса $\overline{\overline{f}}$. Вектор $\overline{b} \in U$ имеет разложение

$$\overline{b} = e_*{}^*b \tag{5.4.4}$$

относительно базиса $\overline{\overline{e}}$.

Так как $\overline{A}$ - линейное отображение, то на основании (5.4.1) следует, что

$$\overline{b} = \overline{A}(\overline{a}) = \overline{A}(f_*{}^*a) = \overline{A}(f)_*{}^*a \tag{5.4.5}$$

$\overline{A}(f_i)$ является вектором ${}_*{}^*D$-векторного пространства $U$ и имеет разложение

$$\overline{A}(\overline{f}_i) = \overline{e}_*{}^*A_i = \overline{e}_j{}^jA_i \tag{5.4.6}$$

относительно базиса $\overline{\overline{e}}$. Комбинируя (5.4.5) и (5.4.6), мы получаем

$$\overline{b} = e_*{}^*A_*{}^*a \tag{5.4.7}$$

(5.4.3) следует из сравнения (5.4.4) и (5.4.7) и теоремы 5.3.3. □

На основании теоремы 5.4.3 мы идентифицируем линейное отображение (5.4.2) ${}_*{}^*D$-векторных пространств и матрицу его представления (5.4.3).

Теорема 5.4.4. *Пусть*

$$\overline{\overline{f}} = (f_i, i \in I)$$

*базис в ${}_*{}^*D$-векторном пространстве $V$,*

$$\overline{\overline{e}} = (e_j, j \in J)$$

*базис в ${}_*{}^*D$-векторном пространстве $U$, и*

$$\overline{\overline{g}} = (g_l, l \in L)$$

*базис в ${}_*{}^*D$-векторном пространстве $W$. Предположим, что мы имеем коммутативную диаграмму отображений*

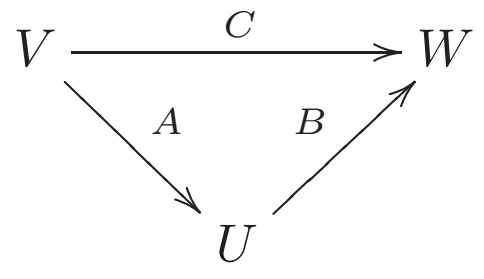

*где линейное отображение $A$ имеет представление*

$$b = A_*{}^*a \tag{5.4.8}$$

*относительно заданных базисов и линейное отображение $B$ имеет представление*

$$c = B_*{}^*b \tag{5.4.9}$$

*относительно заданных базисов. Тогда отображение $C$ является линейным и имеет представление*

$$c = B_*{}^*A_*{}^*a \tag{5.4.10}$$

*относительно заданных базисов.*

Доказательство. Отображение $C$ является линейным, так как

$$\begin{aligned}C_*{}^*(f_*{}^*a) &= (A_*{}^*B)_*{}^*(f_*{}^*a) = B_*{}^*(A_*{}^*(f_*{}^*a))\\ &= B_*{}^*(e_*{}^*(A_*{}^*a)) = g_*{}^*(B_*{}^*(A_*{}^*a))\\ &= g_*{}^*((B_*{}^*A)_*{}^*a) = g_*{}^*(C_*{}^*a)\end{aligned}$$

Равенство (5.4.10) следует из подстановки (5.4.8) в (5.4.9). □

Записывая линейное отображение в форме ${}_*{}^*$-произведения, мы можем переписать (5.4.1) в виде

$$\overline{A}_*{}^*(\overline{a}k) = (\overline{A}_*{}^*\overline{a})k \tag{5.4.11}$$

Утверждение теоремы 5.4.4 мы можем записать в виде

$$\overline{B}_*{}^*(\overline{A}_*{}^*\overline{a}) = (\overline{B}_*{}^*\overline{A})_*{}^*\overline{a} \tag{5.4.12}$$

Равенства (5.4.11) и (5.4.12) представляют собой **закон ассоциативности для линейных отображений ${}_*{}^*D$-векторных пространств**. Это позволяет нам писать подобные выражения не пользуясь скобками.

Равенство (5.4.3) является координатной записью линейного отображения. На основе теоремы 5.4.3 бескоординатная запись также может быть представлена с помощью ${}_*{}^*$-произведения

$$\overline{b} = \overline{A}_*{}^*\overline{a} = \overline{A}_*{}^*\overline{f}_*{}^*a = \overline{e}_*{}^*A_*{}^*a \tag{5.4.13}$$

Если подставить равенство (5.4.13) в теорему 5.4.4, то мы получим цепочку равенств

$$\overline{c} = \overline{B}_*{}^*\overline{b} = \overline{B}_*{}^*\overline{e}_*{}^*b = \overline{g}_*{}^*B_*{}^*b$$

$$\overline{c} = \overline{B}_*{}^*\overline{A}_*{}^*\overline{a} = \overline{B}_*{}^*\overline{A}_*{}^*\overline{f}_*{}^*a = \overline{g}_*{}^*B_*{}^*A_*{}^*a$$

Замечание 5.4.5. На примере линейных отображений легко видеть насколько теорема 3.2.16 облегчает наши рассуждения при изучении морфизма представлений $\Omega$-алгебры. Договоримся в рамках этого замечания теорию линейных отображений называть сокращённой теорией, а теорию, излагаемую в этом замечании, называть расширенной теорией.

Пусть $V$ - ${}_*{}^*S$-векторное пространство. Пусть $U$ - ${}_*{}^*T$-векторное пространство. Пусть

$$r : S \longrightarrow T \qquad \overline{A} : V \longrightarrow U$$

линейное отображение ${}_*{}^*S$-векторного пространства $V$ в ${}_*{}^*T$-векторное пространство $U$. Пусть

$$\overline{\overline{f}} = (f_i, i \in I)$$

базис в ${}_*{}^*S$-векторном пространстве $V$ и

$$\overline{\overline{e}} = (e_j, j \in J)$$

базис в ${}_*{}^*T$-векторном пространстве $U$.

Из определений 5.4.1 и 3.2.2 следует

$$\overline{b} = \overline{A}(\overline{a}) = \overline{A}(f_*{}^*a) = \overline{A}(f)_*{}^*r(a) \tag{5.4.14}$$

$\overline{A}(f_i)$ также вектор векторного пространства $U$ и имеет разложение

$$\overline{A}(\overline{f}_i) = \overline{e}_*{}^*A_i = \overline{e}_j{}^j A_i \tag{5.4.15}$$

относительно базиса $\overline{\overline{e}}$. Комбинируя (5.4.14) и (5.4.15), мы получаем

$$\overline{b}=e_*{}^*A_*{}^*r(a) \tag{5.4.16}$$

Пусть $W$ - ${}_*{}^*D$-векторное пространство. Пусть

$$p:T\longrightarrow D \qquad \overline{B}:U\longrightarrow W$$

линейное отображение ${}_*{}^*T$-векторного пространства $U$ в ${}_*{}^*D$-векторное пространство $W$. Пусть

$$\overline{\overline{g}}=(g_l,l\in L)$$

базис в ${}_*{}^*D$-векторном пространстве $W$. Тогда, согласно (5.4.16), произведение линейного отображения $(r,\overline{A})$ и линейного отображения $(p,\overline{B})$ имеет вид

$$\overline{c}=h_*{}^*B_*{}^*p(A)_*{}^*pr(a) \tag{5.4.17}$$

Сопоставление равенств (5.4.10) и (5.4.17) показывает насколько расширеная теория линейных отображений сложнее сокращённой теории.

При необходимости мы можем пользоваться расширеной теорией, но мы не получим новых результатов по сравнению со случаем сокращённой теорией. В то же время обилие деталей делает картину менее ясной и требует постоянного внимания. □

## 5.5. Система линейных уравнений

Определение 5.5.1. Пусть $V$ - ${}_*{}^*D$-векторное пространство и $\{A_i\in V,i\in I\}$ - множество векторов. **Линейная оболочка в ${}_*{}^*D$-векторном пространстве** - это множество $\mathrm{span}(A_i,i\in I)$ векторов, линейно зависимых от векторов $A_i$. □

Теорема 5.5.2. *Пусть $span(A_i,i\in I)$ - линейная оболочка в ${}_*{}^*D$-векторном пространстве $V$. Тогда $span(A_i,i\in I)$ - подпространство ${}_*{}^*D$-векторного пространства $V$.*

Доказательство. Предположим, что

$$\overline{b}\in\mathrm{span}(A_i,i\in I)$$
$$\overline{c}\in\mathrm{span}(A_i,i\in I)$$

Согласно определению 5.5.1

$$\overline{b}=A_*{}^*b$$
$$\overline{c}=A_*{}^*c$$

Тогда

$$\overline{b}+\overline{c}=A_*{}^*b+A_*{}^*c=A_*{}^*(b+c)\in\mathrm{span}(A_i,i\in I)$$
$$\overline{b}k=(A_*{}^*b)k=A_*{}^*(bk)\in\mathrm{span}(A_i,i\in I)$$

Это доказывает утверждение. □

Пример 5.5.3. Пусть $V$ - ${}_*{}^*D$-векторное пространство и ${}_*$-строка

$$\overline{\overline{A}}=\begin{pmatrix}\overline{A}_1 & \dots & \overline{A}_n\end{pmatrix}=(\overline{A}_i,i\in I)$$

задаёт множество векторов. Чтобы ответить на вопрос, или вектор $\overline{b}\in\mathrm{span}(\overline{A}_i, i\in I)$ , мы запишем линейное уравнение

$$\overline{b}=\overline{A}_*{}^*x \tag{5.5.1}$$

где

$$x=\begin{pmatrix}{}^1x\\ ...\\ {}^nx\end{pmatrix}$$

${}^*$-строка неизвестных коэффициентов разложения. $\overline{b}\in\mathrm{span}(\overline{A}_i, i\in I)$ , если уравнение (5.5.1) имеет решение. Предположим, что $\overline{\overline{f}}=(f_j, j\in J)$ - базис. Тогда векторы $\overline{b}$, $\overline{A}_i$ имеют разложение

$$\overline{b}=f_*{}^*b \tag{5.5.2}$$

$$\overline{A}_i=f_*{}^*A_i \tag{5.5.3}$$

Если мы подставим (5.5.2) и (5.5.3) в (5.5.1), мы получим

$$f_*{}^*b=f_*{}^*A_*{}^*x \tag{5.5.4}$$

Применяя теорему 5.3.3 к (5.5.4), мы получим **систему линейных уравнений**

$$A_*{}^*x=b \tag{5.5.5}$$

Мы можем записать систему линейных уравнений (5.5.5) в одной из следующих форм

$$\begin{pmatrix}{}^1A_1 & ... & {}^1A_n\\ ... & ... & ...\\ {}^mA_1 & ... & {}^mA_n\end{pmatrix}{}_*{}^*\begin{pmatrix}{}^1x\\ ...\\ {}^nx\end{pmatrix}=\begin{pmatrix}{}^1b\\ ...\\ {}^mb\end{pmatrix}$$

$${}^jA_i\ {}^ix={}^jb \tag{5.5.6}$$

$$\begin{array}{cccc}{}^1A_1\ {}^1x & +... & +{}^1A_n\ {}^nx & ={}^1b\\ ... & ... & ... & ...\\ {}^1A_m\ {}^1x & +... & +{}^mA_n\ {}^nx & ={}^mb\end{array}$$

□

Пример 5.5.4. Пусть $V$ - $D_*{}^*$-векторное пространство и ${}^*$-строка

$$\overline{\overline{A}}=\begin{pmatrix}{}^1\overline{A}\\ ...\\ {}^m\overline{A}\end{pmatrix}=({}^j\overline{A}, j\in J)$$

задаёт множество векторов. Чтобы ответить на вопрос, или вектор $\overline{b}\in\mathrm{span}({}^jA, j\in J)$ , мы запишем линейное уравнение

$$\overline{b}=x_*{}^*\overline{A} \tag{5.5.7}$$

где

$$x=\begin{pmatrix} x_1 & ... & x_m \end{pmatrix}$$

${}_*$-строка неизвестных коэффициентов разложения. $\overline{b}\in \mathrm{span}({}^jA, j\in J)$ , если уравнение (5.5.7) имеет решение. Предположим, что $\overline{\overline{f}}=({}^i\overline{f}, i\in I)$ - базис. Тогда векторы $\overline{b}$, ${}^jA$ имеют разложение

$$\overline{b}=b_*{}^*f \tag{5.5.8}$$

$$ {}^jA={}^jA_*{}^*f \tag{5.5.9}$$

Если мы подставим (5.5.8) и (5.5.9) в (5.5.7), мы получим

$$b_*{}^*f=x_*{}^*A_*{}^*f \tag{5.5.10}$$

Применяя теорему 5.3.3 к (5.5.10), мы получим **система линейных уравнений**[5.7]

$$x_*{}^*A=b \tag{5.5.11}$$

Мы можем записать систему линейных уравнений (5.5.11) в одной из следующих форм

$$\begin{pmatrix} x_1 & ... & x_m \end{pmatrix}{}_*{}^*\begin{pmatrix} {}^1A_1 & ... & {}^1A_n\\ ... & ... & ...\\ {}^mA_1 & ... & {}^mA_n \end{pmatrix}=\begin{pmatrix} b_1 & ... & b_n \end{pmatrix}$$

$$x_j\ {}^jA_i=b_i \tag{5.5.12}$$

$$\begin{matrix} x_1\ {}^1A_1 & ... & x_1\ {}^1A_n\\ + & ... & +\\ ... & ... & ...\\ + & ... & +\\ x_m\ {}^mA_1 & ... & x_m\ {}^mA_n\\ = & ... & =\\ b_1 & ... & b_n \end{matrix}$$

□

Чтобы найти решение системы линейных уравнений, мы должны рассмотреть матрицу этой системы. Из примеров 5.3.7, 5.3.9, мы видим, что мы можем рассматривать столбец матрицы как вектор левого или правого векторного пространства. Чтобы сделать утверждения проще, мы будем указывать тип векторного пространства перед словом линейный. Например, утверждение

Столбцы матрицы $D^*{}_*$-линейно зависимы.

означает, что

[5.7] Читая систему ${}_*{}^*D$-линейных уравнений (5.5.5) в ${}_*{}^*D$-векторном пространстве снизу вверх и слева направо, мы получим систему линейных уравнений (5.5.11) в $D_*{}^*$-векторном пространстве.

Столбцы матрицы являются векторами $D^*{}_*$-векторного пространства, и соответствующие векторы линейно зависимы.

В частности, система линейных уравнений (5.5.5) в ${}_*{}^*D$-векторном пространстве называется **системой ${}_*{}^*D$-линейных уравнений** и система линейных уравнений (5.5.11) в $D_*{}^*$-векторном пространстве называется **системой $D_*{}^*$-линейных уравнений**.

Определение 5.5.5. Если $n\times n$ матрица $A$ имеет ${}_*{}^*$-обратную матрицу, мы будем называть такую матрицу **${}_*{}^*$-невырожденной матрицей**. В противном случае, мы будем называть такую матрицу **${}_*{}^*$-вырожденной матрицей**. □

Определение 5.5.6. Предположим, что $A$ - ${}_*{}^*$-невырожденная матрица. Мы будем называть соответствующую систему ${}_*{}^*D$-линейных уравнений

$$A_*{}^*x=b \tag{5.5.13}$$

**невырожденной системой ${}_*{}^*D$-линейных уравнений**. □

Теорема 5.5.7. *Решение невырожденной системы ${}_*{}^*D$-линейных уравнений (5.5.13) определено однозначно и может быть записано в любой из следующих форм*[5.8]

$$x=A^{-1_*{}^*}{}_*{}^*b \tag{5.5.14}$$

$$x=\mathcal{H}\det({}_*{}^*)A_*{}^*b \tag{5.5.15}$$

Доказательство. Умножая обе части равенства (5.5.13) слева на $A^{-1_*{}^*}$, мы получим (5.5.14). Пользуясь определением (2.3.12), мы получим (5.5.15). Решение системы единственно в силу теоремы 2.2.16. □

## 5.6. Ранг матрицы

Определение 5.6.1. Мы будем называть матрицу[5.9] ${}^SA_T$ минорной матрицей порядка $k$. □

Определение 5.6.2. Если минорная матрица ${}^SA_T$ - ${}_*{}^*$-невырожденная матрица, то мы будем говорить, что ${}_*{}^*$-ранг матрицы $A$ не меньше, чем $k$. **${}_*{}^*$-ранг матрицы $A$**

$$\mathrm{rank}_{{}_*{}^*}A$$

- это максимальное значение $k$. Мы будем называть соответствующую минорную матрицу **${}_*{}^*$-главной минорной матрицей**. □

Теорема 5.6.3. *Пусть матрица $A$ - ${}_*{}^*$-вырожденная матрица и минорная матрица ${}^SA_T$ - главный минорная матрица, тогда*

$$ {}^p\det({}_*{}^*)_r\,{}^{S\cup\{p\}}A_{T\cup\{r\}}=0 \tag{5.6.1}$$

[5.8] Мы можем найти решение системы (5.5.13) в теореме [6]-1.6.1. Я повторяю это утверждение, так как я слегка изменил обозначения.

[5.9] Мы делаем следующие предположения в этом разделе

- $i\in M$, $|M|=m$, $j\in N$, $|N|=n$.
- $A=({}^iA_j)$ - произвольная матрица.
- $k$, $s\in S\supseteq M$, $l$, $t\in T\supseteq N$, $k=|S|=|T|$.
- $p\in M\setminus S$, $r\in N\setminus T$.

Доказательство. Чтобы понять, почему минорная матрица

$$B = {}^{S\cup\{p\}}A_{T\cup\{r\}} \tag{5.6.2}$$

не имеет ${}_*{}^*$-обратной матрицы,[5.10] мы предположим, что существует ${}_*{}^*$-обратная матрица $B^{-1_*{}^*}$. Запишем систему линейных уравнений (2.3.2), (2.3.3) полагая $I=\{r\}$, $J=\{p\}$ (в этом случае $[I]=T$, $[J]=S$ )

$${}^SB_T{}_*{}^*{}^TB^{-1_*{}^*}{}_p + {}^SB_r\ {}^rB^{-1_*{}^*}{}_p = 0 \tag{5.6.3}$$

$${}^pB_T{}_*{}^*{}^TB^{-1_*{}^*}{}_p + {}^pB_r\ {}^rB^{-1_*{}^*}{}_p = 1 \tag{5.6.4}$$

и попробуем решить эту систему. Мы умножим (5.6.3) на $({}^SB_T)^{-1_*{}^*}$

$${}^TB^{-1_*{}^*}{}_p + ({}^SB_T)^{-1_*{}^*}{}_*{}^*{}^SB_r\ {}^rB^{-1_*{}^*}{}_p = 0 \tag{5.6.5}$$

Теперь мы можем подставить (5.6.5) в (5.6.4)

$$-\,{}^pB_T{}_*{}^*({}^SB_T)^{-1_*{}^*}{}_*{}^*{}^SB_r\ {}^rB^{-1_*{}^*}{}_p + {}^pB_r\ {}^rB^{-1_*{}^*}{}_p = 1 \tag{5.6.6}$$

Из (5.6.6) следует

$$({}^pB_r - {}^pB_T{}_*{}^*({}^SB_T)^{-1_*{}^*}{}_*{}^*{}^SB_r)\ {}^rB^{-1_*{}^*}{}_p = 1 \tag{5.6.7}$$

Выражение в скобках является квазидетерминантом ${}^p\det({}_*{}^*)_rB$. Подставляя это выражение в (5.6.7), мы получим

$${}^p\det({}_*{}^*)_rB\ {}^rB^{-1_*{}^*}{}_p = 1 \tag{5.6.8}$$

Тем самым мы доказали, что квазидетерминант ${}^p\det({}_*{}^*)_rB$ определён и условие его обращения в 0 необходимое и достаточное условие вырожденности матрицы $B$. Теорема доказана в силу соглашения (5.6.2). □

Теорема 5.6.4. *Предположим, что $A$ - матрица*[5.9]*,*

$$\mathrm{rank}_{{}_*{}^*}A = k < m$$

*и ${}^SA_T$ - ${}_*{}^*$-главная минорная матрица. Тогда ${}_*$-строка ${}^pA$ является $D_*{}^*$-линейной комбинацией ${}_*$-строк ${}^SA$.*

$${}^{M\setminus S}A = R_*{}^*{}^SA \tag{5.6.9}$$

$${}^pA = {}^pR_*{}^*{}^SA \tag{5.6.10}$$

$${}^pA_b = {}^pR_s\ {}^sA_b \tag{5.6.11}$$

Доказательство. Если матрица $A$ имеет $k$ ${}^*$-строк, то, полагая, что ${}_*$-строка ${}^pA$ - $D_*{}^*$-линейная комбинация (5.6.10) ${}_*$-строк ${}^sA$ с коэффициентами ${}^pR_s$, мы получим систему $D_*{}^*$-линейных уравнений (5.6.11). Согласно теореме 5.5.7 система $D_*{}^*$-линейных уравнений (5.6.11) имеет единственное решение[5.11] и это решение нетривиально потому, что все ${}_*{}^*$-квазидетерминанты отличны от 0.

---

[5.10] Естественно ожидать связь между ${}_*{}^*$-вырожденностью матрицы и её ${}_*{}^*$-квазидетерминантом, подобную связи, известной в коммутативном случае. Однако ${}_*{}^*$-квазидетерминант определён не всегда. Например, если ${}_*{}^*$-обратная матрица имеет слишком много элементов, равных 0. Как следует из этой теоремы, не определён ${}_*{}^*$-квазидетерминант также в случае, когда ${}_*{}^*$-ранг матрицы меньше $n-1$.

[5.11] Мы положим, что неизвестные переменные здесь - это $x_s = {}^pR_s$

Остаётся доказать утверждение в случае, когда число ${}^*$-строк матрицы $A$ больше чем $k$. Пусть нам даны ${}_*$-строка ${}^pA$ и ${}^*$-строка $A_r$. Согласно предположению, минорная матрица ${}^{S\cup\{p\}}A_{T\cup\{r\}}$ - ${}_*{}^*$-вырожденная матрица и его ${}_*{}^*$-квазидетерминант

$$ {}^p\det({}_*{}^*)_r\, {}^{S\cup\{p\}}A_{T\cup\{r\}} = 0 \tag{5.6.12}$$

Согласно (2.3.14) равенство (5.6.12) имеет вид

$$ {}^pA_r - {}^pA_{T*}{}^*(({}^SA_T)^{-1_*{}^*})_*{}^{*S}A_r = 0 $$

Матрица

$$ {}^pR = {}^pA_{T*}{}^*(({}^SA_T)^{-1_*{}^*}) \tag{5.6.13}$$

не зависит от $r$. Следовательно, для любых $r \in N \setminus T$

$$ {}^pA_r = {}^pR_*{}^{*S}A_r \tag{5.6.14}$$

Из равенства

$$ (({}^SA_T)^{-1_*{}^*})_*{}^{*S}A_l = {}^T\delta_l $$

следует, что

$$ {}^pA_l = {}^pA_{T*}{}^{*T}\delta_l = {}^pA_{T*}{}^*(({}^SA_T)^{-1_*{}^*})_*{}^{*S}A_l \tag{5.6.15}$$

Подставляя (5.6.13) в (5.6.15), мы получим

$$ {}^pA_l = {}^pR_*{}^{*S}A_l \tag{5.6.16}$$

(5.6.14) и (5.6.16) завершают доказательство. □

Следствие 5.6.5. *Предположим, что $A$ - матрица, $\mathrm{rank}_{*}{}^* A = k < m$. Тогда ${}_*$-строки матрицы $D_*{}^*$-линейно зависимы.*

$$ \lambda_*{}^*A = 0 $$

Доказательство. Предположим, что ${}_*$-строка ${}^pA$ - $D_*{}^*$-линейная комбинация (5.6.10). Мы положим $\lambda_p = -1$, $\lambda_s = {}^pR_s$ и остальные $\lambda_c = 0$. □

Теорема 5.6.6. *Пусть $({}^i\overline{A}, i \in M, |M| = m)$ семейство $D_*{}^*$-линейно независимых векторов. Тогда ${}_*{}^*$-ранг их координатной матрицы равен $m$.*

Доказательство. Пусть $\overline{\overline{e}}$ - базис $D_*{}^*$-векторного пространства. Согласно конструкции, изложенной в примере 5.3.6, координатная матрица семейства векторов $({}^i\overline{A})$ относительно базиса $\overline{\overline{e}}$ состоит из ${}_*$-строк, являющихся координатными матрицами векторов ${}^i\overline{A}$ относительно базиса $\overline{\overline{e}}$. Поэтому ${}_*{}^*$-ранг этой матрицы не может превышать $m$.

Допустим ${}_*{}^*$-ранг координатной матрицы меньше $m$. Согласно следствию 5.6.5, ${}_*$-строки матрицы $D_*{}^*$-линейно зависимы

$$ \lambda_*{}^*A = 0 \tag{5.6.17}$$

Положим $c = \lambda_*{}^*A$. Из равенства (5.6.17) следует, что $D_*{}^*$-линейная комбинация

$$ c_*{}^*\overline{e} = 0 $$

векторов базиса равна 0. Это противоречит утверждению, что векторы $\overline{e}$ образуют базис. Утверждение теоремы доказано. □

Теорема 5.6.7. *Предположим, что $A$ - матрица*[5.9],

$$\mathrm{rank}_{*}{}^{*} A = k < n$$

*и* ${}^S A_T$ - ${}_*{}^*$*-главная минорная матрица. Тогда* ${}^*$*-строка* $A_r$ *является* ${}_*{}^*D$*-линейной композицией* ${}^*$*-строк* $A_t$

$$A_{N\setminus T} = A_T{}_*{}^*R \tag{5.6.18}$$

$$A_r = A_T{}_*{}^*R_r \tag{5.6.19}$$

$$^a A_r = {}^a A_t\ {}^t R_r \tag{5.6.20}$$

Доказательство. Если матрица $A$ имеет $k$ ${}_*$-строк, то, полагая, что ${}^*$-строка $A_r$ - ${}_*{}^*D$-линейная комбинация (5.6.19) ${}^*$-строк $A_t$ с коэффициентами ${}^t R_r$, мы получим систему ${}_*{}^*D$-линейных уравнений (5.6.20). Согласно теореме 5.5.7 система ${}_*{}^*D$-линейных уравнений (5.6.20) имеет единственное решение[5.12] и это решение нетривиально потому, что все ${}_*{}^*$-квазидетерминанты отличны от 0.

Остаётся доказать утверждение в случае, когда число ${}_*$-строк матрицы $A$ больше чем $k$. Пусть нам даны ${}^*$-строка $A_r$ и ${}_*$-строка ${}^p A$. Согласно предположению, минорная матрица ${}^{S\cup\{p\}}A_{T\cup\{r\}}$ - ${}_*{}^*$-вырожденная матрица и его ${}_*{}^*$-квазидетерминант

$$^p \det({}_*{}^*)_r\ {}^{S\cup\{p\}}A_{T\cup\{r\}} = 0 \tag{5.6.21}$$

Согласно (2.3.14) (5.6.21) имеет вид

$$^p A_r - {}^p A_T{}_*{}^*(({}^S A_T)^{-1_*{}^*})_*{}^*{}^S A_r = 0$$

Матрица

$$R_r = (({}^S A_T)^{-1_*{}^*})_*{}^*{}^S A_r \tag{5.6.22}$$

не зависит от $p$, Следовательно, для любых $p \in M \setminus S$

$$^p A_r = {}^p A_T{}_*{}^*R_r \tag{5.6.23}$$

Из равенства

$$^k A_T{}_*{}^*(({}^S A_T)^{-1_*{}^*})_s = {}^k\delta_s$$

следует, что

$$^k A_r = {}^k\delta_S{}_*{}^*{}^S A_r = {}^k A_T{}_*{}^*(({}^S A_T)^{-1_*{}^*})^s{}_*{}^*{}^S A_r \tag{5.6.24}$$

Подставляя (5.6.22) в (5.6.24), мы получим

$$^k A_r = {}^k A_T{}_*{}^*R_r \tag{5.6.25}$$

(5.6.23) и (5.6.25) завершают доказательство. □

Следствие 5.6.8. *Предположим, что $A$ - матрица,* $\mathrm{rank}_{*}{}^{*} A = k < m$. *Тогда* ${}^*$*-строки матрицы* ${}_*{}^*D$*-линейно зависимы.*

$$A_*{}^*\lambda = 0$$

Доказательство. Предположим, что ${}^*$-строка $A_r$ - правая линейная комбинация (5.6.19). Мы положим ${}^r\lambda = -1$, ${}^t\lambda = {}^t R_r$ и остальные ${}^c\lambda = 0$. □

Опираясь на теорему 2.3.9, мы можем записать подобные утверждения для ${}^*{}_*$-ранга матрицы.

[5.12]Мы положим, что неизвестные переменные здесь - это ${}^t x = {}^t R_r$

Теорема 5.6.9. *Предположим, что $A$ - матрица,*

$$\mathrm{rank}_{*}{}_{*}\ A=k<m$$

*и ${}_TA^S$ - ${}^*{}_*$-главная минорная матрица. Тогда ${}_*$-строка $A^p$ является ${}^*{}_*D$-линейной композицией ${}_*$-строк $A^s$*

$$A^{M\setminus S}=A^S{}^*{}_*R \tag{5.6.26}$$

$$A^p=A^S{}^*{}_*R^p \tag{5.6.27}$$

$$ {}_bA^p={}_bA^s\ {}_sR^p \tag{5.6.28}$$

Следствие 5.6.10. *Предположим, что $A$ - матрица,* $\mathrm{rank}_{*}{}_{*}\ A=k<m$. *Тогда ${}_*$-строки матрицы ${}^*{}_*D$-линейно зависимы.*

$$A^*{}_*\lambda=0$$

Теорема 5.6.11. *Предположим, что $A$ - матрица,*

$$\mathrm{rank}_{*}{}_{*}\ A=k<n$$

*и ${}_TA^S$ - ${}^*{}_*$-главная минорная матрица. Тогда ${}^*$-строка ${}_rA$ является $D^*{}_*$-линейной композицией ${}^*$-строк ${}_tA$*

$$ {}_{N\setminus T}A=R^*{}_*{}_TA \tag{5.6.29}$$

$$ {}_rA={}_rR^*{}_*{}_TA \tag{5.6.30}$$

$$ {}_rA^a={}_rR^t\ {}_tA^a \tag{5.6.31}$$

Следствие 5.6.12. *Предположим, что $A$ - матрица,* $\mathrm{rank}_{*}{}_{*}\ a=k<m$. *Тогда ${}^*$-строки матрицы $D^*{}_*$-линейно зависимы.*

$$\lambda^*{}_*A=0$$

## 5.7. Система ${}_*{}^*D$-линейных уравнений

Определение 5.7.1. Предположим, что[5.9] $A$ - матрица системы $D_*{}^*$-линейных уравнений (5.5.12). Мы будем называть матрицу

$$\begin{pmatrix}{}^jA_i\\ b_i\end{pmatrix}=\begin{pmatrix}{}^1A_1&...&{}^1A_n\\ ...&...&...\\ {}^mA_1&...&{}^mA_n\\ b_1&...&b_n\end{pmatrix} \tag{5.7.1}$$

**расширенной матрицей** этой системы. □

Определение 5.7.2. Предположим, что[5.9] $A$ - матрица системы ${}_*{}^*D$-линейных уравнений (5.5.6). Мы будем называть матрицу

$$\begin{pmatrix}{}^jA_i&{}^jb\end{pmatrix}=\begin{pmatrix}{}^1A_1&...&{}^1A_n&{}^1b\\ ...&...&...&...\\ {}^mA_1&...&{}^mA_n&{}^mb\end{pmatrix} \tag{5.7.2}$$

**расширенной матрицей** этой системы. □

Теорема 5.7.3. *Система ${}_*{}^*D$-линейных уравнений* (5.5.6) *имеет решение тогда и только тогда, когда*

$$\mathrm{rank}_{*^*}({}^jA_i) = \mathrm{rank}_{*^*}\begin{pmatrix}{}^jA_i & {}^jb\end{pmatrix} \tag{5.7.3}$$

Доказательство. Пусть ${}^SA_T$ - ${}_*{}^*$-главная минорная матрица матрицы $A$.

Пусть система ${}_*{}^*D$-линейных уравнений (5.5.6) имеет решение ${}^ix = {}^id$. Тогда

$$A_*{}^*d = b \tag{5.7.4}$$

Уравнение (5.7.4) может быть записано в форме

$$A_{T*}{}^{*T}d + A_{N\setminus T*}{}^{*N\setminus T}d = b \tag{5.7.5}$$

Подставляя (5.6.18) в (5.7.5), мы получим

$$A_{T*}{}^{*T}d + A_{T*}{}^*R_*{}^{*N\setminus T}d = b \tag{5.7.6}$$

Из (5.7.6) следует, что ${}_*$-строка $b$ является ${}_*{}^*D$-линейной комбинацией ${}_*$-строк $A_T$

$$A_{T*}{}^*({}^Td + R_*{}^{*N\setminus T}d) = b$$

Это эквивалентно уравнению (5.7.3).

Нам осталось доказать, что существование решения ${}_*{}^*D$-системы линейных уравнений (5.5.6) следует из (5.7.3). Истинность (5.7.3) означает, что ${}^SA_T$ - так же ${}_*{}^*$-главная минорная матрица расширенной матрицы. Из теоремы 5.6.7 следует, что ${}_*$-строка $b$ является ${}_*{}^*D$-линейной композицией ${}_*$-строк $A_T$

$$b = A_{T*}{}^{*T}R$$

Полагая ${}^rR = 0$ , мы получим

$$b = A_*{}^*R$$

Следовательно, мы нашли, по крайней мере, одно решение системы ${}_*{}^*D$-линейных уравнений (5.5.6). □

Теорема 5.7.4. *Предположим, что* (5.5.6) *- система ${}_*{}^*D$-линейных уравнений, удовлетворяющих* (5.7.3). *Если* $\mathrm{rank}_{*^*} A = k \le m$, *то решение системы зависит от произвольных значений $m-k$ переменных, не включённых в ${}_*{}^*$-главную минорную матрицу.*

Доказательство. Пусть ${}^SA_T$ - ${}_*{}^*$-главную минорную матрицу матрицы $a$. Предположим, что

$${}^pA_*{}^*x = {}^pb \tag{5.7.7}$$

уравнение с номером $p$. Применяя теорему 5.6.4 к расширенной матрице (5.7.1), мы получим

$${}^pA = {}^pR_*{}^{*S}A \tag{5.7.8}$$

$${}^pb = {}^pR_*{}^{*S}b \tag{5.7.9}$$

Подставляя (5.7.8) и (5.7.9) в (5.7.7), мы получим

$${}^pR_*{}^{*S}A_*{}^*x = {}^pR_*{}^{*S}b \tag{5.7.10}$$

(5.7.10) означает, что мы можем исключить уравнение (5.7.7) из системы (5.5.6) и новая система эквивалентна старой. Следовательно, число уравнений может быть уменьшено до $k$.

В этом случае у нас есть два варианта. Если число переменных также равно $k$, то согласно теореме 5.5.7 система имеет единственное решение (5.5.15). Если число переменных $m > k$, то мы можем передвинуть $m-k$ переменных, которые не включены в ${}_*{}^*$-главную минорную матрицу в правой части. Присваивая произвольные значения этим переменным, мы определяем значение правой части и для этого значения мы получим единственное решение согласно теореме 5.5.7. □

Следствие 5.7.5. *Система ${}_*{}^*D$-линейных уравнений (5.5.6) имеет единственное решение тогда и только тогда, когда её матрица невырожденная.* □

Теорема 5.7.6. *Решения однородной системы ${}_*{}^*D$-линейных уравнений*

$$A_*{}^*x = 0 \tag{5.7.11}$$

*порождают ${}_*{}^*D$-векторное пространство.*

Доказательство. Пусть $\overline{X}$ - множество решений системы ${}_*{}^*D$-линейных уравнений (5.7.11). Предположим, что $x = ({}^ax) \in \overline{X}$ и $y = ({}^ay) \in \overline{X}$. Тогда

$$x^a\ {}_aa^b = 0$$

$$y^a\ {}_aa^b = 0$$

Следовательно,

$${}^iA_j\ ({}^jx + {}^jy) = {}^iA_j\ {}^jx + {}^iA_j\ {}^jy = 0$$

$$x + y = (\ {}^jx + {}^jy) \in \overline{X}$$

Таким же образом мы видим

$${}^iA_j\ ({}^jxb) = ({}^iA_j\ {}^jx)b = 0$$

$$xb = ({}^jxb) \in \overline{X}$$

Согласно определению 5.1.4 $\overline{X}$ является ${}_*{}^*D$-векторным пространством. □

## 5.8. Невырожденная матрица

Предположим, что нам дана $n\times n$ матрица $A$. Следствия 5.6.5 и 5.6.8 говорят нам, что если $\operatorname{rank}_{{}_*{}^*} A < n$, то ${}_*$-строки $D_*{}^*$-линейно зависимы и ${}^*$-строки ${}_*{}^*D$-линейно зависимы.[5.13]

Теорема 5.8.1. *Пусть $A$ - $n\times n$ матрица и ${}^*$-строка $A_r$ - ${}_*{}^*D$-линейная комбинация других ${}^*$-строк. Тогда $\operatorname{rank}_{{}_*{}^*} A < n$.*

Доказательство. Утверждение, что ${}^*$-строка $A_r$ является ${}_*{}^*D$-линейной комбинацией других ${}^*$-строк, означает, что система ${}_*{}^*D$-линейных уравнений

$$A_r = A_{[r]}{}_*{}^*\lambda$$

имеет, по крайней мере, одно решение. Согласно теореме 5.7.3

$$\operatorname{rank}_{{}_*{}^*} A = \operatorname{rank}_{{}_*{}^*} A_{[r]}$$

Так как число ${}^*$-строк меньше, чем $n$, то $\operatorname{rank}_{{}_*{}^*} A[r] < n$. □

[5.13] Это утверждение похоже на утверждение [7]-1.2.5.

Теорема 5.8.2. *Пусть* $A$ - $n\times n$ *матрица и* ${}_*$*-строка* ${}^pA$ - $D_*{}^*$*-линейная комбинация других* ${}_*$*-строк. Тогда* $\mathrm{rank}_{*}{}^{*} A < n$.

Доказательство. Доказательство утверждения похоже на доказательство теоремы 5.8.1 □

Теорема 5.8.3. *Предположим, что* $A$ *и* $B$ - $n\times n$ *матрицы и*

$$C = A_*{}^*B \tag{5.8.1}$$

$C$ - ${}_*{}^*$*-вырожденная матрица тогда и только тогда, когда либо матрица* $A$*, либо матрица* $B$ - ${}_*{}^*$*-вырожденная матрица.*

Доказательство. Предположим, что матрица $B$ - ${}_*{}^*$-вырожденная. Согласно теореме 5.6.7 ${}^*$-строки матрицы $B$ ${}_*{}^*D$-линейно зависимы. Следовательно,

$$0 = B_*{}^*\lambda \tag{5.8.2}$$

где $\lambda\neq 0$. Из (5.8.1) и (5.8.2) следует, что

$$C_*{}^*\lambda = A_*{}^*B_*{}^*\lambda = 0$$

Согласно теореме 5.8.1 матрица $C$ является ${}_*{}^*$-вырожденной.

Предположим, что матрица $B$ - не ${}_*{}^*$-вырожденная, но матрица $A$ - ${}_*{}^*$-вырожденная. Согласно теореме 5.6.4 ${}^*$-строки матрицы $A$ ${}_*{}^*D$-линейно зависимы. Следовательно,

$$0 = A_*{}^*\mu \tag{5.8.3}$$

где $\mu\neq 0$. Согласно теореме 5.5.7 система

$$B_*{}^*\lambda = \mu$$

имеет единственное решение, где $\lambda\neq 0$. Следовательно,

$$C_*{}^*\lambda = A_*{}^*B_*{}^*\lambda = A_*{}^*\mu = 0$$

Согласно теореме 5.8.1 матрица $C$ является ${}_*{}^*$-вырожденной.

Предположим, что матрица $C$ - ${}_*{}^*$-вырожденная матрица. Согласно теореме 5.6.4 ${}_*$-строки матрицы $C$ ${}_*{}^*D$-линейно зависимы. Следовательно,

$$0 = C_*{}^*\lambda \tag{5.8.4}$$

где $\lambda\neq 0$. Из (5.8.1) и (5.8.4) следует, что

$$0 = A_*{}^*B_*{}^*\lambda$$

Если

$$0 = B_*{}^*\lambda$$

выполнено, то матрица $B$ - ${}_*{}^*$-вырожденная. Предположим, что матрица $B$ не ${}_*{}^*$-вырожденная. Положим

$$\mu = B_*{}^*\lambda$$

где $\mu\neq 0$. Тогда

$$0 = A_*{}^*\mu \tag{5.8.5}$$

Из (5.8.5) следует, что матрица $A$ - ${}_*{}^*$-вырожденная. □

Опираясь на теорему 2.3.9, мы можем записать подобные утверждения для $D^*{}_*$-линейной комбинации ${}_*$-строк или ${}^*{}_*D$-линейной комбинации ${}^*$-строк и ${}^*{}_*$-квазидетерминанта.

Теорема 5.8.4. *Пусть $A$ - $n\times n$ матрица и ${}^*$-строка ${}_rA$ - ${}^*{}_*D$-линейная комбинация других ${}^*$-строк. Тогда* $\mathrm{rank}_{{}^*{}_*} A < n$.

Теорема 5.8.5. *Пусть $A$ - $n\times n$ матрица и ${}_*$-строка $A^p$ - $D^*{}_*$-линейная комбинация других ${}_*$-строк. Тогда* $\mathrm{rank}_{{}^*{}_*} A < n$.

Теорема 5.8.6. *Предположим, что $A$ и $B$ - $n\times n$ матрицы и $C = A^*{}_*B$. $C$ - ${}^*{}_*$-вырожденная матрица тогда и только тогда, когда либо матрица $A$, либо матрица $B$ - ${}^*{}_*$-вырожденная матрица.*

Определение 5.8.7. ${}_*{}^*$**-матричная группа** $GL(n, {}_*{}^*, D)$ - это группа ${}_*{}^*$-невырожденных матриц, где мы определяем ${}_*{}^*$-произведение матриц (2.2.1) и ${}_*{}^*$-обратную матрица $A^{-1^*{}_*}$. □

Определение 5.8.8. ${}^*{}_*$**-матричная группа** $GL(n, {}^*{}_*, D)$ - это группа ${}^*{}_*$-невырожденных матриц где мы определяем ${}^*{}_*$-произведение матриц (2.2.2) и ${}^*{}_*$-обратную матрицу $A^{-1^*{}_*}$. □

Теорема 5.8.9.

$$GL(n, {}_*{}^*, D) \ne GL(n, {}^*{}_*, D)$$

Замечание 5.8.10. Из теоремы 2.3.11 следует, что существуют матрицы, которые ${}^*{}_*$-невырожденны и ${}_*{}^*$-невырожденны. Теорема 5.8.9 означает, что множества ${}^*{}_*$-невырожденных матриц и ${}_*{}^*$-невырожденных матриц не совпадают. Например, существует такая ${}_*{}^*$-невырожденная матрица, которая ${}^*{}_*$-вырожденная матрица. □

Доказательство. Это утверждение достаточно доказать для $n = 2$. Предположим, что каждая ${}^*{}_*$-вырожденная матрица

$$A = \begin{pmatrix} {}^1A_1 & {}^2A_1 \\ {}^1A_2 & {}^2A_2 \end{pmatrix} \tag{5.8.6}$$

является ${}_*{}^*$-вырожденной матрицей. Из теоремы 5.6.7 и теоремы 5.6.7 следует, что ${}_*{}^*$-вырожденная матрица удовлетворяет условию

$$ {}^1A_2 = b\ {}^1A_1 \tag{5.8.7}$$

$$ {}^2A_2 = b\ {}^2A_1 \tag{5.8.8}$$

$$ {}^2A_1 = {}^1A_1c \tag{5.8.9}$$

$$ {}^2A_2 = {}^1A_2c \tag{5.8.10}$$

Если мы подставим (5.8.9) в (5.8.8), мы получим

$$ {}^2A_2 = b\ {}^1A_1\ c$$

$b$ и $c$ - произвольные элементы тела $D$ и ${}_*{}^*$-вырожденная матрица (5.8.6) имеет вид ($d = {}^1A_1$)

$$A = \begin{pmatrix} d & dc \\ bd & bdc \end{pmatrix} \tag{5.8.11}$$

Подобным образом мы можем показать, что ${}^*{}_*$-вырожденная матрица имеет вид

$$A=\begin{pmatrix}d & c'd\\ db' & c'db'\end{pmatrix} \tag{5.8.12}$$

Из предположения следует, что (5.8.12) и (5.8.11) представляют одну и ту же матрицу. Сравнивая (5.8.12) и (5.8.11), мы получим, что для любых $d,c\in D$ существует такое $c'\in D$, которое не зависит от $d$ и удовлетворяет уравнению

$$dc=c'd$$

Это противоречит утверждению, что $D$ - тело. □

Пример 5.8.11. В случае тела кватернионов мы положим $b=1+k$, $c=j$, $d=k$. Тогда мы имеем

$$A=\begin{pmatrix}k & kj\\ (1+k)k & (1+k)kj\end{pmatrix}=\begin{pmatrix}k & -i\\ k-1 & -i-j\end{pmatrix}$$

$$\begin{aligned}{}^2\det({}_*{}^*)_2A&={}^2A_2-{}^2A_1({}^1A_1)^{-1}\ {}^1A_2\\ &=-i-j-(k-1)(k)^{-1}(-i)=-i-j-(k-1)(-k)(-i)\\ &=-i-j-kki+ki=-i-j+i+j=0\end{aligned}$$

$$\begin{aligned}{}_1\det({}^*{}_*)^1A&={}_1A^1-{}_1A^2({}_2A^2)^{-1}\ {}_2A^1\\ &=k-(k-1)(-i-j)^{-1}(-i)=k-(k-1)\frac{1}{2}(i+j)(-i)\\ &=k+\frac{1}{2}((k-1)i+(k-1)j)i=k+\frac{1}{2}(ki-i+kj-j)i\\ &=k+\frac{1}{2}(j-i-i-j)i=k-ii=k+1\end{aligned}$$

$$\begin{aligned}{}_1\det({}^*{}_*)^2A&={}_1A^2-{}_1A^1({}_2A^1)^{-1}\ {}_2A^2\\ &=k-1-k(-i)^{-1}(-i-j)=k-1+ki(i+j)\\ &=k-1+j(i+j)=k-1+ji+jj\\ &=k-1-k-1=-2\end{aligned}$$

$$\begin{aligned}{}_2\det({}^*{}_*)^1A&={}_2A^1-{}_2A^2({}_1A^2)^{-1}\ {}_1A^1\\ &=(-i)-(-i-j)(k-1)^{-1}k=-i+(i+j)\frac{1}{2}(-k-1)k\\ &=-i-\frac{1}{2}(i+j)(k+1)k=-i-\frac{1}{2}(ik+i+jk+j)k\\ &=-i-\frac{1}{2}(-j+i+i+j)k=-i-ik=-i+j\end{aligned}$$

$$\begin{aligned}{}_2\det({}^*{}_*)^2A&={}_2A^2-{}_2A^1({}_1A^1)^{-1}\ {}_1A^2\\ &=-i-j-(-i)(k)^{-1}(k-1)=-i-j+i(-k)(k-1)\\ &=-i-j+j(k-1)=-i-j+jk-j=-i-j+i-j=-2j\end{aligned}$$

Система ${}_*{}^*D$-линейных уравнений

$$\begin{pmatrix}k & -i\\ k-1 & -i-j\end{pmatrix}{}_*{}^*\begin{pmatrix}{}^1x\\ {}^2x\end{pmatrix}=\begin{pmatrix}{}^1b\\ {}^2b\end{pmatrix}\tag{5.8.13}$$

имеет ${}_*{}^*$-вырожденную матрицу. Мы можем записать систему ${}_*{}^*D$-линейных уравнений (5.8.13) в виде

$$\begin{cases}k\ {}^1x - i\ {}^2x={}^1b\\ (k-1)\ {}^1x-(i+j)\ {}^2x={}^2b\end{cases}$$

Система ${}^*{}_*D$-линейных уравнений

$$\begin{pmatrix}k & -i\\ k-1 & -i-j\end{pmatrix}{}^*{}_*\begin{pmatrix}{}_1x & {}_2x\end{pmatrix}=\begin{pmatrix}{}_1b & {}_2b\end{pmatrix}\tag{5.8.14}$$

имеет ${}^*{}_*$-невырожденную матрицу. Мы можем записать систему ${}^*{}_*D$-линейных уравнений (5.8.14) в виде

$$\begin{cases}k\ {}_1x & - & i\ {}_1x\\ +(k-1)\ {}_2x & - & (i+j)\ {}_2x\\ = & {}_1b & = & {}_2b\end{cases}\qquad\begin{cases}k\ {}_1x+(k-1)\ {}_2x={}_1b\\ -i\ {}_1x-(i+j)\ {}_2x={}_2b\end{cases}$$

Система $D_*{}^*$-линейных уравнений

$$\begin{pmatrix}x_1 & x_2\end{pmatrix}{}_*{}^*\begin{pmatrix}k & -i\\ k-1 & -i-j\end{pmatrix}=\begin{pmatrix}b_1 & b_2\end{pmatrix}\tag{5.8.15}$$

имеет ${}_*{}^*$-вырожденную матрицу. Мы можем записать систему $D_*{}^*$-линейных уравнений (5.8.15) в виде

$$\begin{cases}x_1k & -x_1i\\ +x_2(k-1) & -x_2(i+j)\\ =b_1 & =b_2\end{cases}\qquad\begin{cases}x_1k\ +x_2(k-1)=b_1\\ -x_1i\ -x_2(i+j)\ =b_2\end{cases}$$

Система $D^*{}_*$-линейных уравнений

$$\begin{pmatrix}x^1\\ x^2\end{pmatrix}{}^*{}_*\begin{pmatrix}k & -i\\ k-1 & -i-j\end{pmatrix}=\begin{pmatrix}{}^1b\\ {}^2b\end{pmatrix}\tag{5.8.16}$$

имеет ${}^*{}_*$-невырожденную матрицу. Мы можем записать систему $D^*{}_*$-линейных уравнений (5.8.16) в виде

$$\begin{cases}k\ {}^1x - i\ {}^2x={}^1b\\ (k-1)\ {}^1x-(i+j)\ {}^2x={}^2b\end{cases}$$

□

## 5.9. Размерность ${}_*{}^*D$-векторного пространства

Теорема 5.9.1. *Пусть $V$ - ${}_*{}^*D$-векторное пространство. Предположим, что $V$ имеет базисы $\overline{\overline{e}}=(e_i, i\in I)$ и $\overline{\overline{g}}=(g_j, j\in J)$. Если $|I|$ и $|J|$ - конечные числа, то $|I|=|J|$.*

Доказательство. Предположим, что $|I|=m$ и $|J|=n$. Предположим, что

(5.9.1) $$m<n$$

Так как $\overline{\overline{e}}$ - базис, любой вектор $g_j,\ j\in J$ имеет разложение

$$g_j=e_*{}^*A_j$$

Так как $\overline{\overline{g}}$ - базис,

(5.9.2) $$\lambda=0$$

следует из

$$g_*{}^*\lambda=e_*{}^*A_*{}^*\lambda=0$$

Так как $\overline{\overline{e}}$ - базис, мы получим

(5.9.3) $$A_*{}^*\lambda=0$$

Согласно (5.9.1) $\mathrm{rank}_{{}_*{}^*}A\le m$ и система (5.9.3) имеет больше переменных, чем уравнений. Согласно теореме 5.7.4, $\lambda\ne 0$. Это противоречит утверждению (5.9.2). Следовательно, утверждение $m<n$ неверно.

Таким же образом мы можем доказать, что утверждение $n<m$ неверно. Это завершает доказательство теоремы. □

Определение 5.9.2. Мы будем называть **размерностью ${}_*{}^*D$-векторного пространства** число векторов в базисе □

Теорема 5.9.3. *Координатная матрица ${}_*{}^*D$-базиса $\overline{\overline{g}}$ относительно ${}_*{}^*D$-базиса $\overline{\overline{e}}$ векторного пространства $V$ является ${}_*{}^*$-невырожденной матрицей.*

Доказательство. Согласно лемме 5.6.6 ${}_*{}^*D$-ранг координатной матрицы базиса $\overline{\overline{g}}$ относительно базиса $\overline{\overline{e}}$ равен размерности векторного пространства, откуда следует утверждение теоремы. □

Определение 5.9.4. Мы будем называть взаимно однозначное отображение

$$A:V\to W$$

**изоморфизмом ${}_*{}^*D$-векторных пространств**, если это отображение - линейное отображение ${}_*{}^*D$-векторных пространств. □

Определение 5.9.5. **Автоморфизм ${}_*{}^*D$-векторного пространства** $V$ - это изоморфизм $A:V\to V$. □

Теорема 5.9.6. *Предположим, что $\overline{\overline{f}}$ - ${}_*{}^*D$-базис в векторном пространстве $V$. Тогда любой автоморфизм $\overline{A}$ ${}_*{}^*D$-векторного пространства $V$ имеет вид*

(5.9.4) $$v'=A_*{}^*v$$

*где $A$ - ${}_*{}^*$-невырожденная матрица.*

Доказательство. (5.9.4) следует из теоремы 5.4.3. Так как $\overline{A}$ - изоморфизм, то для каждого вектора $v'$ существует единственный вектор $v$ такой, что $v' = v_*{}^*\overline{A}$. Следовательно, система ${}_*{}^*D$-линейных уравнений (5.9.4) имеет единственное решение. Согласно следствию 5.7.5 матрица $A$ невырожденна. □

Теорема 5.9.7. *Автоморфизмы ${}_*{}^*D$-векторного пространства порождают группу $GL(n, {}_*{}^*, D)$.*

Доказательство. Если даны два автоморфизма $\overline{A}$ и $\overline{B}$, то мы можем записать

$$v' = A_*{}^*v$$

$$v'' = B_*{}^*v' = B_*{}^*A_*{}^*v$$

Следовательно, результирующий автоморфизм имеет матрицу $B_*{}^*A$. □

Глава 6

# Многообразие базисов

## 6.1. Линейное представление группы

Пусть $V$ - ${}_*{}^*D$-векторное пространство. Мы доказали в теореме 5.9.6, что любой автоморфизм ${}_*{}^*D$-векторного пространства $V$ можно отождествить с некоторой матрицей. Мы доказали в теореме 5.9.7, что автоморфизмы ${}_*{}^*D$-векторного пространства порождают группу. Поэтому, когда мы изучаем представления в ${}_*{}^*D$-векторном пространстве, нас интересуют линейные отображения ${}_*{}^*D$-векторного пространства.

Определение 6.1.1. Пусть $V$ - ${}_*{}^*D$-векторное пространство.[6.1] Мы будем называть левостороннее представление[6.2] $f$ группы $G$ в ${}_*{}^*D$-векторном пространстве $V$ **линейным $G*$-представлением**. □

Теорема 6.1.2. *Автоморфизмы ${}_*{}^*D$-векторного пространства порождают линейное эффективное $GL(n, {}_*{}^*, D)*$-представление.*

Доказательство. Если даны два автоморфизма $\overline{A}$ и $\overline{B}$, то мы можем записать

$$v' = A_*{}^*v$$

$$v'' = B_*{}^*v' = B_*{}^*A_*{}^*v$$

Следовательно, результирующий автоморфизм имеет матрицу $B_*{}^*A$.

Нам осталось показать, что ядро неэффективности состоит только из единичного элемента. Тождественное преобразование удовлетворяет равенству

$$^iv = {}^iA_j\ {}^jv$$

Выбирая значения координат в форме ${}^iv = {}^i\delta_k$, где $k$ задано, мы получим

$$^i\delta_k = {}^iA_j\ {}^j\delta_k \tag{6.1.1}$$

Из (6.1.1) следует

$$^i\delta_k = {}^iA_k$$

Поскольку $k$ произвольно, мы приходим к заключению $A = \delta$. □

Теорема 6.1.3. *Пусть $\overline{\overline{f}}$, $\overline{\overline{g}}$ - базисы ${}_*{}^*D$-векторного пространства $\overline{V}$. Пусть $\overline{\overline{e}}$, $\overline{\overline{h}}$ - базисы ${}_*{}^*D$-векторного пространства $\overline{W}$. Пусть $A_1$ - матрица линейного отображения*

$$A : V \to W \tag{6.1.2}$$

[6.1] Изучая представление группы в ${}_*{}^*D$-векторном пространстве, мы будем следовать соглашению, описанному в замечании 2.2.15.

[6.2] Согласно теореме 5.9.7, линейное отображение правого векторного пространства является левосторонним преобразованием.

*относительно базисов* $\overline{\overline{f}}$ *и* $\overline{\overline{e}}$ *и* $A_2$ *- матрица линейного отображения* (6.1.2) *относительно базисов* $\overline{\overline{g}}$ *и* $\overline{\overline{h}}$. *Если базис* $\overline{\overline{f}}$ *имеет координатную матрицу* $B$ *относительно базиса* $\overline{\overline{g}}$

$$\overline{\overline{f}}=\overline{\overline{g}}_*{}^*B$$

*и* $\overline{\overline{e}}$ *имеет координатную матрицу* $C$ *относительно базиса* $\overline{\overline{h}}$

$$\overline{\overline{e}}=\overline{\overline{h}}_*{}^*C \tag{6.1.3}$$

*то матрицы* $A_1$ *и* $A_2$ *связаны соотношением*

$$A_1=C^{-1_*{}^*}{}_*{}^*A_{2*}{}^*B \tag{6.1.4}$$

Доказательство. Вектор $\overline{a}\in V$ имеет разложение

$$\overline{a}=f_*{}^*a=g_*{}^*B_*{}^*a$$

относительно базисов $\overline{\overline{f}}$ и $\overline{\overline{g}}$. Так как $A$ - линейное отображение, то мы можем записать его в форме

$$\overline{b}=e_*{}^*A_{1*}{}^*a \tag{6.1.5}$$

относительно базисов $\overline{\overline{f}}$ и $\overline{\overline{e}}$ и

$$\overline{b}=h_*{}^*A_{2*}{}^*B_*{}^*a \tag{6.1.6}$$

относительно базисов $\overline{\overline{g}}$ и $\overline{\overline{h}}$. В силу теоремы 5.9.3 матрица $C$ имеет ${}_*{}^*$-обратную и из равенства (6.1.3) следует

$$\overline{\overline{h}}=\overline{\overline{e}}_*{}^*C^{-1_*{}^*} \tag{6.1.7}$$

Подставив (6.1.7) в равенство (6.1.6), получим

$$\overline{b}=e_*{}^*C^{-1_*{}^*}{}_*{}^*A_{2*}{}^*B_*{}^*a \tag{6.1.8}$$

В силу теоремы 5.3.3 сравнение равенств (6.1.5) и (6.1.8) даёт

$$A_{1*}{}^*a=C^{-1_*{}^*}{}_*{}^*A_{2*}{}^*B_*{}^*a \tag{6.1.9}$$

Так как вектор $a$ - произвольный вектор, из теоремы 6.1.2 и равенства (6.1.9) следует утверждение теоремы. □

Теорема 6.1.4. *Пусть* $\overline{A}$ *- автоморфизм* ${}_*{}^*D$*-векторного пространства. Пусть* $A_1$ *- матрица этого автоморфизма, заданная относительно базиса* $\overline{\overline{f}}$, *и* $A_2$ *- матрица того же автоморфизма, заданная относительно базиса* $\overline{\overline{g}}$. *Если базис* $\overline{\overline{f}}$ *имеет координатную матрицу* $B$ *относительно базиса* $\overline{\overline{g}}$

$$\overline{\overline{f}}=\overline{\overline{g}}_*{}^*B$$

*то матрицы* $A_1$ *и* $A_2$ *связаны соотношением*

$$A_1=B^{-1_*{}^*}{}_*{}^*A_{2*}{}^*B$$

Доказательство. Утверждение следует из теоремы 6.1.3, так как в данном случае $C=B$. □

## 6.2. Многообразие базисов ${}_*{}^*D$-векторного пространства

Теорема 6.2.1. *Автоморфизм $A$, действуя на каждый вектор базиса в ${}_*{}^*D$-векторном пространстве, отображает базис в другой базис.*

Доказательство. Пусть $\overline{\overline{e}}$ - базис в ${}_*{}^*D$-векторном пространстве $V$. Согласно теореме 5.9.6, вектор $e_a$ отображается в вектор $e'_a$

$$e'_a = A_*{}^*e_a \tag{6.2.1}$$

Если векторы $e'_a$ линейно зависимы, то в линейной комбинации

$$e'{}_*{}^*\lambda = 0 \tag{6.2.2}$$

$\lambda \neq 0$. Из равенств (6.2.1) и (6.2.2) следует, что

$$A^{-1_*{}^*}{}_*{}^*e'{}_*{}^*\lambda = e_*{}^*\lambda = 0$$

и $\lambda \neq 0$. Это противоречит утверждению, что векторы $e_a$ линейно независимы. Следовательно, векторы $e'_a$ линейно независимы и порождают базис. □

Таким образом, мы можем распространить линейное $GL(n, {}_*{}^*, D)*$-представление на множество базисов. Мы будем называть преобразование этого левостороннего представления **активным преобразованием** потому, что линейное отображение векторного пространства породило это преобразование ([5]). Активное преобразование не является линейным преобразованием, так как на множестве базисов не определена линейная операция. Соответственно определению мы будем записывать действие активного преобразования $A \in GL(n, {}_*{}^*, D)$ на базис $\overline{\overline{e}}$ в форме $A_*{}^*\overline{\overline{e}}$. Мы будем называть гомоморфизм группы $G$ в группу $GL(n, {}_*{}^*, D)$ активных преобразований **активным $*G$-представлением**.

Теорема 6.2.2. *Активное $GL(n, {}_*{}^*, D)*$-представление на множестве базисов однотранзитивно.*

Доказательство. Чтобы доказать теорему, достаточно показать, что для любых двух базисов определено по крайней мере одно преобразование левостороннего представления и это преобразование единственно. Гомоморфизм $A$, действуя на базис $\overline{\overline{e}}$ имеет вид

$$g_i = A_*{}^*e_i$$

где $g_i$ - координатная матрица вектора $\overline{g}_i$ и $e_i$ - координатная матрица вектора $\overline{e}_i$ относительно базиса $\overline{\overline{h}}$. Следовательно, координатная матрица образа базиса равна ${}_*{}^*$-произведению координатной матрицы исходного базиса и матрицы автоморфизма

$$g = A_*{}^*e$$

В силу теоремы 5.9.3, матрицы $g$ и $e$ невырождены. Следовательно, матрица

$$A = g_*{}^*e^{-1_*{}^*}$$

является матрицей автоморфизма, отображающего базис $\overline{\overline{e}}$ в базис $\overline{\overline{g}}$.

Допустим элементы $g_1$, $g_2$ группы $G$ и базис $\overline{\overline{e}}$ таковы, что

$$g_{1*}{}^*e = g_{2*}{}^*e \tag{6.2.3}$$

В силу теорем 5.9.3 и 2.2.16 справедливо равенство $g_1 = g_2$. Отсюда следует утверждение теоремы. □

Если на векторном пространстве $V$ определена дополнительная структура, не всякое линейное отображение сохраняет свойства заданной структуры. В этом случае нас интересует подгруппа $G$ группы $GL(n, {}_*{}^*, D)$, которая порождает линейные отображения, сохраняющие свойства заданной структуры. Мы обычно будем называть группу $G$ **группой симметрии**. Не нарушая общности, мы будем отождествлять элемент $g$ группы $G$ с соответствующим преобразованием представления и записывать его действие на вектор $v \in V$ в виде $g_*{}^*v$.

Не всякие два базиса могут быть связаны преобразованием группы симметрии потому, что не всякое невырожденное линейное преобразование принадлежит представлению группы $G$. Таким образом, множество базисов можно представить как объединение орбит группы $G$.

Определение 6.2.3. Мы будем называть орбиту $G_*{}^*\overline{\overline{e}}$ выбранного базиса $\overline{\overline{e}}$ **многообразием базисов** $\mathcal{B}(\overline{V}, G)$ ${}_*{}^*D$**-векторного пространства** $V$. □

Теорема 6.2.4. *Активное $G*$-представление на многообразии базисов однотранзитивно.*

Доказательство. Это следствие теоремы 6.2.2 и определения 6.2.3. □

Из теоремы 6.2.4 следует, что многообразие базисов $\mathcal{B}(V, G)$ является однородным пространством группы $G$. Согласно теореме 4.2.13, на многообразии базисов существует $*G$-представление, перестановочное с активным. Как мы видим из замечания 4.2.14 преобразование $*G$-представления отличается от активного преобразования и не может быть сведено к преобразованию пространства $V$. Чтобы подчеркнуть различие, это преобразование называется **пассивным преобразованием** многообразия базисов $\mathcal{B}(\overline{V}, G)$, а $*G$-представление называется **пассивным** $*G$**-представлением**. Согласно определению мы будем записывать пассивное преобразование базиса $\overline{\overline{e}}$, порождённое элементом $A \in G$, в форме $\overline{\overline{e}}_*{}^*A$.

Замечание 6.2.5. Я привёл примеры пассивных и активных представлений в таблице 6.2.1. □

Таблица 6.2.1. Активное и пассивное представления

| векторное пространство | группа представления | активное представление | пассивное представление |
|---|---|---|---|
| ${}_*{}^*D$-векторное пространство | $GL(n, {}_*{}^*, D)$ | левостороннее | правостороннее |
| ${}^*{}_*D$-векторное пространство | $GL(n, {}^*{}_*, D)$ | левостороннее | правостороннее |
| $D_*{}^*$-векторное пространство | $GL(n, {}_*{}^*, D)$ | правостороннее | левостороннее |
| $D^*{}_*$-векторное пространство | $GL(n, {}^*{}_*, D)$ | правостороннее | левостороннее |

Согласно теореме 4.2.7 мы можем определить на $\mathcal{B}(V, G)$ две формы координат, определённые на группе $G$. Так как мы определили два представления группы $G$ на $\mathcal{B}(V, G)$, то мы для определения координат пользуемся пассивным $*G$-представлением. Наш выбор основан на следующей теореме.

Теорема 6.2.6. *Координатная матрица базиса $\overline{\overline{g}}$ относительно базиса $\overline{\overline{e}}$ ${}_*{}^*D$-векторного пространства $V$ совпадает с матрицей пассивного преобразования, отображающего базис $\overline{\overline{e}}$ в базис $\overline{\overline{g}}$.*

Доказательство. Согласно конструкции, изложенной в примере 5.3.9, координатная матрица базиса $\overline{\overline{g}}$ относительно базиса $\overline{\overline{e}}$ состоит из ${}^*$-строк, являющихся координатными матрицами векторов $g_i$ относительно базиса $\overline{\overline{e}}$. Следовательно,

$$\overline{g}_i = \overline{e}_*{}^* g_i \tag{6.2.4}$$

В тоже время пассивное преобразование $A$, связывающее два базиса, имеет вид

$$\overline{g}_i = \overline{e}_*{}^* A_i \tag{6.2.5}$$

В силу теоремы 5.3.3,

$$g_i = A_i$$

для любого $i$. Это доказывает теорему. □

Координаты представления называются **стандартными координатами базиса**. Эта точка зрения позволяет определить два типа координат для элемента $g$ группы $G$. Мы можем либо пользоваться координатами, определёнными на группе, либо определить координаты как элементы матрицы соответствующего преобразования. Первая форма координат более эффективна, когда мы изучаем свойства группы $G$. Вторая форма координат содержит избыточную информацию, но бывает более удобна, когда мы изучаем представление группы $G$. Мы будем называть вторую форму координат **координатами представления**.

## 6.3. Геометрический объект в ${}_*{}^*D$-векторном пространстве

Активное преобразование изменяет базисы и векторы согласовано и координаты вектора относительно базиса не меняются. Пассивное преобразование меняет только базис, и это ведёт к изменению координат вектора относительно базиса.

Допустим пассивное преобразование $A \in G$ отображает базис $\overline{\overline{e}} \in \mathcal{B}(V, G)$ в базис $\overline{\overline{e}}' \in \mathcal{B}(V, G)$

$$e' = e_*{}^* A \tag{6.3.1}$$

Допустим вектор $v \in V$ имеет разложение

$$\overline{v} = e_*{}^* v \tag{6.3.2}$$

относительно базиса $\overline{\overline{e}}$ и имеет разложение

$$\overline{v} = e'{}_*{}^* v' \tag{6.3.3}$$

относительно базиса $\overline{\overline{e}}'$. Из (6.3.1) и (6.3.3) следует, что

$$\overline{v} = e_*{}^* A_*{}^* v' \tag{6.3.4}$$

Сравнивая (6.3.2) и (6.3.4) получаем, что

$$v = A_*{}^* v' \tag{6.3.5}$$

Так как $A$ - ${}_*{}^*$-невырожденная матрица, то из (6.3.5) следует

$$v' = A^{-1_*{}^*}{}_*{}^* v \tag{6.3.6}$$

Преобразование координат (6.3.6) не зависит от вектора $\overline{v}$ или базиса $\overline{\overline{e}}$, а определенно исключительно координатами вектора $\overline{v}$ относительно базиса $\overline{\overline{e}}$.

Теорема 6.3.1. *Преобразования координат* (6.3.6) *порождают эффективное линейное* $G*$*-представление, называемое* **координатным представлением в ${}_*{}^*D$-векторном пространстве**.

Доказательство. Допустим мы имеем два последовательных пассивных преобразования $A$ и $B$. Преобразование координат (6.3.6) соответствует пассивному преобразованию $A$. Преобразование координат

$$v'' = B^{-1_*{}^*}{}_*{}^*v' \tag{6.3.7}$$

соответствует пассивному преобразованию $B$. Произведение преобразований координат (6.3.6) и (6.3.7) имеет вид

$$v'' = B^{-1_*{}^*}{}_*{}^*A^{-1_*{}^*}{}_*{}^*v = (A_*{}^*B)^{-1_*{}^*}{}_*{}^*v \tag{6.3.8}$$

и является координатным преобразованием, соответствующим пассивному преобразованию $A_*{}^*B$. Это доказывает, что преобразования координат порождают линейное $G*$-представление.

Если координатное преобразование не изменяет векторы $\delta_k$, то ему соответствует единица группы $G$, так как пассивное представление однотранзитивно. Следовательно, координатное представление эффективно. □

Предположим, что отображение группы $G$ в группу пассивных преобразований ${}_*{}^*D$-векторного пространства $N$ согласовано с группой симметрий ${}_*{}^*D$-векторного пространства $V$.[6.3] Это означает, что пассивному преобразованию $a$ ${}_*{}^*D$-векторного пространства $V$ соответствует пассивное преобразование $A(a)$ ${}_*{}^*D$-векторного пространства $N$.

$$e'_N = e_{N*}{}^*A(a) \tag{6.3.9}$$

Тогда координатное преобразование в $N$ принимает вид

$$w' = A(a^{-1_*{}^*})_*{}^*w = A(a)^{-1_*{}^*}{}_*{}^*w \tag{6.3.10}$$

Определение 6.3.2. Мы будем называть орбиту

$$(A(G)^{-1_*{}^*}{}_*{}^*w, \overline{\overline{e}}_{V*}{}^*G)$$

**геометрическим объектом в координатном представлении, определённом в ${}_*{}^*D$-векторном пространстве** $V$. Для любого базиса $\overline{\overline{e}}'_V = \overline{\overline{e}}_{V*}{}^*A$ соответствующая точка (6.3.10) орбиты определяет **координаты геометрического объекта в координатном ${}_*{}^*D$-векторном пространстве** относительно базиса $\overline{\overline{e}}'_V$. □

Определение 6.3.3. Мы будем называть орбиту

$$(A(G)^{-1_*{}^*}{}_*{}^*w, \overline{\overline{e}}_{N*}{}^*A(G), \overline{\overline{e}}_{V*}{}^*G)$$

**геометрическим объектом, определённым в ${}_*{}^*D$-векторном пространстве** $V$. Для любого базиса $\overline{\overline{e}}'_V = a_*{}^*\overline{\overline{e}}_V$ соответствующая точка (6.3.10) орбиты определяет **координаты геометрического объекта в ${}_*{}^*D$-векторном пространстве** относительно базиса $\overline{\overline{e}}'_V$ и соответствующий вектор

$$\overline{w} = e'_{N*}{}^*w'$$

[6.3] Мы пользуемся одним и тем же символом типа векторного пространства для векторных пространств $N$ и $V$. Тем не менее они могут иметь различный тип.

называется **представителем геометрического объекта в $D_*{}^*$-векторном пространстве** $V$ в базисе $\overline{\overline{e}}'_V$. □

Так как геометрический объект - это орбита представления, то согласно теореме 4.1.9 определение геометрического объекта корректно.

Мы будем также говорить, что $\overline{w}$ - это **геометрический объект типа** $A$

Определение 6.3.2 строит геометрический объект в координатном пространстве. Определение 6.3.3 предполагает, что мы выбрали базис в векторном пространстве $W$. Это позволяет использовать представитель геометрического объекта вместо его координат.

Вопрос как велико разнообразие геометрических объектов хорошо изучен в случае векторных пространств. Однако он не столь очевиден в случае ${}_*{}^*D$-векторных пространств. Как видно из таблицы 6.2.1 $D_*{}^*$-векторное пространство и ${}_*{}^*D$-векторное пространство имеют общую группу симметрии $GL(n, {}_*{}^*, D)$. Это позволяет, рассматривая пассивное представление в ${}_*{}^*D$-векторное пространство, изучать геометрический объект в $D_*{}^*$-векторном пространстве. Можем ли мы одновременно изучать геометрический объект в ${}^*{}_*D$-пространстве? На первый взгляд, в силу теоремы 5.8.9 ответ отрицательный. Однако, равенство 2.2.6 устанавливает искомое линейное отображение между $GL(n, {}_*{}^*, D)$ и $GL(n, {}^*{}_*, D)$.

Теорема 6.3.4 (принцип инвариантности). *Представитель геометрического объекта не зависит от выбора базиса* $\overline{\overline{e}}'_V$.

Доказательство. Чтобы определить представителя геометрического объекта, мы должны выбрать базис $\overline{\overline{e}}_V$, базис $\overline{\overline{e}}_W$ и координаты геометрического объекта $w^\alpha$. Соответствующий представитель геометрического объекта имеет вид

$$\overline{w} = e_W{}_*{}^*w$$

Предположим базис $\overline{\overline{e}}'_V$ связан с базисом $\overline{\overline{e}}_V$ пассивным преобразованием

$$e'_V = e_V{}_*{}^*A$$

Согласно построению это порождает пассивное преобразование (6.3.9) и координатное преобразование (6.3.10). Соответствующий представитель геометрического объекта имеет вид

$$\overline{w}' = e'_W{}_*{}^*w' = e_W{}_*{}^*A(a){}_*{}^*A(a)^{-1}{}_*{}^*w = e_W{}_*{}^*w = \overline{w}$$

Следовательно, представитель геометрического объекта инвариантен относительно выбора базиса. □

Определение 6.3.5. Пусть

$$\overline{w}_1 = e_W{}_*{}^*w_1$$

$$\overline{w}_2 = e_W{}_*{}^*w_2$$

геометрические объекты одного и того же типа, определённые в ${}_*{}^*D$-векторном пространстве $V$. Геометрический объект

$$\overline{w} = e_W{}_*{}^*(w_1 + w_2)$$

называется **суммой**

$$\overline{w} = \overline{w}_1 + \overline{w}_2$$

**геометрических объектов** $\overline{w}_1$ и $\overline{w}_2$. □

Определение 6.3.6. Пусть

$$\overline{w}_1 = e_{W*}{}^*w_1$$

геометрический объект, определённый в ${}_*{}^*D$-векторном пространстве $V$. Геометрический объект

$$\overline{w}_2 = e_{W*}{}^*(kw_1)$$

называется **произведением**

$$\overline{w}_2 = k\overline{w}_1$$

**геометрического объекта** $\overline{w}_1$ **и константы** $k \in D$. □

Теорема 6.3.7. *Геометрические объекты типа $A$, определённые в ${}_*{}^*D$-векторном пространстве $V$, образуют ${}_*{}^*D$-векторное пространство.*

Доказательство. Утверждение теоремы следует из непосредственной проверки свойств векторного пространства. □

Глава 7

# Линейное отображение ${}_*{}^*D$-векторного пространства

## 7.1. Линейное отображение ${}_*{}^*D$-векторного пространства

В этом разделе мы положим, что $V$, $W$ - это ${}_*{}^*D$-векторные пространства.

Определение 7.1.1. Обозначим $\mathcal{L}({}_*{}^*D;V;W)$ множество линейных отображений

$$A:V\rightarrow W$$

${}_*{}^*D$-векторного пространства $V$ в ${}_*{}^*D$-векторное пространство $W$. Обозначим $\mathcal{L}(D_*{}^*;V;W)$ множество линейных отображений

$$A:V\rightarrow W$$

$D_*{}^*$-векторного пространства $V$ в $D_*{}^*$-векторное пространство $W$. □

Мы можем рассматривать тело $D$ как одномерное ${}_*{}^*D$-векторное пространство. Соответственно мы можем рассматривать множества $\mathcal{L}({}_*{}^*D;D;W)$ и $\mathcal{L}({}_*{}^*D;V;D)$.

Теорема 7.1.2. *Предположим, что $V$, $W$ - ${}_*{}^*D$-векторные пространства. Тогда множество $\mathcal{L}({}_*{}^*D;V;W)$ является абелевой группой относительно закона композиции*

$$(A+B)_*{}^*x=A_*{}^*x+B_*{}^*x \tag{7.1.1}$$

*Линейное отображение $A+B$ называется* **суммой отображений** *$A$ и $B$.*

Доказательство. Нам надо показать, что отображение

$$A+B:V\rightarrow W$$

определённое равенством (7.1.1), - это линейное отображение ${}_*{}^*D$-векторных пространств. Согласно определению 5.4.2

$$A_*{}^*(x_*{}^*a)=(A_*{}^*x)_*{}^*a$$
$$B_*{}^*(x_*{}^*a)=(B_*{}^*x)_*{}^*a$$

Мы видим, что

$$\begin{aligned}(A+B)_*{}^*(x_*{}^*a)&=A_*{}^*(x_*{}^*a)+B_*{}^*(x_*{}^*a)\\&=(A_*{}^*x)_*{}^*a+(B_*{}^*x)_*{}^*a\\&=(A_*{}^*x+B_*{}^*x)_*{}^*a\\&=((A+B)_*{}^*x)_*{}^*a\end{aligned}$$

Нам надо показать так же, что эта операция коммутативна.

$$\begin{aligned}(A+B)_*{}^*x &= A_*{}^*x + B_*{}^*x\\ &= B_*{}^*x + A_*{}^*x\\ &= (B+A)_*{}^*x\end{aligned}$$

□

Теорема 7.1.3. *Пусть* $\overline{\overline{e}}_V = (e_{V\cdot\boldsymbol{i}}, \boldsymbol{i}\in\boldsymbol{I})$ *- базис в* ${}_*{}^*D$*-векторном пространстве* $V$ *и* $\overline{\overline{e}}_W = (e_{W\cdot\boldsymbol{j}}, \boldsymbol{j}\in\boldsymbol{J})$ *- базис в* ${}_*{}^*D$*-векторном пространстве* $W$*. Пусть* $A = ({}^{\boldsymbol{j}}A_{\boldsymbol{i}})$, $\boldsymbol{i}\in\boldsymbol{I}$, $\boldsymbol{j}\in\boldsymbol{J}$, *- произвольная матрица. Тогда отображение*

$$\overline{A}: V\to W$$

*определённое равенством*

$$\overline{a} = e_V{}_*{}^*a \to \overline{A}_*{}^*\overline{a} = e_W{}_*{}^*A_*{}^*a \tag{7.1.2}$$

*относительно выбранных базисов, является линейным отображением* ${}_*{}^*D$*-векторных пространств.*

Доказательство. Теорема 7.1.3 является обратным утверждением теореме 5.4.3. Из равенства (7.1.2) следует, что

$$\overline{A}_*{}^*(\overline{v}_*{}^*a) = e_W{}_*{}^*A_*{}^*v_*{}^*a = (\overline{A}_*{}^*\overline{v})_*{}^*a$$

□

Теорема 7.1.4. *Пусть* $\overline{\overline{e}}_V = (e_{V\cdot\boldsymbol{i}}, \boldsymbol{i}\in\boldsymbol{I})$ *- базис в* ${}_*{}^*D$*-векторном пространстве* $V$ *и* $\overline{\overline{e}}_W = (e_{W\cdot\boldsymbol{j}}, \boldsymbol{j}\in\boldsymbol{J})$ *- базис в* ${}_*{}^*D$*-векторном пространстве* $W$*. Предположим линейное отображение* $\overline{A}$ *имеет матрицу* $A = ({}^{\boldsymbol{j}}A_{\boldsymbol{i}})$, $\boldsymbol{i}\in\boldsymbol{I}$, $\boldsymbol{j}\in\boldsymbol{J}$, *относительно выбранных базисов. Пусть* $m\in D$*. Тогда матрица*

$${}^{\boldsymbol{j}}(mA)_{\boldsymbol{i}} = m\ {}^{\boldsymbol{j}}A_{\boldsymbol{i}}$$

*определяет линейное отображение*

$$m\overline{A}: V\to W$$

*которое мы будем называть* **левосторонним произведением отображения** $A$ **на скаляр**.

Доказательство. Утверждение теоремы является следствием теоремы 7.1.3. □

Теорема 7.1.5. *Множество* $\mathcal{L}({}_*{}^*D;V;W)$ *является* $D_*{}^*$*-векторным пространством.*

Доказательство. Теорема 7.1.2 утверждает, что $\mathcal{L}({}_*{}^*D;V;W)$ - абелевая группа. Из теоремы 7.1.4 следует, что элемент тела $D$ порождает левостороннее преобразование на абелевой группе $\mathcal{L}({}_*{}^*D;V;W)$. Из теорем 7.1.3, 5.1.1 и 5.1.3 следует, что множество $\mathcal{L}({}_*{}^*D;V;W)$ является $D\star$-векторным пространством.

Выписывая элементы базиса $D\star$-векторного пространства $\mathcal{L}({}_*{}^*D;V;W)$ в виде ${}_*$-строки или ${}^*$-строки, мы представим $D\star$-векторное пространство $\mathcal{L}({}_*{}^*D;V;W)$ как $D^*{}_*$- или $D_*{}^*$-векторное пространство. При этом надо иметь в виду, что выбор между $D_*{}^*$- и $D^*{}_*$-линейной зависимостью в $D\star$-векторном пространстве $\mathcal{L}({}_*{}^*D;V;W)$ не зависит от типа векторных пространств $V$ и $W$.

Для того, чтобы выбрать тип векторного пространства $\mathcal{L}({}_*{}^*D;V;W)$ мы обратим внимание на следующее обстоятельство. Допустим $V$ и $W$ - ${}_*{}^*D$-векторные пространства. Допустим $\mathcal{L}({}_*{}^*D;V;W)$ - $D_*{}^*$-векторное пространство. Тогда действие ${}^*$-строки ${}_*{}^*D$-линейных отображений ${}^jA$ на ${}_*$-строку векторов $f_i$ можно представить в виде матрицы

$$\begin{pmatrix}{}^1A\\ ...\\ {}^mA\end{pmatrix}{}_*{}^*\begin{pmatrix}f_1 & ... & f_n\end{pmatrix}=\begin{pmatrix}{}^1A_*{}^*f_1 & ... & {}^1A_*{}^*f_n\\ ... & ... & ...\\ {}^mA_*{}^*f_1 & ... & {}^mA_*{}^*f_n\end{pmatrix}$$

Эта запись согласуется с матричной записью действия $D_*{}^*$-линейной комбинации $a_*{}^*A$ ${}_*{}^*D$-линейных отображений $A$. □

Мы можем также определить правостороннее произведение ${}_*{}^*D$-линейного отображения $A$ на скаляр. Однако, вообще говоря, это левостороннее представление тела $D$ в $D\star$-векторном пространстве $\mathcal{L}({}_*{}^*D;V;W)$ не может быть перенесено в ${}_*{}^*D$-векторное пространство $W$. Действительно, в случае правостороннего произведения мы имеем

$$(Am)_*{}^*v=(Am)_*{}^*f_*{}^*v=e_*{}^*(Am)_*{}^*v$$

Поскольку произведение в теле не коммутативно, мы не можем выразить полученное выражение как произведение $A_*{}^*v$ на скаляр $m$.

Следствием теоремы 7.1.5 является неоднозначность записи

$$m_*{}^*A_*{}^*v$$

Мы можем предположить, что смысл этой записи ясен из контекста. Однако желательно неоднозначность избежать. Мы будем пользоваться для этой цели скобками. Выражение

$$w=[m_*{}^*A]_*{}^*v$$

означает, что $D_*{}^*$-линейная комбинация ${}_*{}^*D$-линейных отображений ${}^iA$ отображает вектор $v$ в вектор $w$. Выражение

$$w=B_*{}^*A_*{}^*v$$

означает, что ${}_*{}^*$-произведение линейных отображений $A$ и $B$ отображает вектор $v$ в вектор $w$.

## 7.2. 1-форма на ${}_*{}^*D$-векторном пространстве

Определение 7.2.1. **1-форма** на ${}_*{}^*D$-векторном пространстве $V$ - это линейное отображение

$$b:V\to D \tag{7.2.1}$$

□

Мы можем значение 1-формы $b$, определённое для вектора $a$, записывать в виде

$$b(a)=<b,a>$$

Теорема 7.2.2. *Множество $\mathcal{L}({}_*{}^*D;V;D)$ является $D\star$-векторным пространством.*

Доказательство. ${}_*{}^*D$-векторное пространство размерности 1 эквивалентно телу $D$ □

Теорема 7.2.3. *Пусть* $\overline{\overline{e}}$ *- базис в* ${}_*{}^*D$*-векторном пространстве* $V$*. 1-форма* $\overline{b}$ *имеет представление*

$$<\overline{b},\overline{a}>=b_*{}^*a \tag{7.2.2}$$

*относительно выбранного базиса, где вектор* $a$ *имеет разложение*

$$\overline{a}=e_*{}^*a \tag{7.2.3}$$

*и*

$$b_i=<\overline{b},e_i> \tag{7.2.4}$$

Доказательство. Так как $\overline{b}$ - 1-форма, из (7.2.3) следует, что

$$<\overline{b},\overline{a}>=<\overline{b},e_*{}^*a>=<\overline{b},e_i>\ {}^ia \tag{7.2.5}$$

(7.2.2) следует из (7.2.5) и (7.2.4). □

Теорема 7.2.4. *Пусть* $\overline{\overline{e}}$ *- базис в* ${}_*{}^*D$*-векторном пространстве* $V$*. 1-форма*

$$\overline{b}:V\to D$$

*однозначно определена значениями* (7.2.4)*, в которые 1-форма* $\overline{b}$ *отображает вектора базиса.*

Доказательство. Утверждение является следствием теорем 7.2.3 и 5.3.3. □

Теорема 7.2.5. *Пусть* $\overline{\overline{e}}$ *- базис в* ${}_*{}^*D$*-векторном пространстве* $V$*. Множество 1-форм* ${}^id$ *таких, что*

$$<{}^jd,e_i>={}^j\delta_i \tag{7.2.6}$$

*является базисом* $\overline{\overline{d}}$ $D_*{}^*$*-векторного пространства* $\mathcal{L}({}_*{}^*D;V;D)$*.*

Доказательство. 1-форма ${}^jd=\overline{b}$ существует для данного $j$ согласно теореме 7.2.4, если положить $b_i={}^j\delta_i$ . Если предположить, что существует 1-форма

$$\overline{c}=c_*{}^*d=0$$

то

$$c_j<{}^jd,e_i>=0$$

Согласно равенству (7.2.6)

$$c_i=c_j\ {}^j\delta_i=0$$

Следовательно, 1-формы ${}^jd$ линейно независимы. □

Определение 7.2.6. Пусть $V$ - ${}_*{}^*D$-векторное пространство. $D_*{}^*$-векторное пространство

$$V^*=\mathcal{L}({}_*{}^*D;V;D)$$

называется **дуальным пространством к** ${}_*{}^*D$**-векторному пространству** $V$. Пусть $\overline{\overline{e}}$ - базис в ${}_*{}^*D$-векторном пространстве $V$. Базис $\overline{\overline{d}}$ $D_*{}^*$-векторного пространства $V^*$, удовлетворяющий равенству (7.2.6), называется **базисом, дуальным базису** $\overline{\overline{e}}$. □

Теорема 7.2.7. *Пусть $A$ - пассивное преобразование многообразия базисов $\mathcal{B}(V, GL(n, {}_*{}^*, D))$. Допустим базис*

$$\overline{\overline{e}}' = \overline{\overline{e}}_*{}^* A \tag{7.2.7}$$

*является образом базиса $\overline{\overline{e}}$. Пусть $B$ - пассивное преобразование многообразия базисов $\mathcal{B}(V^*, GL(n, {}_*{}^*, D))$ такое, что базис*

$$\overline{\overline{d}}' = B_*{}^* \overline{\overline{d}} \tag{7.2.8}$$

*дуален базису. Тогда*

$$B = A^{-1_*{}^*} \tag{7.2.9}$$

Доказательство. Из равенств (7.2.6), (7.2.7), (7.2.8) следует

$$\begin{aligned} {}^j\delta_i &= < {}^j d', e'_i > \\ &= {}^j B_l \ < {}^l d, e_k > \ {}^k A_i \\ &= {}^j B_l \ {}^l\delta_k \ {}^k A_i \\ &= {}^j B_k \ {}^k A_i \end{aligned} \tag{7.2.10}$$

Равенство (7.2.9) следует из равенства (7.2.10). □

## 7.3. Парные представления тела

Теорема 7.3.1. *В любом ${}_*{}^*D$-векторном пространстве можно определить структуру $D*$-векторного пространства, определив правостороннее произведение вектора на скаляр равенством*

$$mv = m\delta_*{}^* v$$

Доказательство. Непосредственная проверка доказывает, что отображение

$$f : D \to V^*$$

определённое равенством

$$f(m) = m\delta$$

определяет $D*$-представление. □

Мы можем иначе сформулировать теорему 7.3.1.

Теорема 7.3.2. *Если мы определили эффективное правостороннее представление $f$ тела $D$ на абелевой группе $V$, то мы можем однозначно определить эффективное левостороннее представление $h$ тела $D$ на абелевой группе $V$ такое, что диаграмма*

$$\begin{array}{ccc} V & \xrightarrow{h(a)} & V \\ \Big\downarrow{\scriptstyle f(b)} & & \Big\downarrow{\scriptstyle f(b)} \\ V & \xrightarrow[h(a)]{} & V \end{array}$$

*коммутативна для любых $a$, $b \in D$.* □

Мы будем называть представления $f$ и $h$ **парными представлениями тела** $D$.

Теорема 7.3.3. *В векторном пространстве $V$ над телом $D$ мы можем определить $D\star$-произведение и $\star D$-произведение вектора на скаляр. Согласно теореме 7.3.2 эти операции удовлетворяют равенству*

$$(am)b = a(mb) \tag{7.3.1}$$

Равенство (7.3.1) представляет **закон ассоциативности для парных представлений**. Это позволяет нам писать подобные выражения не пользуясь скобками.

Доказательство. В разделе 7.1 дано определение правостороннего произведения линейного отображения $A$ ${}_*{}^*D$-векторного пространства на скаляр. Согласно теореме 7.3.1, $D*$-представление в абелевой группе $\mathcal{L}({}_*{}^*D; V; V)$ может быть перенесено в ${}_*{}^*D$-векторное пространство $V$ согласно правилу

$$[mA]_*{}^*v = [m\delta]_*{}^*(A_*{}^*v) = m(A_*{}^*v)$$

□

Аналогия с векторными пространствами над полем заходит столь далеко, что мы можем предположить существование концепции базиса, который годится для левостороннего и правостороннего произведения вектора на скаляр.

Теорема 7.3.4. *Многообразие базисов $D_*{}^*$-векторного пространства $V$ и многообразие базисов ${}^*{}_*D$-векторного пространства $V$ отличны*

$$\mathcal{B}(V, D_*{}^*) \ne \mathcal{B}(V, {}^*{}_*D)$$

Доказательство. При доказательстве этой теоремы мы будем пользоваться стандартным представлением матрицы. Не нарушая общности, мы проведём доказательство в координатном векторном пространстве $D^n$.

Пусть $\overline{\overline{e}} = (e_i = (\delta_i^j), i, j \in I, |I| = n)$ множество векторов векторного пространства $D^n$. Очевидно, $\overline{\overline{e}}$ является одновременно базисом $D_*{}^*$-векторного пространства и базисом ${}^*{}_*D$-векторного пространства. Для произвольного множества векторов $(f_i, i \in I, |I| = n)$ $D_*{}^*$-координатная матрица

$$f = \begin{pmatrix} f_1^1 & ... & f_1^n \\ ... & ... & ... \\ f_n^1 & ... & f_n^n \end{pmatrix} \tag{7.3.2}$$

относительно базиса $\overline{\overline{e}}$ совпадает с ${}^*{}_*D$-координатной матрицей относительно базиса $\overline{\overline{e}}$.

Если множество векторов $(f_i, i \in I, |I| = n)$ - базис ${}_*{}^*D$-векторного пространства, то согласно 5.9.3 матрица (7.3.2) - ${}_*{}^*$-невырожденная матрица.

Если множество векторов $(f_i, i \in I, |I| = n)$ - базис ${}^*{}_*D$-векторного пространства, то согласно теореме 5.9.3 матрица (7.3.2) - ${}^*{}_*$-невырожденная матрица.

Следовательно, если множество векторов $(f_i, i \in I, |I| = n)$ порождают базис $D_*{}^*$-векторного пространства и базис ${}^*{}_*D$-векторного пространства, их координатная матрица (7.3.2) является ${}_*{}^*$-невырожденной и ${}^*{}_*$-невырожденной матрицей. Утверждение следует из теоремы 5.8.9. □

Из теоремы 7.3.4 следует, что существует базис $\overline{\overline{e}}$ $D_*{}^*$-векторного пространства $V$, который не является базисом ${}^*{}_*D$-векторного пространства $V$.

## 7.4. $D$-векторное пространство

При изучении многих задач мы вполне можем ограничиться рассмотрением $D\star$-векторного пространства либо $\star D$-векторного пространства. Однако есть задачи, в которых мы вынуждены отказаться от простой модели и одновременно рассматривать обе структуры векторного пространства.

Парное представление тела $D$ в абелевой группе $V$ определяет структуру $D$ **-векторного пространства**. Из этого определения следует, что если $v, w \in V$, то для любых $a, b, c, d \in D$

$$avb + cwd \in V \tag{7.4.1}$$

Выражение (7.4.1) называется **линейной комбинацией векторов** $D$-векторного пространства $V$.

Определение 7.4.1. Векторы $a_i$, $i \in I$, $D$-векторного пространства $V$ **линейно независимы**, если из равенства

$$b^i a_i c^i = 0$$

следует $b^i c^i = 0$, $i = 1, ..., n$, В противном случае, векторы $a_i$ **линейно зависимы**. □

Определение 7.4.2. Мы называем множество векторов

$$\overline{\overline{e}} = (e_i, i \in I)$$

**базисом $D$-векторного пространства**, если векторы $e_i$ линейно независимы и, добавив к этой системе любого другого вектора, мы получим новую систему, которая является линейно зависимой. □

Теорема 7.4.3. *Базис $D$-векторного пространства $V$ является базисом $D^*{}_*$-векторного пространства $V$.*

Доказательство. Пусть векторы $a_i$, $i \in I$ линейно зависимы в $D^*{}_*$-векторном пространстве. Тогда существует ${}_*$-строка

$$\begin{pmatrix} c^i & i \in I \end{pmatrix} \neq 0 \tag{7.4.2}$$

такая, что

$$c^*{}_* a = 0 \tag{7.4.3}$$

Из равенства (7.4.3), следует, что

$$c^i a_i 1 = 0 \tag{7.4.4}$$

Согласно определению 7.4.1, векторы $a_i$ линейно зависимы в $D$-векторном пространстве.

Пусть $D^*{}_*$-размерность векторного пространства $V$ равна $n$. Пусть $\overline{\overline{e}}$ - базис $D^*{}_*$-векторного пространства. Пусть

$$\begin{pmatrix} 1 & ... & 0 \\ ... & ... & ... \\ 0 & ... & 1 \end{pmatrix}$$

координатная матрица базиса $\overline{\overline{e}}$. Очевидно, что мы не можем представить вектор $e_n$ как $D$-линейную комбинацию векторов $e_1$, ..., $e_{n-1}$. □

Для нас несущественно, является ли множество векторов $(e_i, i \in I)$ базисом ${}_*{}^*D$-векторного пространства или базисом $D^*{}_*$-векторного пространства. Поэтому при записи координат вектора мы индексы будем записывать в стандартном формате. Рассматривая базис $D$-векторного пространства $V$ мы, опираясь на теорему 7.4.3, будем рассматривать его как базис $D^*{}_*$-векторного пространства $V$.

Глава 8

# Произведение представлений

## 8.1. Бимодуль

Определение 8.1.1. $\overline{V}$ - **$(S\star, \star T)$-бимодуль**, если на множестве $\overline{V}$ определены структуры $S\star$-векторного пространства и $\star T$-векторного пространства.[8.1] □

Мы также пользуемся записью ${}_S\overline{V}_T$, когда мы хотим сказать, что $\overline{V}$ - $(S, T)$-бимодуль.

Пример 8.1.2. В разделе 7.4 мы рассмотрели $D$-векторное пространство, в котором можно определить структуру $(D^*{}_*, {}_*{}^*D)$-бимодуля. □

Пример 8.1.3. Множество $n\times m$ матриц порождает $(D\star, \star D)$-бимодуль. $D\star$ -базис можно представить в виде множества матриц ${}_i{}^je = (\delta^l_i\delta^j_k)$. □

Пример 8.1.4. Чтобы увидеть пример $(S_*{}^*, {}_*{}^*T)$-бимодуля, мы можем взять множество матриц

$$\begin{pmatrix} (A_1, {}^1B) & ... & (A_1, {}^mB)\\ ... & ... & ...\\ (A_n, {}^1B) & ... & (A_n, {}^mB) \end{pmatrix} \tag{8.1.1}$$

Формально мы можем представить эту матрицу в виде

$$\begin{pmatrix} (A_1, {}^1B) & ... & (A_1, {}^mB)\\ ... & ... & ...\\ (A_n, {}^1B) & ... & (A_n, {}^mB) \end{pmatrix} = \begin{pmatrix} A_1\\ ...\\ A_n \end{pmatrix} {}_*{}^* \begin{pmatrix} {}^1B & ... & {}^mB \end{pmatrix}$$

Это представление ясно показывает, что $S_*{}^*$-размерность $(S_*{}^*, {}_*{}^*T)$-бимодуля равна $n$ и ${}_*{}^*T$-размерность равна $m$. Однако в отличии от примера 8.1.2 мы не можем построить $S_*{}^*$-базис или ${}_*{}^*T$-базис, который порождал бы $(S_*{}^*, {}_*{}^*T)$-бимодуль. На самом деле базис $(S_*{}^*, {}_*{}^*T)$-бимодуля имеет вид $({}^a\overline{e}, \overline{f}_b)$, и любой вектор $(S_*{}^*, {}_*{}^*T)$-бимодуля имеет разложение

$$(\overline{A}, \overline{B}) = A_a({}^a\overline{e}, \overline{f}_b){}^bB \tag{8.1.2}$$

коэффициенты которого порождают матрицу (8.1.1). □

[8.1] Напомню, что в наших обозначениях $S\star$-векторное пространство означает левое $S$-векторное пространство и $\star T$-векторное пространство означает правое $T$-векторное пространство.

Мы можем записать равенство (8.1.2) в виде

$$(\overline{A},\overline{B}) = (A_*{}^*\overline{e},\overline{f}_*{}^*B) \tag{8.1.3}$$

Равенство (8.1.3) выглядит необычно, но оно становится понятным, если мы положим

$$(\overline{A},\overline{B}) = (\overline{A},\overline{0}) + (\overline{0},\overline{B})$$

Соглашение 5.2.6 справедливо для представления векторов ${}_*{}^*T$-векторного пространства независимо от представления векторов $S_*{}^*$-векторного пространства. Мы можем записать это равенство в виде

$$(\overline{A},\overline{B}) = (A^*{}_*\overline{e},B_*{}^*\overline{f}) \tag{8.1.4}$$

При этом координаты вектора $(\overline{A},\overline{B})$ не порождают матрицы. Однако это для нас сейчас не имеет значения.

Мы хотим записать компоненты $A$ и $B$ в равенстве (8.1.2) с одной стороны по отношению к векторам $({}^a\overline{e},\overline{f}_b)$. На первый взгляд это кажется невозможной задачей. Однако мы примем соглашение записывать это равенство в виде

$$(\overline{A},\overline{B}) = (A_a,B^b)({}^ae,{}_bf) = (A_*{}^*,B^*{}_*)(e,f)$$

Эта конструкция легко подаётся обобщению.

## 8.2. Прямое произведение тел

Определение 8.2.1. Пусть $\mathcal{A}$ - категория. Пусть $\{B_i, i\in I\}$ - множество объектов из $\mathcal{A}$. Объект

$$P=\prod_{i\in I}B_i$$

и множество морфизмов

$$\{\ f_i:P\longrightarrow B_i\ ,i\in I\}$$

называется **произведением объектов** $\{B_i, i\in I\}$ **в категории** $\mathcal{A}$[8.2], если для любого объекта $R$ и множество морфизмов

$$\{\ g_i:R\longrightarrow B_i\ ,i\in I\}$$

существует единственный морфизм

$$h:R\longrightarrow P$$

такой, что диаграмма

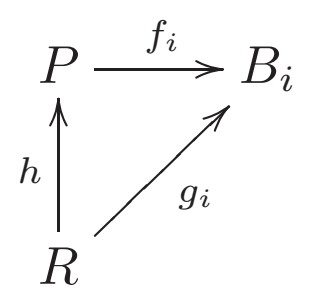

$$f_i\circ h=g_i$$

коммутативна для всех $i\in I$.

Если $|I|=n$, то для произведения объектов $\{B_i, i\in I\}$ в $\mathcal{A}$ мы так же будем пользоваться записью

$$G=\prod_{i=1}^nB_i=B_1\times ...\times B_n$$

□

[8.2] Определение дано согласно [1], страница 45.

Пример 8.2.2. Пусть $\mathcal{S}$ - категория множеств.[8.3] Согласно определению 8.2.1, декартово произведение

$$A = \prod_{i \in I} A_i$$

семейства множеств $(A_i, i \in I)$ и семейство проекций на $i$-й множитель

$$p_i : A \to A_i$$

являются произведением в категории $\mathcal{S}$. □

Теорема 8.2.3. *Пусть произведение существует в категории $\mathcal{A}$ $\Omega$-алгебр. Пусть $\Omega$-алгебра $A$ и семейство морфизмов*

$$p_i : A \to A_i \quad i \in I$$

*является произведением в категории $\mathcal{A}$. Тогда*

8.2.3.1: *Множество $A$ является декартовым произведением семейства множеств $(A_i, i \in I)$*

8.2.3.2: *Гомоморфизм $\Omega$-алгебры*

$$p_i : A \to A_i$$

*является проекцией на $i$-й множитель.*

8.2.3.3: *Любое $A$-число $a$ может быть однозначно представлено в виде кортежа $(p_i(a), i \in I)$ $A_i$-чисел.*

8.2.3.4: *Пусть $\omega \in \Omega$ - $n$-арная операция. Тогда операция $\omega$ определена покомпонентно*

(8.2.1) $$a_1 ... a_n \omega = (a_{1i} ... a_{ni} \omega, i \in I)$$

*где* $a_1 = (a_{1i}, i \in I)$, ..., $a_n = (a_{ni}, i \in I)$ .

Доказательство. Пусть

$$A = \prod_{i \in I} A_i$$

декартово произведение семейства множеств $(A_i, i \in I)$ и, для каждого $i \in I$, отображение

$$p_i : A \to A_i$$

является проекцией на $i$-й множитель. Рассмотрим диаграмму морфизмов в категории множеств $\mathcal{S}$

(8.2.2)

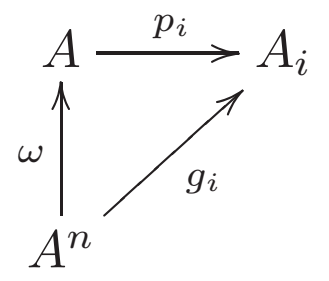

$$p_i \circ \omega = g_i$$

где отображение $g_i$ определено равенством

$$g_i(a_1, ..., a_n) = p_i(a_1) ... p_i(a_n) \omega$$

Согласно определению 8.2.1, отображение $\omega$ определено однозначно из множества диаграмм (8.2.2)

(8.2.3) $$a_1 ... a_n \omega = (p_i(a_1) ... p_i(a_n) \omega, i \in I)$$

Равенство (8.2.1) является следствием равенства (8.2.3). □

[8.3] Смотри также пример в [1], страница 45.

Определение 8.2.4. Пусть произведение существует в категории $\Omega$-алгебр. Если $\Omega$-алгебра $A$ и семейство морфизмов

$$p_i : A \to A_i \quad i \in I$$

является произведением в категории $\mathcal{A}$, то $\Omega$-алгебра $A$ называется **прямым** или **декартовым произведением $\Omega$-алгебр** $(A_i, i \in I)$ . □

Пример 8.2.5. Пусть $\{G_i, i \in I\}$ - множество групп. Пусть

$$G = \prod_{i \in I} G_i$$

декартово произведение множеств $G_i$, $i \in I$. Мы определим на $G$ групповую структуру посредством покомпонентного умножения. Если $x = (x_i, i \in I) \in D$ и $y = (y_i, i \in I) \in D$, мы определим их произведение

$$xy = (x_i y_i, i \in I)$$

Если $x = (x_i, i \in I) \in D$, мы определим обратный элемент

$$x^{-1} = (x_i^{-1}, i \in I)$$

Множество $G$ называется **декартовым произведением групп**[8.4] $G_i$, $i \in I$.

Если $|I| = n$, то для декартова произведения групп $G_1$, ..., $G_n$ мы так же будем пользоваться записью

$$G = \prod_{i=1}^{n} G_i = G_1 \times ... \times G_n$$

□

Теорема 8.2.6. *Декартово произведение абелевых групп является абелевой группой.*

Доказательство. Пусть $x = (x_i, i \in I) \in G$ и $y = (y_i, i \in I) \in G$. Тогда

$$x + y = (x_i + y_i, i \in I) = (y_i + x_i, i \in I) = y + x$$

□

Пример 8.2.7. Пусть $\{D_i, i \in I\}$ - множество тел. Пусть

$$D = \prod_{i \in I} D_i$$

произведение аддитивных групп тел $D_i$, $i \in I$. Если $x = (x_i, i \in I) \in D$ и $y = (y_i, i \in I) \in D$, мы определим их произведение покомпонентно

$$xy = (x_i y_i, i \in I)$$

Мультипликативная единица - это $e = (e_i, i \in I) \in D$, где $e_i$, $i \in I$ - мультипликативная единица тела $D_i$. Множество $D$ называется **прямым произведением тел**[8.5] $D_i$, $i \in I$.

Если $|I| = n$, то для прямого произведения тел $D_1$, ..., $D_n$ мы так же будем пользоваться записью

$$D = \prod_{i=1}^{n} D_i = D_1 \times ... \times D_n$$

[8.4] Определение дано согласно [1], страницы 47, 48.

[8.5] Определение дано согласно предложению 1 из [1], стр. 79

□

Прямое произведение тел $D_i$, $i \in I$, вообще говоря, не является телом. Например, пусть $x_1 \in D_1$, $x_1 \neq 0$, $x_2 \in D_2$, $x_2 \neq 0$. Тогда

$$(x_1, 0)(0, x_2) = (0, 0)$$

Однако прямое произведение тел является кольцом. Следовательно, прямое произведение не определенно в категории тел, но прямое произведение тел определено в категории колец.

Теорема 8.2.8. *Пусть множество $A$ является декартовым произведением множеств $(A_i, i \in I)$ и множество $B$ является декартовым произведением множеств $(B_i, i \in I)$. Для каждого $i \in I$, пусть*

$$f_i : A_i \to B_i$$

*является отображением множества $A_i$ в множество $B_i$. Для каждого $i \in I$, рассмотрим коммутативную диаграмму*

(8.2.4)

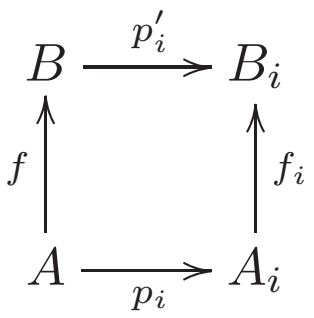

*где отображения $p_i$, $p'_i$ являются проекцией на $i$-й множитель. Множество коммутативных диаграмм (8.2.4) однозначно определяет отображение*

$$f : A \to B$$
$$f(a_i, i \in I) = (f_i(a_i), i \in I)$$

Доказательство. Для каждого $i \in I$, рассмотрим коммутативную диаграмму

(8.2.5)

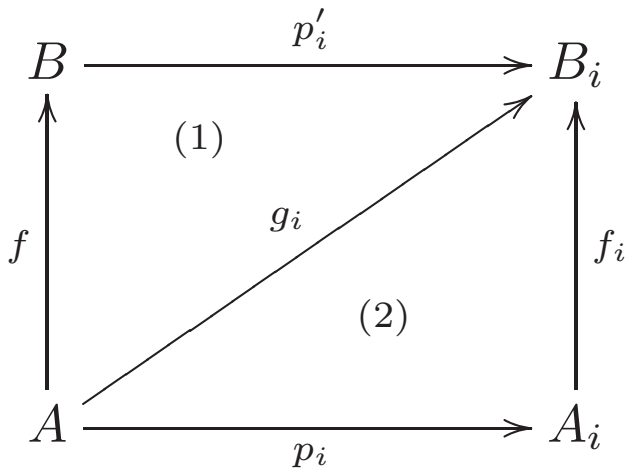

Пусть $a \in A$. Согласно утверждению 8.2.3.3, $A$-число $a$ может быть представлено в виде кортежа $A_i$-чисел

$$a = (a_i, i \in I) \quad a_i = p_i(a) \in A_i \tag{8.2.6}$$

Пусть

$$b = f(a) \in B \tag{8.2.7}$$

Согласно утверждению 8.2.3.3, $B$-число $b$ может быть представлено в виде кортежа $B_i$-чисел

$$b = (b_i, i \in I) \quad b_i = p'_i(b) \in B_i \tag{8.2.8}$$

Из коммутативности диаграммы (1) и из равенств (8.2.7), (8.2.8) следует, что

$$b_i = g_i(b) \tag{8.2.9}$$

Из коммутативности диаграммы (2) и из равенства (8.2.6) следует, что

$$b_i = f_i(a_i)$$

□

Теорема 8.2.9. *Пусть Ω-алгебра $A$ является декартовым произведением Ω-алгебр $(A_i, i \in I)$ и Ω-алгебра $B$ является декартовым произведением Ω-алгебр $(B_i, i \in I)$. Для каждого $i \in I$, пусть отображение*

$$f_i : A_i \to B_i$$

*является гомоморфизмом Ω-алгебры. Тогда отображение*

$$f : A \to B$$

*определённое равенством*

$$f(a_i, i \in I) = (f_i(a_i), i \in I) \tag{8.2.10}$$

*является гомоморфизмом Ω-алгебры.*

Доказательство. Пусть $\omega \in \Omega$ - n-арная операция. Пусть $a_1 = (a_{1i}, i \in I)$, ..., $a_n = (a_{ni}, i \in I)$ и $b_1 = (b_{1i}, i \in I)$, ..., $b_n = (b_{ni}, i \in I)$. Из равенств (8.2.1), (8.2.10) следует, что

$$\begin{aligned} f(a_1...a_n\omega) &= f(a_{1i}...a_{ni}\omega, i \in I) \\ &= (f_i(a_{1i}...a_{ni}\omega), i \in I) \\ &= ((f_i(a_{1i}))...(f_i(a_{ni})), i \in I) \\ &= (b_{1i}...b_{ni}\omega, i \in I) \end{aligned}$$

$$f(a_1)...f(a_n)\omega = b_1...b_n\omega \qquad = (b_{1i}...b_{ni}\omega, i \in I)$$

□

Замечание 8.2.10. Утверждение, обратное утверждениям теорем 8.2.8, 8.2.9, вообще говоря, неверно. Рассмотрим, например, $R$-векторное пространство $V$ размерности 2 и линейное отображение

$$f : V \to V$$

$$x'^1 = x^1$$

$$x'^2 = x^1 + x^2$$

Векторное пространство $V$ - это абелева группа, которая является декартовым произведением абелевых групп

$$V = R \times R$$

Отображение $f$ является гомоморфизмом абелевой группы $V$. Однако соответствие

$$x^2 \to x'^2$$

зависит от значения переменной $x^1$ и поэтому не является отображением. □

## 8.3. Прямое произведение $D_*{}^*$-векторных пространств

Лемма 8.3.1. Пусть
$$A=\prod_{i\in I}A_i$$
декартово произведение семейства $\Omega_1$-алгебр $(A_i,i\in I)$ . Для каждого $i\in I$, пусть множество ${}^*A_i$ является $\Omega_2$-алгеброй. Тогда множество

(8.3.1)
$$^\circ A=\{f\in {}^*A:f(a_i,i\in I)=(f_i(a_i),i\in I)\}$$

является декартовым произведением $\Omega_2$-алгебр ${}^*A_i$.

Доказательство. Согласно определению (8.3.1), мы можем представить отображение $f\in {}^\circ A$ в виде кортежа
$$f=(f_i,i\in I)$$
отображений $f_i\in {}^*A_i$. Согласно определению (8.3.1),
$$(f_i,i\in I)(a_i,i\in I)=(f_i(a_i),i\in I)$$
Пусть $\omega\in\Omega$ - n-арная операция. Мы определим операцию $\omega$ на множестве $^\circ A$ равенством
$$((f_{1i},i\in I)...(f_{ni},i\in I)\omega)(a_i,i\in I)=((f_{1i}(a_i))...(f_{ni}(a_i))\omega,i\in I)$$
□

Теорема 8.3.2. *Допустим категория $\mathcal{A}_1$ $\Omega_1$-алгебр имеет произведение. Допустим категория $\mathcal{A}_2$ $\Omega_2$-алгебр имеет произведение. Тогда в категории*[8.6] *$(\mathcal{A}_1*)\mathcal{A}_2$ существует произведение однотранзитивных левосторонних представлений $\Omega_1$-алгебры в $\Omega_2$-алгебре.*

Доказательство. Для $j=1,\ 2,$ пусть
$$P_j=\prod_{i\in I}B_{ji}$$
произведение семейства $\Omega_j$-алгебр $\{B_{ji},i\in I\}$ и для любого $i\in I$ отображение
$$t_{ji}:P_j\longrightarrow B_{ji}$$
является проекцией на множитель $i$. Для каждого $i\in I$, пусть
$$h_i:B_{1i}\overset{*}{\longrightarrow} B_{2i}$$
однотранзитивное $B_{1i}*$-представление в $\Omega_2$-алгебре $B_{2i}$.

Пусть $b_1\in P_1$. Согласно утверждению 8.2.3.3, $P_1$-число $b_1$ может быть представлено в виде кортежа $B_{1i}$-чисел

(8.3.2)
$$b_1=(b_{1i},i\in I)\quad b_{1i}=t_{1i}(b_1)\in B_{1i}$$

Пусть $b_2\in P_2$. Согласно утверждению 8.2.3.3, $P_2$-число $b_2$ может быть представлено в виде кортежа $B_{2i}$-чисел

(8.3.3)
$$b_2=(b_{2i},i\in I)\quad b_{2i}=t_{2i}(b_2)\in B_{2i}$$

[8.6] Смотри определение 3.2.13

Лемма 8.3.3. Для каждого $i\in I$, рассмотрим диаграмму отображений

(8.3.4)

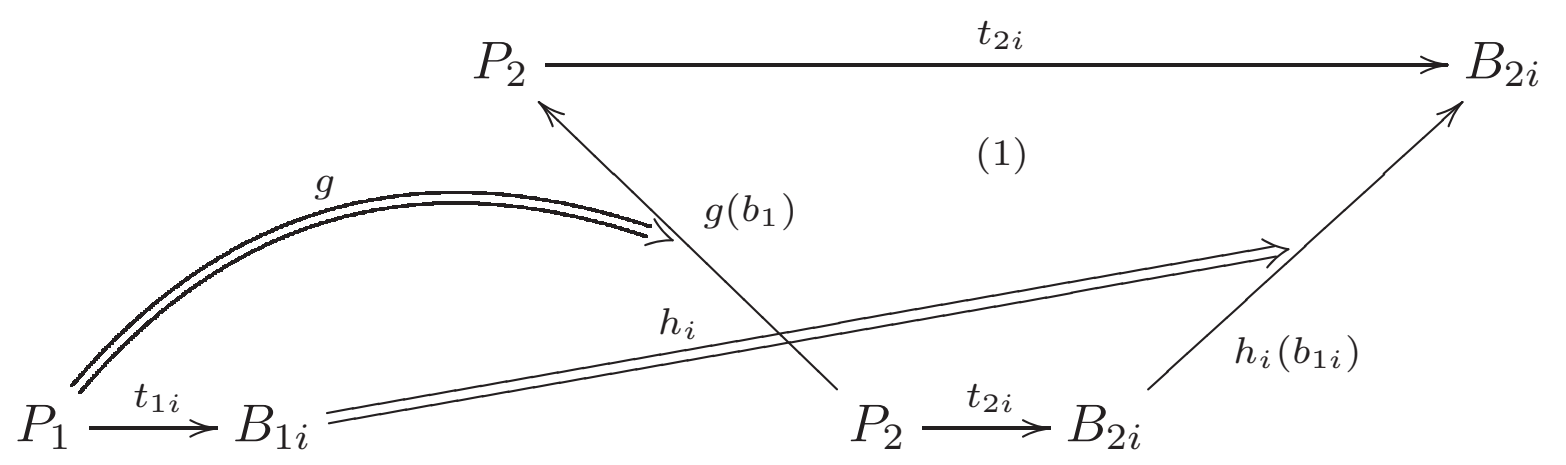

Пусть отображение

$$g:P_1\to {}^*P_2$$

определено равенством

$$g(b_1)(b_2)=(h_i(b_{1i})(b_{2i}),i\in I) \tag{8.3.5}$$

Тогда отображение $g$ является однотранзитивным $P_1*$-представлением в $\Omega_2$-алгебре $P_2$

$$g:P_1 \xrightarrow{*} P_2$$

Отображение $(t_{1i},t_{2i})$ является морфизмом представления $g$ в представление $h_i$.

Доказательство.

8.3.3.1: Согласно определениям 3.1.1, 3.1.2, отображение $h_i(b_{1i})$ является гомоморфизмом $\Omega_2$-алгебры $B_{2i}$. Согласно теореме 8.2.9, из коммутативности диаграммы (1) для каждого $i\in I$, следует, что отображение

$$g(b_1):P_2\to P_2$$

определённое равенством (8.3.5) является гомоморфизмом $\Omega_2$-алгебры $P_2$.

8.3.3.2: Согласно определению 3.1.2, множество ${}^*B_{2i}$ является $\Omega_1$-алгеброй. Согласно лемме 8.3.1, множество ${}^\circ P_2\subseteq {}^*P_2$ является $\Omega_1$-алгеброй.

8.3.3.3: Согласно определению 3.1.2, отображение

$$h_i:B_{1i}\to {}^*B_{2i}$$

является гомоморфизмом $\Omega_1$-алгебры. Согласно теореме 8.2.9, отображение

$$g:P_1\to {}^*P_2$$

определённое равенством

$$g(b_1)=(h_i(b_{1i}),i\in I) \tag{8.3.6}$$

является гомоморфизмом $\Omega_1$-алгебры.

Согласно утверждениям 8.3.3.1, 8.3.3.3 и определению 3.1.2, отображение $g$ является $P_1*$-представлением в $\Omega_2$-алгебре $P_2$.

Пусть $b_{21},b_{22}\in P_2$. Согласно утверждению 8.2.3.3, $P_2$-числа $b_{21},b_{22}$ может быть представлено в виде кортежей $B_{2i}$-чисел

$$\begin{array}{ll} b_{21}=(b_{21i},i\in I) & b_{21i}=t_{2i}(b_{21})\in B_{2i}\\ b_{22}=(b_{22i},i\in I) & b_{22i}=t_{2i}(b_{22})\in B_{2i}\end{array} \tag{8.3.7}$$

Согласно теореме 3.1.9, поскольку представление $h_i$ однотранзитивно, то существует единственное $B_{1i}$-число $b_{1i}$ такое, что

$$b_{22i} = h_i(b_{1i})(b_{21i})$$

Согласно определениям (8.3.2), (8.3.5), (8.3.7), существует единственное $P_1$-число $b_1$ такое, что

$$b_{22} = g(b_1)(b_{21})$$

Согласно теореме 3.1.9, представление $g$ однотранзитивно.

Из коммутативности диаграммы (1) и определения 3.2.2 следует, что отображение $(t_{1i}, t_{2i})$ является морфизмом представления $g$ в представление $h_i$. $\odot$

Пусть

(8.3.8) $$d_2 = g(b_1)(b_2) \quad d_2 = (d_{2i}, i \in I)$$

Из равенств (8.3.5), (8.3.8) следует, что

(8.3.9) $$d_{2i} = h_i(b_{1i})(b_{2i})$$

Для $j = 1, 2,$ пусть $R_j$ другой объект категории $\mathcal{A}_j$. Для любого $i \in I$, пусть

$$r_{1i} : R_1 \longrightarrow B_{1i}$$

морфизм из $\Omega_1$-алгебра $R_1$ в $\Omega_1$-алгебру $B_{1i}$. Согласно определению 8.2.1, существует единственный морфизм $\Omega_1$-алгебры

$$s_1 : R_1 \longrightarrow P_1$$

такой, что коммутативна диаграмма

(8.3.10)

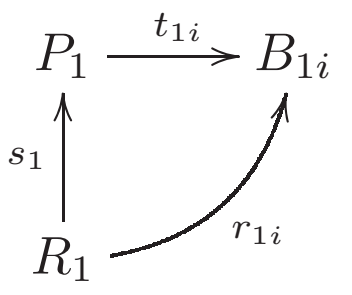

$$t_{1i} \circ s_1 = r_{1i}$$

Пусть $a_1 \in R_1$. Пусть

(8.3.11) $$b_1 = s_1(a_1) \in P_1$$

Из коммутативности диаграммы (8.3.10) и утверждений (8.3.11), (8.3.2) следует, что

(8.3.12) $$b_{1i} = r_{1i}(a_1)$$

Пусть

$$f : R_1 \overset{*}{\longrightarrow} R_2$$

однотранзитивное $R_1*$-представление в $\Omega_2$-алгебре $R_2$. Согласно теореме 3.2.10, морфизм $\Omega_2$-алгебры

$$r_{2i} : R_2 \longrightarrow B_{2i}$$

такой, что отображение $(r_{1i}, r_{2i})$ является морфизмом представлений из $f$ в $h_i$, определён однозначно с точностью до выбора образа $R_2$-числа $a_2$. Согласно замечанию 3.2.6, в диаграмме отображений

(8.3.13)

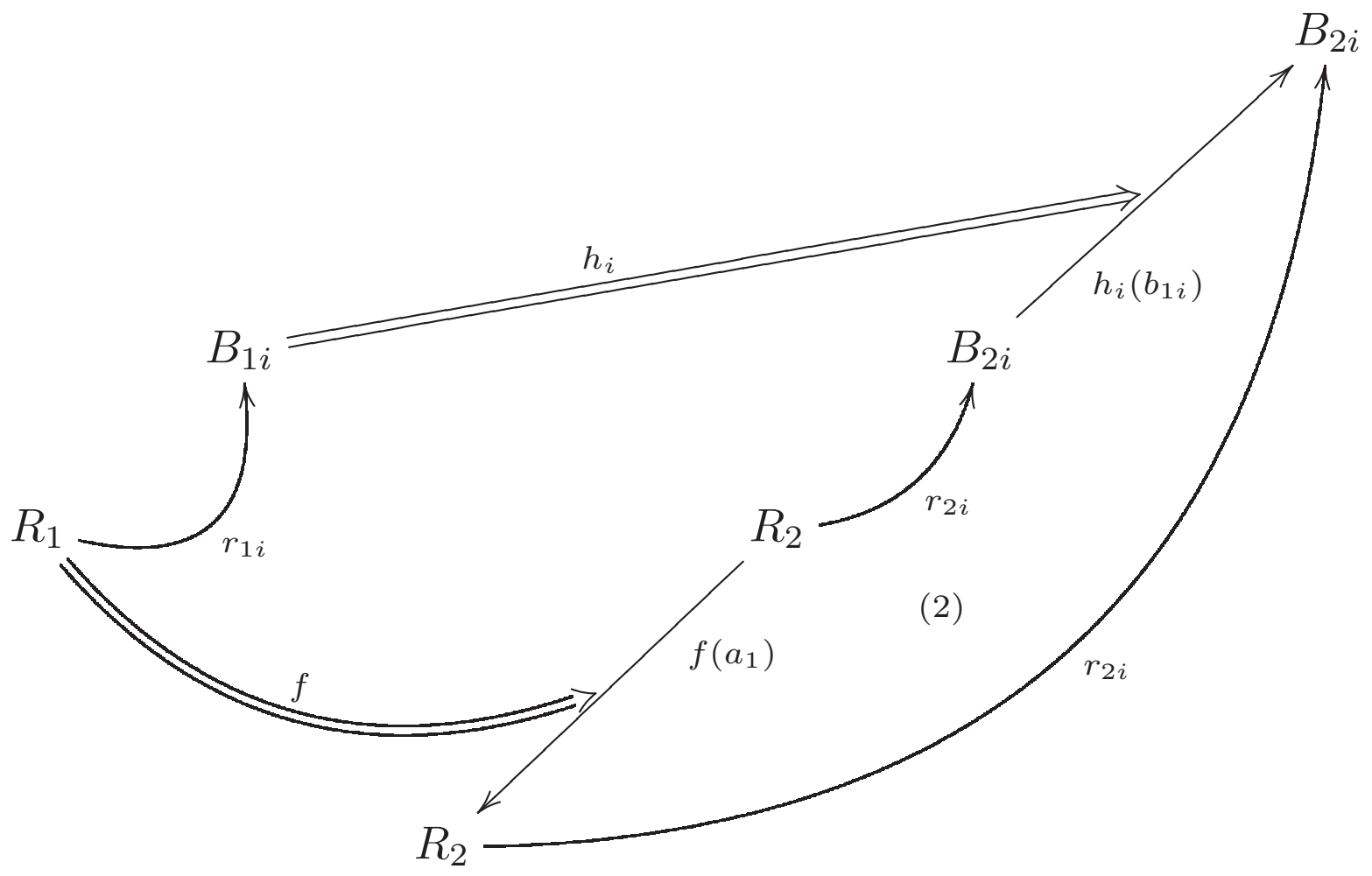

диаграмма (2) коммутативна. Согласно определению 8.2.1, существует единственный морфизм $\Omega_2$-алгебры

$$s_2 : R_2 \longrightarrow P_2$$

такой, что коммутативна диаграмма

(8.3.14)

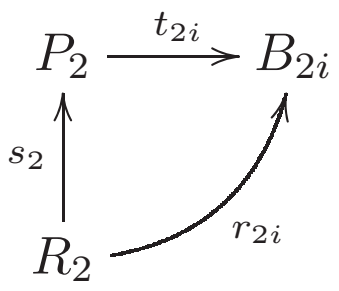

$$t_{2i} \circ s_2 = r_{2i}$$

Пусть $a_2 \in R_2$. Пусть

$$b_2 = s_2(a_2) \in P_2 \tag{8.3.15}$$

Из коммутативности диаграммы (8.3.14) и утверждений (8.3.15), (8.3.3) следует, что

$$b_{2i} = r_{2i}(a_2) \tag{8.3.16}$$

Пусть

$$c_2 = f(a_1)(a_2) \tag{8.3.17}$$

Из коммутативности диаграммы (2) и равенств (8.3.9), (8.3.16), (8.3.17) следует, что

$$d_{2i} = r_{2i}(c_2) \tag{8.3.18}$$

Из равенств (8.3.9), (8.3.18) следует, что

$$d_2 = s_2(c_2) \tag{8.3.19}$$

что согласуется с коммутативносью диаграмы (8.3.14).

Для каждого $i \in I$, мы объединим диаграммы отображений (8.3.4), (8.3.10), (8.3.14), (8.3.13)

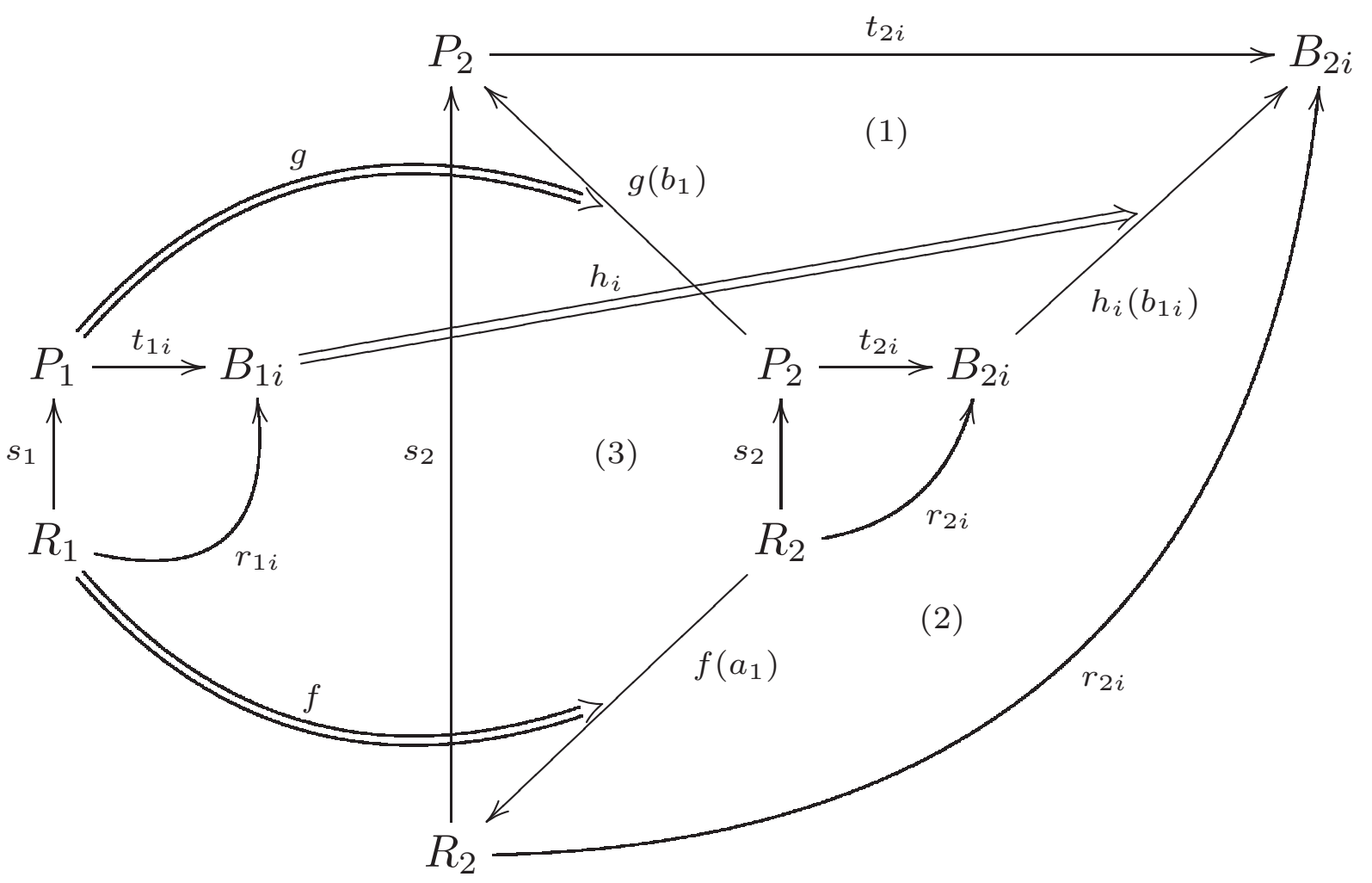

Из равенств (8.3.8) (8.3.15) и из равенств (8.3.17), (8.3.19), следует коммутативность диаграммы (3). Следовательно, отображение $(s_1, s_2)$ является морфизмом представлений из $f$ в $g$., Согласно теореме 3.2.10, морфизм $(s_1, s_2)$ определён однозначно, так как мы требуем (8.3.19).

Согласно определению 8.2.1, представление $g$ и семейство морфизмов представления $((t_{1i}, t_{2i}), i \in I)$ является произведением в категории $(\mathcal{A}_1*)\mathcal{A}_2$ . □

Теорема 8.3.4. *В категории $\mathcal{G}*$ левосторонних представлений группы существует произведение эффективных левосторонних представлений группы.*

Доказательство. Утверждение теоремы является следствием теорем 8.3.2 и 4.2.10 □

Определение 8.3.5. Пусть $\{V_i, i \in I\}$ - множество ${}_*{}^*D_i$-векторных пространств. Левостороннее представление кольца $D = \prod_{i\in I} D_i$ в абелевой группе $V = \prod_{i\in I} V_i$ называется **прямым произведением**

$$V = \prod_{i\in I} V_i$$

**${}_*{}^*D_i$-векторных пространств**[8.7] $V_i$, $i \in I$. □

Если $|I| = n$, то для прямого произведения ${}_*{}^*D_i$-векторных пространств $V_1$, ..., $V_n$ мы так же будем пользоваться записью

$$V = \prod_{i=1}^{n} V_i = V_1 \times ... \times V_n$$

[8.7] Определение дано согласно [1], стр. 98

Согласно определению представление кольца $D$ в абелевой группе $V$ определенно покомпонентно. Если $a=(a_i,i\in I)\in D$ и $v=(v_i,i\in I)\in V$, мы определим правостороннее представление, соответствующее элементу $a$

$$va=(v_ia_i,i\in I)$$

Мы рассматриваем соглашение, описанное в замечании 2.2.15, покомпонентно. Линейную комбинацию ${}_*{}^*D$-векторов $\overline{v}_k$ мы будем записывать в виде

$$\overline{v}_*{}^*a=(\overline{v}_{i*}{}^*a_i,i\in I)=(v_{i\cdot k}\ {}^ka_i,i\in I)$$

Мы так же будем пользоваться записью

$$\overline{v}_*{}^*a=(v_{i\cdot *},i\in I)\ ({}^*a_i,i\in I)=(v_{i\cdot k},i\in I)\ ({}^ka_i,i\in I) \tag{8.3.20}$$

Если $|I|=n$, то равенство (8.3.20) можно записать в виде

$$\overline{v}_*{}^*a=(v_{1\cdot *},...,v_{n\cdot *})\ ({}^*a_1,...,{}^*a_n)=(v_{1\cdot k_1},...,v_{n\cdot k_n})\ ({}^{k_1}a_1,...,{}^{k_n}a_n)$$

Так как прямое произведение тел не является телом, то прямое произведение ${}_*{}^*D_i$-векторных пространств, вообще говоря, является модулем. Однако структура этого модуля близка к структуре ${}_*{}^*D$-векторного пространства.

Теорема 8.3.6. *Пусть $V_1$, ..., $V_n$ - множество ${}_*{}^*D_i$-векторных пространств. Пусть $\overline{\overline{e}}_i$ - ${}_*{}^*D_i$-базис векторного пространства $V_i$. Множество векторов $(e_{1\cdot i_1},...,e_{n\cdot i_n})$ порождает базис прямого произведения*

$$V=V_1\times...\times V_n$$

Доказательство. Мы докажем утверждение теоремы индукцией по $n$.

При $n=1$ утверждение очевидно.

Допустим утверждение справедливо при $n=k-1$. Мы можем представить пространство $V$ в виде[8.8]

$$V=V_1\times...\times V_k=(V_1\times...\times V_{k-1})\times V_k$$

Соответственно, произвольный вектор $(\overline{v}_1,...,\overline{v}_k)\in V$ можно представить в виде

$$(\overline{v}_1,...,\overline{v}_k)=((\overline{v}_1,...,\overline{v}_{k-1}),\overline{v}_k)=((\overline{v}_1,...,\overline{v}_{k-1}),0)+(0,\overline{v}_k) \tag{8.3.21}$$

$(\overline{v}_1,...,\overline{v}_{k-1})\in V_1\times...\times V_{k-1}$ и согласно предположению индукции

$$(\overline{v}_1,...,\overline{v}_{k-1})=(e_{1\cdot *},...,e_{k-1\cdot *})\ ({}^*v_1,...,{}^*v_{k-1}) \tag{8.3.22}$$

Из равенств (8.3.21) и (8.3.22) следует

$$\begin{aligned}(\overline{v}_1,...,\overline{v}_k)&=((e_{1\cdot *},...,e_{k-1\cdot *})\ ({}^*v_1,...,{}^*v_{k-1}),0)+(0,e_{k*}{}^*v_k)\\&=(e_{1\cdot *},...,e_{k-1\cdot *},e_k)\ ({}^*v_1,...,{}^*v_{k-1},0)+(e_1,...,e_{k-1},e_{k\cdot *})\ (0,...,{}^*v_k)\\&=(e_{1\cdot *},...,e_{k-1\cdot *},e_{k\cdot *})\ ({}^*v_1,...,{}^*v_{k-1},{}^*v_k)\end{aligned}$$

Следовательно утверждение справедливо при $n=k$. □

Если $V_1$, ..., $V_n$ - $D$-векторные пространства, то прямое произведение называется **прямым произведением $D$-векторных пространств**.

---

[8.8] Вообще говоря, мы должны рассматривать более широкий класс модулей вида

$$V_1\times...\times V_k$$

включая ассоциативность произведения.

Теорема 8.3.7. *Допустим категория $\mathcal{A}$ $\Omega$-алгебр имеет произведение. Тогда в категории $A*$*[8.9]*, где $A$ - $\Omega$-алгебра, существует произведение однотранзитивных $A*$-представлений.*

Доказательство. Доказательство теоремы похоже на доказательство теоремы 8.3.2 с той разницей, что мы пользуемся диаграммой

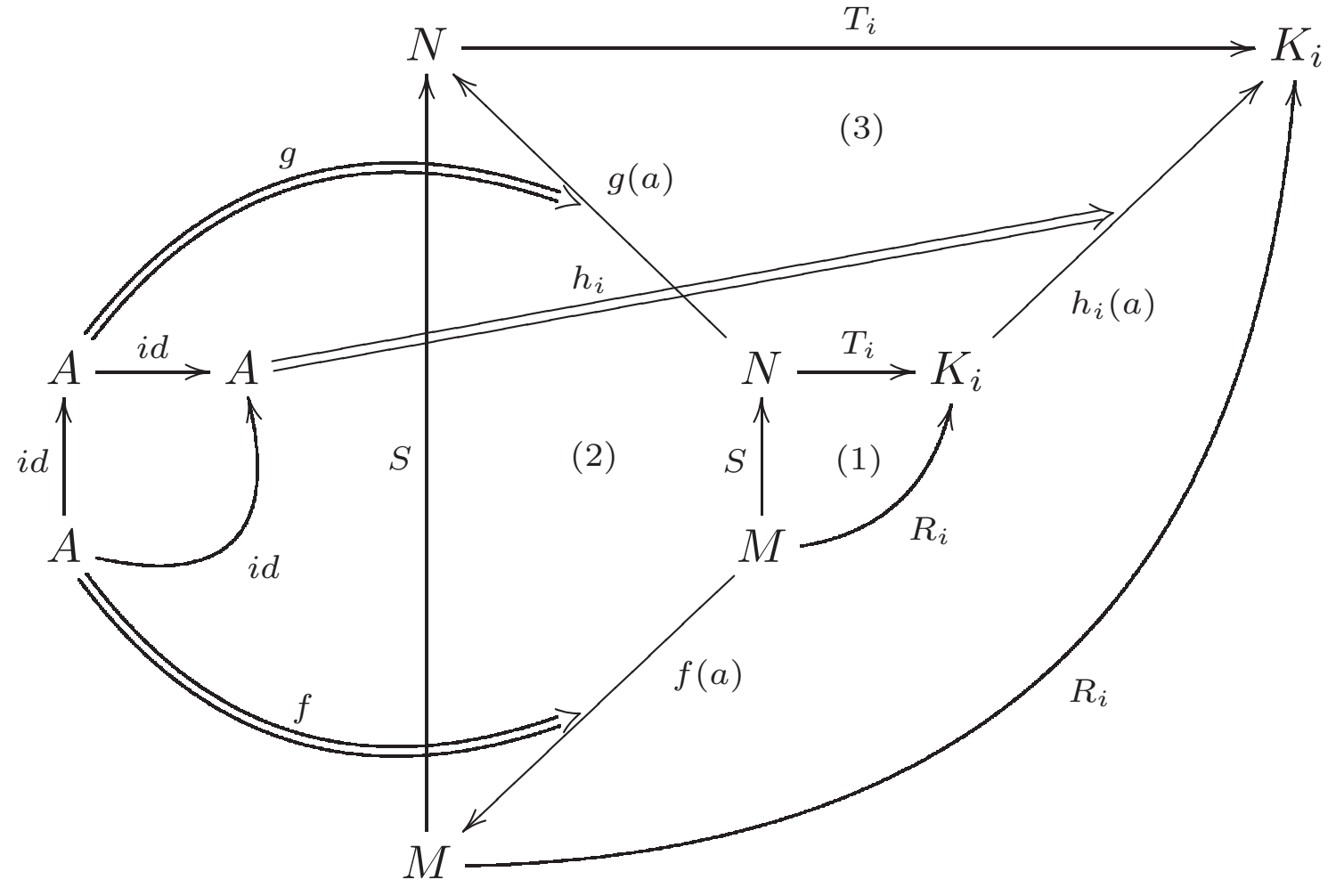

□

Рассмотрим связь между произведением однотранзитивных левосторонних представлений $\Omega$-алгебры $A$ в категории $\mathcal{A}*$ и произведением однотранзитивных левосторонних представлений $\Omega$-алгебры $A$ в категории $A*$. Пусть для каждого $i \in I$ определено представление $\Omega$-алгебры $A$ на множестве $K_i$. Для каждого $i \in I$ рассмотрим диаграмму где $g$ - отображение алгебры $A$ на диагональ в $\prod\limits_{i\in I} A$ и $f_i$ - проекция на множитель $i$. Из диаграммы следует, что отображения $g$ и $f_i$ инъективны. Следовательно оба произведения эквивалентны.

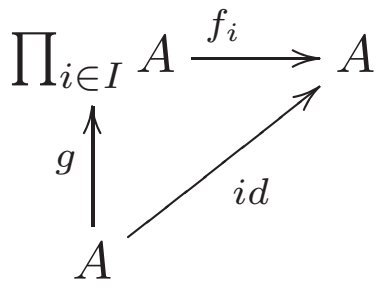

## 8.4. Морфизмы прямого произведения ${}_*{}^*D$-векторных пространств

Пусть $V_i$, $i \in I$, $W_j$, $j \in J$ - ${}_*{}^*D$-векторные пространства. Положим

$$V = \prod_{i\in I} V_i$$

$$W = \prod_{j\in J} W_j$$

[8.9] Смотри определение 3.2.20

Линейное отображение

$$f : V \longrightarrow W \tag{8.4.1}$$

сохраняет линейные операции. Так как операции в модуле $W$ мы выполняем покомпонентно, то мы можем представить отображение $f$ в виде

$$f = (\ f_j : V \longrightarrow W_j, j \in J\ )$$

Пусть

$$I_i : V_i \longrightarrow V \qquad i \in I$$

вложение $V_i$ в $V$. Тогда для любого вектора $v = (v_i, i \in I) \in V$, справедливо

$$f_j(v) = \sum_{i \in I} f_{j*}{}^* I_{i*}{}^* v_i$$

Следовательно, отображение (8.4.1) можно представить в виде матрицы ${}_*{}^*D$-линейных отображений

$$f = (f_{ij} = f_{j*}{}^* I_i : \ V_i \longrightarrow W_j\ , i \in I, j \in J)$$

Глава 9

# Геометрия тела

## 9.1. Центр тела

Определение 9.1.1. Пусть $D$ - кольцо.[9.1] Множество $Z(D)$ элементов $a \in D$ таких, что

(9.1.1) $$ax = xa$$

для всех $x \in D$, называется **центром кольца** $D$. □

Теорема 9.1.2. *Центр $Z(D)$ кольца $D$ является подкольцом кольца $D$.*

Доказательство. Непосредственно следует из определения 9.1.1. □

Определение 9.1.3. Пусть $D$ - кольцо с единицей $e$.[9.2] Отображение

$$l : Z \to D$$

для которого $l(n) = ne$ будет гомоморфизмом колец, и его ядро является идеалом $(n)$, порождённым целым числом $n \geq 0$. Канонический инъективный гомоморфизм

$$Z/nZ \to D$$

является изоморфизмом между $Z/nZ$ и подкольцом в $D$. Если $nZ$ - простой идеал, то у нас возникает два случая.

- $n = 0$. $D$ содержит в качестве подкольца кольцо, изоморфное $Z$ и часто отождествляемое с $Z$. В этом случае мы говорим, что $D$ имеет **характеристику** $0$.
- $n = p$ для некоторого простого числа $p$. $D$ имеет **характеристику** $p$, и $D$ содержит изоморфный образ $F_p = Z/pZ$.

□

Теорема 9.1.4. *Пусть $D$ - кольцо характеристики $0$ и пусть $d \in D$. Тогда любое целое число $n \in Z$ коммутирует с $d$.*

Доказательство. Утверждение теоремы доказывается по индукции. При $n = 0$ и $n = 1$ утверждение очевидно. Допустим утверждение справедливо при $n = k$. Из цепочки равенств

$$(k+1)d = kd + d = dk + d = d(k+1)$$

следует очевидность утверждения при $n = k + 1$. □

Теорема 9.1.5. *Пусть $D$ - кольцо характеристики $0$. Тогда кольцо целых чисел $Z$ является подкольцом центра $Z(D)$ кольца $D$.*

---

9.1[1], стр. 84.

9.2Определение дано согласно определению из [1], стр. 84, 85.

Доказательство. Следствие теоремы 9.1.4. □

Пусть $D$ - тело. Если $D$ имеет характеристику 0, $D$ содержит в качестве подполя изоморфный образ поля $Q$ рациональных чисел. Если $D$ имеет характеристику $p$, $D$ содержит в качестве подполя изоморфный образ $F_p$. В обоих случаях это подполе будет называться простым полем. Так как простое поле является наименьшим подполем в $D$, содержащим 1 и не имеющим автоморфизмов, кроме тождественного, его обычно отождествляют с $Q$ или $F_p$, в зависимости от того, какой случай имеет место.

Теорема 9.1.6. *Центр $Z(D)$ тела $D$ является подполем тела $D$.*

Доказательство. Согласно теореме 9.1.2 достаточно проверить, что $a^{-1} \in Z(D)$, если $a \in Z(D)$. Пусть $a \in Z(D)$. Многократно применяя равенство (9.1.1), мы получим цепочку равенств

$$aa^{-1}x = x = xaa^{-1} = axa^{-1} \tag{9.1.2}$$

Из (9.1.2) следует

$$a^{-1}x = xa^{-1}$$

Следовательно, $a^{-1} \in Z(D)$.

□

Теорема 9.1.7. *Пусть $D$ - тело характеристики 0 и пусть $d \in D$. Тогда для любого целого числа $n \in Z$*

$$n^{-1}d = dn^{-1} \tag{9.1.3}$$

Доказательство. Согласно теореме 9.1.4 справедлива цепочка равенств

$$n^{-1}dn = nn^{-1}d = d \tag{9.1.4}$$

Умножив правую и левую части равенства (9.1.4) на $n^{-1}$, получим

$$n^{-1}d = n^{-1}dnn^{-1} = dn^{-1} \tag{9.1.5}$$

Из (9.1.5) следует (9.1.3). □

Теорема 9.1.8. *Пусть $D$ - тело характеристики 0 и пусть $d \in D$. Тогда любое рациональное число $p \in Q$ коммутирует с $d$.*

Доказательство. Мы можем представить рациональное число $p \in Q$ в виде $p = mn^{-1}$, $m$, $n \in Z$. Утверждение теоремы следует из цепочки равенств

$$pd = mn^{-1}d = n^{-1}dm = dmn^{-1} = dp$$

основанной на утверждении теоремы 9.1.4 и равенстве (9.1.3). □

Теорема 9.1.9. *Пусть $D$ - тело характеристики 0. Тогда поле рациональных чисел $Q$ является подполем центра $Z(D)$ тела $D$.*

Доказательство. Следствие теоремы 9.1.8. □

## 9.2. Геометрия тела над полем

Мы можем рассматривать тело $D$ как векторное пространство над полем $F \subset Z(D)$. Так как $F$ - поле, то мы можем писать все индексы справа от корневой буквы. При этом мы будем пользоваться следующими соглашением 1.5.5.

Замечание 9.2.1. Пусть $\overline{\overline{e}}$ - базис тела $D$ над полем $F$. Тогда произвольный элемент $a \in D$ можно представить в виде

$$a = e_i a^i \quad a^i \in F \tag{9.2.1}$$

Если размерность тела $D$ над полем $F$ бесконечна, то базис может быть либо счётным, либо его мощность может быть не меньше, чем мощность континуума. Если базис счётный, то на коэффициенты $a^i$ разложения (9.2.1) накладываются определённые ограничения. Если мощность множества $I$ континуум, то предполагается, что на множестве $I$ определена мера и сумма в разложении (9.2.1) является интегралом по этой мере. □

Замечание 9.2.2. Поскольку в теле $D$ определена операция произведения, то мы можем рассматривать тело как алгебру над полем $F \subset Z(D)$. Для элементов базиса мы положим

$$e_i e_j = e_k C^k_{ij} \tag{9.2.2}$$

Коэффициенты $C^k_{ij}$ разложения (9.2.2) называются **структурными константами** тела $D$ над полем $F$.

Из равенств (9.2.1), (9.2.2) следует

$$ab = e_k C^k_{ij} a^i b^j \tag{9.2.3}$$

Из равенства (9.2.3) следует

$$(ab)c = e_k C^k_{ij} (ab)^i c^j = e_k C^k_{ij} C^i_{mn} a^m b^n c^j \tag{9.2.4}$$

$$a(bc) = e_k C^k_{ij} a^i (bc)^j = e_k C^k_{ij} a^i C^j_{mn} b^m c^n \tag{9.2.5}$$

Из ассоциативности произведения

$$(ab)c = a(bc)$$

и равенств (9.2.4) и (9.2.5) следует

$$e_k C^k_{ij} C^i_{mn} a^m b^n c^j = e_k C^k_{ij} a^i C^j_{mn} b^m c^n \tag{9.2.6}$$

Так как векторы $a$, $b$, $c$ произвольны, а векторы $e_k$ линейно независимы, то из равенства (9.2.6) следует

$$C^k_{jn} C^j_{im} = C^k_{ij} C^j_{mn} \tag{9.2.7}$$

Теорема 9.2.3. *Координаты $a^j$ вектора $a$ являются тензором*

$$a^j = A^j_i a'^i \tag{9.2.8}$$

Доказательство. Пусть $\overline{\overline{e}}'$ - другой базис. Пусть

$$e'_i = e_j A^j_i \tag{9.2.9}$$

преобразование, отображающее базис $\overline{\overline{e}}$ в базис $\overline{\overline{e}}'$. Так как вектор $a$ не меняется, то

$$a = e'_i a'^i = e_j a^j \tag{9.2.10}$$

Из равенств (9.2.9) и (9.2.10) следует

$$e_j a^j = e'_i a'^i = e_j A^j_i a'^i \tag{9.2.11}$$

Так как векторы $e_j$ линейно независимы, то равенство (9.2.8) следует из равенства (9.2.11). Следовательно, компоненты вектора являются тензором. □

Теорема 9.2.4. *Структурные константы тела $D$ над полем $F$ являются тензором*

$$A^l_k C'^k_{ij} A^{-1\cdot i}_{\ \ \ n} A^{-1\cdot j}_{\ \ \ m} = C^l_{nm} \tag{9.2.12}$$

Доказательство. Рассмотрим аналогичным образом преобразование произведения. Равенство (9.2.3) в базисе $\overline{\overline{e}}'$ имеет вид

$$ab = e'_k C'^k_{ij} a'^i b'^j \tag{9.2.13}$$

Подставив (9.2.8) и (9.2.9) в (9.2.13), получим

$$ab = e_l A^l_k C'^k_{ij} a^n A^{-1\cdot i}_{\ \ \ n} b^m A^{-1\cdot j}_{\ \ \ m} \tag{9.2.14}$$

Из (9.2.3) и (9.2.14) следует

$$e_l A^l_k C'^k_{ij} A^{-1\cdot i}_{\ \ \ n} a^n A^{-1\cdot j}_{\ \ \ m} b^m = e_l C^l_{nm} a^n b^m \tag{9.2.15}$$

Так как векторы $a$ и $b$ произвольны, а векторы $e_l$ линейно независимы, то равенство (9.2.12) следует из равенства (9.2.15). Следовательно, структурные константы являются тензором. □

Глава 10

# Линейное отображение тела

## 10.1. Линейное отображение тела

Согласно замечанию 9.2.2, мы рассматриваем тело $D$ как алгебру над полем $F \subset Z(D)$.

Определение 10.1.1. Пусть $D_1$, $D_2$ - тела. Пусть $F$ - поле такое, что $F \subset Z(D_1)$, $F \subset Z(D_2)$. Линейное отображение

$$f : D_1 \to D_2$$

$F$-векторного пространства $D_1$ в $F$-векторное пространство $D_2$ называется **линейным отображением тела $D_1$ в тело $D_2$**. □

Согласно определению 10.1.1, линейное отображение $f$ тела $D_1$ в тело $D_2$ удовлетворяет свойству

$$\begin{aligned} f(a+b) &= f(a)+f(b) \quad & a,b \in D \\ f(pa) &= pf(a) & p \in F \end{aligned}$$

Теорема 10.1.2. *Пусть отображение*

$$f : D_1 \to D_2$$

*является линейным отображением тела $D_1$ характеристики 0 в тело $D_2$ характеристики 0. Тогда*[10.1]

$$f(nx) = nf(x)$$

*для любого целого $n$.*

Доказательство. Мы докажем теорему индукцией по $n$. При $n = 1$ утверждение очевидно, так как

$$f(1x) = f(x) = 1f(x)$$

Допустим уравнение справедливо при $n = k$. Тогда

$$f((k+1)x) = f(kx+x) = f(kx)+f(x) = kf(x)+f(x) = (k+1)f(x)$$

□

[10.1] Допустим

$$f : D_1 \to D_2$$

является линейным отображением тела $D_1$ характеристики 2 в тело $D_2$ характеристики 3. Тогда для любого $a \in D_1$, $2a = 0$, хотя $2f(a) \neq 0$. Следовательно, если мы предположим, что характеристика тела $D_1$ больше 0, то мы должны потребовать, чтобы характеристика тела $D_1$ равнялась характеристике тела $D_2$.

Теорема 10.1.3. *Пусть отображение*

$$f : D_1 \to D_2$$

*является линейным отображением тела* $D_1$ *характеристики* 0 *в тело* $D_2$ *характеристики* 0. *Тогда*

$$f(ax) = af(x)$$

*для любого рационального* $a$.

Доказательство. Запишем $a$ в виде $a = \frac{p}{q}$. Положим $y = \frac{1}{q}x$. Тогда согласно теореме 10.1.2

$$f(x) = f(qy) = qf(y) = qf\left(\frac{1}{q}x\right) \qquad (10.1.1)$$

Из равенства (10.1.1) следует

$$\frac{1}{q}f(x) = f\left(\frac{1}{q}x\right) \qquad (10.1.2)$$

Из равенства (10.1.2) следует

$$f\left(\frac{p}{q}x\right) = pf\left(\frac{1}{q}x\right) = \frac{p}{q}f(x)$$

□

Мы не можем распространить утверждение теоремы 10.1.3 на произвольное подполе центра $Z(D)$ тела $D$.

Теорема 10.1.4. *Пусть тело* $D$ *является алгеброй над полем* $F \subset Z(D)$. *Если* $F \neq Z(D)$, *то существует линейное отображение*

$$f : D \to D$$

*которое нелинейно над полем* $Z(D)$.

Доказательство. Для доказательства теоремы достаточно рассмотреть поле комплексных чисел $C$ так как $C = Z(C)$. Так как поле комплексных чисел является алгеброй над полем вещественных чисел, то функция

$$z \to \overline{z}$$

линейна. Однако равенство

$$\overline{az} = a\overline{z}$$

неверно. □

Из теоремы 10.1.4 возникает вопрос. Для чего мы рассматриваем линейные отображения над полем $F \neq Z(D)$, если это приводит к резкому расширению множества линейных отображений? Ответом на этот вопрос служит богатый опыт теории функции комплексного переменного.

Теорема 10.1.5. *Пусть отображения*

$$f : D_1 \to D_2$$

$$g : D_1 \to D_2$$

*являются линейными отображениями тела* $D_1$ *в тело* $D_2$. *Тогда отображение* $f + g$ *также является линейным.*

Доказательство. Утверждение теоремы следует из цепочки равенств

$$\begin{aligned}(f+g)(x+y)=&f(x+y)+g(x+y)=f(x)+f(y)+g(x)+g(y)\\=&(f+g)(x)+(f+g)(y)\\(f+g)(px)=&f(px)+g(px)=pf(x)+pg(x)=p(f(x)+g(x))\\=&p(f+g)(x)\end{aligned}$$

□

Теорема 10.1.6. *Пусть отображение*

$$f:D_1\to D_2$$

*является линейным отображением тела* $D_1$ *в тело* $D_2$*. Тогда отображения* $af$, $fb$, $a,\ b\in R_2$, *также являются линейными.*

Доказательство. Утверждение теоремы следует из цепочки равенств

$$\begin{aligned}(af)(x+y)=&a(f(x+y))=a(f(x)+f(y))=af(x)+af(y)\\=&(af)(x)+(af)(y)\\(af)(px)=&a(f(px))=a(pf(x))=p(af(x))\\=&p(af)(x)\\(fb)(x+y)=&(f(x+y))b=(f(x)+f(y))b=f(x)b+f(y)b\\=&(fb)(x)+(fb)(y)\\(fb)(px)=&(f(px))b=(pf(x))b=p(f(x)b)\\=&p(fb)(x)\end{aligned}$$

□

Определение 10.1.7. Обозначим $\mathcal{L}(D_1;D_2)$ множество линейных отображений

$$f:D_1\to D_2$$

тела $D_1$ в тело $D_2$. □

Теорема 10.1.8. *Мы можем представить линейное отображение*

$$f:D_1\to D_2$$

*тела* $D_1$ *в тело* $D_2$ *в виде*

$$f(x)=f_{k\cdot s_k\cdot 0}\ G_k(x)\ f_{k\cdot s_k\cdot 1}\tag{10.1.3}$$

*где* $(G_k,k\in K)$ *- множество аддитивных отображений тела* $D_1$ *в тело* $D_2$*.*[10.2] *Выражение* $f_{k\cdot s_k\cdot p}$, $p=0,\ 1,$ *в равенстве* (10.1.3) *называется* **компонентой линейного отображения** $f$.

Доказательство. Утверждение теоремы следует из теорем 10.1.5 и 10.1.6. □

[10.2] Здесь и в дальнейшем мы будем предполагать сумму по индексу, который встречается в произведении несколько раз. Равенство (10.1.3) является рекурсивным определением и есть надежда, что мы можем его упростить.

Если в теореме 10.1.8 $|K| = 1$, то равенство (10.1.3) имеет вид

$$f(x) = f_{s\cdot 0}\ G(x)\ f_{s\cdot 1} \tag{10.1.4}$$

и отображение $f$ называется **линейным отображением, порождённым отображением** $G$. Отображение $G$ мы будем называть **образующей линейного отображения**.

Теорема 10.1.9. *Пусть $D_1$, $D_2$ - тела характеристики* 0. *Пусть $F$, $F \subset Z(D_1)$, $F \subset Z(D_2)$, - поле. Пусть $G$ - линейное отображение. Пусть $\overline{\overline{e}}$ - базис тела $D_2$ над полем $F$.* **Стандартное представление линейного отображения** (10.1.4) *имеет вид*[10.3]

$$f(x) = f_G^{ij}\ e_i\ G(x)\ e_j \tag{10.1.5}$$

*Выражение $f_G^{ij}$ в равенстве* (10.1.5) *называется* **стандартной компонентой линейного отображения** $f$.

Доказательство. Компоненты линейного отображения $f$ имеют разложение

$$f_{s\cdot p} = f^{i}_{s\cdot p}\ e_i \tag{10.1.6}$$

относительно базиса $\overline{\overline{e}}$. Если мы подставим (10.1.6) в (10.1.4), мы получим

$$f(x) = f^{i}_{s\cdot 0}\ e_i\ G(x)\ f^{j}_{s\cdot 1}\ e_j \tag{10.1.7}$$

Подставив в равенство (10.1.7) выражение

$$f_G^{ij} = f^{i}_{s\cdot 0}\ f^{j}_{s\cdot 1}$$

мы получим равенство (10.1.5). □

Теорема 10.1.10. *Пусть $D_1$, $D_2$ - тела характеристики* 0. *Пусть $F$, $F \subset Z(D_1)$, $F \subset Z(D_2)$, - поле. Пусть $G$ - линейное отображение. Пусть $\overline{\overline{e}}_1$ - базис тела $D_1$ над полем $F$. Пусть $\overline{\overline{e}}_2$ - базис тела $D_2$ над полем $F$. Пусть $C_{2\cdot}{}^{p}_{kl}$ - структурные константы тела $D_2$. Тогда линейное отображение* (10.1.4), *порождённое линейным отображением $G$, можно записать в виде*

$$\begin{array}{ll} f(a) = e_{2\cdot j}\ f^j_i a^i & f^j_k \in F \\ a = e_{1\cdot i}\ a^i & a^i \in F \quad a \in D_1 \end{array} \tag{10.1.8}$$

$$f^j_i = G^l_i f^{kr}_G C_{2\cdot}{}^{p}_{kl} C_{2\cdot}{}^{j}_{pr} \tag{10.1.9}$$

Доказательство. Выберем отображение

$$\begin{array}{c} G : D_1 \to D_2 \quad a = e_{1\cdot i} a^i \to G(a) = e_{2\cdot j} G^j_i a^i \\ a^i \in F \quad G^j_i \in F \end{array} \tag{10.1.10}$$

Согласно теореме 5.4.3 линейное отображение $f(a)$ относительно базисов $\overline{\overline{e}}_1$ и $\overline{\overline{e}}_2$ принимает вид (10.1.8). Из равенств (10.1.5) и (10.1.10) следует

$$f(a) = G^l_i a^i f^{kj}_G e_{2\cdot k} e_{2\cdot l} e_{2\cdot j} \tag{10.1.11}$$

[10.3] Представление линейного отображения тела с помощью компонент линейного отображения неоднозначно. Чисто алгебраическими методами мы можем увеличить либо уменьшить число слагаемых. Если размерность тела $D_2$ над полем $F$ конечна, то стандартное представление линейного отображения гарантирует конечность множества слагаемых в представлении отображения.

Из равенств (10.1.8) и (10.1.11) следует

$$e_{2\cdot j}\ f^{j}_{i}\ a^{i}=G^{l}_{i}a^{i}f^{kr}_{G}e_{2\cdot k}e_{2\cdot l}e_{2\cdot r}=G^{l}_{i}a^{i}f^{kr}_{G}C_{2\cdot}{}^{p}_{kl}C_{2\cdot}{}^{j}_{pr}\ e_{2\cdot j}\tag{10.1.12}$$

Так как векторы $e_{2\cdot r}$ линейно независимы над полем $F$ и величины $a^{k}$ произвольны, то из равенства (10.1.12) следует равенство (10.1.9). □

При рассмотрении отображений

$$f:D\to D\tag{10.1.13}$$

мы будем полагать $G(x)=x$.

Теорема 10.1.11. *Пусть $D$ является телом характеристики 0. Линейное отображение (10.1.13) имеет вид*

$$f(x)=f_{s\cdot 0}\ x\ f_{s\cdot 1}\tag{10.1.14}$$

Теорема 10.1.12. *Пусть $D$ - тело характеристики 0. Пусть $\overline{\overline{e}}$ - базис тела $D$ над полем $F\subset Z(D)$. Стандартное представление линейного отображения (10.1.14) тела имеет вид*

$$f(x)=f^{ij}\ e_{i}xe_{j}\tag{10.1.15}$$

Теорема 10.1.13. *Пусть $D$ является телом характеристики 0. Пусть $\overline{\overline{e}}$ - базис тела $D$ над полем $F\subset Z(D)$. Тогда линейное отображение (10.1.13) можно записать в виде*

$$\begin{aligned}f(a)&=e_{j}f^{j}_{i}a^{i} && f^{j}_{k}\in F\\ a&=e_{i}a^{i} && a^{i}\in F\quad a\in D\end{aligned}\tag{10.1.16}$$

$$f^{j}_{i}=f^{kr}C^{p}_{ki}C^{j}_{pr}\tag{10.1.17}$$

Теорема 10.1.14. *Рассмотрим матрицу*

$$\mathcal{C}=\left(\mathcal{C}^{\cdot j}_{i\cdot kr}\right)=\left(C^{p}_{ki}C^{j}_{pr}\right)\tag{10.1.18}$$

*строки которой пронумерованы индексом ${}^{\cdot j}_{i}$ и столбцы пронумерованы индексом ${}_{\cdot kr}$ Если $\det\mathcal{C}\neq 0$, то для заданных координат линейного преобразования $f^{j}_{i}$ система линейных уравнений (10.1.17) относительно стандартных компонент этого преобразования $f^{kr}$ имеет единственное решение. Если $\det\mathcal{C}=0$, то условием существования решения системы линейных уравнений (10.1.17) является равенство*

$$\operatorname{rank}\left(\mathcal{C}^{\cdot j}_{i\cdot kr}\quad f^{j}_{i}\right)=\operatorname{rank}\mathcal{C}\tag{10.1.19}$$

*В этом случае система линейных уравнений (10.1.17) имеет бесконечно много решений и существует линейная зависимость между величинами $f^{j}_{i}$.*

Доказательство. Равенство (10.1.14) является частным случаем равенства (10.1.4) при условии $G(x)=x$. Теорема 10.1.12 является частным случаем теоремы 10.1.9 при условии $G(x)=x$. Теорема 10.1.13 является частным случаем теоремы 10.1.10 при условии $G(x)=x$. Утверждение теоремы 10.1.14 является следствием теории линейных уравнений над полем. □

Теорема 10.1.15. *Стандартные компоненты тождественного отображения имеют вид*

$$f^{kr}=\delta^{k}_{0}\delta^{r}_{0}\tag{10.1.20}$$

Доказательство. Равенство (10.1.20) является следствием равенства

$$x = e_0\ x\ e_0$$

Убедимся, что стандартные компоненты (10.1.20) линейного преобразования удовлетворяют уравнению

$$\delta_i^j = f^{kr} C_{ki}^p C_{pr}^j \tag{10.1.21}$$

которое следует из уравнения (10.1.17) если $f = \delta$. Из равенств (10.1.20), (10.1.21) следует

$$\delta_i^j = C_{0i}^p C_{p0}^j \tag{10.1.22}$$

Равенство (10.1.22) верно, так как из равенств

$$e_j e_0 = e_0\ e_j = e_j$$

следует

$$C_{0r}^j = \delta_r^j \quad C_{r0}^j = \delta_r^j$$

Если $\det \mathcal{C} \neq 0$, то решение (10.1.20) единственно. Если $\det \mathcal{C} = 0$, то система линейных уравнений (10.1.21) имеет бесконечно много решений. Однако нас интересует по крайней мере одно. □

Теорема 10.1.16. *Если* $\det \mathcal{C} \neq 0$, *то стандартные компоненты нулевого отображения*

$$z : A \to A \quad z(x) = 0$$

*определены однозначно и имеют вид* $z^{ij} = 0$. *Если* $\det \mathcal{C} = 0$, *то множество стандартных компонент нулевого отображения порождает векторное пространство.*

Доказательство. Теорема верна, поскольку стандартные компоненты $z^{ij}$ являются решением однородной системы линейных уравнений

$$0 = z^{kr} C_{ki}^p C_{pr}^j$$

□

Замечание 10.1.17. Рассмотрим равенство

$$a^{kr}\ e_k\ x\ e_r = b^{kr}\ e_k\ x\ e_r \tag{10.1.23}$$

Из теоремы 10.1.16 следует, что только при условии $\det \mathcal{C} \neq 0$ из равенства (10.1.23) следует

$$a^{kr} = b^{kr} \tag{10.1.24}$$

В противном случае мы должны предполагать равенство

$$a^{kr} = b^{kr} + z^{kr} \tag{10.1.25}$$

Несмотря на это, мы и в случае $\det \mathcal{C} = 0$ будем пользоваться стандартным представлением так как в общем случае указать множество линейно независимых векторов - задача достаточно сложная. Если мы хотим определить операцию над линейными отображениями, записанными в стандартном представлении, то также как в случае теоремы 10.1.15 мы будем выбирать представитель из множества возможных представлений. □

Теорема 10.1.18. *Выражение*

$$f_k^r = f^{ij} C_{ik}^p C_{pj}^r$$

*является тензором над полем F*

$$f'^j_i = A^k_i f^l_k A^{-1 \cdot j}_l \tag{10.1.26}$$

Доказательство. $D$-линейное отображение относительно базиса $\overline{\overline{e}}$ имеет вид (10.1.16). Пусть $\overline{\overline{e}}'$ - другой базис. Пусть

$$e'_i = e_j A^j_i \tag{10.1.27}$$

преобразование, отображающее базис $\overline{\overline{e}}$ в базис $\overline{\overline{e}}'$. Так как линейное отображение $f$ не меняется, то

$$f(x) = e'_l f'^l_k x'^k \tag{10.1.28}$$

Подставим (9.2.8), (10.1.27) в равенство (10.1.28)

$$f(x) = e_j A^j_l f'^l_k A^{-1 \cdot k}_i x^i \tag{10.1.29}$$

Так как векторы $e_j$ линейно независимы и компоненты вектора $x^i$ произвольны, то равенство (10.1.26) следует из равенства (10.1.29). Следовательно, выражение $f_k^r$ является тензором над полем $F$. □

Определение 10.1.19. Множество

$$\ker f = \{x \in D_1 : f(x) = 0\}$$

называется **ядром линейного отображения**

$$f : D_1 \to D_2$$

тела $D_1$ в тело $D_2$. □

Теорема 10.1.20. *Ядро линейного отображения*

$$f : D_1 \to D_2$$

*является подгруппой аддитивной группы тела $D_1$.*

Доказательство. Пусть $a$, $b \in \ker f$. Тогда

$$f(a) = 0$$
$$f(b) = 0$$
$$f(a+b) = f(a) + f(b) = 0$$

Следовательно, $a + b \in \ker f$. □

Определение 10.1.21. Линейное отображение

$$f : D_1 \to D_2$$

тела $D_1$ в тело $D_2$ называется **вырожденным**, если

$$\ker f \neq \{0\}$$

□

Теорема 10.1.22. *Пусть $D$ - тело характеристики 0. Пусть $\overline{\overline{e}}$ - базис тела $D$ над центром $Z(D)$ тела $D$. Пусть*

$$f: D \to D \quad f(x) = f_{s\cdot 0}\ x\ f_{s\cdot 1} \tag{10.1.30}$$

$$= f^{ij}\ e_i x e_j \tag{10.1.31}$$

$$g: D \to D \quad g(x) = g_{t\cdot 0}\ x\ g_{t\cdot 1} \tag{10.1.32}$$

$$= g^{ij}\ e_i x e_j \tag{10.1.33}$$

*линейные отображения тела $D$. Отображение*

$$h(x) = gf(x) = g(f(x)) \tag{10.1.34}$$

*является линейным отображением*

$$h(x) = h_{ts\cdot 0}\ x\ h_{ts\cdot 1} \tag{10.1.35}$$

$$= h^{pr}\ e_p\ x\ e_r \tag{10.1.36}$$

*где*

$$h_{ts\cdot 0} = g_{t\cdot 0}\ f_{s\cdot 0} \tag{10.1.37}$$

$$h_{ts\cdot 1} = f_{s\cdot 1}\ g_{t\cdot 1} \tag{10.1.38}$$

$$h^{pr} = g^{ij} f^{kl} C^p_{ik} C^r_{lj} \tag{10.1.39}$$

Доказательство. Отображение (10.1.34) линейно так как

$$h(x+y) = g(f(x+y)) = g(f(x)+f(y)) = g(f(x)) + g(f(y)) = h(x) + h(y)$$

$$h(ax) = g(f(ax)) = g(af(x)) = ag(f(x)) = ah(x)$$

Если мы подставим (10.1.30) и (10.1.32) в (10.1.34), то мы получим

$$h(x) = g_{t\cdot 0}\ f(x)\ g_{t\cdot 1} = g_{t\cdot 0}\ f_{s\cdot 0}\ x\ f_{s\cdot 1}\ g_{t\cdot 1} \tag{10.1.40}$$

Сравнивая (10.1.40) и (10.1.35), мы получим (10.1.37), (10.1.38).

Если мы подставим (10.1.31) и (10.1.33) в (10.1.34), то мы получим

$$\begin{aligned} h(x) =& g^{ij}\ e_i\ f(x)\ e_j \\ =& g^{ij}\ e_i\ f^{kl}\ e_k\ x\ e_l e_j \\ =& g^{ij} f^{kl} C^p_{ik} C^r_{lj}\ e_p\ x\ e_r \end{aligned} \tag{10.1.41}$$

Сравнивая (10.1.41) и (10.1.36), мы получим (10.1.39). □

## 10.2. Полилинейное отображение тела

Определение 10.2.1. Пусть $R_1$, ..., $R_n$, $P$ - кольца характеристики 0. Пусть $S$ - модуль над кольцом $P$. Пусть $F$ - коммутативное кольцо, которое для любого $i$ является подкольцом центра кольца $R_i$. Отображение

$$f: R_1 \times ... \times R_n \to S$$

называется **полилинейным над коммутативным кольцом** $F$, если

$$f(p_1, ..., p_i + q_i, ..., p_n) = f(p_1, ..., p_i, ..., p_n) + f(p_1, ..., q_i, ..., p_n)$$

$$f(a_1, ..., ba_i, ..., a_n) = bf(a_1, ..., a_i, ..., a_n)$$

для любого $i$, $1 \le i \le n$, и любых $p_i$, $q_i \in R_i$, $b \in F$. Обозначим $\mathcal{L}(R_1, ..., R_n; S)$ множество полилинейных отображений колец $R_1$, ..., $R_n$ в модуль $S$. □

Теорема 10.2.2. *Пусть $D$ - тело характеристики 0. Полилинейное отображение*

(10.2.1) $$f : D^n \to D, d = f(d_1, ..., d_n)$$

*имеет вид*

(10.2.2) $$d = f^n_{s\cdot 0}\ \sigma_s(d_1)\ f^n_{s\cdot 1}\ ...\ \sigma_s(d_n)\ f^n_{s\cdot n}$$

*$\sigma_s$ - перестановка множества переменных $\{d_1, ..., d_n\}$*

$$\sigma_s = \begin{pmatrix} d_1 & ... & d_n \\ \sigma_s(d_1) & ... & \sigma_s(d_n) \end{pmatrix}$$

Доказательство. Мы докажем утверждение индукцией по $n$.

При $n = 1$ доказываемое утверждение является следствием теоремы 10.1.11. При этом мы можем отождествить[10.4]

$$f^1_{s\cdot p} = f_{s\cdot p} \quad p = 0, 1$$

Допустим, что утверждение теоремы справедливо при $n = k - 1$. Тогда отображение (10.2.1) можно представить в виде

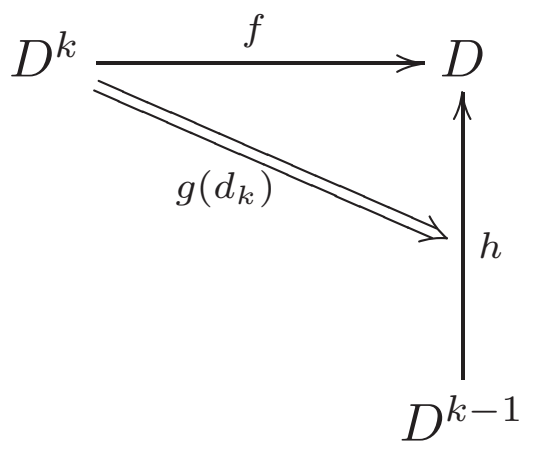

$$d = f(d_1, ..., d_k) = g(d_k)(d_1, ..., d_{k-1})$$

Согласно предположению индукции полилинейное отображение $h$ имеет вид

$$d = h^{k-1}_{t\cdot 0}\ \sigma_t(d_1)\ h^{k-1}_{t\cdot 1}\ ...\ \sigma_t(d_{k-1})\ h^{k-1}_{t\cdot k-1}$$

Согласно построению $h = g(d_k)$. Следовательно, выражения $h_{t\cdot p}$ являются функциями $d_k$. Поскольку $g(d_k)$ - линейная функция $d_k$, то только одно выражение $h_{t\cdot p}$ является линейной функцией переменной $d_k$, и остальные выражения ${}_{t\cdot q}h$ не зависят от $d_k$.

Не нарушая общности, положим $p = 0$. Согласно равенству (10.1.14) для заданного $t$

$$h^{k-1}_{t\cdot 0} = g_{tr\cdot 0}\ d_k\ g_{tr\cdot 1}$$

Положим $s = tr$ и определим перестановку $\sigma_s$ согласно правилу

$$\sigma_s = \sigma_{tr} = \begin{pmatrix} d_k & d_1 & ... & d_{k-1} \\ d_k & \sigma_t(d_1) & ... & \sigma_t(d_{k-1}) \end{pmatrix}$$

[10.4] В представлении (10.2.2) мы будем пользоваться следующими правилами.

- Если область значений какого-либо индекса - это множество, состоящее из одного элемента, мы будем опускать соответствующий индекс.
- Если $n = 1$, то $\sigma_s$ - тождественное преобразование. Это преобразование можно не указывать в выражении.

Положим

$$f^k_{tr\cdot q+1} = h^{k-1}_{t\cdot q} \quad q = 1, ..., k-1$$
$$f^k_{tr\cdot q} = g_{tr\cdot q} \quad q = 0, 1$$

Мы доказали шаг индукции. □

Определение 10.2.3. Выражение $f^n_{s\cdot p}$ в равенстве (10.2.2) называется **компонентой полилинейного отображения** $f$. □

Теорема 10.2.4. *Пусть $D$ - тело характеристики 0. Допустим $\overline{\overline{e}}$ - базис тела $D$ над полем $F \subset Z(D)$.* **Стандартное представление полилинейного отображения** *тела имеет вид*

(10.2.3) $$f(d_1, ..., d_n) = f^{i_0...i_n}_t \ e_{i_0} \ \sigma_t(d_1) \ e_{i_1} \ ... \ \sigma_t(d_n) \ e_{i_n}$$

*Индекс $t$ нумерует всевозможные перестановки $\sigma_t$ множества переменных $\{d_1, ..., d_n\}$. Выражение $f^{i_0...i_n}_t$ в равенстве (10.2.3) называется* **стандартной компонентой полилинейного отображения** $f$.

Доказательство. Компоненты полилинейного отображения $f$ имеют разложение

(10.2.4) $$f^n_{s\cdot p} = e_i f^{ni}_{s\cdot p}$$

относительно базиса $\overline{\overline{e}}$. Если мы подставим (10.2.4) в (10.2.2), мы получим

(10.2.5) $$d = f^{nj_1}_{s\cdot 0} \ e_{j_1} \ \sigma_s(d_1) \ f^{nj_2}_{s\cdot 1} \ e_{j_2} \ ... \ \sigma_s(d_n) \ f^{nj_n}_{s\cdot n} \ e_{j_n}$$

Рассмотрим выражение

(10.2.6) $$f^{j_0...j_n}_t = f^{nj_1}_{s\cdot 0} \ ... f^{nj_n}_{s\cdot n}$$

В правой части подразумевается сумма тех слагаемых с индексом $s$, для которых перестановка $\sigma_s$ совпадает. Каждая такая сумма будет иметь уникальный индекс $t$. Подставив в равенство (10.2.5) выражение (10.2.6) мы получим равенство (10.2.3). □

Теорема 10.2.5. *Пусть $\overline{\overline{e}}$ - базис тела $D$ над полем $F \subset Z(D)$. Полилинейное отображение (10.2.1) можно представить в виде $D$-значной формы степени $n$ над полем $F \subset Z(D)$*[10.5]

(10.2.7) $$f(a_1, ..., a_n) = a^{i_1}_1 ... a^{i_n}_n f_{i_1...i_n}$$

*где*

$$a_j = e_i a^i_j$$

(10.2.8) $$f_{i_1...i_n} = f(e_{i_1}, ..., e_{i_n})$$

*и величины $f_{i_1...i_n}$ являются координатами $D$-значного ковариантнго тензора над полем $F$.*

Доказательство. Согласно определению 10.2.1 равенство (10.2.7) следует из цепочки равенств

$$f(a_1, ..., a_n) = f(e_{i_1} a^{i_1}_1, ..., e_{i_n} a^{i_n}_n) = a^{i_1}_1 ... a^{i_n}_n f(e_{i_1}, ..., e_{i_n})$$

Пусть $\overline{\overline{e}}'$ - другой базис. Пусть

(10.2.9) $$e'_i = e_j A^j_i$$

[10.5] Теорема доказана по аналогии с теоремой в [3], с. 107, 108

преобразование, отображающее базис $\overline{\overline{e}}$ в базис $\overline{\overline{e}}'$. Из равенств (10.2.9) и (10.2.8) следует

$$
\begin{aligned}
f'_{i_1...i_n} &= f(e'_{i_1}, ..., e'_{i_n})\\
&= f(e_{j_1}A^{j_1}_{i_1}, ..., e'_{j_n}A^{j_n}_{i_n})\\
&= A^{j_1}_{i_1}...A^{j_n}_{i_n}f(e_{j_1}, ..., e_{j_n})\\
&= A^{j_1}_{i_1}...A^{j_n}_{i_n}f_{j_1...j_n}
\end{aligned} \tag{10.2.10}
$$

Из равенства (10.2.10) следует тензорный закон преобразования координат полилинейного отображения. Из равенства (10.2.10) и теоремы 9.2.3 следует, что значение отображения $f(a_1, ..., a_n)$ не зависит от выбора базиса. □

Полилинейное отображение (10.2.1) **симметрично**, если

$$f(d_1, ..., d_n) = f(\sigma(d_1), ..., \sigma(d_n))$$

для любой перестановки $\sigma$ множества $\{d_1, ..., d_n\}$.

Теорема 10.2.6. *Если полилинейное отображение $f$ симметрично, то*

$$f_{i_1,...,i_n} = f_{\sigma(i_1),...,\sigma(i_n)} \tag{10.2.11}$$

Доказательство. Равенство (10.2.11) следует из равенства

$$
\begin{aligned}
a_1^{i_1}...a_n^{i_n}f_{i_1...i_n} =&f(a_1, ..., a_n)\\
=&f(\sigma(a_1), ..., \sigma(a_n))\\
=&a_1^{i_1}...a_n^{i_n}f_{\sigma(i_1)...\sigma(i_n)}
\end{aligned}
$$

□

Полилинейное отображение (10.2.1) **косо симметрично**, если

$$f(d_1, ..., d_n) = |\sigma|f(\sigma(d_1), ..., \sigma(d_n))$$

для любой перестановки $\sigma$ множества $\{d_1, ..., d_n\}$. Здесь

$$|\sigma| = \begin{cases} 1 & \text{перестановка } \sigma \text{ чётная}\\ -1 & \text{перестановка } \sigma \text{ нечётная}\end{cases}$$

Теорема 10.2.7. *Если полилинейное отображение $f$ косо симметрично, то*

$$f_{i_1,...,i_n} = |\sigma|f_{\sigma(i_1),...,\sigma(i_n)} \tag{10.2.12}$$

Доказательство. Равенство (10.2.12) следует из равенства

$$
\begin{aligned}
a_1^{i_1}...a_n^{i_n}f_{i_1...i_n} =&f(a_1, ..., a_n)\\
=&|\sigma|f(\sigma(a_1), ..., \sigma(a_n))\\
=&a_1^{i_1}...a_n^{i_n}|\sigma|f_{\sigma(i_1)...\sigma(i_n)}
\end{aligned}
$$

□

Теорема 10.2.8. *Координаты полилинейного над полем $F$ отображения (10.2.1) и его компоненты относительно базиса $\overline{\overline{e}}$ удовлетворяют равенству*

$$f_{j_1...j_n} = f_t^{i_0...i_n} C^{k_1}_{i_0\sigma_t(j_1)} C^{l_1}_{k_1 i_1} ... B^{k_n}_{l_{n-1}\sigma_t(j_n)} C^{l_n}_{k_n i_n} e_{l_n} \tag{10.2.13}$$

$$f^p_{j_1...j_n} = f_t^{i_0...i_n} C^{k_1}_{i_0\sigma_t(j_1)} C^{l_1}_{k_1 i_1} ... C^{k_n}_{l_{n-1}\sigma_t(j_n)} C^{p}_{k_n i_n} \tag{10.2.14}$$

Доказательство. В равенстве (10.2.3) положим

$$d_i = e_{j_i} d_i^{j_i}$$

Тогда равенство (10.2.3) примет вид

$$\begin{aligned} f(d_1, ..., d_n) &= f_t^{i_0...i_n}\ e_{i_0}\sigma_t(d_1^{j_1} e_{j_1}) e_{i_1} ... \sigma_t(d_n^{j_n} e_{j_n}) e_{i_n} \\ &= d_1^{j_1} ... d_n^{j_n} f_t^{i_0...i_n} e_{i_0}\sigma_t(e_{j_1}) e_{i_1} ... \sigma_t(e_{j_n}) e_{i_n} \\ &= d_1^{j_1} ... d_n^{j_n} f_t^{i_0...i_n} C^{k_1}_{i_0\sigma_t(j_1)} C^{l_1}_{k_1 i_1} \\ &\quad ... C^{k_n}_{l_{n-1}\sigma_t(j_n)} C^{l_n}_{k_n i_n} e_{l_n} \end{aligned} \tag{10.2.15}$$

Из равенства (10.2.7) следует

$$f(a_1, ..., a_n) = e_p f^p_{i_1...i_n} a_1^{i_1} ... a_n^{i_n} \tag{10.2.16}$$

Равенство (10.2.13) следует из сравнения равенств (10.2.15) и (10.2.7). Равенство (10.2.14) следует из сравнения равенств (10.2.15) и (10.2.16). □

Глава 11

# Алгебра кватернионов

## 11.1. Линейная функция комплексного поля

Теорема 11.1.1 (Уравнения Коши-Римана). *Рассмотрим поле комплексных чисел $C$ как двумерную алгебру над полем действительных чисел. Положим*

$$e_{C\cdot 0} = 1 \quad e_{C\cdot 1} = i \tag{11.1.1}$$

*базис алгебры $C$. Тогда в этом базисе произведение имеет вид*

$$e_{C\cdot 1}^2 = -e_{C\cdot 0} \tag{11.1.2}$$

*и структурные константы имеют вид*

$$\begin{aligned} C_{C\cdot}{}^{0}_{00} &= 1 & C_{C\cdot}{}^{1}_{01} &= 1 \\ C_{C\cdot}{}^{1}_{10} &= 1 & C_{C\cdot}{}^{0}_{11} &= -1 \end{aligned} \tag{11.1.3}$$

*Матрица линейной функции*

$$y^i = x^j f^i_j$$

*поля комплексных чисел над полем действительных чисел удовлетворяет соотношению*

$$f^0_0 = f^1_1 \tag{11.1.4}$$

$$f^1_0 = -f^0_1 \tag{11.1.5}$$

Доказательство. Равенства (11.1.2) и (11.1.3) следуют из равенства $i^2 = -1$. Пользуясь равенством (10.1.17) получаем соотношения

$$f^0_0 = f^{kr} C_{C\cdot}{}^{p}_{k0} C_{C\cdot}{}^{0}_{pr} = f^{0r} C_{C\cdot}{}^{0}_{00} C_{C\cdot}{}^{0}_{0r} + f^{1r} C_{C\cdot}{}^{1}_{10} C_{C\cdot}{}^{0}_{1r} = f^{00} - f^{11} \tag{11.1.6}$$

$$f^1_0 = f^{kr} C_{C\cdot}{}^{p}_{k0} C_{C\cdot}{}^{1}_{pr} = f^{0r} C_{C\cdot}{}^{0}_{00} C_{C\cdot}{}^{1}_{0r} + f^{1r} C_{C\cdot}{}^{1}_{10} C_{C\cdot}{}^{1}_{1r} = f^{01} + f^{10} \tag{11.1.7}$$

$$f^0_1 = f^{kr} C_{C\cdot}{}^{p}_{k1} C_{C\cdot}{}^{0}_{pr} = f^{0r} C_{C\cdot}{}^{1}_{01} C_{C\cdot}{}^{0}_{1r} + f^{1r} C_{C\cdot}{}^{0}_{11} C_{C\cdot}{}^{0}_{0r} = -f^{01} - f^{10} \tag{11.1.8}$$

$$f^1_1 = f^{kr} C_{C\cdot}{}^{p}_{k1} C_{C\cdot}{}^{1}_{pr} = f^{0r} C_{C\cdot}{}^{1}_{01} C_{C\cdot}{}^{1}_{1r} + f^{1r} C_{C\cdot}{}^{0}_{11} C_{C\cdot}{}^{1}_{0r} = f^{00} - f^{11} \tag{11.1.9}$$

Из равенств (11.1.6) и (11.1.9) следует (11.1.4). Из равенств (11.1.7) и (11.1.8) следует (11.1.5). □

## 11.2. Алгебра кватернионов

В этой статье я рассматриваю множество кватернионных алгебр, определённых в [14].

Определение 11.2.1. Пусть $F$ - поле. Расширение $F(i,j,k)$ поля $F$ называется **алгеброй** $E(F,a,b)$ **кватернионов над полем** $F$[11.1], если произведение в алгебре $E$ определено согласно правилам

(11.2.1)

| | $i$ | $j$ | $k$ |
|---|---|---|---|
| $i$ | $a$ | $k$ | $aj$ |
| $j$ | $-k$ | $b$ | $-bi$ |
| $k$ | $-aj$ | $bi$ | $-ab$ |

где $a$, $b\in F$, $ab\neq 0$. □

Элементы алгебры $E(F,a,b)$ имеют вид

$$x=x^0+x^1i+x^2j+x^3k$$

где $x^i\in F$, $i=0,1,2,3$. Кватернион

$$\overline{x}=x^0-x^1i-x^2j-x^3k$$

называется сопряжённым кватерниону $x$. Мы определим **норму кватерниона** $x$ равенством

$$|x|^2=x\overline{x}=(x^0)^2-a(x^1)^2-b(x^2)^2+ab(x^3)^2 \tag{11.2.2}$$

Из равенства (11.2.2) следует, что $E(F,a,b)$ является алгеброй с делением только когда $a<0$, $b<0$. Тогда мы можем пронормировать базис так, что $a=-1$, $b=-1$.

Мы будем обозначать символом $E(F)$ алгебру $E(F,-1,-1)$ кватернионов с делением над полем $F$. Мы будем полагать $E(R,-1,-1)$. Произведение в алгебре кватернионов $E(F)$ определено согласно правилам

(11.2.3)

| | $i$ | $j$ | $k$ |
|---|---|---|---|
| $i$ | $-1$ | $k$ | $-j$ |
| $j$ | $-k$ | $-1$ | $i$ |
| $k$ | $j$ | $-i$ | $-1$ |

В алгебре $E(F)$ норма кватерниона имеет вид

$$|x|^2=x\overline{x}=(x^0)^2+(x^1)^2+(x^2)^2+(x^3)^2 \tag{11.2.4}$$

При этом обратный элемент имеет вид

$$x^{-1}=|x|^{-2}\overline{x} \tag{11.2.5}$$

Внутренний автоморфизм алгебры кватернионов $H$[11.2]

$$\begin{gathered}p\rightarrow qpq^{-1}\\ q(ix+jy+kz)q^{-1}=ix'+jy'+kz'\end{gathered} \tag{11.2.6}$$

[11.1] Я буду следовать определению из [14].

[11.2] См. [15], с. 643.

описывает вращение вектора с координатами $x$, $y$, $z$. Если $q$ записан в виде суммы скаляра и вектора

$$q=\cos\alpha+(ia+jb+kc)\sin\alpha \quad a^2+b^2+c^2=1$$

то (11.2.6) описывает вращение вектора $(x,y,z)$ вокруг вектора $(a,b,c)$ на угол $2\alpha$.

## 11.3. Линейная функция алгебры кватернионов

Теорема 11.3.1. *Положим*

$$e_0=1 \quad e_1=i \quad e_2=j \quad e_3=k \tag{11.3.1}$$

*базис алгебры кватернионов $H$. Тогда в базисе* (11.3.1) *структурные константы имеют вид*

$$\begin{array}{llll}
C^0_{00}=1 & C^1_{01}=\ \ 1 & C^2_{02}=\ \ 1 & C^3_{03}=\ \ 1\\
C^1_{10}=1 & C^0_{11}=-1 & C^3_{12}=\ \ 1 & C^2_{13}=-1\\
C^2_{20}=1 & C^3_{21}=-1 & C^0_{22}=-1 & C^1_{23}=\ \ 1\\
C^3_{30}=1 & C^2_{31}=\ \ 1 & C^1_{32}=-1 & C^0_{33}=-1
\end{array}$$

Доказательство. Значение структурных констант следует из таблицы умножения (11.2.3). □

Так как вычисления в этом разделе занимают много места, я собрал в одном месте ссылки на теоремы в этом разделе.

**Теорема 11.3.2, страница 129:** определение координат линейного отображения алгебры кватернионов $H$ через стандартные компоненты этого отображения.

**Равенство** (11.3.22)**, страница 132:** матричная форма зависимости координат линейного отображения алгебры кватернионов $H$ от стандартных компонент этого отображения.

**Равенство** (11.3.23)**, страница 133:** матричная форма зависимости стандартных компонент линейного отображения алгебры кватернионов $H$ от координат этого отображения.

**Теорема 11.3.4, страница 135:** зависимость стандартных компонент линейного отображения алгебры кватернионов $H$ от координат этого отображения.

Теорема 11.3.2. *Стандартные компоненты линейной функции алгебры кватернионов $H$ относительно базиса* (11.3.1) *и координаты соответствующего линейного преобразования удовлетворяют соотношениям*

$$\left\{\begin{aligned}
f^0_0&=f^{00}-f^{11}-f^{22}-f^{33}\\
f^1_1&=f^{00}-f^{11}+f^{22}+f^{33}\\
f^2_2&=f^{00}+f^{11}-f^{22}+f^{33}\\
f^3_3&=f^{00}+f^{11}+f^{22}-f^{33}
\end{aligned}\right. \tag{11.3.2}$$

$$\begin{cases} f_0^1 = & f^{01}+f^{10}+f^{23}-f^{32} \\ f_1^0 = & -f^{01}-f^{10}+f^{23}-f^{32} \\ f_2^3 = & -f^{01}+f^{10}-f^{23}-f^{32} \\ f_3^2 = & f^{01}-f^{10}-f^{23}-f^{32} \end{cases} \tag{11.3.3}$$

$$\begin{cases} f_0^2 = & f^{02}-f^{13}+f^{20}+f^{31} \\ f_1^3 = & f^{02}-f^{13}-f^{20}-f^{31} \\ f_2^0 = & -f^{02}-f^{13}-f^{20}+f^{31} \\ f_3^1 = & -f^{02}-f^{13}+f^{20}-f^{31} \end{cases} \tag{11.3.4}$$

$$\begin{cases} f_0^3 = & f^{03}+f^{12}-f^{21}+f^{30} \\ f_1^2 = & -f^{03}-f^{12}-f^{21}+f^{30} \\ f_2^1 = & f^{03}-f^{12}-f^{21}-f^{30} \\ f_3^0 = & -f^{03}+f^{12}-f^{21}-f^{30} \end{cases} \tag{11.3.5}$$

Доказательство. Пользуясь равенством (10.1.17) получаем соотношения

$$\begin{aligned} f_0^0 &= f^{kr}C_{k0}^pC_{pr}^0 \\ &= f^{00}C_{00}^0C_{00}^0+f^{11}C_{10}^1C_{11}^0+f^{22}C_{20}^2C_{22}^0+f^{33}C_{30}^3C_{33}^0 \\ &= f^{00}-f^{11}-f^{22}-f^{33} \end{aligned} \tag{11.3.6}$$

$$\begin{aligned} f_0^1 &= f^{kr}C_{k0}^pC_{pr}^1 \\ &= f^{01}C_{00}^0C_{01}^1+f^{10}C_{10}^1C_{10}^1+f^{23}C_{20}^2C_{23}^1+f^{32}C_{30}^3C_{32}^1 \\ &= f^{01}+f^{10}+f^{23}-f^{32} \end{aligned} \tag{11.3.7}$$

$$\begin{aligned} f_0^2 &= f^{kr}C_{k0}^pC_{pr}^2 \\ &= f^{02}C_{00}^0C_{02}^2+f^{13}C_{10}^1C_{13}^2+f^{20}C_{20}^2C_{20}^2+f^{31}C_{30}^3C_{31}^2 \\ &= f^{02}-f^{13}+f^{20}+f^{31} \end{aligned} \tag{11.3.8}$$

$$\begin{aligned} f_0^3 &= f^{kr}C_{k0}^pC_{pr}^3 \\ &= f^{03}C_{00}^0C_{03}^3+f^{12}C_{10}^1C_{12}^3+f^{21}C_{20}^2C_{21}^3+f^{30}C_{30}^3C_{30}^3 \\ &= f^{03}+f^{12}-f^{21}+f^{30} \end{aligned} \tag{11.3.9}$$

$$\begin{aligned}f_1^0 &= f^{kr}C_{k1}^pC_{pr}^0\\ &= f^{01}C_{01}^1C_{11}^0+f^{10}C_{11}^0C_{00}^0+f^{23}C_{21}^3C_{33}^0+f^{32}C_{31}^2C_{22}^0\\ &= -f^{01}-f^{10}+f^{23}-f^{32}\end{aligned} \tag{11.3.10}$$

$$\begin{aligned}f_1^1 &= f^{kr}C_{k1}^pC_{pr}^1\\ &= f^{00}C_{01}^1C_{10}^1+f^{11}C_{11}^0C_{01}^1+f^{22}C_{21}^3C_{32}^1+f^{33}C_{31}^2C_{23}^1\\ &= f^{00}-f^{11}+f^{22}+f^{33}\end{aligned} \tag{11.3.11}$$

$$\begin{aligned}f_1^2 &= f^{kr}C_{k1}^pC_{pr}^2\\ &= f^{03}C_{01}^1C_{13}^2+f^{12}C_{11}^0C_{02}^2+f^{21}C_{21}^3C_{31}^2+f^{30}C_{31}^2C_{20}^2\\ &= -f^{03}-f^{12}-f^{21}+f^{30}\end{aligned} \tag{11.3.12}$$

$$\begin{aligned}f_1^3 &= f^{kr}C_{k1}^pC_{pr}^3\\ &= f^{02}C_{01}^1C_{12}^3+f^{13}C_{11}^0C_{03}^3+f^{20}C_{21}^3C_{30}^3+f^{31}C_{31}^2C_{21}^3\\ &= f^{02}-f^{13}-f^{20}-f^{31}\end{aligned} \tag{11.3.13}$$

$$\begin{aligned}f_2^0 &= f^{kr}C_{k2}^pC_{pr}^0\\ &= f^{02}C_{02}^2C_{22}^0+f^{13}C_{12}^3C_{33}^0+f^{20}C_{22}^0C_{00}^0+f^{31}C_{32}^1C_{11}^0\\ &= -f^{02}-f^{13}-f^{20}+f^{31}\end{aligned} \tag{11.3.14}$$

$$\begin{aligned}f_2^1 &= f^{kr}C_{k2}^pC_{pr}^1\\ &= f^{03}C_{02}^2C_{23}^1+f^{12}C_{12}^3C_{32}^1+f^{21}C_{22}^0C_{01}^1+f^{30}C_{32}^1C_{10}^1\\ &= f^{03}-f^{12}-f^{21}-f^{30}\end{aligned} \tag{11.3.15}$$

$$\begin{aligned}f_2^2 &= f^{kr}C_{k2}^pC_{pr}^2\\ &= f^{00}C_{02}^2C_{20}^2+f^{11}C_{12}^3C_{31}^2+f^{22}C_{22}^0C_{02}^2+f^{33}C_{32}^1C_{13}^2\\ &= f^{00}+f^{11}-f^{22}+f^{33}\end{aligned} \tag{11.3.16}$$

$$\begin{aligned}f_2^3 &= f^{kr}C_{k2}^pC_{pr}^3\\ &= f^{01}C_{02}^2C_{21}^3+f^{10}C_{12}^3C_{30}^3+f^{23}C_{22}^0C_{03}^3+f^{32}C_{32}^1C_{12}^3\\ &= -f^{01}+f^{10}-f^{23}-f^{32}\end{aligned} \tag{11.3.17}$$

$$\begin{aligned}f_3^0 &= f^{kr}C_{k3}^pC_{pr}^0\\ &= f^{03}C_{03}^3C_{33}^0+f^{12}C_{13}^2C_{22}^0+f^{21}C_{23}^1C_{11}^0+f^{30}C_{33}^0C_{00}^0\\ &= -f^{03}+f^{12}-f^{21}-f^{30}\end{aligned}\tag{11.3.18}$$

$$\begin{aligned}f_3^1 &= f^{kr}C_{k3}^pC_{pr}^1\\ &= f^{02}C_{03}^3C_{32}^1+f^{13}C_{13}^2C_{23}^1+f^{20}C_{23}^1C_{10}^1+f^{31}C_{33}^0C_{01}^1\\ &= -f^{02}-f^{13}+f^{20}-f^{31}\end{aligned}\tag{11.3.19}$$

$$\begin{aligned}f_3^2 &= f^{kr}C_{k3}^pC_{pr}^2\\ &= f^{01}C_{03}^3C_{31}^2+f^{10}C_{13}^2C_{20}^2+f^{23}C_{23}^1C_{13}^2+f^{32}C_{33}^0C_{02}^2\\ &= f^{01}-f^{10}-f^{23}-f^{32}\end{aligned}\tag{11.3.20}$$

$$\begin{aligned}f_3^3 &= f^{kr}C_{k3}^pC_{pr}^3\\ &= f^{00}C_{03}^3C_{30}^3+f^{11}C_{13}^2C_{21}^3+f^{22}C_{23}^1C_{12}^3+f^{33}C_{33}^0C_{03}^3\\ &= f^{00}+f^{11}+f^{22}-f^{33}\end{aligned}\tag{11.3.21}$$

Уравнения (11.3.6), (11.3.11), (11.3.16), (11.3.21) формируют систему линейных уравнений (11.3.2).

Уравнения (11.3.7), (11.3.10), (11.3.17), (11.3.20) формируют систему линейных уравнений (11.3.3).

Уравнения (11.3.8), (11.3.13), (11.3.14), (11.3.19) формируют систему линейных уравнений (11.3.4).

Уравнения (11.3.9), (11.3.12), (11.3.15), (11.3.18) формируют систему линейных уравнений (11.3.5). □

Теорема 11.3.3. *Рассмотрим алгебру кватернионов $H$ с базисом (11.3.1). Стандартные компоненты линейной функции над полем $F$ и координаты этой функции над полем $F$ удовлетворяют соотношениям*

$$\begin{pmatrix}f_0^0 & f_1^0 & f_2^0 & f_3^0\\ f_1^1 & -f_1^0 & f_3^1 & -f_2^1\\ f_2^2 & -f_3^2 & -f_0^2 & f_1^2\\ f_3^3 & f_2^3 & -f_1^3 & -f_0^3\end{pmatrix} = \begin{pmatrix}1 & -1 & -1 & -1\\ 1 & -1 & 1 & 1\\ 1 & 1 & -1 & 1\\ 1 & 1 & 1 & -1\end{pmatrix}\begin{pmatrix}f^{00} & -f^{01} & -f^{02} & -f^{03}\\ f^{11} & f^{10} & f^{13} & -f^{12}\\ f^{22} & -f^{23} & f^{20} & f^{21}\\ f^{33} & f^{32} & -f^{31} & f^{30}\end{pmatrix}\tag{11.3.22}$$

$$\begin{pmatrix}f^{00}&-f^{01}&-f^{02}&-f^{03}\\ f^{11}&f^{10}&f^{13}&-f^{12}\\ f^{22}&-f^{23}&f^{20}&f^{21}\\ f^{33}&f^{32}&-f^{31}&f^{30}\end{pmatrix}$$

$$=\frac{1}{4}\begin{pmatrix}1&1&1&1\\ -1&-1&1&1\\ -1&1&-1&1\\ -1&1&1&-1\end{pmatrix}\begin{pmatrix}f_0^0&f_1^0&f_2^0&f_3^0\\ f_1^1&-f_1^0&f_3^1&-f_2^1\\ f_2^2&-f_3^2&-f_0^2&f_1^2\\ f_3^3&f_2^3&-f_1^3&-f_0^3\end{pmatrix}\tag{11.3.23}$$

*где*

$$\begin{pmatrix}1&-1&-1&-1\\ 1&-1&1&1\\ 1&1&-1&1\\ 1&1&1&-1\end{pmatrix}^{-1}=\frac{1}{4}\begin{pmatrix}1&1&1&1\\ -1&-1&1&1\\ -1&1&-1&1\\ -1&1&1&-1\end{pmatrix}$$

Доказательство. Запишем систему линейных уравнений (11.3.2) в виде произведения матриц

$$\begin{pmatrix}f_0^0\\ f_1^1\\ f_2^2\\ f_3^3\end{pmatrix}=\begin{pmatrix}1&-1&-1&-1\\ 1&-1&1&1\\ 1&1&-1&1\\ 1&1&1&-1\end{pmatrix}\begin{pmatrix}f^{00}\\ f^{11}\\ f^{22}\\ f^{33}\end{pmatrix}\tag{11.3.24}$$

Запишем систему линейных уравнений (11.3.3) в виде произведения матриц

$$\begin{pmatrix}f_0^1\\ f_1^0\\ f_2^3\\ f_3^2\end{pmatrix}=\begin{pmatrix}1&1&1&-1\\ -1&-1&1&-1\\ -1&1&-1&-1\\ 1&-1&-1&-1\end{pmatrix}\begin{pmatrix}f^{01}\\ f^{10}\\ f^{23}\\ f^{32}\end{pmatrix}\tag{11.3.25}$$

Из равенства (11.3.25) следует

$$\begin{pmatrix} -f^{0}_{1} \\ f^{0}_{1} \\ f^{3}_{2} \\ -f^{2}_{3} \end{pmatrix} = \begin{pmatrix} -1 & -1 & -1 & 1 \\ -1 & -1 & 1 & -1 \\ -1 & 1 & -1 & -1 \\ -1 & 1 & 1 & 1 \end{pmatrix} \begin{pmatrix} f^{01} \\ f^{10} \\ f^{23} \\ f^{32} \end{pmatrix}$$

$$= \begin{pmatrix} 1 & -1 & 1 & 1 \\ 1 & -1 & -1 & -1 \\ 1 & 1 & 1 & -1 \\ 1 & 1 & -1 & 1 \end{pmatrix} \begin{pmatrix} -f^{01} \\ f^{10} \\ -f^{23} \\ f^{32} \end{pmatrix}$$

$$\begin{pmatrix} f^{0}_{1} \\ -f^{0}_{1} \\ -f^{2}_{3} \\ f^{3}_{2} \end{pmatrix} = \begin{pmatrix} 1 & -1 & -1 & -1 \\ 1 & -1 & 1 & 1 \\ 1 & 1 & -1 & 1 \\ 1 & 1 & 1 & -1 \end{pmatrix} \begin{pmatrix} -f^{01} \\ f^{10} \\ -f^{23} \\ f^{32} \end{pmatrix} \tag{11.3.26}$$

Запишем систему линейных уравнений (11.3.4) в виде произведения матриц

$$\begin{pmatrix} f^{2}_{0} \\ f^{3}_{1} \\ f^{0}_{2} \\ f^{1}_{3} \end{pmatrix} = \begin{pmatrix} 1 & -1 & 1 & 1 \\ 1 & -1 & -1 & -1 \\ -1 & -1 & -1 & 1 \\ -1 & -1 & 1 & -1 \end{pmatrix} \begin{pmatrix} f^{02} \\ f^{13} \\ f^{20} \\ f^{31} \end{pmatrix} \tag{11.3.27}$$

Из равенства (11.3.27) следует

$$\begin{pmatrix} -f^{2}_{0} \\ -f^{3}_{1} \\ f^{0}_{2} \\ f^{1}_{3} \end{pmatrix} = \begin{pmatrix} -1 & 1 & -1 & -1 \\ -1 & 1 & 1 & 1 \\ -1 & -1 & -1 & 1 \\ -1 & -1 & 1 & -1 \end{pmatrix} \begin{pmatrix} f^{02} \\ f^{13} \\ f^{20} \\ f^{31} \end{pmatrix}$$

$$= \begin{pmatrix} 1 & 1 & -1 & 1 \\ 1 & 1 & 1 & -1 \\ 1 & -1 & -1 & -1 \\ 1 & -1 & 1 & 1 \end{pmatrix} \begin{pmatrix} -f^{02} \\ f^{13} \\ f^{20} \\ -f^{31} \end{pmatrix}$$

$$
\begin{pmatrix} f_{2}^{0} \\ f_{3}^{1} \\ -f_{0}^{2} \\ -f_{1}^{3} \end{pmatrix} = \begin{pmatrix} 1 & -1 & -1 & -1 \\ 1 & -1 & 1 & 1 \\ 1 & 1 & -1 & 1 \\ 1 & 1 & 1 & -1 \end{pmatrix} \begin{pmatrix} -f^{02} \\ f^{13} \\ f^{20} \\ -f^{31} \end{pmatrix} \tag{11.3.28}
$$

Запишем систему линейных уравнений (11.3.5) в виде произведения матриц

$$
\begin{pmatrix} f_{0}^{3} \\ f_{1}^{2} \\ f_{2}^{1} \\ f_{3}^{0} \end{pmatrix} = \begin{pmatrix} 1 & 1 & -1 & 1 \\ -1 & -1 & -1 & 1 \\ 1 & -1 & -1 & -1 \\ -1 & 1 & -1 & -1 \end{pmatrix} \begin{pmatrix} f^{03} \\ f^{12} \\ f^{21} \\ f^{30} \end{pmatrix} \tag{11.3.29}
$$

Из равенства (11.3.29) следует

$$
\begin{pmatrix} -f_{0}^{3} \\ f_{1}^{2} \\ -f_{2}^{1} \\ f_{3}^{0} \end{pmatrix} = \begin{pmatrix} -1 & -1 & 1 & -1 \\ -1 & -1 & -1 & 1 \\ -1 & 1 & 1 & 1 \\ -1 & 1 & -1 & -1 \end{pmatrix} \begin{pmatrix} f^{03} \\ f^{12} \\ f^{21} \\ f^{30} \end{pmatrix}
= \begin{pmatrix} 1 & 1 & 1 & -1 \\ 1 & 1 & -1 & 1 \\ 1 & -1 & 1 & 1 \\ 1 & -1 & -1 & -1 \end{pmatrix} \begin{pmatrix} -f^{03} \\ -f^{12} \\ f^{21} \\ f^{30} \end{pmatrix}
$$

$$
\begin{pmatrix} f_{3}^{0} \\ -f_{2}^{1} \\ f_{1}^{2} \\ -f_{0}^{3} \end{pmatrix} = \begin{pmatrix} 1 & -1 & -1 & -1 \\ 1 & -1 & 1 & 1 \\ 1 & 1 & -1 & 1 \\ 1 & 1 & 1 & -1 \end{pmatrix} \begin{pmatrix} -f^{03} \\ -f^{12} \\ f^{21} \\ f^{30} \end{pmatrix} \tag{11.3.30}
$$

Мы объединяем равенства (11.3.24), (11.3.26), (11.3.28), (11.3.30) в равенстве (11.3.22). □

Теорема 11.3.4. *Стандартные компоненты линейной функции алгебры кватернионов $H$ относительно базиса (11.3.1) и координаты соответствующего линейного преобразования удовлетворяют соотношениям*

$$
\begin{cases} 4f^{00} = f_{0}^{0} + f_{1}^{1} + f_{2}^{2} + f_{3}^{3} \\ 4f^{11} = -f_{0}^{0} - f_{1}^{1} + f_{2}^{2} + f_{3}^{3} \\ 4f^{22} = -f_{0}^{0} + f_{1}^{1} - f_{2}^{2} + f_{3}^{3} \\ 4f^{33} = -f_{0}^{0} + f_{1}^{1} + f_{2}^{2} - f_{3}^{3} \end{cases} \tag{11.3.31}
$$

(11.3.32)
$$\begin{cases} 4f^{10}= -f^0_1 + f^0_1 - f^2_3 + f^3_2 \\ 4f^{01}= -f^0_1 + f^0_1 + f^2_3 - f^3_2 \\ 4f^{32}= -f^0_1 - f^0_1 - f^2_3 - f^3_2 \\ 4f^{23}= f^0_1 + f^0_1 - f^2_3 - f^3_2 \end{cases}$$

(11.3.33)
$$\begin{cases} 4f^{20}= -f^0_2 + f^1_3 + f^2_0 - f^3_1 \\ 4f^{31}= f^0_2 - f^1_3 + f^2_0 - f^3_1 \\ 4f^{02}= -f^0_2 - f^1_3 + f^2_0 + f^3_1 \\ 4f^{13}= -f^0_2 - f^1_3 - f^2_0 - f^3_1 \end{cases}$$

(11.3.34)
$$\begin{cases} 4f^{30}= -f^0_3 - f^1_2 + f^2_1 + f^3_0 \\ 4f^{21}= -f^0_3 - f^1_2 - f^2_1 - f^3_0 \\ 4f^{12}= f^0_3 - f^1_2 - f^2_1 + f^3_0 \\ 4f^{03}= -f^0_3 + f^1_2 - f^2_1 + f^3_0 \end{cases}$$

Доказательство. Системы линейных уравнений (11.3.31), (11.3.32), (11.3.33), (11.3.34) получены в результате перемножения матриц в равенстве (11.3.23). □

Теорема 11.3.5. *Мы можем отождествить кватернион*

(11.3.35)
$$a = a^0 + a^1 i + a^2 j + a^3 k$$

*и матрицу*

(11.3.36)
$$J_a = \begin{pmatrix} a^0 & -a^1 & -a^2 & -a^3 \\ a^1 & a^0 & -a^3 & a^2 \\ a^2 & a^3 & a^0 & -a^1 \\ a^3 & -a^2 & a^1 & a^0 \end{pmatrix}$$

Доказательство. Произведение кватернионов (11.3.35) и

$$x = x^0 + x^1 i + x^2 j + x^3 k$$

имеет вид

$$\begin{aligned} ax = a^0x^0 - a^1x^1 - a^2x^2 - a^3x^3 + (a^0x^1 + a^1x^0 + a^2x^3 - a^3x^2)i \\ + (a^0x^2 + a^2x^0 + a^3x^1 - a^1x^3)j + (a^0x^3 + a^3x^0 + a^1x^2 - a^2x^1)k \end{aligned}$$

Следовательно, функция $f_a(x) = ax$ имеет матрицу Якоби (11.3.36). Очевидно, что $f_a \circ f_b = f_{ab}$. Аналогичное равенство верно для матриц

$$\begin{pmatrix} a^0 & -a^1 & -a^2 & -a^3 \\ a^1 & a^0 & -a^3 & a^2 \\ a^2 & a^3 & a^0 & -a^1 \\ a^3 & -a^2 & a^1 & a^0 \end{pmatrix} \begin{pmatrix} b^0 & -b^1 & -b^2 & -b^3 \\ b^1 & b^0 & -b^3 & b^2 \\ b^2 & b^3 & b^0 & -b^1 \\ b^3 & -b^2 & b^1 & b^0 \end{pmatrix}$$

$$= \begin{pmatrix} \begin{matrix} a^0b^0 - a^1b^1 \\ -a^2b^2 - a^3b^3 \end{matrix} & \begin{matrix} -a^0b^1 - a^1b^0 \\ -a^2b^3 + a^3b^2 \end{matrix} & \begin{matrix} -a^0b^2 + a^1b^3 \\ -a^2b^0 - a^3b^1 \end{matrix} & \begin{matrix} -a^0b^3 - a^1b^2 \\ +a^2b^1 - a^3b^0 \end{matrix} \\ \\ \begin{matrix} a^0b^1 + a^1b^0 \\ +a^2b^3 - a^3b^2 \end{matrix} & \begin{matrix} a^0b^0 - a^1b^1 \\ -a^2b^2 - a^3b^3 \end{matrix} & \begin{matrix} -a^0b^3 - a^1b^2 \\ +a^2b^1 - a^3b^0 \end{matrix} & \begin{matrix} a^0b^2 - a^1b^3 \\ +a^2b^0 + a^3b^1 \end{matrix} \\ \\ \begin{matrix} a^0b^2 - a^1b^3 \\ +a^2b^0 + a^3b^1 \end{matrix} & \begin{matrix} a^0b^3 + a^1b^2 \\ -a^2b^1 + a^3b^0 \end{matrix} & \begin{matrix} a^0b^0 - a^1b^1 \\ -a^2b^2 - a^3b^3 \end{matrix} & \begin{matrix} -a^0b^1 - a^1b^0 \\ -a^2b^3 - a^3b^2 \end{matrix} \\ \\ \begin{matrix} a^0b^3 + a^1b^2 \\ -a^2b^1 + a^3b^0 \end{matrix} & \begin{matrix} -a^0b^2 + a^1b^3 \\ -a^2b^0 - a^3b^1 \end{matrix} & \begin{matrix} a^0b^1 + a^1b^0 \\ +a^2b^3 - a^3b^2 \end{matrix} & \begin{matrix} a^0b^0 - a^1b^1 \\ -a^2b^2 - a^3b^3 \end{matrix} \end{pmatrix}$$

Следовательно, мы можем отождествить кватернион $a$ и матрицу $J_a$. $\square$

Глава 12

# Линейное отображение $D$-векторных пространств

## 12.1. Линейное отображение $D$-векторных пространств

При изучении линейного отображения $D$-векторных пространств мы будем рассматривать тело $D$ как конечно мерную алгебру над полем $F$.

Определение 12.1.1. Пусть поле $F$ является подкольцом центра $Z(D)$ тела $D$. Пусть $V$ и $W$ - $D$-векторные пространства. Мы будем называть отображение

$$A:V\to W$$

$D$-векторного пространства $V$ в $D$-векторное пространство $W$ **линейным отображением**, если

$$A(x+y)=A(x)+A(y)\qquad x,y\in V$$

$$A(px)=pA(x)\qquad\qquad p\in F$$

Обозначим $\mathcal{L}(D;V;W)$ множество линейных отображений

$$A:V\to W$$

$D$-векторного пространства $V$ в $D$-векторное пространство $W$. □

Очевидно, что линейное отображение $D_*{}^*$-векторного пространства, так же как линейное отображение ${}^*{}_*D$-векторного пространства являются линейными отображениями. Множество морфизмов $D$-векторного пространства шире, чем множество морфизмов $D_*{}^*$-векторного пространства. Чтобы рассмотреть линейное отображение векторных пространств, мы будем следовать методике, предложенной в разделе 5.4.

Теорема 12.1.2. *Пусть $D$ - тело характеристики 0. Пусть $V$, $W$ - $D$-векторные пространства. Линейное отображение*

$$A:V\to W$$

*относительно базиса $\overline{\overline{e}}_V$ ${}_*{}^*D$-векторного пространства $V$ и базиса $\overline{\overline{e}}_W$ ${}_*{}^*D$-векторного пространства $W$ имеет вид*

$$A(v)=e_{W\cdot j}\ A^j_i(v^i)\qquad v=e_{V*}{}^*v \tag{12.1.1}$$

*где $A^j_i(v^i)$ линейно зависит от одной переменной $v^i$ и не зависит от остальных координат вектора $v$.*

Доказательство. Согласно определению 12.1.1

$$A(v)=A(e_{V*}{}^*v)=A\left(\sum_i e_{V\cdot i}\ v^i\right)=\sum_i A(e_{V\cdot i}\ v^i) \tag{12.1.2}$$

Для любого заданного $i$ вектор $A(e_{V\cdot i}\ v^i)\in W$ имеет единственное разложение

$$A(e_{V\cdot i}\ v^i)=e_{W\cdot j}\ A^j_i(v^i)\qquad A(e_{V\cdot i}\ v^i)=e_{W*}{}^*A_i(v^i)\tag{12.1.3}$$

относительно базиса $\overline{\overline{e}}_W$ ${}_*{}^*D$-векторного пространства. Подставив (12.1.3) в (12.1.2), мы получим (12.1.1). □

Определение 12.1.3. Линейное отображение

$$A^j_i:D\to D$$

называется **частным линейным отображением** переменной $v^i$. □

Мы можем записать линейное отображение в виде произведения матриц

$$A(v)=\begin{pmatrix}e_{W\cdot 1}&...&e_{W\cdot m}\end{pmatrix}{}_*{}^*\begin{pmatrix}A^1_i(v^i)\\...\\A^m_i(v^i)\end{pmatrix}\tag{12.1.4}$$

Определим произведение матриц

$$\begin{pmatrix}A^1_1&...&A^1_n\\...&...&...\\A^m_1&...&A^m_n\end{pmatrix}\circ^\circ\begin{pmatrix}v^1\\...\\v^n\end{pmatrix}=\begin{pmatrix}A^1_i(v^i)\\...\\A^m_i(v^i)\end{pmatrix}\tag{12.1.5}$$

где $A=\begin{pmatrix}A^j_i\end{pmatrix}$ - матрица частных линейных отображений. Используя равенство (12.1.5), мы можем записать равенство (12.1.4) в виде

$$A(v)=\begin{pmatrix}e_{W\cdot 1}&...&e_{W\cdot m}\end{pmatrix}{}_*{}^*\begin{pmatrix}A^1_1&...&A^1_n\\...&...&...\\A^m_1&...&A^m_n\end{pmatrix}\circ^\circ\begin{pmatrix}v^1\\...\\v^n\end{pmatrix}\tag{12.1.6}$$

Теорема 12.1.4. *Пусть $D$ - тело характеристики $0$. Линейное отображение*

$$A:v\in V\to w\in W\quad w=A(v)\tag{12.1.7}$$

*относительно базиса $\overline{\overline{e}}_V$ в $D$-векторном пространстве $V$ и базиса $\overline{\overline{e}}_W$ в $D$-векторном пространстве $W$*[12.1] *имеет вид*

$$v=e_{V*}{}^*v$$

$$w=e_{W*}{}^*w\tag{12.1.8}$$

$$w^j=A^j_i(v^i)=A_{s\cdot 0\cdot}{}^j_i\ v^i\ A_{s\cdot 1\cdot}{}^j_i\tag{12.1.9}$$

[12.1] Координатная запись отображения (12.1.7) зависит от выбора базиса. Равенства изменят свой вид, если например мы выберем ${}^*{}_*D$-базис $\overline{\overline{r}}_*$ в $D$-векторном пространстве $W$.

Доказательство. Согласно теореме 12.1.2 линейное отображение $A(v)$ можно записать в виде (12.1.1). Так как для заданных индексов $i$, $j$ частное линейное отображение $A^j_i(v^i)$ линейно по переменной $v^i$, то согласно (10.1.14) выражение $A^j_i(v^i)$ можно представить в виде

$$A^j_i(v^i) = A_{s\cdot 0\cdot}{}^j_i \; v^i \; A_{s\cdot 1\cdot}{}^j_i \tag{12.1.10}$$

где индекс $s$ нумерует слагаемые. Множество значений индекса $s$ зависит от индексов $i$ и $j$. Комбинируя равенства (12.1.2) и (12.1.10), мы получим

$$\begin{aligned} A(v) &= e_{W\cdot j}A^j_i(v^i) &&= e_{W\cdot j}A_{s\cdot 0\cdot}{}^j_i \; v^i \; A_{s\cdot 1\cdot}{}^j_i \\ A(v) &= e_{W*}{}^*A_i(v^i) &&= e_{W*}{}^*(A_{s\cdot 0\cdot i} \; v^i \; A_{s\cdot 1\cdot i}) \end{aligned} \tag{12.1.11}$$

В равенстве (12.1.11) мы суммируем также по индексу $i$. Равенство (12.1.9) следует из сравнения равенств (12.1.8) и (12.1.11). □

Определение 12.1.5. Выражение $A_{s\cdot p\cdot}{}^j_i$ в равенстве (12.1.3) называется **компонентой линейного отображения** $A$. □

Теорема 12.1.6. *Пусть $D$ - тело характеристики 0. Пусть $\overline{\overline{e}}_V$ - базис в $D$-векторном пространстве $V$, $\overline{\overline{e}}_U$ - базис в $D$-векторном пространстве $U$, и $\overline{\overline{e}}_W$ - базис в $D$-векторном пространстве $W$. Предположим, что мы имеем коммутативную диаграмму отображений*

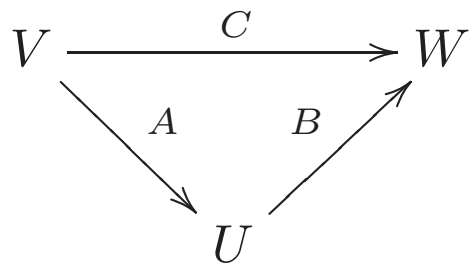

*где линейное отображение $A$ имеет представление*

$$u = A(v) = e_{U\cdot j}A^j_i(v^i) = e_{U\cdot j}A_{s\cdot 0\cdot}{}^j_i \; v^i \; A_{s\cdot 1\cdot}{}^j_i \tag{12.1.12}$$

*относительно заданных базисов и линейное отображение $B$ имеет представление*

$$w = B(u) = e_{W\cdot k}B^k_j(u^j) = e_{W\cdot k}B_{t\cdot 0\cdot}{}^k_j \; u^j \; B_{t\cdot 1\cdot}{}^k_j \tag{12.1.13}$$

*относительно заданных базисов. Тогда отображение $C$ линейно и имеет представление*

$$w = C(v) = e_{W\cdot k}C^k_i(v^i) = e_{W\cdot k}C_{u\cdot 0\cdot}{}^k_i \; v^i \; C_{u\cdot 1\cdot}{}^k_i \tag{12.1.14}$$

*относительно заданных базисов, где*[12.2]

$$\begin{aligned} C^k_i(v^i) && &= B^k_j(A^j_i(v^i)) \\ C_{u\cdot 0\cdot}{}^k_i &= C_{ts\cdot 0\cdot}{}^k_i && = B_{t\cdot 0\cdot}{}^k_j \; A_{s\cdot 0\cdot}{}^j_i \\ C_{u\cdot 1\cdot}{}^k_i &= C_{ts\cdot 1\cdot}{}^k_i && = A_{s\cdot 1\cdot}{}^j_i \; B_{t\cdot 1\cdot}{}^k_j \end{aligned} \tag{12.1.15}$$

[12.2] Индекс $u$ оказался составным индексом, $u = st$. Однако не исключено, что некоторые слагаемые в (12.1.15) могут быть объединены вместе.

Доказательство. Отображение $C$ линейно, так как

$$\begin{aligned}
C(a+b) &= B(A(a+b)) \\
&= B(A(a)+A(b)) \\
&= B(A(a))+B(A(b)) \\
&= C(a)+C(b) \\
C(ab) &= B(A(ab)) = B(aA(b)) \\
&= aB(A(b)) = aC(b) \\
& a, b \in V \\
& a \in F
\end{aligned}$$

Равенство (12.1.14) следует из подстановки (12.1.12) в (12.1.13). □

Теорема 12.1.7. *Для линейного отображения $A$ существует линейное отображение $B$ такое, что*

$$A(axb) = B(x)$$

$$B_{s\cdot 0\cdot}{}^{j}_{i} = A_{s\cdot 0\cdot}{}^{j}_{i}\ a$$

$$B_{s\cdot 1\cdot}{}^{j}_{i} = b\ A_{s\cdot 1\cdot}{}^{j}_{i}$$

Доказательство. Отображение $B$ линейно, так как

$$\begin{aligned}
B(x+y) &= A(a(x+y)b) = A(axb+ayb) = A(axb)+A(ayb) = B(x)+B(y) \\
B(cx) &= A(acxb) = A(caxb) = cA(axb) = cB(x) \\
& x, y \in V, c \in F
\end{aligned}$$

Согласно равенству (12.1.9)

$$B_{s\cdot 0\cdot}{}^{j}_{i}\ v^i\ B_{s\cdot 1\cdot}{}^{j}_{i} = A_{s\cdot 0\cdot}{}^{j}_{i}\ (a\ v^i\ b)A_{s\cdot 1\cdot}{}^{j}_{i} = (A_{s\cdot 0\cdot}{}^{j}_{i}\ a)v^i\ (b\ A_{s\cdot 1\cdot}{}^{j}_{i})$$

□

Теорема 12.1.8. *Пусть $D$ - тело характеристики* 0. *Пусть*

$$A : V \to W$$

*линейное отображение $D$-векторного пространства $V$ в $D$-векторное пространство $W$. Тогда $A(0) = 0$.*

Доказательство. Следствие равенства

$$A(a+0) = A(a) + A(0)$$

□

Определение 12.1.9. Множество

$$\ker f = \{x \in V : f(x) = 0\}$$

называется **ядром линейного отображения**

$$A : V \to W$$

$D$-векторного пространства $V$ в $D$-векторное пространство $W$. □

Определение 12.1.10. Линейное отображение

$$A:V\to W$$

$D$-векторного пространства $V$ в $D$-векторное пространство $W$ называется **вырожденным**, если

$$\ker A=V$$

□

## 12.2. Полилинейное отображение $D$-векторного пространства

Определение 12.2.1. Пусть поле $F$ является подкольцом центра $Z(D)$ тела $D$ характеристики 0. Пусть $V_1$, ..., $V_n$, $W_1$, ..., $W_m$ - $D$-векторные пространства. Мы будем называть отображение

$$\begin{aligned}A:V_1\times...\times V_n&\to W_1\times...\times W_m\\ w_1\times...\times w_m&=A(v_1,...,v_n)\end{aligned}\tag{12.2.1}$$

**полилинейным отображением** $\times$-$D$-векторного пространства $V_1\times...\times V_n$ в $\times$-$D$-векторное пространство $W_1\times...\times W_m$, если

$$A(p_1,...,p_i+q_i,...,p_n)=A(p_1,...,p_i,...,p_n)+A(p_1,...,q_i,...,p_n)$$
$$A(p_1,...,ap_i,...,p_n)=aA(p_1,...,p_i,...,p_n)$$

$$1\le i\le n\quad p_i,q_i\in V_i\quad a\in F$$

□

Определение 12.2.2. Обозначим $\mathcal{L}(D;V_1,...,V_n;W_1,...,W_m)$ множество полилинейных отображений $\times$-$D$-векторного пространства $V_1\times...\times V_n$ в $\times$-$D$-векторное пространство $W_1\times...\times W_m$. □

Теорема 12.2.3. *Пусть $D$ - тело характеристики 0. Для каждого $k\in K=[1,n]$ допустим $\overline{\overline{e}}_{V_k}$ - базис в $D$-векторном пространстве $V_k$ и*

$$v_k={v_k}^*{}_*e_{V_k}\quad v_k\in V_k$$

*Для каждого $l$, $1\le l\le m$, допустим $\overline{\overline{e}}_{W_l}$ - базис в $D$-векторном пространстве $W_l$ и*

$$w_l={w_l}^*{}_*e_{W_l}\quad w_l\in W_l$$

*Полилинейное отображение (12.2.1) относительно базисов $\overline{\overline{e}}_{V_1}$, ..., $\overline{\overline{e}}_{V_n}$, $\overline{\overline{e}}_{W_1}$, ..., $\overline{\overline{e}}_{W_m}$ имеет вид*

$$\begin{aligned}w_l^{\boldsymbol{j}}&=A_{l\cdot}{}^{\boldsymbol{j}}_{\boldsymbol{i}_1...\boldsymbol{i}_n}(v_1^{\boldsymbol{i}_1},...,v_n^{\boldsymbol{i}_n})\\&=A^n_{s\cdot0\cdot l\cdot}{}^{\boldsymbol{j}}_{\boldsymbol{i}_1...\boldsymbol{i}_n}\ \sigma_s(v_1^{i_1})\ A^n_{s\cdot1.l\cdot}{}^{\boldsymbol{j}}_{\boldsymbol{i}_1...\boldsymbol{i}_n}\ ...\ \sigma_s(v_n^{i_n})\ A^n_{s\cdot n\cdot l\cdot}{}^{\boldsymbol{j}}_{\boldsymbol{i}_1...\boldsymbol{i}_n}\end{aligned}\tag{12.2.2}$$

*Область значений $S$ индекса $s$ зависит от значений индексов $\boldsymbol{i}_1$, ..., $\boldsymbol{i}_n$. $\sigma_s$ - перестановка множества переменных $\{v_1^{\boldsymbol{i}_1},...,v_n^{\boldsymbol{i}_n}\}$.*

Доказательство. Так как отображение $A$ в $\times$-$D$-векторное пространство $W_1\times...\times W_m$ можно рассматривать покомпонентно, то мы можем ограничиться рассмотрением отображения

$$A_l:V_1\times...\times V_n\to W_l\quad w_l=A_l(v_1,...,v_n)\tag{12.2.3}$$

Мы докажем утверждение индукцией по $n$.

При $n = 1$ доказываемое утверждение является утверждением теоремы 12.1.4. При этом мы можем отождествить[12.3]

$$A^1_{s\cdot p\cdot}{}^{j}_{i} = A_{s\cdot p\cdot}{}^{j}_{i} \quad p = 0, 1$$

Допустим, что утверждение теоремы справедливо при $n = k - 1$. Тогда отображение (12.2.3) можно представить в виде

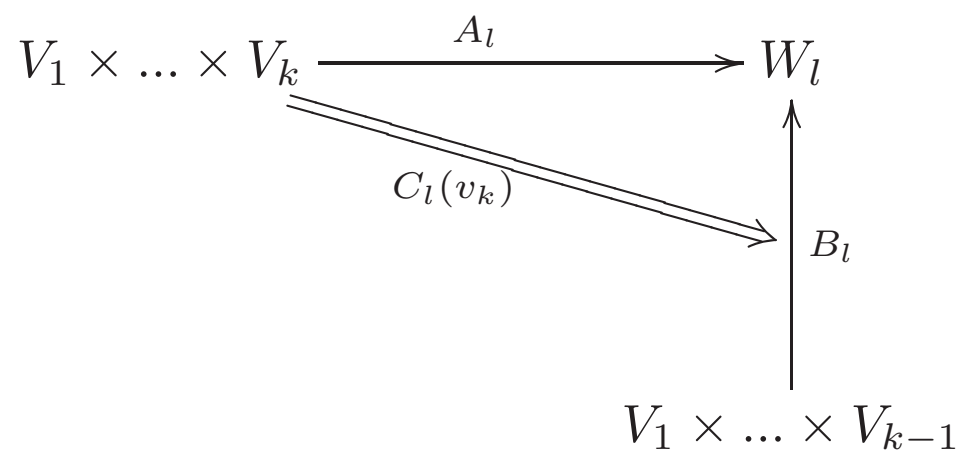

$$w_l = A_l(v_1, ..., v_k) = C_l(v_k)(v_1, ..., v_{k-1})$$

Согласно предположению индукции полилинейное отображение $B_l$ имеет вид

$$w_l^{j} = B^{k-1}_{t\cdot 0\cdot l\cdot}{}^{j}_{i_1...i_{k-1}}\ \sigma_t(v_1^{i_1})\ B^{k-1}_{t\cdot 1.l\cdot}{}^{j}_{i_1...i_{k-1}}\ ...\ \sigma_t(v_{k-1}^{i_{k-1}})\ B^{k-1}_{t\cdot k-1.l\cdot}{}^{j}_{i_1...i_{k-1}}$$

Согласно построению $B_l = C_l(v_k)$. Следовательно, выражения $B^{k-1}_{t\cdot p\cdot l\cdot}{}^{j}_{i_1...i_{k-1}}$ являются функциями $v_k$. Поскольку $C_l(v_k)$ - линейная функция $v_k$, то только одно выражение $B^{k-1}_{t\cdot p\cdot l\cdot}{}^{j}_{i_1...i_{k-1}}$ является линейной функцией $v_k$, и остальные выражения $B^{k-1}_{t\cdot q\cdot l\cdot}{}^{j}_{i_1...i_{k-1}}$ не зависят от $v_k$.

Не нарушая общности, положим $p = 0$. Согласно теореме 12.1.4

$$B^{k-1}_{t\cdot 0\cdot l\cdot}{}^{j}_{i_1...i_{k-1}} = C^k_{tr\cdot 0\cdot l\cdot}{}^{j}_{i_k i_1...i_{k-1}}\ v_k^{i_k}\ C^k_{tr\cdot 1.l\cdot}{}^{j}_{i_k i_1...i_{k-1}}$$

Положим $s = tr$ и определим перестановку $\sigma_s$ согласно правилу

$$\sigma_s = \sigma(tr) = \begin{pmatrix} v_k^{i_k} & v_1^{i_1} & ... & v_{k-1}^{i_{k-1}} \\ v_k^{i_k} & \sigma_t(v_1^{i_1}) & ... & \sigma_t(v_{k-1}^{i_{k-1}}) \end{pmatrix}$$

Положим

$$A^k_{tr\cdot q+1.l\cdot}{}^{j}_{i_k i_1...i_{k-1}} = B^{k-1}_{t\cdot q\cdot l\cdot}{}^{j}_{i_1...i_{k-1}} \qquad q = 1, ..., k-1$$
$$A^k_{tr\cdot q\cdot l\cdot}{}^{j}_{i_k i_1...i_{k-1}} = C^k_{tr\cdot q\cdot l\cdot}{}^{j}_{i_k i_1...i_{k-1}} \qquad q = 0, 1$$

Мы доказали шаг индукции. □

Определение 12.2.4. Выражение $A^n_{s\cdot p\cdot l\cdot}{}^{j}_{i_1...i_n}$ в равенстве (12.2.2) называется **компонентой полилинейного отображения** $A$. □

[12.3] В представлении (12.2.2) мы будем пользоваться следующими правилами.

- Если область значений какого-либо индекса - это множество, состоящее из одного элемента, мы будем опускать соответствующий индекс.
- Если $n = 1$, то $\sigma_s$ - тождественное преобразование. Это преобразование можно не указывать в выражении.

Глава 13

# Тензорное произведение

## 13.1. Тензорное произведение колец

Теорема 13.1.1. *Пусть $G$ - полугруппа. Пусть $M$ - свободный модуль с базисом $G$ над коммутативным кольцом $F$. Тогда на $M$ определена структура кольца таким образом, что $F$ является подкольцом центра $Z(M)$ кольца $M$.*

Доказательство. Произвольные векторы $a$, $b$, $c \in M$ имеют единственное разложение

$$a = a^i g_i \qquad b = b^j g_j \qquad c = c^k g_k$$

относительно базиса $G$. $i \in I$, $j \in J$, $k \in K$, где $I$, $J$, $K$ - конечные множества. Не нарушая общности, мы можем положить $i$, $j$, $k \in I \cup J \cup K$.

Для $a$ и $b$ определим произведение равенством

$$ab = (a^i b^j)(g_i g_j) \tag{13.1.1}$$

Из цепочки равенств

$$\begin{aligned}
(a+b)c &= ((a^i + b^i)g_i)(c^j g_j)\\
&= ((a^i + b^i)c^j)(g_i g_j)\\
&= (a^i c^j + b^i c^j)(g_i g_j)\\
&= (a^i c^j)(g_i g_j) + (b^i c^j)(g_i g_j)\\
&= (a^i g_i)(c^j g_j) + (b^i g_i)(c^j g_j)\\
&= ac + bc
\end{aligned}$$

следует, что предложенное произведение дистрибутивно справа по отношению к сложению. Аналогично доказывается, что предложенное произведение дистрибутивно слева по отношению к сложению. Из цепочки равенств

$$\begin{aligned}
(ab)c &= ((a^i g_i)(b^j g_j))(c^k g_k)\\
&= ((a^i b^j)(g_i g_j))(c^k g_k)\\
&= ((a^i b^j)c^k)((g_i g_j)g_k)\\
&= (a^i (b^j c^k))(g_i (g_j g_k))\\
&= (a^i g_i)((b^j c^k)(g_j g_k))\\
&= (a^i g_i)((b^j g_j)(c^k g_k))\\
&= a(bc)
\end{aligned}$$

следует, что предложенное произведение ассоциативно.

Пусть $e \in G$ - единица полугруппы $G$. Из равенства

$$ae = (a^i g_i)e = a^i(g_i e) = a^i g_i = a$$

следует, что $e$ является единицей кольца $M$.

Отображение

$$f : p \in F \rightarrow pe \in M$$

является гомоморфным вложением кольца $F$ в кольцо $M$, откуда следует, что кольцо $F$ содержится в центре $Z(M)$ кольца $M$. □

Пусть $R_1$, ..., $R_n$ - кольца характеристики 0.[13.1] Пусть $F$ - максимальное коммутативное кольцо, являющееся подкольцом центра $Z(R_i)$ для любого i, $i = 1,...,n$. Рассмотрим категорию $\mathcal{A}$ объектами которой являются полилинейные над коммутативным кольцом $F$ отображения

$$f : R_1 \times ... \times R_n \longrightarrow S_1 \qquad g : R_1 \times ... \times R_n \longrightarrow S_2$$

где $S_1$, $S_2$ - модули над кольцом $F$. Мы определим морфизм $f \rightarrow g$ как линейное над коммутативным кольцом $F$ отображение $h : S_1 \rightarrow S_2$, для которого коммутативна диаграмма

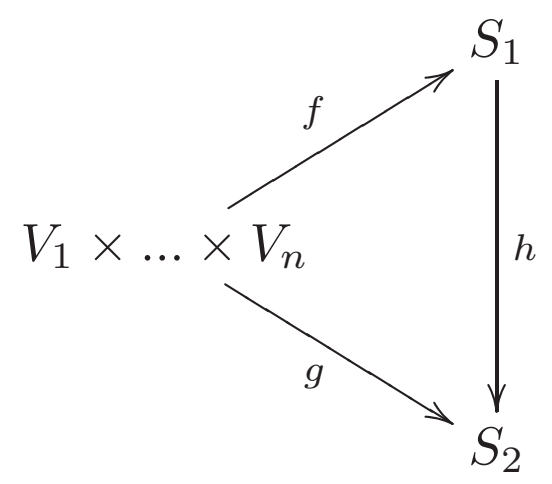

Универсальный объект $R_1 \otimes ... \otimes R_n$ категории $\mathcal{A}$ называется **тензорным произведением колец $R_1$, ..., $R_n$ над коммутативным кольцом $F$.**

Теорема 13.1.2. *Тензорное произведение колец существует.*

Доказательство. Пусть $M$ - модуль над кольцом $F$, порождённый произведением $R_1 \times ... \times R_n$ мультипликативных полугрупп колец $R_1$, ..., $R_n$. Инъекция

$$i : R_1 \times ... \times R_n \longrightarrow M$$

определена по правилу

(13.1.2) $$i(d_1, ..., d_n) = (d_1, ..., d_n)$$

и является гомоморфизмом мультипликативной полугруппы кольца $R_1 \times ... \times R_n$ на базис модуля $M$. Следовательно, произведение векторов базиса определено покомпонентно

(13.1.3) $$(d_1, ..., d_n)(c_1, ..., c_n) = (d_1 c_1, ..., d_n c_n)$$

Пусть $N \subset M$ - подмодуль, порождённый элементами вида

(13.1.4) $$\begin{gathered}(d_1, ..., d_i + c_i, ..., d_n) - (d_1, ..., d_i, ..., d_n) - (d_1, ..., c_i, ..., d_n)\\(d_1, ..., ad_i, ..., d_n) - a(d_1, ..., d_i, ..., d_n)\end{gathered}$$

[13.1] Я определяю тензорное произведение колец по аналогии с определением в [1], с. 456 - 458.

где $d_i, c_i \in R_i$, $a \in F$. Пусть

$$j : M \to M/N$$

каноническое отображение на фактормодуль. Рассмотрим коммутативную диаграмму

(13.1.5)

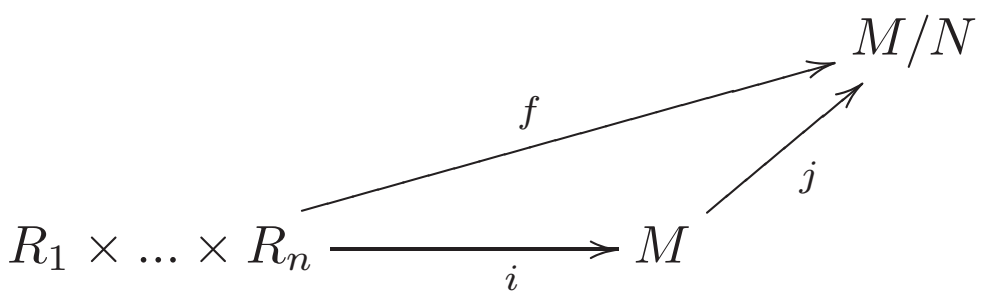

Поскольку элементы (13.1.4) принадлежат ядру линейного отображения $j$, то из равенства (13.1.2) следует

$$\begin{aligned} f(d_1, ..., d_i + c_i, ..., d_n) &= f(d_1, ..., d_i, ..., d_n) + f(d_1, ..., c_i, ..., d_n) \\ f(d_1, ..., ad_i, ..., d_n) &= af(d_1, ..., d_i, ..., d_n) \end{aligned} \tag{13.1.6}$$

Из равенств (13.1.6) следует, что отображение $f$ полилинейно над полем $F$. Поскольку $M$ - модуль с базисом $R_1 \times ... \times R_n$, то согласно теореме [1]-1, на с. 104 для любого модуля $V$ и любого полилинейного над $F$ отображения

$$g : R_1 \times ... \times R_n \longrightarrow V$$

существует единственный гомоморфизм $k : M \to V$, для которого коммутативна следующая диаграмма

(13.1.7)

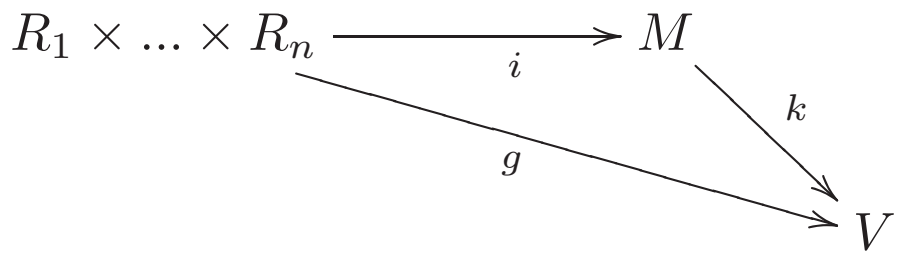

Так как $g$ - полилинейно над $F$, то $\ker k \subseteq N$. Согласно утверждению на с. [1]-94 отображение $j$ универсально в категории гомоморфизмов векторного пространства $M$, ядро которых содержит $N$. Следовательно, определён гомоморфизм

$$h : M/N \to V$$

для которого коммутативна диаграмма

(13.1.8)

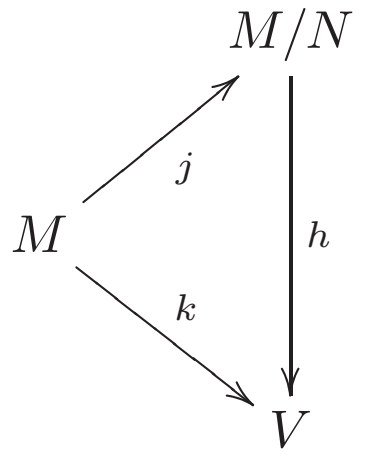

Объединяя диаграммы (13.1.5), (13.1.7), (13.1.8), получим коммутативную диаграмму

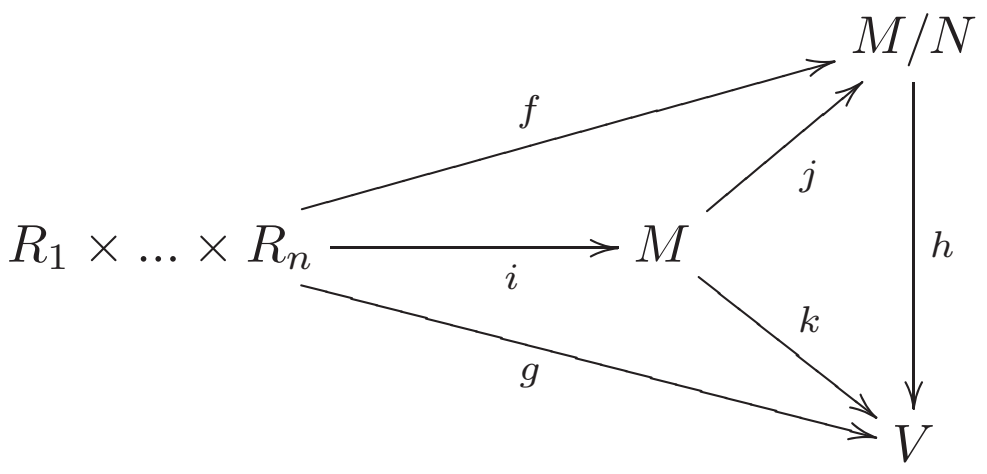

Так как $\mathrm{Im} f$ порождает $M/N$, то отображение $h$ однозначно определено. □

Согласно доказательству теоремы 13.1.2

$$R_1 \otimes ... \otimes R_n = M/N$$

Для $d_i \in R_i$ будем записывать

$$d_1 \otimes ... \otimes d_n = j(d_1, ..., d_n) \tag{13.1.9}$$

Из равенств (13.1.2) и (13.1.9) следует

$$f(d_1, ..., d_n) = d_1 \otimes ... \otimes d_n$$

Равенства (13.1.6) можно записать в виде

$$\begin{aligned} &d_1 \otimes ... \otimes (d_i + c_i) \otimes ... \otimes d_n \\ =&d_1 \otimes ... \otimes d_i \otimes ... \otimes d_n + d_1 \otimes ... \otimes c_i \otimes ... \otimes d_n \\ &d_1 \otimes ... \otimes (ad_i) \otimes ... \otimes d_n = a(d_1 \otimes ... \otimes d_i \otimes ... \otimes d_n) \end{aligned} \tag{13.1.10}$$

Теорема 13.1.3. *Тензорное произведение $R_1 \otimes ... \otimes R_n$ колец $R_1$, ..., $R_n$ характеристики 0 над кольцом $F$ является кольцом.*

Доказательство. Согласно доказательству теоремы 13.1.2 базис модуля $M$ над кольцом $F$ является полугруппой. Согласно теореме 13.1.1 в модуле $M$ определена структура кольца. Подмодуль $N$ является двусторонним идеалом кольца $M$. Каноническое отображение на фактормодуль

$$j : M \to M/N$$

является также каноническим гомоморфизмом в факторкольцо. Следовательно, в модуле $M/N = R_1 \otimes ... \otimes R_n$ определена структура кольца. Из равенства (13.1.3) следует

$$(d_1 \otimes ... \otimes d_n)(c_1 \otimes ... \otimes c_n) = (d_1c_1) \otimes ... \otimes (d_nc_n)$$

Из равенства (13.1.10), следует

$$\begin{aligned} &d_1 \otimes ... \otimes d_i \otimes ... \otimes d_n \\ =&d_1 \otimes ... \otimes (d_i + 0) \otimes ... \otimes d_n \\ =&d_1 \otimes ... \otimes d_i \otimes ... \otimes d_n + d_1 \otimes ... \otimes 0 \otimes ... \otimes d_n \end{aligned}$$

Следовательно, мы можем отождествить тензор $d_1 \otimes ... \otimes 0 \otimes ... \otimes d_n$ с нулём $0 \otimes ... \otimes 0$. В силу условия теоремы произведение равно $0 \otimes ... \otimes 0$ только, если один из сомножителей равен $0 \otimes ... \otimes 0$. Единицей произведения служит тензор $e_1 \otimes ... \otimes e_n$. □

## 13.2. Тензорное произведение тел

Пусть $D_1$, ..., $D_n$ - тела характеристики 0. Пусть $F$ - поле, являющееся подкольцом центра $Z(D_i)$ для любого i, $i = 1,...,n$. Тензор

$$a_1 \otimes ... \otimes a_n \in D_1 \otimes ... \otimes D_n$$

имеет обратный

$$(a_1 \otimes ... \otimes a_n)^{-1} = (a_1)^{-1} \otimes ... \otimes (a_n)^{-1}$$

Однако **тензорное произведение** $D_1 \otimes ... \otimes D_n$ **тел** $D_1$, **...,** $D_n$, вообще говоря, не является телом, так как мы не можем сказать, является ли обратим элемент

$$p(a_1 \otimes ... \otimes a_n) + q(b_1 \otimes ... \otimes b_n)$$

Замечание 13.2.1. Представление тензора в виде

$$a^s\ d_{s\cdot 1} \otimes ... \otimes d_{s\cdot n} \tag{13.2.1}$$

неоднозначно. Если $r \in F$, то

$$(d_1 r) \otimes d_2 = d_1 \otimes (r d_2)$$

Чисто алгебраическими методами мы можем увеличить либо уменьшить число слагаемых. К числу допустимых преобразований относится

$$\begin{aligned} &d_1 \otimes d_2 + c_1 \otimes c_2 \\ =&d_1 \otimes (d_2 - c_2 + c_2) + c_1 \otimes c_2 \\ =&d_1 \otimes (d_2 - c_2) + d_1 \otimes c_2 + c_1 \otimes c_2 \\ =&d_1 \otimes (d_2 - c_2) + (d_1 + c_1) \otimes c_2 \end{aligned}$$

□

В качестве иллюстрации замечания 13.2.1 рассмотрим следующую теорему.

Теорема 13.2.2. *Если в теле $D$ характеристики 0 существуют элементы $a$, $b$, $c$, произведение которых не коммутативно, то существует нетривиальная запись нулевого тензора.*

Доказательство. Доказательство теоремы следует из цепочки равенств

$$\begin{aligned} 0 \otimes 0 =&a \otimes a - a \otimes a + b \otimes b - b \otimes b \\ =&(a + c) \otimes a - c \otimes a - a \otimes (a - c) - a \otimes c \\ &+ (b - c) \otimes b + c \otimes b - b \otimes (b + c) + b \otimes c \\ =&(a + c) \otimes a - a \otimes (a - c) + (b - c) \otimes b \\ &+ c \otimes (b - a) - b \otimes (b + c) + (b - a) \otimes c \end{aligned}$$

□

Неоднозначность представления тензора приводит к тому, что мы должны найти каноническое представление тензора. Мы можем найти решение поставленной задачи в случае тензорного произведения тел. Допустим тело $D_i$ является векторным пространством над полем $F \subset Z(D_i)$. Допустим это векторное пространство имеет конечный базис $e_{i\cdot \boldsymbol{s}_i}$, $\boldsymbol{s}_i \in \boldsymbol{S}_i$, $|\boldsymbol{S}_i| = m_i$. Следовательно, любой элемент $f_i$ тела $D_i$ можно представить в виде

$$f_i = e_{i\cdot \boldsymbol{s}_i} f_i^{\boldsymbol{s}_i}$$

В этом случае тензор (13.2.1) можно записать в виде

$$a^s(f_{s\cdot 1}^{\ \ s_1} e_{1\cdot s_1})\otimes ...\otimes (f_{s\cdot n}^{\ \ s_n} e_{n\cdot s_n}) \tag{13.2.2}$$

где $a^s$, ${}_s f_1^{s_1} \in F$. Мы можем упростить выражение (13.2.2)

$$a^s f_{s\cdot 1}^{\ \ s_1} ... f_{s\cdot n}^{\ \ s_n} e_{1\cdot s_1}\otimes ...\otimes e_{n\cdot s_n} \tag{13.2.3}$$

Положим

$$a^s f_{s\cdot 1}^{\ \ s_1} ... f_{s\cdot n}^{\ \ s_n} = f^{s_1...s_n}$$

Тогда равенство (13.2.3) примет вид

$$f^{s_1...s_n} e_{1\cdot s_1}\otimes ...\otimes e_{n\cdot s_n} \tag{13.2.4}$$

Выражение (13.2.4) определенно однозначно с точностью до выбора базиса. Мы будем называть выражение $f^{s_1...s_n}$ **стандартной компонентой тензора**.

## 13.3. Тензорное произведение $D_*{}^*$-векторных пространств

Пусть $D_1$, ..., $D_n$ - тела характеристики 0.[13.2] Пусть $V_i$ - $D_{i*}{}^*$-векторное пространство, $i = 1$, ..., $n$. Рассмотрим категорию $\mathcal{A}$, объектами которой являются полилинейные отображения

$$f : V_1 \times ... \times V_n \longrightarrow W_1 \qquad g : V_1 \times ... \times V_n \longrightarrow W_2$$

где $W_1$, $W_2$ - $D\star$-модули над кольцом $D_1 \otimes ... \otimes D_n$. Мы определим морфизм $f \to g$ как аддитивное отображение

$$h : W_1 \to W_2$$

для которого коммутативна диаграмма

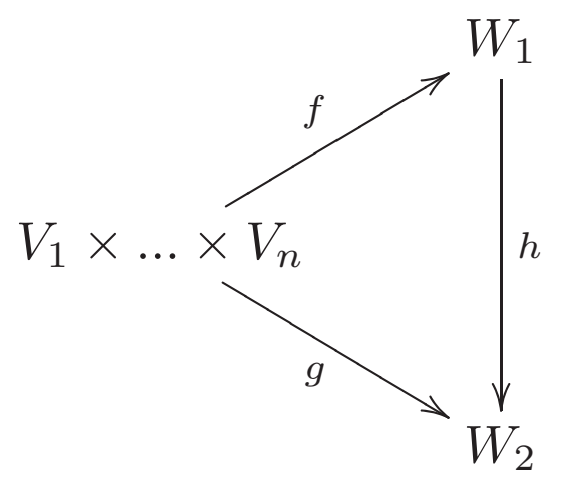

Универсальный объект $V_1 \otimes ... \otimes V_n$ категории $\mathcal{A}$ называется **тензорным произведением $D\star$-векторных пространств** $V_1$, ..., $V_n$.

Теорема 13.3.1. *Тензорное произведение $D\star$-векторных пространств существует.*

Доказательство. Пусть $F$ - поле, являющееся подкольцом центра $Z(D_i)$ для любого i, $i = 1,...,n$.

Пусть $D$ - свободное векторное пространство над кольцом $F$, порождённое произведением $D_1 \times ... \times D_n$ мультипликативных полугрупп тел $D_1$, ..., $D_n$. Инъекция

$$i' : D_1 \times ... \times D_n \longrightarrow D$$

определена по правилу

$$i'(d_1, ..., d_n) = (d_1, ..., d_n) \tag{13.3.1}$$

[13.2] Я определяю тензорное произведение $D_*{}^*$-векторных пространств по аналогии с определением в [1], с. 456 - 458.

и является гомоморфизмом мультипликативной полугруппы кольца $D_1 \times ... \times D_n$ на базис модуля $D$. Следовательно, произведение векторов базиса определено покомпонентно

(13.3.2) $$(d_1, ..., d_n)(c_1, ..., c_n) = (d_1c_1, ..., d_nc_n)$$

Согласно теореме 13.1.1 на векторном пространстве $D$ определена структура кольца.

Рассмотрим прямое произведение $F \times D$ поля $F$ и кольца $D$. Мы будем отождествлять элемент $(f, e)$ с элементом $f \in F$ и элемент $(1, d)$ с элементом $d \in D$.

Пусть $M$ - свободный модуль над кольцом $F \times D$, порождённый декартовым произведением $V_1 \times ... \times V_n$ . Если $v_1 \in V_1$, ..., $v_n \in V_n$, то соответствующий вектор из $M$ мы будем обозначать $(v_1, ..., v_n)$. Пусть

$$i : V_1 \times ... \times V_n \longrightarrow M$$

инъекция, определённая по правилу

(13.3.3) $$i(v_1, ..., v_n) = (v_1, ..., v_n)$$

Пусть $N \subset M$ - векторное подпространство, порождённое элементами вида

(13.3.4) $$\begin{matrix}(v_1, ..., v_i + w_i, ..., v_n) - (v_1, ..., v_i, ..., v_n) - (v_1, ..., w_i, ..., v_n)\\(v_1, ..., av_i, ..., v_n) - a(v_1, ..., v_i, ..., v_n)\end{matrix}$$

где $v_i$, $w_i \in V_i$, $a \in F$. Пусть

$$j : M \to M/N$$

каноническое отображение на фактормодуль. Рассмотрим коммутативную диаграмму

(13.3.5)

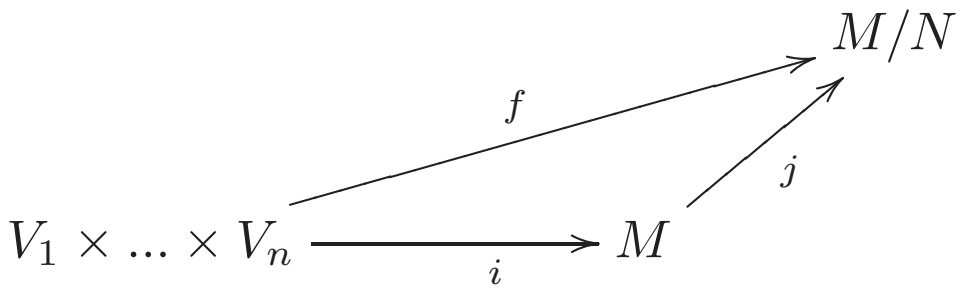

Поскольку элементы (13.3.4) принадлежат ядру линейного отображения $j$, то из равенства (13.3.3) следует

(13.3.6) $$\begin{aligned}f(v_1, ..., v_i + w_i, ..., v_n) =& f(v_1, ..., v_i, ..., v_n) + f(v_1, ..., w_i, ..., v_n)\\f(v_1, ..., av_i, ..., v_n) =& af(v_1, ..., v_i, ..., v_n)\end{aligned}$$

Из равенств (13.3.6) следует, что отображение $f$ полилинейно над полем $F$. Поскольку $M$ - модуль с базисом $V_1 \times ... \times V_n$, то согласно теореме [1]-1, на с. 104 для любого модуля $V$ и любого полилинейного над $F$ отображения

$$g : V_1 \times ... \times V_n \longrightarrow V$$

существует единственный гомоморфизм $k : M \to V$, для которого коммутативна следующая диаграмма

(13.3.7)

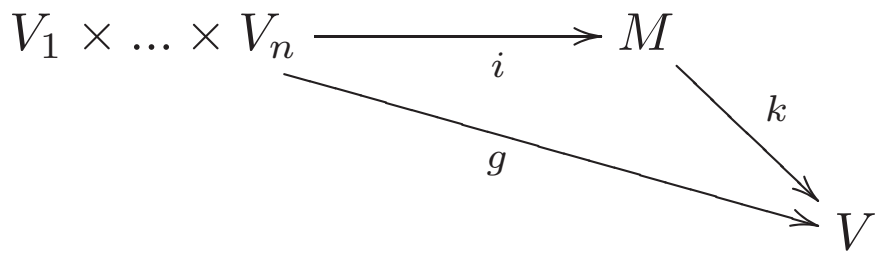

Так как $g$ - полилинейно над $F$, то $\ker k \subseteq N$. Согласно утверждению на с. [1]-94 отображение $j$ универсально в категории гомоморфизмов модуля $M$, ядро которых содержит $N$. Следовательно, определён гомоморфизм

$$h : M/N \to V$$

для которого коммутативна диаграмма

(13.3.8)

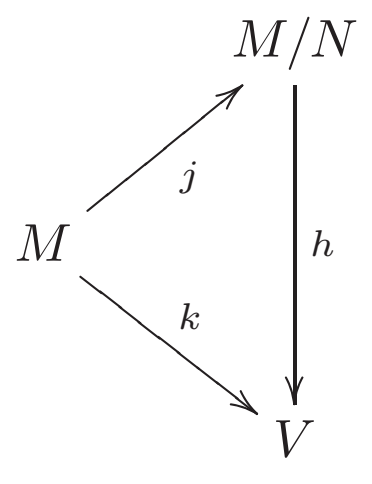

Объединяя диаграммы (13.3.5), (13.3.7), (13.3.8), получим коммутативную диаграмму

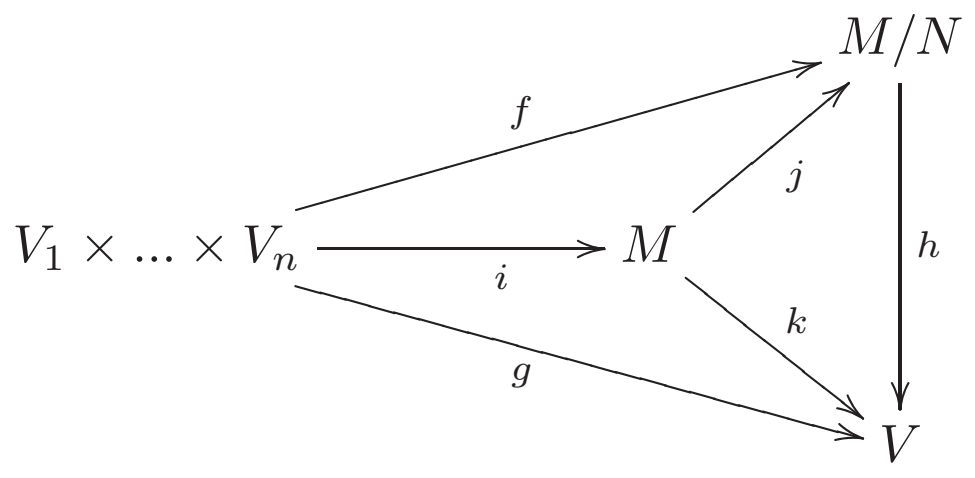

Так как $\mathrm{Im} f$ порождает $M/N$, то отображение $h$ однозначно определено. □

Согласно доказательству теоремы 13.3.1

$$V_1 \otimes ... \otimes V_n = M/N$$

Для $v_i \in V_i$ будем записывать

$$v_1 \otimes ... \otimes v_n = j(v_1, ..., v_n) \tag{13.3.9}$$

Из равенств (13.3.3) и (13.3.9) следует

$$f(v_1, ..., v_n) = v_1 \otimes ... \otimes v_n$$

Равенства (13.3.6), можно записать в виде

$$\begin{aligned} & v_1 \otimes ... \otimes (v_i + w_i) \otimes ... \otimes v_n \\ =& v_1 \otimes ... \otimes v_i \otimes ... \otimes v_n + v_1 \otimes ... \otimes w_i \otimes ... \otimes v_n \\ & v_1 \otimes ... \otimes (a v_i) \otimes ... \otimes v_n = a(v_1 \otimes ... \otimes v_i \otimes ... \otimes v_n) \end{aligned}$$

Теорема 13.3.2. *Тензорное произведение $V_1 \otimes ... \otimes V_n$ является модулем над тензорным произведением $D_1 \otimes ... \otimes D_n$.*

Доказательство. Для доказательства утверждения теоремы необходимо доказать, что действие

$$(d_1 \otimes ... \otimes d_n)(v_1 \otimes ... \otimes v_n) = (d_1, ..., d_n)(v_1 \otimes ... \otimes v_n)$$

кольца $D_1\otimes...\otimes D_n$ на модуле $V_1\otimes...\otimes V_n$ определено корректно. Это следует из цепочек равенств

$$\begin{aligned}(d_1\otimes...\otimes(ad_i)\otimes...\otimes d_n)(v_1\otimes...\otimes v_n)&=(d_1,...,ad_i,...,d_n)(v_1\otimes...\otimes v_n)\\&=(d_1v_1)\otimes...\otimes(ad_iv_i)\otimes...\otimes(d_nv_n)\\&=a(d_1v_1)\otimes...\otimes(d_nv_n)\\&=a((d_1,...,d_n)(v_1\otimes...\otimes v_n))\\&=a((d_1\otimes...\otimes d_n)(v_1\otimes...\otimes v_n))\\&=(a(d_1\otimes...\otimes d_n))(v_1\otimes...\otimes v_n)\end{aligned}$$

$$\begin{aligned}&(d_1\otimes...\otimes(d_i+c_i)\otimes...\otimes d_n)(v_1\otimes...\otimes v_n)\\=&(d_1,...,d_i+c_i,...,d_n)(v_1\otimes...\otimes v_n)\\=&(d_1v_1)\otimes...\otimes((d_i+c_i)v_i)\otimes...\otimes(d_nv_n)\\=&(d_1v_1)\otimes...\otimes(d_iv_i+c_iv_i)\otimes...\otimes(d_nv_n)\\=&(d_1v_1)\otimes...\otimes(d_iv_i)\otimes...\otimes(d_nv_n)+(d_1v_1)\otimes...\otimes(c_iv_i)\otimes...\otimes(d_nv_n)\\=&(d_1v_1)\otimes...\otimes(d_iv_i)\otimes...\otimes(d_nv_n)+(d_1v_1)\otimes...\otimes(c_iv_i)\otimes...\otimes(d_nv_n)\\=&(d_1,...,d_i,...,d_n)(v_1\otimes...\otimes v_n)+(d_1,...,c_i,...,d_n)(v_1\otimes...\otimes v_n)\\=&(d_1\otimes...\otimes d_i\otimes...\otimes d_n)(v_1\otimes...\otimes v_n)+(d_1\otimes...\otimes c_i\otimes...\otimes d_n)(v_1\otimes...\otimes v_n)\\=&((d_1\otimes...\otimes d_i\otimes...\otimes d_n)+(d_1\otimes...\otimes c_i\otimes...\otimes d_n))(v_1\otimes...\otimes v_n)\end{aligned}$$

□

## 13.4. Тензорное произведение $D$-векторных пространств

Пусть $D_1$, ..., $D_n$ - тела характеристики 0.[13.3] Пусть $V_i$ - $D_i$-векторное пространство, $i=1$, ..., $n$. Рассмотрим категорию $\mathcal{A}$ объектами которой являются полилинейные отображения

$$f:V_1\times...\times V_n\longrightarrow W_1\qquad g:V_1\times...\times V_n\longrightarrow W_2$$

где $W_1$, $W_2$ - модули над кольцом $D_1\otimes...\otimes D_n$. Мы определим морфизм $f\to g$ как линейное отображение

$$h:W_1\to W_2$$

для которого коммутативна диаграмма

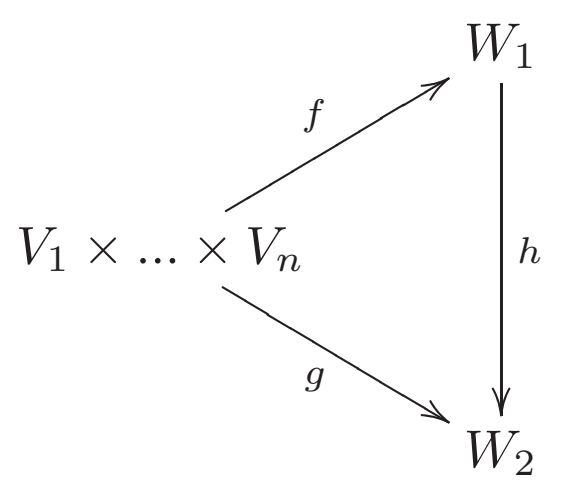

Универсальный объект $V_1\otimes...\otimes V_n$ категории $\mathcal{A}$ называется **тензорным произведением $D$-векторных пространств** $V_1$, ..., $V_n$.

[13.3] Я определяю тензорное произведение $D$-векторных пространств по аналогии с определением в [1], с. 456 - 458.

Теорема 13.4.1. *Тензорное произведение $D$-векторных пространств существует.*

Доказательство. Пусть $F$ - поле, являющееся подкольцом центра $Z(D_i)$ для любого i, $i = 1,...,n$.

Пусть $D$ - свободное векторное пространство над кольцом $F$, порождённое произведением $D_1 \times ... \times D_n$ мультипликативных полугрупп тел $D_1$, ..., $D_n$. Инъекция

$$i' : D_1 \times ... \times D_n \longrightarrow D$$

определена по правилу

$$i'(d_1, ..., d_n) = (d_1, ..., d_n) \tag{13.4.1}$$

и является гомоморфизмом мультипликативной полугруппы кольца $D_1 \times ... \times D_n$ на базис модуля $D$. Следовательно, произведение векторов базиса определено покомпонентно

$$(d_1, ..., d_n)(c_1, ..., c_n) = (d_1c_1, ..., d_nc_n) \tag{13.4.2}$$

Согласно теореме 13.1.1 на векторном пространстве $D$ определена структура кольца.

Рассмотрим прямое произведение $F \times D$ поля $F$ и кольца $D$. Мы будем отождествлять элемент $(f, e)$ с элементом $f \in F$ и элемент $(1, d)$ с элементом $d \in D$.

Пусть $M$ - свободный модуль над кольцом $F \times D$, порождённый декартовым произведением $V_1 \times ... \times V_n$ . Если $v_1 \in V_1$, ..., $v_n \in V_n$, то соответствующий вектор из $M$ мы будем обозначать $(v_1, ..., v_n)$. Пусть

$$i : V_1 \times ... \times V_n \longrightarrow M$$

инъекция, определённая по правилу

$$i(v_1, ..., v_n) = (v_1, ..., v_n) \tag{13.4.3}$$

Пусть $N \subset M$ - векторное подпространство, порождённое элементами вида

$$\begin{gathered}(v_1, ..., v_i + w_i, ..., v_n) - (v_1, ..., v_i, ..., v_n) - (v_1, ..., w_i, ..., v_n)\\(v_1, ..., av_i, ..., v_n) - a(v_1, ..., v_i, ..., v_n)\end{gathered} \tag{13.4.4}$$

где $v_i$, $w_i \in V_i$, $a \in F$. Пусть

$$j : M \to M/N$$

каноническое отображение на фактормодуль. Рассмотрим коммутативную диаграмму

(13.4.5)

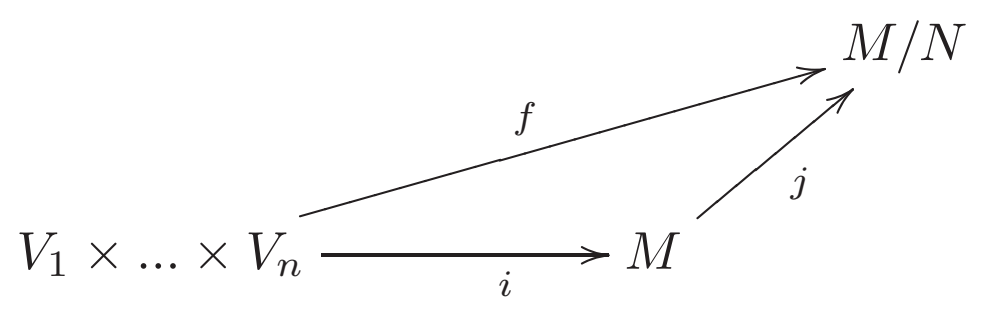

Поскольку элементы (13.4.4) принадлежат ядру линейного отображения $j$, то из равенства (13.4.3) следует

$$\begin{aligned}f(v_1, ..., v_i + w_i, ..., v_n) &= f(v_1, ..., v_i, ..., v_n) + f(v_1, ..., w_i, ..., v_n)\\f(v_1, ..., av_i, ..., v_n) &= af(v_1, ..., v_i, ..., v_n)\end{aligned} \tag{13.4.6}$$

Из равенств (13.4.6) следует, что отображение $f$ полилинейно над полем $F$. Поскольку $M$ - модуль с базисом $V_1 \times ... \times V_n$, то согласно теореме [1]-1, на с. 104 для любого модуля $V$ и любого полилинейного над $F$ отображения

$$g : V_1 \times ... \times V_n \longrightarrow V$$

существует единственный гомоморфизм $k : M \to V$, для которого коммутативна следующая диаграмма

(13.4.7)

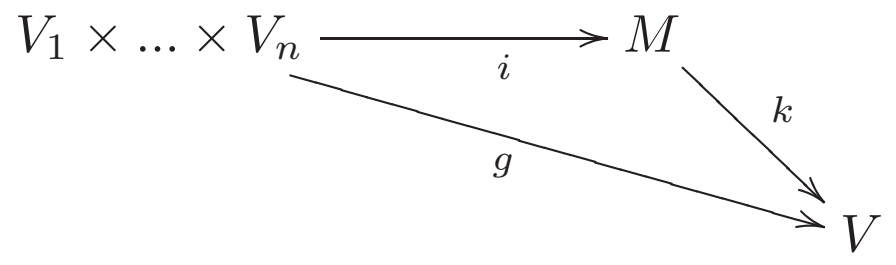

Так как $g$ - полилинейно над $F$, то $\ker k \subseteq N$. Согласно утверждению на с. [1]-94 отображение $j$ универсально в категории гомоморфизмов модуля $M$, ядро которых содержит $N$. Следовательно, определён гомоморфизм

$$h : M/N \to V$$

для которого коммутативна диаграмма

(13.4.8)

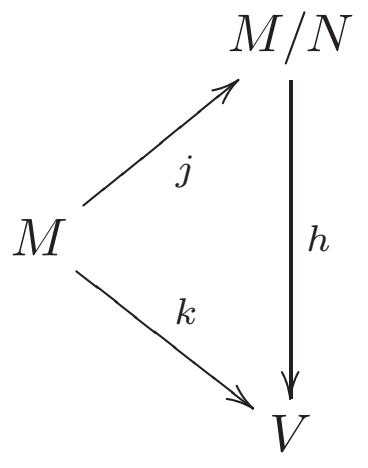

Объединяя диаграммы (13.4.5), (13.4.7), (13.4.8), получим коммутативную диаграмму

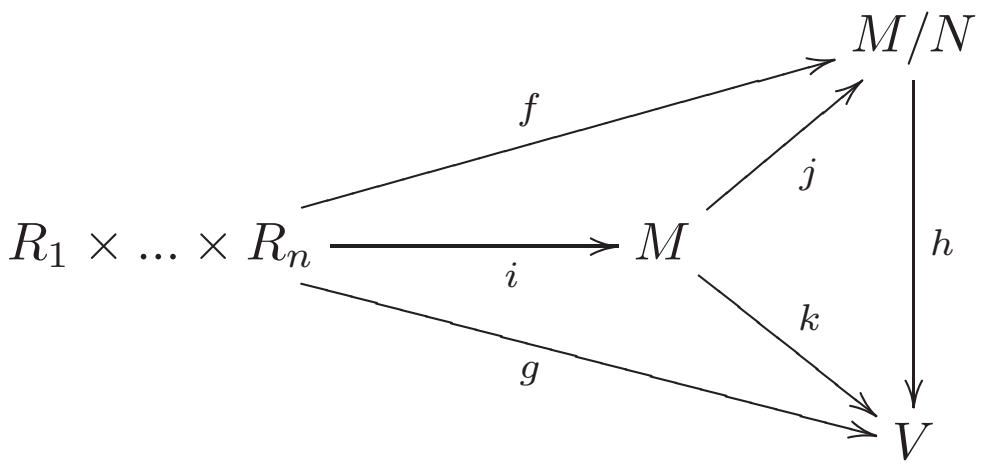

Так как $\mathrm{Im} f$ порождает $M/N$, то отображение $h$ однозначно определено. □

Согласно доказательству теоремы 13.4.1

$$V_1 \otimes ... \otimes V_n = M/N$$

Для $v_i \in V_i$ будем записывать

(13.4.9)
$$v_1 \otimes ... \otimes v_n = j(v_1, ..., v_n)$$

Из равенств (13.4.3) и (13.4.9) следует

$$f(v_1, ..., v_n) = v_1 \otimes ... \otimes v_n$$

Равенства (13.4.6) можно записать в виде

$$
\begin{aligned}
&v_1 \otimes ... \otimes (v_i + w_i) \otimes ... \otimes v_n \\
=&v_1 \otimes ... \otimes v_i \otimes ... \otimes v_n + v_1 \otimes ... \otimes w_i \otimes ... \otimes v_n \\
&v_1 \otimes ... \otimes (av_i) \otimes ... \otimes v_n = a(v_1 \otimes ... \otimes v_i \otimes ... \otimes v_n)
\end{aligned}
$$

Теорема 13.4.2. *Тензорное произведение $V_1 \otimes ... \otimes V_n$ является бимодулем над тензорным произведением $D_1 \otimes ... \otimes D_n$.*

Доказательство. Доказательство утверждения теоремы аналогично доказательству утверждения теоремы 13.3.2. □

Глава 14

# Список литературы

Глава 15

# Предметный указатель

Глава 16

# Специальные символы и обозначения